\documentclass[a4paper]{book}

\usepackage{amssymb,eucal,amsmath,amsthm,graphicx,pst-all,makeidx,fancyhdr,mathrsfs}
\usepackage[latin1]{inputenc}

\makeindex

\newcommand{\mt}[1]{\text{\rm #1}}    
\newcommand{\comment}[1]{}            
\newcommand{\Proof}{\begin{proof}[\textsc{Proof}] \;}
\newcommand{\set}[1]{\{#1\}}
\newcommand{\suchthat}{\;|\;}         
\newcommand{\restrict}{\,|}           
\newcommand{\compose}{\circ}          
\newcommand{\define}{\;{\rm :=}\;}
\newcommand{\invdef}{\;{\rm =:}\;}

\newcommand{\without}{\mathord{\setminus}}
\newcommand{\blank}{\text{\textvisiblespace}}

\newcommand{\id}{\mt{id}}             
\newcommand{\incl}{\mt{incl}}         
\newcommand{\pr}{\mt{pr}}             
\newcommand{\R}{\mathbb{R}}           
\newcommand{\N}{\mathbb{N}}           
\newcommand{\Z}{\mathbb{Z}}           
\newcommand{\leer}{\varnothing}       
\newcommand{\mfbd}{\partial}          
\newcommand{\CP}{\mathbb{CP}}         
\newcommand{\RP}{\mathbb{RP}}         
\newcommand{\eval}[2]{\langle #1,#2 \rangle}     
\newcommand{\abs}[1]{\left\lvert#1\right\rvert}  
\newcommand{\norm}[1]{\left\lVert#1\right\rVert} 
\newcommand{\diag}{\mt{diag}}         
\newcommand{\Gr}{G}                   
\newcommand{\eps}{\varepsilon}        
\newcommand{\ot}{\text{$\gets$}}      
\newcommand{\oointerval}[2]{(#1,#2)}  
\newcommand{\ocinterval}[2]{(#1,#2]}  
\newcommand{\cointerval}[2]{[#1,#2)}  
\newcommand{\Sob}[2]{H^{#2,#1}}       
\newcommand{\Sobzero}[2]{H^{#2,#1}_0} 

\newcommand{\FILL}{\mspace{30mu}}

\swapnumbers
\theoremstyle{definition}
\newtheorem{remark}{Remark}[section]
\newtheorem{remarks}[remark]{Remarks}

\newtheorem{remarksconventions}[remark]{Remarks and conventions}
\newtheorem{remarkdefinition}[remark]{Remark and definition}
\newtheorem{definition}[remark]{Definition}

\newtheorem{notation}[remark]{Notation}
\newtheorem{facts}[remark]{Facts}
\newtheorem{fact}[remark]{Fact}

\newtheorem{example}[remark]{Example}
\newtheorem{examples}[remark]{Examples}

\newtheorem{theplainproblem}[remark]{The plain problem}
\newtheorem{thedistributionproblem}[remark]{The distribution problem}
\newtheorem{theorthogonalproblem}[remark]{The orthogonal problem}
\newtheorem{thetimedistributionproblem}[remark]{The time distribution problem}
\newtheorem{thespacedistributionproblem}[remark]{The space distribution problem}
\newtheorem{thediffeotopyproblem}[remark]{The diffeotopy class problem}
\newtheorem{thetimediffeotopyproblem}[remark]{The time diffeotopy class problem}
\newtheorem{thespacediffeotopyproblem}[remark]{The space diffeotopy class problem}
\newtheorem{thehomotopyproblem}[remark]{The homotopy class problem}

\newtheorem{frameformulae}[remark]{ON frame formulae}

\newtheorem{finalremark}[remark]{Final remark}
\theoremstyle{plain}
\newtheorem{theorem}[remark]{Theorem}
\newtheorem{trichotomy}[remark]{The trichotomy theorem}
\newtheorem{lemma}[remark]{Lemma}
\newtheorem{corollary}[remark]{Corollary}
\newtheorem{proposition}[remark]{Proposition}
\newtheorem{conjecture}[remark]{Conjecture}
\newtheorem{escconjecture}[remark]{The esc Conjecture}

\DeclareMathOperator{\Torsion}{\mt{Torsion}}
\DeclareMathOperator{\im}{\mt{im}}

\newcommand{\End}{\mt{End}}
\DeclareMathOperator{\Mon}{\mt{Mon}} 

\newcommand{\Hmat}{\binom{0\;1}{1\;0}}

\DeclareMathOperator{\Sym}{\mt{Sym}} 

\DeclareMathOperator{\SO}{SO}
\DeclareMathOperator{\GL}{GL}
\DeclareMathOperator{\OO}{O}
\DeclareMathOperator{\U}{U}

\newcommand{\Diff}{\mt{Diff}}

\DeclareMathOperator{\vol}{vol}            
\DeclareMathOperator{\volume}{volume}      
\DeclareMathOperator{\scal}{scal}
\DeclareMathOperator{\Ric}{Ric}
\DeclareMathOperator{\grad}{grad}
\DeclareMathOperator{\Hess}{Hess}
\DeclareMathOperator{\divergence}{div}
\DeclareMathOperator{\laplace}{\Delta}
\DeclareMathOperator{\SecondFF}{\varPi}    
\DeclareMathOperator{\secondFF}{\mspace{2mu}\hat{\mspace{-2mu}\SecondFF}}
\DeclareMathOperator{\ind}{index}

\DeclareMathOperator{\rank}{rank}
\DeclareMathOperator{\spann}{span}

\DeclareMathOperator{\Lin}{Lin}

\newcommand{\Poincare}{Poincar\'{e}}
\newcommand{\Eliasson}{El{\'{\i}}asson}
\newcommand{\BerardBergery}{B\'{e}rard Bergery}

\newcommand{\LeviCivita}{Levi-Civita}
\newcommand{\Frechet}{Fr\'{e}chet}
\newcommand{\Gateaux}{G\^{a}teaux}
\newcommand{\Theoreme}{Th\'{e}or\`{e}me}

\renewcommand{\labelenumi}{(\roman{enumi})}
\renewcommand{\bar}{\overline}

\newcommand{\Qurv}{Q}               
\DeclareMathOperator{\qual}{qual}   
\newcommand{\fscal}{\mt{scal}^{\mt{fol}}}
\newcommand{\flaplace}{\laplace^{\mt{fol}}}

\newcommand{\lot}{\text{l.o.t.}}    

\newcommand{\conf}{\text{\tt conform}}
\newcommand{\stre}{\text{\tt stretch}}
\newcommand{\change}{\text{\tt change}}

\newcommand{\switch}{\text{\tt switch}}

\newcommand{\connsum}{\text{\rm\#}}
\newcommand{\Compl}{\text{\rm Compl}}
\newcommand{\Uni}{\mathscr{U}}     

\DeclareMathOperator{\Distr}{Distr}
\DeclareMathOperator{\Twist}{Twist}
\DeclareMathOperator{\Metr}{Metr}

\newcommand{\moebius}{\mathfrak{M}}     
\newcommand{\klein}{\mathfrak{K}}       

\newcommand{\Time}{\mathfrak{Time}}
\newcommand{\Space}{\mathfrak{Space}}
\newcommand{\Timifier}{\mathfrak{Timifier}}
\newcommand{\Spacifier}{\mathfrak{Spacifier}}

\newcommand{\tdc}{\mt{tdc}}
\newcommand{\sdc}{\mt{sdc}}
\newcommand{\TDC}{\mt{TDC}}
\newcommand{\SDC}{\mt{SDC}}
\newcommand{\tmc}{\mt{tmc}}
\newcommand{\smc}{\mt{smc}}
\newcommand{\TMC}{\mt{TMC}}
\newcommand{\SMC}{\mt{SMC}}
\newcommand{\cdc}{\mt{cdc}}
\newcommand{\CDC}{\mt{CDC}}

\newcommand{\comp}{\text{\rm comp}}

\newcommand{\PD}{\Upsilon}

\begin{document}

\pagenumbering{roman}


\begin{titlepage}

\begin{centering}
\vspace*{1.5cm}
{\Huge \bf Pseudo-Riemannian metrics}\\
\vspace*{5mm}
{\Huge \bf with prescribed scalar curvature}\\
\vspace*{2cm}
\Large Marc Nardmann\\
\end{centering} 
\vspace*{2cm}

\begin{center}
{\bf Abstract}
\end{center}

We consider the following generalisation of a well-known problem in Riemannian geometry: When is a smooth real-valued function $s$ on a given compact $n$-dimensional manifold $M$ (with or without boundary) the scalar curvature of some smooth pseudo-Riemannian metric of index $q\in\set{1,\dots,n-1}$ on $M$? We prove that this is the case for every $s$ if $3\leq q\leq n-3$, provided $M$ admits a metric of index $q$ at all. In fact, if $3\leq q\leq n-3$, then each connected component of the space of all pseudo-Riemannian metrics of index $q$ on $M$ contains a metric with scalar curvature $s$. We prove several theorems for pseudo-Riemannian metrics of index $1$ or $2$ as well. For instance, we show that on a compact orientable connected $4$-dimensional manifold $M$ with nonempty boundary, for every function $s\in C^\infty(M,\R)$ which is positive in at least one point and for every connected component $\mathscr{C}$ of the space of time-orientable Lorentzian (i.e.\ \mbox{$-$+++}) metrics on $M$, there exists one metric in $\mathscr{C}$ whose scalar curvature is $s$.
\medskip

The present work is essentially (i.e.\ up to one minor supplement) my PhD thesis, submitted to and accepted by the faculty of mathematics and computer science at the University of Leipzig (Germany).

\end{titlepage}


\chapter*{Preface\markboth{Preface}{}}
\addtocontents{toc}{\protect\contentsline {chapter}{Preface}{\roman{page}}} 

A classical problem in Riemannian geometry is the \emph{prescribed scalar curvature problem}: Given a smooth manifold $M$ and a smooth real-valued function $s$ on $M$, is there a Riemannian metric on $M$ whose scalar curvature is $s$?
\smallskip\\
The same question can be asked for pseudo-Riemannian metrics of any given index, for example for Lorentzian metrics. This pseudo-Riemannian generalisation of the prescribed scalar curvature problem is the topic of the present thesis.
\smallskip\\
The Riemannian case has been solved largely by J.\ L.\ Kazdan and F.\ W.\ Warner in 1974--75 (cf.\ Appendix \ref{D1} for a list of the main results). It seems that no work has been done on the pseudo-Riemannian analogue of the prescribed scalar curvature problem. The reason is apparently that a promising approach to a solution has been missing: The main tool in the Riemannian case is an elliptic partial differential equation whose obvious generalisation to the pseudo-Riemannian case is no longer elliptic and therefore virtually useless for solving the problem globally on a manifold with possibly complicated topology.
\medskip\\
The aim of the present thesis is to remedy this situation by the construction of a new PDE which generalises the PDE from the Riemannian case and is still \emph{elliptic}, in contrast to the naive generalisation mentioned above. Solutions of this new PDE yield solutions of the pseudo-Riemannian prescribed scalar curvature problem. The construction of the PDE is sketched and explained in Section \ref{ONETWO}, and carried out in detail in Chapter \ref{FOUR}. The second part of the thesis is concerned with solving the PDE.
\medskip\\
We prove the existence of solutions mainly via the method of sub- and supersolutions and via the method that Kazdan and Warner developed for the Riemannian prescribed scalar curvature problem; the latter employs the implicit function theorem for Banach spaces and a theorem which describes the $L^p$-closure of the orbit of a function under the right action by the diffeomorphism group. For the case of Lorentzian metrics on a $2$-manifold with nonempty boundary, we use direct methods in the calculus of variations.
\smallskip\\
These analytic techniques require that the manifold is \emph{compact} (with or without boundary). Although there exist methods to deal with noncompact manifolds, we will restrict ourselves to the compact case in the present work. In this way, we avoid complicated analytic estimates and can focus on the geometric and topological phenomena.
\smallskip\\
However, we clearly cannot do entirely without estimates. We need at least \emph{some} knowledge about the coefficient functions which appear in our PDE in order to prove that solutions exist. These coefficient functions are determined by ``background'' data which we have to put into the construction of our elliptic PDE: a Riemannian metric (analogous to the background metric which appears in the Riemannian prescribed scalar curvature problem, that is, in the PDE of the Yamabe problem) and a distribution (in the differential topological sense, i.e.\ a sub vector bundle of the tangent bundle). If we choose these background data carefully, we get coefficient functions with the desired properties and can solve our PDE.
\medskip\\
In particular, we have to establish the existence of distributions which satisfy a certain nonintegrability condition. (Contact structures in odd dimensions are examples which have this property; but in dimensions $>3$, the contact condition is stronger than the property we need.) This is a purely differential topological problem which can be solved by M.\ Gromov's h-principle techniques, in particular the convex integration method --- provided the manifold under consideration has dimension $\geq4$. In dimension $3$ (where the partial differential relation in question is not ample in the sense of Gromov), we can refer to well-known results about existence of contact structures, due to J.\ Martinet, R.\ Lutz, and Y.\ Eliashberg.
\smallskip\\
The prescribed scalar curvature problem turns out to be harder in the case of metrics with index $1$ (i.e.\ Lorentzian metrics) or $2$ than in the case of metrics with higher index. In order to deal with those higher-index metrics, we just have to choose appropriately nonintegrable background distributions; whereas for instance in the Lorentzian case, also the choice of the background Riemannian metric is crucial. (The reason for this difference is that a certain coefficient function in our PDE vanishes automatically in the Lorentzian case because line distributions are always integrable.)
\smallskip\\
Since the existence of good background \emph{metrics} is much harder to establish than the existence of good background \emph{distributions}, the main theorems in this thesis contain only those facts whose proofs avoid subtle constructions of suitable Riemannian metrics. However, I state a general conjecture which claims that such metrics do almost always exist, and I give a rough idea of how to prove that. This conjecture would solve the Lorentzian case of the prescribed scalar curvature problem more or less completely on manifolds of dimension $\geq4$.
\medskip\\
Here are the main results of this thesis (all manifolds and metrics are of class $C^\infty$).
{\it \begin{itemize}
\item
Let $M$ be a compact $n$-dimensional manifold (with or without boundary), let $q\in\set{3,\dots,n-3}$, and let $s\in C^\infty(M,\R)$. Then every connected component of the space of pseudo-Riemannian metrics of index $q$ on $M$ contains a metric whose scalar curvature is $s$. (In particular, if $M$ admits a pseudo-Riemannian metric of index $q$, then it does also admit such a metric with scalar curvature $s$.)
\item
Let $M$ be a connected compact $n$-manifold (with or without boundary) where $n\geq5$, let $q\in\set{1,2}$, and let $s\in C^\infty(M,\R)$ be somewhere positive (i.e.\ positive in at least one point). Then every connected component of the space of pseudo-Riemannian metrics of index $q$ on $M$ contains a metric with scalar curvature $s$.
\item
Let $M$ be a connected compact orientable $4$-manifold with nonempty boundary, let $s\in C^\infty(M,\R)$ be somewhere positive, and let $\mathscr{C}$ be a connected component of the space of Lorentzian (i.e.\ index-$1$) metrics on $M$; we assume that $\mathscr{C}$ contains a time-orientable metric\footnote{Then all metrics in $\mathscr{C}$ are time-orientable.}. Then $\mathscr{C}$ contains a metric with scalar curvature $s$.
\item
Let $M$ be a connected compact orientable $3$-manifold (with or without boundary), and let $s\in\! C^\infty(M,\R)$ be somewhere positive. Then every connected component of the space of Lorentzian metrics on $M$ contains a metric with scalar curvature $s$.
\item
Let $M$ be a compact connected $2$-manifold with nonempty boundary, and let $s\in C^\infty(M,\R)$. Then there is a Lorentzian metric on $M$ with scalar curvature $s$. 
\item
Let $M$ be either the $2$-dimensional torus or the Klein bottle (these are the only closed nonempty connected $2$-manifolds which admit a Lorentzian metric), and let $s\in C^\infty(M,\R)$. Then there is a Lorentzian metric on $M$ with scalar curvature $s$ if and only if $s$ is the constant $0$ or changes its sign (i.e.\ is positive somewhere and negative somewhere else).
\end{itemize} }
Let me refer you to Section \ref{ONETHREE} for a complete list of results.
\medskip\\
Obviously, the theorems above leave many questions open. Some of them are discussed and supplemented with conjectures and speculations in Chapters \ref{FIVE}--\ref{SEVEN}. Chapter \ref{ONE} contains an extensive introduction to the pseudo-Riemannian prescribed scalar curvature problem, raising additional, more refined questions which are not answered by the statements above. Certainly a lot of work remains to be done on the prescribed scalar curvature problem.
\medskip\\
It is not necessary for a reader of the present thesis to have background knowledge in all the fields which are touched upon; i.e.\ Riemannian and pseudo-Riemannian geometry, the theory of partial differential equations, differential and algebraic topology. Since this work does not build on too many previous results, it should be possible to explain it to a not too specialised audience, and I have tried to do so. Most of the prerequisites can be found in the appendices.

\subsection*{Acknowledgements}

This work was done while I was a fellow of the Graduiertenkolleg \emph{Analysis, Geometry and their Interaction with the Natural Sciences} in the faculty of mathematics and computer sciences at the University of Leipzig, supported by the Deutsche Forschungsgemeinschaft (DFG) and the state of Sachsen.
\medskip\\
Profs.\ Matthias Günther, Hans-Bert Rademacher, Matthias Schwarz, and Stephan Luckhaus gave me the opportunity to write my dissertation on this topic in the Graduiertenkolleg. I am grateful to Profs.\ Rademacher and Felix Finster for several discussions about my work.
\medskip\\
Countless remarks by my colleagues in the Graduiertenkolleg have influenced me during all stages of the project. In addition to encouragement and black humour, I received much relevant mathematical information, in particular from Michael Holicki, Frank Klinker, Matthias Kurzke, Mario Listing, Olaf Müller, and Kai Zehmisch. I have to emphasise the patience with which Matthias and Kai listened to and commented on my respective idea of the day.

\subsection*{How the arXiv version differs from my thesis}

A few typos have been corrected throughout the text, a few words changed for stylistic reasons. The condition in Theorems \ref{INTROmain3b} and \ref{INTROmain4b} is now stated in terms of the signature of the $4$-manifold; correspondingly, parts of Subsection \ref{threeplanefour} have been rewritten. Because the main body of the text has not changed, the present work is still called \emph{thesis} throughout.


\cleardoublepage
\addtocontents{toc}{\protect\contentsline {chapter}{Contents}{\roman{page}}}
\tableofcontents


\chapter*{Conventions and notations\markboth{Conventions and notations}{}}
\addtocontents{toc}{\protect\contentsline {chapter}{Conventions and notations}{\roman{page}}}

\begin{itemize}
\item
The word \emph{smooth} means $C^\infty$. This is the category we usually work in; i.e., all our manifolds, functions, vector fields etc.\ are assumed to be smooth, except explicitly stated otherwise.
\smallskip\\
Our definition of \emph{manifold} includes the Hausdorff property and paracompactness. (Every generalised --- that is, Hausdorff but not necessarily paracompact --- manifold which admits a semi-Riemannian metric is paracompact; cf.\ Corollary 2 in \cite{Marathe1972}.) Manifolds are assumed to be finite-dimensional (except when the word \emph{manifold} is preceded by \emph{\Frechet} or \emph{Banach} or \emph{Hilbert}) and pure-dimensional (that is, all connected components have the same dimension). An \emph{$n$-manifold} is an $n$-dimensional manifold.
\smallskip\\
\emph{Manifolds} may have a boundary. (Most of the results in this thesis hold for manifolds with nonempty boundary as well as for manifolds without boundary.) \emph{Closed} manifolds are compact manifolds without boundary. \emph{Open} manifolds are manifolds all of whose connected components are not closed. (When we talk about a subset of a manifold, we distinguish carefully whether it is open/closed as a topological \emph{subset} or open/closed as a \emph{manifold}.)
\item
The \emph{index} of a nondegenerate symmetric bilinear form on a finite-dimensional $\R$-vector space is the number of negative entries in any of its diagonalisations.
\smallskip\\
Let $M$ be an $n$-manifold. A semi-Riemannian metric on $M$ is a section in the bundle of nondegenerate symmetric bilinear forms on the tangent bundle $TM$. A semi-Riemannian metric $g$ on $M$ \emph{has index $q$} if and only if the bilinear form $g_x$ has index $q$ for every $x\in M$. A \emph{pseudo-Riemannian metric} on $M$ is a semi-Riemannian metric on $M$ which has index $\in\set{1,\dots,n-1}$. A \emph{Riemannian metric} on $M$ is a semi-Riemannian metric on $M$ which has index $0$. A \emph{Lorentzian metric} on $M$ is a semi-Riemannian metric on $M$ which has index $1$.
\smallskip\\
According to a different convention, a Lorentzian metric on $M$ is a semi-Riemannian metric on $M$ which has index $n-1$. The choice of convention is (at least slightly) relevant to the prescribed scalar curvature problem, because the scalar curvature of a semi-Riemannian metric $g$ (of index $q$) is related to the scalar curvature of the metric $-g$ (of index $n-q$) by the formula $\scal_{-g} = -\scal_g$ (cf.\ \ref{scalsigninversion}).
\item
The word \emph{distribution} will only be used in the sense of differential topology (cf.\ Subsection \ref{distributionproblem}): a \emph{$q$-plane distribution} on a manifold $M$ is a smooth rank-$q$ sub vector bundle of the tangent bundle $TM$.
\item
Take note of our usage of the word \emph{conformal}; cf.\ Definition \ref{confdef}.
\item
We denote the space of $C^r$ sections in a (smooth) fibre bundle $E\to M$ by $C^r(M\ot E)$; cf.\ also Notation \ref{Crnotation}. (I don't know who invented this notation; I learned it from \cite{KrieglMichor}.)
\end{itemize}

\begin{tabular*}{163mm}{@{}p{43mm}@{}p{120mm}@{}}
$\N$ & $=\set{0,1,2,\dots}$\\
$\oointerval{a}{b}$, \;\;$\ocinterval{a}{b}$, \;\;$\cointerval{a}{b}$ & open resp.\ half-open interval\\
$\N_{\geq2}$, \;$\R_{>0}$, \;etc. & $\set{n\in\N \suchthat n\geq2}$, \;$\set{r\in\R \suchthat r>0}$, \;etc.\\
$R_g$, $\Ric_g$, $\scal_g$ & Riemann (cf.\ \ref{Riemanntensor}) resp.\ Ricci resp.\ scalar curvature of the metric $g$\\
$\delta_{ij}$ & Kronecker symbol: $\delta_{ij}=1$ if $i=j$, and $\delta_{ij}=0$ if $i\neq j$\\
$\id_M$ & identity map on the set $M$\\
$f\restrict M$ & restriction of the function $f$ to the subset $M$ of its domain\\
$df(X)\equiv \partial_Xf$ & derivative of the function $f$ in the direction of the vector field $X$\\
$Df \equiv Tf$ & derivative of a map $f$ between (Banach) manifolds\\
$D_xf\equiv T_xf$ & value of the derivative $Df\equiv Tf$ in the point $x$
\end{tabular*}

\begin{tabular*}{163mm}{@{}p{43mm}@{}p{120mm}@{}}
$[v,w]$ & Lie bracket of vector fields\\
$S^n$ & $n$-dimensional sphere\\
$T^n$ & $n$-dimensional torus $S^1\times\ldots\times S^1$ ($n$ factors)\\
$\RP^n$, $\CP^n$ & real resp.\ complex $n$-dimensional projective space\\
$\moebius$, $\klein$ & Möbius strip, Klein bottle\\
$M\connsum N$ & connected sum of manifolds\\
$\beta_M$, $\sigma_M$ & intersection form of the $4$-manifold $M$, resp.\ signature of $M$\\
$\chi(M) \equiv \chi_M$ & Euler characteristic of $M$\vspace{0.012cm}\\
$\Sob{p}{k}(M,\R)$, $\Sobzero{p}{k}(M,\R)$ & Sobolev spaces (cf.\ Appendix \ref{sobolevappendix})\\
$H^k(M;G)$, $H_k(M;G)$ & $k$th cohomology resp.\ homology group of the space $M$ with coefficients in $G$\\
$\pi_0(X)$ & set of path-connected components of the topological space $X$\\
$\pi_k(X;x_0)$, $\pi_k(X)$ & $k$th homotopy group of the space $X$ (with base point $x_0$)\vspace{0.1cm}\\
$C^r(M,N)$ & space of $r$-times continuously differentiable maps from $M$ to $N$\\
$C^r(M\ot E)$ & space of $C^r$ sections in the fibre bundle $E\to M$\\
$J^r(M,N)$, $J^r(M\ot E)$ & jet spaces (cf.\ Appendix \ref{jets})\\
$j^rf$, $j^r_xf$ & $r$-jet prolongation of the map/section $f$, resp.\ its evaluation in $x$ (cf.\ \ref{jets})\\
$J^1_{\bot W}E$ & perp-jet space of the bundle $E$ (cf.\ Appendix \ref{hprincipleONE})\vspace{0.13cm}\\
$\Diff(M)$ & group of all diffeomorphisms from $M$ to $M$\\
$\Diff^0(M)$ & path component of the identity in $\Diff(M)$\\
$\ind(g)$ & index of a nondeg. symm. bilinear form resp.\ semi-Riemannian metric\\
$\Metr_q(M)$ & space of semi-Riemannian metrics of index $q$ on $M$; cf.\ \ref{metrdef}\\
$\Distr_q(M)$ & space of $q$-plane distributions on $M$; cf.\ \ref{distrdef}\vspace{0.015cm}\\
$\Lin(V,W)$ & Banach space of continuous linear maps between the real Banach spaces $V,W$, resp.\ vector bundle of linear maps between the real vector bundles $V,W$\\
$\Compl(V)$ & affine space resp.\ affine bundle of complementary subspaces (cf.\ Appendix \ref{Grassmannappendix})\\
$\Gr_q(E)$ & $q$th Grassmann manifold resp.\ Grassmann bundle (cf.\ Appendix \ref{Grassmannappendix})\\
$\Sym(E)$ & vector space/bundle of symmetric bilinear forms on the v.\ space/bundle $E$\\
$\Sym_q(E)$ & subset of $\Sym(E)$ consisting of nondegenerate bilinear forms of index $q$\\
$\Time$, $\Space$ & cf.\ \ref{TimeSpacedef}\\
$\Timifier$, $\Spacifier$ & cf.\ \ref{TimifierSpacifierdef}\\
$\tdc$, $\sdc$, $\TDC$, $\SDC$ & cf.\ \ref{distrcompdef}, \ref{TDCSDC}\\
$\tmc$, $\smc$, $\TMC$, $\SMC$ & cf.\ \ref{metrcompdef}, \ref{TMCSMC}\\
$\cdc$, $\CDC$ & cf.\ \ref{cdcdef}, \ref{CDCdef}\smallskip\\
$E^\ast$ & dual vector space resp.\ dual vector bundle\\
$\bot_gU$ & $g$-orthonormal bundle: cf.\ \ref{maximally}, \ref{musical}\\
$\bot U$ & normal bundle: cf.\ \ref{normalbundle}\\
$\flat_g$, $\sharp_g$ & musical isomorphisms: cf.\ \ref{musical}\\
$i_U$ & inclusion map: cf.\ \ref{musical}\\
$\pr^U_g$ & orthogonal projection: cf.\ \ref{projdef}\\
$\eval{\alpha}{\beta}_g$ & scalar product of cotangent vectors resp.\ $1$-forms: cf.\ \ref{musical}\\
$\eval{\alpha}{\beta}_{g,U}$ & cf.\ \ref{scalproddef}\\
$\abs{\alpha}_g$ & $\equiv \sqrt{\eval{\alpha}{\alpha}_g}$ (only used if $g$ is positive definite); similarly $\abs{\alpha}_{g,U} \equiv \sqrt{\eval{\alpha}{\alpha}_{g,U}}$\\
$\nabla$, \;$\nabla^{(g)}$ & covariant derivative with respect to the \LeviCivita\ connection of $g$\\
$i:U$, $i,j:U$, $i:\bot U$ & notation for $U$-adapted ON frames: cf.\ \ref{colonnotation}\\
$\eps_i$ & notation for ON frames: cf.\ \ref{ONframedef}\\
$\Gamma^k_{ij}$ & ON frame Christoffel symbol: cf.\ \ref{ONChristoffel}\\
$\varGamma^k_{ij}$ & coordinate system Christoffel symbol (will rarely occur; no danger of confusion)\vspace{0.12cm}\\
$\switch(g,V)$ & cf.\ \ref{switchdefinition}\\
$\stre(g,f,V)$ & cf.\ \ref{stretchdefinition}\\
$\conf(g,f)$ & cf.\ \ref{confdefinition}\\
$\change(g,f,\kappa,V)$ & cf.\ \ref{changedefinition}\vspace{0.13cm}\\
$\divergence_g(X)$ & divergence of the vector field $X$ with respect to the semi-Riemannian metric $g$\\
$\divergence^U_g(X)$ & cf.\ \ref{divUXdef}
\end{tabular*}

\begin{tabular*}{163mm}{@{}p{43mm}@{}p{120mm}@{}}
$\divergence^U_g$ & cf.\ \ref{scalprodexamples}\\
$\eval{\divergence^U_g}{df}_{g,\bot U}$ & cf.\ \ref{scalprodexamples}\\
$\eval{\divergence^U_g}{\divergence^U_g}_{g,\bot U}$ & cf.\ \ref{scalprodexamples}\\
$\laplace_g(f)$ & Laplacian of the function $f$ with respect to the metric $g$\\
$\laplace^V_{g,W}$, \;$\laplace^V_g$ & cf.\ \ref{laplacedef}\\
$\sigma_{g,V}$, $\tau_{g,V}$ & cf.\ \ref{defsigmatau}\\
$\scal^{V,W}_g$ & cf.\ \ref{scaldef}\\
$\Qurv^V_g$ & cf.\ \ref{Qurvdef}\\
$\qual^V_g$ & cf.\ \ref{qualdef}\\
$\xi_{g,V}$ & cf.\ \ref{xidef}\\
$\chi_{g,V}$ & cf.\ \ref{chidef}\\
$\nabla^{(g)}_VV$ & cf.\ \ref{lineobjects}\\
$\divergence_g(V)^2$, $\partial_V\divergence_g(V)$ & cf.\ \ref{lineobjects}\\
$\flaplace_{g,H}$ & cf.\ \ref{flaplacedef}\\
$\fscal_{g,H}$ & scalar curvature of a foliation: cf.\ \ref{fscal}\\
$\PD_{g,f,V}$ & our favourite elliptic operator: cf.\ \ref{PDdef}\smallskip\\
$A^\top$ & transpose of the matrix $A$\\
$\int_{(M,g)}f$ & the integral of a measurable function $f\colon M\to \R$ with respect to the measure induced by the semi-Riemannian metric $g$ on $M$
\end{tabular*}

\bigskip
The total contraction of a tensor field $T\otimes T$ (i.e.\ contraction of $T\otimes T$ in all corresponding indices) with respect to a metric $g$ is denoted by $\abs{T}_g^2$ --- but only if $g$ is Riemannian, not in the pseudo-Riemannian case. (As a matter of taste, I avoid a notation of the form $\abs{\ldots}^2$ in a situation where it might represent a negative number or function and where thus the notation $\abs{\ldots}$ has no unambiguous meaning.)

\medskip
The notion \emph{$g$-good} is defined in \ref{ggood}.


\chapter{Introduction} \label{ONE}

\pagenumbering{arabic}
A classical problem in Riemannian geometry is to determine which real-valued functions on a given manifold $M$ can be represented as scalar curvatures of Riemannian metrics on $M$. The results of J.\ L.\ Kazdan and F.\ W.\ Warner from the mid-1970s solve this problem largely; cf.\ Appendix \ref{D1} for a review of some facts which might help to put the present work into context.
\smallskip\\
This thesis addresses the analogous problem for pseudo-Riemannian metrics --- with an emphasis on global solutions on manifolds with possibly complicated topology. Our aim is \emph{not} to give the sharpest possible results for the problem on open subsets of $\R^n$, or on open neighbourhoods of an initial data hypersurface in a Lorentzian manifold. Rather, we would like to obtain any results at all for metrics on manifolds which are not diffeomorphic to a product manifold.
\smallskip\\
To my knowledge, there has been no prior work in this direction, not even in the case of Lorentzian metrics. However, even the much harder problem of prescribed \emph{Ricci} curvature for Lorentzian metrics has been investigated as an initial value problem on small open neighbourhoods of hypersurfaces; cf.\ D.\ DeTurck's Theorem 3.6 in \cite{DeTurck1983}.
\smallskip\\
But what exactly is the pseudo-Riemannian prescribed scalar curvature problem on arbitrary manifolds? Actually, there are \emph{several} semi-Riemannian generalisations of the Riemannian problem, which answer different natural questions in the pseudo-Riemannian case. We start this chapter with a detailed exposition of these different versions, in Section \ref{ONEONE}.
\smallskip\\
Section \ref{ONETWO} explains how we are going to solve the pseudo-Riemannian prescribed scalar curvature problem, Section \ref{ONETHREE} contains a complete list of results, and Section \ref{ONEFOUR} gives a short overview of the structure of the thesis.
\medskip\\
Why should we care about the scalar curvature of pseudo-Riemannian metrics at all? Let me conclude these introductory remarks with a short motivation (albeit limited to the Lorentzian case).

\subsubsection*{Connections to general relativity}

The physical theory of general relativity models spacetime as a Lorentzian manifold of dimension $4$. On every open subset of spacetime which contains no matter (i.e.\ a \emph{vacuum region}), A.\ Einstein's equation
\begin{equation} \label{EINSTEIN}
\Ric_g -\frac{1}{2}\scal_gg +\Lambda g = 0 \end{equation}
holds; here $g$ is the Lorentzian spacetime metric, and $\Lambda\in\R$ is the \emph{cosmological constant}. (We could also discuss non-vacuum regions, but let's keep things simple.)
\smallskip\\
A natural question arises: Given $\Lambda\in\R$ and a manifold $M$ (of dimension $4$, say), is there a Lorentzian metric $g$ on $M$ which satisfies \eqref{EINSTEIN}? In other words, can $M$ be a vacuum spacetime for the given cosmological constant $\Lambda$? (By rescaling $g$ with a constant, one sees easily that the answer depends only on the sign --- positive, zero, or negative --- of $\Lambda$.)
\smallskip\\
This is certainly a very difficult question in general. Even the famous analogous question in Riemannian geometry --- which has been studied extensively --- is still far from being solved (cf.\ \cite{LeBrunMcKenzie} for several overview articles on the subject). So one should probably try to solve a simpler problem first.
\smallskip\\
Recall in this context that $\Ric_g -\frac{1}{2}\scal_gg = 0$ holds for every $2$-dimensional semi-Riemannian manifold $(M,g)$. For every semi-Riemannian manifold of dimension $n>2$, the vacuum Einstein equation $\Ric_g -\frac{1}{2}\scal_gg +\Lambda g = 0$ implies that $g$ is an Einstein metric: contracting this tensor equation yields $\scal_g = \frac{2n}{n-2}\Lambda$, hence $\Ric_g = \frac{2}{n-2}\Lambda g$. Conversely, if the metric $g$ is Einstein, i.e.\ $\Ric_g = \lambda g$ for some $\lambda\in\R$, then the vacuum Einstein equation holds with cosmological constant $\Lambda = \frac{n-2}{2}\lambda$. In particular, the constant scalar curvature of a vacuum Einstein metric in dimension $n>2$ has the same sign (positive, zero, or negative) as the cosmological constant.
\smallskip\\
This suggests that one should simplify our question from above as follows: Given $\Lambda\in\R$ (without loss of generality $\Lambda\in\set{1,0,-1}$) and a manifold $M$ of dimension $>2$, is there a Lorentzian metric $g$ on $M$ with constant scalar curvature $\Lambda$? If the answer is \emph{no} for some pair $(\Lambda,M)$, then the answer to the question from above is negative for this pair as well.
\smallskip\\
The first obstruction to a Lorentzian metric on $M$ with constant scalar curvature $\Lambda$ is of course non-existence of Lorentzian metrics on $M$. The results of this thesis suggest that there are no further obstructions if $n\geq4$. If this conjecture is true, then there are no nontrivial obstructions \emph{on the scalar curvature level} to the existence of solutions of \eqref{EINSTEIN}, in sharp contrast to the situation in the Riemannian case, where obstructions to positive or zero scalar curvature metrics exist and yield obstructions to Einstein metrics with these signs.
\medskip\\
We follow here a philosophy which says that it is interesting to understand the global geometry of Lorentzian manifolds with complicated topology. In Lorentzian geometry and general relativity, much is known about \emph{causal structure} and \emph{geodesic incompleteness} (cf.\ e.g.\ the singularity theorems of R.\ Penrose, S.\ Hawking, and others), or about the curvature of \emph{globally hyperbolic} spacetimes. General information on the \emph{curvature} of manifolds which are not globally hyperbolic seems to be rare. The results of the present thesis could be seen as a small step toward a better understanding of the curvature of such manifolds.\footnote{Some physicists might object to this motivation of the Lorentzian prescribed scalar curvature problem because they think that spacetimes which contain closed timelike curves are forbidden for physical-philosophical reasons (cf.\ e.g.\ \cite{HawkingEllis}, p.~189). Demanding the non-existence of closed timelike curves forces the topology of spacetime to be quite trivial: all closed manifolds are excluded (cf.\ e.g.\ \cite{HawkingEllis}, Proposition 6.4.2), and compact Lorentz cobordisms (cf.\ \cite{Yodzis1}, \cite{Yodzis2}) have to be topological products $[0,1]\times N$ then (cf.\ Theorem 2 in \cite{Geroch1967}). Thus the Lorentzian metrics which we are going to construct in the present thesis will usually admit closed timelike curves. Depending on the philosophical viewpoint, people might have different opinions as to whether this fact reduces the motivation one can draw from physics for the prescribed scalar curvature problem. Personally, I think there is no experimental evidence against closed timelike curves in the universe we live in, so I see no problem.} \footnote{There are several related questions which I cannot discuss here in detail, for example the problem of \emph{topology change} in general relativity; cf.\ e.g.\ \S3.3 in \cite{Tod1999}. In particular, to what extent remains Tipler's Theorem 4 in \cite{Tipler1977} true if one omits the so-called \emph{generic condition}? When is topology change by compact constant scalar curvature spacetimes possible (being a necessary condition for topology change by compact vacuum spacetimes)? Cf.\ Remark \ref{boundaryvalue} for a first result in this direction.}

\section{Statement of the problem} \label{ONEONE}

We state the three most important versions of the prescribed scalar curvature problem: the plain problem \ref{plainproblem}, the distribution problem \ref{distributionproblem}, and the homotopy class problem \ref{homotopyclassproblem}. On a quick reading, it suffices to read just the statements of these three versions because everything else in this section is explanations and background.

\subsection{The plain version} \label{ONEONEONE}

In its simplest version, the problem we want to solve is the following:
\begin{theplainproblem} \label{plainproblem}
Let $n$ be an integer $\geq2$, and let $q\in\set{0,\dots,n}$. We are given an $n$-dimensional smooth manifold $M$ which admits a semi-Riemannian metric of index $q$; and, in addition, a smooth real-valued function $s$ on $M$. Is there a smooth semi-Riemannian metric of index $q$ on $M$ whose scalar curvature is $s$?
\end{theplainproblem}

\begin{remarksconventions}{\ } \label{conventions}
\begin{enumerate}
\item
If $M$ is a manifold of dimension $n\leq1$, only the constant function $0$ can occur as the scalar curvature of some semi-Riemannian metric on $M$. That's why we restrict our considerations to manifolds of dimension $\geq2$. Of course we could also restrict ourselves to nonempty connected manifolds, without loss of generality.
\item
An $n$-dimensional manifold admits a semi-Riemannian metric of index $q$ if and only if it admits a semi-Riemannian metric of index $n-q$ (because of the correspondence $g \leftrightarrow -g$). Every manifold admits a Riemannian metric. A connected manifold admits a Lorentzian metric if and only if it is either closed with zero Euler characteristic, or open; cf.\ Proposition 5.37 in \cite{ONeill}.\footnote{Note that the reference cited in the proof of this proposition deals only with one direction of the claimed equivalence. For orientable manifolds, this equivalence is proved in \S12 of \cite{MilnorStasheff}. I do not know a good reference for the nonorientable case.} The existence of ultra-Lorentzian metrics\footnote{i.e., pseudo-Riemannian metrics of index $\in\set{2,\dots,n-2}$ on an $n$-dimensional manifold} is harder to characterise; cf.\ \cite{HirzebruchHopf}, \cite{Thomas1967a}, \cite{Thomas1967b} for some results\footnote{These articles are concerned with distributions (i.e.\ sub vector bundles of the tangent bundle) instead of pseudo-Riemannian metrics. Note that a manifold admits a $q$-plane distribution if and only if it admits a semi-Riemannian metric of index $q$, by \ref{propeleven} and \ref{baumtheorem} below.}. The topological assumption \emph{which admits a semi-Riemannian metric of index $q$} in our statement of the plain problem is just the trivial necessary condition for the existence of a solution.
\item
We have assumed that the manifold $M$ and the function $s$ are smooth, and asked for a smooth semi-Riemannian metric. This is quite natural from a geometer's viewpoint. But we could consider lower regularity as well, of course; for example, we could ask: \emph{Let $r\geq0$, assume that $M$ is (at least\footnote{Since, for $r\in\N_{\geq1}$, every $C^r$ manifold admits a compatible $C^\infty$ structure (cf.\ \cite{Hirsch}, Theorem 2.2.9), and since this $C^\infty$ structure is unique up to $C^\infty$ diffeomorphism, we can assume without loss of generality that $M$ is smooth. In fact, we can even assume without loss of generality that $M$ is real-analytic; cf.\ Lemma \ref{realanalyticatlas} in Appendix \ref{diffbackground}.}) $C^{r+2}$, and that $s$ is $C^r$; is there a $C^{r+2}$ semi-Riemannian metric with scalar curvature $s$?} Alternatively (and, from a PDE theorist's viewpoint, more conveniently), we might consider $C^{r,\alpha}$ (Hölder) regularity or $\Sob{p}{r}$ (Sobolev) regularity (where $\alpha\in\oointerval{0}{1}$ and $p\in\cointerval{1}{\infty}$\;). We derive results along these lines naturally in this thesis by applying standard PDE techniques, but I have made no effort to state them in maximal generality and sharpness. The aim is to find smooth solutions.
\smallskip\\
However, we will discuss what happens if $M$ is equipped with a real-analytic structure (this we can assume without loss of generality) and $s$ is real-analytic. We will see that whenever we can prove that a smooth solution metric $g$ exists, we can even choose $g$ real-analytic.
\item
Without loss of generality, we could restrict ourselves to the consideration of the case $0\leq q\leq \frac{n}{2}$, by taking $\scal_{-g} = -\scal_g$ into account. Note that the problem has an obvious symmetry in the case $q=\frac{n}{2}$: If the function $s\in C^\infty(M,\R)$ can be realised as the scalar curvature of some index-$q$ metric, then the function $-s$ can also be realised as the scalar curvature of some index-$q$ metric.
\item
The main theorems in this thesis are stated for \emph{compact} manifolds. According to our convention, a \emph{manifold} may have a nonempty boundary. We allow this possibility because the pseudo-Riemannian prescribed scalar curvature problem for compact manifolds with nonempty boundary can in general not be reduced to the case of closed manifolds, in contrast to the Riemannian problem.
\smallskip\\
Even if one deals only with compact manifolds, there is no need to mention manifolds with nonempty boundary at all in the \emph{Riemannian} prescribed scalar curvature problem: Every compact manifold with nonempty boundary $M$ can be imbedded as an open subset into a closed manifold $\bar{M}$ (the so-called \emph{double of $M$}, for instance); and every smooth prescribed function $s$ on $M$ can be extended to a smooth function on $\bar{M}$.\footnote{This extendability is is tautological if we define smoothness accordingly, and it is nontrivial if we define smoothness by just demanding that all derivatives on the manifold interior can be continuously extended to the boundary. However, R.\ T.\ Seeley's extension theorem \cite{Seeley1964} for functions defined on a half-space in $\R^n$ (or a related result by H.\ Whitney) implies that extendability holds even in this case. (Seeley's result for $\R$-valued functions generalises to sections in a fibre bundle over an arbitrary manifold, by standard arguments in differential topology.)} Thus the results for closed manifolds solve the problem for manifolds with nonempty boundary as well: One extends $s$ to $\bar{M}$ in a suitable way\footnote{A \emph{suitable} extension would be one that is negative in some point of $\bar{M}\without M$.}, solves the prescribed scalar curvature problem for the extended function on the closed manifold $\bar{M}$, and restricts the solution metric to $M$.
\smallskip\\
This does not always work in the pseudo-Riemannian version. Consider the case where $M$ is a closed $2$-manifold of sufficiently high genus (a connected sum $T^2\connsum T^2$ would do) with an open ball removed. Then $M$ admits a Lorentzian metric, but it cannot be imbedded into any closed $2$-manifold which admits a Lorentzian metric (the only closed $2$-manifolds which admit a Lorentzian metric are the torus and the Klein bottle). It would therefore cause a loss of generality if we restricted ourselves to closed manifolds, while there would be no technical advantage in doing so; all we need is compactness.
\end{enumerate}
\end{remarksconventions}

\begin{remark}[completeness] \label{completenessremark}
Let us consider the Riemannian special case of the plain problem, i.e.\ the case $q=0$. If $M$ is noncompact, a natural sharper version of the plain problem is to look for \emph{complete} solution metrics. This makes the problem much harder (cf.\ Subsection \ref{D14}).
\smallskip\\
Things are even worse in the pseudo-Riemannian case since there is no pseudo-Riemannian analogue of the Hopf/Rinow theorem, which gives several equivalent descriptions of completeness of Riemannian metrics; cf.\ Chapters 1 and 6 in \cite{BeemEhrlichEasley}. One useful notion of completeness in the general case is geodesic completeness (which can be divided into timelike, spacelike and lightlike geodesic completeness). But recall that even \emph{closed} Lorentzian manifolds might be timelike, spacelike and lightlike geodesically incomplete (cf.\ \cite{Wald}, p.~242, or \cite{ONeill}, Example 7.16 and Exercise 9.12, for an example metric on the $2$-dimensional torus).
\smallskip\\
It is not clear under which conditions the metrics we are going to construct will be geodesically complete. We shall not deal with this question in the present thesis.
\end{remark}

\subsection{The distribution version(s)} \label{ONEONETWO}

In contrast to Riemannian geometry, Lorentzian geometry (or, more generally, pseudo-Riemannian geometry) is concerned with the \emph{causal properties} of metrics, i.e.\ with the question which directions in each tangent space are timelike resp.\ lightlike resp.\ spacelike. Our plain problem does not reflect this; if we solve it, we have a priori no information about the causal behaviour of the solutions. We will now state several versions of the problem which remove this failure by prescribing the causal character of the solution metrics. This subsection explains the sharpest versions one can hope to solve in useful generality.
\medskip\\
Let us recall some notions from differential topology. A \emph{distribution of rank $q$} (synonymously, a \emph{$q$-plane distribution}) on a manifold $M$ is a sub vector bundle of $TM$ which has rank $q$. A \emph{line distribution} is just a $1$-plane distribution. Two distributions $V$ and $H$ on $M$ are \emph{complementary} if and only if $TM$ is the internal direct sum of $V$ and $H$. (In this situation, we have obviously $\rank(V)+\rank(H)=\dim(M)$.)
\smallskip\\
Let $W$ be a finite-dimensional real vector space equipped with a symmetric nondegenerate bilinear form $g$. Recall that a vector $v\in W$ is called \emph{timelike} [resp.\ \emph{spacelike} resp.\ \emph{lightlike}] if and only if $v\neq0$ and $g(v,v)<0$ [resp.\ $g(v,v)>0$ resp.\ $g(v,v)=0$]. A sub vector space $V$ of $W$ is \emph{timelike} resp. \emph{spacelike} resp.\ \emph{lightlike} if and only if every non-zero vector in $V$ has the respective property.
\smallskip\\
Let $(M,g)$ be a semi-Riemannian manifold. A distribution is \emph{timelike} resp. \emph{spacelike} resp.\ \emph{lightlike} if and only if every non-zero vector in the distribution has the respective property.

\begin{thedistributionproblem} \label{distributionproblem}
Let $n$ be an integer $\geq2$, and let $q\in\set{0,\dots,n}$. We are given \begin{itemize}
\item
an $n$-dimensional smooth manifold $M$;
\item
a smooth real-valued function $s$ on $M$;
\item
two complementary distributions\footnote{The letters $V$ and $H$ remind us that in (special or general) relativistic diagrams, timelike directions correspond usually to the \emph{v}ertical axis, spacelike directions to the \emph{h}orizontal axis.} $V$ and $H$ on $M$, where $V$ has rank $q$.
\end{itemize}
Is there a smooth semi-Riemannian metric $g$ of index $q$ on $M$ whose scalar curvature is $s$, such that the distribution $V$ is timelike and the distribution $H$ is spacelike with respect to $g$?
\end{thedistributionproblem}

\begin{remarks}{\ } \label{distributionremarks}
\begin{enumerate}
\item
The assumption \emph{which admits a semi-Riemannian metric of index $q$} from the plain problem is superfluous in the distribution problem; cf.\ Subsection \ref{ONEONETHREE} below.
\item
We could have made the following elementary remark already when we were discussing the plain version of the prescribed scalar curvature problem; it applies, mutatis mutandis, to all versions.
\smallskip\\
If the distribution problem can be solved for the function $s\in C^\infty(M,\R)$ with respect to given distributions $V$ and $H$, then it can be solved, with respect to the same distributions $V$ and $H$, for every function of the form $cs$, where $c$ is a positive constant. Namely, if $s=\scal_g$ for some metric $g$ which makes $V$ timelike and $H$ spacelike, then $cs=\scal_{\tilde{g}}$ for the metric $\tilde{g}=c^{-1}g$ which makes $V$ timelike and $H$ spacelike.
\item
If we want to prescribe the causal behaviour of the desired metric, why not even prescribe the lightcones? I.e., why don't we fix a set (which satisfies some necessary properties) in each tangent space and demand that this set becomes the lightcone of our solution? Or, only slightly more drastic, why don't we prescribe for each nonzero vector in the tangent bundle whether it should become timelike, spacelike, or lightlike? The reason is as follows.
\smallskip\\
When two pseudo-Riemannian metrics $g_0,g_1$ on a manifold $M$ induce the same causal structure on $TM$ (i.e.\ each tangent vector is timelike/spacelike/lightlike with respect to $g_0$ if and only if it is timelike/spacelike/lightlike with respect to $g_1$), then $g_0$ and $g_1$ are conformal, i.e., there is a function $\lambda\in C^\infty(M,\R_{>0})$ such that $g_1 = \lambda g_0$. When two pseudo-Riemannian metrics $g_0,g_1$ on $M$ have the same index $q$ and the same lightcones (i.e.\ each tangent vector is $g_0$-lightlike if and only if it is $g_1$-lightlike), then $g_0$ and $g_1$ are conformal as well, provided $q\neq\frac{n}{2}$; in the case $q=\frac{n}{2}$, there is at least a function $\lambda\in C^\infty(M,\R_{\neq0})$ such that $g_1 = \lambda g_0$ (in other words: on each connected component of $M$, the metric $g_1$ is either conformal to $g_0$ or conformal to $-g_0$). These standard facts can be found in \cite{BeemEhrlichEasley}, Section 2.3 (in particular Lemma 2.1), for instance.
\smallskip\\
The scalar curvature problem with prescribed lightcones is therefore essentially the \emph{conformal problem}: \emph{Given a semi-Riemannian metric $g$ on $M$ and a function $s\in C^\infty(M,\R)$, can we find a function $\lambda\in C^\infty(M,\R_{>0})$ such that the metric $\lambda g$ has scalar curvature $s$?} This is a rigid problem: we have to solve a partial differential equation (second-order, nonlinear but quasilinear) for the function $\lambda$. (We will write down the equation explicitly later on in this chapter.) In contrast, the plain and distribution problems are far from being rigid: a naive count of scalar functions yields $n(n+1)/2$ free parameters (for a metric on an $n$-manifold) but just one equation, $\scal_g=s$.
\smallskip\\
The conformal problem has been investigated intensively in the Riemannian case; cf.\ Subsection \ref{D1} below. It is one of the approaches Kazdan and Warner used in their work on the Riemannian prescribed scalar curvature problem. The main reason why it was possible to get useful results for the Riemannian conformal problem is that the PDE one has to solve is \emph{elliptic} in that case. The PDE corresponding to the pseudo-Riemannian conformal problem is not elliptic, however; it is hyperbolic in the Lorentzian case, and ultrahyperbolic in the ultra-Lorentzian case.
\smallskip\\
This non-ellipticity is the reason why we don't deal with the conformal problem: No PDE technique that I am aware of could produce global solutions of this nonlinear hyperbolic (or ultrahyperbolic) equation on, e.g., a closed manifold. Nonlinear hyperbolic equations are already quite hard to solve on domains in $\R^n$, but closed timelike curves produce additional difficulties. Lorentzian metrics on manifolds like $S^3$ or $T^4\connsum T^4\connsum \CP^2\connsum \CP^2$ do always admit closed timelike curves, and the patterns in which the timelike curves close might be quite complicated. There seems to be no chance to find useful criteria for the solvability of such equations.
\end{enumerate}
\end{remarks}

We can demand even more from a solution metric than we did in the distribution problem:
\begin{theorthogonalproblem} \label{orthogonalproblem}
We are given the same data $n$, $q$, $M$, $s$, $V$, $H$ as in the distribution problem. Is there a metric $g$ which solves the distribution problem for these data such that $V$ and $H$ are $g$-orthogonal?
\end{theorthogonalproblem}

Or we can demand only one half of what we demanded in the distribution problem:
\begin{thetimedistributionproblem} \label{timedistributionproblem}
We are given $n$, $q$, $M$, $s$ with the same properties as in the distribution problem, and we are given a $q$-plane distribution $V$ on $M$. Is there a smooth semi-Riemannian metric $g$ of index $q$ on $M$ whose scalar curvature is $s$, such that $V$ is timelike with respect to $g$?
\end{thetimedistributionproblem}

\begin{thespacedistributionproblem} \label{spacedistributionproblem}
We are given $n$, $q$, $M$, $s$ with the same properties as in the distribution problem, and we are given an $(n-q)$-plane distribution $H$ on $M$. Is there a smooth semi-Riemannian metric $g$ of index $q$ on $M$ whose scalar curvature is $s$, such that $H$ is spacelike with respect to $g$?
\end{thespacedistributionproblem}

\subsection{The topology of metrics and distributions} \label{ONEONETHREE}

For the reader's convenience, Appendix \ref{AppendixC} contains an extensive discussion of the connection between $q$-plane resp.\ $(n-q)$-plane distributions on an $n$-manifold $M$ and semi-Riemannian metrics of index $q$ on $M$. This connection is important for the homotopy class problem below. Here we repeat some of the basic facts and definitions from Appendix \ref{AppendixC}.

\begin{definition} \label{maximally}
Let $n\in\N$ and $q\in\set{0,\dots,n}$, let $(M,g)$ be an $n$-dimensional semi-Riemannian manifold of index $q$. A distribution on $M$ is called \emph{maximally timelike} if and only if it is timelike and has rank $q$. A distribution on $M$ is called \emph{maximally spacelike} if and only if it is spacelike and has rank $n-q$. If $V$ is any distribution on $M$, then the \emph{$g$-orthogonal distribution of $V$}, denoted by $\bot_gV$, is the distribution whose fibre over each point $x\in M$ is $\bot_gV_x = \set{w\in T_xM \suchthat g(v,w)=0 \text{ for all $v\in V_x$}}$ (where $V_x$ denotes the fibre of $V$ over $x$).
\end{definition}

\begin{definition}[$\Metr_q(M)$, $\Distr_q(M)$]
Let $M$ be an $n$-manifold, let $q\in\set{0,\dots,n}$. We denote the set of all $q$-plane distributions on $M$ by $\Distr_q(M)$ and equip it with the compact-open $C^\infty$ topology. $\Metr_q(M)$ denotes the set of all semi-Riemannian metrics of index $q$ on $M$, also equipped with the compact-open $C^\infty$-topology.
\end{definition}

\begin{definition}[$i$th factor distribution] \label{factordistribution}
For $i\in\set{1,2}$, let $M_i$ be an $n_i$-manifold. We define the \emph{first-factor distribution on $M_1\times M_2$} to be the $n_1$-plane distribution on $M_1\times M_2$ which is everywhere tangential to $M_1$, i.e.\ whose value in each point $(x_1,x_2)\in M_1\times M_2$ is the vector space $(T_{x_1}M_1)\oplus\set{0} \subseteq (T_{x_1}M_1)\oplus(T_{x_2}M_2) = T_{(x_1,x_2)}(M_1\times M_2)$. We define the \emph{second-factor distribution on $M_1\times M_2$} to be the $n_2$-plane distribution on $M_1\times M_2$ which is everywhere tangential to $M_2$.
\end{definition}

\begin{definition}
As usual, $\pi_0(X)$ denotes the set of (path-)connected components of a topological space $X$. Let $M$ be an $n$-manifold, let $q\in\set{0,\dots,n}$. There is a canonical bijection $\TMC\colon \pi_0(\Distr_q(M)) \to \pi_0(\Metr_q(M))$ which sends the connected component of each $V\in\Distr_q(M)$ to the connected component of $\Metr_q(M)$ which contains all metrics which make $V$ timelike. Analogously, there is a canonical bijection $\SMC\colon \pi_0(\Distr_{n-q}(M)) \to \pi_0(\Metr_q(M))$ which sends the connected component of each $H\in\Distr_{n-q}(M)$ to the connected component of $\Metr_q(M)$ which contains all metrics which make $H$ spacelike.
\smallskip\\
We define $\TDC\colon \pi_0(\Metr_q(M)) \to \pi_0(\Distr_q(M))$ to be the inverse of $\TMC$. We define $\SDC\colon \pi_0(\Metr_q(M))$ $\to \pi_0(\Distr_{n-q}(M))$ to be the inverse of $\SMC$. The map $\CDC \define \SDC\compose\TMC = \TDC\compose\SMC \colon \pi_0(\Distr_q(M))$ $\to \pi_0(\Distr_{n-q}(M))$ sends the connected component of each $V\in\Distr_q(M)$ to the connected component of $\Distr_{n-q}(M)$ that contains all distributions which are complementary to $V$.
\end{definition}

\subsection{The homotopy class version} \label{ONEONEFOUR}

Recall that the expressions \emph{homotopy class of $\sigma$} and \emph{path component of $\sigma$ (in $C^\infty(M\ot E)$)} are used synonymously for a section $\sigma$ in a fibre bundle $E\to M$. That explains the name of the following version of the prescribed scalar curvature problem.

\begin{thehomotopyproblem} \label{homotopyclassproblem}
Let $n$ be an integer $\geq2$, and let $q\in\set{0,\dots,n}$. We are given an $n$-dimensional smooth manifold $M$, a smooth real-valued function $s$ on $M$, and either
\begin{description}
\item[(metric component variant)]
a path component $C_0$ of $\Metr_q(M)$; or \item[(time component variant)]
a path component $C_1$ of $\Distr_q(M)$; or \item[(space component variant)]
a path component $C_2$ of $\Distr_{n-q}(M)$. \end{description}
Is there a smooth semi-Riemannian metric $g$ with index $q$ on $M$ whose scalar curvature is $s$, such that
\begin{description}
\item[(metric component variant)]
$g\in C_0$? \;\;resp.\ \item[(time component variant)]
$g$ makes some element of $C_1$ timelike? \;\;resp.\ \item[(space component variant)]
$g$ makes some element of $C_2$ spacelike? \end{description}
\end{thehomotopyproblem}

\begin{remarks}{\ }
\begin{enumerate}
\item
All three variants are equivalent in the following sense:
\smallskip\\
A metric $g$ solves the metric variant with respect to $C_0\in\pi_0(\Metr_q(M))$ if and only if it solves the time component variant with respect to $C_1=\TDC(C_0)\in\pi_0(\Distr_q(M))$. (Put another way, $g$ solves the time component variant with respect to $C_1\in\pi_0(\Distr_q(M))$ if and only if it solves the metric variant with respect to $C_0=\TMC(C_1)\in\pi_0(\Metr_q(M))$.)
\smallskip\\
A metric $g$ solves the metric variant with respect to $C_0\in\pi_0(\Metr_q(M))$ if and only if it solves the space component variant with respect to $C_2=\SDC(C_0)\in\pi_0(\Distr_{n-q}(M))$. (Put another way, $g$ solves the space component variant with respect to $C_2\in\pi_0(\Distr_{n-q}(M))$ if and only if it solves the metric variant with respect to $C_0=\SMC(C_2)\in\pi_0(\Metr_q(M))$.)
\smallskip\\
The equivalence of the time and space variants can also be expressed directly: A metric $g$ solves the time component variant with respect to $C_1\in\pi_0(\Distr_q(M))$ if and only if it solves the space component variant with respect to $C_2=\CDC(C_1)\in\pi_0(\Distr_{n-q}(M))$.
\item
Let $n,q,M,s$ satisfy the assumptions of the homotopy class problem. If $g\in\Metr_q(M)$ is a solution of the time [resp.\ space] distribution problem for a given $V\in\Distr_q(M)$ [resp.\ $H\in\Distr_{n-q}(M)$], then $g$ is a solution of the time [space] variant of the homotopy class problem with respect to $C_1$ being the homotopy class of $V$ [$C_2$ being the homotopy class of $H$]. Put more sloppily, the homotopy problem makes the (time/space) distribution problem easier in the sense that we do not want to make a \emph{fixed} given distribution timelike [spacelike], but any distribution which lies in the same homotopy class as a fixed given distribution. (One can easily visualise the differences between the plain and the homotopy class and the distribution problems in the case where $M$ is the $2$-torus.)
\end{enumerate}
\end{remarks}

\subsection{The diffeotopy class version} \label{ONEONEFIVE}

Let us discuss one more version of the prescribed scalar curvature problem. Its degree of difficulty lies between those of the distribution problem and the homotopy class problem.
\medskip\\
Let $M,N$ be $n$-manifolds, and let $\varphi\colon M\to N$ be a diffeomorphism. Recall that if $V$ is a $q$-plane distribution on $N$, then $\varphi^\ast(V)$ is the $q$-plane distribution on $M$ defined by $\varphi^\ast(V)_x = (\varphi_\ast)_x^{-1}(V_{\varphi(x)})$ for all $x\in M$ (where a subscript term denotes evaluation at the corresponding point, and $(\varphi_\ast)_x \colon T_xM \to T_{\varphi(x)}N$ is a vector space isomorphism).
\smallskip\\
If $\varphi\colon M\to N$ is any map and $s\in C^\infty(N,\R)$, then $\varphi^\ast s\in C^\infty(M,\R)$ is defined by $\varphi^\ast s = s\compose\varphi$. \smallskip\\
Let $M$ be a manifold. Recall that a smooth map $\Phi\colon M\times[0,1]\to M$ is a \emph{diffeotopy} if and only if (A) the map $\Phi_t\colon M\to M$ given by $x\mapsto \Phi(x,t)$ is a diffeomorphism for every $t\in[0,1]$, and (B) $\Phi_0$ is the identity. As usual, $\Diff^0(M)$ denotes the set of all diffeomorphisms $\varphi\colon M\to M$ such that there is a diffeotopy $\Phi\colon M\times[0,1]\to M$ with $\varphi = \Phi_1$. In other words, $\Diff^0(M)$ is the path component of the identity in the diffeomorphism group $\Diff(M)$ with respect to (the $C^\infty(M,M)$ subspace topology of) the compact-open $C^\infty$-topology (or, equivalently, any compact-open $C^r$-topology, where $r\in\N$). $\Diff^0(M)$ is a subgroup of $\Diff(M)$.
\smallskip\\
We call the orbits of the right group action $\Distr_q(M)\times\Diff^0(M) \to \Distr_q(M)$, $(V,\varphi)\mapsto\varphi^\ast(V)$ the \emph{diffeotopy classes} of $q$-plane distributions on $M$. This explains the name of the following problem.

\begin{thediffeotopyproblem} \label{diffeotopyclassproblem}
Let $n$ be an integer $\geq2$, and let $q\in\set{0,\dots,n}$. We are given an $n$-dimensional smooth manifold $M$, a smooth real-valued function $s$ on $M$, and two complementary distributions $V$ and $H$ on $M$, where $V$ is of rank $q$. Are there
\begin{itemize}
\item
a smooth semi-Riemannian metric $g$ of index $q$ on $M$ whose scalar curvature is $s$, and
\item
a diffeomorphism $\varphi\in\Diff^0(M)$ such that the distribution $\varphi^\ast(V)$ is timelike and the distribution $\varphi^\ast(H)$ is spacelike with respect to $g$?
\end{itemize}
\end{thediffeotopyproblem}

\begin{remarks}{\ } \label{diffeotopyremarks}
\begin{enumerate}
\item
The distribution problem and the diffeotopy class problem have the same data as input. If, for given input $s,V,H$, the distribution problem is solvable --- i.e., a suitable metric $g$ exists ---, then the diffeotopy class problem is solvable for this input, by the same metric $g$ and the diffeomorphism $\varphi = \id_M$. In this sense, the distribution problem is harder than the diffeotopy class problem.
\smallskip\\
But if the function $s$ is constant, then the converse holds as well, that is, both problems are equivalent: If, for given $V,H$ and constant $s$, there exist a metric $g$ and a diffeomorphism $\varphi$ which solve the diffeotopy class problem, then the metric $(\varphi^{-1})^\ast g$ solves the distribution problem for the same input data. Namely, we have then $\scal_{(\varphi^{-1})^\ast g} = (\varphi^{-1})^\ast \scal_g = (\varphi^{-1})^\ast s = s$ (since $s$ is constant); $V = (\varphi^{-1})^\ast \varphi^\ast V$ is timelike with respect to $(\varphi^{-1})^\ast g$ since $\varphi^\ast V$ is timelike with respect to $g$; and, analogously, $H$ is spacelike with respect to $(\varphi^{-1})^\ast g$.
\smallskip\\
To put it more sloppily: The idea behind the diffeotopy class problem is to make the distribution problem easier by allowing that the given distributions are $\Diff^0(M)$-modified relative to the prescribed function $s$; or, equivalently, that $s$ is $\Diff^0(M)$-modified relative to the distributions. If $s$ is constant, then $\Diff^0(M)$-modifications do not change it, so we are still left with the distribution problem. (Observe that the homotopy class problem is easier than the distribution problem even if $s$ is constant.)
\item
If the distributions $V_0,V_1\in\Distr_q(M)$ lie in the same diffeotopy class, then they lie in the same homotopy class: any diffeotopy $\Phi\colon[0,1]\times M\to M$ with $\Phi_1^\ast(V_0) = V_1$ yields a path $V_{\blank} \colon [0,1]\to\Distr_q(M)$ from $V_0$ to $V_1$, via $V_t \define \Phi_t^\ast(V_0)$.
\smallskip\\
The homotopy class problem is therefore easier than the diffeotopy class problem, in the following sense: If $g$ is a solution of the diffeotopy class problem for given distributions $V\in\Distr_q(M)$ and $H\in\Distr_{n-q}(M)$, then it is a solution of the time component variant [space component variant] of the homotopy class problem with respect to $C_0$ being the homotopy class of $V$ [$C_1$ being the homotopy class of $H$].
\item
The main additional information that a solution of the diffeotopy class problem provides compared to an arbitrary solution of the homotopy class problem lies in the \emph{twistedness} of the prescribed distribution (cf.\ Subsection \ref{twistsection}). For instance, we might ask whether we can find a Lorentzian metric with scalar curvature $s$ \emph{which admits a maximally spacelike foliation}, i.e.\ which admits a maximally spacelike distribution which is \emph{integrable}. An arbitrary solution of the homotopy class problem might not have this property. But if $H$ is integrable, then any distribution of the form $\varphi^\ast(H)$, where $\varphi\in\Diff(M)$, is integrable, too; so if we can solve the diffeotopy class problem for an integrable $H$, then our solution metric admits a maximally spacelike foliation.
\smallskip\\
Similarly, we could consider a distribution $H$ which is everywhere twisted (cf.\ \ref{twistsection}), i.e.\ nowhere integrable. Then every solution of the diffeotopy class problem with respect to this $H$ admits a maximally spacelike distribution which is everywhere twisted.
\item
One might have the idea to create even more versions of the prescribed scalar curvature problem by allowing that $V$ and $H$ in the diffeotopy class problem get modified by two different elements of $\Diff^0(M)$; or by allowing that they get modified by an arbitrary diffeomorphism of $M$ (not necessarily contained in $\Diff^0(M)$). I think that these versions have no practical or theoretical value: they do not make any of the mentioned problems easier to solve in practice, and their solutions would not provide any interesting information beyond what one can get out of the problems we have already introduced.
\end{enumerate}
\end{remarks}

Like the distribution problem, we can split the diffeotopy class problem into a time half and a space half:

\begin{thetimediffeotopyproblem} \label{timediffeotopyclassproblem}
We are given $n$, $q$, $M$, $s$ with the same properties as in the diffeotopy class problem, and we are given a $q$-plane distribution $V$ on $M$. Are there a smooth semi-Riemannian metric $g$ of index $q$ on $M$ whose scalar curvature is $s$, and a diffeomorphism $\varphi\in\Diff^0(M)$ such that $\varphi^\ast(V)$ is timelike with respect to $g$?
\end{thetimediffeotopyproblem}

\begin{thespacediffeotopyproblem} \label{spacediffeotopyclassproblem}
We are given $n$, $q$, $M$, $s$ with the same properties as in the diffeotopy class problem, and we are given an $(n-q)$-plane distribution $H$ on $M$. Are there a smooth semi-Riemannian metric $g$ of index $q$ on $M$ whose scalar curvature is $s$, and a diffeomorphism $\varphi\in\Diff^0(M)$ such that $\varphi^\ast(H)$ is spacelike with respect to $g$?
\end{thespacediffeotopyproblem}

The Remarks \ref{diffeotopyremarks} generalise mutatis mutandis to the time/space diffeotopy class problem.

\subsection{Summary and outlook} \label{ONEONESIX}

\begin{figure}[tb]
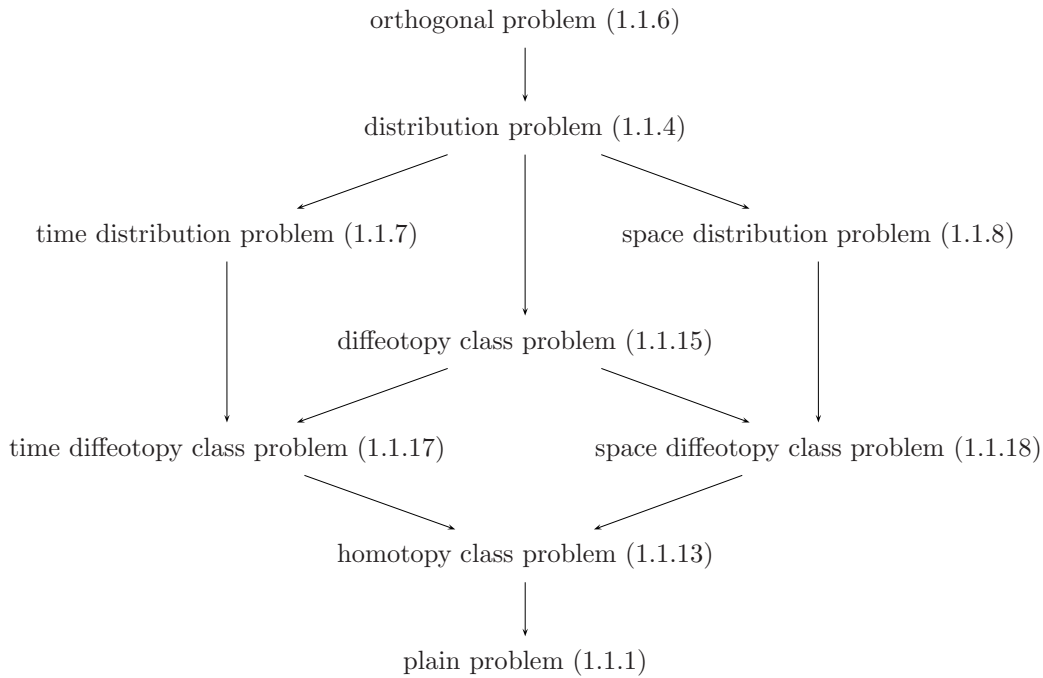

\begin{center}
\newcommand{\centre}[1]{\multicolumn{2}{c}{#1}}
\begin{tabular}{c@{\hspace{2cm}}c}
\centre{\Rnode{p0}{orthogonal problem (\ref{orthogonalproblem})}}\\[1cm]
\centre{\Rnode{p1}{distribution problem (\ref{distributionproblem})}}\\[1cm]
\Rnode{p2l}{time distribution problem (\ref{timedistributionproblem})} &\Rnode{p2r}{space distribution problem (\ref{spacedistributionproblem})}\\[1cm]
\centre{\Rnode{p3}{diffeotopy class problem (\ref{diffeotopyclassproblem})}}\\[1cm]
\Rnode{p4l}{time diffeotopy class problem (\ref{timediffeotopyclassproblem})} &\Rnode{p4r}{space diffeotopy class problem (\ref{spacediffeotopyclassproblem})}\\[1cm]
\centre{\Rnode{p5}{homotopy class problem (\ref{homotopyclassproblem})}}\\[1cm]
\centre{\Rnode{p6}{plain problem (\ref{plainproblem})}}
\end{tabular}
\psset{arrows=->,nodesep=5pt,linewidth=0.3pt}
\ncline{p0}{p1} \ncline{p1}{p2l} \ncline{p1}{p2r} \ncline{p1}{p3} \ncline{p3}{p4l} \ncline{p3}{p4r} \ncline{p2l}{p4l} \ncline{p2r}{p4r} \ncline{p4l}{p5} \ncline{p4r}{p5} \ncline{p5}{p6}
\end{center}
\caption{Diagram of the versions of the semi-Riemannian prescribed scalar curvature problem.} \label{problemdiagram}
\end{figure}

The different versions of the semi-Riemannian prescribed scalar curvature problem are summarised in Figure \ref{problemdiagram}, whose arrows generate the order relation \emph{solutions of \dots\ yield solutions of \dots} (i.e.\ the order relation \emph{\dots\ is harder than\ \dots}).
\smallskip

All these versions reduce to the same problem in the Riemannian case, because there the only possible timelike distribution is the unique distribution of rank $0$, and the only possible spacelike distribution is the whole tangent bundle. In my opinion, the most natural generalisation of the Riemannian prescribed scalar curvature problem is the homotopy class problem. We will therefore concentrate on solving that one.
\smallskip\\
This turns out to be possible in many cases on manifolds of dimension $n\geq3$. In the $2$-dimensional case, however, we solve only the plain problem.
\smallskip\\
For metrics of index $\in\set{3,\dots,n-3}$, we will even be able to solve the orthogonal problem for a set of pairs $(V,H)$ which is $C^0$-dense and $C^1$-open in $\set{(V,H)\in\Distr_q(M)\times\Distr_{n-q}(M) \suchthat \text{$V$ is complementary to $H$}}$.

\smallskip
We will discuss the solvability of the diffeotopy class problem in the case where $V$ and $H$ are the first- resp.\ second-factor distributions on a product manifold (cf.\ \ref{productexample}).


\section{How to solve the problem? The main idea} \label{ONETWO}

\newcommand{\gm}{\mathfrak{g}}
\newcommand{\hm}{\mathfrak{h}}

\emph{This section explains on an informal level the idea of our approach to the prescribed scalar curvature problem for pseudo-Riemannian metrics. It contains neither definitions nor facts which are used later.}
\medskip\\
It is not completely obvious how to generalise the methods which Kazdan and Warner applied in the Riemannian case (cf.\ Appendix \ref{D1}) to the pseudo-Riemannian case of the prescribed scalar curvature problem. The direct approach would be to try again to find a pseudo-Riemannian metric $h$ which is conformal to a given metric $g$ of index $q$ and satisfies the equation $\scal_h = s$. As I already mentioned in Remark \ref{distributionremarks}, this does not work since the resulting equation
\begin{equation} \label{conformalequation}
0 = 2(n-1)\laplace_g(\kappa) -\frac{n(n-1)}{\kappa}\abs{d\kappa}_g^2 +\kappa\scal_g -\frac{1}{\kappa}s
\end{equation}
(the precise analogue of the Riemannian case equation \eqref{yamabe1}) is not elliptic. If $g$ is Lorentzian, the equation is hyperbolic; if $g$ is ultra-Lorentzian, it is ultrahyperbolic.
\medskip\\
It would be nice if we could construct an \emph{elliptic} equation whose solvability implies the solvability of the prescribed curvature problem. That is indeed possible. The equation (cf.\ Theorem \ref{THEEQUATION}) will --- at least in a general situation --- have a much more complicated form than Equation \eqref{conformalequation}, but this is more than compensated by the virtue of ellipticity.
\medskip\\
The purpose of this section is to sketch in a simple (but not too simple) special case the construction of the elliptic equation. It might be surprising that a global problem in pseudo-Riemannian geometry can be solved via an \emph{elliptic} equation at all; that is presumably not what one would expect. Therefore I would like to convince the reader that this is indeed possible, and that the ellipticity is not caused by some hypothetical sign error in the quite long (but straightforward) calculations in Chapter 3 and 4 where we determine the explicit form of the equation. The basic idea is quite simple, as we will see.

\subsection{The equation in a special case} \label{ONETWOONE}

We consider a simple special case of the Lorentzian distribution problem, in which $M$ is a product manifold: Let $N$ be any manifold of dimension $n-1$, where $n\geq2$; let $M$ be the $n$-manifold $S^1\times N$; let $V$ be the first-factor distribution (cf.\ Definition \ref{factordistribution}) on $S^1\times N$; and let $H$ be the second-factor distribution on $S^1\times N$.
\smallskip\\
Given any function $s\in C^\infty(M,\R)$, we want to find a Lorentzian metric $h$ on $M$ which makes $V$ timelike and $H$ spacelike, and whose scalar curvature is $s$. (Note that in spite of the special form of $M,V,H$, this problem is not easy, except for very special functions $s$. We want time to move in circles, and that would cause difficult periodicity conditions for an initial value approach via a hyperbolic equation, for instance.)
\smallskip\\
We choose a Riemannian metric $g_N$ on $N$ and try to find a suitable metric $h$ in the form
\[
h(\kappa,f) \define \frac{1}{\kappa^2}\Big((-\frac{1}{f^2}dt^2)\oplus g_N\Big) \equiv \frac{1}{\kappa^2}\Big(-\frac{1}{f^2}\pi_{S^1}^\ast(dt^2) +\pi_{N}^\ast(g_N)\Big) \;\;,
\]
where $dt^2$ is the standard Riemannian metric on $S^1$, the maps $\pi_{S^1}\colon M\to S^1$ and $\pi_N\colon M\to N$ are the obvious projections, and $\kappa,f\in C^\infty(M,\R_{>0})$. For all functions $\kappa,f\in C^\infty(M,\R_{>0})$, this $h(\kappa,f)$ is obviously a Lorentzian metric which makes $V$ timelike, and makes $H$ orthogonal to $V$ and thus spacelike.
\smallskip\\
Now we have to compute the scalar curvature of $h(\kappa,f)$ explicitly in terms of $\kappa$ and $f$.

\subsubsection{First step}

Denoting the metric $h(1,f)$ by $h(f)$ for the sake of brevity, we get \begin{equation} \label{babyconform}
\scal_{h(\kappa,f)} = \kappa^2\scal_{h(f)} +2(n-1)\kappa\laplace_{h(f)}(\kappa) -n(n-1)\eval{d\kappa}{d\kappa}_{h(f)} \;\;,
\end{equation}
by the well-known formula for the conformal change of a semi-Riemannian metric; cf.\ e.g.\ \cite{Besse}, 1.159(f).\footnote{Note that Besse uses a sign convention for the Laplacian which is opposite to ours.} (We rederive this formula as a special case of a much more general one in Chapter 3; cf.\ Section \ref{conformalsection}.)
\smallskip\\
Here we used the following notation: $\laplace_h(u)\in C^\infty(M,\R)$ is the Laplacian of $u\in C^\infty(M,\R)$ with respect to a semi-Riemannian metric $h$. For $u_0,u_1\in C^\infty(M,\R)$, the function $\eval{du_0}{du_1}_h \in C^\infty(M,\R)$ is the pointwise scalar product (with respect to $h$) of the $1$-forms $du_0,du_1\in C^\infty(M\ot T^\ast M)$; in other words, $\eval{du_0}{du_1}_h = h(\grad_hu_0,\grad_hu_1)$, where $\grad_hu_i$ is the gradient vector field of $u_i$ with respect to $h$.

\subsubsection{Second step}

For every function $u\in C^\infty(M,\R)$, we define the function $\laplace_{g_N}(u)\in C^\infty(M,\R)$ in a natural way: For each $t\in S^1$, we consider the function $u_t\in C^\infty(N,\R)$ given by $x\mapsto u(t,x)$. We define the value of $\laplace_{g_N}(u)$ in the point $(t,x)\in S^1\times N$ to be the value of $\laplace_{g_N}(u_t)\in C^\infty(N,\R)$ in the point $x$.
\smallskip\\
For all functions $u,v\in C^\infty(M,\R)$, we define the function $\eval{du}{dv}_{g_N} \in C^\infty(M,\R)$ in an equally natural way: we assign to each $(t,x)\in M$ the value of $\eval{du_t}{dv_t}_{g_N}\in C^\infty(N,\R)$ in $x$. Finally, we consider $\scal_{g_N}$ as a function on $M$ whose value in $(t,x)$ is just the value of $\scal_{g_N}\in C^\infty(N,\R)$ in $x$. With these notations, we get for every function $u\in C^\infty(M,\R)$:
\begin{subequations} \label{babystretch} \begin{align}
\scal_{h(f)} &= \scal_{g_N} +\frac{2}{f}\laplace_{g_N}(f) -\frac{4}{f^2}\eval{df}{df}_{g_N}^2 \;\;, \label{babya}\\
\laplace_{h(f)}(u) &= \laplace_{g_N}(u) -f^2\partial_t\partial_tu -\frac{1}{f}\eval{df}{du}_{g_N} -f(\partial_tf)(\partial_tu) \;\;, \label{babyb}\\
\eval{du}{du}_{h(f)} &= \eval{du}{du}_{g_N} -f^2(\partial_tu)^2 \;\;. \label{babyc} \end{align} \end{subequations}
Here $\partial_t$ denotes a partial derivative with respect to the canonical (local) coordinate on $S^1 = \R/\Z$. Note that $\partial_t\partial_tu = \laplace_{dt^2}(u)$ and $(\partial_tf)(\partial_tu) = \eval{df}{du}_{dt^2}$ (we use a similar notation here for the metric $dt^2$ on the first factor $S^1$ as for the metric $g_N$ on the second factor $N$).
\medskip\\
The formulae \eqref{babystretch} are special cases of more general formulae which we will prove in Chapter 3; cf.\ Theorem \ref{stretchformulae}. They are, from a different direction, also special cases of O'Neill's equations for semi-Riemannian submersions (cf.\ e.g.\ \S\S9.B-D in \cite{Besse}), because the projection $\pi_N \colon (M,h(f)) \to (N,g_N)$ is a semi-Riemannian submersion. In particular, $h(f)$ is a warped product metric if the function $f$ on $S^1\times N$ does only depend on the $N$ factor: $M$ equipped with $h(f)$ is the warped product $N\times_{1/f}S^1$ then. In that case, the equations \eqref{babystretch} coincide with the well-known warped product formulae; cf. e.g.\ \cite{ONeill}, Exercise 7.13.
\smallskip\\
However, I don't know a reference from which \eqref{babystretch} could be read off immediately; and I cannot refer to Chapter 3 because I promised to explain --- without using the computation that we will do there --- why there is an approach to the prescribed scalar curvature problem via an elliptic equation. Luckily, all we need to understand this approach is those summands on the right-hand sides of \eqref{babystretch} which contain \emph{second} derivatives of the functions $f$ and $u$; i.e.\ the summands containing $\laplace_{g_N}(f)$, $\laplace_{g_N}(u)$, or $\partial_t\partial_tu$. So let us just check that the equations \eqref{babystretch} are correct up to terms of (derivative) order less than $2$.
\medskip\\
We will verify that the local coordinate expressions for the left-hand sides coincide with those for the right-hand sides, up to lower order terms. Let $(x_1,\dots,x_{n-1})$ be local coordinates on $N$, let $\gm$ be the $(n-1)\times(n-1)$-matrix which represents the Riemannian metric $g_N$ with respect to these local coordinates, and let $x_0=t$ be the standard (local) coordinate on $S^1$. Then the Lorentzian metric $h(f)$ is represented, with respect to the local coordinates $(x_0,\dots,x_{n-1})$ on $M$, by the block matrix-valued function
\[
\hm \define \begin{pmatrix} -1/f^2 &0\\ 0 &\gm\end{pmatrix} \;\;.
\]
Its pointwise inverse (whose components will be denoted by $\hm^{ij}$, as usual) is
\[
\hm^{-1} = \begin{pmatrix} -f^2 &0\\ 0 &\gm^{-1}\end{pmatrix} \;\;.
\]
By well-known formulae (involving the Riemann curvature tensor and Christoffel symbols of $h(f)$ in the intermediate steps; cf.\ e.g.\ \cite{ONeill}, pp.\ 88, 76, 62, 87), we have
\[ \begin{split}
\scal_{h(f)} &= \sum_{i,j,k}\hm^{ij}R^k_{ijk}\\[2mm]
&= \sum_{i,j,k}\hm^{ij}\big(\partial_k\varGamma^k_{ji} -\partial_j\varGamma^k_{ki}\big) +\lot\\[2mm]
&= \frac{1}{2}\sum_{i,j,k,m}\hm^{ij}\bigg(\partial_k\Big(\hm^{km}(\partial_j\hm_{im} +\partial_i\hm_{jm} -\partial_m\hm_{ji})\Big) -\partial_j\Big(\hm^{km}(\partial_k\hm_{im} +\partial_i\hm_{km} -\partial_m\hm_{ki})\Big)\bigg) +\lot\\[2mm]
&= \frac{1}{2}\sum_{i,j,k,m}\hm^{ij}\hm^{km}\bigg((\partial_k\partial_j\hm_{im} +\partial_k\partial_i\hm_{jm} -\partial_k\partial_m\hm_{ji}) -(\partial_j\partial_k\hm_{im} +\partial_j\partial_i\hm_{km} -\partial_j\partial_m\hm_{ki})\bigg) +\lot
\end{split} \]

\[ \begin{split}
&= \frac{1}{2}\bigg(\sum_{i,j,k,m}\hm^{ij}\hm^{km}\partial_k\partial_i\hm_{jm} -\!\!\sum_{i,j,k,m}\hm^{ij}\hm^{km}\partial_k\partial_m\hm_{ji} -\!\!\sum_{i,j,k,m}\hm^{ij}\hm^{km}\partial_j\partial_i\hm_{km} +\!\!\sum_{i,j,k,m}\hm^{ij}\hm^{km}\partial_j\partial_m\hm_{ki}\bigg)\\
&\mspace{701mu} +\lot\\
&= \frac{1}{2}\bigg(\sum_{i,k}\hm^{i0}\hm^{k0}\partial_k\partial_i\hm_{00} -\sum_{k,m}\hm^{00}\hm^{km}\partial_k\partial_m\hm_{00} -\sum_{i,j}\hm^{ij}\hm^{00}\partial_j\partial_i\hm_{00} +\sum_{j,m}\hm^{0j}\hm^{0m}\partial_j\partial_m\hm_{00}\bigg) +\lot\\
&= \hm^{00}\hm^{00}\partial_0\partial_0\hm_{00} -\hm^{00}\sum_{i,j}\hm^{ij}\partial_i\partial_j\hm_{00} +\lot\\
&= -\hm^{00}\sum_{i,j=1}^{n-1}\hm^{ij}\partial_i\partial_j\hm_{00} +\lot\\
&= -f^2\sum_{i,j=1}^{n-1}\gm^{ij}\partial_i\partial_j(1/f^2) +\lot\\
&= -f^2\laplace_{g_N}(1/f^2) +\lot\\
&= \frac{2}{f}\laplace_{g_N}(f) +\lot \;\;,
\end{split} \]
where all sums run from $0$ to $n-1$, except where indicated; $\partial_i$ denotes a partial derivative with respect to the coordinate $x_i$; and $\lot$ denotes summands (not the same ones in each step) which contain no second derivatives of $\hm_{00}$ or $\hm^{00}$ and thus no second derivatives of $f$.
\smallskip\\
This proves that \eqref{babya} is correct at least up to terms of order less than $2$. In a similar way, we compute
\[ \begin{split}
\laplace_{h(f)}(u) &= \sum_{i,j}\hm^{ij}\partial_i\partial_ju +\lot\\
&= -f^2\partial_0\partial_0u +\sum_{i,j=1}^{n-1}\gm^{ij}\partial_i\partial_ju +\lot\\
&= -f^2\partial_t\partial_tu +\laplace_{g_N}(u) +\lot \;\;,
\end{split} \]
which proves that \eqref{babyb} is correct at least up to terms of order less than $2$ in $f$ or $u$. Equation \eqref{babyc} is true up to terms of order less than $2$ for the simple reason that neither the left-hand nor the right-hand side contains terms of order $\geq2$. (And one can easily check \eqref{babyc} directly.)

\subsubsection{Third step}

Equations \eqref{babyconform} and \eqref{babystretch} together yield
\begin{equation} \label{babymix} \begin{split}
\scal_{h(\kappa,f)} &= \kappa^2\bigg(\scal_{g_N} +\frac{2}{f}\laplace_{g_N}(f) -\frac{4}{f^2}\eval{df}{df}_{g_N}^2\bigg) -n(n-1)\bigg(\eval{d\kappa}{du\kappa}_{g_N} -f^2(\partial_t\kappa)^2\bigg)\\ &\mspace{20mu}+2(n-1)\kappa\bigg(\laplace_{g_N}(\kappa) -f^2\partial_t\partial_t\kappa -\frac{1}{f}\eval{df}{d\kappa}_{g_N} -f(\partial_tf)(\partial_t\kappa)\bigg)\\
&= \frac{2\kappa^2}{f}\laplace_{g_N}(f) +2(n-1)\kappa\laplace_{g_N}(\kappa) -2(n-1)\kappa f^2\partial_t\partial_t\kappa +\lot \;\;.
\end{split} \end{equation}
We have not checked that the equations \eqref{babystretch} are really correct (be sure, they are), but we have checked the second-order terms and hence can conclude that \eqref{babymix} is true at least up to terms of order less than $2$ in $\kappa$ or $f$, which is enough for our present purposes.
\smallskip\\
Now we consider the Riemannian product metric $g \define dt^2\oplus g_N$ on $M$. Because the $g$-Laplacian of a function $u\in C^\infty(M,\R)$ is $\laplace_g(u) = \laplace_{g_N}(u) +\partial_t\partial_tu$, we can write \eqref{babymix} in the form
\begin{equation} \label{babymix2}
\scal_{h(\kappa,f)} = \frac{2\kappa^2}{f}\laplace_g(f) -\frac{2\kappa^2}{f}\partial_t\partial_tf +2(n-1)\kappa\laplace_g(\kappa) -2(n-1)\kappa(1+f^2)\partial_t\partial_t\kappa +\lot \;\;.
\end{equation}
This completes our calculation of $\scal_{h(\kappa,f)}$ in terms of $\kappa$ and $f$ (up to lower order terms).

\subsubsection{The main idea}

For every fixed $f\in C^\infty(M,\R_{>0})$, the right-hand side of \eqref{babymix2} is evidently a \emph{hyperbolic} differential operator in $\kappa$. For every fixed $\kappa\in C^\infty(M,\R_{>0})$, the right-hand side of \eqref{babymix2} is a \emph{parabolic} differential operator in $f$. (Here we use a general definition of parabolicity which does not exclude degenerate cases: we say that a linear second-order differential operator $P\colon C^\infty(M,\R) \to C^\infty(M,\R)$ is parabolic if and only if the value of its symbol $\sigma(P)\in C^\infty(M\ot\Sym(T^\ast M))$ in each point of $M$ is a symmetric bilinear form whose diagonalisations have $n-1$ positive diagonal entries, and one diagonal entry equal to $0$. Parabolicity of nonlinear operators is defined in terms of parabolicity of linearisations.)
\smallskip\\
So we have, loosely speaking, a two-($C^\infty(M,\R_{>0})$-)parameter family of functions which contains a hyperbolic one-parameter subfamily and a parabolic one-parameter subfamily. It seems natural to ask whether it contains also an \emph{elliptic} one-parameter subfamily.
\smallskip\\
We try to find such a subfamily by considering $\scal_{h(K\compose f,f)}$ for some function $K\in C^\infty(\R_{>0},\R_{>0})$. Using the notation $K(f)\define K\compose f\in C^\infty(M,\R_{>0})$, we get
\[ \begin{split}
\scal_{h(K(f),f)} &= \frac{2K(f)^2}{f}\laplace_g(f) -\frac{2K(f)^2}{f}\partial_t\partial_tf\\
&\mspace{20mu}+2(n-1)K(f)K'(f)\laplace_g(f) -2(n-1)K(f)(1+f^2)K'(f)\partial_t\partial_tf +\lot \;\;,
\end{split} \]
where $\lot$ denotes summands which do not contain second derivatives of $f$. Let us choose the function $K\in C^\infty(\R_{>0},\R_{>0})$ given by
\[
K(x) = \Big(\frac{1+x^2}{x^2}\Big)^{\frac{1}{2(n-1)}} \;\;.
\]
Then
\[
K'(x) = \frac{1}{2(n-1)}\cdot\frac{x^2}{1+x^2}K(x)\cdot\frac{-2}{x^3} = -\frac{1}{(n-1)x(1+x^2)}K(x) \;\;,
\]
and thus
\[ \begin{split}
\scal_{h(K(f),f)} &= K(f)^2\bigg(\frac{2}{f}\laplace_g(f) -\frac{2}{f}\partial_t\partial_tf -\frac{2}{f(1+f^2)}\laplace_g(f) +\frac{2(1+f^2)}{f(1+f^2)}\partial_t\partial_tf\bigg) +\lot\\
&= \frac{2fK(f)^2}{1+f^2}\laplace_g(f) +\lot \;\;.
\end{split} \]
So the map $C^\infty(M,\R_{>0}) \to C^\infty(M,\R)$ given by $f\mapsto \scal_{h(K(f),f)}$ is a second-order \emph{elliptic} differential operator.
\medskip\\
Note that no miracle had to occur to produce this elliptic PDE: It was clear already before our three-step calculation above that functions $C_0,C_1 \colon \R_{>0}\times\R_{>0}\times\R \to \R$ exist such that for all $K\in C^\infty(\R_{>0},\R_{>0})$, the scalar curvature of $h(K(f),f)$ has the form
\[
\scal_{h(K(f),f)} = C_0(f,K(f),K'(f))\laplace_g(f) + C_1(f,K(f),K'(f))\partial_t\partial_tf +\lot \;\;;
\]
this form can be deduced from general properties of the map which assigns to each metric its scalar curvature. The point is that $C_0,C_1$ are such that a function $K\in C^\infty(\R_{>0},\R)$ exists with the following three properties:
\begin{enumerate}
\item
$K$ solves the ordinary differential equation $C_1(x,K(x),K'(x))=0$;
\item
$K$ takes indeed values in $\R_{>0}$;
\item
the function $\R_{>0}\to\R$ given by $x\mapsto C_0(x,K(x),K'(x))$ has no zeroes.
\end{enumerate}
It is not really surprising that one can find a $K$ with property (i), and that this $K$ is unique up to multiplication with a constant. That some such $K$ satisfies properties (ii) and (iii) is nice, but not exactly a miracle. (In fact, it is clear that properties (ii) and (iii) hold for the restriction $K\restrict I$ of a solution $K$ of (i) to some nonempty open interval $I\subseteq \R_{>0}$. If this interval $I$ could not be taken to be the whole interval $\R_{>0}$, then we could still consider the elliptic equation, but we had to prove that it admits a solution $f$ whose values are contained in $I$.)

\subsubsection{The elliptic equation}

Realising our prescribed function $s\in C^\infty(M,\R)$ as the scalar curvature of a Lorentzian metric of the form $h(K(f),f)$ is equivalent to finding a solution $f\in C^\infty(M,\R_{>0})$ of the second-order semilinear\footnote{Recall that a $k$th-order partial differential operator is \emph{semilinear} if and only if it is the sum of a linear $k$th-order partial differential operator and a partial differential operator of order $<k$. A semilinear PDE is a PDE of the form $P(f)=0$, where $P$ is a semilinear operator.} elliptic PDE
\[
\frac{1+f^2}{fK(f)^2}s = \frac{1+f^2}{fK(f)^2}\scal_{h(K(f),f)} \;\;,
\]
whose symbol is just the symbol of $2\laplace_g$. Including all lower-order terms, the explicit form of this PDE is
\begin{equation} \begin{split} \label{babyequation}
0 &= 2\laplace_g(f) -\frac{4f^4+\alpha(n)}{f^3(1+f^2)}\abs{df}_g^2 +\frac{\alpha(n)}{f^3}(\partial_tf)^2 +\frac{1+f^2}{f}\scal_g -f\Big(\frac{1+f^2}{f^2}\Big)^{\alpha(n)}s \;\;,
\end{split} \end{equation}
where $\alpha(n)\define\frac{n-2}{n-1}$. (Observe that $\scal_g = \scal_{g_N}$.)
\medskip\\
Our discussion of the special case $S^1\times N$ ends here. The result is as follows:
\smallskip\\
\emph{If the semilinear elliptic second-order equation \eqref{babyequation} has a solution $f\in C^\infty(M,\R_{>0})$, then there is a Lorentzian metric on $S^1\times N$ (namely the metric $h(K(f),f)$) whose scalar curvature is $s$ and which, moreover, solves the orthogonal problem \ref{orthogonalproblem} with respect to the factor distributions $V$ and $H$.}

\subsection{Some remarks} \label{ONETWOTWO}

The approach from the special case above can be generalised: to arbitrary manifolds $M$ (of dimension $\geq2$) with arbitrary prescribed complementary distribution $V,H$. The idea is simple: We choose Riemannian metrics on the vector bundles $V,H$ and consider the Riemannian metric $g= g_V\oplus g_H$ on $M$ and the semi-Riemannian metric $h(\kappa,f) = \kappa^{-2}((-f^{-2}g_V) \oplus g_H)$, whose index is $\rank(V)$.
\smallskip\\
The main differences are that the computations become much longer (recall that, even in the simple example from above, we wrote down the calculations only up to lower order terms) and that the function $K$ is given by
\[
K(x) = \Big(\frac{1+x^2}{x^{2q}}\Big)^{\frac{1}{2(n-1)}}
\]
when $V$ is an arbitrary $q$-plane distribution, i.e.\ when we want to find a metric of index $q$.

\subsubsection{The Riemannian case as a special case}

In the Riemannian case $q=0$, the resulting elliptic equation is \emph{almost} equivalent to the well-known equation \eqref{conformalequation}, which is related to a conformal change of the metric. \emph{Almost} means that our elliptic equation admits a solution $f\in C^\infty(M,\R_{>0})$ if and only if \eqref{conformalequation} admits a solution $\kappa\in C^\infty(M,\R_{>1})$. (Note that \eqref{conformalequation} is defined for all $\kappa\in C^\infty(M,\R_{>0})$.) The reason for this restriction of the allowed functions $\kappa$ is of course that our approach uses conformal factors of the form $\kappa = K(f) = (1+f^2)^{1/(2n-2)}$ in the case $q=0$.
\smallskip\\
Nevertheless, one loses nothing essential by restricting the allowed functions $\kappa$ in this way: If \eqref{conformalequation} admits a solution $\in C^\infty(M,\R_{>0})$, then there is a large set of functions $\lambda\in C^\infty(M,\R_{>0})$ such that Equation \eqref{conformalequation}, with the metric $g$ replaced by $\lambda g$, admits a solution $\in C^\infty(M,\R_{>1})$. Since the background metric $g$ is not fixed but has to be chosen suitably in the proofs anyway, it is justified to say that our approach to the semi-Riemannian prescribed scalar curvature problem is a generalisation of the approach which Kazdan and Warner chose\footnote{in \cite{KazdanWarner2a}, \cite{KazdanWarner2b}, \cite{KazdanWarner0}, \cite{KazdanWarner1}; the method in \cite{KazdanWarner2} is different} to solve the Riemannian case.

\subsubsection{How unique is our approach to the construction of an elliptic equation?}

This question is a bit too vague and general to allow a useful answer. Let us first consider the special case of the prescribed scalar curvature problem from the previous section, and let us ask a more specific question: \emph{For given functions $A,B\in C^\infty(\R_{>0},\R_{>0})$, when precisely is the map $P \colon C^\infty(M,\R_{>0}) \to C^\infty(M,\R)$ which is defined by $f\mapsto \scal_{h(A(f),B(f))}$ an elliptic operator?} (Note that the target of the functions $A,B$ has to be $\R_{>0}$, and that we can assume the domain to be $\R_{>0}$ without loss of generality.)
\smallskip\\
The answer is simple: From Equation \eqref{babymix}, we see that the second-order terms of $\frac{1}{2}\scal_{h(A(f),B(f))}$ have the form
\[
\bigg(\frac{A(f)^2}{B(f)}B'(f) +(n-1)A(f)A'(f)\bigg)\laplace_{g_N}(f) -(n-1)A(f)B(f)^2A'(f)\partial_t\partial_tf \;\;.
\]
Let $C_0\define A^2B'/B +(n-1)AA' \in C^\infty(\R_{>0},\R)$ denote the coefficient of $\laplace_{g_N}$, and let $C_1\define -(n-1)AB^2A'\in C^\infty(\R_{>0},\R)$ denote the coefficient of $\partial_t\partial_tf$.
\smallskip\\
Our operator $P$ is elliptic if and only if $C_0$ has the same sign as $C_1$ everywhere, i.e.\ if and only if $A'$ vanishes nowhere and there exists a function $E\in C^\infty(\R_{>0},\R_{>0})$ such that $C_0 = EC_1$. There are clearly lots of pairs of functions $(A,B)$ which satisfy this condition.
\medskip\\
The picture changes if we ask for which $(A,B)$ the operator $P \colon f\mapsto \scal_{h(A(f),B(f))}$ has the form $P(f) = \beta(f)Q(f)$, where $\beta\colon \R_{>0}\to\R$ is a nowhere vanishing function and $Q$ is a \emph{uniformly} elliptic operator; the uniform condition here refers to the Riemannian metric $g = dt^2\oplus g_N$. Since uniform ellipticity is important for the proofs (e.g.\ for the method of sub- and supersolutions), this question is what we are really interested in. The nowhere vanishing factor $\beta(f)$ might destroy the uniform ellipticity, but that makes no problem since we get a uniformly elliptic equation after division by $\beta(f)$.
\smallskip\\
Clearly $P$ has the form $\beta Q$ with nowhere vanishing $\beta$ and $g$-uniformly elliptic $Q$ if and only if $A'$ vanishes nowhere and there exist $a,b\in\R_{>0}$ with $a\leq b$ and a function $E\in C^\infty(\R_{>0},[a,b])$ such that $C_0 = EC_1$.
\smallskip\\
The equation $C_0 = EC_1$ is equivalent to
\[
\frac{B'}{B(1+EB^2)} = -(n-1)\frac{A'}{A} \;\;.
\]
Because of
\[
\frac{1}{2}\Big(\log\frac{B^2}{1+bB^2}\Big)' = \frac{B'}{B(1+bB^2)} \leq \frac{B'}{B(1+EB^2)} \leq \frac{B'}{B(1+aB^2)} = \frac{1}{2}\Big(\log\frac{B^2}{1+aB^2}\Big)' \;\;,
\]
this implies that there exist constants $d_0,d_1\in\R$ (which can be written as functions of $a,b,n,A(1),B(1)$) such that
\[
\frac{1}{2}\log\frac{B^2}{1+bB^2} +d_0 \leq -(n-1)\log(A) \leq \frac{1}{2}\log\frac{B^2}{1+aB^2} +d_1 \;\;.
\]
In other words, there exist constants $c_0,c_1\in\R_{>0}$ (which can be written as functions of $a,b,n,A(1),B(1)$) such that
\begin{equation} \label{strictA}
c_0\bigg(\frac{1+aB^2}{B^2}\bigg)^{\frac{1}{2(n-1)}} \leq A \leq c_1\bigg(\frac{1+bB^2}{B^2}\bigg)^{\frac{1}{2(n-1)}} \;\;.
\end{equation}
This means that the elliptic PDE that we have constructed is essentially unique: we would gain nothing if we replaced $f$ by $B(f)$, and the additional freedom from replacing the function $K$ that we chose in the preceding section by a more general function $A$ is quite limited in view of \eqref{strictA}. The solvability properties of the resulting PDE are essentially the same for all possible $A$; our $K$ is just the simplest choice.
\medskip\\
These considerations generalise mutatis mutandis to the case of arbitrary given distributions $V,H$ on an arbitrary manifold $M$, and arbitrary metrics of index $q$.
\medskip\\
Let us now return to our general initial question: How unique is our approach to the construction of an elliptic equation?
\smallskip\\
One could ask the same question for the \emph{Riemannian} prescribed scalar curvature problem. Are there alternatives to an approach via conformal deformations? The method from \cite{KazdanWarner2} (which uses results from \mbox{J.-P.} Bourguignon, and A.\ E.\ Fischer and J.\ E.\ Marsden; cf.\ \cite{Bourguignon1975}, \cite{FischerMarsden1975}) does not employ conformal deformations, nor even a PDE for real-valued functions on $M$; instead it works directly on the space of Riemannian metrics on $M$. This method cannot be applied in the pseudo-Riemannian case since the relevant fourth-order differential operator is not elliptic then. When we want to use a PDE for real-valued functions, there seems to be no alternative to conformal deformation. This is not just an ellipticity issue; there is simply no other global deformation of a metric which works for every manifold. Of course, one can always deform a metric inside a given chart in many ways, but that does not help to realise a prescribed function globally as a scalar curvature. On special manifolds, there are also other global parameters one can get a grip on: on a \emph{parallelisable} $n$-manifold for instance, one can control globally each of the $n(n+1)/2$ real-valued functions which define the metric. But on an arbitrary $n$-manifold, there is in general only one function (out of the $n(n+1)/2$ local ones) which one can influence \emph{globally} and thus employ for a proof: the conformal factor. Luckily, this unique possible deformation yields an \emph{elliptic} equation.
\smallskip\\
In the pseudo-Riemannian case, the situation is better. On every manifold $M$ which admits an index-$q$ metric at all, there exist at least \emph{two} controllable parameters since the tangent bundle splits off a vector bundle $V$ of rank $q$. We can thus decompose any given index-$q$ metric $g$ on $M$ which makes $V$ timelike, into a metric on $V$ and a metric on the $g$-orthogonal complement of $V$; these metrics can be rescaled independently.
\smallskip\\
It is precisely this advantage of the pseudo-Riemannian situation which we used to construct an elliptic PDE, by suitable mixing of the two parameters. Analogously to the Riemannian problem, these are the only parameters which one can get a grip on globally on every manifold.
\medskip\\
In this sense, our approach seems to be unique: we use the only two parameters which work on every manifold, and we mix them in an essentially unique way in order to get a uniformly elliptic equation.

\subsubsection{How to prove that the elliptic equation has solutions?}

There exist several standard techniques for elliptic equations which we can now apply in order to solve the prescribed scalar curvature problem. Apart from special cases, the best one in our situation seems to be the method of sub- and supersolutions. This is not obvious from looking at Equation \eqref{babyequation}, since the case we considered in the previous section is just one of those special cases where the sub-/supersolution method is \emph{not} optimal. (This special case will be solved in Section \ref{SIXwhatever}.)
\smallskip\\
Namely, in the situation of the previous section, our given distributions $V,H$ were \emph{integrable}. In the general situation, where $V$ and $H$ are not necessarily integrable, our elliptic PDE (cf.\ Theorem \ref{THEEQUATION}) contains additional summands. One might think naively that these terms make the equation even more complicated. But actually, these terms make it much \emph{easier} to find sub- and supersolutions, especially if they vanish \emph{nowhere} on $M$.
\smallskip\\
It turns out that the set where these terms vanish does not depend on the Riemannian background metric $g$ which occurs in the equation; instead, their zeroes (but not the terms themselves) are determined by the integrability properties of the distributions $V$ and $H$. That's why Chapter \ref{FIVE} is concerned with the construction of distributions (in an arbitrary given homotopy class of distributions) which are \emph{nowhere integrable} in a suitable sense. If $V$ and $H$ have this property, then we can choose constant functions as our sub- and supersolutions. This is how we solve the homotopy class version of the prescribed scalar curvature problem in many cases.
\smallskip\\
The difficulty in the case of Lorentzian metrics is that line distributions are always integrable. So we might still be able to make the $(n-1)$-plane distribution $H$ nowhere integrable and hence get a subsolution of our elliptic equation, but we have to find a supersolution in a different way than before. This is still easy if our prescribed function $s$ is everywhere positive. Otherwise, we have to work harder: then we must find a background metric $g$ with certain properties.
\smallskip\\
The prescribed scalar curvature problem for Lorentzian metrics on $2$-manifolds is even more different, since in this case, $V$ \emph{and} $H$ are line distributions and thus integrable. The sub-/supersolution method does not work there, so we apply other techniques; cf.\ Chapter \ref{SEVEN} for details.

\subsubsection{Perspectives}

The technique explained above --- to construct an elliptic one-parameter family in a multi-parameter family of equations --- might apply to other underdetermined partial differential problems. In fact, although I do not know examples where it had been used before, I would be surprised if this idea was entirely new.
\smallskip\\
Finally, I would like to emphasise that I do not regard our elliptic equation (cf.\ Theorem \ref{THEEQUATION}) as an artificial construction or just a trick which solves the problem it is supposed to solve. Despite its complicated appearance, it is just as natural an equation as, say, the equation \eqref{conformalequation} which describes the behaviour of scalar curvature under conformal change; and the uniform ellipticity requirement singles it out from all equations describing the behaviour of scalar curvature under modifications of the metric, making it just as unique as the conformal equation \eqref{conformalequation} is in the Riemannian special case. Therefore it deserves, in my opinion, to be studied as carefully as the conformal equation has been studied in Riemannian geometry.


\section{Summary of the results proved in this thesis} \label{ONETHREE}

\begin{theorem} \label{INTROmain1}
Let $M$ be a compact $n$-manifold, let $q\in\set{3,\dots,n-3}$, and let $s\in C^\infty(M,\R)$. Then every connected component of $\Metr_q(M)$ contains a metric with scalar curvature $s$.
\smallskip\\
Moreover, let $V$ be a $q$-plane distribution on $M$, let $H$ be an $(n-q)$-plane distribution which is complementary to $V$, let $\mathscr{V}\subseteq\Distr_q(M)$ be a $C^0$-neighbourhood of $V$, and let $\mathscr{H}\subseteq\Distr_{n-q}(M)$ be a $C^0$-neighbourhood of $H$. Then there is a pseudo-Riemannian metric of index $q$ on $M$ with scalar curvature $s$ which makes some element of $\mathscr{V}$ timelike and makes some element of $\mathscr{H}$ spacelike.
\end{theorem}
\Proof
Theorem \ref{main1}.
\end{proof}

\begin{theorem}[metrics of index $2$ in dimension $\geq5$, no restriction on $s$] \label{INTROmain2a}
Let $n\geq5$, let $M$ be a compact $n$-manifold, let $s\in C^\infty(M,\R)$, let $V$ be an orientable $2$-plane distribution on $M$, let $H$ be an $(n-2)$-plane distribution which is complementary to $V$ and admits a nowhere vanishing section. Let $\mathscr{V}\subseteq\Distr_2(M)$ be a $C^0$-neighbourhood of $V$, and let $\mathscr{H}\subseteq\Distr_{n-2}(M)$ be a $C^0$-neighbourhood of $H$. Then there is a pseudo-Riemannian metric of index $2$ on $M$ with scalar curvature $s$ which makes some element of $\mathscr{V}$ timelike and makes some element of $\mathscr{H}$ spacelike. There is a pseudo-Riemannian metric of index $n-2$ on $M$ with scalar curvature $s$ which makes some element of $\mathscr{V}$ spacelike and makes some element of $\mathscr{H}$ timelike.
\end{theorem}
\Proof
Theorem \ref{main2a}.
\end{proof}

\begin{theorem}[metrics of index $2$ in dimension $4$, no restriction on $s$] \label{INTROmain2b}
Let $M$ be a compact $4$-manifold, let $s\in C^\infty(M,\R)$, let $V,H$ be complementary $2$-plane distributions on $M$ which are trivial as vector bundles (so in particular $M$ is parallelisable). Let $\mathscr{V}\subseteq\Distr_2(M)$ be a $C^0$-neighbourhood of $V$, and let $\mathscr{H}\subseteq\Distr_2(M)$ be a $C^0$-neighbourhood of $H$. Then there is a pseudo-Riemannian metric of index $2$ on $M$ with scalar curvature $s$ which makes some element of $\mathscr{V}$ timelike and makes some element of $\mathscr{H}$ spacelike.
\end{theorem}
\Proof
Theorem \ref{main2b}.
\end{proof}

\begin{theorem}[metrics of index $1$ or $2$ in dimension $\geq5$, everywhere positive $s$] \label{INTROmain3a}
Let $n\geq5$, let $M$ be a compact $n$-manifold, let $q\in\set{1,\dots,n-3}$, let $s\in C^\infty(M,\R)$ be everywhere positive. Then every connected component of $\Metr_q(M)$ contains a metric with scalar curvature $s$.
\smallskip\\
Moreover, let $V$ be a $q$-plane distribution on $M$, let $H$ be an $(n-q)$-plane distribution which is complementary to $V$, and let $\mathscr{H}\subseteq\Distr_{n-q}(M)$ be a $C^0$-neighbourhood of $H$. Then there is a pseudo-Riemannian metric of index $q$ on $M$ with scalar curvature $s$ which makes $V$ timelike and makes some element of $\mathscr{H}$ spacelike.
\end{theorem}
\Proof
Theorem \ref{main3a}.
\end{proof}

\begin{theorem}[Lorentzian metrics in dimension $4$, everywhere positive $s$] \label{INTROmain3b}
Let $M$ be a compact connected orientable $4$-manifold which either has nonempty boundary, or is closed with $\sigma_M\not\equiv 2\mod4$.\footnote{Cf.\ Notation \ref{intersectionform}.} Let $s\in C^\infty(M,\R)$ be everywhere positive. Then every connected component of $\Metr_1(M)$ which consists of time-orientable metrics contains a metric with scalar curvature $s$.
\smallskip\\
Moreover, let $V$ be an orientable line distribution on $M$, let $H$ be a $3$-plane distribution which is complementary to $V$, and let $\mathscr{H}\subseteq\Distr_3(M)$ be a $C^0$-neighbourhood of $H$. Then there is a Lorentzian metric on $M$ with scalar curvature $s$ which makes $V$ timelike and makes some element of $\mathscr{H}$ spacelike.
\end{theorem}
\Proof
Theorem \ref{main3b}.
\end{proof}

\begin{theorem}[metrics of index $2$ in dimension $4$, everywhere positive $s$] \label{INTROmain3d}
Let $M$ be a compact $4$-manifold, let $s\in C^\infty(M,\R)$ be everywhere positive, let $V$ be a $2$-plane distribution on $M$ which admits a nowhere vanishing section, let $H$ be an orientable $2$-plane distribution which is complementary to $V$, and let $\mathscr{H}\subseteq\Distr_2(M)$ be a $C^0$-neighbourhood of $H$. Then there is a pseudo-Riemannian metric of index $2$ on $M$ with scalar curvature $s$ which makes $V$ timelike and makes some element of $\mathscr{H}$ spacelike.
\end{theorem}
\Proof
Theorem \ref{main3d}.
\end{proof}

\begin{theorem}[metrics of index $1$ or $2$ in dimension $\geq5$, somewhere positive $s$] \label{INTROmain4a}
Let $n\geq5$, let $M$ be a compact connected $n$-manifold, let $q\in\set{1,2}$, and let $s\in C^\infty(M,\R)$ be somewhere positive. Then every connected component of $\Metr_q(M)$ contains a metric with scalar curvature $s$.
\end{theorem}
\Proof
Theorem \ref{main4a}.
\end{proof}

\begin{theorem}[Lorentzian metrics in dimension $4$, somewhere positive $s$] \label{INTROmain4b}
Let $M$ be a compact connected orientable $4$-manifold which either has nonempty boundary, or is closed with $\sigma_M\not\equiv 2\mod4$. Let $s\in C^\infty(M,\R)$ be somewhere positive. Then every connected component of $\Metr_1(M)$ which consists of time-orientable metrics contains a metric with scalar curvature $s$.
\end{theorem}
\Proof
Theorem \ref{main4b}.
\end{proof}

\begin{theorem}[Lorentzian metrics in dimension $3$, somewhere positive $s$] \label{INTROmain4c}
Let $M$ be a compact connected orientable $3$-manifold, let $s\in C^\infty(M,\R)$ be somewhere positive. Then every connected component of $\Metr_1(M)$ contains a metric with scalar curvature $s$.
\end{theorem}
\Proof
Theorem \ref{main4c}. Cf.\ also Theorem \ref{main3c}.
\end{proof}

\begin{theorem}[metrics of index $2$ in dimension $4$, somewhere positive $s$] \label{INTROmain4d}
Let $M$ be a compact connected $4$-manifold, let $s\in C^\infty(M,\R)$ be somewhere positive, let $\mathscr{C}$ be a connected component of $\Metr_2(M)$ consisting of space-orientable metrics, such that the elements of $\TDC(\mathscr{C})$ admit a nowhere vanishing section. Then $\mathscr{C}$ contains a metric with scalar curvature $s$.
\end{theorem}
\Proof
Theorem \ref{main4d}.
\end{proof}

\begin{theorem}[real-analytic versions] \label{INTROmain1ra}
All the theorems above hold in the real-analytic category: If the manifold $M$ is equipped with a real-analytic atlas, if the function $s$ is real-analytic, and if in the theorems \ref{INTROmain3a}, \ref{INTROmain3b}, \ref{INTROmain3d} the distribution $V$ is real-analytic, then we can find a metric with the stated properties which is real-analytic.
\end{theorem}
\Proof
Theorem \ref{main1ra}, Remark \ref{ranaremark}.
\end{proof}

\begin{theorem}[product manifolds, somewhere negative $s$] \label{INTROmain5}
Let $q,m\in\N_{\geq1}$, let $B$ be a closed connected $q$-manifold, let $N$ be a closed connected $m$-manifold, let $M$ be the product manifold $B\times N$. Let $V,H$ denote the first-factor resp.\ second-factor distribution on $M$, and let $s\in C^\infty(M,\R)$ be a somewhere negative function. If $m\geq3$, or if $m=2$ and $\chi(N)<0$, or if $q=m=2$ and $\chi(B)>0$, or if $q\geq3$ and $B$ admits a Riemannian metric with positive scalar curvature, then there exists a pseudo-Riemannian metric $h$ of index $q$ on $M$ with scalar curvature $s$, and there exists a diffeomorphism $\varphi\in\Diff^0(M)$ such that $\varphi^\ast(V)$ is $h$-timelike and $\varphi^\ast(H)$ is $h$-spacelike.
\end{theorem}

\subsection{The $2$-dimensional case}

\begin{theorem}[Lorentzian metrics on closed $2$-manifolds] \label{introclosedtwo}
Let $M$ be either the $2$-dimensional torus or the Klein bottle (these are the only closed nonempty connected $2$-manifolds which admit a Lorentzian metric), and let $s\in C^\infty(M,\R)$. Then there is a Lorentzian metric on $M$ with scalar curvature $s$ if and only if $s$ is either identically zero or changes its sign (i.e.\ is positive somewhere and negative somewhere else).
\end{theorem}
\Proof
Theorem \ref{LorentzGaussBonnet}, Theorem \ref{closedtwotheorem}.
\end{proof}

\begin{theorem} \label{introopentwo}
Let $M$ be a connected compact $2$-manifold with nonempty boundary, let $s\in C^\infty(M,\R)$. Then there is a Lorentzian metric on $M$ with scalar curvature $s$.
\end{theorem}
\Proof
Theorem \ref{twoopenmain}.
\end{proof}

\subsection{Byproducts}

Some results of this thesis might be interesting independent of the prescribed scalar curvature problem. We will not list them here but point out where they can be found.
\medskip\\
The definition of the \emph{twistedness} $\Twist_V\in C^\infty(M\ot \Lambda^2(V^\ast)\otimes\bot V)$ of a distribution $V$ on a manifold $M$ is reviewed in Subsection \ref{twistsection}. In Chapter 5, we prove several theorems which say that under certain conditions on the dimension $n$ of $M$ and the rank $q$ of $V$, the set of distributions which are \emph{everywhere twisted} (i.e.\ whose twistedness, being a section in a vector bundle, vanishes nowhere on $M$) is dense in the space $\Distr_q(M)$ of all $q$-plane distributions on $M$ with respect to the fine $C^0$-topology. For instance, the condition $(n-q)(q-2)\geq2$ suffices. Cf.\ the theorems \ref{FIVEMAIN}, \ref{FIVEMAINfourdim}, \ref{FIVEMAINtwoplane}.
\smallskip\\
Moreover, we show that under the same conditions \emph{integrable} distributions on a \emph{compact} manifold can be approximated by everywhere twisted distributions even with respect to the $C^\infty$-topology; cf.\ Section \ref{twistCinfty}.
\medskip\\
It might also be interesting to have a look at Theorem \ref{realanalyticRiemann1} in Appendix \ref{D3}, which deals with the real-analytic Riemannian prescribed scalar curvature problem, quite analogous to our Theorem \ref{INTROmain1ra}. The statement follows easily from standard theorems and is therefore probably well-known, but I have never seen it in the literature.

\subsection{A conjecture}

The following conjecture does of course not count as a \emph{result} of the thesis, but it indicates quite clearly what remains to be done on the prescribed scalar curvature problem. The message is that there are probably no obstructions in dimension $\geq4$ to the existence of Lorentzian metrics with certain functions as scalar curvatures; i.e., there is probably no Lorentzian analogue to the well-known obstructions against Riemannian metrics with positive scalar curvature. We will discuss the conjecture briefly in Section \ref{SIXesc}; a strategy for a proof is outlined there. The letters \emph{esc} suggest that ``\emph{e}very function is a \emph{s}calar \emph{c}urvature'' (which might remind you of the notion of \emph{psc manifolds}, i.e.\ manifolds which admit a Riemannian metric of \emph{p}ositive \emph{s}calar \emph{c}urvature).

\begin{escconjecture} \label{esc}
Let $M$ be a compact manifold of dimension $\geq4$, and let $s\in C^\infty(M,\R)$. Then every connected component of the space of Lorentzian metrics on $M$ contains a metric with scalar curvature $s$.
\end{escconjecture}


\section{Overview of the further chapters} \label{ONEFOUR}

The rest of the thesis is arranged in a logical order: each chapter uses only facts from previous chapters.

\begin{description}
\item[{\sc Chapter 2.}]
Here we introduce all the tensor fields and functions which play a role in the rest of the thesis, and state some of their basic properties.

\item[{\sc Chapter 3.}]
We define three simple ways to modify a given pseudo-Riemannian metric, and we calculate how these manipulations affect the scalar curvature. The tedious computational details can be skipped by everyone who trusts the results.

\item[{\sc Chapter 4.}]
Using the results of Chapter 3, we construct the elliptic equation which forms the core of our approach to the pseudo-Riemannian prescribed scalar curvature problem. Again, the details consist of nothing but high school algebra and can safely be skipped.

\item[{\sc Chapter 5.}]
As an application of M.~Gromov's convex integration technique, we prove several existence theorems for $q$-plane distributions which are \emph{everywhere twisted} (cf.\ Subsection \ref{twistsection} for a definition of this notion), i.e.\ nonintegrable in a weak sense. The contents of this chapter are purely differential-topological and can be read independently of the rest of the thesis. They are the key ingredients to most of the results in Chapter \ref{SIX}.

\item[{\sc Chapter 6.}]
We prove the main existence theorems for solutions of our elliptic PDE, employing the method of sub- and supersolutions and the technique that Kazdan and Warner developed for their solution of the Riemannian prescribed scalar curvature problem.

\item[{\sc Chapter 7.}]
Here we deal with the case of Lorentzian metrics on $2$-manifolds. This chapter is independent of the results of Chapters 5 and 6. The analytic techniques here are the Kazdan/Warner method (on closed manifolds) and direct methods in the calculus of variations (on open manifolds).
\end{description}

A fast way to read the thesis would be as follows: Skim the definitions in Chapter 2, then jump to Theorem \ref{THEEQUATION} at the end of Chapter 4 where the elliptic PDE is written down. If you are interested in Lorentz surfaces, continue with Chapter 7. If you are interested in manifolds of dimension $\geq3$ instead, continue with Chapter 6. There you will see what everywhere twisted distributions are good for. After that, you can go back to Chapter 5 where these distributions are constructed.

\bigskip
We conclude this overview with a list of all tools and theorems that we apply in the present work (without proving them here); it might give an idea how elementary its results are:

\bigskip
\textsc{Standard tools (theorems applicable to a wide range of problems).}
\begin{itemize}
\item
Standard tools from analysis, in particular elliptic PDE theory: the method of sub- and supersolutions, elliptic regularity, solving a PDE by variational methods, the $L^p$ approximation theorem of Kazdan and Warner (cf.\ Theorem \ref{KWapproximation1} in Appendix \ref{D2}) and the corresponding implicit function technique.
\item
The most basic result of obstruction theory (cf.\ Appendix \ref{obstructionappendix}). (needed for the problem in dimensions $\geq5$)
\item
M.\ Gromov's convex integration technique for solving partial differential relations on first-order jet bundles (cf.\ Appendix \ref{hprinciple}). (needed for the problem in dimensions $\geq4$)
\end{itemize}

\textsc{More specialised theorems.}
\begin{itemize}
\item
The classical existence theorem for almost-complex structures on $4$-manifolds, due to F.~Hirzebruch and H.~Hopf (cf.\ \cite{HirzebruchHopf}). (only needed for the $4$-dimensional case)
\item
The classical existence theorem for contact structures on orientable $3$-manifolds, due to J.~Martinet, R.~Lutz, and Y.~Eliashberg (cf.\ \cite{Geigessurvey}, Theorem 3.1, and \cite{Eliashberg1989}). (only needed for the $3$-dimensional case)
\item
P.~Percell's existence theorem for Riemannian metrics with parallel vector fields on compact manifolds with nonempty boundary (cf.\ \cite{Percell1981}). (only needed for the $2$-dimensional case)
\item
W.~Thurston's existence theorem for codimension-$1$ foliations (cf.\ \cite{Thurston1976}). (occurs only in the discussion of the esc Conjecture in Section \ref{SIXesc})
\end{itemize}


\chapter{Basic differential geometry} \label{TWO}

Given a pseudo-Riemannian manifold $(M,g)$, a distribution $V$ on $M$, and perhaps a function $f$ on $M$, we can form several other functions on $M$. They appear in the elliptic equation which we will construct in Chapter \ref{FOUR} and which we will use in Chapters \ref{SIX} and \ref{SEVEN} to solve the pseudo-Riemannian prescribed scalar curvature problem. The aim of the present chapter is to define all these functions and to state their basic properties, in particular the expressions which describe them with respect to certain local orthonormal frames of the tangent bundle $TM$.
\smallskip\\
As far as I know, there exists no established name or notation for most of the objects we are going to define. So if they have been introduced somewhere else after all, our names and notations will most likely differ.
\smallskip\\
To be specific, we will define the following functions on $M$:
\begin{quote}
$\eval{\divergence^V_g}{df}_{g,\bot V}$, \;\;$\eval{\divergence^V_g}{\divergence^V_g}_{g,\bot V}$, \;\;$\eval{df}{df}_{g,V}$, \;\;$\laplace^V_{g,V}(f)$
\smallskip\\
$\sigma_{g,V}$, \;\;$\tau_{g,V}$, \;\;$\scal^{V,V}_g$, \;\;$\qual^V_g$, \;\;$\xi_{g,V}$, \;\;$\chi_{g,V}$ \;\;.
\end{quote}

\section{Preliminaries} \label{TWOONE}

This section collects some well-known facts, which are repeated here mainly for the reader's convenience.

\subsection{Linear algebra}

We review some definitions and elementary facts about the linear algebra of semi-Riemannian metrics --- i.e.\ about those aspects of semi-Riemannian metrics which can be explained on the level of vector spaces, without derivatives or curvature ---, and we introduce a convenient language for the definitions in the next section.

\begin{definition}
A \emph{semi-Riemannian vector bundle} [\emph{of index $q$}] is a vector bundle $E\to M$ together with a section $g\in C^\infty(M\ot\Sym(E))$ such that for all $x\in M$, the symmetric bilinear form $g_x$ on the fibre $E_x$ is nondegenerate [and has index $q$].
\end{definition}

\begin{notation} \label{musical}
Let $(E\to M,\,g)$ be a semi-Riemannian vector bundle. We use the following notation: $\flat_g \colon E \to E^\ast$ is the vector bundle isomorphism induced by $g$; i.e., $\flat_g(v)(w) = g(v,w)$ for all $x\in M$ and $v,w\in E_x$. The vector bundle isomorphism $\sharp_g\colon E^\ast\to E$ is the inverse of $\flat_g$. For any sub vector bundle $U$ of $E$, $i_U\colon U\to E$ is the inclusion, and $i_U^\ast\colon E^\ast \to U^\ast$ is its dual bundle morphism. We denote the $g$-induced metric on $E^\ast$ by $\eval{.}{.}_g$; i.e., $\eval{\alpha}{\beta}_g \define g(\sharp_g\alpha,\sharp_g\beta)$ for all $\alpha,\beta$ which lie in the same fibre of $E^\ast$. If $g$ is a Riemannian metric, we define $\abs{\alpha}_g\define \sqrt{\eval{\alpha}{\alpha}_g}$.
\smallskip\\
The \emph{$g$-orthogonal bundle of $U$}, denoted by $\bot_gU$ (cf.\ Definition \ref{maximally}), is the sub vector bundle $\ker(i_U^\ast\compose\flat_g)$ of $E\to M$; it consists of all vectors $w\in E$ such that $g(w,u)=0$ holds for all $u\in U$ in the fibre of $w$.
\end{notation}

\begin{fact} \label{nondegfacts}
Let $(E\to M,\,g)$ be a semi-Riemannian vector bundle, and let $U$ be a sub vector bundle of $E$. Recall (cf.\ e.g.\ \cite{ONeill}, Lemma 2.22, 2.23) that $\rank(U) +\rank(\bot_g U) = \rank(E)$ and $\bot_g\bot_gU = U$, and that the following statements are equivalent:
\begin{enumerate}
\item
The restriction of $g$ to $U$ is nondegenerate\footnote{We call a section $g\in C^\infty(M\ot\Sym(E))$ \emph{nondegenerate} [resp.\ \emph{positive/negative definite}] if it has this property pointwise.}.
\item
The restriction of $g$ to $\bot_gU$ is nondegenerate.
\item
$U$ and $\bot_gU$ are complementary subbundles of $E$, i.e., $E$ is the internal direct sum of $U$ and $\bot_gU$.
\end{enumerate}
A sufficient condition for $g\restrict U$ to be nondegenerate is of course that $g\restrict U$ is positive or negative definite.
\end{fact}

\begin{definition}[$g$-good] \label{ggood}
Let $(E\to M,\,g)$ be a semi-Riemannian vector bundle. We call a sub vector bundle $U$ of $E$ \emph{$g$-good} if and only if the restriction of $g$ to $U$ is nondegenerate.
\end{definition}

\begin{remark} \label{normalbundle}
Let $(E\to M,\,g)$ be a semi-Riemannian vector bundle, and let $U$ be a sub vector bundle of $E$. We will sometimes denote the bundle $\bot_gU$ simply by $\bot U$ when this is not likely to cause confusion. When we use --- later in the thesis --- the notation $\bot U$ in a context where no metric is specified, then we mean the quotient bundle $E/U$. If a metric $g$ on $E$ is given and $U$ is $g$-good, then there is of course a canonical vector bundle isomorphism $\bot_g U \to E/U$, namely the composition of the inclusion $\bot_g U\to E$ and the projection $E\to E/U$.
\end{remark}

\begin{notation} \label{projdef}
Let $(E\to M,\,g)$ be a semi-Riemannian vector bundle, and let $U$ be a $g$-good subbundle of $E$. Since $E = U\oplus\bot_gU$, there is a well-defined orthogonal projection map $E\to U$, given by $w\mapsto u$ whenever $w = u+v$ with $u\in U$, $v\in\bot_gU$. We denote this map by $\pr^U_g$.
\end{notation}

\begin{remark} \label{commdiag}
Let $(E\to M,\,g)$ be a semi-Riemannian vector bundle, and let $U$ be a $g$-good subbundle of $E$. Then the following diagrams commute:
\vspace{2mm} \[
\setlength{\arraycolsep}{1cm} \begin{array}{cc}
\setlength{\arraycolsep}{1cm} \begin{array}{cc}
\Rnode{p00}{E} & \Rnode{p01}{U} \\[1cm]
\Rnode{p10}{E^\ast} & \Rnode{p11}{U^\ast} \end{array}
\psset{arrows=->,nodesep=5pt,linewidth=0.3pt} \everypsbox{\scriptstyle} \ncline{p00}{p01}\Aput{\pr^U_g} \ncline{p10}{p11}\Aput{i_U^\ast} \ncline{p00}{p10}\Bput{\flat_g}\Aput{\cong} \ncline{p01}{p11}\Aput{\flat_{g\restrict U}}\Bput{\cong} %
&\begin{array}{cc} \setlength{\arraycolsep}{1cm} \Rnode{q00}{E^\ast\otimes E^\ast} & \Rnode{q01}{U^\ast\otimes U^\ast} \\[1cm]
\Rnode{q10}{E^\ast\otimes E} & \Rnode{q11}{U^\ast\otimes U} \end{array}
\psset{arrows=->,nodesep=5pt,linewidth=0.3pt} \everypsbox{\scriptstyle} \ncline{q00}{q01}\Aput{i_U^\ast\otimes i_U^\ast} \ncline{q10}{q11}\Aput{i_U^\ast\otimes\pr^U_g} \ncline{q00}{q10}\Bput{\id\otimes\sharp_g}\Aput{\cong} \ncline{q01}{q11}\Aput{\id\otimes\sharp_{g\restrict U}}\Bput{\cong} \end{array}
\]
\end{remark}
\Proof
Every $w\in E$ has a unique decomposition $w=u+v$, where $u\in U$ and $v\in\bot_gU$. This yields $i_U^\ast(\flat_g(w))(u') = (\flat_g(w)\compose i_U)(u') = g(w,i_U(u')) = g(u,u')+g(v,u') = g(u,u')$ for all $u'\in U$, and, on the other hand, $\flat_{g\restrict U}(\pr^U_g(w))(u') = \flat_{g\restrict U}(u)(u') = g(u,u')$. This proves the commutativity of the diagram on the left. The commutativity of the other diagram follows now, since $\sharp$ is the inverse of $\flat$.
\end{proof}

\begin{definition}[contractions]
Let $M$ be a manifold, let $k\in\N_{\geq2}$, and let $V_1,\dots,V_k$ be vector bundles over $M$. We consider the tensor bundle $V_1^\ast\otimes\ldots\otimes V_k^\ast$ over $M$, and a section $T$ in this tensor bundle. Since we are in a finite-dimensional situation, we can and will identify\footnote{Observe that the canonical isomorphism $\iota\colon V\to V^{\ast\ast}$ (which is defined by $\iota(v)(\lambda)=\lambda(v)$ for all $(v,\lambda)\in V\times V^\ast$) is equal to $\flat_{g^\ast}\compose\flat_g$ when $g$ is a semi-Riemannian metric on $V$ and $g^\ast=\eval{.}{.}_g$ is the induced metric on $V^\ast$.} $V_i^{\ast\ast}$ with $V_i$, and elements of $V_1^\ast\otimes\ldots\otimes V_k^\ast$ with multilinear forms $V_1\times\ldots\times V_k\to\R$.
\smallskip\\
Let $(E\to M,\,g)$ be a semi-Riemannian vector bundle. Assume that $V_i=V_j=E$ for some $i,j\in\set{1,\dots,k}$ with $i\neq j$. Then we define the \emph{$g$-contraction of $T$ in the $i$th and $j$th index} as usual, namely as the (pointwise) contraction in the $i$th and $j$th index of the tensor field $T^{\sharp i}\in C^\infty(M\;\ot\; V_1^\ast\otimes\ldots\otimes V_{i-1}^\ast\otimes E\otimes V_{i+1}^\ast\otimes\ldots\otimes V_k^\ast)$ given by
\[
T^{\sharp i}(v_1,\dots,v_{i-1},\lambda,v_{i+1},\dots,v_k) = T(v_1,\dots,v_{i-1},\sharp_g(\lambda),v_{i+1},\dots,v_k) \;\;.
\]
In other words, $T^{\sharp i}$ arises from $T$ by pulling up the $i$th index via the metric $g$. Note that the contraction in the $i$th and $j$th index of the tensor field $T^{\sharp i}$ is equal to the contraction in the $i$th and $j$th index of $T^{\sharp j}$.
\medskip\\
Now we generalise this definition: Let $m\in\N_{\geq1}$, let $(E_1,g_1),\dots,(E_m,g_m)$ be semi-Riemannian vector bundles over $M$, let $i(1),j(1),\dots,i(m),j(m)$ be distinct elements of $\set{1,\dots,k}$, and assume that $V_{i(\mu)}=V_{j(\mu)}=E_\mu$ for all $\mu\in\set{1,\dots,m}$. (The definition above was the case $m=1$.)
\smallskip\\
Then we define the \emph{$(g_1,\dots,g_m)$-contraction of $T$ in the $i(1)$th and $j(1)$th, the $i(2)$th and $j(2)$th, \dots, and the $i(m)$th and $j(m)$th index} by the following straightforward generalisation of the definition above: We perform $m$ consecutive contractions, where the $\mu$th contraction is done with respect to the metric $g_\mu$, and in those indices which were the indices $i(\mu)$ and $j(\mu)$ before we deleted several indices by the earlier $\mu-1$ contractions. Note that this $m$-fold contraction is independent of the order in which we perform the $m$ contractions.
\end{definition}

\begin{remark}
Let $(E\to M,\,g)$ be a semi-Riemannian vector bundle, let $U$ be a $g$-good subbundle of $E$, and let $T$ be a section in $E^\ast\otimes E^\ast$. We define the section $T_U$ in $U^\ast\otimes U^\ast$ to be the restriction of $T$ to $U$. The commutativity of the diagram on the right in \ref{commdiag} shows that we can compute the $(g\restrict U)$-contraction of $T_U$ in two ways: Either we restrict $T$ to a section in $U^\ast\otimes U^\ast$, pull up one index by the metric $g\restrict U$, and contract; or we pull up one index of $T$ by the metric $g$, turn the resulting section in $E^\ast\otimes E$ into a section in $U^\ast\otimes U$ by restriction and $g$-projection, and contract.
\end{remark}

\begin{remark} \label{multilinear}
Recall the following standard method for defining tensor fields: Let $E\to M$ be a vector bundle, and let $\tilde{T}\colon C^\infty(M\ot E)\to C^\infty(M,\R)$ be a $C^\infty(M,\R)$-linear map. Then there is a unique section $T$ in the vector bundle $E^\ast\to M$ such that for all $x\in M$, $v\in E_x$ and $\tilde{v}\in C^\infty(M\ot E)$ with $\tilde{v}(x)=v$, we have $T(v)=\tilde{T}(\tilde{v})(x)$. The statement generalises to $C^\infty(M,\R)$-multilinear maps $\tilde{T}\colon C^\infty(M\ot E_1)\times\ldots\times C^\infty(M\ot E_k) \to C^\infty(M,\R)$ and tensor fields $T\in C^\infty(M\ot E_1^\ast\otimes\dots\otimes E_k^\ast)$.
\end{remark}

\subsection{Adapted orthonormal frames and ON Christoffel symbols}

We will define some tensor fields and functions in the next section, and we will calculate them with respect to suitable local orthonormal frames. Here we introduce resp. review the relevant notions.

\begin{definition} \label{ONframedef}
An \emph{orthonormal frame} (abbreviated: \emph{ON frame}) of a semi-Riemannian vector bundle $(E\to M,g)$ of rank $k$ is a $k$-frame $(e_1,\dots,e_k)$ of $E\to M$ (i.e., the values of the sections $e_i\in C^\infty(M\ot E)$ in each point $x\in M$ form a basis of $E_x$) such that $g(e_i,e_j)=0$ for all $i,j\in\set{1,\dots,k}$ with $i\neq j$, and, for each $i\in\set{1,\dots,k}$, the function $\eps_i\define g(e_i,e_i)$ is either the constant $1$ or the constant $-1$. (Clearly, if $g$ has index $q$, then there are exactly $q$ elements $i\in\set{1,\dots,k}$ with $g(e_i,e_i)=-1$, and thus exactly $k-q$ elements $i\in\set{1,\dots,k}$ with $g(e_i,e_i)=1$.) Whenever we use the symbol $\eps_i$ in a context where an ON frame is given, it denotes the constant $g(e_i,e_i)$.
\end{definition}

\begin{remark} \label{ONcontract}
Let $(E\to M,\;g)$ be a semi-Riemannian vector bundle of rank $k$, and let $(e_1,\dots,e_k)$ be an orthonormal frame for it. Let $T$ be a section in $E^\ast\otimes E^\ast$. Then the $g$-contraction of $T$ is the real-valued function $\sum_{i=1}^k\eps_iT(e_i,e_i)$ on $M$. (By definition, the $g$-contraction of $T$ is the function $\sum_{i=1}^kT^{\sharp2}(e_i,e_i^\ast) = \sum_{i=1}^kT(e_i,\sharp_g(e_i^\ast))$, where $(e_1^\ast,\dots,e_k^\ast)$ is the algebraic dual frame of $(e_1,\dots,e_k)$, which is given by $e_i^\ast(e_j)=\delta_{ij}$ for all $i,j\in\set{1,\dots,k}$. Since $g(\eps_ie_i,e_j)=\delta_{ij}=e_i^\ast(e_j)$ for all $i,j$, we have $\flat_g(\eps_ie_i)=e_i^\ast$ and thus $\eps_ie_i=\sharp_g(e_i^\ast)$, which proves the claimed formula.)
\smallskip\\
This fact generalises in an obvious way to the situation where $E_1,\dots,E_m$ are vector bundles over $M$ with $E=E_i=E_j$ for distinct $i,j\in\set{1,\dots,m}$, where $g$ is a semi-Riemannian metric on $E$, and $T$ is a section in $E_1^\ast\otimes\dots\otimes E_m^\ast$, which we contract in the $i$th and $j$th index.
\end{remark}

\begin{definition}[$U$-adapted ON frame]
Let $(M,g)$ be a semi-Riemannian $n$-manifold, let $q\in\set{0,\dots,n}$, and let $U$ be a $g$-good $q$-plane distribution. A $g$-orthonormal frame $(e_1,\dots,e_n)$ of the tangent bundle $TM$ is called \emph{$U$-adapted} if and only if there is a subset $\varrho \subseteq \set{1,\dots,n}$ of cardinality $q$ such that $e_j(x)\in U_x$ for all $i\in\varrho$ and $x\in M$.
\smallskip\\
In other words, $(e_1,\dots,e_n)$ is $U$-adapted if and only if there is a subset $\varrho \subseteq \set{1,\dots,n}$ such that $U = \spann\set{e_i \suchthat i\in\varrho}$. (This set $\varrho$ is obviously unique, and, moreover, satisfies $\bot U = \spann\set{e_i \suchthat i\in\set{1,\dots,n}\without\varrho}$.)
\end{definition}

\begin{notation} \label{colonnotation}
Let $U$ be a $g$-good distribution on a semi-Riemannian $n$-manifold $(M,g)$, let $(e_1,\dots,e_n)$ be a $U$-adapted orthonormal frame, and let $\varrho$ be the unique subset of $\set{1,\dots,n}$ with $U = \spann\set{e_i \suchthat i\in\varrho}$. We use the notation $i:U$ as an abbreviation for the expression $i\in\varrho$, and we use the notation $i:\bot U$ for the expression $i\in\set{1,\dots,n}\without\varrho$. Expressions like $i,j:U$ have to be interpreted as $i:U, \;j:U$. A summation sign $\sum_i$ denotes a sum from $1$ to $n$. (Sometimes we use indices $0,\dots,n-1$ instead of $1,\dots,n$. The notations generalise to this situation in an obvious way.)
\end{notation}

\begin{proposition}
Let $(E\to M,g)$ be a semi-Riemannian vector bundle. Then there exists, for every $x\in M$, an open neighbourhood $N$ of $x$ such that $E\restrict N$ admits a $g$-orthonormal frame.
\end{proposition}
\Proof
We prove this by induction over the rank $k$ of $E$. For $k=0$, the statement of the proposition is true. If $k\geq1$ and $x\in M$, then there is a vector $v_x\in E_x$ with $g(v_x,v_x)\neq0$. We choose any section $v\in C^\infty(N_0\ot E)$ on some neighbourhood $N_0$ of $x$ such that $v(x)=v_x$. By continuity, there is a connected neighbourhood $N_1\subseteq N_0$ of $x$ such that $g(v,v)$ vanishes nowhere on $N_1$. The vector field $e_k\in C^\infty(N_1\ot E)$ defined by $e_k \define v/\sqrt{\abs{g(v,v)}}$ satisfies $\abs{g(e_k,e_k)}=1$. Let $E_{k-1}\to N_1$ be the $g$-orthogonal bundle of the line bundle $\R e_k$. Since $(E_{k-1},g\restrict E_{k-1})$ is a semi-Riemannian vector bundle of rank $k-1$ (note that $\R e_k$ is $g$-good since $g\restrict \R e_k$ is positive or negative definite), we can apply the induction hypothesis and find a neighbourhood $N\subseteq N_1$ of $x$ which admits a $(g\restrict E_{k-1})$-orthonormal frame $(e_1,\dots,e_{k-1})$. Because $(e_1,\dots,e_{k-1},e_k\restrict N)$ is a $g$-orthonormal frame of $E\restrict N$, the induction proof is complete.
\end{proof}

\begin{corollary} \label{adaptedframe}
Let $(M,g)$ be a semi-Riemannian $n$-manifold equipped with a $g$-good distribution $U$. Then there exists, for every $x\in M$, an open neighbourhood $N$ of $x$ such that $TN$ admits a $g$-orthonormal frame which is $U$-adapted.
\end{corollary}
\Proof
Since $U$ is $g$-good, the semi-Riemannian vector bundle $(TM,g)$ is the internal $g$-orthogonal direct sum of the semi-Riemannian vector bundles $(U,g\restrict U)$ and $(\bot_gU,g\restrict\bot_gU)$. Both of them admit local orthonormal frames by the preceding proposition. If $k=\rank(U)$, and $(e_1,\dots,e_k)$ and $(e_{k+1},\dots,e_n)$ are such local ON frames, then the restriction of $(e_1,\dots,e_n)$ to the intersection of the neighbourhoods over which these two frames are defined is a $U$-adapted local ON frame of $(TM,g)$.
\end{proof}

The preceding corollary is the main reason why we use local ON frames instead of local coordinates for the calculations in Chapter \ref{THREE}: $U$-adapted ON frames do always exist locally, $U$-adapted local coordinates exist (if and) only if the distribution $U$ is integrable. In the situation of Subsection \ref{ONETWOONE} of the introduction, we used local coordinates because they made the computation very simple. In the general situation, however, they would instead complicate it considerably.
\medskip\\
From now on, we will frequently use the standard term \emph{[$U$-adapted] local orthonormal frame} for an ON frame which is defined on an open subset of the manifold under consideration [and is adapted to the restriction of $U$ to this subset].

\begin{remark} \label{projformula}
Let $(M,g)$ be a semi-Riemannian $n$-manifold, and let $(e_1,\dots,e_n)$ be an orthonormal frame on $(M,g)$. Then every $v\in TM$ satisfies $v = \sum_{i=1}^n\eps_ig(v,e_i)e_i$. If $U$ is a $g$-good distribution on $M$ and $(e_1,\dots,e_n)$ is $U$-adapted, then $\pr^U_g(v) = \sum_{i:U}\eps_ig(v,e_i)e_i$.
\end{remark}

\begin{definition}[ON Christoffel symbols] \label{ONChristoffel}
Let $(M,g)$ be a semi-Riemannian $n$-manifold, and let $(e_1,\dots,e_n)$ be a local orthonormal frame on $(M,g)$. We define the \emph{orthonormal Christoffel symbols} by
\[
\boxed{\Gamma^k_{ij} \define g(\nabla_{e_i}e_j,e_k)}
\]
for $i,j,k\in\set{1,\dots,n}$, where $\nabla$ denotes the \LeviCivita\ connection of $g$. (The metric $g$ will always be clear from the context and is therefore suppressed in our notation.)
\begin{center} \fbox{\parbox{160mm}{
\emph{In Chapters \ref{TWO}--\ref{FOUR}, Christoffel symbols appear only in this form, that is, only with respect to local orthonormal frames; we will not use Christoffel symbols with respect to coordinate systems.}}}
\end{center}
(Coordinate Christoffel symbols --- denoted by $\varGamma$ instead of $\Gamma$ --- occurred in one computation in the Second Step of Subsection \ref{ONETWOONE}. The rest of the thesis is free of any Christoffel symbols.)
\end{definition}

\begin{remark} \label{Christoffel}
Let $(M,g)$ be a semi-Riemannian $n$-manifold, and let $(e_1,\dots,e_n)$ be a local orthonormal frame on $(M,g)$. Since $0 = \partial_{e_i}g(e_j,e_k) = g(\nabla_{e_i}e_j,e_k)+g(e_j,\nabla_{e_i}e_k)$, the equation
\[
\boxed{\Gamma^k_{ij} = -\Gamma^j_{ik}}
\]
holds for all $i,j,k\in\set{1,\dots,n}$. In particular,
\[
\boxed{\Gamma^k_{ik} = 0} \;\;.
\]
\emph{We will use these rules from now on without further mention.}
\smallskip\\
The Koszul formula (cf.\ e.g.\ \cite{ONeill}, p.~61) yields for all $i,j,k\in\set{1,\dots,n}$
\[
2\Gamma^k_{ij} = g([e_i,e_j],e_k) +g([e_k,e_i],e_j) +g([e_k,e_j],e_i) \;\;.
\]
\end{remark}


\section{Tensor fields defined by a distribution on a manifold} \label{TWOTWO}

\emph{Throughout this section, we consider a semi-Riemannian manifold $(M,g)$ of dimension $n\in\N$.}
\medskip\\
Now we are ready to introduce some tensor fields which can be defined using $(M,g)$ and some $g$-good distribution on $M$. Since we are dealing with scalar curvature (as opposed to Riemann or Ricci curvature), we are mostly interested not in the tensor fields themselves but in their total contractions, which are functions on $M$ --- namely those listed at the beginning of this chapter.
\smallskip\\
The philosophy behind most of the following definitions is that usual tensor fields or functions on $M$ --- e.g.\ the scalar curvature, or the Laplacian of $f$, where $f\in C^\infty(M,\R)$ --- can be split up into several summands if we take the decomposition $TM = U\oplus\bot U$ with respect to a given $g$-good distribution $U$ into account.
\smallskip\\
$\nabla$ denotes the (\LeviCivita) covariant derivative with respect to the metric $g$ in this section. In the rest of the thesis, we will sometimes have to specify the metric $g$ explicitly in our notation; we will then write $\nabla^{(g)}$ instead of just $\nabla$.

\subsection{Laplacians and divergences}

\begin{definition}[$\divergence^U_g(X)$] \label{divUXdef}
Let $U$ be a $g$-good distribution. Let $X$ be a vector field on $M$. We define $\nabla^{[U]}X$ as the section in $U^\ast\otimes U$ which is the image of the section $\nabla X$ in $T^\ast M \otimes TM$ under the map $T^\ast M \otimes TM \to U^\ast\otimes U$ induced by $i_U^\ast \colon T^\ast M \to U^\ast$ and $\pr^U_g \colon TM \to U$.
\smallskip\\
We define $\divergence^U_g(X)\in C^\infty(M,\R)$ to be the contraction of $\nabla^{[U]}X\in C^\infty(M\ot\,U^\ast\otimes U)$.
\end{definition}

\begin{frameformulae} \label{divUformulae}
Let $U$ be a $g$-good distribution. Let $(e_1,\dots,e_n)$ be a local $U$-adapted orthonormal frame. For every function $h\in C^\infty(M,\R)$ and every vector field $X$ on $M$, we have by \ref{projformula}:
\[ \begin{split}
\divergence^U_g(X) &= \sum_{i:U}\eps_ig(\nabla_{e_i}X,e_i) \;\;,\\
\divergence^U_g(hX) &= \sum_{i:U}\eps_ig(\nabla_{e_i}(hX),e_i) = h\sum_{i:U}\eps_ig(\nabla_{e_i}X,e_i)+\sum_{i:U}\eps_idh(e_i)g(X,e_i)\\
&= h\divergence^U_g(X) +dh(\pr^U_g(X)) \;\;.
\end{split} \]
The equation $\divergence^U_g(hX) = h\divergence^U_g(X) +dh(\pr^U_g(X))$ holds globally on $M$, since local $U$-adapted ON frames exist around each point in $M$.
\end{frameformulae}

\begin{definition}[$\eval{\alpha}{\beta}_{g,U}$] \label{scalproddef}
Let $U$ be a $g$-good distribution. For any two sections $\alpha,\beta$ in $U^\ast$, we define the function $\eval{\alpha}{\beta}_{g,U}\in C^\infty(M,\R)$ to be the scalar product $\eval{\alpha}{\beta}_{g\restrict U}$ (cf.\ \ref{musical}), i.e.\ the $(g\restrict U)$-contraction of $\alpha\otimes\beta\in C^\infty(M\ot\,U^\ast\otimes U^\ast)$. If $g\restrict U$ is positive definite, we sometimes use the notation $\abs{\alpha}_g \define \sqrt{\eval{\alpha}{\alpha}_{g,U}}$.
\par
If $\alpha$ is a section in $T^\ast M$ or in $U^\ast$, and if $\beta$ is a section in $T^\ast M$ or in $U^\ast$, we define $\eval{\alpha}{\beta}_{g,U}$ to be $\eval{\alpha\restrict U}{\beta\restrict U}_{g,U}$, where $\alpha\restrict U = i_U^\ast(\alpha)$ and $\beta\restrict U = i_U^\ast(\beta)$ are the restrictions of $\alpha$ resp.\ $\beta$ to sections in $U^\ast$.
\end{definition}

\begin{examples} \label{scalprodexamples}
Let $U$ be a $g$-good distribution. Note that we can interpret $\divergence^{\bot U}_g$ as a section in the vector bundle $U^\ast$, since $\divergence^{\bot U}_g(fv) = f\divergence^{\bot U}_g(v) +df(\pr^{\bot U}_g(v)) = f\divergence^{\bot U}_g(v)$ for every $f\in C^\infty(M,\R)$ and $v\in C^\infty(M\ot U)$; cf.\ Remark \ref{multilinear}.
\smallskip\\
Hence $\eval{\divergence^{\bot U}_g}{\divergence^{\bot U}_g}_{g,U}$ is a well-defined smooth function on $M$. For all functions $f,h\in C^\infty(M,\R)$, we can also consider the functions $\eval{df}{dh}_{g,U}$ and $\eval{\divergence^{\bot U}_g}{df}_{g,U}$.
\end{examples}

\begin{frameformulae} \label{moreONformulae}
Let $U$ be a $g$-good distribution, let $(e_1,\dots,e_n)$ be a local $U$-adapted orthonormal frame, and let $f,h\in C^\infty(M,\R)$. Then (cf.\ \ref{ONChristoffel} for the definition of the functions $\Gamma^k_{ij}$):
\[ \begin{split}
\eval{\divergence^{\bot U}_g}{\divergence^{\bot U}_g}_{g,U} &= \sum_{i:U}\eps_i\divergence^{\bot U}_g(e_i)\divergence^{\bot U}_g(e_i) = \sum_{i:U}\sum_{j,k:\bot U}\eps_i\eps_j\eps_k\Gamma^j_{ji}\Gamma^k_{ki} \;\;,\\
\eval{\divergence^{\bot U}_g}{df}_{g,U} &= \sum_{i:U}\eps_i\divergence^{\bot U}_g(e_i)df(e_i) = \sum_{i:U}\sum_{j:\bot U}\eps_i\eps_j\Gamma^j_{ji}df(e_i) \;\;,\\
\eval{df}{dh}_{g,U} &= \sum_{i:U}\eps_idf(e_i)dh(e_i)\;\;.
\end{split} \]
\end{frameformulae}

\begin{definition}[$\laplace^U_{g,W}(f)$, $\laplace^U_g(f)$] \label{laplacedef}
Let $U$ be a $g$-good distribution, and let $f\in C^\infty(M,\R)$. For $W\in\set{U,\bot U}$, we define the function
\[
\laplace^U_{g,W}(f) \define \divergence^U_g(\pr^W_g(\grad_g(f))) \;\;;
\]
and we define the function
\[
\laplace^U_g(f) \define \divergence^U_g(\grad_g(f)) = \laplace^U_{g,U}(f) +\laplace^U_{g,\bot U}(f) \;\;.
\]
\end{definition}

\begin{remark}
Let $U$ be a $g$-good distribution. Let $f\in C^\infty(M,\R)$. The Hessian $\Hess_g(f)$ is a symmetric $(0,2)$-tensor field on $M$, i.e.\ a symmetric section in the vector bundle $T^\ast M \otimes T^\ast M$. We define $\Hess^U_g(f)$ to be the restriction of $\Hess_g(f)$ to a section in the vector bundle $U^\ast\otimes U^\ast$. Then $\laplace^U_g(f)$ is the $(g\restrict U)$-contraction of $\Hess^U_g(f)$, as the following formulae prove.
\end{remark}

\begin{frameformulae} \label{laplacianformulae}
Let $U$ be a $g$-good distribution, let $(e_1,\dots,e_n)$ be a local $U$-adapted ON frame, and let $f\in C^\infty(M,\R)$. Since the Hessian of $f$ is given by $\Hess_g(f)(v,w) = \partial_v\partial_wf -df(\nabla_vw)$, the $(g\restrict U)$-contraction of $\Hess^U_g(f)$ is
\[ \begin{split}
\sum_{i:U}\eps_i\Hess_g(e_i,e_i) &= \sum_{i:U}\eps_i\partial_{e_i}\partial_{e_i}f -\sum_{i:U}\sum_k\eps_i\eps_kg(\nabla_{e_i}e_i,e_k)df(e_k)\\
&= \sum_{i:U}\eps_i\partial_{e_i}\partial_{e_i}f +\sum_k\eps_k\divergence^U_g(e_k)df(e_k) \;\;.
\end{split} \]
For $W\in\set{U,\bot U}$, we compute (using \ref{divUformulae})
\[ \begin{split}
\laplace^U_{g,W}(f) &= \divergence^U_g\Big(\sum_{i:W}\eps_idf(e_i)e_i\Big) = \sum_{i:W}\eps_idf(e_i)\divergence^U_g(e_i) +\sum_{i:W}\eps_id(df(e_i))(\pr^U_g(e_i))\\
&= \sum_{i:W\cap U}\eps_i\partial_{e_i}\partial_{e_i}f +\sum_{i:W}\eps_i\divergence^U_g(e_i)df(e_i)
\end{split} \]
(note that $\sum_{i:\bot U\cap U}$ is the empty sum since $\bot U\cap U$ is a $0$-plane distribution). Hence
\[ \begin{split}
\laplace^U_g(f) &= \laplace^U_{g,U}(f) +\laplace^U_{g,\bot U}(f) = \sum_{i:U}\eps_i\partial_{e_i}\partial_{e_i}f +\sum_{i:U}\eps_i\divergence^U_g(e_i)df(e_i) +\sum_{i:\bot U}\eps_i\divergence^U_g(e_i)df(e_i)\\
&= \sum_{i:U}\eps_i\partial_{e_i}\partial_{e_i}f +\sum_i\eps_i\divergence^U_g(e_i)df(e_i) \;\;,
\end{split} \]
which is equal to the $(g\restrict U)$-contraction of $\Hess^U_g(f)$.
\end{frameformulae}

\begin{remark} \label{laplaciansplit}
Let $U$ be a $g$-good distribution, let $f\in C^\infty(M,\R)$. By calculating contractions with respect to local $U$-adapted $g$-orthonormal frames, we get
\[
\laplace_g(f) = \laplace^U_g(f) +\laplace^{\bot U}_g(f) \;\;.
\]
\end{remark}

\begin{remark} \label{laplacedivergence}
Let $U$ be a $g$-good distribution, let $f\in C^\infty(M,\R)$. Then
\[
\laplace^U_{g,\bot U}(f) = \eval{\divergence^U_g}{df}_{g,\bot U} \;\;,
\]
since for every local $U$-adapted ON frame $(e_1,\dots,e_n)$, we have by \ref{moreONformulae} and \ref{laplacianformulae}
\[
\laplace^U_{g,\bot U}(f) = \sum_{i:\bot U}\eps_i\divergence^U_g(e_i)df(e_i) = \eval{\divergence^U_g}{df}_{g,\bot U} \;\;.
\]
In particular, $\laplace^U_{g,\bot U}(f) \colon C^\infty(M,\R) \to C^\infty(M,\R)$ is a \emph{first}-order differential operator. This is not really surprising, but nonetheless crucial for the present work: The only second-order terms which occur anywhere in the computations in Chapter \ref{THREE} have the form $\laplace^U_{g,W}(f)$, where $V$ is a given distribution $V$ and $U,W\in\set{V,\bot V}$. So we have to deal with only two second-order terms $\laplace^V_{g,V}(f)$ and $\laplace^{\bot V}_{g,\bot V}(f)$, which are analogous to the functions $\partial_t\partial_tf$ and $\laplace_{g_N}(f)$ that we have seen in the discussion of a special case in Subsection \ref{ONETWOONE}. Hence the general situation is, up to first and zeroth order terms, not very different from the special case: we can construct an elliptic equation by the same recipe that we have discussed in Subsection \ref{ONETWOONE}. If e.g.\ $\laplace^{\bot V}_{g,V}(f)$ included second-order terms which vanished only in the special case, we would have a problem; but that is not the case.
\end{remark}

\begin{remark} \label{laplacecompose}
Let $U$ be a $g$-good distribution, let $W\in\set{U,\bot U}$, let $f\in C^\infty(M,\R)$, let $K\in C^\infty(\R,\R)$. Then
\[
\laplace^U_{g,W}(K\compose f) = (K'\compose f)\laplace^U_{g,W}(f) +(K''\compose f)\eval{df}{df}_{g,U\cap W} \;\;.
\]
\end{remark}
\Proof
With respect to every local $U$-adapted ON frame $(e_1,\dots,e_n)$, we have by \ref{divUformulae}
\[ \begin{split}
&\laplace^U_{g,W}(K\compose f) = \divergence^U_g\Big(\sum_{i:W}\eps_id(K\compose f)(e_i)e_i\Big) = \divergence^U_g\Big((K'\compose f)\sum_{i:W}\eps_idf(e_i)e_i\Big)\\
&\mspace{30mu}= (K'\compose f)\divergence^U_g\Big(\sum_{i:W}\eps_idf(e_i)e_i\Big) +d(K'\compose f)\Big(\sum_{i:U\cap W}\eps_idf(e_i)e_i\Big)\\
&\mspace{30mu}= (K'\compose f)\laplace^U_{g,W}(f) +(K''\compose f)\sum_{i:U\cap W}\eps_idf(e_i)df(e_i) = (K'\compose f)\laplace^U_{g,W}(f) +(K''\compose f)\eval{df}{df}_{g,U\cap W} \;.\hfill\qedhere
\end{split} \]
\end{proof}

\begin{remark} \label{laplacemultiply}
Let $U$ be a $g$-good distribution, let $W\in\set{U,\bot U}$, let $f_0,f_1\in C^\infty(M,\R)$. Then
\[
\laplace^U_{g,W}(f_0f_1) = f_0\laplace^U_{g,W}(f_1) +f_1\laplace^U_{g,W}(f_0) +2\eval{df_0}{df_1}_{g,W\cap U} \;\;.
\]
\end{remark}
\Proof
With respect to every local $U$-adapted ON frame $(e_1,\dots,e_n)$, we have by \ref{laplacianformulae}
\[ \begin{split}
\laplace^U_{g,W}(f_0f_1) &= \sum_{i:W\cap U}\eps_i\partial_{e_i}\partial_{e_i}(f_0f_1) +\sum_{i:W}\eps_i\divergence^U_g(e_i)d(f_0f_1)(e_i)\\
&= \sum_{i:W\cap U}\eps_i\partial_{e_i}\Big(f_0df_1(e_i) +f_1df_0(e_i)\Big) +\sum_{i:W}\eps_i\divergence^U_g(e_i)\Big(f_0df_1(e_i) +f_1df_0(e_i)\Big)\\
&= \sum_{i:W\cap U}\eps_idf_0(e_i)df_1(e_i) +f_0\sum_{i:W\cap U}\eps_i\partial_{e_i}\partial_{e_i}f_1 +\sum_{i:W\cap U}\eps_idf_1(e_i)df_0(e_i) +f_1\sum_{i:W\cap U}\eps_i\partial_{e_i}\partial_{e_i}f_0\\
&\mspace{20mu}+f_0\sum_{i:W}\eps_i\divergence^U_g(e_i)df_1(e_i) +f_1\sum_{i:W}\eps_i\divergence^U_g(e_i)df_0(e_i)\\
&= f_0\laplace^U_{g,W}(f_1) +f_1\laplace^U_{g,W}(f_0) +2\eval{df_0}{df_1}_{g,W\cap U} \;\;.\hfill\qedhere
\end{split} \]
\end{proof}

\subsection{The functions $\sigma_{g,U}$, $\tau_{g,U}$}

Given a distribution $U$ and the $1$-jet of the metric $g$ (i.e.\ the first derivatives of $g$), we define two functions $\sigma_{g,U}$, $\tau_{g,U}$. Their relation to the integrability properties of $U$ will be discussed in Subsection \ref{twistsection}.

\begin{definition} \label{secondff}
Let $U$ be a $g$-good distribution. We define a section $\mathscr{T}^U_g$ in $T^\ast M \otimes U^\ast \otimes (\bot U)^\ast$ by $(w,u,v)\mapsto g(\nabla_wu,v)$. To see that this section is well-defined, cf.\ Remark \ref{multilinear} and note that the defining map is $C^\infty(M,\R)$-linear in the second argument since $g(\nabla_w(fu),v) = f\,g(\nabla_wu,v) +df(w)g(u,v) = f\,g(\nabla_wu,v)$ for every section $u$ in $U$ and every function $f$ on $M$.
\end{definition}

\begin{remark} \label{alternatingindices}
We have $\mathscr{T}^{\bot U}_g(w,v,u) = -\mathscr{T}^U_g(w,u,v)$, since $g(\nabla_wu,v)+g(\nabla_wv,u) = \partial_wg(u,v) = 0$.
\end{remark}

\begin{remark}[relation to the second fundamental form]
Let us consider the case in which our $g$-good distribution $U$ is integrable. Then it induces a foliation of $M$, and at each point $x\in M$, the leaf through $x$ determines a submanifold germ $N$ at the point $x$. (In other words, for each $x\in M$ there is a submanifold $N$ of $M$ [of dimension $\rank(U)$] which contains $x$ and satisfies $U_y = T_yN$ for all $y\in N$. This submanifold is uniquely determined up to restriction to a smaller neighbourhood of $x$. Note that the whole leaf through $x$ is in general not a submanifold of $M$.)
\smallskip\\
The second fundamental form $\SecondFF$ of $N$ is a section in the vector bundle $T^\ast N\otimes T^\ast N\otimes \bot(TN)$; i.e., $\SecondFF\in C^\infty(N\ot\;U^\ast\otimes U^\ast\otimes \bot U)$. It is given by $\SecondFF(u,v) = \pr^{\bot U}_g(\nabla_uv)$ for all $u,v\in C^\infty(N\ot U) = C^\infty(N\ot TN)$; recall that $\nabla$ is the covariant derivative on $(M,g)$. We can pull down the upper index of the tensor field $\SecondFF$ via the metric $g$, thereby defining a section $\secondFF\in C^\infty(N\ot\,U^\ast\otimes U^\ast\otimes(\bot U)^\ast)$; i.e., $\secondFF(u,v,w) = g(\nabla_uv,w)$ for all $u,v\in C^\infty(N\ot U)$, $w\in C^\infty(N\ot\bot U)$.
\smallskip\\
The preceding definition of the second fundamental form of submanifold germs yields also a definition of the second fundamental form of a foliation: If $U$ is integrable, then $\secondFF_U\in C^\infty(M\ot U^\ast\otimes U^\ast\otimes(\bot U)^\ast)$ assigns to each point $x\in M$ the second fundamental form in $x$ of (the germ of) the leaf through $x$. \smallskip\\
In this situation, the (pointwise) restriction of the tensor field $\mathscr{T}^U_g\in C^\infty(M\ot\,T^\ast M\otimes U^\ast \otimes (\bot U)^\ast)$ to a section in $U^\ast\otimes U^\ast\otimes(\bot U)^\ast$ is equal to the second fundamental form $\secondFF_U$. In this sense, $\mathscr{T}^U_g$ is a generalisation of the second fundamental form to arbitrary (not necessarily integrable) distributions.
\end{remark}

\begin{definition}[$\sigma_{g,U}, \tau_{g,U}$] \label{defsigmatau}
Let $U$ be a $g$-good distribution, and let $W\in\set{U,\bot U}$. We define the section $\mathscr{D}^U_{g,W}$ in the vector bundle $W^\ast \otimes U^\ast \otimes (\bot U)^\ast \otimes W^\ast\otimes U^\ast\otimes(\bot U)^\ast$ as the (pointwise) restriction of the section $\mathscr{T}^U_g \otimes \mathscr{T}^U_g$ in $T^\ast M\otimes U^\ast\otimes(\bot U)^\ast\otimes T^\ast M\otimes U^\ast\otimes(\bot U)^\ast$.
\smallskip\\
By taking suitable contractions of $\mathscr{D}^U_{g,W}$, we could now define many tensor fields. Because functions is all we need in the following, we restrict ourselves to the consideration of total contractions of $\mathscr{D}^U_{g,W}$.
\smallskip\\
We define the function $\sigma_{g,U}$ as the contraction of $\mathscr{D}^{\bot U}_{g,U} \in C^\infty(M\ot\;U^\ast\otimes(\bot U)^\ast\otimes U^\ast\otimes U^\ast\otimes(\bot U)^\ast\otimes U^\ast)$ in the first and fourth, the second and fifth, and the third and sixth index. The contractions are of course taken with respect to the metrics $g\restrict U$ and $g\restrict(\bot U)$ on the vector bundles $U$ and $\bot U$, respectively.
\smallskip\\
We define the function $\tau_{g,U}$ as the contraction of $\mathscr{D}^{\bot U}_{g,U} \in C^\infty(M\ot\;U^\ast\otimes(\bot U)^\ast\otimes U^\ast\otimes U^\ast\otimes(\bot U)^\ast\otimes U^\ast)$ in the first and sixth, the third and fourth, and the second and fifth index.
\end{definition}

\begin{remark}
In Riemannian geometry, it would be customary to use a notation like $\abs{T}^2_g$ for $\sigma_{g,U}$, where $T$ is the restriction of $\mathscr{T}^{\bot U}_g$ to a section in $U^\ast\otimes(\bot U)^\ast\otimes U^\ast$.
\end{remark}

\begin{frameformulae} \label{sigmatauframe}
Let $U$ be a $g$-good distribution, let $(e_1,\dots,e_n)$ be a local $U$-adapted ON frame. Then
\[ \begin{split}
\sigma_{g,U} &= \sum_{i,k:U}\sum_{j:\bot U}\eps_i\eps_j\eps_k\Gamma^k_{ij}\Gamma^k_{ij} = \sum_{i,k:U}\sum_{j:\bot U}\eps_i\eps_j\eps_k\Gamma^j_{ik}\Gamma^j_{ik} \;\;,\\
\tau_{g,U} &= \sum_{i,k:U}\sum_{j:\bot U}\eps_i\eps_j\eps_k\Gamma^k_{ij}\Gamma^i_{kj} = \sum_{i,k:U}\sum_{j:\bot U}\eps_i\eps_j\eps_k\Gamma^j_{ik}\Gamma^j_{ki} \;\;.
\end{split} \]
\end{frameformulae}

\begin{remark}[why no other total contractions exist]
There are several other possibilities to form total contractions of $\mathscr{D}^{\bot U}_{g,U}$ or $\mathscr{D}^U_{g,U}$, but all of them can be expressed by functions that we know already:
\smallskip\\
Since $\mathscr{D}^U_{g,U}(u_0,v_0,w_0,u_1,v_1,w_1) = \mathscr{D}^{\bot U}_{g,U}(u_0,w_0,v_0,u_1,w_1,v_1)$ (cf.\ Remark \ref{alternatingindices}), considering only total contractions of $\mathscr{D}^{\bot U}_{g,U}$ is no loss of generality. In addition to the contractions which define $\sigma_{g,U}$ and $\tau_{g,U}$, exactly one other contraction of $\mathscr{D}^{\bot U}_{g,U}$ is possible: contraction in the first and third, the fourth and sixth, and the second and fifth index. This contraction is equal to the function $\eval{\divergence^U_g}{\divergence^U_g}_{g,\bot U}$.
\end{remark}

\subsection{Scalar curvatures}

Given a distribution $U$ and the $2$-jet of the metric $g$, we can define more functions than the $\sigma$s and $\tau$s from the preceding subsection. The sum of these functions $\scal^{U,U}_g$, $\scal^{U,\bot U}_g$, $\scal^{\bot U,U}_g$, $\scal^{\bot U,\bot U}_g$ which we are going to introduce now is the scalar curvature of $g$.

\begin{frameformulae}[Riemann curvature] \label{Riemanntensor}
Let $R_g\in C^\infty(M\ot\,T^\ast M\otimes T^\ast M\otimes T^\ast M\otimes T^\ast M)$ be the $(0,4)$-tensor version of the Riemann curvature of $(M,g)$. Our sign convention\footnote{Our sign convention agrees e.g.\ with \cite{KobayashiNomizu1}, \cite{Petersen}, \cite{Klingenberg}, \cite{KolarMichorSlovak}, \cite{Lang}; and it is opposite to e.g.\ \cite{Besse}, \cite{GallotHulinLafontaine}, \cite{ONeill}. However, all conventions agree on the definition of the Ricci curvature and the scalar curvature.} yields for vector fields $u,v,w,z\in C^\infty(M\ot TM)$:
\[
R_g(u,v,w,z) = g\big(\nabla_u\nabla_vw-\nabla_v\nabla_uw-\nabla_{[u,v]}w,z\big) \;\;.
\]
Let $(e_1,\dots,e_n)$ be a local $g$-orthonormal frame of $TM$. Then we have for all $i,j,k,l\in\set{1,\dots,n}$:
\[ \begin{split}
R_g(e_i,e_j,e_k,e_l) &= \partial_{e_i}\Gamma^l_{jk} -\partial_{e_j}\Gamma^l_{ik} +\sum_\mu\eps_\mu\bigg(\Gamma^l_{i\mu}\Gamma^\mu_{jk} -\Gamma^l_{j\mu}\Gamma^\mu_{ik} -(\Gamma^\mu_{ij}-\Gamma^\mu_{ji})\Gamma^l_{\mu k}\bigg) \;\;,\\
\scal_g &= -\bigg(2\sum_i\eps_i\partial_{e_i}\divergence_g(e_i) +\sum_i\eps_i\divergence_g(e_i)^2 +\sum_{i,j,k}\eps_i\eps_j\eps_k\Gamma^k_{ij}\Gamma^k_{ji}\bigg) \;\;.
\end{split} \]
\end{frameformulae}
\renewcommand{\FILL}{\mspace{20mu}}
\Proof
We compute $R_g(e_i,e_j,e_k,e_l)$:
\[ \begin{split}
&g\Big(\nabla_{e_i}\nabla_{e_j}e_k -\nabla_{e_j}\nabla_{e_i}e_k
-\nabla_{[e_i,e_j]}e_k, \;e_l\Big)
\\
&\FILL= g\bigg(\nabla_{e_i}\Big(\sum_\mu\eps_\mu g(\nabla_{e_j}e_k,e_\mu)e_\mu\Big) -\nabla_{e_j}\Big(\sum_\mu\eps_\mu g(\nabla_{e_i}e_k,e_\mu)e_\mu\Big) -\nabla_{\sum_\mu\eps_\mu g([e_i,e_j],e_\mu)e_\mu}e_k, \;e_l\bigg)\\
&\FILL= \sum_\mu\eps_\mu g\Big(\nabla_{e_i}\Big(\Gamma^\mu_{jk}e_\mu\Big),e_l\Big) -\sum_\mu\eps_\mu g\Big(\nabla_{e_j}\Big(\Gamma^\mu_{ik}e_\mu\Big),e_l\Big) -\sum_\mu\eps_\mu g\Big([e_i,e_j],e_\mu\Big)g\Big(\nabla_{e_\mu}e_k,e_l\Big)\\
&\FILL= \sum_\mu\eps_\mu\bigg(\Gamma^\mu_{jk}g(\nabla_{e_i}e_\mu,e_l) +(\partial_{e_i}\Gamma^\mu_{jk})g(e_\mu,e_l) -\Gamma^\mu_{ik}g(\nabla_{e_j}e_\mu,e_l) -(\partial_{e_j}\Gamma^\mu_{ik})g(e_\mu,e_l) -(\Gamma^\mu_{ij}-\Gamma^\mu_{ji})\Gamma^l_{\mu k}\bigg)\\
&\FILL= \partial_{e_i}\Gamma^l_{jk} -\partial_{e_j}\Gamma^l_{ik} +\sum_\mu\eps_\mu\bigg(\Gamma^l_{i\mu}\Gamma^\mu_{jk} -\Gamma^l_{j\mu}\Gamma^\mu_{ik} -(\Gamma^\mu_{ij}-\Gamma^\mu_{ji})\Gamma^l_{\mu k}\bigg) \;\;.
\end{split} \]
Contraction yields (note that $\divergence_g(e_i) = \sum_k\eps_kg(\nabla_{e_k}e_i,e_k) = \sum_k\eps_k\Gamma^k_{ki}$):
\[ \begin{split}
\scal_g &= \sum_{i,j}\eps_i\eps_j R_g(e_i,e_j,e_j,e_i)\\
&= \sum_{i,j}\eps_i\eps_j\bigg(\partial_{e_i}\Gamma^i_{jj}-\partial_{e_j}\Gamma^i_{ij} +\sum_k\eps_k(\Gamma^i_{ik}\Gamma^k_{jj} -\Gamma^i_{jk}\Gamma^k_{ij}) -\sum_k\eps_k(\Gamma^k_{ij}-\Gamma^k_{ji})\Gamma^i_{kj}\bigg)\\
&= -2\sum_i\eps_i\partial_{e_i}\Big(\sum_j\eps_j\Gamma^j_{ji}\Big) -\sum_k\eps_k\Big(\sum_i\eps_i\Gamma^i_{ik}\Big)\Big(\sum_j\eps_j\Gamma^j_{jk}\Big)\\
&\mspace{20mu}+\sum_{i,j,k}\eps_i\eps_j\eps_k\Gamma^k_{ji}\Gamma^k_{ij} -\sum_{i,j,k}\eps_i\eps_j\eps_k\Gamma^j_{ik}\Gamma^j_{ki} -\sum_{i,j,k}\eps_i\eps_j\eps_k\Gamma^i_{jk}\Gamma^i_{kj}\\
&= -\bigg(2\sum_i\eps_i\partial_{e_i}\divergence_g(e_i) +\sum_i\eps_i\divergence_g(e_i)^2 +\sum_{i,j,k}\eps_i\eps_j\eps_k\Gamma^k_{ij}\Gamma^k_{ji}\bigg) \;\;.\hfill\qedhere
\end{split} \]
\end{proof}

\begin{definition}[$\scal^{U,W}_g$] \label{scaldef}
Let $U$ be a $g$-good distribution. We define $\Ric^U_g\in C^\infty(M\ot\,T^\ast M\otimes T^\ast M)$ to be the $(g\restrict U)$-contraction of the section in $U^\ast\otimes T^\ast M\otimes T^\ast M\otimes U^\ast$ which we get as a restriction of the Riemann tensor $R_g$, in the first and fourth index. The tensor field $\Ric^U_g$ is symmetric because of the symmetries of $R_g$. We define the function $\scal^U_g\in C^\infty(M,\R)$ to be its $g$-contraction.
\smallskip\\
For $W\in\set{U,\bot U}$, we define the function $\scal^{U,W}_g\in C^\infty(M,\R)$ to be the $(g\restrict W)$-contraction of the restriction of $\Ric^U_g\in C^\infty(M\ot T^\ast M\otimes T^\ast M)$ to a section in the vector bundle $W^\ast\otimes W^\ast$.
\end{definition}

\begin{remark} \label{scaldecomposition}
Obviously, we have $\Ric_g = \Ric^U_g +\Ric^{\bot U}_g$, hence $\scal_g = \scal^U_g +\scal^{\bot U}_g$, and, moreover, $\scal^U_g = \scal^{U,U}_g +\scal^{U,\bot U}_g$. The well-known symmetry $R_g(u,v,v,u) = R_g(v,u,u,v)$ of the curvature tensor implies $\scal^{U,\bot U}_g = \scal^{\bot U,U}_g$.
\end{remark}

\begin{frameformulae} \label{curvatureON}
Let $U$ be a $g$-good distribution, let $(e_1,\dots,e_n)$ be a local $U$-adapted ON frame, and let $W\in\set{U,\bot U}$. By formula \ref{Riemanntensor}, we have
\[ \begin{split}
\scal^{U,W}_g &= \sum_{i:U}\sum_{k:W}\eps_i\eps_kR_g(e_i,e_k,e_k,e_i)\\
&= \sum_{i:U}\sum_{k:W}\eps_i\eps_k\bigg(\partial_{e_i}\Gamma^i_{kk} -\partial_{e_k}\Gamma^i_{ik} +\sum_\mu\eps_\mu\Big(\Gamma^i_{i\mu}\Gamma^\mu_{kk} -\Gamma^i_{k\mu}\Gamma^\mu_{ik} -(\Gamma^\mu_{ik}-\Gamma^\mu_{ki})\Gamma^i_{\mu k}\Big)\bigg) \;\;,
\end{split} \]
in particular
\[ \label{scalON} \begin{split}
\scal^{U,U}_g &= -\sum_{i,k:U}\eps_i\eps_k\partial_{e_i}\Gamma^k_{ki} -\sum_{i,k:U}\eps_i\eps_k\partial_{e_k}\Gamma^i_{ik} -\sum_{i,k:U}\sum_\mu\eps_i\eps_k\eps_\mu\Gamma^i_{i\mu}\Gamma^k_{k\mu} +\sum_{i,k:U}\sum_\mu\eps_i\eps_k\eps_\mu\Gamma^\mu_{ki}\Gamma^\mu_{ik}\\
&\mspace{20mu} -\sum_{i,k:U}\sum_\mu\eps_i\eps_k\eps_\mu\Gamma^k_{i\mu}\Gamma^k_{\mu i} -\sum_{i,k:U}\sum_\mu\eps_i\eps_k\eps_\mu\Gamma^i_{k\mu}\Gamma^i_{\mu k}\\
&= -2\sum_{j:U}\eps_j\partial_{e_j}\divergence^U_g(e_j) -\sum_j\eps_j\divergence^U_g(e_j)^2 -\sum_{i,j,k:U}\eps_i\eps_j\eps_k\Gamma^k_{ij}\Gamma^k_{ji} +\tau_{g,U} -2\sum_{i,k:U}\sum_{j:\bot U}\eps_i\eps_j\eps_k\Gamma^k_{ij}\Gamma^k_{ji}
\end{split} \]
and \[ \begin{split}
\scal^{U,\bot U}_g &= -\sum_{i:U}\sum_{k:\bot U}\eps_i\eps_k\partial_{e_i}\Gamma^k_{ki} -\sum_{i:U}\sum_{k:\bot U}\eps_i\eps_k\partial_{e_k}\Gamma^i_{ik} -\sum_{i:U}\sum_{k:\bot U}\sum_\mu\eps_i\eps_k\eps_\mu\Gamma^i_{i\mu}\Gamma^k_{k\mu}\\
&\mspace{20mu}+\sum_{i:U}\sum_{k:\bot U}\sum_\mu\eps_i\eps_k\eps_\mu\Gamma^\mu_{ki}\Gamma^\mu_{ik} -\sum_{i:U}\sum_{k:\bot U}\sum_\mu\eps_i\eps_k\eps_\mu\Gamma^k_{i\mu}\Gamma^k_{\mu i} -\sum_{i:U}\sum_{k:\bot U}\sum_\mu\eps_i\eps_k\eps_\mu\Gamma^i_{k\mu}\Gamma^i_{\mu k}\\
&= -\sum_{i:U}\eps_i\partial_{e_i}\divergence^{\bot U}_g(e_i) -\sum_{k:\bot U}\eps_k\partial_{e_k}\divergence^U_g(e_k) -\sum_\mu\eps_\mu\divergence^U_g(e_\mu)\divergence^{\bot U}_g(e_\mu)\\
&\mspace{20mu} +\sum_{i:U}\sum_{j,k:\bot U}\eps_i\eps_j\eps_k\Gamma^k_{ij}\Gamma^k_{ji} +\sum_{i,k:U}\sum_{j:\bot U}\eps_i\eps_j\eps_k\Gamma^k_{ij}\Gamma^k_{ji} -\sum_{i:U}\sum_{j,k:\bot U}\eps_i\eps_j\eps_k\Gamma^k_{ij}\Gamma^k_{ji}\\
&\mspace{20mu} -\sum_{i,j:U}\sum_{k:\bot U}\eps_i\eps_j\eps_k\Gamma^k_{ij}\Gamma^k_{ji} -\sum_{k:U}\sum_{i,j:\bot U}\eps_i\eps_j\eps_k\Gamma^k_{ij}\Gamma^k_{ji} -\sum_{j,k:U}\sum_{i:\bot U}\eps_i\eps_j\eps_k\Gamma^k_{ij}\Gamma^k_{ji}\\
&= -\sum_{j:U}\eps_j\partial_{e_j}\divergence^{\bot U}_g(e_j) -\sum_{j:\bot U}\eps_j\partial_{e_j}\divergence^U_g(e_j) -\sum_j\eps_j\divergence^U_g(e_j)\divergence^{\bot U}_g(e_j) -\tau_{g,U} -\tau_{g,\bot U} \;\;.
\end{split} \]
\end{frameformulae}

\subsection{Qualar curvatures}

Until now, we have decomposed the scalar curvature of $g$ into the sum $\scal^{U,U}_g +\scal^{\bot U,\bot U}_g +2\scal^{U,\bot U}_g$. Now we will decompose the summand $\scal^{U,\bot U}_g$ even further, namely into a sum $-(\qual^U_g +\qual^{\bot U}_g)$. For lack of a better name, I call these summands \emph{qualar curvatures} (where \emph{qualar} is a crude mixture of \emph{quasi} and \emph{scalar}; alternatively, it might remind you of \emph{quarter of scalar}). We start with the definition of the corresponding ``quasi'' version of Riemann curvature:

\begin{definition}[$\Qurv^U_g$] \label{Qurvdef}
Let $U$ be a $g$-good distribution. We define the section $\Qurv^U_g$ in the vector bundle $T^\ast M\otimes T^\ast M\otimes U^\ast\otimes (\bot U)^\ast$ by mapping sections $u,v\in C^\infty(M\ot TM)$, $w\in C^\infty(M\ot U)$, $z\in C^\infty(M\ot\bot U)$ to the function
\[
\Qurv^U_g(u,v,w,z) \define g(\nabla_v\nabla_uw,z) -g(\nabla_{\nabla_vu}w,z) +g\big(\pr^U_g(\nabla_vz),\pr^U_g(\nabla_uw)\big) +g\big(\pr^U_g(\nabla_uz),\pr^U_g(\nabla_vw)\big) \;\;.
\]
In order to prove that we really define a section in $T^\ast M\otimes T^\ast M\otimes U^\ast\otimes (\bot U)^\ast$ in this way, we have to check that our map $C^\infty(M\ot TM)\times C^\infty(M\ot TM)\times C^\infty(M\ot U)\times C^\infty(M\ot\bot U) \to C^\infty(M,\R)$ is $C^\infty(M,\R)$-quadrilinear. $C^\infty(M,\R)$-linearity in the second argument (i.e.\ $v$) is obvious. Note that $g(w,z)=0$ and thus $g(\nabla_\eta w,z)+g(\nabla_\eta z,w)=0$ for all $\eta\in C^\infty(M\ot TM)$. Moreover, $g(\pr^U_g(\eta),\pr^U_g(w)) = g(\eta,w)$. Using this, we compute for every $f\in C^\infty(M,\R)$:
\[ \begin{split}
\Qurv^U_g(fu,v,w,z) &= f\Qurv^U_g(u,v,w,z) +df(v)g(\nabla_uw,z) -df(v)g(\nabla_uw,z) = f\Qurv^U_g(u,v,w,z) \;\;,\\
\Qurv^U_g(u,v,fw,z) &= g\big(\nabla_v(f\nabla_uw),z\big) +g\big(\nabla_v\big((\partial_uf)w\big),z\big)\\
&\mspace{20mu}-g\big(\nabla_{\nabla_vu}(fw),z\big) +g\big(\pr^U_g(\nabla_vz),\pr^U_g(\nabla_u(fw))\big) +g\big(\pr^U_g(\nabla_uz),\pr^U_g(\nabla_v(fw))\big)\\
&= f\,g(\nabla_v\nabla_uw,z) +df(v)g(\nabla_uw,z) +df(u)g(\nabla_vw,z) +(\partial_v\partial_uf)g(w,z)\\
&\mspace{20mu}-f\,g(\nabla_{\nabla_vu}w,z) -df(\nabla_vu)g(w,z) +f\,g\big(\pr^U_g(\nabla_vz),\pr^U_g(\nabla_uw)\big)\\
&\mspace{20mu}+df(u)g\big(\pr^U_g(\nabla_vz),\pr^U_g(w)\big) +f\,g\big(\pr^U_g(\nabla_uz),\pr^U_g(\nabla_vw)\big) +df(v)g\big(\pr^U_g(\nabla_uz),\pr^U_g(w)\big)\\
&= f\Qurv^U_g(u,v,w,z) \;\;,\\
\Qurv^U_g(u,v,w,fz) &= f\Qurv^U_g(u,v,w,z) +df(v)g\big(\pr^U_g(z),\pr^U_g(\nabla_uw)\big) +df(u)g\big(\pr^U_g(z),\pr^U_g(\nabla_vw)\big)\\
&= f\Qurv^U_g(u,v,w,z) \;\;.
\end{split} \]
This completes our verification that $\Qurv^U_g\in C^\infty(M\ot\;T^\ast M\otimes T^\ast M\otimes U^\ast\otimes (\bot U)^\ast)$ is well-defined.
\end{definition}

\begin{remark} \label{remarkantisym}
Let $U$ be a $g$-good distribution. Then we have for all $x\in M$, $u,v\in T_xM$, $w\in U_x$, and $z\in \bot U_x$:
\[
\Qurv^U_g(u,v,w,z) = -\Qurv^{\bot U}_g(u,v,z,w) \;\;.
\]
\end{remark}
\Proof
All $u,v\in C^\infty(M\ot TM)$, $w\in C^\infty(M\ot U)$, $z\in C^\infty(M\ot\bot U)$ satisfy
\[
\Qurv^U_g(u,v,w,z) = \partial_vg(\nabla_uw,z) -g(\nabla_{\nabla_vu}w,z) -g\big(\pr^{\bot U}_g(\nabla_vz),\pr^{\bot U}_g(\nabla_uw)\big) +g\big(\pr^U_g(\nabla_uz),\pr^U_g(\nabla_vw)\big)
\]
and $g(\nabla_\eta w,z) = -g(\nabla_\eta z,w)$ for every $\eta\in C^\infty(M\ot TM)$; hence
\[ \begin{split}
\Qurv^U_g(u,v,w,z) &= \partial_vg(\nabla_uw,z) -g(\nabla_{\nabla_vu}w,z) -g\big(\pr^{\bot U}_g(\nabla_vz),\pr^{\bot U}_g(\nabla_uw)\big) +g\big(\pr^U_g(\nabla_uz),\pr^U_g(\nabla_vw)\big)\\
&= -\partial_vg(\nabla_uz,w) +g(\nabla_{\nabla_vu}z,w) +g\big(\pr^U_g(\nabla_vw),\pr^U_g(\nabla_uz)\big) -g\big(\pr^{\bot U}_g(\nabla_uw),\pr^{\bot U}_g(\nabla_vz)\big)\\
&= -\Qurv^{\bot U}_g(u,v,z,w) \;\;.\qedhere
\end{split} \]
\end{proof}

\begin{remark}[relation to Riemann curvature] \label{RiemannQ}
Let $U$ be a $g$-good distribution. Then we have for all $x\in M$, $u,v\in T_xM$, $w\in U_x$, and $z\in \bot U_x$:
\[
R_g(u,v,w,z) = \Qurv^U_g(v,u,w,z) -\Qurv^U_g(u,v,w,z) \;\;.
\]
\end{remark}
\Proof
\[ \begin{split}
R_g(u,v,w,z) &= g(\nabla_u\nabla_vw,z) -g(\nabla_v\nabla_uw,z) -g(\nabla_{[u,v]}w,z)\\
&= -g(\nabla_v\nabla_uw,z) +g(\nabla_u\nabla_vw,z) +g(\nabla_{\nabla_vu}w,z) -g(\nabla_{\nabla_uv}w,z) -g\big(\pr^U_g(\nabla_vz),\pr^U_g(\nabla_uw)\big)\\
&\mspace{20mu}+g\big(\pr^U_g(\nabla_uz),\pr^U_g(\nabla_vw)\big) -g\big(\pr^U_g(\nabla_uz),\pr^U_g(\nabla_vw)\big) +g\big(\pr^U_g(\nabla_vz),\pr^U_g(\nabla_uw)\big)\\
&= -\Qurv^U_g(u,v,w,z) +\Qurv^U_g(v,u,w,z) \;\;.\qedhere
\end{split} \]
\end{proof}

\begin{definition}[$\qual^U_g$] \label{qualdef}
Let $U$ be a $g$-good distribution. We consider the (pointwise) restriction of $\Qurv^U_g$ to a section in $(\bot U)^\ast\otimes U^\ast\otimes U^\ast\otimes (\bot U)^\ast$, and we define the function $\qual^U_g\in C^\infty(M,\R)$ to be its contraction in the first and fourth, and in the second and third index. We call $\qual^U_g$ the \emph{$U$-qualar curvature} of $(M,g)$.
\end{definition}

\begin{remark}
The functions $\qual^U_g$ and $\qual^{\bot U}_g$ are up to a sign the only total contractions of (pointwise) restrictions of $\Qurv^U_g$ or $\Qurv^{\bot U}_g$ (by restrictions we mean restrictions to multilinear forms on $U$ or $\bot U$). This follows immediately from Remark \ref{remarkantisym}.
\end{remark}

\renewcommand{\FILL}{\mspace{10mu}}
\begin{frameformulae} \label{qON}
Let $U$ be a $g$-good distribution. Let $(e_1,\dots,e_n)$ be a local $U$-adapted orthonormal frame. Then we have (by the first formula in the proof of \ref{remarkantisym}):
\[ \begin{split}
\Qurv^U_g(u,v,w,z) &= \partial_vg(\nabla_uw,z) -\sum_k\eps_kg(\nabla_vu,e_k)g(\nabla_{e_k}w,z)\\
&\mspace{20mu}-\sum_{k:\bot U}\eps_kg(\nabla_vz,e_k)g(\nabla_uw,e_k) +\sum_{k:U}\eps_kg(\nabla_uz,e_k)g(\nabla_vw,e_k) \;\;,
\end{split} \]
\[ \begin{split}
\qual^U_g &= \sum_{i:\bot U}\sum_{j:U}\eps_i\eps_j\partial_{e_j}g(\nabla_{e_i}e_j,e_i) -\sum_{i:\bot U}\sum_{j:U}\sum_k\eps_i\eps_j\eps_kg(\nabla_{e_j}e_i,e_k)g(\nabla_{e_k}e_j,e_i)\\
&\mspace{20mu}-\sum_{i,k:\bot U}\sum_{j:U}\eps_i\eps_j\eps_kg(\nabla_{e_j}e_i,e_k)g(\nabla_{e_i}e_j,e_k) +\sum_{i:\bot U}\sum_{j,k:U}\eps_i\eps_j\eps_kg(\nabla_{e_i}e_i,e_k)g(\nabla_{e_j}e_j,e_k)\\
&= \sum_{i:\bot U}\sum_{j:U}\eps_i\eps_j\partial_{e_j}\Gamma^i_{ij} -\sum_{i,k:\bot U}\sum_{j:U}\eps_i\eps_j\eps_k\Gamma^k_{ji}\Gamma^i_{kj} -\sum_{i:\bot U}\sum_{j,k:U}\eps_i\eps_j\eps_k\Gamma^k_{ji}\Gamma^i_{kj}\\
&\mspace{20mu} -\sum_{i,k:\bot U}\sum_{j:U}\eps_i\eps_j\eps_k\Gamma^k_{ji}\Gamma^k_{ij} +\sum_{i:\bot U}\sum_{j,k:U}\eps_i\eps_j\eps_k\Gamma^k_{ii}\Gamma^k_{jj}\\
&= \sum_{j:U}\eps_j\partial_{e_j}\divergence^{\bot U}_g(e_j) +\sum_{i,k:\bot U}\sum_{j:U}\eps_i\eps_j\eps_k\Gamma^i_{jk}\Gamma^i_{kj} +\sum_{i:\bot U}\sum_{j,k:U}\eps_i\eps_j\eps_k\Gamma^i_{jk}\Gamma^i_{kj}\\
&\mspace{20mu} -\sum_{i,k:\bot U}\sum_{j:U}\eps_i\eps_j\eps_k\Gamma^i_{jk}\Gamma^i_{kj} +\sum_{k:U}\eps_k\divergence^U_g(e_k)\divergence^{\bot U}_g(e_k)\\
&= \sum_{j:U}\eps_j\partial_{e_j}\divergence^{\bot U}_g(e_j) +\sum_{j:U}\eps_j\divergence^U_g(e_j)\divergence^{\bot U}_g(e_j) +\tau_{g,U} \;\;.
\end{split} \]
\end{frameformulae}

\begin{fact} \label{scalqual}
Let $U$ be a $g$-good distribution. Then
\[
\scal^{U,\bot U}_g = \scal^{\bot U,U}_g = -(\qual^U_g +\qual^{\bot U}_g) \;\;.
\]
\end{fact}
\begin{proof}[\textsc{First proof.}]
By \ref{curvatureON} and \ref{qON}, we get for every local $U$-adapted ON frame $(e_1,\dots,e_n)$:
\[ \begin{split}
-(\qual^U_g +\qual^{\bot U}_g) &= -\sum_{j:U}\eps_j\partial_{e_j}\divergence^{\bot U}_g(e_j) -\sum_{j:U}\eps_j\divergence^U_g(e_j)\divergence^{\bot U}_g(e_j) -\tau_{g,U}\\
&\mspace{20mu}-\sum_{j:\bot U}\eps_j\partial_{e_j}\divergence^U_g(e_j) -\sum_{j:\bot U}\eps_j\divergence^{\bot U}_g(e_j)\divergence^U_g(e_j) -\tau_{g,\bot U}\\
&= \scal^{U,\bot U}_g \;\;.\hfill\qedhere
\end{split} \]
\end{proof}
\begin{proof}[\textsc{Second proof.}]
By \ref{remarkantisym} and \ref{RiemannQ}, we get
\[ \begin{split}
\scal^{U,\bot U}_g &= \sum_{i:U}\sum_{j:\bot U}\eps_i\eps_jR_g(e_i,e_j,e_j,e_i)\\
&= \sum_{i:U}\sum_{j:\bot U}\eps_i\eps_jQ^U_g(e_j,e_i,e_j,e_i) -\sum_{i:U}\sum_{j:\bot U}\eps_i\eps_jQ^U_g(e_i,e_j,e_j,e_i)\\
&= -\sum_{i:U}\sum_{j:\bot U}\eps_i\eps_jQ^{\bot U}_g(e_j,e_i,e_i,e_j) -\sum_{i:U}\sum_{j:\bot U}\eps_i\eps_jQ^U_g(e_i,e_j,e_j,e_i)\\
&= -\qual^{\bot U}_g -\qual^U_g \;\;.\qedhere
\end{split} \]
\end{proof}

\subsection{The functions $\xi_{g,V}$ and $\chi_{g,V}$}

\begin{definition}[$\xi_{g,V}$] \label{xidef}
Let $V$ be a $g$-good distribution, and let $H\define\bot V$. We define the function $\xi_{g,V}\in C^\infty(M,\R)$ by
\[
\xi_{g,V} \define \eval{\divergence^H_g}{\divergence^H_g}_{g,V} -\eval{\divergence^V_g}{\divergence^V_g}_{g,H} +\frac{\sigma_{g,H}+\tau_{g,H}}{2} -\frac{\sigma_{g,V}+\tau_{g,V}}{2} -\scal^{V,V}_g +2\qual^V_g \;\;.
\]
\end{definition}

This function deserves a name because it will appear as a coefficient in our elliptic PDE \eqref{THEEQUATION}. It is especially relevant in the Lorentzian case of the prescribed scalar curvature problem.

\begin{definition}[$\chi_{g,V}$] \label{chidef}
Let $V$ be a $g$-good distribution, and let $H\define\bot V$. We define the function $\chi_{g,V}\in C^\infty(M,\R)$ by
\[ \begin{split}
\chi_{g,V} &\define \scal_g +\xi_{g,V} +\frac{\sigma_{g,H}-\tau_{g,H}}{2}\\
&\phantom{\rm :}= \scal^{H,H}_g -2\qual^H_g +\eval{\divergence^H_g}{\divergence^H_g}_{g,V} -\eval{\divergence^V_g}{\divergence^V_g}_{g,H} +\sigma_{g,H} -\frac{\sigma_{g,V}+\tau_{g,V}}{2} \;\;.
\end{split} \]
\end{definition}

In Section \ref{SIXesc}, we will discuss the relation of $\chi_{g,V}$ to the esc Conjecture.

\subsection{The line distribution case}

The functions that we have defined above can be interpreted in a more direct way when $V$ is a line distribution, because local orthonormal frames of $V$ are essentially unique in that case. To make this interpretation explicit, we introduce the following objects.

\begin{definition}[$\partial_V\divergence_g(V)$, $\divergence_g(V)^2$, $\nabla^{(g)}_VV$, $\eps_V$] \label{lineobjects}
Let $V$ be a $g$-good line distribution.
\smallskip\\
We define a function $\partial_V\divergence_g(V)\in C^\infty(M,\R)$ as follows. For every $x\in M$, there exist an open neighbourhood $N$ of $x$ and a section $X$ of $V\restrict N$ with $\abs{g(X,X)}=1$. (In other words, there exists a $1$-tuple $(X)$ which is a local ON frame around $x$ of the line bundle $V$. If the line bundle $V$ is orientable and thus trivial, then there is even a \emph{global} ON frame of $V$.) Moreover, if $X_0$ and $X_1$ are two sections with this property, then their germs at $x$ are either equal or equal up to a sign; i.e., there is a neighbourhood $N'$ of $x$ on which both sections are defined, such that either $X_0\restrict N' = X_1\restrict N'$ or $X_0\restrict N' = -X_1\restrict N'$.
\smallskip\\
We define the value of $\partial_V\divergence_g(V)$ in $x$ to be the value of the function $\partial_X\divergence_g(X)$ in $x$. This yields a globally well-defined function, because $\partial_X\divergence_g(X)$ is equal to $\partial_{-X}\divergence_g(-X)$ and the value of $\partial_X\divergence_g(X)$ in $x$ depends only on the germ of $X$ at $x$.
\smallskip\\
Analogously, we define a function $\divergence_g(V)^2$ by declaring that it be locally the function $\divergence_g(X)^2$, where $(X)$ is again a local ON frame of $V$. Since $\divergence_g(X)^2 = \divergence_g(-X)^2$, we get a globally well-defined function in this way.
\smallskip\\
We define a vector field $\nabla^{(g)}_VV\in C^\infty(M\ot TM)$ by declaring that it be locally the vector field $\nabla_XX$, where $(X)$ is a local ON frame of $V$. This is a globally well-defined vector field since $\nabla_XX = \nabla_{-X}(-X)$.
\smallskip\\
We define $\eps_V\define g(X,X)\in\set{1,-1}$, where $(X)$ is a local ON frame of $V$; in other words, $\eps_V=1$ [resp.\ $\eps_V=-1$] if and only if the line distribution $V$ is spacelike [timelike].
\end{definition}

\begin{remark}[$\partial_V\divergence^V_g(V) = 0 = \divergence^V_g(V)^2$]
Let $V$ be a $g$-good line distribution, let $H$ denote the $g$-orthogonal distribution of $V$. In addition to the functions defined above, we could also define functions $\partial_V\divergence^V_g(V)$, $\partial_V\divergence^H_g(V)$, $\divergence^V_g(V)^2$, $\divergence^H_g(V)^2$, in a completely analogous way: we would just replace $\divergence_g(X)$ by $\divergence^V_g(X)$ resp.\ $\divergence^H_g(X)$ in the definitions. However, there is no point in doing so since $\divergence^V_g(X)=0$, which implies $\divergence^H_g(X) = \divergence_g(X)$ and thus $\partial_V\divergence^V_g(V) = 0$, $\divergence^V_g(V)^2 = 0$, $\partial_V\divergence^H_g(V) = \partial_V\divergence_g(V)$, $\divergence^H_g(V)^2 = \divergence_g(V)^2$.
\smallskip\\
In order to see why $\divergence^V_g(X)=0$, consider a local $V$-adapted ON frame $(e_0,\dots,e_{n-1})$ on $M$ such that, without loss of generality, $X=e_0$. Then $\divergence^V_g(X) = \divergence^V_g(e_0) = \sum_{i:V}\eps_i\Gamma^i_{i0} = \eps_0\Gamma^0_{00} = 0$, as claimed.
\end{remark}

The relation of the objects $\partial_V\divergence_g(V)$, $\divergence_g(V)^2$, $\nabla^{(g)}_VV$ to the functions which we defined for distributions of arbitrary rank is as follows.

\begin{fact} \label{lineformulae}
Let $V$ be a $g$-good line distribution, let $H$ denote the $g$-orthogonal distribution of $V$. Then
\begin{gather*}
\sigma_{g,V} = \tau_{g,V} = \eval{\divergence^V_g}{\divergence^V_g}_{g,H} = g(\nabla_VV,\nabla_VV) \;\;,\\
\eval{\divergence^H_g}{\divergence^H_g}_{g,V} = \eps_V\divergence_g(V)^2 \;\;,\\
\scal^{V,V}_g = 0 \;\;,\\
\qual^V_g = \eps_V\partial_V\divergence_g(V) +\sigma_{g,V} \;\;,\\
\qual^H_g = -\eps_V\divergence_g(\nabla_VV) -\sigma_{g,V} +\tau_{g,H} \;\;,\\
\xi_{g,V} = 2\eps_V\partial_V\divergence_g(V) +\eps_V\divergence_g(V)^2 +\frac{\sigma_{g,H}+\tau_{g,H}}{2} \;\;,
\end{gather*}
where we have used the abbreviation $\nabla_VV\define \nabla^{(g)}_VV$.
\end{fact}
\Proof
Let $(e_0,\dots,e_{n-1})$ be a $V$-adapted local ON frame of $TM$ such that, without loss of generality, $e_0$ is a section in $V$. Then $\nabla_VV = \nabla_{e_0}e_0 = \sum_i\eps_i\Gamma^i_{00}e_i$, and we thus get (by the preceding remark and \ref{moreONformulae}, \ref{sigmatauframe}, \ref{curvatureON}, \ref{qON}):
\begin{gather*}
\sum_{i:H}\eps_i\Gamma^i_{00}\Gamma^i_{00} = \sigma_{g,V} = \tau_{g,V} = \eval{\divergence^V_g}{\divergence^V_g}_{g,H} = g(\nabla_VV,\nabla_VV) \;\;,\\
\eval{\divergence^V_g}{\divergence^V_g}_{g,H} = \eps_0\divergence_g(e_0)^2 = \eps_V\divergence_g(V)^2 \;\;,\\
\scal^{V,V}_g = -2(\eps_0)^2\partial_{e_0}\Gamma^0_{00} -\sum_i\eps_i(\Gamma^0_{0i})^2 -\eps_0(\Gamma^0_{00})^2 +\tau_{g,V} -2\sum_{i:H}\eps_i\Gamma^0_{0i}\Gamma^0_{i0} = -\sum_{i:H}\eps_i(\Gamma^0_{0i})^2 +\tau_{g,V} = 0 \;\;,\\
\qual^V_g = \eps_0\partial_{e_0}\divergence^H_g(e_0) +(\eps_0)^2\Gamma^0_{00}\divergence^H_g(e_0) +\tau_{g,V} = \eps_V\partial_{e_0}\divergence_g(e_0) +\tau_{g,V} = \eps_V\partial_V\divergence_g(V) +\sigma_{g,V} \;\;,\\
\qual^H_g = \sum_{j:H}\eps_j\partial_{e_j}\divergence^V_g(e_j) +\sum_{j:H}\eps_j\divergence^H_g(e_j)\divergence^V_g(e_j) +\tau_{g,H} \;\;,
\end{gather*}
\[ \begin{split}
-\eps_V\divergence_g(\nabla_VV) &= -\eps_V\sum_{i:H}\eps_i\divergence_g(\Gamma^i_{00}e_i) = -\eps_V\sum_{i:H}\eps_i\Gamma^i_{00}\divergence_g(e_i) -\eps_V\sum_{i:H}\eps_i\partial_{e_i}\Gamma^i_{00}\\
&= \sum_{i:H}\eps_i\eps_0\Gamma^0_{0i}\Big(\divergence^V_g(e_i)+\divergence^H_g(e_i)\Big) +\sum_{i:H}\eps_i\eps_0\partial_{e_i}\Gamma^0_{0i}\\
&= \sum_{i:H}\eps_i\divergence^V_g(e_i)\divergence^H_g(e_i) +\sum_{i:H}\eps_i\divergence^V_g(e_i)^2 +\sum_{i:H}\eps_i\partial_{e_i}\divergence^V_g(e_i)\\
&= \qual^H_g -\tau_{g,H} +\sigma_{g,V} \;\;;
\end{split} \]
here we used the fact that $\divergence_g(fX) = f\divergence_g(X) +df(X)$ for every function $f$ and every vector field $X$ (cf.\ \ref{divUformulae} for $U=TM$). Finally,
\[ \begin{split}
\xi_{g,V} &= \eval{\divergence^H_g}{\divergence^H_g}_{g,V} -\eval{\divergence^V_g}{\divergence^V_g}_{g,H} +\frac{\sigma_{g,H}+\tau_{g,H}}{2} -\frac{\sigma_{g,V}+\tau_{g,V}}{2} -\scal^{V,V}_g +2\qual^V_g\\
&= \eps_V\divergence_g(V)^2 +\frac{\sigma_{g,H}+\tau_{g,H}}{2} +2\eps_V\partial_V\divergence_g(V) \;\;.\qedhere
\end{split} \]
\end{proof}


\section{Integrability properties of distributions}

\subsection{$\sigma-\tau$ and the twistedness tensor field} \label{twistsection}

Recall some facts and notions that we used already in Chapter \ref{ONE}: \smallskip\\
Let $M$ be an $n$-manifold. A $q$-plane distribution $H$ on $M$ is called \emph{integrable} if and only if for each $x\in M$ there exist local coordinates $(x_1,\dots,x_n)$ on a neighbourhood $U$ of $x$ such that $H\restrict U$ is the span of the vector fields $\frac{\partial}{\partial x_1},\dots,\frac{\partial}{\partial x_q}$ on $U$. By the Frobenius theorem (cf.\ e.g.\ \cite{Conlon}, Theorem 4.5.5), a distribution $H$ is integrable if and only if the sections in $H$ form a sub Lie algebra of the Lie algebra of all vector fields on $M$, that is, if and only if the Lie bracket of every two sections in $H$ is again a section in $H$. In particular, every line distribution is integrable. There is a natural bijective correspondence between integrable $q$-plane distributions on $M$ and $q$-dimensional foliations of $M$; it maps every foliation to its tangent distribution (cf.\ e.g.\ \cite{Conlon}, \S4.5).
\medskip\\
Moreover, recall the following closely related standard construction. Its name seems to be less standardised; I follow W.\ Thurston (cf.\ \cite{Thurston}, p.~176) in calling it the \emph{twistedness} of a distribution.

\begin{definition}[twistedness] \label{twistednessdef}
Let $M$ be an $n$-manifold, and let $H$ be a $q$-plane distribution on $M$. Denoting the fibrewise projection (which is a vector bundle morphism) $TM\to TM/H \invdef \bot H$ by $\pr^{\bot H}$, we define a section $\Twist_H$ in the vector bundle $\Lambda^2(H^\ast)\otimes\bot H$ via
\[
\Twist_H(v,w) \define \pr^{\bot H}([v,w])
\]
for all $v,w\in C^\infty(M\ot H)$. This yields indeed a well-defined section in $\Lambda^2(H^\ast)\otimes\bot H$ because $(v,w)\mapsto \pr^{\bot H}([v,w])$ is an alternating $C^\infty(M,\R)$-bilinear map $C^\infty(M\ot H)\times C^\infty(M\ot H)\to C^\infty(M\ot\bot H)$ (cf.\ Remark \ref{multilinear}): for every $f\in C^\infty(M,\R)$, we have $\pr^{\bot H}([fv,w]) = \pr^{\bot H}(f[v,w]-df(w)v) = f\pr^{\bot H}([v,w])$ and $[v,w]=-[w,v]$.
\smallskip\\
When $g$ is a semi-Riemannian metric on $M$ such that $H$ is $g$-good, then we can identify $\bot H$ with $\bot_gH$ in such a way that $\pr^{\bot H}$ becomes $\pr^{\bot H}_g$ under the identification; cf.\ Remark \ref{normalbundle}. In this situation, we consider $\Twist_H$ as a section in $\Lambda^2(H^\ast)\otimes\bot_gH$, given by $\Twist_H(v,w) = \pr^{\bot H}_g([v,w])$.
\smallskip\\
If $x$ is a point in $M$, we call $H$ \emph{untwisted} [resp.\ \emph{twisted}] \emph{at $x$} if and only if the twistedness of $H$ vanishes [resp.\ does not vanish] in $x$. We call $H$ \emph{everywhere twisted} if and only if $\Twist_H\in C^\infty(M\ot\Lambda^2(H^\ast)\otimes\bot H)$ vanishes nowhere. (The analogous concept \emph{everywhere untwisted} is the same as \emph{integrable}.)
\end{definition}

\begin{remark}
When I mentioned a \emph{suitable nonintegrability condition} (for some distribution) in the preface and in Chapter \ref{ONE}, I always meant everywhere twistedness. This is a very weak nonintegrability condition, compared for instance to the contact condition for a hyperplane distribution on an odd-dimensional manifold: everywhere twistedness demands only that for each $x\in M$, there are at least two vectors $v,w\in H_x$ such that $\Twist_H(v,w)\neq0$. When the rank of $H$ is large, the contact condition demands much more; cf.\ Chapter \ref{FIVE} for a detailed discussion.
\end{remark}

\begin{fact} \label{sti}
Let $(M,g)$ be a semi-Riemannian manifold, let $H$ be a $g$-good distribution on $M$. We consider the section $T$ in $H^\ast\otimes H^\ast\otimes(\bot H)^\ast\otimes H^\ast\otimes H^\ast\otimes(\bot H)^\ast$ which is defined by
\[
T(u_0,v_0,w_0,u_1,v_1,w_1) \define g(\Twist_H(u_0,v_0),w_0)g(\Twist_H(u_0,v_0),w_1) \;\;,
\]
and we consider the contraction $t\in C^\infty(M,\R)$ of $T$ in the first and fourth, the second and fifth, and the third and sixth index (via the metrics $g\restrict H$ and $g\restrict\bot H$, of course). I.e., we consider, in classical Riemannian notation, the total contraction $t = 2\abs{\Twist_H}^2_g$. (The factor $2$ appears because we defined $\Twist_H$ as a section in $\Lambda^2(H^\ast)\otimes\bot H$ instead of $H^\ast\otimes H^\ast\otimes\bot H$.) Then we have
\[
t = 2(\sigma_{g,H}-\tau_{g,H}) \;\;.
\]
\end{fact}
\Proof
Let $x$ in $M$. By \ref{adaptedframe}, we can choose an $H$-adapted local orthonormal frame $(e_1,\dots,e_n)$ on a neighbourhood of $x$; here $n$ is the dimension of $M$. By \ref{sigmatauframe}, the following equation holds on this neighbourhood:
\[ \begin{split}
t &= \sum_{i,j:H}\sum_{k:\bot H}\eps_i\eps_j\eps_kg(\Twist_H(e_i,e_j),e_k)^2 = \sum_{i,j:H}\sum_{k:\bot H}\eps_i\eps_j\eps_kg(\nabla_{e_i}e_j -\nabla_{e_j}e_i,e_k)^2\\
&= \sum_{i,j:H}\sum_{k:\bot H}\eps_i\eps_j\eps_k(\Gamma^k_{ij}-\Gamma^k_{ji})^2 = 2\sum_{i,j:H}\sum_{k:\bot H}\eps_i\eps_j\eps_k(\Gamma^k_{ij})^2 -2\sum_{i,j:H}\sum_{k:\bot H}\eps_i\eps_j\eps_k\Gamma^k_{ij}\Gamma^k_{ji} = 2(\sigma_{g,H}-\tau_{g,H}) \;.\qedhere
\end{split} \]
\end{proof}
\
\smallskip\\
The following fact is the reason why everywhere twisted distributions play a crucial role in the pseudo-Riemannian prescribed scalar curvature problem.

\begin{proposition} \label{sigmatauintegrable}
Let $(M,g)$ be a Riemannian manifold, let $H$ be a distribution on $M$. Then $\sigma_{g,H}\geq\abs{\tau_{g,H}}$. Moreover, for every $x\in M$, the following statements are equivalent:
\begin{enumerate}
\item
$\sigma_{g,H}(x) = \tau_{g,H}(x)$.
\item
$H$ is untwisted at $x$.
\end{enumerate}
In particular, $H$ is integrable if and only if $\sigma_{g,H} = \tau_{g,H}$.
\end{proposition}
\Proof
By \ref{sti}, $\sigma_{g,H}-\tau_{g,H} = \abs{\Twist_H}^2_g \geq0$ (the inequality is valid since $g$ is Riemannian). Equality holds exactly in those points where $\Twist_H$ vanishes. It remains to prove $\sigma_{g,H} \geq -\tau_{g,H}$. This follows from an ON frame calculation like that in the proof of \ref{sti}:
\[ \begin{split}
0\leq \sum_{i,j:H}\sum_{k:\bot H}(\Gamma^k_{ij}+\Gamma^k_{ji})^2 = 2\sum_{i,j:H}\sum_{k:\bot H}(\Gamma^k_{ij})^2 +2\sum_{i,j:H}\sum_{k:\bot H}\Gamma^k_{ij}\Gamma^k_{ji} = 2(\sigma_{g,H}+\tau_{g,H}) \;\;.\qedhere
\end{split} \]
\end{proof}

\begin{remarks}{\ }
\begin{enumerate}
\item
The proposition shows that the zeroes of the nonnegative function $\sigma_{g,H} -\tau_{g,H}$ do not depend on the Riemannian metric $g$, but only on the distribution $H$. However, the zeroes of the nonnegative function $\sigma_{g,H}+\tau_{g,H}$ depend on $g$.
\item
We need Proposition \ref{sigmatauintegrable} only in the Riemannian version we stated, but it can be generalised:
\smallskip\\
Let $(M,g)$ be any semi-Riemannian manifold, and let $H$ be a maximally timelike distribution on $M$. Then all conclusions of the proposition are still true. This is easy to prove: In the ON frame calculation in the proofs of \ref{sti} and \ref{sigmatauintegrable}, we have  $\eps_i=\eps_j=-1$ and $\eps_k=1$ because $H$ is timelike and $\bot_gH$ is spacelike. Hence $\eps_i\eps_j\eps_k=1$ in each case, so the proof of \ref{sigmatauintegrable} still works.
\smallskip\\
If $H$ is not maximally timelike but maximally spacelike, then $\eps_i\eps_j\eps_k=-1$; so the equivalence of (i) and (ii) in the proposition remains true, but $\sigma_{g,H}\geq\abs{\tau_{g,H}}$ has to be replaced by $\sigma_{g,H}\leq-\abs{\tau_{g,H}}$.
\end{enumerate}
\end{remarks}

\begin{center} \fbox{\parbox{160mm}{
\emph{The rest of this chapter can be skipped: It contains facts and definitions which are not used in the proofs of the main results of this thesis, but are just cited later in the discussion of the esc Conjecture and in several side remarks.}}}
\end{center}

\subsection{The Laplacian with respect to a foliation}

\begin{definition} \label{flaplacedef}
Let $(M,g)$ be a semi-Riemannian manifold, let $H$ be an integrable $g$-good distribution on $M$, and let $f\in C^\infty(M,\R)$. Then we define $\flaplace_{g,H}(f)\in C^\infty(M,\R)$ to be the Laplacian of the restriction of $f$ to the leaves of the foliation corresponding to $H$, with respect to the restriction of the metric on the leaves. I.e., for each $x\in M$, the value of $\flaplace_{g,H}(f)$ in $x$ is the value in $x$ of $\laplace_{g\restrict F_x}(f\restrict F_x)\in C^\infty(F_x,\R)$, where $F_x$ denotes the (germ of the) leaf through $x$ of the foliation defined by $H$.
\end{definition}

Recall the following elementary fact:
\begin{fact} \label{foliationnabla}
Let $(M,g)$ be a semi-Riemannian manifold, let $F$ be a foliation on $M$ such that the distribution $TF$ on $M$ is $g$-good, let $\nabla$ be the \LeviCivita\ connection on $(M,g)$. Let $\nabla^F$ be the \LeviCivita\ connection on $(F,g\restrict F)$; i.e., for all $x\in M$, $u\in T_xF$ and each section $v\in C^\infty(M\ot TF)$, the vector $\nabla^F_uv\in T_xF$ is by definition the vector $\nabla^N_u(v\restrict N)$, where $\nabla^N$ is the \LeviCivita\ connection on the (germ of the) leaf $N=F_x$ through $x$ with respect to the metric $g\restrict N$.
\smallskip\\
Then for all $x\in M$, all $u,w\in T_xF$ and every section $v\in C^\infty(M\ot TF)$, we have $g(\nabla_uv,w) = g(\nabla^F_uv,w)$; in other words, $\pr^{TF}_g(\nabla_uv) = \nabla^F_uv$.
\end{fact}
\Proof
We extend the vectors $u,w$ to sections in $TF\to M$. The Koszul formula yields
\[
2g(\nabla^F_uv,w) = \partial_ug(v,w) +\partial_vg(u,w) -\partial_wg(u,v) +g([u,v],w) +g([w,u],v) +g([w,v],u) = 2g(\nabla_uv,w) \;.\qedhere
\]
\end{proof}

From this, we deduce that the function $\flaplace_{g,H}(f)$ is equal to another one that we knew before:

\begin{lemma}
Let $(M,g)$ be a semi-Riemannian manifold, let $H$ be an integrable $g$-good distribution on $M$, and let $f\in C^\infty(M,\R)$. Then
\[
\flaplace_{g,H}(f) = \laplace^H_{g,H}(f) \;\;.
\]
\end{lemma}
\Proof
Let $(e_1,\dots,e_n)$ be an $H$-adapted local $g$-orthonormal frame of $TM$. We denote the orthonormal Christoffel symbols with respect to the \LeviCivita\ connection of $g$ by $\Gamma^k_{ij}$, as usual. Since the set of all $e_i$ with $i:H$ is a $(g\restrict H)$-orthonormal frame of $H$, we get by \ref{laplacianformulae}:
\[
\flaplace_{g,H}(f) = \sum_{i:H}\epsilon_i\partial_{e_i}\partial_{e_i}f + \sum_{i,j:H}\eps_i\eps_j\Gamma^j_{ji}df(e_i) = \laplace^H_{g,H}(f) \;\;;
\]
the first equality holds because of \ref{foliationnabla}.
\end{proof}

\subsection{The scalar curvature of a foliation}

\begin{definition}[$\fscal_{g,H}$] \label{fscal}
Let $(M,g)$ be a semi-Riemannian manifold, and let $H$ be an integrable $g$-good distribution on $M$. Then we define $\fscal_{g,H}\in C^\infty(M,\R)$ to be the scalar curvature of the leaves of the foliation corresponding to $H$. I.e., for each $x\in M$, the value of $\fscal_{g,H}$ in $x$ is the value in $x$ of the scalar curvature of the (germ of the) leaf $F_x$ through $x$ with respect to the semi-Riemannian metric $g\restrict H$ on the tangent bundle $TF_x = H\restrict F_x$.
\end{definition}

$\fscal_{g,H}$ does obviously not depend on the whole metric $g$ but only on its restriction to $H$. It is related to the functions from Section \ref{TWOTWO} by the following formula:

\begin{lemma} \label{fscallemma}
Let $(M,g)$ be a semi-Riemannian manifold, and let $H$ be an integrable $g$-good distribution on $M$. Then
\[
\fscal_{g,H} = \scal_g^{H,H} +\eval{\divergence^H_g}{\divergence^H_g}_{g,\bot H} -\sigma_{g,H} \;\;.
\]
\end{lemma}
\Proof
It suffices to check the formula locally, so let $(e_1,\dots,e_n)$ be any local $H$-adapted $g$-orthonormal frame. Since $H$ is integrable, we have $\sigma_{g,H}=\tau_{g,H}$ and $\Gamma^k_{ij}=\Gamma^k_{ji}$ for all $i,j:H$ and $k:\bot H$ (because $0 = g(\Twist_H(e_i,e_j),e_k) = \Gamma^k_{ij}-\Gamma^k_{ji}$ in this case). This implies for every $j:\bot H$:
\[
\sum_{i,k:H}\eps_i\eps_k\Gamma^k_{ij}\Gamma^k_{ji} = \sum_{i,k:H}\eps_i\eps_k\Gamma^j_{ik}\Gamma^i_{jk} = \sum_{i,k:H}\eps_i\eps_k\Gamma^j_{ki}\Gamma^i_{jk} = -\sum_{i,k:H}\eps_i\eps_k\Gamma^i_{kj}\Gamma^i_{jk} = -\sum_{i,k:H}\eps_i\eps_k\Gamma^k_{ij}\Gamma^k_{ji} \;\;,
\]
i.e.\ $\sum_{i,k:H}\eps_i\eps_k\Gamma^k_{ij}\Gamma^k_{ji}=0$. Hence \ref{curvatureON} yields
\[
\scal^{H,H}_g = -2\sum_{j:H}\eps_j\partial_{e_j}\divergence^H_g(e_j) -\sum_j\eps_j\divergence^H_g(e_j)^2 -\sum_{i,j,k:H}\eps_i\eps_j\eps_k\Gamma^k_{ij}\Gamma^k_{ji} +\sigma_{g,H} \;\;.
\]
Consider (a germ of) a fixed leaf $F$ of the foliation defined by $H$. Since the set of all $e_i\restrict F$ with $i:H$ is a $(g\restrict F)$-orthonormal frame of $TF$, the function $\fscal_{g,H}$ has, by \ref{Riemanntensor}, the following form:
\[
\fscal_{g,H} = -2\sum_{i:H}\eps_i\partial_{e_i}\divergence^H_g(e_i) -\sum_{i:H}\eps_i\divergence^H_g(e_i)^2 -\sum_{i,j,k:H}\eps_i\eps_j\eps_k\Gamma^k_{ij}\Gamma^k_{ji} \;\;.
\]
(Here we used again Fact \ref{foliationnabla}.) Thus
\[
\scal^{H,H}_g = \fscal_{g,H} -\sum_{j:\bot H}\eps_j\divergence^H_g(e_j)^2 +\sigma_{g,H} = \fscal_{g,H} -\eval{\divergence^H_g}{\divergence^H_g}_{g,\bot H} +\sigma_{g,H} \;\;.\qedhere
\]
\end{proof}


\section{Partial integration formulae}

The following lemma generalises the well-known equation
\[
\int_{(M,g)}\eval{df}{dh}_{g} = -\int_{(M,g)}\laplace_g(f)\cdot h
\]
for functions $f,h$ on a compact semi-Riemannian manifold (in the sense that this equation is the special case of the lemma in which $H=TM$).

\begin{lemma} \label{adjointlaplace}
Let $(M,g)$ be a compact semi-Riemannian manifold, let $V$ be a $g$-good distribution on $M$, let $f,h\in C^\infty(M,\R)$ with $h\restrict\mfbd M \equiv 0$, and let $H$ denote the $g$-orthogonal distribution of $V$. Then
\[
\int_{(M,g)}\eval{df}{dh}_{g,H} = -\int_{(M,g)}\Big(\laplace^H_{g,H}(f)+\eval{\divergence^V_g}{df}_{g,H}\Big)h \;\;.
\]
\end{lemma}
\Proof
Using a partition of unity, we can decompose $h$ into a finite sum of functions each summand of which has support in a contractible subset of $M$. If the statement of the lemma holds for each summand, then it does also hold for $h$. Therefore it suffices to prove the lemma in the case in which $M$ is orientable. We choose an orientation.
\smallskip\\
We define the vector field $X\define\pr^H_g(\grad_g(f))$. By Definition \ref{laplacedef} and Remark \ref{laplacedivergence}, we have
\[
\divergence_g(X) = \divergence^H_g(X) +\divergence^V_g(X) = \laplace^H_{g,H}(f) +\laplace^V_{g,H}(f) = \laplace^H_{g,H}(f) +\eval{\divergence^V_g}{df}_{g,H} \;\;.
\]
Let $\vol_g$ denote the volume form of $g$ with respect to the orientation. The Lie derivative $L_X\omega$ of the $n$-form $\omega\define h\vol_g$ satisfies $L_X\omega = i_X(d\omega) +d(i_X\omega) = d(i_X\omega)$. Since $\omega$ and thus $i_X\omega$ vanish on the boundary of $M$, we get by Stokes' theorem
\[
0 = \int_ML_X\omega \;\;.
\]
From the well-known equation $L_X\vol_g = \divergence_g(X)\vol_g$, we infer
\[ \begin{split}
L_X\omega &= L_X(h)\vol_g +hL_X\vol_g = dh(X)\vol_g +hL_X\vol_g\\
&= \bigg(dh\Big(\pr^H_g(\grad_g(f))\Big) +h\divergence_g(X)\bigg)\vol_g\\
&= \bigg(\eval{df}{dh}_{g,H} +\Big(\laplace^H_{g,H}(f) +\eval{\divergence^V_g}{df}_{g,H}\Big)h\bigg)\vol_g \;\;.
\end{split} \]
This implies the statement of the lemma.
\end{proof}

\begin{remark} \label{dfdffunctional}
The preceding lemma shows in particular that the function $\laplace^H_{g,H}(f)+\eval{\divergence^V_g}{df}_{g,H}$ has variational form, i.e., $\laplace^H_{g,H}(f)+\eval{\divergence^V_g}{df}_{g,H} = 0$ is the Euler/Lagrange equation of some functional: If $B$ is a suitable Banach space of real-valued functions on $M$ with zero (Dirichlet) boundary values (e.g.\ the Sobolev space $H^{1,2}_0(M,\R)$), and if $f_0$ is a suitable real-valued function on $M$ (not necessarily vanishing on the boundary), then the functional $E\colon f_0+B \to \R$ given by
\[
E(f) \define \int_{(M,g)}\eval{df}{df}_{g,H}
\]
is well-defined, \Frechet\ differentiable, and its \Frechet\ derivative $D_fE \colon B\to \R$ in the point $f$ is given by
\[
(D_fE)(h) = -2\int_{(M,g)}\Big(\laplace^H_{g,H}(f)+\eval{\divergence^V_g}{df}_{g,H}\Big)h \;\;.
\]
In particular, $f$ is a critical point of $E$ if and only if $\laplace^H_{g,H}(f)+\eval{\divergence^V_g}{df}_{g,H} = 0$.
\smallskip\\
(This is just a side remark, so we may be vague about the allowed choices of $B$ and $f_0$, and we may omit the proof.)
\end{remark}

\begin{lemma} \label{adjointdivergence}
Let $(M,g)$ be a compact semi-Riemannian manifold, let $V$ be a $g$-good distribution on $M$, let $u\in C^\infty(M,\R)$ with $u\restrict\mfbd M \equiv 0$, and let $H$ denote the $g$-orthogonal distribution of $V$. Then
\[
\int_{(M,g)}\eval{\divergence^V_g}{du}_{g,H} = -\int_{(M,g)}\Big(\qual^H_g +\eval{\divergence^V_g}{\divergence^V_g}_{g,H} -\tau_{g,H}\Big)u \;\;.
\]
\end{lemma}
\Proof
There is an open cover of $M$ each element of which admits a $V$-adapted $g$-orthonormal frame. We choose a finite subcover and a subordinate partition of unity. Using this, we can write $u$ as a finite sum of functions each of which vanishes on $\mfbd M$ and has support in some element of the open cover. Since it suffices to prove the lemma for each summand, we can assume that $M$ admits a global $V$-adapted $g$-orthonormal frame $(e_1,\dots,e_n)$.
\smallskip\\
By \ref{moreONformulae} and \ref{qON}, we get
\[ \begin{split}
&\eval{\divergence^V_g}{du}_{g,H} +\Big(\qual^H_g +\eval{\divergence^V_g}{\divergence^V_g}_{g,H} -\tau_{g,H}\Big)\,u\\
&= \sum_{i:H}\eps_i\bigg(\divergence^V_g(e_i)du(e_i) +\partial_{e_i}\divergence^V_g(e_i)\,u +\divergence^V_g(e_i)\divergence^H_g(e_i)\,u +\divergence^V_g(e_i)\divergence^V_g(e_i)\,u\bigg) \;\;.
\end{split} \]
Therefore we just have to prove that
\[
\int_{(M,g)}\bigg(\partial_{e_i}\Big(\divergence^V_g(e_i)u\Big) +\divergence^V_g(e_i)u\divergence_g(e_i)\bigg) = 0
\]
for each $i:H$. This follows like in the proof of Lemma \ref{adjointlaplace}: Let $\vol_g$ denote the volume form of $g$ with respect to some orientation of $M$ (recall that we assume $M$ to be parallelisable and hence orientable), let $\omega$ denote the $n$-form $\divergence^V_g(e_i)u\vol_g$, and let $X$ denote the vector field $e_i$. Then $L_X\omega = d(i_X\omega)$, and since $i_X\omega$ vanishes on the boundary of $M$, we get by Stokes' theorem
\[
0 = \int_ML_X\omega \;\;.
\]
But
\[
L_X\omega = L_X\Big(\divergence^V_g(e_i)u\Big)\vol_g +\divergence^V_g(e_i)u\,L_X\vol_g = \bigg(\partial_{e_i}\Big(\divergence^V_g(e_i)u\Big) +\divergence^V_g(e_i)u\divergence_g(e_i)\bigg)\vol_g \;\;,
\]
so the proof is complete.
\end{proof}


\chapter{Modifications of the metric} \label{THREE}

In Subsection \ref{ONETWOONE} of the introduction, we have seen how one can modify a given Riemannian or Lorentzian product metric on a product manifold $M = S^1\times N$, and how these modifications change the scalar curvature. We verified the leading order terms in the corresponding formulae. In the present chapter, we are going to generalise those modifications to the situation where $M$ is an arbitrary (not necessarily product) manifold, and where the first-factor (line) distribution on $S^1\times N$ is replaced by an arbitrary $q$-plane distribution on $M$.
\smallskip\\
While the definitions generalise in an easy and obvious way, the --- completely straightforward --- computations of the formulae describing the change of scalar curvature become much longer. In fact, they occupy the whole chapter, since they will be spelled out in such detail that the reader can check them without doing separate auxiliary calculations.

\section{Definition of the $\switch$, $\stre$, $\conf$ operations}

\subsection{Switching}

Our first construction turns Riemannian metrics into pseudo-Riemannian metrics and vice versa:

\begin{definition}[switched metric] \label{switchdefinition}
Let $(M,g)$ be a semi-Riemannian manifold, let $V$ be a $g$-good distribution on $M$. Then we define the semi-Riemannian metric $\switch(g,V)$ on $M$ as follows. Let $H$ be the $g$-orthogonal distribution of $V$. Then $TM = V\oplus H$, and $g$ has, with respect to this decomposition, the form $g_V\oplus g_H$, where $g_V$ is a semi-Riemannian metric on the vector bundle $V$ and $g_H$ is a semi-Riemannian metric on $H$. We define $\switch(g,V)$ to be the metric $(-g_V)\oplus g_H$ on the vector bundle $V\oplus H = TM$. (This is indeed a semi-Riemannian metric on $M$ since $-g_V$ and $g_H$ are semi-Riemannian metrics on the vector bundles $V$ and $H$, respectively.)
\smallskip\\
In other words, we define $\switch(g,V)$ by
\begin{equation}
\switch(g,V)(w,z) = -g(w_V,z_V) +g(w_H,z_H)
\end{equation}
for all $x\in M$ and $z,w\in T_xM$; here $u_V\define\pr^V_g(u)$ denotes the $g$-orthogonal projection of $u\in TM$ to $V$, and $u_H\define\pr^H_g(u)$ denotes the $g$-orthogonal projection of $u$ to $H$.
\end{definition}

\begin{remark} \label{switchtrivialities}
Let $(M,g)$ be a semi-Riemannian manifold, let $V$ be a $g$-good distribution on $M$. Then $V$ is $\switch(g,V)$-good, and the $g$-orthogonal distribution of $V$ is equal to the $\switch(g,V)$-orthogonal distribution of $V$.
\end{remark}
\Proof
Since the restriction of $g$ to $V$ is nondegenerate, the restriction of $-g$ to $V$ --- i.e.\ the restriction of $\switch(g,V)$ to $V$ --- is nondegenerate, too. In other words, $V$ is $\switch(g,V)$-good.
\smallskip\\
Let $H\define\bot_gV$. A vector $w\in TM$ is contained in $\bot_{\switch(g,V)}V$ if and only if $\switch(g,V)(v,w)=0$ for all $v\in V$, i.e.\ if and only if $0=g(v_V,w_V)$ for all $v\in V$, i.e.\ (since $g\restrict V$ is nondegenerate) if and only if $w_V=0$, that is, if and only if $w\in H$. Hence the $g$-orthogonal distribution of $V$ is equal to the $\switch(g,V)$-orthogonal distribution of $V$.
\end{proof}

The only special cases we are interested in are those when $g$ is Riemannian or $V$ is maximally timelike:

\begin{remark}
If $g$ is a Riemannian metric and $V$ is a $q$-plane distribution on $M$, then $\switch(g,V)$ is a semi-Riemannian metric of index $q$ which makes $V$ (maximally) timelike. If $g$ is a semi-Riemannian metric of index $q$ and $V$ is a (maximally) timelike $q$-plane distribution on $M$, then $\switch(g,V)$ is a Riemannian metric.
\end{remark}

\begin{remark}[idempotency]
Let $(M,g)$ be a semi-Riemannian manifold, let $V$ be a $g$-good distribution on $M$. Then $\switch(\switch(g,V),V)=g$. (That's why we call it the ``switch'' operation.) If $V$ is the unique $0$-plane distribution on $M$, then $\switch(g,V)=g$.
\end{remark}

We will make use of the preceding elementary remarks without further mention.
\begin{remark}[adapted ON frames] \label{switchONframe}
Let $(M,g)$ be a semi-Riemannian $n$-manifold, let $V$ be a $g$-good distribution on $M$, and let $(e_1,\dots,e_n)$ be a (local) $V$-adapted $g$-orthonormal frame of $TM$. We use the abbreviations $\bar{g}\define\switch(g,V)$ and $H\define \bot_gV = \bot_{\bar{g}}V$. Then $(e_1,\dots,e_n)$ is obviously also a (local) $V$-adapted $\bar{g}$-orthonormal frame of $TM$; more precisely, $\bar{g}(e_i,e_i) = g(e_i,e_i)$ if $i:H$, whereas $\bar{g}(e_i,e_i) = -g(e_i,e_i)$ if $i:V$.
\end{remark}

\subsection{Stretching}

The second modification ``stretches'' a metric along a given distribution, by an amount which is specified by a function on the manifold. This construction is a generalisation of warped product metrics.

\begin{definition}[stretched metric] \label{stretchdefinition}
Let $(M,g)$ be a semi-Riemannian manifold of index $q$, let $V$ be a $g$-good distribution on $M$, and let $f\in C^\infty(M,\R_{>0})$. Using the same notation as in the definition of the $\switch$ operation, we define the semi-Riemannian metric $\stre(g,f,V)$ (of index $q$) on $M$ to be $(f^{-2}g_V)\oplus g_H$. (This is indeed a semi-Riemannian metric of index $q$ on $M$ since $f^{-2}g_V$ and $g_H$ are semi-Riemannian metrics on the vector bundles $V$ and $H$, respectively, and $f^{-2}g_V$ has the same index as $g_V$.) In other words, $\stre(g,f,V)$ is given by
\begin{equation}
\stre(g,f,V)(w,z) = \frac{1}{f^2}g(w_V,z_V) +g(w_H,z_H) \;\;.
\end{equation}
\end{definition}

\begin{remark} \label{stretchtrivialities}
Let $(M,g)$ be a semi-Riemannian manifold, let $V$ be a $g$-good distribution on $M$, and let $f\in C^\infty(M,\R_{>0})$. Then $V$ is $\stre(g,f,V)$-good, and the $g$-orthogonal distribution of $V$ is equal to the $\stre(g,f,V)$-orthogonal distribution of $V$.
\end{remark}

\begin{remark}
Let $(M,g)$ be a semi-Riemannian manifold, let $V$ be a $g$-good distribution on $M$, and let $f,f_0,f_1\in C^\infty(M,\R_{>0})$. Then we have obviously
\[ \begin{split}
\stre(\stre(g,f_0,V),f_1,V) &= \stre(g,f_0f_1,V) \;\;,\\
\stre(\switch(g,V),f,V) &= \switch(\stre(g,f,V),V) \;\;.
\end{split} \]
If $V$ is the unique $0$-plane distribution on $M$ or $f$ is the constant $1$, then $\stre(g,f,V)=g$.
\end{remark}

\begin{remark}[warped products as a special case]
Let $(B,g_B)$ and $(F,g_F)$ be semi-Riemannian manifolds, and let $f\in C^\infty(B,\R_{>0})$. Recall that the \emph{warped product} $B\times_fF$ is the manifold $M\define B\times F$ equipped with the semi-Riemannian metric $g_B\oplus f^2g_F = \pi_B^\ast(g_B) +(f\compose\pi_B)^2\pi_F^\ast(g_F)$ (where $\pi_B \colon M\to B$ and $\pi_F\colon M\to F$ denote the obvious projections).
\smallskip\\
In this situation, let $g$ denote the semi-Riemannian product metric $g_B\oplus g_N = \pi_B^\ast(g_B)+\pi_F^\ast(g_F)$ on $M$, and let $V$ denote the second-factor distribution on $M = B\times F$. Then the warped product metric of $B\times_fF$ is equal to $\stre(g,1/(f\compose\pi_B),V)$.
\smallskip\\
Warped product metrics are thus a special case of our stretch metrics. Note that even in the case where $M$ is a product manifold and $V$ is the second-factor distribution on $M$, stretch metrics are much more general than warped products because the stretch factor is a function which may depend on the whole manifold $M = B\times F$, while the warp factor is a function which depends only on $B$.
\end{remark}

\begin{remark}[adapted ON frames] \label{stretchONframe}
Let $(M,g)$ be a semi-Riemannian $n$-manifold, let $V$ be a $g$-good distribution on $M$, let $f\in C^\infty(M,\R_{>0})$, and let $(e_1,\dots,e_n)$ be a (local) $V$-adapted $g$-orthonormal frame of $TM$. We use the abbreviations $\bar{g}\define\stre(g,f,V)$ and $H\define \bot_gV = \bot_{\bar{g}}V$. Then the tuple $(\bar{e}_1,\dots,\bar{e}_n)$, where $\bar{e}_i\define e_i$ if $i:H$, and $\bar{e}_i\define fe_i$ if $i:V$, is a (local) $V$-adapted $\stre(g,f,V)$-orthonormal frame of $TM$ such that $g(e_i,e_i) = \bar{g}(\bar{e}_i,\bar{e}_i)$ for all $i\in\set{1,\dots,n}$.
\end{remark}
\Proof
This is obvious from
\[
\bar{g}(\bar{e}_i,\bar{e}_j) =
\left.\begin{cases}
  \stre(g,f,V)(fe_i,fe_j)    &\text{if $i,j:V$}\\
  \stre(g,f,V)(e_i,e_j)      &\text{if $i,j:H$}\\
  \stre(g,f,V)(fe_i,e_j) = 0 &\text{if $i:V$ and $j:H$}
\end{cases}\right\} = g(e_i,e_j) \;\;.\qedhere
\]
\end{proof}

\subsection{Conformal deformation}

This way of modifying a semi-Riemannian metric is well-known. It is a special case of stretching a metric.

\begin{definition}[conformally deformed metric] \label{confdefinition}
Let $(M,g)$ be a semi-Riemannian manifold of index $q$, let $\kappa\in C^\infty(M,\R_{>0})$. Then we define the semi-Riemannian metric $\conf(g,\kappa)$ (of index $q$) on $M$ to be the conformally deformed metric $\kappa^{-2}g$. In other words, $\conf(g,\kappa) = \stre(g,\kappa,TM)$.
\end{definition}

\begin{remark}
Let $(M,g)$ be a semi-Riemannian manifold, let $V$ be a $g$-good distribution on $M$, and let $\kappa\in C^\infty(M,\R_{>0})$. Then $V$ is $\conf(g,\kappa)$-good, and the $g$-orthogonal distribution of $V$ is equal to the $\conf(g,\kappa)$-orthogonal distribution of $V$. With respect to this orthogonal distribution $H$, we have $\conf(g,\kappa) = \stre(\stre(g,\kappa,H),\kappa,V)$ (since $\kappa^{-2}g = (\kappa^{-2}g_V)\oplus(\kappa^{-2}g_H)$).
\end{remark}

\begin{remark}
Let $(M,g)$ be a semi-Riemannian manifold, let $V$ be a $g$-good distribution on $M$, and let $f,\kappa\in C^\infty(M,\R_{>0})$. Then
\[ \begin{split}
\conf(\stre(g,f,V),\kappa) &= \stre(\conf(g,\kappa),f,V) \;\;,\\
\conf(\switch(g,V),\kappa) &= \switch(\conf(g,\kappa),V) \;\;.
\end{split} \]
\end{remark}

Now we can state the aim of this chapter more precisely: Let $g$ be a Riemannian metric on some manifold $M$, let $V$ be a $q$-plane distribution on $M$, and let $f,\kappa\in C^\infty(M,\R_{>0})$. We want to derive a formula for the scalar curvature of the index-$q$ semi-Riemannian metric $\conf(\stre(\switch(g,V),f,V),\kappa)$.
\smallskip\\
This formula contains all the functions that we have defined in Chapter \ref{TWO} (i.e. $\qual^V_g, \eval{\divergence_g^{\bot V}}{df}_{g,V}$, etc.). We could compute it in one (huge) step, but I prefer to state in separate formulae how the scalar curvature and all those other functions behave under the $\switch$, the $\stre$, and the $\conf$ operations. This has the advantage of splitting the computation into smaller portions, and of providing more information.
\smallskip\\
Note that the metric $\conf(\stre(\switch(g,V),f,V),\kappa)$ specialises in the situation of Subsection \ref{ONETWOONE} to the metric $h(\kappa,f)$ we considered there (recall that $g$ was the Riemannian product metric $dt^2\oplus g_N$, and $V$ was the first-factor distribution on $M=S^1\times N$). So in this chapter, we will in particular generalise Equation \eqref{babymix} and prove that generalisation. (Recall that we had already proved Equation \eqref{babymix} up to terms of order less than $2$ in $\kappa$ and $f$.)

\section{Formulae for switching}

\emph{Throughout this section, $M$ is an $n$-dimensional manifold, $V$ is a $q$-plane distribution on $M$, $g$ is a Riemannian metric on $M$, $h$ denotes the semi-Riemannian metric $\switch(g,V)$ (of index $q$), and $H$ denotes the $g$-orthogonal distribution of $V$ (which is also the $h$-orthogonal distribution of $V$; cf.\ Remark \ref{switchtrivialities}).}
\medskip\\
Our aim is to express functions which are defined by the semi-Riemannian metric $h$ --- i.e.\ the functions $\scal_h$, $\qual^V_h$, $\sigma_{h,H}$, etc.\ --- in terms of functions which are defined by the Riemannian metric $g$ (i.e.\ $\scal_g$, $\qual^V_g$, $\sigma_{g,H}$, etc.). We will summarise the results of our calculations in Theorem \ref{switchformulae} at the end of this section.

\begin{remark} \label{switchgeneral}
We assume that $g$ is a Riemannian (instead of an arbitrary semi-Riemannian) metric just for simplicity; it is the only case we will need later. In this way, we avoid additional $\eps_i$s in the computations: For the general case, we had to introduce numbers $\bar{\eps}_1,\dots,\bar{\eps}_n\in\set{1,-1}$ by $\bar{\eps}_i = \eps_i$ if $i:H$, and $\bar{\eps}_i = -\eps_i$ if $i:V$, and put them into the calculations below. (The $\eps_i$s are defined by $\eps_i = \switch(g,V)(e_i,e_i)$, and the $\bar{\eps}_i$s would be defined by $\bar{\eps}_i = g(e_i,e_i)$. In the Riemannian case, $\bar{\eps}_i=1$ for all $i$.) By making these modifications, we could prove that the formulae from Theorem \ref{switchformulae} below hold also in the general semi-Riemannian case.
\end{remark}

\begin{remark} \label{wickedtrick}
The computations of the formulae for the switch and stretch modifications could be unified; i.e., there is a general computation which has these two as special cases: We just had to allow in the stretch computation that the (nonvanishing) function $f$ takes values in the purely imaginary complex numbers. The case $f\equiv\sqrt{-1}$ would then correspond to the switch operation. This formal trick (one could call it a \emph{Wick rotation}) can easily be justified --- i.e.\ shown to yield the correct results ---, but it might be a bit confusing. Therefore we will do two separate computations.
\end{remark}

\begin{remark}
One might ask why we compute the switch and stretch formulae for the scalar curvature in a direct way instead of systematically computing first the Riemann tensor, then the Ricci tensor, and finally the scalar curvature. The answer is of course that the computation becomes much shorter this way (shortness is relative), and that we do not need the Riemann and Ricci tensors.
\end{remark}

Now we start with the computation. By Corollary \ref{adaptedframe}, each point in $M$ has an open neighbourhood which admits a $V$-adapted $g$-orthonormal frame $(e_1,\dots,e_n)$. This frame is also a $V$-adapted $h$-orthonormal frame by Remark \ref{switchONframe}, and the numbers $\eps_i\define h(e_i,e_i) \in \set{1,-1}$ satisfy
\[
\eps_i = \begin{cases} -1 &\text{if $i:V$}\\ \phantom{-}1 &\text{if $i:H$} \end{cases} \;\;.
\]
We denote the \LeviCivita\ connections of $g$ and $h$ by $\nabla^{(g)}$ and $\nabla^{(h)}$, respectively, and we define the ON Christoffel symbols (with respect to the ON frame $(e_1,\dots,e_n)$) by
\[ \begin{split}
G^k_{ij} &\define g(\nabla^{(g)}_{e_i}e_j,e_k) \;\;,\\
H^k_{ij} &\define h(\nabla^{(h)}_{e_i}e_j,e_k)
\end{split} \]
for $i,j,k\in\set{1,\dots,n}$ (cf.\ Definition \ref{ONChristoffel}, and recall the facts from Remark \ref{Christoffel}).
\smallskip\\
For each $i\in\set{1,\dots,n}$, we define $\delta_{iV}\in\set{0,1}$ by
\[
\delta_{iV} = \begin{cases} 1 &\text{if $i:V$}\\ 0 &\text{if $i:H$} \end{cases} \;\;.
\]
Every local vector field $v$ on $M$ satisfies obviously
\[
h(v,e_i) =
\left.\begin{cases}
  -g(v,e_i)           &\text{if $i:V$}\\
  \phantom{-}g(v,e_i) &\text{if $i:H$}
\end{cases}\right\} = (1-2\delta_{iV})g(v,e_i) \;\;.
\]

\subsection{The orthonormal Christoffel symbols}

We compute the ON Christoffel symbols $H^k_{ij}$ of the metric $h$ in terms of the ON Christoffel symbols $G^k_{ij}$ of the metric $g$. The Koszul formula yields
\[ \begin{split}
2G^k_{ij} &= g([e_i,e_j],e_k)+g([e_k,e_i],e_j)+g([e_k,e_j],e_i) \;\;,\\
2H^k_{ij} &= h([e_i,e_j],e_k)+h([e_k,e_i],e_j)+h([e_k,e_j],e_i)\\
&= g([e_i,e_j],e_k)+g([e_k,e_i],e_j)+g([e_k,e_j],e_i)\\
&\mspace{20mu}-2\bigg(\delta_{kV}g([e_i,e_j],e_k) +\delta_{jV}g([e_k,e_i],e_j) +\delta_{iV}g([e_k,e_j],e_i)\bigg) \;\;,
\end{split} \]
hence (using $g([e_a,e_b],e_c) = g(\nabla_{e_a}e_b -\nabla_{e_b}e_a,e_c) = G^c_{ab}-G^c_{ba}$):
\begin{equation} \label{switchChristoffel} \begin{split}
H^k_{ij} &= G^k_{ij} -\delta_{kV}(G^k_{ij}-G^k_{ji}) -\delta_{jV}(G^j_{ki}-G^j_{ik}) -\delta_{iV}(G^i_{kj}-G^i_{jk})\\
&= \begin{cases}
  G^k_{ij} &\text{if $i,j,k:H$}\\
  G^k_{ij} -G^k_{ij}+G^k_{ji} = G^k_{ji} &\text{if $i,j:H$, $k:V$}\\
  G^k_{ij} -G^j_{ki}+G^j_{ik} = G^i_{kj} &\text{if $i,k:H$, $j:V$}\\
  G^k_{ij} -G^i_{kj}+G^i_{jk} = G^k_{ij}+G^j_{ki}+G^i_{jk} &\text{if $j,k:H$, $i:V$}\\
  G^k_{ij} -G^j_{ki}+G^j_{ik} -G^i_{kj}+G^i_{jk} = -G^k_{ji} &\text{if $k:H$, $i,j:V$}\\
  G^k_{ij} -G^k_{ij}+G^k_{ji} -G^i_{kj}+G^i_{jk} = G^j_{ki} &\text{if $j:H$, $i,k:V$}\\
  G^k_{ij} -G^k_{ij}+G^k_{ji} -G^j_{ki}+G^j_{ik} = G^k_{ji}+G^i_{kj}+G^j_{ik} &\text{if $i:H$, $j,k:V$}\\
  G^k_{ij} -G^k_{ij}+G^k_{ji} -G^j_{ki}+G^j_{ik} -G^i_{kj}+G^i_{jk} = -G^k_{ij} &\text{if $i,j,k:V$}
\end{cases} \;\;.
\end{split} \end{equation}

\subsection{Divergences, $\sigma$s, $\tau$s}

In particular, \eqref{switchChristoffel} implies for all $j,k$:
\[
\eps_kH^k_{kj} = \left.\begin{cases}
   G^k_{kj}            &\text{if $j,k:H$}\\
  -G^j_{kk} = G^k_{kj} &\text{if $j:H$, $k:V$}\\
   G^k_{kj}            &\text{if $j:V$, $k:H$}\\
  -(-G^k_{kj})         &\text{if $j,k:V$}
\end{cases}\right\} = G^k_{kj} \;\;;
\]
and thus (cf.\ \ref{divUformulae}):
\begin{equation} \label{switchdivergence} \begin{split}
\divergence_h(e_j) = \sum_k\eps_kH^k_{kj} = \sum_kG^k_{kj} = \divergence_g(e_j) \;\;,\\
\divergence^V_h(e_j) = \sum_{k:V}\eps_kH^k_{kj} = \sum_{k:V}G^k_{kj} = \divergence^ V_g(e_j) \;\;,\\
\divergence^H_h(e_j) = \sum_{k:H}\eps_kH^k_{kj} = \sum_{k:H}G^k_{kj} = \divergence^H_g(e_j) \;\;.
\end{split} \end{equation}

This yields (cf.\ \ref{moreONformulae}):
\begin{equation} \begin{split}
\eval{\divergence^V_h}{\divergence^V_h}_{h,H} &= \sum_{i:H}\eps_i\divergence^V_h(e_i)^2 = \sum_{i:H}\divergence^V_g(e_i)^2 = \eval{\divergence^V_g}{\divergence^V_g}_{g,H} \;\;,\\
\eval{\divergence^H_h}{\divergence^H_h}_{h,V} &= \sum_{i:V}\eps_i\divergence^H_h(e_i)^2 = -\sum_{i:V}\divergence^H_g(e_i)^2 = -\eval{\divergence^H_g}{\divergence^H_g}_{g,V} \;\;.
\end{split} \end{equation}

From \ref{sigmatauframe} and \eqref{switchChristoffel}, we get
\begin{equation} \label{switchsigmatau} \begin{split}
\sigma_{h,H} &= -\sum_{i,j:H}\sum_{k:V}H^k_{ij}H^k_{ij} = -\sum_{i,j:H}\sum_{k:V}G^k_{ji}G^k_{ji} = -\sigma_{g,H} \;\;,\\
\sigma_{h,V} &= \sum_{i,j:V}\sum_{k:H}H^k_{ij}H^k_{ij} = \sum_{i,j:V}\sum_{k:H}G^k_{ji}G^k_{ji} = \sigma_{g,V} \;\;,\\
\tau_{h,H} &= -\sum_{i,j:H}\sum_{k:V}H^k_{ij}H^k_{ji} = -\sum_{i,j:H}\sum_{k:V}G^k_{ji}G^k_{ij} = -\tau_{g,H} \;\;,\\
\tau_{h,V} &= \sum_{i,j:V}\sum_{k:H}H^k_{ij}H^k_{ji} = \sum_{i,j:V}\sum_{k:H}G^k_{ji}G^k_{ij} = \tau_{g,V} \;\;.
\end{split} \end{equation}

\subsection{Scalar curvatures and qualar curvatures}

The next equations follow from \eqref{switchdivergence}, \eqref{switchsigmatau}, and \ref{qON}:
\pagebreak \begin{subequations}
\begin{equation} \label{switchqualV} \begin{split}
\qual^V_h &= \sum_{i:V}\eps_i\partial_{e_i}\divergence^H_h(e_i) +\sum_{i:V}\eps_i\divergence^V_h(e_i)\divergence^H_h(e_i) +\tau_{h,V}\\
&= -\sum_{i:V}\partial_{e_i}\divergence^H_g(e_i) -\sum_{i:V}\divergence^V_g(e_i)\divergence^H_g(e_i) +\tau_{g,V}\\
&= -\qual^V_g +2\tau_{g,V} \;\;,
\end{split} \end{equation}

\begin{equation} \label{switchqualH} \begin{split}
\qual^H_h &= \sum_{i:H}\eps_i\partial_{e_i}\divergence^V_h(e_i) +\sum_{i:H}\eps_i\divergence^H_h(e_i)\divergence^V_h(e_i) +\tau_{h,H}\\
&= \sum_{i:H}\partial_{e_i}\divergence^V_g(e_i) +\sum_{i:H}\divergence^H_g(e_i)\divergence^V_g(e_i) -\tau_{g,H}\\
&= \qual^H_g -2\tau_{g,H} \;\;.
\end{split} \end{equation}
\end{subequations}

\begin{subequations}
The formulae \ref{scalON}, \eqref{switchChristoffel}, \eqref{switchdivergence}, \ref{moreONformulae}, and \ref{sigmatauframe} yield
\begin{equation} \label{switchscalVV} \begin{split}
\scal^{V,V}_h &= -2\sum_{i:V}\eps_i\partial_{e_i}\divergence^V_h(e_i) -\sum_i\eps_i\divergence^V_h(e_i)^2 -\sum_{i,j,k:V}\eps_i\eps_j\eps_kH^k_{ij}H^k_{ji} +\tau_{h,V} -2\sum_{i,k:V}\sum_{j:H}\eps_i\eps_j\eps_kH^k_{ij}H^k_{ji}\\
&= 2\sum_{i:V}\partial_{e_i}\divergence^V_g(e_i) +\sum_{i:V}\divergence^V_g(e_i)^2 -\sum_{i:H}\divergence^V_g(e_i)^2\\
&\mspace{20mu}+\sum_{i,j,k:V}G^k_{ij}G^k_{ji} +\tau_{g,V} -2\sum_{i,k:V}\sum_{j:H}G^j_{ki}(G^k_{ij}+G^j_{ki}+G^i_{jk})\\
&= \Big(2\sum_{i:V}\partial_{e_i}\divergence^V_g(e_i) +\sum_i\divergence^V_g(e_i)^2 +\sum_{i,j,k:V}G^k_{ij}G^k_{ji} -\tau_{g,V} +2\sum_{i,k:V}\sum_{j:H}G^k_{ij}G^k_{ji}\Big)\\
&\mspace{20mu}-2\sum_{i:H}\divergence^V_g(e_i)^2 +2\tau_{g,V} -2\sum_{i,k:V}\sum_{j:H}G^k_{ij}G^k_{ji} -2\sum_{i,k:V}\sum_{j:H}G^j_{ki}(G^k_{ij}+G^j_{ki}+G^i_{jk})\\
&= -\scal^{V,V}_g -2\eval{\divergence^V_g}{\divergence^V_g}_{g,H} +2\tau_{g,V} -2\sum_{i,k:V}\sum_{j:H}G^k_{ij}G^k_{ji}\\
&\mspace{20mu}+2\sum_{i,k:V}\sum_{j:H}G^j_{ki}G^j_{ik} -2\sum_{i,k:V}\sum_{j:H}G^j_{ki}G^j_{ki} +2\sum_{i,k:V}\sum_{j:H}G^i_{kj}G^i_{jk}\\
&= -\scal^{V,V}_g -2\eval{\divergence^V_g}{\divergence^V_g}_{g,H} +4\tau_{g,V} -2\sigma_{g,V} \;\;,
\end{split} \end{equation}

\begin{equation} \label{switchscalHH} \begin{split}
\scal^{H,H}_h &= -2\sum_{i:H}\eps_i\partial_{e_i}\divergence^H_h(e_i) -\sum_i\eps_i\divergence^H_h(e_i)^2 -\sum_{i,j,k:H}\eps_i\eps_j\eps_kH^k_{ij}H^k_{ji} +\tau_{h,H} -2\sum_{i,k:H}\sum_{j:V}\eps_i\eps_j\eps_kH^k_{ij}H^k_{ji}\\
&= -2\sum_{i:H}\partial_{e_i}\divergence^H_g(e_i) +\sum_{i:V}\divergence^H_g(e_i)^2 -\sum_{i:H}\divergence^H_g(e_i)^2\\ &\mspace{20mu}-\sum_{i,j,k:H}G^k_{ij}G^k_{ji} -\tau_{g,H} +2\sum_{i,k:H}\sum_{j:V}G^i_{kj}(G^k_{ji}+G^i_{kj}+G^j_{ik})\\
&= \Big(-2\sum_{i:H}\partial_{e_i}\divergence^H_g(e_i) -\sum_i\divergence^H_g(e_i)^2 -\sum_{i,j,k:H}G^k_{ij}G^k_{ji} +\tau_{g,H} -2\sum_{i,k:H}\sum_{j:V}G^k_{ij}G^k_{ji}\Big)\\
&\mspace{20mu} +2\sum_{i:V}\divergence^H_g(e_i)^2 -2\tau_{g,H} +2\sum_{i,k:H}\sum_{j:V}G^k_{ij}G^k_{ji} +2\sum_{i,k:H}\sum_{j:V}G^i_{kj}(G^k_{ji}+G^i_{kj}+G^j_{ik})\\
&= \scal^{H,H}_g +2\eval{\divergence^H_g}{\divergence^H_g}_{g,V} -2\tau_{g,H} +2\sum_{i,k:H}\sum_{j:V}G^k_{ij}G^k_{ji}\\
&\mspace{20mu}-2\sum_{i,k:H}\sum_{j:V}G^i_{kj}G^i_{jk} +2\sum_{i,k:H}\sum_{j:V}G^i_{kj}G^i_{kj} -2\sum_{i,k:H}\sum_{j:V}G^j_{ki}G^j_{ik}\\
&= \scal^{H,H}_g +2\eval{\divergence^H_g}{\divergence^H_g}_{g,V} +2\sigma_{g,H} -4\tau_{g,H} \;\;.
\end{split} \end{equation}

From \ref{scalqual} and \eqref{switchqualV}, \eqref{switchqualH}, we get
\begin{equation} \label{switchscalVH} \begin{split}
\scal^{V,H}_h = \scal^{H,V}_h &= -(\qual^V_h +\qual^H_h)\\
&= \qual^V_g -2\tau_{g,V} -\qual^H_g +2\tau_{g,H}\\
&= \scal^{V,H}_g +2\qual^V_g -2\tau_{g,V} +2\tau_{g,H} \;\;.
\end{split} \end{equation}
\end{subequations}

We combine \eqref{switchscalVV}, \eqref{switchscalHH}, \eqref{switchscalVH}: \begin{equation} \begin{split}
\scal_h &= \scal^{V,V}_h +\scal^{H,H}_h +2\scal^{V,H}_h\\[1ex]
&= -\scal^{V,V}_g -2\eval{\divergence^V_g}{\divergence^V_g}_{g,H} +4\tau_{g,V} -2\sigma_{g,V} +\scal^{H,H}_g +2\eval{\divergence^H_g}{\divergence^H_g}_{g,V}\\
&\mspace{20mu}+2\sigma_{g,H} -4\tau_{g,H} +2\scal^{V,H}_g +4\qual^V_g -4\tau_{g,V} +4\tau_{g,H}\\[1ex]
&= \scal_g -2\scal^{V,V}_g +4\qual^V_g -2\eval{\divergence^V_g}{\divergence^V_g}_{g,H} +2\eval{\divergence^H_g}{\divergence^H_g}_{g,V} -2\sigma_{g,V} +2\sigma_{g,H} \;\;.
\end{split} \end{equation}

\subsection{Laplacians}

For every $u\in C^\infty(M,\R)$, we calculate (cf.\ \ref{laplacianformulae}, \ref{laplacedivergence}):
\begin{equation} \label{switchlaplace} \begin{split}
\laplace^V_{h,V}(u) &= \sum_{i:V}\eps_i\partial_{e_i}\partial_{e_i}u +\sum_{i:V}\eps_i\divergence^V_h(e_i)du(e_i)\\
&= -\sum_{i:V}\partial_{e_i}\partial_{e_i}u -\sum_{i:V}\divergence^V_g(e_i)du(e_i)\\
&= -\laplace^V_{g,V}(u) \;\;,\\[2ex]
\eval{\divergence^H_h}{du}_{h,V} = \laplace^H_{h,V}(u) &= \sum_{i:V}\eps_i\divergence^H_h(e_i)du(e_i)\\
&= -\sum_{i:V}\divergence^H_g(e_i)du(e_i)\\
&= -\eval{\divergence^H_g}{du}_{g,V} = -\laplace^H_{g,V}(u) \;\;,\\[2ex]
\eval{\divergence^V_h}{du}_{h,H} = \laplace^V_{h,H}(u) &= \sum_{i:H}\eps_i\divergence^V_h(e_i)du(e_i)\\
&= \sum_{i:H}\divergence^V_g(e_i)du(e_i)\\
&= \eval{\divergence^V_g}{du}_{g,H} = \laplace^V_{g,H}(u) \;\;,\\[2ex]
\laplace^H_{h,H}(u) &= \sum_{i:H}\eps_i\partial_{e_i}\partial_{e_i}u +\sum_{i:H}\eps_i\divergence^H_h(e_i)du(e_i)\\
&= \sum_{i:H}\partial_{e_i}\partial_{e_i}u +\sum_{i:H}\divergence^H_g(e_i)du(e_i)\\
&= \laplace^H_{g,H}(u) \;\;. \end{split} \end{equation}

Consequently (cf.\ \ref{laplacedef}, \ref{laplaciansplit}):
\begin{equation} \begin{split}
\laplace_h(u) &= \laplace^V_{h,V}(u) +\laplace^V_{h,H}(u) +\laplace^H_{h,V}(u) +\laplace^H_{h,H}(u)\\
&= -\laplace^V_{g,V}(u) +\laplace^V_{g,H}(u) -\laplace^H_{g,V}(u) +\laplace^H_{g,H}(u)\\
&= \laplace_g(u) -2\laplace^V_{g,V}(u) -2\eval{\divergence^H_g}{du}_{g,V} \;\;.
\end{split} \end{equation}

\pagebreak

\subsection{Summary of the results}

\begin{theorem} \label{switchformulae}
Let $(M,g)$ be a Riemannian manifold, let $V$ be a $q$-plane distribution on $M$, let $h$ be the semi-Riemannian metric $\switch(g,V)$ of index $q$, let $u\in C^\infty(M,\R)$. Then the following formulae hold (where $H$ denotes the $g$-orthogonal distribution of $V$):
\[ \begin{split}
\eval{\divergence^V_h}{\divergence^V_h}_{h,H} &= \eval{\divergence^V_g}{\divergence^V_g}_{g,H} \;\;,\\
\eval{\divergence^H_h}{\divergence^H_h}_{h,V} &= -\eval{\divergence^H_g}{\divergence^H_g}_{g,V} \;\;,\\
\sigma_{h,H} &= -\sigma_{g,H} \;\;,\\
\sigma_{h,V} &= \sigma_{g,V} \;\;,\\
\tau_{h,H} &= -\tau_{g,H} \;\;,\\
\tau_{h,V} &= \tau_{g,V} \;\;,\\
\qual^V_h &= -\qual^V_g +2\tau_{g,V} \;\;,\\
\qual^H_h &= \qual^H_g -2\tau_{g,H} \;\;,\\
\scal^{V,H}_h &= \scal^{V,H}_g +2\qual^V_g -2\tau_{g,V} +2\tau_{g,H} \;\;,\\
\scal^{V,V}_h &= -\scal^{V,V}_g -2\eval{\divergence^V_g}{\divergence^V_g}_{g,H} +4\tau_{g,V} -2\sigma_{g,V} \;\;,\\
\scal^{H,H}_h &= \scal^{H,H}_g +2\eval{\divergence^H_g}{\divergence^H_g}_{g,V} +2\sigma_{g,H} -4\tau_{g,H} \;\;,\\
\scal_h &= \scal_g -2\scal^{V,V}_g +4\qual^V_g -2\eval{\divergence^V_g}{\divergence^V_g}_{g,H} +2\eval{\divergence^H_g}{\divergence^H_g}_{g,V} -2\sigma_{g,V} +2\sigma_{g,H} \;\;,\\
\laplace^V_{h,V}(u) &= -\laplace^V_{g,V}(u) \;\;,\\
\laplace^H_{h,H}(u) &= \laplace^H_{g,H}(u) \;\;,\\
\laplace_h(u) &= \laplace_g(u) -2\laplace^V_{g,V}(u) -2\eval{\divergence^H_g}{du}_{g,V} \;\;,\\
\eval{\divergence^V_h}{du}_{h,H} &= \eval{\divergence^V_g}{du}_{g,H} \;\;,\\
\eval{\divergence^H_h}{du}_{h,V} &= -\eval{\divergence^H_g}{du}_{g,V} \;\;.
\end{split} \]
\end{theorem}
\Proof
This has been proved in the previous subsections.
\end{proof}

\begin{corollary} \label{scalsigninversion}
Let $(M,h)$ be a semi-Riemannian manifold. Then $\scal_{-h} = -\scal_h$.
\end{corollary}
\Proof
This (well-known) statement can easily be verified in a direct way, but let us deduce it from the preceding theorem. By Theorem \ref{baumtheorem}, the tangent bundle $TM$ splits into an $h$-orthogonal sum $V\oplus H$, where $V$ is an $h$-timelike distribution on $M$ and $H$ is an $h$-spacelike distribution on $M$. We consider the Riemannian metric $g\define\switch(h,V)$. Since $h=\switch(g,V)$ and $-h=\switch(g,H)$, we get
\[ \begin{split}
\scal_h &= \scal_g -2\scal^{V,V}_g +4\qual^V_g -2\eval{\divergence^V_g}{\divergence^V_g}_{g,H} +2\eval{\divergence^H_g}{\divergence^H_g}_{g,V} -2\sigma_{g,V} +2\sigma_{g,H} \;\;,\\
\scal_{-h} &= \scal_g -2\scal^{H,H}_g +4\qual^H_g -2\eval{\divergence^H_g}{\divergence^H_g}_{g,V} +2\eval{\divergence^V_g}{\divergence^V_g}_{g,H} -2\sigma_{g,H} +2\sigma_{g,V} \;\;.
\end{split} \]
By \ref{scaldecomposition} and \ref{scalqual}, we have $2\scal_g = 2\scal^{V,V}_g +2\scal^{H,H}_g -4\qual^V_g -4\qual^H_g$. Hence
\[
\scal_g -2\scal^{H,H}_g +4\qual^H_g = -\scal_g +2\scal^{V,V}_g -4\qual^V_g
\]
and thus $\scal_{-h} = -\scal_h$.
\end{proof}

\newcommand{\Triv}{\textbf{0}}
\begin{remark}
As we mentioned in Remark \ref{switchgeneral}, the formulae of Theorem \ref{switchformulae} remain true if $g$ is an arbitrary semi-Riemannian metric (provided we assume that $V$ is $g$-good). Hence we can deduce Corollary \ref{scalsigninversion} also via $-h = \switch(h,TM)$ and the obvious fact that $\scal^{TM,TM}_h = \scal_h$ and $0 = \qual^{TM}_g = \eval{\divergence^{TM}_g}{\divergence^{TM}_g}_{g,\Triv} = \eval{\divergence^\Triv_g}{\divergence^\Triv_g}_{g,TM} = \sigma_{g,TM} = \sigma_{g,\Triv}$, where $\Triv$ denotes the unique $0$-plane distribution on $M$. Note that the latter functions vanish since their definitions involve contractions over this trivial distribution.
\end{remark}


\section{Formulae for stretching}

\emph{Throughout this section, let $(M,g)$ be an $n$-dimensional semi-Riemannian manifold, let $V$ be a $g$-good $q$-plane distribution on $M$, let $f\in C^\infty(M,\R_{>0})$, let $\bar{g}$ denote the semi-Riemannian metric $\stre(g,f,V)$, and let $H$ denote the $g$-orthogonal distribution of $V$ (which is also the $\bar{g}$-orthogonal distribution of the $\bar{g}$-good distribution $V$; cf.\ Remark \ref{stretchtrivialities}).}
\medskip\\
Our aim is to express functions which are defined with respect to the metric $\bar{g}$ --- e.g.\ the functions $\scal_{\bar{g}}$, $\qual^V_{\bar{g}}$, $\sigma_{\bar{g},H}$ --- in terms of functions which are defined with respect to $g$ (e.g.\ $\scal_g$, $\qual^V_g$, $\sigma_{g,H}$). We will summarise the results of our calculations in Theorem \ref{stretchformulae} at the end of this section.
\medskip\\
As usual, we can do the computations locally: By Corollary \ref{adaptedframe}, each point in $M$ has an open neighbourhood which admits a $V$-adapted $g$-orthonormal frame $(e_1,\dots,e_n)$. The local frame $(\bar{e}_1,\dots,\bar{e}_n)$ which is defined by
\[
\bar{e}_i = \begin{cases} f\,e_i &\text{if $i:V$}\\ \phantom{f}\,e_i &\text{if $i:H$}\\ \end{cases}
\]
for all $i\in\set{1,\dots,n}$ is a $V$-adapted $\bar{g}$-orthonormal frame such that $\eps_i \define g(e_i,e_i) = \bar{g}(\bar{e}_i,\bar{e}_i) \in \set{1,-1}$ for all $i\in\set{1,\dots,n}$; cf.\ Remark \ref{stretchONframe}.
\smallskip\\
We denote the \LeviCivita\ connections of $g$ and $\bar{g}$ by $\nabla$ and $\bar{\nabla}$, respectively. The ON Christoffel symbols of $g$ with respect to the ON frame $(e_1,\dots,e_n)$ are denoted by
\[ \begin{split}
\Gamma^k_{ij} &\define g(\nabla_{e_i}e_j,e_k) \;\;,
\end{split} \]
whereas the ON Christoffel symbols of $\bar{g}$ with respect to the ON frame $(\bar{e}_1,\dots,\bar{e}_n)$ are denoted by
\[ \begin{split}
\bar{\Gamma}^k_{ij} &\define \bar{g}(\bar{\nabla}_{\bar{e}_i}\bar{e}_j,\bar{e}_k)
\end{split} \]
for $i,j,k\in\set{1,\dots,n}$. (Recall that the ON Christoffel symbols satisfy the equations from Remark \ref{Christoffel}.)
\smallskip\\
Every local vector field $v$ on $M$ satisfies obviously
\begin{equation} \label{stretchprep}
\bar{g}(v,\bar{e}_i) = \begin{cases} f\bar{g}(v,e_i) = \frac{1}{f}g(v,e_i) &\text{if $i:V$}\\ \phantom{f}\bar{g}(v,e_i) = \phantom{\frac{1}{f}}g(v,e_i) &\text{if $i:H$} \end{cases} \;\;.
\end{equation}

\subsection{The orthonormal Christoffel symbols}

Using Equation \eqref{stretchprep}, we compute the following cases (the intermediate step which occurs in the third line is omitted in the calculations of the lines below it because it should be obvious then):
\[
\bar{g}([\bar{e}_i,\bar{e}_j],\bar{e}_k) = \begin{cases} g([e_i,e_j],e_k) &\text{if $i,j,k:H$}\\
\frac{1}{f}g([e_i,e_j],e_k) &\text{if $i,j:H$, $k:V$}\\
g([e_i,f\,e_j],e_k) = f\,g([e_i,e_j],e_k) +df(e_i)g(e_j,e_k) = f\,g([e_i,e_j],e_k) &\text{if $i,k:H$, $j:V$}\\
g([f\,e_i,e_j],e_k) = f\,g([e_i,e_j],e_k) &\text{if $j,k:H$, $i:V$}\\
\frac{1}{f}g([e_i,f\,e_j],e_k) = g([e_i,e_j],e_k) +\eps_k\delta_{jk}\frac{1}{f}df(e_i) &\text{if $i:H$, $j,k:V$}\\
\frac{1}{f}g([f\,e_i,e_j],e_k) = g([e_i,e_j],e_k) -\eps_k\delta_{ik}\frac{1}{f}df(e_j) &\text{if $j:H$, $i,k:V$}\\
g([f\,e_i,f\,e_j],e_k) = f^2g([e_i,e_j],e_k) &\text{if $k:H$, $i,j:V$}\\
\frac{1}{f}\,g([f\,e_i,f\,e_j],e_k) = f\,g([e_i,e_j],e_k) +\eps_k\delta_{jk}df(e_i) -\eps_k\delta_{ik}df(e_j) &\text{if $i,j,k:V$}\\
\end{cases} \;\;.
\]
(By the antisymmetry of the Lie bracket, the fourth case formula can be read off directly from the third case, and the sixth case formula can be read off from the fifth case.) Note that in the case where $V$ has rank $1$, the seventh and eighth case yield $0$, because $e_i=e_j$ holds then.
\smallskip\\
Hence $\bar{\Gamma}^k_{ij} = \Gamma^k_{ij}$ if $i,j,k:H$. If $i,j:H$ and $k:V$, we get (using $g([e_a,e_b],e_c) = g(\nabla_{e_a}e_b -\nabla_{e_b}e_a,e_c) = \Gamma^c_{ab}-\Gamma^c_{ba}$):
\[\begin{split}
\bar{\Gamma}^k_{ij} &= \frac{1}{2}\Big(\bar{g}([\bar{e}_i,\bar{e}_j],\bar{e}_k) +\bar{g}([\bar{e}_k,\bar{e}_i],\bar{e}_j) +\bar{g}([\bar{e}_k,\bar{e}_j],\bar{e}_i)\Big)\\
&= \frac{1}{2}\Big(\frac{1}{f}g([e_i,e_j],e_k) +f\,g([e_k,e_i],e_j) +f\,g([e_k,e_j],e_i)\Big)\\
&= \frac{1}{2}\Big(\frac{1}{f}(\Gamma^k_{ij}-\Gamma^k_{ji}) +f(\Gamma^j_{ki}-\Gamma^j_{ik} +\Gamma^i_{kj}-\Gamma^i_{jk})\Big)\\
&= \frac{1}{2}\Big(\frac{1}{f}(\Gamma^k_{ij}-\Gamma^k_{ji}) +f(\Gamma^k_{ij}+\Gamma^k_{ji})\Big) \;\;. \end{split} \]

This implies the formula for the case $i,k:H$ and $j:V$: \[
\bar{\Gamma}^k_{ij} = -\bar{\Gamma}^j_{ik} = -\frac{1}{2}\Big(\frac{1}{f}(\Gamma^j_{ik}-\Gamma^j_{ki}) +f(\Gamma^j_{ik}+\Gamma^j_{ki})\Big) \;\;.
\]

In the case $j,k:H$ and $i:V$, we get
\[\begin{split}
\bar{\Gamma}^k_{ij} &= \frac{1}{2}\Big(\bar{g}([\bar{e}_i,\bar{e}_j],\bar{e}_k) +\bar{g}([\bar{e}_k,\bar{e}_i],\bar{e}_j) +\bar{g}([\bar{e}_k,\bar{e}_j],\bar{e}_i)\Big)\\
&= \frac{1}{2}\Big(f\,g([e_i,e_j],e_k) +f\,g([e_k,e_i],e_j) +\frac{1}{f}g([e_k,e_j],e_i)\Big)\\
&= \frac{1}{2}\Big(f(\Gamma^k_{ij}-\Gamma^k_{ji} +\Gamma^j_{ki}-\Gamma^j_{ik}) +\frac{1}{f}(\Gamma^i_{kj}-\Gamma^i_{jk})\Big)\\
&= f\Gamma^k_{ij} +\frac{1}{2}(\frac{1}{f}-f)(\Gamma^i_{kj}-\Gamma^i_{jk}) \;\;.
\end{split} \]

In the case $i:H$ and $j,k:V$, we obtain:
\[ \begin{split}
\bar{\Gamma}^k_{ij} &= \frac{1}{2}\Big(\bar{g}([\bar{e}_i,\bar{e}_j],\bar{e}_k) +\bar{g}([\bar{e}_k,\bar{e}_i],\bar{e}_j) +\bar{g}([\bar{e}_k,\bar{e}_j],\bar{e}_i)\Big)\\
&= \frac{1}{2}\Big(g([e_i,e_j],e_k) +\eps_k\delta_{jk}\frac{1}{f}df(e_i) +g([e_k,e_i],e_j) -\eps_j\delta_{kj}\frac{1}{f}df(e_i) +f^2g([e_k,e_j],e_i)\Big)\\
&= \frac{1}{2}\Big((\Gamma^k_{ij}-\Gamma^k_{ji} +\Gamma^j_{ki}-\Gamma^j_{ik}) +f^2(\Gamma^i_{kj}-\Gamma^i_{jk})\Big)\\ &= \Gamma^k_{ij} -\frac{1}{2}(1-f^2)(\Gamma^i_{kj}-\Gamma^i_{jk}) \;\;;
\end{split} \]
we have used $\eps_k\delta_{jk} = \eps_j\delta_{kj}$ here.
\smallskip\\
In the case $k:H$ and $i,j:V$, we get
\[ \begin{split}
\bar{\Gamma}^k_{ij} &= \frac{1}{2}\Big(\bar{g}([\bar{e}_i,\bar{e}_j],\bar{e}_k) +\bar{g}([\bar{e}_k,\bar{e}_i],\bar{e}_j) +\bar{g}([\bar{e}_k,\bar{e}_j],\bar{e}_i)\Big)\\
&=  \frac{1}{2}\Big(f^2g([e_i,e_j],e_k) +g([e_k,e_i],e_j) +\eps_j\delta_{ij}\frac{1}{f}df(e_k) +g([e_k,e_j],e_i) +\eps_i\delta_{ji}\frac{1}{f}df(e_k)\Big)\\
&= \frac{1}{2}\Big(f^2(\Gamma^k_{ij}-\Gamma^k_{ji}) +(\Gamma^j_{ki}-\Gamma^j_{ik} +\Gamma^i_{kj}-\Gamma^i_{jk})\Big) +\eps_j\delta_{ij}\frac{1}{f}df(e_k)\\
&= \frac{1}{2}\Big((\Gamma^k_{ij}+\Gamma^k_{ji}) +f^2(\Gamma^k_{ij}-\Gamma^k_{ji})\Big) +\eps_j\delta_{ij}\frac{1}{f}df(e_k) \;\;.
\end{split} \]

This implies the formula for the case $j:H$ and $i,k:V$:
\[ \begin{split}
\bar{\Gamma}^k_{ij} &= -\bar{\Gamma}^j_{ik} = -\frac{1}{2}\Big((\Gamma^j_{ik}+\Gamma^j_{ki}) +f^2(\Gamma^j_{ik}-\Gamma^j_{ki})\Big) -\eps_k\delta_{ik}\frac{1}{f}df(e_j) \;\;.
\end{split} \]

In the case $i,j,k:V$, we obtain
\[ \begin{split}
\bar{\Gamma}^k_{ij} &= \frac{1}{2}\Big(\bar{g}([\bar{e}_i,\bar{e}_j],\bar{e}_k) +\bar{g}([\bar{e}_k,\bar{e}_i],\bar{e}_j) +\bar{g}([\bar{e}_k,\bar{e}_j],\bar{e}_i)\Big)\\
&= \frac{1}{2}\Big(f\,g([e_i,e_j],e_k) +\eps_k\delta_{jk}df(e_i) -\eps_k\delta_{ik}df(e_j) +f\,g([e_k,e_i],e_j)\\
&\mspace{45mu} +\eps_j\delta_{ij}df(e_k) -\eps_j\delta_{kj}df(e_i) +f\,g([e_k,e_j],e_i) +\eps_i\delta_{ji}df(e_k) -\eps_i\delta_{ki}df(e_j)\Big)\\
&= f\Gamma^k_{ij} -\eps_i\delta_{ik}df(e_j) +\eps_i\delta_{ij}df(e_k) \;\;.
\end{split} \]

To summarise, the orthonormal Christoffel symbols of the metric $\bar{g} = \stre(g,f,V)$ are given by
\begin{center} \fbox{\parbox{160mm}{
\begin{equation} \label{stretchChristoffel}
\bar{\Gamma}^k_{ij} = \begin{cases}
\Gamma^k_{ij} &\text{if $i,j,k:H$}\\[1ex]
\frac{1}{2}\Big(\frac{1}{f}(\Gamma^k_{ij}-\Gamma^k_{ji}) +f(\Gamma^k_{ij}+\Gamma^k_{ji})\Big) &\text{if $i,j:H$, $k:V$}\\[1ex]
-\frac{1}{2}\Big(\frac{1}{f}(\Gamma^j_{ik}-\Gamma^j_{ki}) +f(\Gamma^j_{ik}+\Gamma^j_{ki})\Big) &\text{if $i,k:H$, $j:V$}\\[1ex]
f\Gamma^k_{ij} +\frac{1}{2}(\frac{1}{f}-f)(\Gamma^i_{kj}-\Gamma^i_{jk}) &\text{if $j,k:H$, $i:V$}\\[1ex]
\Gamma^k_{ij} -\frac{1}{2}(1-f^2)(\Gamma^i_{kj}-\Gamma^i_{jk}) &\text{if $i:H$, $j,k:V$}\\[1ex]
-\frac{1}{2}\Big((\Gamma^j_{ik}+\Gamma^j_{ki}) +f^2(\Gamma^j_{ik}-\Gamma^j_{ki})\Big) -\eps_k\delta_{ik}\frac{1}{f}df(e_j) &\text{if $j:H$, $i,k:V$}\\[1ex]
\frac{1}{2}\Big((\Gamma^k_{ij}+\Gamma^k_{ji}) +f^2(\Gamma^k_{ij}-\Gamma^k_{ji})\Big) +\eps_j\delta_{ij}\frac{1}{f}df(e_k) &\text{if $k:H$, $i,j:V$}\\[1ex]
f\Gamma^k_{ij} -\eps_i\delta_{ik}df(e_j) +\eps_i\delta_{ij}df(e_k) &\text{if $i,j,k:V$}
\end{cases} \;\;.
\end{equation}}}
\end{center}

\subsection{Divergences, $\sigma$s, $\tau$s}

\begin{subequations} \label{threefourteen}
For $i:H$, we infer from \eqref{stretchChristoffel} and \ref{divUformulae}:
\begin{equation} \label{threefourteena} \begin{split}
\divergence^V_{\bar{g}}(\bar{e}_i) &= -\sum_{k:V}\eps_k\bar{\Gamma}^i_{kk}\\
&= -\sum_{k:V}\eps_k\bigg(\frac{1}{2}\Big((\Gamma^i_{kk}+\Gamma^i_{kk}) +f^2(\Gamma^i_{kk}-\Gamma^i_{kk})\Big) +\eps_k\delta_{kk}\frac{1}{f}df(e_i)\bigg)\\
&= \divergence^V_g(e_i) -\frac{q}{f}df(e_i) \;\;; \end{split} \end{equation}
moreover,
\begin{equation} \label{threefourteenb}
\divergence^H_{\bar{g}}(\bar{e}_i) = \sum_{k:H}\eps_k\bar{\Gamma}^k_{ki} = \sum_{k:H}\eps_k\Gamma^k_{ki} = \divergence^H_g(e_i) \;\;.
\end{equation}

For $i:V$, we get
\begin{equation} \label{threefourteenc} \begin{split}
\divergence^V_{\bar{g}}(\bar{e}_i) &= \sum_{k:V}\eps_k\bar{\Gamma}^k_{ki} = \sum_{k:V}\eps_k(f\Gamma^k_{ki} -\eps_k\delta_{kk}df(e_i) +\eps_k\delta_{ki}df(e_k))\\
&= f\divergence^V_g(e_i) -(q-1)df(e_i)
\end{split} \end{equation}
and
\begin{equation} \label{threefourteend} \begin{split}
\divergence^H_{\bar{g}}(\bar{e}_i) &= -\sum_{k:H}\eps_k\bar{\Gamma}^i_{kk} = -\sum_{k:H}\frac{\eps_k}{2}\Big(\frac{1}{f}(\Gamma^i_{kk}-\Gamma^i_{kk}) +f(\Gamma^i_{kk}+\Gamma^i_{kk})\Big)\\
&= f\divergence^H_g(e_i) \;\;.
\end{split} \end{equation}
\end{subequations}

From \eqref{stretchChristoffel}, we infer
\begin{subequations} \label{threefifteen}
\begin{equation} \label{threefifteena} \begin{split}
\sigma_{\bar{g},H} &= \sum_{i,j:H}\sum_{k:V}\eps_i\eps_j\eps_k(\bar{\Gamma}^k_{ij})^2\\[0.9ex]
&= \frac{1}{4}\sum_{i,j:H}\sum_{k:V}\eps_i\eps_j\eps_k \Big(\frac{1}{f}(\Gamma^k_{ij}-\Gamma^k_{ji}) +f(\Gamma^k_{ij}+\Gamma^k_{ji})\Big)^2\\[0.9ex]
&= \frac{1}{4}\sum_{i,j:H}\sum_{k:V}\eps_i\eps_j\eps_k\Big((f+\frac{1}{f})\Gamma^k_{ij} +(f-\frac{1}{f})\Gamma^k_{ji}\Big)^2\\[0.9ex]
&= \frac{1}{4}\sum_{i,j:H}\sum_{k:V}\eps_i\eps_j\eps_k \Big((f^2+\frac{1}{f^2}+2)(\Gamma^k_{ij})^2 +(f^2+\frac{1}{f^2}-2)(\Gamma^k_{ji})^2 +2(f^2-\frac{1}{f^2})\Gamma^k_{ij}\Gamma^k_{ji}\Big)\\[0.9ex]
&= \frac{1}{2}(f^2+\frac{1}{f^2})\sum_{i,j:H}\sum_{k:V}\eps_i\eps_j\eps_k(\Gamma^k_{ij})^2 +\frac{1}{2}(f^2-\frac{1}{f^2})\sum_{i,j:H}\sum_{k:V}\eps_i\eps_j\eps_k \Gamma^k_{ij}\Gamma^k_{ji}\\[0.9ex]
&= \frac{1}{2}(f^2+\frac{1}{f^2})\sigma_{g,H} +\frac{1}{2}(f^2-\frac{1}{f^2})\tau_{g,H} \;\;,
\end{split} \end{equation}

\begin{equation} \label{threefifteenb} \begin{split}
\tau_{\bar{g},H} &= \sum_{i,j:H}\sum_{k:V}\eps_i\eps_j\eps_k\bar{\Gamma}^k_{ij}\bar{\Gamma}^k_{ji}\\[0.9ex]
&= \frac{1}{4}\sum_{i,j:H}\sum_{k:V}\eps_i\eps_j\eps_k \Big(\frac{1}{f}(\Gamma^k_{ij}-\Gamma^k_{ji}) +f(\Gamma^k_{ij}+\Gamma^k_{ji})\Big) \Big(\frac{1}{f}(\Gamma^k_{ji}-\Gamma^k_{ij}) +f(\Gamma^k_{ji}+\Gamma^k_{ij})\Big)\\[0.9ex]
&= \frac{1}{4}\sum_{i,j:H}\sum_{k:V}\eps_i\eps_j\eps_k \Big(f^2(\Gamma^k_{ij}+\Gamma^k_{ji})^2 -\frac{1}{f^2}(\Gamma^k_{ij}-\Gamma^k_{ji})^2\Big)\\[0.9ex]
&= \frac{1}{4}\sum_{i,j:H}\sum_{k:V}\eps_i\eps_j\eps_k \Big((f^2-\frac{1}{f^2})(\Gamma^k_{ij})^2 +(f^2-\frac{1}{f^2})(\Gamma^k_{ji})^2 +2(f^2+\frac{1}{f^2})\Gamma^k_{ij}\Gamma^k_{ji}\Big)\\[0.91ex]
&= \frac{1}{2}(f^2+\frac{1}{f^2})\tau_{g,H} +\frac{1}{2}(f^2-\frac{1}{f^2})\sigma_{g,H} \;\;,
\end{split} \end{equation}

\begin{equation} \label{threefifteenc} \begin{split}
\sigma_{\bar{g},V} &= \sum_{i,j:V}\sum_{k:H}\eps_i\eps_j\eps_k(\bar{\Gamma}^k_{ij})^2\\[0.9ex]
&= \sum_{i,j:V}\sum_{k:H}\eps_i\eps_j\eps_k \bigg(\frac{1}{2}\Big((\Gamma^k_{ij}+\Gamma^k_{ji}) +f^2(\Gamma^k_{ij}-\Gamma^k_{ji})\Big) +\eps_j\delta_{ij}\frac{1}{f}df(e_k)\bigg)^2\\[0.9ex]
&= \frac{1}{4}\sum_{i,j:V}\sum_{k:H}\eps_i\eps_j\eps_k\Big((\Gamma^k_{ij}+\Gamma^k_{ji}) +f^2(\Gamma^k_{ij}-\Gamma^k_{ji})\Big)^2 +\frac{1}{f^2}\sum_{i,j:V}\sum_{k:H}\eps_i\eps_j\eps_k\delta_{ij}df(e_k)^2\\[0.9ex]
&\mspace{20mu}+\frac{1}{f}\sum_{i,j:V}\sum_{k:H}\eps_i\eps_j\eps_k \Big((\Gamma^k_{ij}+\Gamma^k_{ji}) +f^2(\Gamma^k_{ij}-\Gamma^k_{ji})\Big)\eps_j\delta_{ij}df(e_k)\\[0.9ex]
&= \frac{1}{4}\sum_{i,j:V}\sum_{k:H}\eps_i\eps_j\eps_k\Big((1+f^2)^2(\Gamma^k_{ij})^2 +(1-f^2)^2(\Gamma^k_{ji})^2 +2(1-f^4)\Gamma^k_{ij}\Gamma^k_{ji}\Big)\\[0.9ex]
&\mspace{20mu}+\frac{q}{f^2}\sum_{k:H}\eps_k df(e_k)^2 +\frac{2}{f}\sum_{i:V}\sum_{k:H}\eps_i\eps_k\Gamma^k_{ii}df(e_k)\\[0.9ex]
&= \frac{1}{2}(1+f^4)\sigma_{g,V} +\frac{1}{2}(1-f^4)\tau_{g,V} +\frac{q}{f^2}\eval{df}{df}_{g,H} -\frac{2}{f}\eval{\divergence^V_g}{df}_{g,H} \;\;.
\end{split} \end{equation}

Moreover, using $\eps_j\delta_{ij} = \eps_i\delta_{ji}$, we get:
\begin{equation} \label{threefifteend} \begin{split}
\tau_{\bar{g},V} &= \sum_{i,j:V}\sum_{k:H}\eps_i\eps_j\eps_k\bar{\Gamma}^k_{ij}\bar{\Gamma}^k_{ji}\\
&= \sum_{i,j:V}\sum_{k:H}\eps_i\eps_j\eps_k \bigg(\frac{1}{2}\Big((\Gamma^k_{ij}+\Gamma^k_{ji}) +f^2(\Gamma^k_{ij}-\Gamma^k_{ji})\Big) +\eps_j\delta_{ij}\frac{1}{f}df(e_k)\bigg)\cdot\\
&\mspace{135mu}\bigg(\frac{1}{2}\Big((\Gamma^k_{ji}+\Gamma^k_{ij}) +f^2(\Gamma^k_{ji}-\Gamma^k_{ij})\Big) +\eps_i\delta_{ji}\frac{1}{f}df(e_k)\bigg)\\
&= \frac{1}{4}\sum_{i,j:V}\sum_{k:H}\eps_i\eps_j\eps_k\Big((\Gamma^k_{ij}+\Gamma^k_{ji}) +f^2(\Gamma^k_{ij}-\Gamma^k_{ji})\Big)\Big((\Gamma^k_{ji}+\Gamma^k_{ij}) +f^2(\Gamma^k_{ji}-\Gamma^k_{ij})\Big)\\
&\mspace{20mu}+\frac{1}{f^2}\sum_{i,j:V}\sum_{k:H}\eps_i\eps_j\eps_k \delta_{ij}df(e_k)^2\\
&\mspace{20mu}+\frac{1}{2f}\sum_{i,j:V}\sum_{k:H}\eps_i\eps_j\eps_k \Big((\Gamma^k_{ij}+\Gamma^k_{ji}) +f^2(\Gamma^k_{ij}-\Gamma^k_{ji}) +(\Gamma^k_{ji}+\Gamma^k_{ij}) +f^2(\Gamma^k_{ji}-\Gamma^k_{ij})\Big)\eps_j\delta_{ij}df(e_k)\\
&= \frac{1}{4}\sum_{i,j:V}\sum_{k:H}\eps_i\eps_j\eps_k\Big((\Gamma^k_{ij}+\Gamma^k_{ji})^2 -f^4(\Gamma^k_{ij}-\Gamma^k_{ji})^2\Big) +\frac{q}{f^2}\sum_{k:H}\eps_kdf(e_k)^2 +\frac{2}{f}\sum_{i:V}\sum_{k:H}\eps_i\eps_k\Gamma^k_{ii}df(e_k)\\
&= \frac{1}{2}(1+f^4)\tau_{g,V} +\frac{1}{2}(1-f^4)\sigma_{g,V} +\frac{q}{f^2}\eval{df}{df}_{g,H} -\frac{2}{f}\eval{\divergence^V_g}{df}_{g,H} \;\;.
\end{split} \end{equation}
\end{subequations}

\subsection{Laplacians}

For every $u\in C^\infty(M,\R)$, we calculate (cf.\ \ref{laplacianformulae}, \ref{laplacedivergence}):
\begin{subequations} \label{threesixteen}
\begin{equation} \label{threesixteena} \begin{split}
\laplace^V_{\bar{g},V}(u) &= \sum_{i:V}\eps_i\partial_{\bar{e}_i}\partial_{\bar{e}_i}u +\sum_{i:V}\eps_i\divergence^V_{\bar{g}}(\bar{e}_i)du(\bar{e}_i)\\
&= \sum_{i:V}\eps_if\partial_{e_i}(f\partial_{e_i}u) +\sum_{i:V}\eps_i\Big(f\divergence^V_g(e_i)-(q-1)df(e_i)\Big)f\,du(e_i)\\
&= f^2\sum_{i:V}\eps_i\partial_{e_i}\partial_{e_i}u +f\eval{df}{du}_{g,V} +f^2\sum_{i:V}\eps_i\divergence^V_g(e_i)du(e_i) -(q-1)f\eval{df}{du}_{g,V}\\
&= f^2\laplace^V_{g,V}(u) -(q-2)f\eval{df}{du}_{g,V} \;\;,
\end{split} \end{equation}

\begin{equation} \label{threesixteenb} \begin{split}
\laplace^H_{\bar{g},H}(u) = \sum_{i:H}\eps_i\partial_{\bar{e}_i}\partial_{\bar{e}_i}u +\sum_{i:H}\eps_i\divergence^H_{\bar{g}}(\bar{e}_i)du(\bar{e}_i) = \sum_{i:H}\eps_i\partial_{e_i}\partial_{e_i}u +\sum_{i:H}\eps_i\divergence^H_g(e_i)du(e_i) = \laplace^H_{g,H}(u) \;\;,
\end{split} \end{equation}

\begin{equation} \label{threesixteenc} \begin{split}
\eval{\divergence^H_{\bar{g}}}{du}_{\bar{g},V} = \laplace^H_{\bar{g},V}(u) &= \sum_{i:V}\eps_i\divergence^H_{\bar{g}}(\bar{e}_i)du(\bar{e}_i) = \sum_{i:V}\eps_if\divergence^H_g(e_i)f\,du(e_i)\\
&= f^2\eval{\divergence^H_g}{du}_{g,V} = f^2\laplace^H_{g,V}(u) \;\;,
\end{split} \end{equation}

\begin{equation} \label{threesixteend} \begin{split}
\eval{\divergence^V_{\bar{g}}}{du}_{\bar{g},H} = \laplace^V_{\bar{g},H}(u) &= \sum_{i:H}\eps_i\divergence^V_{\bar{g}}(\bar{e}_i)du(\bar{e}_i)\\
&= \sum_{i:H}\eps_i\Big(\divergence^V_g(e_i)-\frac{q}{f}df(e_i)\Big)du(e_i)\\
&= \eval{\divergence^V_g}{du}_{g,H} -\frac{q}{f}\eval{df}{du}_{g,H} = \laplace^V_{g,H}(u) -\frac{q}{f}\eval{df}{du}_{g,H}\;\;.
\end{split} \end{equation}
\end{subequations}

Hence \begin{equation} \label{threeseventeen} \begin{split}
\laplace_{\bar{g}}(u) &= \laplace^V_{\bar{g},V}(u) +\laplace^H_{\bar{g},H}(u) +\laplace^H_{\bar{g},V}(u) +\laplace^V_{\bar{g},H}(u)\\[1ex]
&= f^2\laplace^V_{g,V}(u) -(q-2)f\eval{df}{du}_{g,V} +\laplace^H_{g,H}(u)\\
&\mspace{20mu}+f^2\eval{\divergence^H_g}{du}_{g,V} +\eval{\divergence^V_g}{du}_{g,H} -\frac{q}{f}\eval{df}{du}_{g,H} \;\;.
\end{split} \end{equation}

\subsection{Scalar curvatures and qualar curvatures}

Now we compute the four functions $\sum_{j:A}\eps_j\divergence^B_{\bar{g}}(\bar{e}_j)^2$, where $A,B\in\set{V,H}$. Note that such a function has no invariant meaning (i.e.\ it depends on the choice of orthonormal frame) if $A=B$. However, these functions appear as summands in the scalar curvature formulae \ref{curvatureON}, so we compute them here as a preparation for the scalar curvatures.
\smallskip\\
By \eqref{threefourteen}, we get
\begin{equation} \label{threeeighteen} \begin{split}
\sum_{i:V}\eps_i\divergence^V_{\bar{g}}(\bar{e}_i)^2 &= \sum_{i:V}\eps_i\Big(f\divergence^V_g(e_i) -(q-1)df(e_i)\Big)^2\\
&= f^2\sum_{i:V}\eps_i\divergence^V_g(e_i)^2 +(q-1)^2\eval{df}{df}_{g,V}\\
&\mspace{20mu}-2(q-1)f\sum_{i:V}\eps_i\divergence^V_g(e_i)df(e_i) \;\;,\\[1.8ex]
\eval{\divergence^V_{\bar{g}}}{\divergence^V_{\bar{g}}}_{\bar{g},H} = \sum_{i:H}\eps_i\divergence^V_{\bar{g}}(\bar{e}_i)^2 &= \sum_{i:H}\eps_i\Big(\divergence^V_g(e_i) -\frac{q}{f}df(e_i)\Big)^2\\
&= \eval{\divergence^V_g}{\divergence^V_g}_{g,H} +\frac{q^2}{f^2}\eval{df}{df}_{g,H} -\frac{2q}{f}\eval{\divergence^V_g}{df}_{g,H} \;\;,\\[1.8ex]
\eval{\divergence^H_{\bar{g}}}{\divergence^H_{\bar{g}}}_{\bar{g},V} = \sum_{i:V}\eps_i\divergence^H_{\bar{g}}(\bar{e}_i)^2 &= f^2\sum_{i:V}\eps_i\divergence^H_g(e_i)^2 = f^2\eval{\divergence^H_g}{\divergence^H_g}_{g,V} \;\;,\\[1.8ex]
\sum_{i:H}\eps_i\divergence^H_{\bar{g}}(\bar{e}_i)^2 &= \sum_{i:H}\eps_i\divergence^H_g(e_i)^2 \;\;.
\end{split} \end{equation}

We continue our preparations for the scalar curvature calculations by computing the following terms (via \eqref{threefourteen}):
\begin{equation} \label{threenineteen} \begin{split}
\sum_{i:V}\eps_i\partial_{\bar{e}_i}\divergence^V_{\bar{g}}(\bar{e}_i) &= f\sum_{i:V}\eps_i\partial_{e_i}\Big(f\divergence^V_g(e_i) -(q-1)df(e_i)\Big)\\
&= f^2\sum_{i:V}\eps_i\partial_{e_i}\divergence^V_g(e_i) +f\sum_{i:V}\eps_i\divergence^V_g(e_i)df(e_i) -(q-1)f\sum_{i:V}\eps_i\partial_{e_i}\partial_{e_i}f \;\;,\\[1.8ex]
\sum_{i:V}\eps_i\partial_{\bar{e}_i}\divergence^H_{\bar{g}}(\bar{e}_i) &= f\sum_{i:V}\eps_i\partial_{e_i}\Big(f\divergence^H_g(e_i)\Big) = f^2\sum_{i:V}\eps_i\partial_{e_i}\divergence^H_g(e_i) +f\eval{\divergence^H_g}{df}_{g,V} \;\;,\\[1.8ex]
\sum_{i:H}\eps_i\partial_{\bar{e}_i}\divergence^V_{\bar{g}}(\bar{e}_i) &= \sum_{i:H}\eps_i\partial_{e_i}\Big(\divergence^V_g(e_i) -\frac{q}{f}df(e_i)\Big)\\
&= \sum_{i:H}\eps_i\partial_{e_i}\divergence^V_g(e_i) -\frac{q}{f}\sum_{i:H}\eps_i\partial_{e_i}\partial_{e_i}f +\frac{q}{f^2}\eval{df}{df}_{g,H} \;\;,\\[1.8ex]
\sum_{i:H}\eps_i\partial_{\bar{e}_i}\divergence^H_{\bar{g}}(\bar{e}_i) &= \sum_{i:H}\eps_i\partial_{e_i}\divergence^H_g(e_i) \;\;.
\end{split} \end{equation}

From \ref{qON}, \eqref{threenineteen}, \eqref{threefourteenc}, \eqref{threefourteend}, \eqref{threefifteend}, we obtain:
\begin{subequations}
\begin{equation} \begin{split}
\qual^V_{\bar{g}} &= \sum_{i:V}\eps_i\partial_{\bar{e}_i}\divergence^H_{\bar{g}}(\bar{e}_i) +\sum_{i:V}\eps_i\divergence^V_{\bar{g}}(\bar{e}_i)\divergence^H_{\bar{g}}(\bar{e}_i) +\tau_{\bar{g},V}\\[1.4ex]
&= f^2\sum_{i:V}\eps_i\partial_{e_i}\divergence^H_g(e_i) +f\eval{\divergence^H_g}{df}_{g,V} +\sum_{i:V}\eps_i\Big(f\divergence^V_g(e_i) -(q-1)df(e_i)\Big)f\divergence^H_g(e_i)\\[0.8ex]
&\mspace{20mu}+\frac{1}{2}(1+f^4)\tau_{g,V} +\frac{1}{2}(1-f^4)\sigma_{g,V} +\frac{q}{f^2}\eval{df}{df}_{g,H} -\frac{2}{f}\eval{\divergence^V_g}{df}_{g,H}\\[1.4ex]
&= f^2\qual^V_g -(q-2)f\eval{\divergence^H_g}{df}_{g,V}\\[0.8ex]
&\mspace{20mu}+\frac{1}{2}(1+f^4-2f^2)\tau_{g,V} +\frac{1}{2}(1-f^4)\sigma_{g,V} +\frac{q}{f^2}\eval{df}{df}_{g,H} -\frac{2}{f}\eval{\divergence^V_g}{df}_{g,H} \;\;.\\[1ex]
\end{split} \end{equation}

By \ref{qON}, \eqref{threenineteen}, \eqref{threefourteena}, \eqref{threefourteenb}, \eqref{threefifteenb}, \ref{laplacianformulae}, we get:
\begin{equation} \begin{split}
\qual^H_{\bar{g}} &= \sum_{i:H}\eps_i\partial_{\bar{e}_i}\divergence^V_{\bar{g}}(\bar{e}_i) +\sum_{i:H}\eps_i\divergence^H_{\bar{g}}(\bar{e}_i)\divergence^V_{\bar{g}}(\bar{e}_i) +\tau_{\bar{g},H}\\[1.4ex]
&= \sum_{i:H}\eps_i\partial_{e_i}\divergence^V_g(e_i) -\frac{q}{f}\sum_{i:H}\eps_i\partial_{e_i}\partial_{e_i}f +\frac{q}{f^2}\eval{df}{df}_{g,H}\\[0.8ex]
&\mspace{20mu}+\sum_{i:H}\eps_i\divergence^H_g(e_i)\Big(\divergence^V_g(e_i) -\frac{q}{f}df(e_i)\Big) +\frac{1}{2}(f^2+\frac{1}{f^2})\tau_{g,H} +\frac{1}{2}(f^2-\frac{1}{f^2})\sigma_{g,H}\\[1.4ex]
&= \qual^H_g -\frac{q}{f}\laplace^H_{g,H}(f) +\frac{q}{f^2}\eval{df}{df}_{g,H} +\frac{1}{2}(f^2+\frac{1}{f^2}-2)\tau_{g,H} +\frac{1}{2}(f^2-\frac{1}{f^2})\sigma_{g,H} \;\;.\\[1ex]
\end{split} \end{equation}
\end{subequations}

The following terms do also appear in the formulae for the scalar curvatures, so we compute them as a preparation. \eqref{stretchChristoffel} yields:
\begin{subequations} \label{threetwentyone}
\begin{equation} \label{threetwentyonea} \begin{split}
\sum_{i,j,k:V}\eps_i\eps_j\eps_k\bar{\Gamma}^k_{ij}\bar{\Gamma}^k_{ji} &= \sum_{i,j,k:V}\eps_i\eps_j\eps_k\Big(f\Gamma^k_{ij} -\eps_i\delta_{ik}df(e_j) +\eps_i\delta_{ij}df(e_k)\Big)\Big(f\Gamma^k_{ji} -\eps_j\delta_{jk}df(e_i) +\eps_j\delta_{ji}df(e_k)\Big)\\[1.4ex]
&= f^2\sum_{i,j,k:V}\eps_i\eps_j\eps_k\Gamma^k_{ij}\Gamma^k_{ji} -f\sum_{i,j,k:V}\eps_i\eps_k\Gamma^k_{ij}\delta_{jk}df(e_i) +f\sum_{i,j,k:V}\eps_i\eps_k\Gamma^k_{ij}\delta_{ji}df(e_k)\\[0.8ex]
&\mspace{20mu}-f\sum_{i,j,k:V}\eps_j\eps_k\delta_{ik}df(e_j)\Gamma^k_{ji} +\sum_{i,j,k:V}\eps_k\delta_{ik}df(e_j)\delta_{jk}df(e_i) -\sum_{i,j,k:V}\eps_k\delta_{ik}df(e_j)\delta_{ji}df(e_k)\\[0.8ex]
&\mspace{20mu}+f\sum_{i,j,k:V}\eps_j\eps_k\delta_{ij}df(e_k)\Gamma^k_{ji} -\sum_{i,j,k:V}\eps_k\delta_{ij}df(e_k)\delta_{jk}df(e_i) +\sum_{i,j,k:V}\eps_k\delta_{ij}df(e_k)\delta_{ji}df(e_k)\\[1.4ex]
&= f^2\sum_{i,j,k:V}\eps_i\eps_j\eps_k\Gamma^k_{ij}\Gamma^k_{ji} +0 -f\sum_{i,k:V}\eps_i\eps_k\Gamma^i_{ik}df(e_k)\\[0.8ex]
&\mspace{20mu}+0 +\sum_{k:V}\eps_kdf(e_k)df(e_k) -\sum_{k:V}\eps_kdf(e_k)df(e_k)\\[0.8ex]
&\mspace{20mu}-f\sum_{i,k:V}\eps_i\eps_kdf(e_k)\Gamma^i_{ik} -\sum_{k:V}\eps_kdf(e_k)df(e_k) +q\sum_{k:V}\eps_kdf(e_k)df(e_k)\\[1.4ex]
&= f^2\sum_{i,j,k:V}\eps_i\eps_j\eps_k\Gamma^k_{ij}\Gamma^k_{ji} -2f\sum_{k:V}\eps_k\divergence^V_g(e_k)df(e_k) +(q-1)\eval{df}{df}_{g,V} \;\;.
\end{split} \end{equation}

From \eqref{stretchChristoffel}, \ref{sigmatauframe}, and by some obvious exchanges of the summation indices $i,k$, we obtain:
\begin{equation} \label{threetwentyoneb} \begin{split}
&\sum_{i,k:V}\sum_{j:H}\eps_i\eps_j\eps_k\bar{\Gamma}^k_{ij}\bar{\Gamma}^k_{ji}\\
&= \sum_{i,k:V}\sum_{j:H}\eps_i\eps_j\eps_k \Big(-\frac{1}{2}\Big((\Gamma^j_{ik}+\Gamma^j_{ki}) +f^2(\Gamma^j_{ik}-\Gamma^j_{ki})\Big) -\eps_k\delta_{ik}\frac{1}{f}df(e_j)\Big)\Big(\Gamma^k_{ji} -\frac{1}{2}(1-f^2)(\Gamma^j_{ki}-\Gamma^j_{ik})\Big)\\
&= -\frac{1}{2}\sum_{i,k:V}\sum_{j:H}\eps_i\eps_j\eps_k (\Gamma^j_{ik}+\Gamma^j_{ki})\Gamma^k_{ji} +\frac{1}{4}(1-f^2)\sum_{i,k:V}\sum_{j:H}\eps_i\eps_j\eps_k (\Gamma^j_{ik}+\Gamma^j_{ki})(\Gamma^j_{ki}-\Gamma^j_{ik})\\
&\mspace{20mu}-\frac{f^2}{2}\sum_{i,k:V}\sum_{j:H}\eps_i\eps_j\eps_k (\Gamma^j_{ik}-\Gamma^j_{ki})\Gamma^k_{ji} +(1-f^2)\frac{f^2}{4}\sum_{i,k:V}\sum_{j:H}\eps_i\eps_j\eps_k (\Gamma^j_{ik}-\Gamma^j_{ki})(\Gamma^j_{ki}-\Gamma^j_{ik})\\
&\mspace{20mu}-\frac{1}{f}\sum_{i,k:V}\sum_{j:H}\eps_i\eps_j \delta_{ik}df(e_j)\Gamma^k_{ji} +\frac{1}{2f}(1-f^2)\sum_{i,k:V}\sum_{j:H}\eps_i\eps_j \delta_{ik}df(e_j)(\Gamma^j_{ki}-\Gamma^j_{ik})\\
&= \frac{1}{2}\sum_{i,k:V}\sum_{j:H}\eps_i\eps_j\eps_k(\Gamma^k_{ij}\Gamma^k_{ji} -\Gamma^i_{kj}\Gamma^i_{jk}) +\frac{1}{4}(1-f^2)\sum_{i,k:V}\sum_{j:H}\eps_i\eps_j\eps_k((\Gamma^j_{ki})^2 -(\Gamma^j_{ik})^2)\\
&\mspace{20mu}+\frac{f^2}{2}\sum_{i,k:V}\sum_{j:H}\eps_i\eps_j\eps_k (\Gamma^k_{ij}\Gamma^k_{ji} +\Gamma^i_{kj}\Gamma^i_{jk}) -(1-f^2)\frac{f^2}{4}\sum_{i,k:V}\sum_{j:H}\eps_i\eps_j\eps_k((\Gamma^j_{ik})^2 +(\Gamma^j_{ki})^2 -2\Gamma^j_{ik}\Gamma^j_{ki}) -0 +0\\
&= f^2\sum_{i,k:V}\sum_{j:H}\eps_i\eps_j\eps_k\Gamma^k_{ij}\Gamma^k_{ji} -\frac{1}{2}f^2(1-f^2)\sigma_{g,V} +\frac{1}{2}f^2(1-f^2)\tau_{g,V} \;\;,
\end{split} \end{equation}
\begin{equation} \label{threetwentyonec} \begin{split}
\sum_{i,k:H}\sum_{j:V}\eps_i\eps_j\eps_k\bar{\Gamma}^k_{ij}\bar{\Gamma}^k_{ji} &= -\frac{1}{2}\sum_{i,k:H}\sum_{j:V}\eps_i\eps_j\eps_k \Big(\frac{1}{f}(\Gamma^j_{ik}-\Gamma^j_{ki}) +f(\Gamma^j_{ik}+\Gamma^j_{ki})\Big)\Big(f\Gamma^k_{ji} +\frac{1}{2}(\frac{1}{f}-f)(\Gamma^j_{ki}-\Gamma^j_{ik})\Big)\\
&= -\frac{1}{2}\sum_{i,k:H}\sum_{j:V}\eps_i\eps_j\eps_k (\Gamma^j_{ik}-\Gamma^j_{ki})\Gamma^k_{ji} -\frac{1}{4f}(\frac{1}{f}-f)\!\!\sum_{i,k:H}\sum_{j:V}\eps_i\eps_j\eps_k (\Gamma^j_{ik}-\Gamma^j_{ki})(\Gamma^j_{ki}-\Gamma^j_{ik})\\
&\mspace{20mu} -\frac{f^2}{2}\!\sum_{i,k:H}\sum_{j:V}\eps_i\eps_j\eps_k (\Gamma^j_{ik}+\Gamma^j_{ki})\Gamma^k_{ji} -\frac{f}{4}(\frac{1}{f}-f)\!\!\sum_{i,k:H}\sum_{j:V}\eps_i\eps_j\eps_k (\Gamma^j_{ik}+\Gamma^j_{ki})(\Gamma^j_{ki}-\Gamma^j_{ik})\\
&= \frac{1}{2}\sum_{i,k:H}\sum_{j:V}\eps_i\eps_j\eps_k(\Gamma^k_{ij}\Gamma^k_{ji} +\Gamma^i_{kj}\Gamma^i_{jk})\\
&\mspace{20mu}+\frac{1}{4f}(\frac{1}{f}-f)\sum_{i,k:H}\sum_{j:V}\eps_i\eps_j\eps_k ((\Gamma^j_{ik})^2 +(\Gamma^j_{ki})^2 -2\Gamma^j_{ik}\Gamma^j_{ki})\\
&\mspace{20mu} +\frac{f^2}{2}\sum_{i,k:H}\sum_{j:V}\eps_i\eps_j\eps_k(\Gamma^k_{ij}\Gamma^k_{ji} -\Gamma^i_{kj}\Gamma^i_{jk}) -\frac{f}{4}(\frac{1}{f}-f)\sum_{i,k:H}\sum_{j:V}\eps_i\eps_j\eps_k((\Gamma^j_{ki})^2 -(\Gamma^j_{ik})^2)\\
&= \sum_{i,k:H}\sum_{j:V}\eps_i\eps_j\eps_k\Gamma^k_{ij}\Gamma^k_{ji} +\frac{1}{2f^2}(1-f^2)\sigma_{g,H} -\frac{1}{2f^2}(1-f^2)\tau_{g,H} \;\;,
\end{split} \end{equation}
\begin{equation} \label{threetwentyoned} \begin{split}
\sum_{i,j,k:H}\eps_i\eps_j\eps_k\bar{\Gamma}^k_{ij}\bar{\Gamma}^k_{ji} &= \sum_{i,j,k:H}\eps_i\eps_j\eps_k\Gamma^k_{ij}\Gamma^k_{ji} \;\;.
\end{split} \end{equation}
\end{subequations}

\begin{subequations} \label{threetwentytwo}
\begin{equation} \label{threetwentytwoa} \begin{split}
\scal^{V,H}_{\bar{g}} &= -(\qual^V_{\bar{g}} +\qual^H_{\bar{g}})\\
&= -f^2\qual^V_g +(q-2)f\eval{\divergence^H_g}{df}_{g,V}\\
&\mspace{20mu}-\frac{1}{2}(1-f^2)^2\tau_{g,V} -\frac{1}{2}(1-f^4)\sigma_{g,V} -\frac{2q}{f^2}\eval{df}{df}_{g,H} +\frac{2}{f}\eval{\divergence^V_g}{df}_{g,H}\\
&\mspace{20mu}-\qual^H_g +\frac{q}{f}\laplace^H_{g,H}(f) -\frac{1}{2f^2}(1-f^2)^2\tau_{g,H} +\frac{1}{2f^2}(1-f^2)(1+f^2)\sigma_{g,H} \;\;.
\end{split} \end{equation}

By \ref{curvatureON}, \eqref{threenineteen}, \eqref{threeeighteen}, \eqref{threetwentyoned}, \eqref{threefifteenb}, and \eqref{threetwentyonec}, we get:
\begin{equation} \label{threetwentytwob} \begin{split}
\scal^{H,H}_{\bar{g}} &= -2\sum_{i:H}\eps_i\partial_{\bar{e}_i}\divergence^H_{\bar{g}}(\bar{e}_i) -\sum_i\eps_i\divergence^H_{\bar{g}}(\bar{e}_i)^2 -\sum_{i,j,k:H}\eps_i\eps_j\eps_k\bar{\Gamma}^k_{ij}\bar{\Gamma}^k_{ji} +\tau_{{\bar{g}},H} -2\sum_{i,k:H}\sum_{j:V}\eps_i\eps_j\eps_k\bar{\Gamma}^k_{ij}\bar{\Gamma}^k_{ji}\\
&= -2\sum_{i:H}\eps_i\partial_{e_i}\divergence^H_g(e_i) -\sum_{i:H}\eps_i\divergence^H_g(e_i)^2 -f^2\eval{\divergence^H_g}{\divergence^H_g}_{g,V} -\sum_{i,j,k:H}\eps_i\eps_j\eps_k\Gamma^k_{ij}\Gamma^k_{ji}\\
&\mspace{20mu}+\frac{1}{2}(f^2+\frac{1}{f^2})\tau_{g,H} +\frac{1}{2}(f^2-\frac{1}{f^2})\sigma_{g,H}\\
&\mspace{20mu}-2\sum_{i,k:H}\sum_{j:V}\eps_i\eps_j\eps_k\Gamma^k_{ij}\Gamma^k_{ji} -\frac{1}{f^2}(1-f^2)\sigma_{g,H} +\frac{1}{f^2}(1-f^2)\tau_{g,H}\\
&= \scal^{H,H}_g +(1-f^2)\eval{\divergence^H_g}{\divergence^H_g}_{g,V}\\
&\mspace{20mu}-\tau_{g,H} +\frac{1}{2f^2}(f^4+1)\tau_{g,H} +\frac{1}{2f^2}(f^4-1)\sigma_{g,H} -\frac{1}{f^2}(1-f^2)\sigma_{g,H} +\frac{1}{f^2}(1-f^2)\tau_{g,H}\\
&= \scal^{H,H}_g +(1-f^2)\eval{\divergence^H_g}{\divergence^H_g}_{g,V}\\ &\mspace{20mu}+\frac{1}{2f^2}(f^4-1-2(1-f^2))\sigma_{g,H} +\frac{1}{2f^2}(f^4+1-2f^2 +2(1-f^2))\tau_{g,H}\\
&= \scal^{H,H}_g +(1-f^2)\eval{\divergence^H_g}{\divergence^H_g}_{g,V} -\frac{1}{2f^2}(1-f^2)(3+f^2)\sigma_{g,H} +\frac{1}{2f^2}(1-f^2)(3-f^2)\tau_{g,H} \;\;,
\end{split} \end{equation}

From \ref{curvatureON}, \eqref{threenineteen}, \eqref{threeeighteen}, \eqref{threetwentyonea}, \eqref{threefifteend}, \eqref{threetwentyoneb}, and \ref{laplacianformulae}, we obtain:
\begin{equation} \label{threetwentytwoc} \begin{split}
\scal^{V,V}_{\bar{g}} &= -2\sum_{i:V}\eps_i\partial_{\bar{e}_i}\divergence^V_{\bar{g}}(\bar{e}_i) -\sum_i\eps_i\divergence^V_{\bar{g}}(\bar{e}_i)^2 -\sum_{i,j,k:V}\eps_i\eps_j\eps_k\bar{\Gamma}^k_{ij}\bar{\Gamma}^k_{ji} +\tau_{{\bar{g}},V} -2\sum_{i,k:V}\sum_{j:H}\eps_i\eps_j\eps_k\bar{\Gamma}^k_{ij}\bar{\Gamma}^k_{ji}\\
&= -2f^2\sum_{i:V}\eps_i\partial_{e_i}\divergence^V_g(e_i) -2f\sum_{i:V}\eps_i\divergence^V_g(e_i)df(e_i) +2(q-1)f\sum_{i:V}\eps_i\partial_{e_i}\partial_{e_i}f\\
&\mspace{20mu} -f^2\sum_{i:V}\eps_i\divergence^V_g(e_i)^2 -(q-1)^2\eval{df}{df}_{g,V} +2(q-1)f\sum_{i:V}\eps_i\divergence^V_g(e_i)df(e_i)\\
&\mspace{20mu}-\eval{\divergence^V_g}{\divergence^V_g}_{g,H} -\frac{q^2}{f^2}\eval{df}{df}_{g,H} +\frac{2q}{f}\eval{\divergence^V_g}{df}_{g,H}\\
&\mspace{20mu} -f^2\sum_{i,j,k:V}\eps_i\eps_j\eps_k\Gamma^k_{ij}\Gamma^k_{ji} +2f\sum_{k:V}\eps_k\divergence^V_g(e_k)df(e_k) -(q-1)\eval{df}{df}_{g,V}\\
&\mspace{20mu}+\frac{1}{2}(1+f^4)\tau_{g,V} +\frac{1}{2}(1-f^4)\sigma_{g,V} +\frac{q}{f^2}\eval{df}{df}_{g,H} -\frac{2}{f}\eval{\divergence^V_g}{df}_{g,H}\\
&\mspace{20mu} -2f^2\sum_{i,k:V}\sum_{j:H}\eps_i\eps_j\eps_k\Gamma^k_{ij}\Gamma^k_{ji} +f^2(1-f^2)\sigma_{g,V} -f^2(1-f^2)\tau_{g,V}\\
&= f^2\scal^{V,V}_g +2(q-1)f\laplace^V_{g,V}(f) -q(q-1)\eval{df}{df}_{g,V} -q(q-1)\frac{1}{f^2}\eval{df}{df}_{g,H}\\
&\mspace{20mu}+(f^2-1)\eval{\divergence^V_g}{\divergence^V_g}_{g,H} +2(q-1)\frac{1}{f}\eval{\divergence^V_g}{df}_{g,H}\\
&\mspace{20mu}-f^2\tau_{g,V} +\frac{1}{2}(1+f^4)\tau_{g,V} +\frac{1}{2}(1-f^4)\sigma_{g,V} +(f^2-f^4)\sigma_{g,V} +(f^4-f^2)\tau_{g,V}\\
&= f^2\scal^{V,V}_g +2(q-1)f\laplace^V_{g,V}(f) -q(q-1)\eval{df}{df}_{g,V} -q(q-1)\frac{1}{f^2}\eval{df}{df}_{g,H}\\
&\mspace{20mu}+2(q-1)\frac{1}{f}\eval{\divergence^V_g}{df}_{g,H} -(1-f^2)\eval{\divergence^V_g}{\divergence^V_g}_{g,H}\\
&\mspace{20mu}+\frac{1}{2}(1-f^2)(1+3f^2)\sigma_{g,V} +\frac{1}{2}(1-f^2)(1-3f^2)\tau_{g,V} \;\;.
\end{split} \end{equation}
\end{subequations}

\eqref{threetwentytwoa}, \eqref{threetwentytwob}, and \eqref{threetwentytwoc} together yield
\begin{equation} \label{stretchscalar} \begin{split}
\scal_{\bar{g}} &= \scal^{V,V}_{\bar{g}} +\scal^{H,H}_{\bar{g}} +2\scal^{V,H}_{\bar{g}}\\
&= f^2\scal^{V,V}_g +2(q-1)f\laplace^V_{g,V}(f) -q(q-1)\eval{df}{df}_{g,V} -q(q-1)\frac{1}{f^2}\eval{df}{df}_{g,H}\\
&\mspace{20mu}+2(q-1)\frac{1}{f}\eval{\divergence^V_g}{df}_{g,H} -(1-f^2)\eval{\divergence^V_g}{\divergence^V_g}_{g,H}\\
&\mspace{20mu}+\frac{1}{2}(1-f^2)(1+3f^2)\sigma_{g,V} +\frac{1}{2}(1-f^2)(1-3f^2)\tau_{g,V}\\
&\mspace{20mu}+\scal^{H,H}_g +(1-f^2)\eval{\divergence^H_g}{\divergence^H_g}_{g,V} -\frac{1}{2f^2}(1-f^2)(3+f^2)\sigma_{g,H} +\frac{1}{2f^2}(1-f^2)(3-f^2)\tau_{g,H}\\
&\mspace{20mu}-2f^2\qual^V_g +2(q-2)f\eval{\divergence^H_g}{df}_{g,V}\\
&\mspace{20mu}-(1-f^2)^2\tau_{g,V} -(1-f^4)\sigma_{g,V} -\frac{4q}{f^2}\eval{df}{df}_{g,H} +\frac{4}{f}\eval{\divergence^V_g}{df}_{g,H}\\
&\mspace{20mu}-2\qual^H_g +\frac{2q}{f}\laplace^H_{g,H}(f) -\frac{1}{f^2}(1-f^2)^2\tau_{g,H} +\frac{1}{f^2}(1-f^2)(1+f^2)\sigma_{g,H}\\
&= 2(q-1)f\laplace^V_{g,V}(f) +\frac{2q}{f}\laplace^H_{g,H}(f)\\
&\mspace{20mu}-q(q-1)\eval{df}{df}_{g,V} -\frac{q(q+3)}{f^2}\eval{df}{df}_{g,H} +2(q-2)f\eval{\divergence^H_g}{df}_{g,V} +\frac{2(q+1)}{f}\eval{\divergence^V_g}{df}_{g,H}\\
&\mspace{20mu}+f^2(\scal^{V,V}_g -2\qual^V_g) +(\scal^{H,H}_g -2\qual^H_g) +(1-f^2)\Big(\eval{\divergence^H_g}{\divergence^H_g}_{g,V} -\eval{\divergence^V_g}{\divergence^V_g}_{g,H}\Big)\\
&\mspace{20mu}-\frac{(1-f^2)^2}{2f^2}\sigma_{g,H} +\frac{(1-f^2)(1+f^2)}{2f^2}\tau_{g,H} -\frac{(1-f^2)^2}{2}\sigma_{g,V} -\frac{(1-f^2)(1+f^2)}{2}\tau_{g,V} \;\;.
\end{split} \end{equation}


\subsection{Summary of the results}

\begin{theorem} \label{stretchformulae}
Let $(M,g)$ be a semi-Riemannian manifold, let $V$ be a $g$-good $q$-plane distribution on $M$, let $f\in C^\infty(M,\R_{>0})$, let $H$ denote the $g$-orthogonal distribution of $V$, and let $\bar{g}$ denote the semi-Riemannian metric $\stre(g,f,V)$. Then the following formulae hold:
\[ \begin{split}
\eval{\divergence^V_{\bar{g}}}{\divergence^V_{\bar{g}}}_{\bar{g},H} &= \eval{\divergence^V_g}{\divergence^V_g}_{g,H} +\frac{q^2}{f^2}\eval{df}{df}_{g,H} -\frac{2q}{f}\eval{\divergence^V_g}{df}_{g,H} \;\;,\\
\eval{\divergence^H_{\bar{g}}}{\divergence^H_{\bar{g}}}_{\bar{g},V} &= f^2\eval{\divergence^H_g}{\divergence^H_g}_{g,V} \;\;,\\
\sigma_{\bar{g},H} &= \frac{1}{2}(f^2+\frac{1}{f^2})\sigma_{g,H} +\frac{1}{2}(f^2-\frac{1}{f^2})\tau_{g,H} \;\;,\\
\sigma_{\bar{g},V} &= \frac{1}{2}(1+f^4)\sigma_{g,V} +\frac{1}{2}(1-f^4)\tau_{g,V} +\frac{q}{f^2}\eval{df}{df}_{g,H} -\frac{2}{f}\eval{\divergence^V_g}{df}_{g,H} \;\;,\\
\tau_{\bar{g},H} &= \frac{1}{2}(f^2+\frac{1}{f^2})\tau_{g,H} +\frac{1}{2}(f^2-\frac{1}{f^2})\sigma_{g,H} \;\;,\\
\tau_{\bar{g},V} &= \frac{1}{2}(1+f^4)\tau_{g,V} +\frac{1}{2}(1-f^4)\sigma_{g,V} +\frac{q}{f^2}\eval{df}{df}_{g,H} -\frac{2}{f}\eval{\divergence^V_g}{df}_{g,H} \;\;,\\
\laplace^V_{\bar{g},V}(u) &= f^2\laplace^V_{g,V}(u) -(q-2)f\eval{df}{du}_{g,V} \;\;,\\
\laplace^H_{\bar{g},H}(u) &= \laplace^H_{g,H}(u) \;\;,\\
\eval{\divergence^H_{\bar{g}}}{du}_{\bar{g},V} &= f^2\eval{\divergence^H_g}{du}_{g,V} \;\;,\\
\eval{\divergence^V_{\bar{g}}}{du}_{\bar{g},H} &= \eval{\divergence^V_g}{du}_{g,H} -\frac{q}{f}\eval{df}{du}_{g,H} \;\;,\\
\laplace_{\bar{g}}(u) &= f^2\laplace^V_{g,V}(u) -(q-2)f\eval{df}{du}_{g,V} +\laplace^H_{g,H}(u)\\
&\mspace{20mu}+f^2\eval{\divergence^H_g}{du}_{g,V} +\eval{\divergence^V_g}{du}_{g,H} -\frac{q}{f}\eval{df}{du}_{g,H} \;\;,
\end{split} \]

\[ \begin{split}
\qual^V_{\bar{g}} &= f^2\qual^V_g -(q-2)f\eval{\divergence^H_g}{df}_{g,V}\\ &\mspace{20mu}+\frac{(1-f^2)^2}{2}\tau_{g,V} +\frac{(1+f^2)(1-f^2)}{2}\sigma_{g,V} +\frac{q}{f^2}\eval{df}{df}_{g,H} -\frac{2}{f}\eval{\divergence^V_g}{df}_{g,H} \;\;,\\
\qual^H_{\bar{g}} &= \qual^H_g -\frac{q}{f}\laplace^H_{g,H}(f) +\frac{q}{f^2}\eval{df}{df}_{g,H} +\frac{(1-f^2)^2}{2f^2}\tau_{g,H} -\frac{(1+f^2)(1-f^2)}{2f^2}\sigma_{g,H} \;\;,\\
\scal^{V,H}_{\bar{g}} &= -f^2\qual^V_g +(q-2)f\eval{\divergence^H_g}{df}_{g,V} -\frac{1}{2}(1-f^2)^2\tau_{g,V} -\frac{(1+f^2)(1-f^2)}{2}\sigma_{g,V} -\frac{2q}{f^2}\eval{df}{df}_{g,H}\\
&\mspace{20mu}+\frac{2}{f}\eval{\divergence^V_g}{df}_{g,H} -\qual^H_g +\frac{q}{f}\laplace^H_{g,H}(f) -\frac{1}{2f^2}(1-f^2)^2\tau_{g,H} +\frac{1}{2f^2}(1-f^2)(1+f^2)\sigma_{g,H} \;\;,\\
\scal^{V,V}_{\bar{g}} &= f^2\scal^{V,V}_g +2(q-1)f\laplace^V_{g,V}(f) -q(q-1)\eval{df}{df}_{g,V} -q(q-1)\frac{1}{f^2}\eval{df}{df}_{g,H}\\
&\mspace{20mu}+2(q-1)\frac{1}{f}\eval{\divergence^V_g}{df}_{g,H} -(1-f^2)\eval{\divergence^V_g}{\divergence^V_g}_{g,H}\\ &\mspace{20mu}+\frac{1}{2}(1-f^2)(1+3f^2)\sigma_{g,V} +\frac{1}{2}(1-f^2)(1-3f^2)\tau_{g,V} \;\;,\\
\scal^{H,H}_{\bar{g}} &= \scal^{H,H}_g +(1-f^2)\eval{\divergence^H_g}{\divergence^H_g}_{g,V} -\frac{1}{2f^2}(1-f^2)(3+f^2)\sigma_{g,H} +\frac{1}{2f^2}(1-f^2)(3-f^2)\tau_{g,H} \;\;,\\
\scal_{\bar{g}} &= 2(q-1)f\laplace^V_{g,V}(f) +\frac{2q}{f}\laplace^H_{g,H}(f)\\
&\mspace{20mu}-q(q-1)\eval{df}{df}_{g,V} -\frac{q(q+3)}{f^2}\eval{df}{df}_{g,H} +2(q-2)f\eval{\divergence^H_g}{df}_{g,V} +\frac{2(q+1)}{f}\eval{\divergence^V_g}{df}_{g,H}\\
&\mspace{20mu}+f^2(\scal^{V,V}_g -2\qual^V_g) +(\scal^{H,H}_g -2\qual^H_g) +(1-f^2)\Big(\eval{\divergence^H_g}{\divergence^H_g}_{g,V} -\eval{\divergence^V_g}{\divergence^V_g}_{g,H}\Big)\\
&\mspace{20mu}-\frac{(1-f^2)^2}{2f^2}\sigma_{g,H} +\frac{(1-f^2)(1+f^2)}{2f^2}\tau_{g,H} -\frac{(1-f^2)^2}{2}\sigma_{g,V} -\frac{(1-f^2)(1+f^2)}{2}\tau_{g,V} \;\;.\\
\end{split} \]
\end{theorem}

\begin{remark}
As we mentioned in Remark \ref{wickedtrick}, the formulae in Theorem \ref{switchformulae} can formally be obtained as special cases of the formulae in Theorem \ref{stretchformulae}, by setting $f\equiv\sqrt{-1}$. It is not hard to justify this Wick rotation trick, but we do not need to do so.
\end{remark}


\section{Conformal deformation} \label{conformalsection}

We could now compute how all the functions $\sigma_{g,V}$, $\qual^H_g$, $\laplace^V_{g,V}(u)$, etc.\ (where $V$ is a $g$-good $q$-plane distribution and $H=\bot_gV$) behave under conformal changes of the metric.
\smallskip\\
Since $\conf(g,\kappa)=\stre(\stre(g,\kappa,V),\kappa,H)$, this information can easily be obtained from Theorem \ref{stretchformulae}; e.g.,
\[ \begin{split}
\sigma_{\conf(g,\kappa),V} &= \sigma_{\stre(\stre(g,\kappa,V),\kappa,H),V}\\
&= \frac{1}{2}(\kappa^2+\frac{1}{\kappa^2})\sigma_{\stre(g,\kappa,V),V} +\frac{1}{2}(\kappa^2-\frac{1}{\kappa^2})\tau_{\stre(g,\kappa,V),V}\\
&= \frac{1}{2}(\kappa^2+\frac{1}{\kappa^2})\Big(\frac{1}{2}(1+\kappa^4)\sigma_{g,V} +\frac{1}{2}(1-\kappa^4)\tau_{g,V} +\frac{q}{\kappa^2}\eval{d\kappa}{d\kappa}_{g,H} -\frac{2}{\kappa}\eval{\divergence^V_g}{d\kappa}_{g,H}\Big)\\
&\mspace{20mu} +\frac{1}{2}(\kappa^2-\frac{1}{\kappa^2})\Big(\frac{1}{2}(1+\kappa^4)\tau_{g,V} +\frac{1}{2}(1-\kappa^4)\sigma_{g,V} +\frac{q}{\kappa^2}\eval{d\kappa}{d\kappa}_{g,H} -\frac{2}{\kappa}\eval{\divergence^V_g}{d\kappa}_{g,H}\Big)\\
&= \frac{1}{4}\Big(\frac{(\kappa^4+1)^2}{\kappa^2} -\frac{(\kappa^4-1)^2}{\kappa^2}\Big)\sigma_{g,V} +\frac{1}{4}\Big(\frac{(1+\kappa^4)(1-\kappa^4)}{\kappa^2} +\frac{(\kappa^4-1)(\kappa^4+1)}{\kappa^2}\Big)\tau_{g,V}\\
&\mspace{20mu}+\frac{q}{2}\Big(\frac{\kappa^4+1}{\kappa^4} +\frac{\kappa^4-1}{\kappa^4}\Big)\eval{d\kappa}{d\kappa}_{g,H} -\Big(\frac{\kappa^4+1}{\kappa^3} +\frac{\kappa^4-1}{\kappa^3} \Big)\eval{\divergence^V_g}{d\kappa}_{g,H}\\
&= \kappa^2\sigma_{g,V} +q\eval{d\kappa}{d\kappa}_{g,H} -2\kappa\eval{\divergence^V_g}{d\kappa}_{g,H} \;\;.
\end{split} \]

We will not write down all the other formulae, since all we need in the following is information about the scalar curvature. The well-known formula for conformal changes of $\scal_g$ can be rederived in several ways from Theorem \ref{stretchformulae}:

\begin{proposition} \label{confscalar}
Let $(M,g)$ be a semi-Riemannian manifold, and let $\kappa\in C^\infty(M,\R_{>0})$. Then the scalar curvature of the metric $\conf(g,\kappa)$ is given by
\begin{equation} \label{confscalareq} \begin{split}
\scal_{\conf(g,\kappa)} = 2(n-1)\kappa\laplace_g\kappa -n(n-1)\eval{d\kappa}{d\kappa}_g +\kappa^2\scal_g \;\;.
\end{split} \end{equation}
\end{proposition}
\Proof
We could check this using the equation $\conf(g,\kappa)=\stre(\stre(g,\kappa,V),\kappa,H)$, as explained above. However, we can argue much more simply in the case of scalar curvature:
\smallskip\\
Since $\scal_{\conf(g,\kappa)} = \scal_{\stre(g,\kappa,TM)}$, we can apply the scalar curvature formula from Theorem \ref{stretchformulae} in the case $V=TM$, $q=n$, $f=\kappa$. Because all contractions over $H$ yield $0$, we have
\[
0 = \sigma_{g,H} = \tau_{g,H} = \sigma_{g,V} = \tau_{g,V} = \eval{\divergence^V_g}{\divergence^V_g}_{g,H} = \eval{\divergence^H_g}{\divergence^H_g}_{g,V} = \scal^{H,H}_g = \qual^V_g = \qual^H_g
\]
and
\[
0 = \laplace^H_{g,H}(\kappa) = \eval{d\kappa}{d\kappa}_{g,H} = \eval{\divergence^H_g}{d\kappa}_{g,V} = \eval{\divergence^H_g}{d\kappa}_{g,V} \;\;.
\]
In particular, $\scal^{V,V}_g = \scal_g$ and $\laplace^V_{g,V}(u) = \laplace_g(u)$. These facts imply the claimed equation.
\end{proof}


\section{Stretching and conformal deformation}

Let $(M,g)$ be a semi-Riemannian $n$-manifold, let $V$ be a $g$-good $q$-plane distribution on $M$, let $H$ denote the $g$-orthogonal distribution of $V$, and let $\kappa,f\in C^\infty(M,\R_{>0})$. Then

\renewcommand{\FILL}{\mspace{30mu}}
\begin{equation} \label{stretchconform} \begin{split}
&\scal_{\conf(\stre(g,f,V),\kappa)}\\[0.5ex]
&\FILL= 2(n-1)\kappa\laplace_{\stre(g,f,V)}(\kappa) -n(n-1)\eval{d\kappa}{d\kappa}_{\stre(g,f,V)} +\kappa^2\scal_{\stre(g,f,V)}\\[0.5ex]
&\FILL= 2(n-1)\kappa\Big(f^2\laplace^V_{g,V}(\kappa) -(q-2)f\eval{df}{d\kappa}_{g,V} +\laplace^H_{g,H}(\kappa)+f^2\eval{\divergence^H_g}{d\kappa}_{g,V} +\eval{\divergence^V_g}{d\kappa}_{g,H}\\
&\FILL\mspace{105mu}-\frac{q}{f}\eval{df}{d\kappa}_{g,H}\Big) \;\;-n(n-1)\Big(\eval{d\kappa}{d\kappa}_{g,H} +f^2\eval{d\kappa}{d\kappa}_{g,V}\Big)\\
&\FILL\mspace{20mu}+\kappa^2\Big(2(q-1)f\laplace^V_{g,V}(f) +\frac{2q}{f}\laplace^H_{g,H}(f)\\ &\FILL\mspace{65mu}-q(q-1)\eval{df}{df}_{g,V} -\frac{q(q+3)}{f^2}\eval{df}{df}_{g,H} +2(q-2)f\eval{\divergence^H_g}{df}_{g,V} +\frac{2(q+1)}{f}\eval{\divergence^V_g}{df}_{g,H}\\
&\FILL\mspace{65mu}+f^2(\scal^{V,V}_g -2\qual^V_g) +(\scal^{H,H}_g -2\qual^H_g) +(1-f^2)\big(\eval{\divergence^H_g}{\divergence^H_g}_{g,V} -\eval{\divergence^V_g}{\divergence^V_g}_{g,H}\big)\\
&\FILL\mspace{65mu}-\frac{(1-f^2)^2}{2f^2}\sigma_{g,H} +\frac{(1-f^2)(1+f^2)}{2f^2}\tau_{g,H} -\frac{(1-f^2)^2}{2}\sigma_{g,V} -\frac{(1-f^2)(1+f^2)}{2}\tau_{g,V}\Big)\\[0.5ex]
&\FILL= 2(n-1)\kappa\laplace^H_{g,H}(\kappa) +\frac{2q\kappa^2}{f}\laplace^H_{g,H}(f) +2(n-1)\kappa f^2\laplace^V_{g,V}(\kappa) +2(q-1)\kappa^2f\laplace^V_{g,V}(f)\\
&\FILL\mspace{20mu}-n(n-1)\eval{d\kappa}{d\kappa}_{g,H} -\frac{q(q+3)\kappa^2}{f^2}\eval{df}{df}_{g,H} -\frac{2(n-1)q\kappa}{f}\eval{df}{d\kappa}_{g,H}\\
&\FILL\mspace{20mu}-n(n-1)f^2\eval{d\kappa}{d\kappa}_{g,V} -q(q-1)\kappa^2\eval{df}{df}_{g,V} -2(n-1)(q-2)\kappa f\eval{df}{d\kappa}_{g,V}\\
&\FILL\mspace{20mu}+2(n-1)\kappa\eval{\divergence^V_g}{d\kappa}_{g,H} +\frac{2(q+1)\kappa^2}{f}\eval{\divergence^V_g}{df}_{g,H}\\
&\FILL\mspace{20mu}+2(n-1)\kappa f^2\eval{\divergence^H_g}{d\kappa}_{g,V} +2(q-2)\kappa^2f\eval{\divergence^H_g}{df}_{g,V}\\
&\FILL\mspace{20mu}+\kappa^2f^2(\scal^{V,V}_g -2\qual^V_g) +\kappa^2(\scal^{H,H}_g -2\qual^H_g) +\kappa^2(1-f^2)\big(\eval{\divergence^H_g}{\divergence^H_g}_{g,V} -\eval{\divergence^V_g}{\divergence^V_g}_{g,H}\big)\\
&\FILL\mspace{20mu}-\frac{\kappa^2(1-f^2)^2}{2f^2}\sigma_{g,H} +\frac{\kappa^2(1-f^2)(1+f^2)}{2f^2}\tau_{g,H} -\frac{\kappa^2(1-f^2)^2}{2}\sigma_{g,V} -\frac{\kappa^2(1-f^2)(1+f^2)}{2}\tau_{g,V} \;\;.
\end{split} \end{equation}

\section{All modifications in one formula}

\begin{proposition} \label{changeformula}
Let $(M,g)$ be a Riemannian $n$-manifold, let $V$ be a $q$-plane distribution on $M$, let $\kappa,f\in C^\infty(M,\R_{>0})$. Then the following formula holds (where $H\define\bot_gV$):
\renewcommand{\FILL}{\mspace{18mu}}
\[ \begin{split}
&\scal_{\conf(\stre(\switch(g,V),f,V),\kappa)}\\
&\FILL= 2(n-1)\kappa\laplace^H_{g,H}(\kappa) +\frac{2q\kappa^2}{f}\laplace^H_{g,H}(f) -2(n-1)\kappa f^2\laplace^V_{g,V}(\kappa) -2(q-1)\kappa^2f\laplace^V_{g,V}(f)\\
&\FILL\mspace{20mu}-n(n-1)\eval{d\kappa}{d\kappa}_{g,H} -\frac{q(q+3)\kappa^2}{f^2}\eval{df}{df}_{g,H} -\frac{2(n-1)q\kappa}{f}\eval{df}{d\kappa}_{g,H}\\
&\FILL\mspace{20mu}+n(n-1)f^2\eval{d\kappa}{d\kappa}_{g,V} +q(q-1)\kappa^2\eval{df}{df}_{g,V} +2(n-1)(q-2)\kappa f\eval{df}{d\kappa}_{g,V}\\
&\FILL\mspace{20mu}+2(n-1)\kappa\eval{\divergence^V_g}{d\kappa}_{g,H} +\frac{2(q+1)\kappa^2}{f}\eval{\divergence^V_g}{df}_{g,H} -2(n-1)\kappa f^2\eval{\divergence^H_g}{d\kappa}_{g,V}\\
&\FILL\mspace{20mu}-2(q-2)\kappa^2f\eval{\divergence^H_g}{df}_{g,V} +\kappa^2\bigg((1+f^2)\xi_{g,V} +\frac{1+f^2}{2f^2}\abs{\Twist_H}^2_g -\frac{f^2(1+f^2)}{2}\abs{\Twist_V}^2_g +\scal_g\bigg) \;.
\end{split} \]
\end{proposition}

\Proof
From \eqref{stretchconform} and Theorem \ref{switchformulae}, we obtain
\renewcommand{\FILL}{\mspace{30mu}}
\[ \begin{split}
&\scal_{\conf(\stre(\switch(g,V),f,V),\kappa)}\\
&\FILL= 2(n-1)\kappa\laplace^H_{\switch(g,V),H}(\kappa) +\frac{2q\kappa^2}{f}\laplace^H_{\switch(g,V),H}(f)\\
&\FILL\mspace{20mu}+2(n-1)\kappa f^2\laplace^V_{\switch(g,V),V}(\kappa) +2(q-1)\kappa^2f\laplace^V_{\switch(g,V),V}(f) -n(n-1)\eval{d\kappa}{d\kappa}_{\switch(g,V),H}\\
&\FILL\mspace{20mu}-\frac{q(q+3)\kappa^2}{f^2}\eval{df}{df}_{\switch(g,V),H} -\frac{2(n-1)q\kappa}{f}\eval{df}{d\kappa}_{\switch(g,V),H} -n(n-1)f^2\eval{d\kappa}{d\kappa}_{\switch(g,V),V}\\
&\FILL\mspace{20mu}-q(q-1)\kappa^2\eval{df}{df}_{\switch(g,V),V} -2(n-1)(q-2)\kappa f\eval{df}{d\kappa}_{\switch(g,V),V}\\
&\FILL\mspace{20mu} +2(n-1)\kappa\eval{\divergence^V_{\switch(g,V)}}{d\kappa}_{\switch(g,V),H} +\frac{2(q+1)\kappa^2}{f}\eval{\divergence^V_{\switch(g,V)}}{df}_{\switch(g,V),H}\\
&\FILL\mspace{20mu}+2(n-1)\kappa f^2\eval{\divergence^H_{\switch(g,V)}}{d\kappa}_{\switch(g,V),V} +2(q-2)\kappa^2f\eval{\divergence^H_{\switch(g,V)}}{df}_{\switch(g,V),V}\\
&\FILL\mspace{20mu}+\kappa^2f^2(\scal^{V,V}_{\switch(g,V)} -2\qual^V_{\switch(g,V)}) +\kappa^2(\scal^{H,H}_{\switch(g,V)} -2\qual^H_{\switch(g,V)})\\
&\FILL\mspace{20mu}+\kappa^2(1-f^2) \big(\eval{\divergence^H_{\switch(g,V)}}{\divergence^H_{\switch(g,V)}}_{\switch(g,V),V} -\eval{\divergence^V_{\switch(g,V)}}{\divergence^V_{\switch(g,V)}}_{\switch(g,V),H} \big)\\
&\FILL\mspace{20mu}-\frac{\kappa^2(1-f^2)^2}{2f^2}\sigma_{\switch(g,V),H} +\frac{\kappa^2(1-f^2)(1+f^2)}{2f^2}\tau_{\switch(g,V),H}\\
&\FILL\mspace{20mu}-\frac{\kappa^2(1-f^2)^2}{2}\sigma_{\switch(g,V),V} -\frac{\kappa^2(1-f^2)(1+f^2)}{2}\tau_{\switch(g,V),V}\\
&\FILL= 2(n-1)\kappa\laplace^H_{g,H}(\kappa) +\frac{2q\kappa^2}{f}\laplace^H_{g,H}(f) -2(n-1)\kappa f^2\laplace^V_{g,V}(\kappa) -2(q-1)\kappa^2f\laplace^V_{g,V}(f)\\
&\FILL\mspace{20mu}-n(n-1)\eval{d\kappa}{d\kappa}_{g,H} -\frac{q(q+3)\kappa^2}{f^2}\eval{df}{df}_{g,H} -\frac{2(n-1)q\kappa}{f}\eval{df}{d\kappa}_{g,H}\\
&\FILL\mspace{20mu}+n(n-1)f^2\eval{d\kappa}{d\kappa}_{g,V} +q(q-1)\kappa^2\eval{df}{df}_{g,V} +2(n-1)(q-2)\kappa f\eval{df}{d\kappa}_{g,V}\\
&\FILL\mspace{20mu} +2(n-1)\kappa\eval{\divergence^V_g}{d\kappa}_{g,H} +\frac{2(q+1)\kappa^2}{f}\eval{\divergence^V_g}{df}_{g,H}\\
&\FILL\mspace{20mu}-2(n-1)\kappa f^2\eval{\divergence^H_g}{d\kappa}_{g,V} -2(q-2)\kappa^2f\eval{\divergence^H_g}{df}_{g,V}\\
&\FILL\mspace{20mu}+\kappa^2f^2\Big(-\scal^{V,V}_g -2\eval{\divergence^V_g}{\divergence^V_g}_{g,H} +4\tau_{g,V} -2\sigma_{g,V} +2\qual^V_g -4\tau_{g,V}\Big)\\
&\FILL\mspace{20mu}+\kappa^2\Big(\scal^{H,H}_g +2\eval{\divergence^H_g}{\divergence^H_g}_{g,V} +2\sigma_{g,H} -4\tau_{g,H} -2\qual^H_g +4\tau_{g,H}\Big)\\
&\FILL\mspace{20mu}-\kappa^2(1-f^2)\eval{\divergence^H_g}{\divergence^H_g}_{g,V} -\kappa^2(1-f^2)\eval{\divergence^V_g}{\divergence^V_g}_{g,H}\\
&\FILL\mspace{20mu}+\frac{\kappa^2(1-f^2)^2}{2f^2}\sigma_{g,H} -\frac{\kappa^2(1-f^2)(1+f^2)}{2f^2}\tau_{g,H} -\frac{\kappa^2(1-f^2)^2}{2}\sigma_{g,V} -\frac{\kappa^2(1-f^2)(1+f^2)}{2}\tau_{g,V} \;\;.
\end{split} \]

Using the equation $\scal^{H,H}_g-2\qual^H_g = \scal_g -(\scal^{V,V}_g-2\qual^V_g)$ (cf.\ \ref{scaldecomposition} and \ref{scalqual}), the definition of $\xi_{g,V}$ (cf.\ \ref{xidef}), and the equation $\abs{\Twist_U}_g^2 = \sigma_{g,U}-\tau_{g,U}$ (cf.\ \ref{sti}), we can write the derivative-free summands in the form $\kappa^2(\ldots)$, where $(\ldots)$ is the following term:
\renewcommand{\FILL}{\mspace{20mu}}
\[ \begin{split}
&f^2\Big(-\scal^{V,V}_g -2\eval{\divergence^V_g}{\divergence^V_g}_{g,H} -2\sigma_{g,V} +2\qual^V_g\Big)\\
&\FILL\mspace{20mu}+\Big(\scal^{H,H}_g +2\eval{\divergence^H_g}{\divergence^H_g}_{g,V} +2\sigma_{g,H} -2\qual^H_g\Big)\\
&\FILL\mspace{20mu}-(1-f^2)\eval{\divergence^H_g}{\divergence^H_g}_{g,V} -(1-f^2)\eval{\divergence^V_g}{\divergence^V_g}_{g,H}\\
&\FILL\mspace{20mu}+\frac{(1-f^2)^2}{2f^2}\sigma_{g,H} -\frac{(1-f^2)(1+f^2)}{2f^2}\tau_{g,H} -\frac{(1-f^2)^2}{2}\sigma_{g,V} -\frac{(1-f^2)(1+f^2)}{2}\tau_{g,V}\\
&\FILL= \scal_g -(1+f^2)(\scal^{V,V}_g -2\qual^V_g) +(1+f^2)\Big(\eval{\divergence^H_g}{\divergence^H_g}_{g,V} -\eval{\divergence^V_g}{\divergence^V_g}_{g,H}\Big)\\
&\FILL\mspace{20mu}+\frac{(1+f^2)^2}{2f^2}\sigma_{g,H} -\frac{(1-f^2)(1+f^2)}{2f^2}\tau_{g,H} -\frac{(1+f^2)^2}{2}\sigma_{g,V} -\frac{(1-f^2)(1+f^2)}{2}\tau_{g,V}\\
&\FILL= (1+f^2)\Big(\eval{\divergence^H_g}{\divergence^H_g}_{g,V} -\eval{\divergence^V_g}{\divergence^V_g}_{g,H} -\scal^{V,V}_g +2\qual^V_g)\Big) +\scal_g\\
&\FILL\mspace{20mu}+\frac{(1+f^2)f^2}{2f^2}(\sigma_{g,H}+\tau_{g,H}) +\frac{1+f^2}{2f^2}(\sigma_{g,H}-\tau_{g,H}) -\frac{1+f^2}{2}(\sigma_{g,V}+\tau_{g,V}) -\frac{f^2(1+f^2)}{2}(\sigma_{g,V}-\tau_{g,V})\\
&\FILL= (1+f^2)\xi_{g,V} +\frac{1+f^2}{2f^2}\abs{\Twist_H}^2_g -\frac{f^2(1+f^2)}{2}\abs{\Twist_V}^2_g +\scal_g \;\;.
\end{split} \]
This completes the proof.
\end{proof}


\section{The effects on $\chi$ curvature} \label{chisection}

The contents of the present section are not needed for the proofs of the main theorems of this thesis. We will only refer to them when we discuss the esc Conjecture \ref{esc} in Section \ref{SIXesc}.
\medskip\\
Recall the definition of the function $\chi_{g,U}$ from \ref{chidef}.

\begin{lemma} \label{chiformulae}
Let $(M,g)$ be a semi-Riemannian manifold, let $V$ be a $g$-good $q$-plane distribution on $M$, let $f\in C^\infty(M,\R_{>0})$, and let $H$ denote the $g$-orthogonal distribution of $V$. Then
\[ \begin{split}
\chi_{\stre(g,f,V),V} &= \chi_{g,V} +\frac{2q}{f}\laplace^H_{g,H}(f) -\frac{q(q+3)}{f^2}\eval{df}{df}_{g,H} +\frac{2(q+1)}{f}\eval{\divergence^V_g}{df}_{g,H} \;\;,\\
\chi_{\stre(g,f,H),V} &= f^2\chi_{g,V} +2(n-q-1)f\laplace^H_{g,H}(f) -(n-q)(n-q-1)\eval{df}{df}_{g,H}\\
&\mspace{20mu}+2(n-q-2)f\eval{\divergence^V_g}{df}_{g,H} \;\;,\\
\chi_{\conf(g,f),V} &= f^2\chi_{g,V} +2(n-1)f\laplace^H_{g,H}(f) -n(n-1)\eval{df}{df}_{g,H} +2(n-1)f\eval{\divergence^V_g}{df}_{g,H} \;\;.
\end{split} \]
\end{lemma}

\Proof
Let $\bar{g}\define \stre(g,f,V)$. Theorem \ref{stretchformulae} yields
\[ \begin{split}
\chi_{\bar{g},V} &= \scal^{H,H}_{\bar{g}} -2\qual^H_{\bar{g}} +\eval{\divergence^H_{\bar{g}}}{\divergence^H_{\bar{g}}}_{\bar{g},V} -\eval{\divergence^V_{\bar{g}}}{\divergence^V_{\bar{g}}}_{\bar{g},H} +\sigma_{\bar{g},H} -\frac{\sigma_{\bar{g},V}+\tau_{\bar{g},V}}{2}\\
&= \scal^{H,H}_g +(1-f^2)\eval{\divergence^H_g}{\divergence^H_g}_{g,V} -\frac{1}{2f^2}(1-f^2)(3+f^2)\sigma_{g,H} +\frac{1}{2f^2}(1-f^2)(3-f^2)\tau_{g,H}\\
&\mspace{20mu}-2\qual^H_g +\frac{2q}{f}\laplace^H_{g,H}(f) -\frac{2q}{f^2}\eval{df}{df}_{g,H} -\frac{(1-f^2)^2}{f^2}\tau_{g,H} +\frac{(1+f^2)(1-f^2)}{f^2}\sigma_{g,H}\\
&\mspace{20mu}+f^2\eval{\divergence^H_g}{\divergence^H_g}_{g,V} -\eval{\divergence^V_g}{\divergence^V_g}_{g,H} -\frac{q^2}{f^2}\eval{df}{df}_{g,H} +\frac{2q}{f}\eval{\divergence^V_g}{df}_{g,H}\\
&\mspace{20mu}+\frac{1}{2}(f^2+\frac{1}{f^2})\sigma_{g,H} +\frac{1}{2}(f^2-\frac{1}{f^2})\tau_{g,H} -\frac{q}{f^2}\eval{df}{df}_{g,H} +\frac{2}{f}\eval{\divergence^V_g}{df}_{g,H}\\
&\mspace{20mu}-\frac{1}{2}\Big(\frac{1}{2}(1+f^4)\sigma_{g,V} +\frac{1}{2}(1-f^4)\tau_{g,V} +\frac{1}{2}(1+f^4)\tau_{g,V} +\frac{1}{2}(1-f^4)\sigma_{g,V}\Big)
\end{split} \]

\[ \begin{split}
&= \scal^{H,H}_g -2\qual^H_g +\eval{\divergence^H_g}{\divergence^H_g}_{g,V} -\eval{\divergence^V_g}{\divergence^V_g}_{g,H} -\frac{\sigma_{g,V}+\tau_{g,V}}{2}\\
&\mspace{20mu}+\frac{-3+f^4+2f^2 +2(1-f^4) +f^4+1}{2f^2}\sigma_{g,H} +\frac{3+f^4-4f^2 -2(1+f^4-2f^2) +f^4-1}{2f^2}\tau_{g,H}\\
&\mspace{20mu}+\frac{2q}{f}\laplace^H_{g,H}(f) -\frac{q(q+3)}{f^2}\eval{df}{df}_{g,H} +\frac{2(q+1)}{f}\eval{\divergence^V_g}{df}_{g,H}\\
&= \scal^{H,H}_g -2\qual^H_g +\eval{\divergence^H_g}{\divergence^H_g}_{g,V} -\eval{\divergence^V_g}{\divergence^V_g}_{g,H} -\frac{\sigma_{g,V}+\tau_{g,V}}{2} +\sigma_{g,H}\\
&\mspace{20mu}+\frac{2q}{f}\laplace^H_{g,H}(f) -\frac{q(q+3)}{f^2}\eval{df}{df}_{g,H} +\frac{2(q+1)}{f}\eval{\divergence^V_g}{df}_{g,H}\\
&= \chi_{g,V} +\frac{2q}{f}\laplace^H_{g,H}(f) -\frac{q(q+3)}{f^2}\eval{df}{df}_{g,H} +\frac{2(q+1)}{f}\eval{\divergence^V_g}{df}_{g,H} \;\;.
\end{split} \]

For $\hat{g}\define \stre(g,f,H)$, we obtain
\[ \begin{split}
\chi_{\hat{g},V} &= \scal^{H,H}_{\hat{g}} -2\qual^H_{\hat{g}} +\eval{\divergence^H_{\hat{g}}}{\divergence^H_{\hat{g}}}_{\hat{g},V} -\eval{\divergence^V_{\hat{g}}}{\divergence^V_{\hat{g}}}_{\hat{g},H} +\sigma_{\hat{g},H} -\frac{\sigma_{\hat{g},V}+\tau_{\hat{g},V}}{2}\\
&= f^2\scal^{H,H}_g +2(n-q-\!1)f\laplace^H_{g,H}(f) -(n-q)(n-q-\!1)\eval{df}{df}_{g,H} -(n-q)(n-q-\!1)\frac{1}{f^2}\eval{df}{df}_{g,V}\\
&\mspace{20mu}+2(n-q-1)\frac{1}{f}\eval{\divergence^H_g}{df}_{g,V} -(1-f^2)\eval{\divergence^H_g}{\divergence^H_g}_{g,V}\\
&\mspace{20mu}+\frac{1}{2}(1-f^2)(1+3f^2)\sigma_{g,H} +\frac{1}{2}(1-f^2)(1-3f^2)\tau_{g,H} -2f^2\qual^H_g +2(n-q-2)f\eval{\divergence^V_g}{df}_{g,H}\\ &\mspace{20mu}-(1-f^2)^2\tau_{g,H} -(1+f^2)(1-f^2)\sigma_{g,H} -\frac{2(n-q)}{f^2}\eval{df}{df}_{g,V} +\frac{4}{f}\eval{\divergence^H_g}{df}_{g,V}\\
&\mspace{20mu}+\eval{\divergence^H_g}{\divergence^H_g}_{g,V} +\frac{(n-q)^2}{f^2}\eval{df}{df}_{g,V} -\frac{2(n-q)}{f}\eval{\divergence^H_g}{df}_{g,V} -f^2\eval{\divergence^V_g}{\divergence^V_g}_{g,H}\\
&\mspace{20mu}+\frac{1}{2}(1+f^4)\sigma_{g,H} +\frac{1}{2}(1-f^4)\tau_{g,H} +\frac{n-q}{f^2}\eval{df}{df}_{g,V} -\frac{2}{f}\eval{\divergence^H_g}{df}_{g,V} -\frac{1}{2}(f^2\sigma_{g,V} +f^2\tau_{g,V})\\
&= f^2\chi_{g,V} +2(n-q-1)f\laplace^H_{g,H}(f) -(n-q)(n-q-1)\eval{df}{df}_{g,H} +2(n-q-2)f\eval{\divergence^V_g}{df}_{g,H} \;\;.
\end{split} \]

Finally, we get
\renewcommand{\FILL}{\mspace{20mu}}
\[ \begin{split}
&\chi_{\conf(g,f),V}\\
&\FILL= \chi_{\stre(\stre(g,f,V),f,H),V}\\
&\FILL= f^2\chi_{\stre(g,f,V),V} +2(n-q-1)f\laplace^H_{\stre(g,f,V),H}(f)\\
&\FILL\mspace{20mu}-(n-q)(n-q-1)\eval{df}{df}_{\stre(g,f,V),H} +2(n-q-2)f\eval{\divergence^V_{\stre(g,f,V)}}{df}_{\stre(g,f,V),H}\\
&\FILL= f^2\chi_{g,V} +2qf\laplace^H_{g,H}(f) -q(q+3)\eval{df}{df}_{g,H} +2(q+1)f\eval{\divergence^V_g}{df}_{g,H} +2(n-q-1)f\laplace^H_{g,H}(f)\\
&\FILL\mspace{20mu}-(n-q)(n-q-1)\eval{df}{df}_{g,H} +2(n-q-2)f\eval{\divergence^V_g}{df}_{g,H} -2(n-q-2)q\eval{df}{df}_{g,H}\\
&\FILL= f^2\chi_{g,V} +2(n-1)f\laplace^H_{g,H}(f) -n(n-1)\eval{df}{df}_{g,H} +2(n-1)f\eval{\divergence^V_g}{df}_{g,H} \;\;.\qedhere
\end{split} \]
\end{proof}

\begin{remark}
Note that the $\chi_{g,V}$ curvature has some remarkable features: Under stretching with a constant function --- either in $V$ or in $\bot_gV$ direction ---, it just rescales by a constant function. None of our other functions which are determined by a metric and a distribution and involve second derivatives of the metric has this property, as one can see from Theorem \ref{stretchformulae}. It ruins any simple strategy to construct, for a given distribution $V$, a metric $g$ such that $\chi_{g,V}$ becomes everywhere negative; and as we will see in Section \ref{SIXesc}, we would like to do just that.
\smallskip\\
Moreover, the behaviour of $\chi_{g,V}$ under conformal rescaling of $g$ is strikingly similar to the behaviour of $\scal_g$ under conformal rescaling of $g$; compare \ref{confscalar} and \ref{chiformulae}. The difference between the two formulae is just that $\laplace_g(f)$ occurs in the $\scal_g$ formula, while $\laplace^H_{g,H}(f) +\eval{\divergence^V_g}{df}_{g,H}$ occurs in the $\chi_{g,V}$ case, and that $\eval{df}{df}_g$ is replaced by $\eval{df}{df}_{g,H}$. Moreover, this similarity is related to Remark \ref{dfdffunctional}, which says that on a compact manifold $(M,g)$, the function $f$ is a critical point of the functional $f\mapsto \int_{(M,g)}\eval{df}{df}_{g,H}$ if and only if $\laplace^H_{g,H}(f) +\eval{\divergence^V_g}{df}_{g,H} = 0$; whereas $f$ is a critical point of the functional $f\mapsto \int_{(M,g)}\eval{df}{df}_g$ if and only if $\laplace_g(f)=0$.
\end{remark}


\chapter{The elliptic equation} \label{FOUR}

Given a manifold $M$, a function $s\in C^\infty(M,\R)$, a $q$-plane distribution $V$ on $M$, and a Riemannian ``background'' metric $g$ on $M$, we construct a semilinear second-order elliptic partial differential operator $\PD_{g,V,s} \colon C^\infty(M,\R_{>0}) \to C^\infty(M,\R)$ with the following property: If the equation $\PD_{g,V,s}(f)=0$ admits a solution $f\in C^\infty(M,\R_{>0})$, then there is a semi-Riemannian metric $h$ on $M$ with index $q$ whose scalar curvature is $s$, such that the distribution $V$ is timelike with respect to $h$.
\smallskip\\
The computation of $\PD_{g,V,s}$ will be carried out in Section \ref{FOURONE}. (The content of that section is an obvious generalisation of the approach explained in Section \ref{ONETWO} of the introduction. I assume that every reader of Section \ref{FOURONE} is already familiar, from \ref{ONETWO}, with the basic idea behind the construction, so that we can focus on the computation.) The final results, which form the basis of the rest of the thesis, are summarised in Section \ref{FOURTWO}.

\section{Computation of the PDE} \label{FOURONE}

\emph{Throughout this section, we consider a Riemannian manifold $(M,g)$ of dimension $n\in\N_{\geq2}$ (!) and a $q$-plane distribution $V$ on $M$, where $q\in\set{0,\dots,n}$. We denote the $g$-orthogonal distribution of $V$ by $H$.}
\bigskip\\
We introduce the following short notation:

\begin{definition} \label{changedefinition}
Let $f,\kappa\in C^\infty(M,\R_{>0})$. We define a semi-Riemannian metric $\change(g,f,\kappa,V)$ with index $q$ on $M$ by
\[
\change(g,f,\kappa,V) \define \conf(\stre(\switch(g,V),f,V),\kappa) \;\;.
\]
\end{definition}

In the special case where $M$ is a smooth product $S^1\times N$, where $g$ is a product metric, and where $V$ is the first-factor distribution on $M$, the Lorentzian metric $h(\kappa,f)$ which we considered in Subsection \ref{ONETWOONE} is equal to the metric $\change(g,f,\kappa,V)$. So we can use $\change(g,f,\kappa,V)$ to generalise our approach from Subsection \ref{ONETWOONE}, i.e., to construct an elliptic PDE which serves to solve the prescribed scalar problem for pseudo-Riemannian metrics with arbitrary index on arbitrary manifolds. This generalisation will be carried out in the present section.

\subsection{The scalar curvature of $\change(g,f,K(f),V)$}

For every $K\in C^\infty(\R_{>0},\R_{>0})$ and every $f\in C^\infty(M,\R_{>0})$, we can consider the function $K(f)\define K\compose f\in C^\infty(M,\R_{>0})$ and thus the semi-Riemannian metric $\change(g,f,K(f),V)$.
\smallskip\\
Now we compute the scalar curvature of $\change(g,f,K(f),V)$. The first step in the following calculation follows from Proposition \ref{changeformula}. The second step employs \ref{laplacecompose} as well as the rules $\eval{d(K\compose f)}{d(K\compose f)}_{g,U} = K'(f)^2\eval{df}{df}_{g,U}$ and $\eval{\divergence^{\bot U}_g}{d(K\compose f)}_{g,U} = K'(f)\eval{\divergence^{\bot U}_g}{df}_{g,U}$ for $U\in\set{V,H}$, which are obvious from the formulae in \ref{moreONformulae} and the fact that $d(K\compose f) = (K'\compose f)df$.

\renewcommand{\FILL}{\mspace{0mu}}
\[ \begin{split}
&\scal_{\change(g,f,K(f),V)}\\ %
&\FILL= 2(n-1)K(f)\laplace^H_{g,H}(K(f)) +\!\frac{2qK(f)^2}{f}\laplace^H_{g,H}(f) -\!2(n-1)K(f)f^2\!\laplace^V_{g,V}(K(f)) -2(q-1)K(f)^2f\!\laplace^V_{g,V}(f)\\
&\FILL\mspace{20mu}-n(n-1)\eval{d(K\mathord{\compose}f)}{d(K\mathord{\compose}f)}_{g,H} -\frac{q(q+3)K(f)^2}{f^2}\eval{df}{df}_{g,H} -\frac{2(n-1)qK(f)}{f}\eval{df}{d(K\mathord{\compose}f)}_{g,H}\\
&\FILL\mspace{20mu}+n(n-1)f^2\eval{d(K\mathord{\compose}f)}{d(K\mathord{\compose}f)}_{g,V} +q(q-1)K(f)^2\eval{df}{df}_{g,V} +\mspace{-1mu}2(n-1)(q-2)K(f)f\eval{df}{d(K\mathord{\compose}f)}_{g,V}\\
&\FILL\mspace{20mu}+2(n-1)K(f)\eval{\divergence^V_g}{d(K\mathord{\compose}f)}_{g,H} +\frac{2(q+1)K(f)^2}{f}\eval{\divergence^V_g}{df}_{g,H}\\
&\FILL\mspace{20mu}-2(n-1)K(f)f^2\eval{\divergence^H_g}{d(K\mathord{\compose}f)}_{g,V} -2(q-2)K(f)^2f\eval{\divergence^H_g}{df}_{g,V}\\
&\FILL\mspace{20mu}+K(f)^2\bigg((1+f^2)\xi_{g,V} +\frac{1+f^2}{2f^2}\abs{\Twist_H}^2_g -\frac{f^2(1+f^2)}{2}\abs{\Twist_V}^2_g +\scal_g\bigg)\\
&\FILL= 2(n-1)K(f)K'(f)\laplace^H_{g,H}(f) +2(n-1)K(f)K''(f)\eval{df}{df}_{g,H} +\frac{2qK(f)^2}{f}\laplace^H_{g,H}(f)\\
&\FILL\mspace{20mu}-2(n-1)f^2K(f)K'(f)\laplace^V_{g,V}(f) -2(n-1)f^2K(f)K''(f)\eval{df}{df}_{g,V} -2(q-1)K(f)^2f\laplace^V_{g,V}(f)\\
&\FILL\mspace{20mu}-n(n-1)K'(f)^2\eval{df}{df}_{g,H} -\frac{q(q+3)K(f)^2}{f^2}\eval{df}{df}_{g,H} -\frac{2(n-1)qK(f)K'(f)}{f}\eval{df}{df}_{g,H}\\
&\FILL\mspace{20mu}+n(n-1)f^2K'(f)^2\eval{df}{df}_{g,V} +q(q-1)K(f)^2\eval{df}{df}_{g,V} +2(n-1)(q-2)K(f)K'(f)f\eval{df}{df}_{g,V}\\
&\FILL\mspace{20mu}+2(n-1)K(f)K'(f)\eval{\divergence^V_g}{df}_{g,H} +\frac{2(q+1)K(f)^2}{f}\eval{\divergence^V_g}{df}_{g,H}\\
&\FILL\mspace{20mu}-2(n-1)f^2K(f)K'(f)\eval{\divergence^H_g}{df}_{g,V} -2(q-2)K(f)^2f\eval{\divergence^H_g}{df}_{g,V}\\
&\FILL\mspace{20mu}+K(f)^2\bigg((1+f^2)\xi_{g,V} +\frac{1+f^2}{2f^2}\abs{\Twist_H}^2_g -\frac{f^2(1+f^2)}{2}\abs{\Twist_V}^2_g +\scal_g\bigg) \;\;.
\end{split} \]

We define functions $E,F\in C^\infty(\R_{>0},\R)$ by $E\define K'/K$ and $F\define K''/K$. Then we get

\renewcommand{\FILL}{\mspace{15mu}}
\begin{equation} \begin{split} \label{eqfourtwo}
&K(f)^{-2}\scal_{\change(g,f,K(f),V)}\\ %
&\FILL= 2\Big((n-1)E(f) +\frac{q}{f}\Big)\laplace^H_{g,H}(f) -2f\Big((n-1)fE(f) +(q-1)\Big)\laplace^V_{g,V}(f)\\
&\FILL\mspace{20mu}+\Big(2(n-1)F(f) -n(n-1)E(f)^2 -\frac{q(q+3)}{f^2} -\frac{2(n-1)qE(f)}{f}\Big)\eval{df}{df}_{g,H}\\
&\FILL\mspace{20mu}+\Big(-2(n-1)f^2F(f) +n(n-1)f^2E(f)^2 +q(q-1) +2(n-1)(q-2)E(f)f\Big)\eval{df}{df}_{g,V}\\
&\FILL\mspace{20mu}+2\Big((n-1)E(f)+\frac{q+1}{f}\Big)\eval{\divergence^V_g}{df}_{g,H} -2f\Big((n-1)fE(f)+(q-2)\Big)\eval{\divergence^H_g}{df}_{g,V}\\
&\FILL\mspace{20mu}+(1+f^2)\xi_{g,V} +\scal_g +\frac{1+f^2}{2f^2}\abs{\Twist_H}^2_g -\frac{f^2(1+f^2)}{2}\abs{\Twist_V}^2_g\\
&\FILL= 2\Big((n-1)E(f) +\frac{q}{f}\Big)\laplace_g(f) -2\Big((n-1)(1+f^2)E(f) +\frac{(q-1)f^2+q}{f}\Big)\laplace^V_{g,V}(f)\\
&\FILL\mspace{20mu}+\Big(2(n-1)F(f) -n(n-1)E(f)^2 -\frac{q(q+3)}{f^2} -\frac{2(n-1)qE(f)}{f}\Big)\eval{df}{df}_g\\
&\FILL\mspace{20mu}+\Big(-2(n-1)(1+f^2)F(f) +n(n-1)(1+f^2)E(f)^2\\
&\FILL\mspace{50mu}+\frac{q}{f^2}\big((q-1)f^2+(q+3)\big) +\frac{2(n-1)E(f)}{f}\big((q-2)f^2+q\big)\Big)\eval{df}{df}_{g,V}\\
&\FILL\mspace{20mu}+\frac{2}{f}\eval{\divergence^V_g}{df}_{g,H} -2\Big((n-1)(1+f^2)E(f)+\frac{(q-2)f^2+q}{f}\Big)\eval{\divergence^H_g}{df}_{g,V}\\
&\FILL\mspace{20mu}+(1+f^2)\xi_{g,V} +\scal_g +\frac{1+f^2}{2f^2}\abs{\Twist_H}^2_g -\frac{f^2(1+f^2)}{2}\abs{\Twist_V}^2_g \;\;;
\end{split} \end{equation}

here we used the formulae $\laplace_g(f) = \laplace^H_{g,H}(f) +\laplace^V_{g,V}(f) +\eval{\divergence^H_g}{df}_{g,V} +\eval{\divergence^V_g}{df}_{g,H}$ (cf.\ \ref{laplaciansplit}, \ref{laplacedivergence}) and $\eval{df}{df}_g = \eval{df}{df}_{g,H} +\eval{df}{df}_{g,V}$.

\subsection{The correct choice of $K$}

The right hand side of \eqref{eqfourtwo} is the value of a second-order differential operator in $f$. As we explained in Section \ref{ONETWO} for a special case, the idea is to choose the function $K\in C^\infty(\R_{>0},\R_{>0})$ in such a way that this operator becomes elliptic. This happens if the function $E\define K'/K$ has the properties
\[ \begin{split}
0 &= (n-1)(1+f^2)E(f) +\frac{(q-1)f^2+q}{f}\\
\text{and} \mspace{30mu} 0 &< (n-1)E(f) +\frac{q}{f} \;\;;
\end{split} \]
because then the summand involving $\laplace^V_{g,V}(f)$ vanishes in \eqref{eqfourtwo} and the only remaining second-order term $2\big((n-1)E(f)+q/f\big)\laplace_g(f)$ is elliptic since the coefficient has no zeroes and $g$ is a Riemannian metric.
\medskip\\
We claim that the function $K\in C^\infty(\R_{>0},\R_{>0})$ given by
\begin{equation} \label{Kdef}
\boxed{ K(x) = \Big(\frac{1+x^2}{x^{2q}}\Big)^{\frac{1}{2(n-1)}} }
\end{equation}
has the desired properties. Namely, for all $x\in\R_{>0}$,
\begin{subequations}
\begin{equation} \begin{split} \label{eqa}
K'(x) &= \frac{1}{2(n-1)}\frac{x^{2q}}{(1+x^2)}K(x)\Big(-2qx^{-2q-1}+2(1-q)x^{1-2q}\Big)\\
&= -\frac{(q-1)x^2+q}{(n-1)x(1+x^2)}K(x)\\
\end{split} \end{equation}
and thus
\[
\boxed{E(x) = -\frac{(q-1)x^2+q}{(n-1)x(1+x^2)}} \;\;.
\]
This yields indeed for all $x\in\R_{>0}$
\begin{equation} \label{eqc}
(n-1)(1+x^2)E(x) +\frac{(q-1)x^2+q}{x} = 0
\end{equation}
and
\begin{equation} \label{eqd}
(n-1)E(x)+\frac{q}{x} = \frac{-(q-1)x^2-q +q(1+x^2)}{x(1+x^2)} = \frac{x}{1+x^2} > 0 \;\;.
\end{equation}
So the function $K$ given by \eqref{Kdef} is the correct choice for our purposes, as we had claimed. \emph{In the rest of this subsection, $K$ will denote this specific function.}
\medskip\\
Now we compute the coefficients of $\eval{df}{df}_g$ and $\eval{df}{df}_{g,V}$ on the right hand side of \eqref{eqfourtwo} more explicitly. (Note that we do this mostly in order to produce an explicit formula. The concrete form of these squared first-order term coefficients turns out to be irrelevant, except for the $2$-dimensional case of our problem.)
\smallskip\\
From \eqref{eqa}, we get

\[ \begin{split}
K''(x) &= -\frac{(q-1)x^2+q}{(n-1)x(1+x^2)}K'(x) -\frac{2(q-1)x^2(1+x^2) -\big((q-1)x^2+q\big)(1+3x^2)}{(n-1)x^2(1+x^2)^2}K(x)\\
&= \frac{\big((q-1)x^2+q\big)^2 -2(q-1)(n-1)x^2(1+x^2) +\big((q-1)x^2+q\big)(n-1)(1+3x^2)}{(n-1)^2x^2(1+x^2)^2}K(x)\\
&= \frac{(q-1)(n+q-2)x^4 +\big((n-1)(2q+1)+2q(q-1)\big)x^2 +q(n-1+q)}{(n-1)^2x^2(1+x^2)^2}K(x) \;\;,
\end{split} \]
where in the last step, we used the equations
\begin{gather*}
(q-1)^2 -2(q-1)(n-1) +3(q-1)(n-1) = (q-1)(q-1+n-1) = (q-1)(n+q-2) \;\;,\\[0.1ex]
2q(q-1) -2(q-1)(n-1) +(n-1)(q-1) +3q(n-1) = (n-1)(2q+1) +2q(q-1) \;\;.
\end{gather*}
We thus have
\[
\boxed{F(x) = \frac{(q-1)(n+q-2)x^4 +\big((n-1)(2q+1)+2q(q-1)\big)x^2 +q(n-1+q)}{(n-1)^2x^2(1+x^2)^2}} \;\;.
\]

We use this to compute the coefficient of $\eval{df}{df}_g$ in \eqref{eqfourtwo}:
\renewcommand{\FILL}{\mspace{20mu}}
\begin{equation} \begin{split} \label{squarecoefficient}
&2(n-1)F(x) -n(n-1)E(x)^2 -\frac{q(q+3)}{x^2} -\frac{2(n-1)qE(x)}{x}\\[0.3ex]
&\FILL= 2(n-1)\frac{(q-1)(n+q-2)x^4 +\big((n-1)(2q+1)+2q(q-1)\big)x^2 +q(n-1+q)}{(n-1)^2x^2(1+x^2)^2}\\
&\FILL\mspace{20mu}-n(n-1)\Big(\frac{(q-1)x^2+q}{(n-1)x(1+x^2)}\Big)^2 -\frac{q(q+3)}{x^2} +\frac{2(n-1)q\big((q-1)x^2+q\big)}{(n-1)x^2(1+x^2)}\\[0.3ex]
&\FILL= \frac{1}{(n-1)x^2(1+x^2)^2}\bigg(2(q-1)(n+q-2)x^4 +2\big((n-1)(2q+1)+2q(q-1)\big)x^2\\
&\FILL\mspace{200mu}+2q(n-1+q) -n\Big((q-1)^2x^4 +2q(q-1)x^2 +q^2\Big)\\
&\FILL\mspace{200mu}-q(q+3)(n-1)(1+x^2)^2 +2(n-1)q(1+x^2)\big((q-1)x^2+q\big)\bigg)\\[0.3ex]
&\FILL= \frac{1}{(n-1)x^2(1+x^2)^2}\bigg(\Big(2(q-1)(n+q-2) -n(q-1)^2 -q(q+3)(n-1) +2(n-1)q(q-1)\Big)x^4\\
&\FILL\mspace{100mu}+2\Big((n-1)(2q+1)+2q(q-1) -nq(q-1) -q(q+3)(n-1) +(n-1)q(2q-1)\Big)x^2\\
&\FILL\mspace{100mu}+q\Big(2(n-1+q) -nq -(q+3)(n-1) +2(n-1)q\Big)\bigg)\\[0.3ex]
&\FILL= \frac{\Big((q-1)^2-(n-1)(q+3)\Big)x^4 -2(q-1)(n-1-q)x^2 -q(n-1-q)}{(n-1)x^2(1+x^2)^2} \;\;;
\end{split} \end{equation}

in the last step, we used
\renewcommand{\FILL}{\mspace{20mu}}
\[ \begin{split}
&2(q-1)(n+q-2) -n(q-1)^2 -q(q+3)(n-1) +2(n-1)q(q-1)\\[0.1ex]
&\FILL= (2nq+2q^2-4q-2n-2q+4) -(nq^2-2nq+n) -(nq^2+3nq-q^2-\!3q) +(2nq^2-2nq-2q^2+2q)\\[0.1ex]
&\FILL= q^2 -q +4 -nq -3n = q^2-2q+1-nq-3n+q+3 = (q-1)^2 -(n-1)(q+3) \;\;,\\[2ex]
&(n-1)(2q+1) +2q(q-1) -nq(q-1) -q(q+3)(n-1) +(n-1)q(2q-1)\\[0.1ex]
&\FILL= (2nq+n-2q-1) +(2q^2-2q) +(-nq^2+nq) +(-nq^2+q^2-3nq+3q) +(2nq^2-nq-2q^2+q)\\[0.1ex]
&\FILL= n-1+q^2-nq = -(nq-q^2-q-n+q+1) = -(q-1)(n-q-1) \;\;,\\
\intertext{and} &2(n-1+q) -nq -(q+3)(n-1) +2(n-1)q = (2n-2+2q) -nq -(nq-q+3n-3) +(2nq-2q)\\[0.1ex]
&\FILL= -(n-1-q) \;\;.
\end{split} \]

Now we compute the coefficient of $\eval{df}{df}_{g,V}$ in \eqref{eqfourtwo}:
\renewcommand{\FILL}{\mspace{20mu}}
\begin{equation} \begin{split} \label{Vcoefficient}
&(n-1)(1+x^2)\big(-2F(x) +nE(x)^2\big) +\frac{q}{x^2}\big((q-1)x^2+(q+3)\big) +\frac{2(n-1)E(x)}{x}\big((q-2)x^2+q\big)\\[2ex]
&\FILL= \frac{-2(q-1)(n+q-2)x^4 -2\big((n-1)(2q+1)+2q(q-1)\big)x^2 -2q(n-1+q) +n((q-1)x^2+q)^2}{(n-1)x^2(1+x^2)}\\[0.4ex]
&\FILL\mspace{20mu}+\frac{q(n-1)(1+x^2)\big((q-1)x^2+(q+3)\big)}{(n-1)x^2(1+x^2)} -\frac{2(n-1)\big((q-1)x^2+q\big)\big((q-2)x^2+q\big)}{(n-1)x^2(1+x^2)}\\[2ex]
&\FILL= \frac{(q-1)\big(-2(n+q-2)+n(q-1)+q(n-1)-2(n-1)(q-2)\big)x^4}{(n-1)x^2(1+x^2)}\\[0.4ex]
&\FILL\mspace{20mu}+\frac{\big(-2(n-1)(2q+1) -4q(q-1) +2nq(q-1) +q(n-1)(2q+2) -2(n-1)q(2q-3)\big)x^2}{(n-1)x^2(1+x^2)}\\[0.4ex]
&\FILL\mspace{20mu} +\frac{q\big(-2(n-1+q)+nq+(n-1)(q+3)-2(n-1)q\big)}{(n-1)x^2(1+x^2)}\\[2ex]
&\FILL= \frac{(q-1)\big(-2n-2q+4 +nq-n +nq-q -2nq+4n+2q-4\big)x^4}{(n-1)x^2(1+x^2)}\\[0.4ex]
&\FILL\mspace{20mu}+\frac{\big(\!-4nq-\!2n+\!4q+\!2-\!4q^2+\!4q+\!2nq^2 -\!2nq+\!2nq^2+\!2nq-\!2q^2-\!2q-\!4nq^2+\!6nq+\!4q^2 -\!6q\big)x^2}{(n-1)x^2(1+x^2)}\\[0.4ex]
&\FILL\mspace{20mu} +\frac{q\big(-2n+2-2q+nq+nq+3n-q-3-2nq+2q\big)}{(n-1)x^2(1+x^2)}\\[2ex]
&\FILL= \frac{(q-1)(n-q)x^4 +2(q-1)(n-1-q)x^2 +q(n-1-q)}{(n-1)x^2(1+x^2)} \;\;; \end{split} \end{equation}
\end{subequations}

in the last step, we used
\[ \begin{split}
&-4nq-2n+4q+2-4q^2+4q+2nq^2-2nq+2nq^2+2nq-2q^2-2q-4nq^2+6nq+4q^2-6q\\[0.2ex]
&\FILL= 2nq-2q^2-2n+2 = 2(nq-q-q^2-n+1+q) = 2(q-1)(n-1-q) \;\;.
\end{split} \]

\subsection{The equation}

For the function $K$ given by \eqref{Kdef}, we can rewrite \eqref{eqfourtwo}, taking  the equations \eqref{eqc}, \eqref{eqd}, \eqref{squarecoefficient}, \eqref{Vcoefficient} into account (note that we can use \eqref{eqc} also for the coefficient of $\eval{\divergence^H_g}{df}_{g,V}$). This yields

\renewcommand{\FILL}{\mspace{0mu}}
\begin{equation} \begin{split}
&K(f)^{-2}\scal_{\change(g,f,K(f),V)}\\[1.5ex]
&\FILL= \frac{2f}{1+f^2}\laplace_g(f) +\frac{\Big((q-1)^2-(n-1)(q+3)\Big)f^4 -2(q-1)(n-1-q)f^2 -q(n-1-q)}{(n-1)f^2(1+f^2)^2}\eval{df}{df}_g\\[0.9ex]
&\FILL\mspace{20mu}+\frac{(q-1)(n-q)f^4 +\!2(q-1)(n-1-q)f^2 +\!q(n-1-q)}{(n-1)f^2(1+f^2)}\eval{df}{df}_{g,V} +\!\frac{2}{f}\eval{\divergence^V_g}{df}_{g,H} +\!2f\eval{\divergence^H_g}{df}_{g,V}\\[0.9ex]
&\FILL\mspace{20mu}+(1+f^2)\xi_{g,V} +\scal_g +\frac{1+f^2}{2f^2}\abs{\Twist_H}^2_g -\frac{f^2(1+f^2)}{2}\abs{\Twist_V}^2_g \;\;,
\end{split} \end{equation}

hence
\begin{equation} \begin{split} \label{preeq}
0 &= 2\laplace_g(f) +\frac{\Big((q-1)^2-(n-1)(q+3)\Big)f^4 -2(q-1)(n-1-q)f^2 -q(n-1-q)}{(n-1)f^3(1+f^2)}\eval{df}{df}_g\\
&\mspace{20mu}+\frac{(q-1)(n-q)f^4 +2(q-1)(n-1-q)f^2 +q(n-1-q)}{(n-1)f^3}\eval{df}{df}_{g,V}\\
&\mspace{20mu}+\frac{2(1+f^2)}{f^2}\eval{\divergence^V_g}{df}_{g,H} +2(1+f^2)\eval{\divergence^H_g}{df}_{g,V} +\frac{(1+f^2)^2}{f}\xi_{g,V} +\frac{1+f^2}{f}\scal_g\\
&\mspace{20mu}+\frac{(1+f^2)^2}{2f^3}\abs{\Twist_H}^2_g -\frac{f(1+f^2)^2}{2}\abs{\Twist_V}^2_g -K(f)^{-2}\frac{1+f^2}{f}\scal_{\change(g,f,K(f),V)} \;\;,
\end{split} \end{equation}

where $K(f)^{-2}(1+f^2)/f = f^{2q/(n-1)-1}(1+f^2)^{1-1/(n-1)}$.
\smallskip\\
This equation for $\scal_{\change(g,f,K(f),V)}$ is the main result of the computations we did in Chapters \ref{THREE} and \ref{FOUR}. Let us summarise what we have proved.


\section{Summary of the results obtained so far} \label{FOURTWO}

\begin{definition} \label{PDdef}
Let $(M,g)$ be a Riemannian manifold of dimension $n\geq2$, let $q\in\set{0,\dots,n}$, let $V$ be a $q$-plane distribution on $M$, and let $s\in C^\infty(M,\R)$. We denote the $g$-orthogonal distribution of $V$ by $H$.
\smallskip\\
We define the functions $a_{n,q},b_{n,q}\in C^\infty(\R_{>0},\R)$ by
\[ \begin{split}
a_{n,q}(x) &= \frac{\Big((q-1)^2-(n-1)(q+3)\Big)x^4 -2(q-1)(n-1-q)x^2 -q(n-1-q)}{(n-1)x^3(1+x^2)} \;\;,\\
b_{n,q}(x) &= \frac{(q-1)(n-q)x^4 +2(q-1)(n-1-q)x^2 +q(n-1-q)}{(n-1)x^3} \;\;.
\end{split} \]
We define a semilinear elliptic second order differential operator $\PD_{g,V,s} \colon C^\infty(M,\R_{>0})\to C^\infty(M,\R)$ by
\begin{center} \fbox{\parbox{160mm}{
\[ \begin{split}
\PD_{g,V,s}(f) &\define 2\laplace_g(f) +a_{n,q}(f)\abs{df}^2_g +b_{n,q}(f)\abs{df}^2_{g,V} +\frac{2(1+f^2)}{f^2}\eval{\divergence^V_g}{df}_{g,H} +2(1+f^2)\eval{\divergence^H_g}{df}_{g,V}\\
&\mspace{20mu}+\frac{(1+f^2)^2}{2f^3}\abs{\Twist_H}^2_g -\frac{f(1+f^2)^2}{2}\abs{\Twist_V}^2_g +\frac{(1+f^2)^2}{f}\xi_{g,V} +\frac{1+f^2}{f}\scal_g\\
&\mspace{20mu}-f^{\frac{2q}{n-1}-1}(1+f^2)^{1-\frac{1}{n-1}}s \;\;.
\end{split} \]}}
\end{center}
\end{definition}

\begin{theorem} \label{THEEQUATION}
Let $(M,g)$ be a Riemannian manifold of dimension $n\geq2$, let $q\in\set{0,\dots,n}$, let $V$ be a $q$-plane distribution on $M$, let $s\in C^\infty(M,\R)$, let $H$ denote the $g$-orthogonal distribution of $V$. If the elliptic partial differential equation
\[
\PD_{g,V,s}(f) = 0
\]
has a solution $f\in C^\infty(M,\R_{>0})$, then there is a semi-Riemannian metric $h$ of index $q$ on $M$ such that $\scal_h = s$, such that $V$ is timelike with respect to $h$, and such that $H$ is $h$-orthogonal to $V$, and so, in particular, $H$ is spacelike with respect to $h$. Namely, we can choose $h = \change(g,f,K\compose f,V)$, where the function $K\in C^\infty(\R_{>0},\R_{>0})$ is given by
\[
K(x) = \Big(\frac{1+x^2}{x^{2q}}\Big)^{\frac{1}{2(n-1)}}.
\]
\end{theorem}
\Proof
The semi-Riemannian metric \mbox{$\change(g,f,K\!\compose\!f,V)$} has index $q$, makes $V$ timelike, and makes $H$ orthogonal to $V$; these properties follow immediately from the definition of the $\switch$, $\stre$, and $\conf$ operations. By Equation \eqref{preeq}, the function $S\define\scal_{\change(g,f,K(f),V)}$ satisfies $\PD_{g,V,S}(f)=0$. If $f\in C^\infty(M,\R_{>0})$ satisfies $\PD_{g,V,s}(f)=0$, then $\PD_{g,V,s}(f) = \PD_{g,V,S}(f)$, which obviously implies $s=S$.
\end{proof}

\subsection{Some special cases}

We start with a well-known special case of Definition \ref{PDdef}:

\begin{remark}[the Riemannian case $q=0$] \label{equationRiemann}
Let $(M,g)$ be a Riemannian manifold of dimension $n\geq2$, let $V$ be the unique $0$-plane distribution on $M$, let $s\in C^\infty(M,\R)$. Then the $g$-orthogonal distribution $H$ of $V$ is the whole tangent bundle $TM$, and the operator $\PD_{g,V,s} \colon C^\infty(M,\R_{>0})\to C^\infty(M,\R)$ is given by
\begin{equation} \label{riemcase}
f\mapsto 2\laplace_g(f) +\frac{(1-3(n-1))f^2 +2(n-1)}{(n-1)f(1+f^2)}\abs{df}^2_g +\frac{1+f^2}{f}\scal_g -\frac{(1+f^2)^{\frac{n-2}{n-1}}}{f}s \;\;.
\end{equation}
(This follows immediately by substituting $q=0$ in Definition \ref{PDdef}. The terms $\abs{df}^2_{g,V}$, $\eval{\divergence^V_g}{df}_{g,H}$, $\eval{\divergence^H_g}{df}_{g,V}$, $\xi_{g,V}$, $\abs{\Twist_H}^2_g$, $\abs{\Twist_V}^2_g$ vanish since their definitions involve contractions over the bundle $V$, and each such contraction vanishes since $V$ has rank $0$.)
\smallskip\\
Moreover,
\[ \begin{split}
\change(g,f,K(f),V) &= \conf(\stre(\switch(g,V),f,V),K(f))\\
&= \conf(\stre(g,f,V),K(f)) = \conf(g,K(f)) \;\;,
\end{split} \]
so the elliptic equation $\PD_{g,V,s}=0$ is just the well-known equation \eqref{yamabe1} of conformal rescaling, up to substitution of $\kappa$ by $K(f) = (1+f^2)^{\frac{1}{2(n-1)}}$. (This can of course also be checked by a direct calculation.)
\smallskip\\
Equation \eqref{yamabe1} is a standard tool for the solution of the Riemannian prescribed scalar curvature problem; cf.\ the review in Appendix \ref{D1}. It has a solution $\kappa\in C^\infty(M,\R_{>1})$ if and only if \eqref{riemcase} has a solution $f\in C^\infty(M,\R_{>0})$. As we already remarked during the discussion of the special case in Subsection \ref{ONETWOTWO}, our approach to the pseudo-Riemannian prescribed scalar curvature problem is thus --- in a slightly restricted sense, because of the $\kappa>1$ requirement --- a generalisation of the standard approach to the Riemannian problem.
\end{remark}

\begin{remark}[the Lorentzian case $q=1$] \label{lorentzequation}
Let $(M,g)$ be a Riemannian manifold of dimension $n\geq2$, let $V$ be a line distribution on $M$, let $s\in C^\infty(M,\R)$, let $H$ denote the $g$-orthogonal distribution of $V$.
\smallskip\\
With the abbreviation $\alpha(n)\define\frac{n-2}{n-1} \in\cointerval{0}{1}$, the operator $\PD_{g,V,s} \colon C^\infty(M,\R_{>0})\to C^\infty(M,\R)$ is given by
\[ \begin{split}
f &\mapsto 2\laplace_g(f) -\frac{4f^4+\alpha(n)}{f^3(1+f^2)}\abs{df}^2_g +\frac{\alpha(n)}{f^3}\abs{df}^2_{g,V} +\frac{2(1+f^2)}{f^2}\eval{\divergence^V_g}{df}_{g,H} +2(1+f^2)\eval{\divergence^H_g}{df}_{g,V}\\
&\mspace{20mu}+\frac{(1+f^2)^2}{2f^3}\abs{\Twist_H}^2_g +\frac{(1+f^2)^2}{f}\xi_{g,V} +\frac{1+f^2}{f}\scal_g -f\Big(\frac{1+f^2}{f^2}\Big)^{\alpha(n)}s \;\;.
\end{split} \]
(This follows immediately by substituting $q=1$ in Definition \ref{PDdef}. Since $V$ is a line distribution, we have $\abs{\Twist_V}^2_g = 0$.)
\smallskip\\
Recall that we have seen in Fact \ref{lineformulae} how the function $\xi_{g,V}$ looks like when $V$ is a line bundle.
\end{remark}


\chapter{Everywhere twisted distributions} \label{FIVE}

As we will see in Chapter \ref{SIX}, it would be very helpful to know whether a given $n$-manifold admits an everywhere twisted\footnote{Recall that \emph{everywhere twistedness} has been introduced in Definition \ref{twistednessdef}.} $q$-plane distribution, in order to solve the plain version of the prescribed scalar curvature problem. It would be even more helpful --- in order to solve the homotopy class version --- to know whether \emph{each homotopy class} of $q$-plane distributions contains one which is everywhere twisted. Our aim in the present chapter is to prove existence results in this direction.
\smallskip\\
A special case of this problem is existence of contact structures on $3$-manifolds: A $2$-plane distribution on a $3$-manifold is a contact structure if and only if it is everywhere twisted. More generally, every contact structure on a $(2n+1)$-manifold $M$ is an everywhere twisted $2n$-plane distribution on $M$. The converse is not true in general if $n\geq2$ --- the everywhere twisted condition is much weaker than the contact condition then. On manifolds of dimension $2n+1\geq5$, existence of everywhere twisted $2n$-plane distributions does therefore hold in many more cases, and is much easier to prove, than existence of contact structures. In dimension $3$, we can apply the standard existence theorem for contact structures.
\smallskip\\
Another --- less well-known --- special case of our problem is existence of even-contact structures on $4$-manifolds. Even-contact structures are analogues of contact structures on even-dimensional manifolds. Appendix \ref{CONTACT} contains a review of basic facts about contact and even-contact structures. Though we will not use them (except for the $3$-dimensional contact theorem), we refer to them occasionally to put the everywhere twisted results into context.
\smallskip\\
Except in the $3$-dimensional case, M.\ Gromov's convex integration technique applies to the twistedness partial differential relation; i.e., the relation is \emph{ample} in the sense of Gromov. This means that the existence problem for everywhere twisted distributions can be reduced from a \emph{differential} topological problem to an obstruction-theoretic problem in \emph{algebraic} topology; namely to the question whether a certain vector bundle admits a nowhere vanishing section. (D.\ McDuff has verified ampleness in the even-dimensional case even for the stronger even-contact partial differential relation; so it is hardly surprising that the twistedness relation turns out to be ample. But as far as I know, this result is not contained in the literature.)
\smallskip\\
The remaining obstruction-theoretic problem for $q$-plane distributions on an $n$-manifold is completely trivial if $(n-q)(q-2)\geq 2$. This will enable us to solve (in Chapter \ref{SIX}) the homotopy class version of the prescribed scalar curvature problem completely in the case $3\leq q\leq n-3$. But still something can be said if $q=2$ or $(n,q)=(4,3)$.
\smallskip\\
In Section \ref{twistCinfty} at the end of the chapter, we will discuss the question whether integrable distributions can be approximated by everywhere twisted ones \emph{in the $C^\infty$-topology}; the most straightforward application of the convex integration technique yields only $C^0$-approximations. These approximation theorems will not be used for the main theorems of the thesis but just for the discussion of the esc Conjecture in Section \ref{SIXesc}. We will therefore not provide detailed proofs of all statements in Section \ref{twistCinfty}.


\section{Twistedness as a partial differential relation} \label{FIVEONE}

We want to prove existence of everywhere twisted distributions via Gromov's h-principle theorems (cf.\ Appendix \ref{hprinciple}). In order to do that, we have to interpret twistedness as a partial differential relation on the $1$-jet bundle of the Grassmann bundle $\Gr_q(TM)\to M$.\footnote{cf.\ Appendix \ref{Grassmannappendix} for basic facts about Grassmann bundles} This is straightforward, but quite technical because one has to deal with derivatives of sections in the bundle of Grassmannians on the one hand, and has to interpret each section as a bundle and consider derivatives of sections in that bundle on the other hand.
\medskip\\
More precisely, we want to show that Definition \ref{twistjetdefinition} below is well-defined. The point of the definition is this: The twistedness of a distribution $V$ is defined via the Lie bracket of sections in $V$, i.e.\ via \emph{the $1$-jets of sections in $V$}. Interpreting $V$ as a section in a Grassmann bundle, we can consider \emph{the $1$-jet of $V$ itself}. One would guess that this $1$-jet contains enough information about the $1$-jets of sections in $V$ to determine the twistedness of $V$.
\smallskip\\
In order to prove that this is indeed true, we start with some preparations.

\begin{definition}
Let $E\to M$ be a vector bundle, let $q\in\N$. Then we can consider the trivial vector bundle $\R^q\times M\to M$ and the thereby defined vector bundle $\Lin(\R^q,E)$ over $M$ (whose fibre over $x$ consists of the linear maps from $\R^q$ to $E_x$). We denote by $\Mon(\R^q,E)$ the open sub fibre bundle of $\Lin(\R^q,E)$ whose fibre over $x\in M$ is the set $\Mon(\R^q,E_x)$ of all monomorphisms (i.e.\ injective linear maps) from $\R^q$ to $E_x$.
\end{definition}

\begin{lemma} \label{grassmannlift}
Let $M$ be a manifold, let $V$ be a $q$-plane distribution on $M$, let $x\in M$. Then there exist an open neighbourhood $U$ of $x$ in $M$ and a section $\hat{V}\in C^\infty(U\ot\Mon(\R^q,TM))$ such that $\im(\hat{V}(y))=V(y)$ for all $y\in U$.
\end{lemma}
\Proof
Let $n$ denote the dimension of $M$. The sub vector bundle $V$ of $TM$ admits a sub vector bundle chart $\phi\colon TU\to U\times\R^n$ around $x$; i.e., $U$ is an open neighbourhood of $x$, and $\phi$ is a vector bundle chart of $TU$ which maps each fibre $V_y$ to $\set{y}\times\R^q \subseteq \set{y}\times\R^n$. The map $\hat{V}\in C^\infty(U\ot\Mon(\R^q,TM))$ given by $\hat{V}(y)(v) = \phi^{-1}(y,v)$ has the desired properties.
\end{proof}

\begin{lemma} \label{derivative}
Let $U\subseteq\R^n$ be an open neighbourhood of $0$, let $q\in\set{0,\dots,n}$, let $\hat{V}\in C^\infty(U,\Mon(\R^q,\R^n))$, let $V\in C^\infty(U,\Gr_q(\R^n))$ be the map which is pointwise the image of $\hat{V}$. Let $\iota$ denote the vector space isomorphism $\hat{V}(0)\in\Lin(\R^q,V(0))$, and let $\pr\colon \R^n \to \R^n/V(0)$ denote the obvious projection. Then the derivative $D_0V\in \Lin(\R^n,T_{V(0)}\Gr_q(\R^n))$ is via the identification\footnote{cf.\ Appendix \ref{Grassmannappendix}} $T_{V(0)}\Gr_q(\R^n) = \Lin(V(0),\R^n/V(0))$ given by
\[
(D_0V)(q) = \pr\compose(D_0\hat{V})(q)\compose\iota^{-1} \in \Lin(V(0),\R^n/V(0)) \;\;.
\]
\end{lemma}
\Proof
Let $H$ be a complementary sub vector space of $V(0)$ in $\R^n$, and let $\pr_H\colon \R^n = V(0)\oplus H\to H$ denote the obvious projection. We identify $H$ with $\R^n/V(0)$, and thus $\pr_H$ with $\pr$.
\smallskip\\
Since $V(0)$ is complementary to $H$, there exists an open neighbourhood $U'\subseteq U$ of $x$ such that $V(x)$ is complementary to $H$ for all $x\in U'$. Recall that for every $\lambda\in\Lin(V(0),H)$, the vector space $\lambda+V(0)\in\Compl(H)$ is defined to be $\set{v+\lambda(v)\suchthat v\in V(0)}$; and that the map $\varphi\colon \Lin(V(0),H)\to \Gr_q(\R^n)$ given by $\varphi(\lambda) = \lambda+V(0)$ is a diffeomorphism onto an open neighbourhood of $V(0)$ in $\Gr_q(\R^n)$. The canonical vector space isomorphism $\Lin(V(0),H)\to T_{V(0)}\Gr_q(\R^n)$ which appears implicitly in the statement of the lemma is $D_{V(0)}\varphi$. So what we have to prove is, for all $q\in\R^n$, the equality
\[
(D_0V)(q) = (D_{V(0)}\varphi)\Big(\pr_H\compose(D_0\hat{V})(q)\compose\iota^{-1}\Big) \in T_{V(0)}\Gr_q(\R^n) \;\;.
\]
Let $\pr_V \colon \R^n = V(0)\oplus H\to V(0)$ denote the obvious projection. For each $x\in U'$, the map $\iota_x\define \pr_V\compose(\hat{V}(x))\in\Lin(\R^q,V(0))$ is bijective because $\im(\hat{V}(x))+H$ is complementary to $H$. We define $\lambda_x\in\Lin(V(0),H)$ by $\lambda_x\define \pr_H\compose(\hat{V}(x))\compose\iota_x^{-1}$. This $\lambda_x$ is the unique $\lambda\in\Lin(V(0),H)$ such that $V(x)=\lambda+V(0)$ since
\[ \begin{split}
\lambda_x+V(0) &= \set{v+\lambda_x(v)\suchthat v\in V(0)} = \set{\iota_x(w) +\lambda_x(\iota_x(w)) \suchthat w\in\R^q}\\
&= \set{\pr_V\big(\hat{V}(x)(w)\big) +\pr_H\big(\hat{V}(x)(w)\big) \suchthat w\in\R^q} = \im(\hat{V}(x)) = V(x) \;\;.
\end{split} \]
Hence $(\varphi^{-1}\compose V)(x)=\lambda_x$ for all $x\in U'$, and thus we obtain (by the chain rule and the product rule):
\[ \begin{split}
(D_{V(0)}\varphi)^{-1}\compose(D_0V)(q) &= D_0(x\mapsto\lambda_x)(q)\\
&= \pr_H\compose\Big(D_0\big(x\mapsto\hat{V}(x)\compose\iota_x^{-1}\big)(q)\Big)\\
&= \pr_H\compose\Big((D_0\hat{V})(q)\compose\iota_0^{-1} +\hat{V}(0)\compose\Big(D_0\big(x\mapsto \iota_x\big)^{-1}(q)\Big)\Big)\\
&= \pr_H\compose(D_0\hat{V})(q)\compose\iota^{-1} \;\;,
\end{split} \]
since $\pr_H\compose(\hat{V}(0))=0$ and $\iota_0=\iota$. This implies the statement of the lemma.
\end{proof}

\begin{lemma} \label{bracketlemma}
Let $U\subseteq\R^n$ be an open neighbourhood of $0$, let $q\in\set{0,\dots,n}$, let $\hat{V}\in C^\infty(U,\Mon(\R^q,\R^n))$, let $V\in C^\infty(U,\Gr_q(\R^n))$ be the map which is pointwise the image of $\hat{V}$, and let $\iota$ denote the vector space isomorphism $\hat{V}(0)\in\Lin(\R^q,V(0))$. For every $v\in V(0)$, the map $\hat{v}\in C^\infty(U,\R^n)$ given by $\hat{v}(y)=\hat{V}(y)(\iota^{-1}(v))$ is a section in the vector bundle $V$ with $\hat{v}(0)=v$. If $w\in V(0)$, then the value of the Lie bracket of $\hat{v}$ and $\hat{w}$ in the point $0$ (where $\hat{w}\in C^\infty(U,\R^n)$ is given by $\hat{w}(y)=\hat{V}(y)(\iota^{-1}(w))$) is
\[
[\hat{v},\hat{w}](0) = (D_0\hat{V})(v)(\iota^{-1}(w)) -(D_0\hat{V})(w)(\iota^{-1}(v)) \in\R^n \;\;.
\]
\end{lemma}
\Proof
The map $\hat{v}$ is a section in $V$ since $\hat{v}(y) \in \im(\hat{V}(y))=V(y)$ for all $y\in U$. Moreover, $\hat{v}(0) = \hat{V}(0)(\iota^{-1}(v)) = \iota(\iota^{-1}(v)) = v$. The Lie bracket of vector fields $\hat{v},\hat{w}$ on an open subset of $\R^n$ is given by $[\hat{v},\hat{w}](x) = (D_x\hat{w})(\hat{v}(x)) -(D_x\hat{v})(\hat{w}(x))$. Since $(D_0\hat{v})(q) = (D_0\hat{V})(q)(\iota^{-1}(v))$ for all $q\in\R^n$ (and similarly for $\hat{w}$), we obtain
\[
[\hat{v},\hat{w}](0) = (D_0\hat{w})(v) -(D_0\hat{v})(w) = (D_0\hat{V})(v)(\iota^{-1}(w)) -(D_0\hat{V})(w)(\iota^{-1}(v)) \;\;.\qedhere
\]
\end{proof}

\begin{lemma} \label{twistedjet}
Let $U\subseteq\R^n$ be an open neighbourhood of $0$, let $q\in\set{0,\dots,n}$, let $V\in C^\infty(U,\Gr_q(\R^n))$. Then the twistedness $\Twist_V(0)\in\Lambda^2(V(0)^\ast)\otimes(\R^n/V(0))$ of the distribution $V$ in the point $0$ is via the derivative $D_0V\in \Lin(\R^n,\Lin(V(0),\R^n/V(0)))$ given by
\[
\Twist_V(0)(v_0,v_1) = (D_0V)(v_0)(v_1) -(D_0V)(v_1)(v_0) \;\;.
\]
\end{lemma}
\Proof
By Lemma \ref{grassmannlift}, there exist an open neighbourhood $U'\subseteq U$ of $0$ in $\R^n$ and a section $\hat{V}\in C^\infty(U',\Mon(\R^q,\R^n))$ such that $\im(\hat{V}(y))=V(y)$ for all $y\in U'$. Let $\iota$ denote the vector space isomorphism $\hat{V}(0)\in\Lin(\R^q,V(0))$, and let $\pr\colon \R^n\to\R^n/V(0)$ denote the obvious projection.
\smallskip\\
Lemma \ref{bracketlemma} tells us that for each $i\in\set{0,1}$, the map $\hat{v}_i\in C^\infty(U',\R^n)$ given by $\hat{v}_i(y)=\hat{V}(y)(\iota^{-1}(v_i))$ is a section in $V$ with $\hat{v}_i(0)=v_i$. Lemma \ref{bracketlemma} and \ref{derivative} yield:
\[ \begin{split}
\Twist_V(0)(v_0,v_1) &= \pr\Big([\hat{v}_0,\hat{v}_1](0)\Big)\\
&= \pr\Big((D_0\hat{V})(v_0)(\iota^{-1}(v_1)) -(D_0\hat{V})(v_1)(\iota^{-1}(v_0))\Big)\\
&= \Big((D_0V)(v_0)\compose\iota\Big)(\iota^{-1}(v_1)) -\Big((D_0V)(v_1)\compose\iota\Big)(\iota^{-1}(v_0))\\
&= (D_0V)(v_0)(v_1) -(D_0V)(v_1)(v_0) \;\;.\qedhere
\end{split} \]
\end{proof}

\begin{corollary} \label{twistlemma}
Let $M$ be an $n$-manifold, let $x\in M$, $q\in\set{0,\dots,n}$, and let $V,W$ be $q$-plane distributions on $M$; we interpret them as sections in the Grassmann bundle $\Gr_q(TM)$ over $M$. If $j^1_xV = j^1_xW$, i.e.\ the $1$-jets of $V,W$ in the point $x$ are equal, then $\Twist_V(x) = \Twist_W(x)$, i.e., the twistedness of the distribution $V$ in $x$ is equal to the twistedness of $W$ in $x$.
\end{corollary}
\Proof
Since twistedness and jets are defined by local data, it suffices to consider the germs of $V,W\in C^\infty(M\ot\Gr_q(TM))$ at the point $x$. By choosing a suitable manifold chart of $M$ around $x$, we may therefore assume without loss of generality that $M$ is an open subset of $\R^n$, and that $x=0\in M$. Then there is a trivialisation $TM = M\times\R^n$ of the tangent bundle, an induced trivialisation $\Gr_q(TM) = M\times\Gr_q(\R^n)$ of the Grassmann bundle, and $V,W$ can be interpreted as elements of $C^\infty(M,\Gr_q(\R^n))$.
\smallskip\\
The statement $j^1_0V = j^1_0W$ means precisely that $V(0)=W(0)$ and that the derivatives $D_0V,D_0W \colon \R^n \to T_{V(0)}(\Gr_q(\R^n))=\Lin(V(0),\R^n/V(0))$ are equal (here we have identified $T_0\R^n$ with $\R^n$). By Lemma \ref{twistedjet}, this implies for all $v,w\in V(0)=W(0)$:
\[
\Twist_V(0)(v,w) = (D_0V)(v)(w) -(D_0V)(w)(v) = (D_0W)(v)(w) -(D_0W)(w)(v) = \Twist_W(0)(v,w) \;\;.
\]
Hence $\Twist_V(0)=\Twist_W(0)$.
\end{proof}

\begin{definition}[twistedness on the jet level] \label{twistjetdefinition}
Let $M$ be an $n$-manifold, $q\in\set{0,\dots,n}$, and let $\mathscr{V}\in J^1\Gr_q(TM)$. We define $\Twist_{\mathscr{V}}$, the \emph{twistedness of $\mathscr{V}$}, as follows. Let $p^1\colon J^1\Gr_q(TM) \to M$ and $p^1_0 \colon J^1\Gr_q(TM) \to \Gr_q(TM)$ denote the standard bundle projections, and let $x\define p^1(\mathscr{V}) \in M$ and $V\define p^1_0(\mathscr{V}) \in \Gr_q(T_xM)$. Then $\Twist_{\mathscr{V}}$ is the element in the vector space $\Lambda^2(V^\ast)\otimes\bot V$ given by $\Twist_{\mathscr{V}} \define \Twist_{\tilde{V}}(x)$, where $\tilde{V}\in C^\infty(U\ot\Gr_q(TU))$ is any distribution on a neighbourhood $U$ of $x$ in $M$ such that $j^1_x\tilde{V} = \mathscr{V}$ (and thus in particular $\tilde{V}(x)=V$).\footnote{Note that although the $1$-jet bundle $J^1E$ of a fibre bundle $E\to M$ is defined via $1$-jets of local $C^1$ sections in $E$, there is for every $\mathscr{V}\in J^1_xE$ also a local $C^\infty$ section $V$ in $E$ with $\mathscr{V}=j^1_xV$. So our definition makes sense.} Corollary \ref{twistlemma} implies that $\Twist_{\mathscr{V}}$ is well-defined, i.e.\ independent of the choice of $\tilde{V}$.
\end{definition}

\begin{definition}[the twistedness relation $\mathscr{R}_{M,q}$]
Let $M$ be an $n$-manifold, and let $q\in\set{0,\dots,n}$. We define $\mathscr{R}_{M,q}$ to be the subset of $J^1\Gr_q(TM)$ consisting of all $\mathscr{V}$ with $\Twist_{\mathscr{V}}\neq 0$. (Here $0$ denotes the zero element of the vector space $\Lambda^2(V^\ast)\otimes\bot V$, where $V = p^1_0(\mathscr{V})$.) We call $\mathscr{R}_{M,q}$ the \emph{twistedness relation}.
\end{definition}

\begin{lemma} \label{openrelation}
Let $M$ be an $n$-manifold, let $q\in\set{0,\dots,n}$. Then $\mathscr{R}_{M,q}$ is an open subset of $J^1\Gr_q(TM)$.
\end{lemma}
\Proof
We choose a Riemannian metric $g$ on $M$ and consider the map $\Xi_g\colon J^1\Gr_q(TM) \to \R$ which assigns to each $\mathscr{V}$ the $g$-induced norm of $\Twist_{\mathscr{V}}\in\Lambda^2(V^\ast)\otimes\bot V$; here $V\define p^1_0(\mathscr{V})$ is a sub vector space of a fibre of $TM$ and thus inherits a metric, and $\bot V$ can be identified with $\bot_gV\subseteq TM$ and thus inherits a metric, too. It is a routine matter to check (via local trivialisations of $TM$) that $\Xi_g$ is continuous; we omit the details. Now $\mathscr{R}_{M,g}$ is the $\Xi_g$-preimage of the open set $\R\without\set{0}$, thus open.
\end{proof}

\begin{remark}
It is also easy to verify that the relation $\mathscr{R}_{M,g}$ is \emph{diff-invariant} (cf.\ Appendix \ref{openhprinciple}). If $M$ is open, Gromov's h-principle for open diff-invariant relations can thus be applied to prove existence of everywhere twisted $q$-plane distributions on $M$. However, we can do better by applying the h-principle for \emph{ample} relations, which holds also on closed manifolds $M$; and to prove our main theorems for the prescribed scalar curvature problem, we \emph{have} to do better even in the case of open manifolds, since we apply a $C^0$-denseness statement (cf.\ Theorem \ref{FIVEMAIN}) which we could not get out of the diff-invariant h-principle (cf.\ the remarks in \cite{EliashbergMishachev}, Chapter 7).
\end{remark}


\section{Existence of formal solutions} \label{FIVETWO}

\subsection{Formal solutions vs.\ nowhere vanishing sections in $\Lambda^2(V^\ast)\otimes\bot V$}

Let $V$ be a $q$-plane distribution on the manifold $M$, and let $p^1_0\colon J^1\Gr_q(TM)\to \Gr_q(TM)$ denote the standard projection. The aim of this subsection is to prove that the twistedness relation admits a formal solution\footnote{Recall that a \emph{formal solution} of $\mathscr{R}_{M,q}$ is a section in $J^1\Gr_q(TM)\to M$ whose image is contained in $\mathscr{R}_{M,q}$.} $\bar{V}\in C^\infty(M\ot J^1\Gr_q(TM))$ with $p^1_0\compose\bar{V} = V$ if and only if the vector bundle $\Lambda^2(V^\ast)\otimes\bot V$ admits a nowhere vanishing section.

\newcommand{\JJJ}{\mathscr{J}}
\begin{definition}[$p_V \colon \JJJ_V\to M$] \label{JVdef}
Let $M$ be a manifold, let $V$ be a $q$-plane distribution on $M$. Via the section $V\in C^\infty(M\ot\Gr_q(TM))$, we can pull back the affine bundle $p^1_0\colon J^1\Gr_q(TM)\to \Gr_q(TM)$ to an affine bundle over $M$; we denote this bundle by $p_V \colon \JJJ_V\to M$.
\smallskip\\
Let $p\colon\Gr_q(TM)\to M$ denote the bundle projection. Since the affine bundle $p^1_0$ is modelled on the vector bundle $p^\ast(T^\ast M)\otimes\ker(Tp\colon T\Gr_q(TM)\to TM)$ over $\Gr_q(TM)$, the affine bundle $p_V$ is modelled on the $V$-pullback of this vector bundle, i.e.\ on the vector bundle $T^\ast M\otimes V^\ast(\ker(Tp))$ over $M$.
\smallskip\\
The affine bundle $p_V \colon \JJJ_V\to M$ has a canonical section, namely $j^1V\in C^\infty(M\ot J^1\Gr_q(TM))$. Via $j^1V$, we can identify $p_V \colon \JJJ_V\to M$ with the vector bundle $T^\ast M\otimes V^\ast(\ker(Tp))\to M$.
\end{definition}

\newcommand{\Tw}{\mt{Tw}}
\begin{definition}[$\Tw_V\colon \JJJ_V\to \Lambda^2(V^\ast)\otimes\bot V$]
Let $M$ be a manifold, let $V$ be a $q$-plane distribution on $M$. We define a map $\Tw_V\colon \JJJ_V\to \Lambda^2(V^\ast)\otimes\bot V$ as follows. Let $\mathscr{W}\in\JJJ_V$. For $x\define p_V(\mathscr{W})\in M$, consider $\mathscr{V}\define j^1_xV\in p_V^{-1}(\set{x})$. We define $\Tw_V(\mathscr{W})\define \Twist_{\mathscr{W}}-\Twist_{\mathscr{V}} = \Twist_{\mathscr{W}}-\Twist_V(x) \in \Lambda^2(V_x^\ast)\otimes\bot V_x$.
\end{definition}

\begin{fact} \label{Twbundlemap}
Let $M$ be a manifold, let $V$ be a $q$-plane distribution on $M$. The map $\Tw_V\colon \JJJ_V\to \Lambda^2(V^\ast)\otimes\bot V$ is fibre-preserving (with respect to the bundle projections to $M$) and smooth.
\end{fact}
\Proof
$\Tw_V$ is fibre-preserving by definition. It is a routine matter to verify smoothness; we will omit the details. (It would actually suffice for our applications to check continuity.)
\end{proof}

Our aim is now to prove that $\Tw_V$ is surjective and fibrewise linear with respect to the vector bundle structures on $\JJJ_V$ and $\Lambda^2(V^\ast)\otimes\bot V$.

\begin{lemma} \label{twistsurj}
Let $M$ be an $n$-manifold, let $x\in M$, $q\in\set{0,\dots,n}$, let $V_x\in\Gr_q(T_xM)$, and let $\omega\in\Lambda^2(V_x^\ast)\otimes\bot V_x$. Then there exist an open neighbourhood $U$ of $x$ and a $q$-plane distribution $V$ on $U$ such that $V(x)=V_x$ and $\Twist_V(x)=\omega$.
\end{lemma}
\Proof
Since the statement is local, it suffices to prove it for the case $M=\R^n$, $x=0\in\R^n$. We choose a complementary sub vector space $H\subseteq T_xM=\R^n$ of $V_x$ and identify $T_xM/V_x$ with $H$. (Without loss of generality, we could restrict our considerations to the case $V_x=\R^q\times\set{0}$, $H=\set{0}\times\R^{n-q}$, but there is no advantage in doing so.) Let $\pr_V \colon \R^n = V_x\oplus H\to V_x$ and $\pr_H \colon \R^n = V_x\oplus H \to H$ denote the obvious projections.
\smallskip\\
We define the map $\hat{V}\in C^\infty(\R^n,\Lin(V_x,\R^n))$ by $\hat{V}(y)(v)\define v+\frac{1}{2}\omega(\pr_V(y),v)$. Since its value in $0$ is injective --- namely, $\hat{V}(0)$ is the inclusion $V_x\to\R^n$ ---, there is an open neighbourhood $U$ of $0\in\R^n$ such that $\hat{V}\restrict U$ is monomorphism-valued. Thus the map $V$ on $U$ which is defined to be pointwise the image of $\hat{V}$ is a $q$-plane distribution on $U$ with $V(x)=\im(\hat{V}(x))=V_x$.
\smallskip\\
It remains to prove that the twistedness of $V$ in the point $0$ satisfies $\Twist_V(v_0,v_1) = \omega(v_0,v_1)$ for all $v_0,v_1\in V_x\subseteq T_0\R^n$. For each $i\in\set{0,1}$, the function $\hat{v}_i\in C^\infty(U,\R^n)$ given by $\hat{v}_i(y)\define \hat{V}(y)(v_i)$ is a section in the vector bundle $V$ with $\hat{v}_i(0)=v_i$. Thus
\[ \begin{split}
\Twist_V(v_0,v_1) &= \pr_H([\hat{v}_0,\hat{v}_1])(x)\\
&= \pr_H\Big(\,\big[v_0+\tfrac{1}{2}\omega(\pr_V(.),v_0), v_1+\tfrac{1}{2}\omega(\pr_V(.),v_1)\big]\,\Big)(x)\\
&= \pr_H\Big(D_0\big(v_1+\tfrac{1}{2}\omega(\pr_V(.),v_1)\big)(v_0) -D_0\big(v_0+\tfrac{1}{2}\omega(\pr_V(.),v_0)\big)(v_1)\Big)\\
&= \pr_H\Big(\tfrac{1}{2}\omega(\pr_V(v_0),v_1) -\tfrac{1}{2}\omega(\pr_V(v_1),v_0)\Big)\\
&= \tfrac{1}{2}\big(\omega(v_0,v_1)-\omega(v_1,v_0)\big)\\
&= \omega(v_0,v_1) \;\;.\qedhere
\end{split} \]
\end{proof}

\begin{corollary} \label{Twsurjectivity}
Let $M$ be a manifold, let $V$ be a $q$-plane distribution on $M$. Then the map $\Tw_V\colon \JJJ_V\to \Lambda^2(V^\ast)\otimes\bot V$ is surjective.
\end{corollary}
\Proof
Let $x\in M$, let $\omega\in\Lambda^2(V_x^\ast)\otimes\bot V_x$. By the preceding lemma, there exist an open neighbourhood $U$ of $x$ and a $q$-plane distribution $W$ on $U$ such that $W(x)=V(x)$ and $\Twist_W(x) = \omega+\Twist_V(x)$. With the notation $\mathscr{V}\define j^1_xV\in p_V^{-1}(\set{x})$ and $\mathscr{W}\define j^1_xW\in p_V^{-1}(\set{x})$, we get $\Tw_V(\mathscr{W}) = \Twist_{\mathscr{W}} -\Twist_{\mathscr{V}} = \Twist_W(x) -\Twist_V(x) = \omega$.
\end{proof}

\begin{lemma} \label{Twlinearity}
Let $M$ be a manifold, let $V$ be a $q$-plane distribution on $M$. Then the map $\Tw_V\colon \JJJ_V\to \Lambda^2(V^\ast)\otimes\bot V$ is fibrewise linear: for every $x\in M$ and \mbox{$\gamma\in\Lin(T_xM,\ker(T_{V(x)}p\colon T_{V(x)}\Gr_q(TM)\to T_xM))$} and all $v,w\in V(x)$, we have (via the identifications $(p_V)^{-1}(\set{x}) = \Lin(T_xM,\ker(T_{V(x)}p))$ and $\ker(T_{V(x)}p)$ $= T_{V(x)}\Gr_q(T_xM) = \Lin(V(x),T_xM/V(x))$):
\[
\Tw_V(\gamma)(v,w) = \gamma(v)(w) -\gamma(w)(v) \in T_xM/V(x) \;\;.
\]
\end{lemma}
\Proof
The identification $(p_V)^{-1}(\set{x}) = \Lin(T_xM,\ker(T_{V(x)}p))$ comes from Definition \ref{JVdef}; the identification $\ker(T_{V(x)}p\colon T_{V(x)}\Gr_q(TM)\to T_xM)) = T_{V(x)}\Gr_q(T_xM)$ is standard (cf.\ e.g.\ \cite{Saunders}, Lemma 3.1.2); and $T_{V(x)}\Gr_q(T_xM) = \Lin(V(x),T_xM/V(x))$ is explained in Appendix \ref{Grassmannappendix}. Thus both sides of the formula are well-defined (independent of the local trivialisation of $TM$ which we are going to choose).
\smallskip\\
Since the statement of the corollary depends only on the germ of $V$ in $x$, we may assume without loss of generality that $M=\R^n$ and $x=0\in\R^n$. Let $\mathscr{V}\define j^1_0V$ and $\mathscr{W}\define \gamma+\mathscr{V}\in(\JJJ_V)_x$. We choose a map $W\in C^\infty(U,\Gr_q(\R^n))$, defined on a neighbourhood $U$ of $0\in\R^n$, with $j^1_0W = \mathscr{W}$; that is, $W(0)=V(0)$ and $D_0W = \gamma+D_0V$. Using Lemma \ref{twistedjet}, we obtain
\[ \begin{split}
\Tw_V(\gamma)(v,w) &= \Tw_V(\gamma+\mathscr{V})(v,w)\\
&= \Twist_{\gamma+\mathscr{V}}(v,w) -\Twist_{\mathscr{V}}(v,w)\\
&= \Twist_{\mathscr{W}}(v,w) -\Twist_{\mathscr{V}}(v,w)\\
&= \Twist_W(0)(v,w) -\Twist_V(0)(v,w)\\
&= (D_0W)(v)(w) -(D_0W)(w)(v) -(D_0V)(v)(w) +(D_0V)(w)(v)\\
&= \gamma(v)(w) -\gamma(w)(v) \;\;.
\end{split} \]
This proves that the formula is correct. Since it is linear in $\gamma$, the map $\Tw_V\colon \JJJ_V\to \Lambda^2(V^\ast)\otimes\bot V$ is fibrewise linear.
\end{proof}

\emph{Remark.} We could have proved \ref{Twsurjectivity} also as a corollary to the formula in \ref{Twlinearity}, but the proof we gave via Lemma \ref{twistsurj} is more direct.

\begin{lemma} \label{affbundlesect}
Let $f\colon\xi\to\eta$ be a surjective morphism in the category of finite-rank real vector bundles over some manifold\footnote{Of course the lemma remains true if we replace the category of smooth vector bundles over manifolds by topological vector bundles over arbitrary topological spaces.} $M$, and let $s$ be a section in $\eta$. Then the fibre bundle $\xi\without f^{-1}(\im(s))$ (whose fibre over $x\in M$ is $\xi_x\without f^{-1}(\set{s(x)})$) admits a section if and only if $\eta$ admits a nowhere vanishing section.
\end{lemma}
\Proof\footnote{Note that the lemma is completely obvious if $s$ is the zero section. To prove the general statement, we just have to perform translations in the fibres. If you think that this is a banality then you're probably right and should skip the proof.}
$f^{-1}(\im(s))$ is an affine bundle modelled on the vector bundle $\ker(f)$. Like every bundle with contractible fibres, it admits a section $\sigma$. The bundle $\zeta\define -\sigma +f^{-1}(\im(s))$, whose fibre over $x$ is $-\sigma(x)+f^{-1}(\set{s(x)}) = \set{-\sigma(x)+v \suchthat v\in\xi_x,\; f(v)=s(x)}$, is a sub vector bundle of $\xi$. We choose a sub vector bundle $\tilde{\eta}$ of $\xi$ which is complementary to $\zeta$.
\smallskip\\
We define a vector bundle morphism $\phi\colon \xi/\zeta \to \eta$ by $[w]\mapsto f(w)$. This map is well-defined and injective: a vector $w\in\xi_x$ is contained in $\zeta_x$ if and only if $f(w+\sigma(x))=s(x)$, i.e.\ if and only if $f(w)=0$. The map $\phi$ is also surjective since $f$ is surjective. Hence the vector bundle $\tilde{\eta}\cong \xi/\zeta$ is isomorphic to $\eta$.
\smallskip\\
If $\tilde{\eta}$ admits a nowhere vanishing section $\tilde{\sigma}$, then $\sigma+\tilde{\sigma}$ is a section in the bundle $\xi\without f^{-1}(\im(s))$: namely, if $\sigma(x)+\tilde{\sigma}(x)$ were contained in $f^{-1}(\set{s(x)})$ for some $x\in M$, then $\tilde{\sigma}(x)\in\zeta_x$ and thus $\tilde{\sigma}(x) \in \zeta_x\cap\tilde{\eta}_x = \set{0}$; that's a contradiction.
\smallskip\\
Conversely, assume that $\xi\without f^{-1}(\im(s))$ admits a section $\bar{\sigma}$. Let $\pr$ denote the projection $\xi = \zeta\oplus\tilde{\eta} \to \tilde{\eta}$. Then the section $\pr\compose(\bar{\sigma}-\sigma)$ in $\tilde{\eta}$ vanishes nowhere: otherwise we had $\bar{\sigma}(x)-\sigma(x)\in\zeta_x$ and thus $\bar{\sigma}(x)\in f^{-1}(\im(s))$ for some $x\in M$, again a contradiction.
\smallskip\\
This shows that $\tilde{\eta}$ admits a nowhere vanishing section if and only if $\xi\without f^{-1}(\im(s))$ admits a section. Since $\eta$ is isomorphic to $\tilde{\eta}$, the proof is complete.
\end{proof}

Now we have assembled all ingredients for the proof of the main result in this subsection.

\begin{proposition} \label{formalsolution}
Let $M$ be a manifold, let $V$ be a $q$-plane distribution on $M$, let $p^1_0\colon J^1\Gr_q(TM) \to \Gr_q(TM)$ denote the standard projection. Then the following statements are equivalent:
\begin{enumerate}
\item
There is a section $\bar{V}\in C^\infty(M\ot J^1\Gr_q(TM))$ which takes values in $\mathscr{R}_{M,q}$ (i.e., the twistedness relation $\mathscr{R}_{M,q}$ admits a \emph{formal solution} $\bar{V}$) such that $p^1_0\compose\bar{V} = V$.
\item
The vector bundle $\Lambda^2(V^\ast)\otimes\bot V \to M$ admits a nowhere vanishing section.
\end{enumerate}
\end{proposition}
\Proof
We apply the preceding lemma in the case where $\xi$ is the vector bundle $\JJJ_V\to M$, where $\eta$ is the vector bundle $\Lambda^2(V^\ast)\otimes\bot V \to M$, where $s$ is the section $-\Twist_V$ in $\eta$, and where $f\colon\xi\to\eta$ is the surjective vector bundle morphism $\Tw_V$; cf.\ Fact \ref{Twbundlemap}, Corollary \ref{Twsurjectivity}, and Lemma \ref{Twlinearity} for the proof that $\Tw_V$ is indeed a surjective vector bundle morphism. This shows that $\Lambda^2(V^\ast)\otimes\bot V$ admits a nowhere vanishing section if and only if the bundle $\JJJ_V\without f^{-1}(\im(s))$ over $M$ admits a section.
\smallskip\\
For every $x\in M$, the fibre of $f^{-1}(\im(s))$ over $x$ consists precisely of those $\mathscr{W}\in J^1\Gr_q(TM)$ which satisfy $p^1_0(\mathscr{W})=V(x)$ and $\Tw_V(\mathscr{W}) = -\Twist_V(x)$. The latter equation is equivalent to $\Twist_{\mathscr{W}}=0$, by the definition of $\Tw_V$. Therefore the fibre of $\JJJ_V\without f^{-1}(\im(s))$ over $x$ is $\mathscr{R}_{M,q}\cap(p^1_0)^{-1}(V(x))$. Thus a section in $\JJJ_V\without f^{-1}(\im(s))$ is the same as a section $\bar{V}\in C^\infty(M\ot\mathscr{R}_{M,q})$ with $p^1_0\compose\bar{V} = V$.
\end{proof}

\begin{remark}
If $V_0,V_1$ are contained in the same connected component of $C^\infty(M\ot\Gr_q(TM))$, then they are isomorphic as vector bundles. Hence also the vector bundles $\Lambda^2(V_0^\ast)\otimes\bot V_0$ and $\Lambda^2(V_1^\ast)\otimes\bot V_1$ are isomorphic then. This shows that if $V\in C^\infty(M\ot\Gr_q(TM))$ has a $1$-jet prolongation which is a formal solution of the twistedness relation $\mathscr{R}_{M,q}$, then every distribution in the same homotopy class as $V$ has such a $1$-jet prolongation.
\end{remark}

\subsection{The case $(n-q)(q-2)\geq2$} \label{nooriginalname}

Having proved Proposition \ref{formalsolution}, we must now check under which conditions $\Lambda^2(V^\ast)\otimes\bot V$ admits a nowhere vanishing section. Since $\Lambda^2(V^\ast)\otimes\bot V$ is the unique rank-$0$ vector bundle over the $n$-manifold $M$ if $q\in\set{0,1,n}$, it can never have a nowhere vanishing section in that case (unless $M$ is empty).
\smallskip\\
But if $3\leq q\leq n-1$ and $(n,q)\neq(4,3)$, then $\Lambda^2(V^\ast)\otimes\bot V$ admits a nonvanishing section:

\begin{proposition} \label{genericformalsolution}
Let $M$ be an $n$-dimensional manifold, and let $V$ be a $q$-plane distribution on $M$ such that $(n-q)(q-2)\geq2$. Then $\Lambda^2(V^\ast)\otimes\bot V$ admits a nowhere vanishing section.
\end{proposition}
\Proof
The vector bundle $E\define \Lambda^2(V^\ast)\otimes\bot V$ has rank $r\define \frac{q(q-1)}{2}(n-q)$. Because $(n-q)(q-2)\geq2$ implies $n-q\geq1$ and $q\geq3$, we have $r\geq3$.
\smallskip\\
We want to prove that the fibre bundle $F\define E\without(\text{image of the zero section})$ admits a section. The fibres of $F$ are homotopy equivalent to $S^{r-1}$ and thus $(r-2)$-connected, where $r-2\geq1$. If $r-2\geq n-1$, then the fibres are $(n-1)$-connected with $n\geq2$, and thus Theorem \ref{smoothobstruction} tells us that $F$ admits a section. We will show that the condition $\frac{(n-q)q(q-1)}{2}-2\geq n-1$ is equivalent to $(n-q)(q-2)\geq2$; this implies $r-2\geq n-1$ and hence the statement of the proposition.
\[ \begin{split}
\frac{(n-q)q(q-1)}{2}-2 \geq n-1 &\iff nq(q-1) -q^2(q-1) \geq 2(n+1)\\
&\iff n(q^2-q-2) \geq q^3-q^2+2\\
&\iff n(q+1)(q-2) \geq (q+1)(q^2-2q+2)\\
&\iff n(q-2) \geq q(q-2)+2\\
&\iff (n-q)(q-2) \geq 2 \;\;.\qedhere
\end{split} \]
\end{proof}

It remains to consider the cases $q=2$ and $(n,q)=(4,3)$. The latter case has to do with Lorentzian metrics on $4$-manifolds and is thus the most interesting one from a physical point of view. I will therefore concentrate on that one below. Subsection \ref{twoplane} contains a few remarks on the case $q=2$.

\subsection{$3$-plane distributions on $4$-manifolds} \label{threeplanefour}

Let us recall some (basically well-known) statements about orientability:

\begin{remark} \label{lambdaori}
If $\xi$ is any real vector bundle of rank $3$, then the vector bundle $\Lambda^2\xi$ (of rank $3$) is orientable: it has a canonical orientation, fibrewise defined by the ordered basis $(e_1\wedge e_2, e_2\wedge e_3, e_3\wedge e_1)$ of $\Lambda^2\xi_x$ which is induced by any basis $(e_1,e_2,e_3)$ of the fibre $\xi_x$.
\end{remark}
\emph{Remark.} The statement can be generalised: $\Lambda^2\xi$ is canonically oriented whenever the rank of $\xi$ is odd.
\Proof
We have to prove that any two bases $(e_1,e_2,e_3)$ and $(e_1',e_2',e_3')$ of $\xi_x$ induce the same orientation of $\Lambda^2\xi_x$; i.e., we have to show that for all $A\in\GL(\xi_x)$, the determinant of the matrix which describes the basis change from $(e_1\wedge e_2, e_2\wedge e_3, e_3\wedge e_1)$ to $(Ae_1\wedge Ae_2, Ae_2\wedge Ae_3, Ae_3\wedge Ae_1)$ is positive. It is easy to check that this determinant is $\det(A)^2$.
\end{proof}

\begin{lemma} \label{HvsLambda}
Let $H$ be a real vector bundle of rank $3$. Then there is a canonical vector bundle isomorphism $\varphi\colon H\otimes\Lambda^3(H^\ast)\to \Lambda^2(H^\ast)$, given by $\varphi(u\otimes\omega)(v,w)\define \omega(u,v,w)$. In particular, $H$ is isomorphic to $\Lambda^2(H^\ast)$ if and only if $H$ is orientable.
\end{lemma}
\Proof
The map $\varphi$ is clearly a well-defined vector bundle morphism. It is fibrewise injective: For each fibre $H_x$, we can choose a nonvanishing element $\omega$ of the $1$-dimensional vector space $\Lambda^3(H_x^\ast)$. Each element of $H_x\otimes\Lambda^3(H_x^\ast)$ has the form $u\otimes\omega$. If $\varphi(u\otimes\omega)=0$, then $u=0$ since $\omega$ is a nondegenerate $3$-form. So the kernel of $\varphi\colon H_x\otimes\Lambda^3(H_x^\ast)\to \Lambda^2(H_x^\ast)$ is trivial, as claimed. Because $H\otimes\Lambda^3(H^\ast)$ and $\Lambda^2(H^\ast)$ have both rank $3$, $\varphi$ is a vector bundle isomorphism.
\smallskip\\
If $H$ is orientable, then the line bundle $\Lambda^3(H^\ast)$ is trivial. Hence $H$ is isomorphic to $\Lambda^2(H^\ast)$. Conversely, if $H$ is not orientable, then $H$ and $\Lambda^2(H^\ast)$ are not isomorphic since $\Lambda^2(H^\ast)$ is orientable, by Remark \ref{lambdaori}.
\end{proof}

\begin{remark} \label{hhlambda}
Let $H$ be a $3$-plane distribution on a $4$-manifold $M$. Then the vector bundles $\Lambda^2(H^\ast)\otimes\bot H$ and $H\otimes\Lambda^4(T^\ast M)$ are isomorphic. In particular, if $M$ is orientable, then $\Lambda^2(H^\ast)\otimes\bot H \cong H$.
\end{remark}
\Proof
We choose a line distribution $V$ on $M$ which is complementary to $H$. Since $TM = V\oplus H$, we get $\Lambda^4(T^\ast M)\cong \Lambda^4(V\oplus H)\cong \Lambda^1(V)\otimes\Lambda^3(H) \cong \Lambda^3(H^\ast)\otimes\bot H$. The preceding lemma yields $H\otimes\Lambda^4(T^\ast M) \cong H\otimes\Lambda^3(H^\ast)\otimes\bot H \cong \Lambda^2(H^\ast)\otimes\bot H$.
\end{proof}

We will now state our main criterion for existence of formal solutions of the twistedness relation in the case of $3$-plane distributions on a $4$-manifold. The orientability discussion above indicates that the situation is easier to understand if we deal with an orientable manifold. Then $\Lambda^2(H^\ast)\otimes\bot H$ is isomorphic to $H$, i.e.\ to a subbundle of $TM$. Using this fact and assuming in addition that $H$ is orientable, we can relate existence of almost-contact structures on $M$ to existence of nowhere vanishing sections in $\Lambda^2(H^\ast)\otimes\bot H$. To avoid all complications, we will restrict our considerations to this oriented-cooriented case (which corresponds to Lorentzian metrics that are both time-orientable and space-orientable, i.e.\ time-orientable metrics on an orientable manifold).

\begin{notation} \label{intersectionform}
Let $M$ be an oriented closed $4$-manifold. We introduce the abbreviation $X_M$ for the finitely generated free $\Z$-module $H_2(M;\Z)/\Torsion(H_2(M;\Z))$. We denote the intersection form of $M$ by $\beta_M\colon X_M\times X_M\to\Z$, and the signature of $M$ (i.e.\ of $\beta_M$) by $\sigma_M\in\Z$. (Cf.\ e.g.\ Chapter 1 in \cite{GompfStipsicz} for definitions and related information; they use the notation $Q_M$ instead of $\beta_M$ there.) Recall that an element $u\in X_M$ is called \emph{characteristic} if and only if $\beta_M(u,x) \equiv \beta_M(x,x) \mod 2$ for all $x\in X_M$.
\end{notation}

Recall that an \emph{almost-complex structure} on a manifold $M$ is a complex structure on the vector bundle $TM$, i.e.\ a section $J$ in the endomorphism bundle $\End(TM)$ such that $J(x)\compose J(x) = -\id_{T_xM}$ for all $x\in M$.

\begin{proposition} \label{fourthreesection}
Let $M$ be a connected orientable $4$-manifold which admits a line distribution\footnote{Recall that a connected manifold $M$ admits a line distribution if and only if either $M$ is open, or $M$ is closed with zero Euler characteristic.}. Then $\sigma_M\in\Z$ is even\footnote{Both $\beta_M$ and $\sigma_M$ are taken with respect to some fixed orientation of $M$. The relations $\sigma_M\equiv 0\mod2$ and $\sigma_M\equiv 0\mod4$ and $\beta_M(u,u)=3\sigma_M$ do not depend on the choice of orientation since $\beta_M$ and $\sigma_M$ change their signs if the orientation is reversed.} if $M$ is closed; and the following statements are equivalent:
\begin{enumerate}
\item
There is an orientable $3$-plane distribution $H$ on $M$ such that $\Lambda^2(H^\ast)\otimes\bot H$ admits a nowhere vanishing section.
\item
For every orientable $3$-plane distribution $H$ on $M$, $\Lambda^2(H^\ast)\otimes\bot H$ admits a nowhere vanishing section.
\item
$M$ admits an almost-complex structure.
\item
Either $M$ is open; or $M$ is closed and there exists a characteristic element $u\in X_M$ such that $\beta_M(u,u)=3\sigma_M$.\addtocounter{footnote}{-1}\footnotemark
\item
Either $M$ is open; or $M$ is closed and $\sigma_M \equiv 0 \mod 4$.\addtocounter{footnote}{-1}\footnotemark
\end{enumerate}
\end{proposition}
\Proof
If $M$ is closed, then its Euler characteristic vanishes. Thus \Poincare\ duality yields $0 = 2 -2b_1 +b_2^+ +b_2^-$. This implies that $b_2^+ +b_2^-$ and hence $\sigma_M = b_2^+ -b_2^-$ are even. It remains to prove the equivalences.

\medskip
(ii)$\implies$(i): Every manifold which admits a line distribution does also admit an orientable line distribution $H$. Since $M$ is orientable, it thus admits an orientable $3$-plane distribution. Statement (ii) implies that $\Lambda^2(H^\ast)\otimes\bot H$ admits a nowhere vanishing section.
\medskip\\
(i)$\implies$(iii): Let $H$ be an orientable $3$-plane distribution on $M$ such that $\Lambda^2(H^\ast)\otimes\bot H$ admits a nowhere vanishing section. Since $H\cong \Lambda^2(H^\ast)\otimes\bot H$ by Remark \ref{hhlambda}, there exist a trivial line subbundle $L$ and an orientable $2$-plane subbundle $\xi$ of $H$ with $H= L\oplus\xi$. The bundle $\xi$ admits a complex structure since every orientable rank-$2$ vector bundle does.\footnote{This is a consequence of the fact that the Lie groups $\SO(2)$ and $\U(1)$ coincide.} Hence the vector bundle $TM\cong \bot H\oplus L\oplus\xi$ admits a complex structure, i.e., $M$ admits an almost-complex structure.
\medskip\\
(iii)$\implies$(ii): Let $J$ be an almost-complex structure on $M$, and let $H$ be an orientable $3$-plane distribution on $M$. We choose a complementary distribution $V$ of $H$; the line bundle $V$ is orientable and thus admits a nowhere vanishing section. $H$ is isomorphic to $\Lambda^2(H^\ast)\otimes\bot H$ by Remark \ref{hhlambda}, so we just have to prove that $H$ admits a nowhere vanishing section. Since all $3$-plane distributions which are complementary to $V$ are isomorphic to the vector bundle $TM/V$, it suffices to show that one such distribution splits off a trivial line bundle. We choose any $2$-plane distribution $\xi$ which is complementary to $V\oplus J(V)$ (note that the line distribution $J(V)$ is pointwise different from $V$ because of $J^2=-\id_{TM}$). Then the $3$-plane distribution $J(V)\oplus\xi$ has the required property.
\medskip\\
(iii)$\iff$(iv): The closed case follows from Satz 4.6 in \cite{HirzebruchHopf} since the Euler characteristic $\chi_M$ of $M$ is zero for every closed manifold which admits a line distribution. (The general condition is existence of a characteristic element $u\in X_M$ with $\beta_M(u,u) = 3\sigma_M+2\chi_M$.) It remains to prove that every open orientable $4$-manifold admits an almost-complex structure.
\smallskip\\
Satz 2.5 (cf.\ also the remarks at the beginning of 4.6, and the definition of $\xi_1$ in 3.1) in \cite{HirzebruchHopf} states the general obstruction-theoretic criterion for the existence of an almost-complex structure on $M$: it exists if a certain bundle $\xi_1\to M$ admits a section, and this is the case if the characteristic class $W_3(M)\in H^3(M;\Z)$ is zero and, moreover, some secondary obstruction in $H^4(M;\Z)$ vanishes.\footnote{Note that Hirzebruch and Hopf assume that $M$ can be represented as a \emph{finite} cell complex, but this assumption is not necessary. At the time when they wrote their article, obstruction theory was usually formulated only for finite complexes; cf.\ Steenrod's classic text \cite{Steenrod}. Modern treatments like \cite{Whitehead} work for arbitrary CW complexes. The arguments of Hirzebruch and Hopf remain true word by word in this general situation.}
\smallskip\\
If $M$ is open, then every secondary obstruction vanishes since $H^4(M;\Z)\cong\set{0}$.\footnote{Since the triviality of the top \emph{co}homology group for open manifolds seems to be less well-known than the corresponding homology statement, let me outline two proofs. I am grateful to Bruce Westbury and Boudewijn Moonen for pointing out these arguments to me. (Note that the main results of the present work use only the case of compact manifolds with nonempty boundary, which does not require such elaborate arguments.)
\smallskip\\
If $M$ is an open connected $n$-manifold, then $M$ has the homotopy type of an $(n-1)$-dimensional simplicial complex. Namely, we triangulate $M$, obtaining a simplicial complex $K$. We can find a smoothly and properly imbedded tree (i.e.\ graph without cycles) $T\subseteq M$, in such a way that the interior of each $n$-simplex contains exactly one $T$-vertex, that $T$ does not meet the $(n-2)$-skeleton, and that each $T$-edge intersects the $(n-1)$-skeleton in exactly one point. Clearly $M$ is diffeomorphic to $M\without\text{(closed tubular neighbourhood of $T$)}$ and thus to $M\without T$. Let $K'$ be the $(n-1)$-dimensional simplicial complex obtained from $K$ by throwing out all simplices which meet $T$. It is easy to construct a deformation retraction of $M\without T$ to $K'$. So $M$ has the homotopy type of an $(n-1)$-dimensional simplicial complex, as claimed. This implies that $H^n(M;\Z)$ is trivial.
\smallskip\\
Alternatively, one can argue that $H^n(M;\Z)$ is isomorphic to $H^\infty_0(M;\Gamma^M)$, by \Poincare\ duality. (For the definitions of the local system of $\Z$-modules $\Gamma^M$ and the locally finite homology groups $H^\infty_\ast(M;\Gamma)$, cf.\ \cite{Spanier1993}: p.~197, \S7. By the remarks on p.~183, Theorem 10.2, and the remarks on p.~188 there, we get $H_0^\infty(M;\Gamma^M)\cong {}^MH_0^\infty(M;\Gamma^M) \cong \bar{H}^n(M;\Z) \cong H^n(M;\Z)$.) So it remains to show that $H^\infty_0(M;\Gamma^M)$ is trivial if $M$ is open. But this is essentially the statement that every locally finite set of points in $M$ is the start point set of a proper imbedding into $M$ of a disjoint union of copies of $\cointerval{0}{\infty}$.}
\smallskip\\
The exact sequence $\set{0}\to\Z\xrightarrow{\cdot2}\Z\xrightarrow{r}\Z_2\to\set{0}$ induces the long exact cohomology sequence
\[
\setlength{\arraycolsep}{0.5cm}
\begin{array}{cccccc}
\Rnode{p00}{\dots}
&\Rnode{p01}{H^2(M;\Z)}
&\Rnode{p02}{H^2(M;\Z)}
&\Rnode{p03}{H^2(M;\Z_2)}
&\Rnode{p04}{H^3(M;\Z)}
&\Rnode{p05}{\dots}
\end{array}
\psset{arrows=->,nodesep=5pt,linewidth=0.3pt} \everypsbox{\scriptstyle}
\ncline{p00}{p01}
\ncline{p01}{p02}
\ncline{p02}{p03}\Aput{\bar{r}}
\ncline{p03}{p04}\Aput{\delta}
\ncline{p04}{p05}
\]
The primary obstruction $W_3(M)\in H^3(M;\Z)$ is $\delta(w_2(M))$, where $w_2(M)\in H^2(M;\Z_2)$ denotes the second Stiefel/Whitney class of $M$; cf.\ \cite{HirzebruchHopf}. But $w_2(M)$ is contained in the image of $\bar{r}$; cf.\ \cite{GompfStipsicz}, Remark 5.7.5. Thus $W_3(M)=0$, so $M$ admits an almost-complex structure.
\medskip\\
(iv)$\iff$(v): We have to show that if $M$ is closed, then existence of a characteristic element $u\in X_M$ with $\beta_M(u,u)= 3\sigma_M$ is equivalent to $\sigma_M\equiv0 \mod4$.
\smallskip\\
By the classification of indefinite nondegenerate symmetric bilinear forms on finitely generated free $\Z$-modules (cf.\ Theorem 1.2.21 in \cite{GompfStipsicz}; or \cite{MilnorHusemoller}) and S.\ Donaldson's theorem about diagonalisability of definite intersection forms of closed oriented smooth $4$-manifolds (Theorem 1 of \cite{Donaldson1987c}; cf.\ also Remark 2.4.30 in \cite{GompfStipsicz}), the bilinear form $\beta_M$ is equivalent either to $n(1)\oplus m(-1)$ for some $n,m\in\N$, or to $n\,E_8\oplus m\Hmat$ for some $n,m\in\N$ with $m\geq1$.
\smallskip\\
\emph{First we consider the case $\beta_M \cong n(1)\oplus m(-1)$.} We choose a basis $(e_1,\dots,e_{n+m})$ of $X_M$ with respect to which the bilinear form $\beta_M$ is represented by the diagonal matrix $n(1)\oplus m(-1)$. The equation $\beta_M(u,u)=3\sigma_M$ for an element $u=(u_1,\dots,u_{n+m})\in \Z^{n+m}=X_M$ means just that $\sum_{i=1}^nu_i^2-\sum_{i=n+1}^{n+m}u_i^2 = 3(n-m)$.

\smallskip
An element $u\in\Z^{n+m}=X_M$ is characteristic if and only if all numbers $u_1,\dots,u_{n+m}\in\Z$ are odd: ``Only if'' because we have $\beta_M(u,e_i)=\pm u_i$ and $\beta_M(e_i,e_i)=\pm 1$, hence $\beta_M(u,e_i)\equiv\beta_M(e_i,e_i)\mod2$ only if $u_i$ is odd. ``If'' because $\beta_M(u,x)=\sum_i\pm u_ix_i \equiv \sum_i\pm x_i^2 = \beta_M(x,x) \mod2$ when all $u_i$ are odd.
\smallskip\\
Assume that there exists a characteristic element $u\in X_M$ such that $\beta_M(u,u)=3\sigma_M$. Since $u_i^2\equiv1\mod8$ for all $i\in\set{1,\dots,n+m}$, this implies that $n-m\equiv 3(n-m)\mod 8$, i.e.\ $\sigma_M=n-m\equiv 0\mod4$.
\smallskip\\
Conversely, assume that $\sigma_M\equiv0\mod4$. We have to prove that there exist odd numbers $u_1,\dots,u_{n+m}$ with $\sum_{i=1}^nu_i^2-\sum_{i=n+1}^{n+m}u_i^2 = 3(n-m)$. If $n=m$, we can choose all $u_i$ to be $1$. The case $n<m$ follows from the case $n>m$ by reversing the roles of $n$ and $m$.
\smallskip\\
If $n>m$, then $n\geq m+4$ (because $n-m\equiv0\mod4$). In this case, we choose $u_i=1$ for all $i>4$. Every positive integer $\equiv4\mod8$ is the sum of four odd squares (cf.\ e.g.\ \cite{IrelandRosen}, Proposition 17.7.2). In particular, there exist odd numbers $u_1,\dots,u_4\in\Z$ with $\sum_{i=1}^4u_i^2 = 2(n-m)+4 = 3(n-m) -(n-4) +m$, as we wished to prove.
\smallskip\\
This completes the proof of the case $\beta_M \cong n(1)\oplus m(-1)$.
\medskip\\
\emph{Now we consider the case $\beta_M \cong nE_8\oplus m\Hmat$, where $m\geq1$.} Then $\sigma_M = 8n$ is divisible by $4$, so we just have to prove that a characteristic element $u\in X_M$ with $\beta_M(u,u)=24n$ exists. For every $v\in X_M$, the element $2v\in X_M$ is characteristic since $\beta_M(x,x)\in\Z$ is even for all $x\in X_M$. (The converse holds as well: every characteristic element is divisible by $2$. But we don't need that information.) We choose a $\beta_M$-orthogonal decomposition $X_M = Y\oplus\Z^2$ such that the restriction of $\beta_M$ to the $\Z^2$ submodule is represented by the matrix $\Hmat$. The element $v=(0,1,3n)\in Y\oplus\Z\oplus\Z = X_M$ satisfies $\beta_M(2v,2v) = 4\binom{1}{3n}^\top\Hmat\binom{1}{3n} = 24n$. So $2v\in X_M$ has the desired properties. This completes the proof.
\end{proof}

\emph{Remark.} The arguments in the proof of (iii)$\iff$(iv)$\iff$(v) are well-known in $4$-dimensional topology, but I couldn't find an explicit reference for the statement.

\begin{corollary} \label{closedD}
Let $M$ be a closed orientable $4$-manifold, let $H$ be an orientable $3$-plane distribution on $M$, let $D\subseteq M$ be the interior of a closed imbedded ball. Then the line bundle $\Lambda^2(H^\ast)\otimes\bot H$ admits a section all of whose zeroes are contained in $D$.
\end{corollary}
\emph{Remark.} We can easily arrange that the section has at most one zero, but we don't need that.
\Proof
Since the compact manifold $M\without D$ has nonempty boundary, the restriction of $\Lambda^2(H^\ast)\otimes\bot H$ to $M\without D$ admits a nowhere vanishing section $\sigma$ by Proposition \ref{fourthreesection}. We choose any smooth extension of $\sigma$ to $M$ (sections in vector bundles over a paracompact base are always extendible).
\end{proof}

\subsection{$2$-plane distributions} \label{twoplane}

Our general criterion \ref{formalsolution} re-establishes a well-known fact about contact topology in dimension $3$: Recall from the discussion in Appendix \ref{CONTACTONE} that contact structures in dimension $3$, i.e.\ \emph{solutions} of the twistedness relation for $2$-plane distributions on $3$-manifolds, exist in every homotopy class of $2$-plane distributions on \emph{orientable} $3$-manifolds, whereas they never exist on \emph{nonorientable} $3$-manifolds. While the existence part is a hard theorem in the case of closed manifolds, the following side remark shows in an elementary way that \emph{at least formal solutions} of the twistedness relation do always exist on orientable $3$-manifolds; and it proves that \emph{not even formal solutions} exist on nonorientable $3$-manifolds.
\smallskip\\
Our remark below has as a corollary the existence of contact structures in each homotopy class of $2$-plane distributions on \emph{open} connected orientable $3$-manifolds: the twistedness relation is open and diff-invariant, so Gromov's h-principle for such relations on open manifolds (cf.\ Appendix \ref{openhprinciple}) applies. Because these facts are easy to check and well-known anyway, we will not provide further details. (Note that our criterion \ref{formalsolution} for the existence of formal solutions of the contact relation on $3$-manifolds looks --- at least superficially --- not the same as the standard criterion; cf.\ e.g.\ Theorem 3.8 in \cite{Geigeshprinciple} and the remarks preceding it. The key to the translation of both approaches is of course Proposition \ref{twistvswedge}.)

\begin{remark} \label{threecontactremark}
Let $M$ be a $3$-manifold, and let $H$ be a $2$-plane distribution on $M$. Then the line bundle $\Lambda^2(H^\ast)\otimes\bot H$ over $M$ admits a nowhere vanishing section if and only if $M$ is orientable.
\end{remark}
\Proof
Let $V\define\bot H$. The line bundle $\Lambda^3(H\oplus V)$ admits a nowhere vanishing section if and only if $H\oplus V\cong TM$ is orientable, i.e.\ if and only if $M$ is orientable. Because $\rank(H)=2$ and $\rank(V)=1$, we have vector bundle isomorphisms
\[
\Lambda^3(H\oplus V)\cong \bigoplus_{j=1}^3\Lambda^j(H)\otimes\Lambda^{3-j}(V)\cong \Lambda^2(H)\otimes\Lambda^1(V)\cong \Lambda^2(H^\ast)\otimes V \;\;.\qedhere
\]
\end{proof}

On manifolds of dimension $\geq4$, we will only consider the case of \emph{orientable} $2$-plane distributions:

\begin{remark} \label{twopl}
A vector bundle $H$ of rank $2$ is orientable if and only if the line bundle $\Lambda^2(H^\ast)$ is trivial. In particular, if $H$ is an orientable $2$-plane distribution on a manifold $M$, then $\Lambda^2(H^\ast)\otimes\bot H$ admits a nowhere vanishing section if and only if $\bot H$ does. If $H$ is an orientable $2$-plane distribution on an orientable $4$-manifold $M$, then $\Lambda^2(H^\ast)\otimes\bot H$ admits a nowhere vanishing section if and only if $\bot H$ is trivial; if $M$ is closed, this can only happen if the Euler characteristic $\chi_M$ vanishes and the signature $\sigma_M$ is divisible by $4$.
\end{remark}
\Proof
The first two statements are obvious. For the third one, let $H$ be an orientable $2$-plane distribution on an orientable $4$-manifold $M$. Since $\bot H$ has rank $2$ and is orientable, it admits a nowhere vanishing section if and only if it is trivial. If $M$ is closed in that case, $\chi_M$ vanishes because $M$ admits a nowhere vanishing vector field; and $TM = H\oplus\bot H$ admits a complex structure since $H$ and $\bot H$ are orientable rank-$2$ bundles and thus admit a complex structure. Then $\sigma_M\equiv0\mod4$ by the implication (iii)$\implies$(v) of Proposition \ref{fourthreesection}.
\end{proof}


\section{Ampleness and the main results} \label{FIVETHREE}

Recall the definition \ref{ampledeftwo} of \emph{ampleness} of a first-order partial differential relation. The following proposition tells us that the twistedness relation for $q$-plane distributions on an $n$-manifold is ample except in the case corresponding to contact structures on $3$-manifolds. As I mentioned before, at least the physically interesting special case $(n,q)=(4,3)$ of the proposition is not new but due to D.\ McDuff; cf.\ \cite{McDuff1987}, Lemma 2.7. (But McDuff's motivation to consider this problem had nothing to do with general relativity.)

\begin{proposition} \label{ampleness}
Let $M$ be an $n$-manifold, let $q\in\set{0,\dots,n}$. If $(n,q)\neq(3,2)$, then the twistedness relation $\mathscr{R}_{M,q}\subseteq J^1\Gr_q(TM)$ is ample.
\end{proposition}
\Proof
Let $x\in M$, let $W$ be a codimension-$1$ sub vector space of $T_xM$, and let $e\in J^1_{\bot W}\Gr_q(TM)$. We denote by $\mathscr{F}$ the fibre over $e$ of the affine bundle $p^1_{\bot W}\colon J^1\Gr_q(TM) \to J^1_{\bot W}\Gr_q(TM)$, and we write simply $\mathscr{R}$ instead of $\mathscr{R}_{M,q}$. We have to prove that $\mathscr{R}\cap \mathscr{F}$ is an ample subset of $\mathscr{F}$.
\smallskip\\
If $\Twist_{\mathscr{V}}\neq0$ for all $\mathscr{V}\in\mathscr{F}$, then $\mathscr{R}\cap \mathscr{F} = \mathscr{F}$, hence $\mathscr{R}\cap\mathscr{F}$ is ample. So we assume that there is a $\mathscr{V}\in\mathscr{F}$ with $\Twist_{\mathscr{V}}=0$. Let $V\define p^1_0(\mathscr{V})$, where $p^1_0 \colon J^1_x\Gr_q(TM)\to \Gr_q(T_xM)$ is the standard projection.
\smallskip\\
Via the identification $\ker(T_Vp\colon T_V\Gr_q(TM)\to T_xM) = \Lin(V,T_xM/V)$ from Lemma \ref{Twlinearity}, we can consider the set $Q\subseteq\Lin(T_xM,\Lin(V,T_xM/V))$ consisting of all $\gamma\in\Lin(T_xM,\Lin(V,T_xM/V))$ such that $\gamma\restrict W = 0$ and $\Twist_{\gamma+\mathscr{V}}=0$. This $Q$ is a sub vector space of $\Lin(T_xM,\Lin(V,T_xM/V))$ because the map $\gamma\mapsto \Twist_{\gamma+\mathscr{V}}$ is linear: in fact, $\Twist_{\gamma+\mathscr{V}}\in\Lambda^2(V^\ast)\otimes(T_xM/V)$ is given by $\Twist_{\gamma+\mathscr{V}}(v,w) = \gamma(v)(w) -\gamma(w)(v)$. (This equation follows from $\Twist_{\mathscr{V}}=0$ and Lemma \ref{Twlinearity}; to apply the lemma, just extend $V$ locally to a distribution whose $1$-jet is $\mathscr{V}$).
\smallskip\\
By \ref{perpfacts}, the affine space $\mathscr{F}$ is modelled on the vector space $\set{\gamma\in\Lin(T_xM,\Lin(V,T_xM/V)) \suchthat \gamma\restrict W=0}$. Thus $B\define \set{\gamma+\mathscr{V} \suchthat \gamma\in Q} \subseteq J^1_x\Gr_q(TM)$ is actually an affine subspace of $\mathscr{F}$, and we have $\mathscr{R}\cap\mathscr{F} = \mathscr{F}\without B$. The complement of the nonempty affine subspace $B$ is ample in $\mathscr{F}$ if and only if the codimension of $B$ is not $1$; cf.\ \ref{ampleexamples}. Since $\dim(B)=\dim(Q)$, and since the dimension of $\mathscr{F}$ equals the dimension of the fibres of $\Gr_q(TM)\to M$, i.e.\ $q(n-q)$, it just remains to show that $q(n-q)-\dim(Q)\neq1$ if $(n,q)\neq(3,2)$. We distinguish two cases.
\smallskip\\
\emph{First case: $V\subseteq W$.} Then $Q$ is the whole vector space $\set{\gamma\in\Lin(T_xM,\Lin(V,T_xM/V)) \suchthat \gamma\restrict W = 0}$, because each element of that vector space satisfies $\Twist_{\gamma+\mathscr{V}}(v,w) = \gamma(v)(w) -\gamma(w)(v) = 0$ for all $v,w\in V\subseteq W$. Thus $q(n-q)-\dim(Q)=0\neq1$ in this case.
\smallskip\\
\emph{Second case: $V\not\subseteq W$.} Then $V$ is the internal direct sum of $U\define V\cap W$ and a one-dimensional vector space $L$ (a complementary vector space of $V\cap W$ in $V$ has dimension $\leq1$ since $W$ has codimension $1$ in $T_xM$, and it has dimension $\geq1$ by the assumption of the second case). We choose a nonzero vector $l_0\in L$.
\smallskip\\
If $\gamma\in Q$, then $\gamma(w)=0$ for all $w\in W$, and $\gamma(v_0)(v_1)=\gamma(v_1)(v_0)$ for all $v_0,v_1\in V$. Thus $\gamma(u)(v)=0$ for all $(u,v)\in U\times V$, and $\gamma(l,u)=\gamma(u,l)=0$ for all $(u,l)\in U\times L$. This leaves only $\gamma(l_0)(l_0)\in T_xM/V$ arbitrary. Hence $\dim(Q)\leq n-q$.
\smallskip\\
Conversely, every $\gamma\in\Lin(T_xM,\Lin(V,T_xM/V))$ satisfying $\gamma\restrict W=0$ and $\gamma(l,u)=0$ for all $(l,u)\in L\times U$ (but with arbitrary $\gamma(l_0)(l_0)$) is contained in $Q$. Hence $\dim(Q)\geq n-q$.
\smallskip\\
In the case $V\not\subseteq W$, we thus have $q(n-q)-\dim(Q) = (q-1)(n-q)$, and this is not $1$ because of $(n,q)\neq(3,2)$. This completes the proof.
\end{proof}

Now we can state the main results of the present chapter. All $C^0$-denseness statements in the following three theorems refer to the \emph{fine} $C^0$-topology on $C^\infty(M\ot\Gr_q(TM))$. (Fine $C^0$-denseness implies compact-open $C^0$-denseness, of course.)

\begin{theorem} \label{FIVEMAIN}
Let $M$ be a manifold of dimension $n\geq5$, and let $q\in\set{3,\dots,n-1}$. Then the everywhere twisted $q$-plane distributions on $M$ form a $C^0$-dense subset of the space $\Distr_q(M)$ of all $q$-plane distributions on $M$. In particular, each connected component of $\Distr_q(M)$ contains a distribution which is everywhere twisted.
\end{theorem}
\Proof
Let $V\subseteq\Distr_q(M)$. By Proposition \ref{genericformalsolution}, the vector bundle $\Lambda^2(V^\ast)\otimes\bot V$ admits a nowhere vanishing section. Thus Proposition \ref{formalsolution} implies that there is a formal solution $\bar{V}\in C^\infty(M\ot J^1\Gr_q(TM))$ of the twistedness relation $\mathscr{R}_{M,q}$ with $p^1_0\compose\bar{V}=V$. Since $\mathscr{R}_{M,q}$ is an open ample relation (cf.\ Lemma \ref{openrelation} and Proposition \ref{ampleness}), Gromov's h-principle for such relations (cf.\ Theorem \ref{hprincipleone}) tells us that every fine $C^0$-neighbourhood of $V\in C^\infty(M\ot\Gr_q(TM)) = \Distr_q(M)$ contains a distribution which is a solution of $\mathscr{R}_{M,q}$, i.e.\ a distribution with the property that the image of its $1$-jet is contained in $\mathscr{R}_{M,q}$, i.e.\ an everywhere twisted distribution.
\end{proof}

\begin{theorem} \label{FIVEMAINfourdim}
Let $M$ be a connected orientable $4$-manifold, let $\mathscr{C}$ be a connected component of $\Distr_3(M)$ consisting of orientable distributions. If $M$ either is open, or is closed with $\sigma_M\not\equiv 2\mod4$, then the everywhere twisted distributions in $\mathscr{C}$ form a $C^0$-dense subset of $\mathscr{C}$. If $M$ is closed and $D\subseteq M$ is the interior of a closed imbedded $4$-ball, then those distributions in $\mathscr{C}$ which are twisted everywhere outside $D$ form a $C^0$-dense subset of $\mathscr{C}$.
\end{theorem}
\Proof
Let $V\in\mathscr{C}\subseteq C^\infty(M\ot\Gr_3(TM))$ be orientable. If $M$ either is open, or is closed with $\sigma_M\not\equiv 2\mod4$, then $\Lambda^2(V^\ast)\otimes\bot V$ admits a nowhere vanishing section by Proposition \ref{fourthreesection}. Using the same arguments as in the proof of Theorem \ref{FIVEMAIN}, we conclude that every fine $C^0$-neighbourhood of $V$ contains an everywhere twisted distribution.
\smallskip\\
If $M$ is closed and $D\subseteq M$ is the interior of a closed imbedded $4$-ball, then we consider the following first-order partial differential relation $\mathscr{R}_D\subseteq J^1\Gr_3(TM)$: Let $p^1\colon J^1\Gr_3(TM)\to M$ denote the standard projection. We define the intersection of $\mathscr{R}_D$ with $(p^1)^{-1}(M\without D)$ to be the intersection of $\mathscr{R}_{M,3}$ with $(p^1)^{-1}(M\without D)$, and we define the intersection of $\mathscr{R}_D$ with $(p^1)^{-1}(D)$ to be the whole set $(p^1)^{-1}(D)$.
\smallskip\\
This subset $\mathscr{R}_D\subseteq J^1\Gr_3(TM)$ is open because it is the union of the open subsets $\mathscr{R}_{M,3}$ and $(p^1)^{-1}(D)$ of $J^1\Gr_3(TM)$. Moreover, it is ample: Over each point $x\in M\without D$, Proposition \ref{ampleness} implies that $\mathscr{R}_D$ is ample; and over each point $x\in D$, the relation is ample because its complement is empty.
\smallskip\\
By Corollary \ref{closedD}, $\Lambda^2(V^\ast)\otimes\bot V$ admits a section all of whose zeroes are contained in $D$. Hence Proposition \ref{formalsolution} tells us that there is a section $\bar{V}\in C^\infty(M\without D\ot J^1\Gr_3(TM))$ which takes values in $\mathscr{R}_D$ and satisfies $p^1_0\compose\bar{V} = V\restrict(M\without D)$. This $\bar{V}$ gives us a section \emph{over $M\without D$} in the affine bundle $\JJJ_V\to M$. Like every partial section in an affine bundle, this section can be extended to a continuous section over all of $M$. In other words, there is a section $\bar{V}\in C^0(M\ot J^1\Gr_3(TM))$ which takes values in $\mathscr{R}_D$ and satisfies $p^1_0\compose\bar{V} = V$; i.e., $\mathscr{R}_D$ admits a formal solution.
\smallskip\\
Now the h-principle for ample relations shows that every fine $C^0$-neighbourhood of $V$ in $C^\infty(M\ot\Gr_3(TM))$ contains a solution of $\mathscr{R}_D$, i.e.\ a section $\tilde{V}\in C^\infty(M\ot\Gr_3(TM))$ with the property that the image of its $1$-jet is a subset of $\mathscr{R}_D$. In other words, it contains a distribution which is twisted everywhere outside $D$.
\end{proof}

\begin{theorem}[$2$-plane distributions] \label{FIVEMAINtwoplane}
Let $M$ be a manifold of dimension $n\geq4$, let $\mathscr{C}$ be a connected component of $\Distr_2(M)$ consisting of orientable distributions. If one (and thus every) element of the complementary connected component\footnote{cf.\ Definition \ref{CDCdef}} $\CDC(\mathscr{C})\in\pi_0(\Distr_{n-2}(M))$ admits a nowhere vanishing section, then the everywhere twisted distributions in $\mathscr{C}$ form a $C^0$-dense subset of $\mathscr{C}$.
\end{theorem}
\Proof
Let the bundle $V\in\Distr_2(M)$ be orientable, and let the bundle $\bot V\in\Distr_{n-2}(M)$ admit a nowhere vanishing section. Then $\Lambda^2(V^\ast)\otimes\bot V$ admits a nowhere vanishing section (cf.\ Remark \ref{twopl}), so the same argument as in the proof of Theorem \ref{FIVEMAIN} shows that every fine $C^0$ neighbourhood of $V$ contains an everywhere twisted distribution.
\end{proof}

\begin{remark}
Since there is a relative version of the h-principle for ample relations, we could also prove relative versions of the preceding theorems. Here \emph{relative} means: If we are given a distribution $V$ which is already twisted in every point of a neighbourhood of some closed set $K\subseteq M$, then we can find an everywhere twisted distribution whose restriction to $K$ is $V\restrict K$. I leave to the interested reader the task to find the optimal statements in this respect.
\end{remark}

\begin{remark}
Theorem \ref{threecontact} solves the existence problem for everywhere twisted $2$-plane distributions on $3$-manifolds (i.e.\ contact structures on $3$-manifolds). Theorem \ref{FIVEMAIN} shows that many of the obstructions which exist for contact or even-contact structures in higher dimensions do not occur for everywhere twisted codimension-$1$ distributions. For instance, contact structures on a $(2n+1)$-manifold $M$ with odd $n$ can only exist if $M$ is orientable; and contact structures on a $(2n+1)$-manifold $M$ with even $n$ can only exist in homotopy classes of orientable $2n$-plane distributions (cf.\ Appendix \ref{CONTACTONE}). Neither of these conditions yields an obstruction to existence of everywhere twisted distributions if $n\geq2$.
\end{remark}


\section[$C^\infty$-approximation]{$C^\infty$-approximation by everywhere twisted distributions} \label{twistCinfty}

The present section can be skipped by readers who are only interested in the \emph{theorems} of this thesis. However, it is important for our discussion of the esc Conjecture in Section \ref{SIXesc}. For that application, it would suffice to construct $C^2$-approximations (of integrable distributions by everywhere twisted ones), but we get $C^\infty$-approximations --- in fact, $C^\infty$-\emph{deformations} --- without extra work. (Recall that we have already proved that $C^0$-approximations exist in many cases: we got that information for free from the convex integration technique.)
\smallskip\\
Because the results of this section will not be applied in the proof of the main theorems of the thesis, we will below give only sketches whenever the arguments are extremely similar to what he did in the first three sections of this chapter.
\medskip\\
Let $M$ be an $n$-manifold, let $V$ be a $q$-plane distribution on $M$, and let $H$ be an $(n-q)$-plane distribution which is complementary to $V$. Given a section $\lambda$ in the vector bundle $\Lin(V,H)\to M$, we want to find out whether there is a $t_0\in\R_{>0}$ such that for all $t\in\oointerval{0}{t_0}$, the distribution $V+t\lambda \in C^\infty(M\ot\Compl(H))$ (recall the affine structure on $\Compl(H)$ from Definition \ref{affdef}) is everywhere twisted. If this is the case for some $\lambda$, then $V$ can be approximated by everywhere twisted distributions in the compact-open $C^\infty$-topology.
\smallskip\\
To keep things simple, we consider only the case that the distribution $V$ is integrable. As we will see in Section \ref{SIXesc}, this is the most interesting case for the Lorentzian prescribed scalar curvature problem.
\medskip\\
We start by introducing a concept which should be thought of as a directional derivative of twistedness.

\begin{definition}[first-order twistedness]
Let $M$ be a manifold, let $V$ and $H$ be complementary distributions on $M$, and let $\lambda\in C^\infty(M\ot\Lin(V,H))$. Then we define the \emph{first-order twistedness of $V$ in the direction $\lambda$} to be the section $\Twist^1_{V,H,\lambda}\in C^\infty(M\ot \Lambda^2(V^\ast)\otimes H)$ which can be described as follows:
\smallskip\\
Let $\pi_V \colon TM=V\oplus H\to V$ and $\pi_H\colon TM=V\oplus H\to H$ be the obvious projections. For all sections $v,w\in C^\infty(M\ot V)$, we define
\[
\Twist^1_{V,H,\lambda}(v,w) \define \pi_H([\lambda(v),w]) +\pi_H([v,\lambda(w)]) -\lambda(\pi_V([v,w])) \;\;.
\]
This map $\Twist^1_{V,H,\lambda}$ is $C^\infty(M,\R)$-bilinear and alternating and thus yields indeed a well-defined section in the vector bundle $\Lambda^2(V^\ast)\otimes H$: alternatingness is obvious, and for every $f\in C^\infty(M,\R)$, we have
\[ \begin{split}
\Twist^1_{V,H,\lambda}(fv,w) &= f\Twist^1_{V,H,\lambda}(v,w) -df(w)\pi_H(\lambda(v)) -df(\lambda(w))\pi_H(v) +df(w)\lambda(\pi_V(v))\\
&= f\Twist^1_{V,H,\lambda}(v,w) \;\;.
\end{split} \]
\end{definition}

\begin{proposition}
Let $M$ be a compact manifold, let $V$ be an integrable $q$-plane distribution on $M$, let $H$ be a distribution on $M$ which is complementary to $V$, and let $\lambda\in C^\infty(M\ot\Lin(V,H))$. Assume that the first-order twistedness $\Twist^1_{V,H,\lambda}\in C^\infty(M\ot \Lambda^2(V^\ast)\otimes H)$ vanishes nowhere. Then each $C^\infty$-neighbourhood of $V\in C^\infty(M\ot \Gr_q(TM))$ contains a distribution which is everywhere twisted and complementary to $H$. Even more is true: For each $C^\infty$-neighbourhood $\mathscr{U}$ of $V\in C^\infty(M\ot \Gr_q(TM))$, there is a $t_0\in\R_{>0}$ such that for all $t\in\oointerval{0}{t_0}$, the distribution $V+t\lambda\in C^\infty(M\ot\Compl(H))$ is everywhere twisted and contained in $\mathscr{U}$.
\end{proposition}
\Proof
Let $\mathscr{U}$ be a $C^\infty$-neighbourhood of $V$ in $C^\infty(M\ot \Gr_q(TM))$. Since $\Compl(H)\to M$ is an affine bundle, the affine space operation $\R \to C^\infty(M\ot\Compl(H))$ given by $t\mapsto V+t\lambda$ is continuous with respect to the (compact-open) $C^\infty$-topology on the space of sections in $\Compl(H)$. Thus there exists a $t_1\in\R_{>0}$ such that $V+t\lambda\in\mathscr{U}$ for all $t\in[-t_1,t_1]$.
\smallskip\\ {}
[The intuitive idea is now to consider the Taylor expansion of the map $t\mapsto \Twist_{V+t\lambda}$ around the point $0$. (This needs clarification since $\Twist_{V+t\lambda}$ takes values in different bundles for each $t$.) Loosely speaking, the value of this map in $0$ vanishes because $V$ is integrable, and the first derivative in $0$ is the first-order twistedness of $V$ in the direction $\lambda$ (this fact explains the terminology \emph{first-order twistedness}), hence nonvanishing. Thus we expect $\Twist_{V+t\lambda}$ to be nonvanishing for sufficiently small positive values of $t$.] {}
\smallskip\\
Let $\pi_V \colon TM=V\oplus H\to V$ and $\pi_H\colon TM=V\oplus H\to H$ denote the obvious projections. For each $t\in\R$, let $\pr_H^t \colon TM=(V+t\lambda)\oplus H \to H$ denote the projection onto the second summand. Recall that $V+t\lambda = \set{v+t\lambda(v) \suchthat v\in V}$, by definition. Therefore $\pr_H^t(u) = \pi_H(u)-t\lambda(\pi_V(u))$ for all $u\in TM$, because $u = \big(\pi_V(u)+t\lambda(\pi_V(u))\big) + \big(\pi_H(u)-t\lambda(\pi_V(u))\big)$ and $\pi_V(u)+t\lambda(\pi_V(u))\in V+t\lambda$ and $\pi_H(u)-t\lambda(\pi_V(u))\in H$.

\smallskip
Via identification of $TM/(V+t\lambda)$ with $H$, we can consider $\Twist_{V+t\lambda}$ as a section in $\Lambda^2(V+t\lambda)^\ast\otimes H$. Moreover, $\Twist_{V+t\lambda}$ then defines a section $T_t$ in the bundle $\Lambda^2(V^\ast)\otimes H$ by
\[
T_t(v,w) \define \Twist_{V+t\lambda}(v+t\lambda(v),w+t\lambda(w)) \;\;.
\]
$T_t$ vanishes nowhere if and only if $\Twist_{V+t\lambda}$ vanishes nowhere.
\smallskip\\
For all local sections $v,w$ in $V$, we have $\pi_H([v,w])=0$ because $V$ is integrable. Hence, for all $t\in\R$,\pagebreak
\[ \begin{split}
T_t(v,w) &= \Twist_{V+t\lambda}(v+t\lambda(v),w+t\lambda(w))\\
&= \pr_H^t([v+t\lambda(v),w+t\lambda(w)])\\
&= \pr_H^t([v,w]) +t\pr_H^t([\lambda(v),w]) +t\pr_H^t([v,\lambda(w)]) +t^2\pr_H^t([\lambda(v),\lambda(w)])\\
&= \pi_H([v,w])-t\lambda(\pi_V([v,w]))
+t\pi_H([\lambda(v),w])-t^2\lambda(\pi_V([\lambda(v),w]))\\
&\mspace{20mu}+t\pi_H([v,\lambda(w)])-t^2\lambda(\pi_V([v,\lambda(w)]))
+t^2\pi_H([\lambda(v),\lambda(w)])-t^3\lambda(\pi_V([\lambda(v),\lambda(w)]))\\
&= t\Twist^1_{V,H,\lambda}(v,w) -t^2\Big(\lambda(\pi_V([\lambda(v),w])) +\lambda(\pi_V([v,\lambda(w)])) -\pi_H([\lambda(v),\lambda(w)])\Big)\\
&\mspace{20mu}-t^3\lambda(\pi_V([\lambda(v),\lambda(w)])) \;\;.
\end{split} \]
Note that the assignment $(v,w)\mapsto \lambda(\pi_V([\lambda(v),\lambda(w)]))$ for $v,w\in C^\infty(M\ot V)$ is $C^\infty(M,\R)$-bilinear and alternating and thus defines a section in $\Lambda^2(V^\ast)\otimes H$, which we call $\Twist^3_{V,H,\lambda}$. Hence also the assignment
\[
(v,w)\mapsto \lambda(\pi_V([\lambda(v),w])) +\lambda(\pi_V([v,\lambda(w)])) -\pi_H([\lambda(v),\lambda(w)])
\]
defines a section in $\Lambda^2(V^\ast)\otimes H$ (being the sum of three other sections in this bundle), which we call $\Twist^2_{V,H,\lambda}$.
\smallskip\\
Now we fix any Riemannian metric on $M$. This induces (fibrewise) a norm on the vector bundle $\Lambda^2(V^\ast)\otimes H$. We have to prove that there exists a $t_0$ with $0<t_0\leq t_1$ such that for all $t$ with $0<t<t_0$, the function $\norm{T_t}\in C^0(M,\R_{\geq0})$ is nowhere zero.
\smallskip\\
Since by assumption the section $\Twist^1_{V,H,\lambda}$ in $\Lambda^2(V^\ast)\otimes H$ vanishes nowhere and $M$ is compact, there is an $\varepsilon\in\R_{>0}$ such that $\norm{\Twist^1_{V,H,\lambda}}\in C^0(M,\R_{>0})$ is everywhere $\geq\varepsilon$. On the other hand, there is a $C\in\R_{>0}$ such that $\norm{\Twist^2_{V,H,\lambda}}$ and $\norm{\Twist^3_{V,H,\lambda}}$ are everywhere $\leq C$.
\smallskip\\
We define $t_0\define \min\set{t_1,1,\varepsilon/(2C)}$. For every $t\in\R$ with $0<t<t_0$, we obtain
\[ \begin{split}
\norm{T_t} &= \norm{t\Twist^1_{V,H,\lambda} +t^2\Twist^2_{V,H,\lambda} +t^3\Twist^3_{V,H,\lambda}}\\
&\geq t\bigg(\norm{\Twist^1_{V,H,\lambda}} -t\norm{\Twist^2_{V,H,\lambda}} -t^2\norm{\Twist^3_{V,H,\lambda}}\bigg)\\
&\geq t\Big(\varepsilon -tC -t^2C\Big) \geq t\Big(\varepsilon -2tC\Big)\\
&> 0 \;\;.
\end{split} \]
Hence $t_0$ has the desired property, and we are done.
\end{proof}

\begin{remark}
One can in general not expect that the solutions of an open first-order partial differential relation $\mathscr{R}\subseteq J^1E$ (where $E\to M$ is any fibre bundle) form a $C^1$-dense subset of $C^\infty(M\ot E)$, let alone a $C^\infty$-dense subset. (In particular, there is no general h-principle theorem making such a statement.) Namely, if $\mathscr{R}$ is not dense in $J^1E$, and $\sigma\in C^\infty(M\ot E)$ is any section such that not the whole image of $j^1\sigma$ is contained in the closure of $\mathscr{R}$, then $\sigma$ can obviously not be $C^1$-approximated by solutions of $\mathscr{R}$. But in the case of the twistedness relation $\mathscr{R}_{M,q}\in J^1\Gr_q(TM)$, the trivial necessary condition for $C^1$-denseness (or $C^\infty$-denseness) of solutions is satisfied: the closure of $\mathscr{R}_{M,q}$ is the whole $1$-jet manifold $J^1\Gr_q(TM)$.
\end{remark}

Now we have to find a $\lambda\in C^\infty(M\ot\Lin(V,H))$ such that $\Twist^1_{V,H,\lambda}$ vanishes nowhere. This can again be done via the usual h-principle for open ample relations, following the same recipe as in our analysis of the $\Twist_V\neq0$ relation. I give only a sketch of the argument; it is straightforward to fill in the details.

\begin{lemma}
Let $M$ be a manifold, let $V$ and $H$ be complementary distributions on $M$, let $x\in M$, let $\lambda_0,\lambda_1\in C^\infty(M\ot\Lin(V,H))$. If $j^1_x\lambda_0 = j^1_x\lambda_1$, then $\Twist^1_{V,H,\lambda_0}(x) = \Twist^1_{V,H,\lambda_1}(x)$.
\end{lemma}
\Proof
Omitted (very similar to the proof of Corollary \ref{twistlemma}).
\end{proof}

\begin{definition}[$\Twist^1_{V,H}$ on the jet level]
Let $M$ be a manifold, let $V$ and $H$ be complementary distributions on $M$, and let $\varLambda\in J^1\Lin(V,H)$. We define $\Twist^1_{V,H,\varLambda}$, the \emph{first-order twistedness of $\varLambda$}, as follows. Let $p^1\colon J^1\Lin(V,H)\to M$ denote the standard bundle projection, and let $x\define p^1(\varLambda)\in M$. Then $\Twist^1_{V,H,\varLambda}\in \Lambda^2(V_x^\ast)\otimes H_x$ is given by $\Twist^1_{V,H,\varLambda}\define \Twist^1_{V,H,\lambda}(x)$, where $\lambda\in C^\infty(U\ot\Lin(V,H))$ is any local section on a neighbourhood $U$ of $x$ in $M$ such that $j^1_x\lambda=\varLambda$.
\end{definition}

\begin{definition}[the first-order twistedness relation $\mathscr{R}^1_{V,H}$]
Let $M$ be a manifold, let $V$ and $H$ be complementary distributions on $M$. We define $\mathscr{R}^1_{V,H}$ to be the subset of $J^1\Lin(V,H)$ consisting of all $\varLambda$ with $\Twist^1_{V,H,\varLambda}\neq 0\in\Lambda^2(V_x^\ast)\otimes H_x$, where $\varLambda\in J^1_x\Lin(V,H)$. We call $\mathscr{R}^1_{V,H}$ the \emph{first-order twistedness relation}.
\end{definition}

\begin{remark}
$\mathscr{R}^1_{V,H}$ is an open subset of $J^1\Lin(V,H)$ (by the same argument as in \ref{openrelation}). But except in trivial cases (e.g.\ $M=\leer$), it is not diff-invariant; so even on open manifolds, we need the h-principle for ample relations to prove that solutions of $\mathscr{R}^1_{V,H}$ exist.
\end{remark}

\begin{definition}
Let $M$ be a manifold, let $V$ and $H$ be complementary distributions on $M$. Via the zero section in the vector bundle $\Lin(V,H)\to M$, we can pull back the vector bundle $J^1\Lin(V,H)\to \Lin(V,H)$ to a vector bundle over $M$; we denote this pull-back bundle by $\JJJ^1_{V,H}\to M$.
\smallskip\\
We define a morphism $\Tw^1_{V,H}\colon \JJJ^1_{V,H}\to \Lambda^2(V^\ast)\otimes H$ in the category of vector bundles over $M$ by sending each $\varLambda\in(\JJJ^1_{V,H})_x$ to $\Twist^1_{V,H,\varLambda}\in\Lambda^2(V_x^\ast)\otimes H_x$. (Fibrewise linearity is obvious here.)
\end{definition}

\begin{lemma}
Let $M$ be a manifold, let $V$ and $H$ be complementary distributions on $M$. The vector bundle morphism $\Tw^1_{V,H}\colon \JJJ^1_{V,H}\to \Lambda^2(V^\ast)\otimes H$ is surjective.
\end{lemma}
\Proof
Omitted (very similar to the proof of \ref{Twsurjectivity}).
\end{proof}

\begin{proposition}
Let $M$ be a manifold, let $V$ and $H$ be complementary distributions on $M$, and let $p^1_0:$ $J^1\Lin(V,H) \to \Lin(V,H)$ denote the standard projection. Then the following statements are equivalent:
\begin{enumerate}
\item
There is a section $\bar{\lambda}\in C^\infty(M\ot J^1\Lin(V,H))$ which takes values in $\mathscr{R}^1_{V,H}$ (i.e., the twistedness relation $\mathscr{R}^1_{V,H}$ admits a \emph{formal solution} $\bar{\lambda}$) such that $p^1_0\compose\bar{\lambda}\in C^\infty(M\ot\Lin(V,H))$ is the zero section.
\item
The vector bundle $\Lambda^2(V^\ast)\otimes H \to M$ admits a nowhere vanishing section.
\end{enumerate}
\end{proposition}
\Proof
Like \ref{formalsolution}, this is now a straightforward application of Lemma \ref{affbundlesect}.
\end{proof}

This tells us that the criterion for the existence of formal solutions of $\mathscr{R}^1_{V,H}$ is exactly the same as for the existence of formal solutions of $\mathscr{R}_{M,q}$. So we can apply our results from the subsections \ref{nooriginalname}, \ref{threeplanefour}, and \ref{twoplane} again.

\begin{proposition}
Let $M$ be an $n$-manifold, let $V$ and $H$ be complementary distributions on $M$, where $\rank(V)=q$. If $(n,q)\neq(3,2)$, then the first-order twistedness relation $\mathscr{R}^1_{V,H}\subseteq J^1\Lin(V,H)$ is ample.
\end{proposition}
\Proof
Omitted (very similar to the proof of \ref{ampleness}).
\end{proof}

Putting together the preceding results, mutatis mutandis as before in Section \ref{FIVETHREE}, we get theorems which look precisely like Theorems \ref{FIVEMAIN}, \ref{FIVEMAINfourdim}, and \ref{FIVEMAINtwoplane} --- except that we have to assume this time that the manifold $M$ is compact, and can deduce that integrable distributions are not only $C^0$-approximable but even $C^\infty$-approximable by everywhere twisted ones.

\begin{theorem}
Let $M$ be a compact manifold of dimension $n\geq5$, let $q\in\set{3,\dots,n-1}$, and let $V$ be an integrable $q$-plane distribution on $M$. Then every $C^\infty$-neighbourhood of $V$ in $\Distr_q(M)$ contains an everywhere twisted distribution.
\end{theorem}

\begin{theorem} \label{fourCinfty}
Let $M$ be a connected orientable $4$-manifold, and let $V$ be an orientable integrable $3$-plane distribution on $M$. If $M$ either is open, or is closed with $\sigma_M\not\equiv 2\mod4$, then every $C^\infty$-neighbourhood of $V$ in $\Distr_3(M)$ contains an everywhere twisted distribution. If $M$ is closed and $D\subseteq M$ is the interior of a closed imbedded $4$-ball, then every $C^\infty$-neighbourhood of $V$ in $\Distr_3(M)$ contains a distribution which is twisted everywhere outside $D$.
\end{theorem}

\begin{theorem}[$2$-plane distributions]
Let $M$ be a manifold of dimension $n\geq4$, and let $V$ be an orientable integrable $2$-plane distribution on $M$. If $\bot V$ admits a nowhere vanishing section, then every $C^\infty$-neighbourhood of $V$ in $\Distr_2(M)$ contains an everywhere twisted distribution.
\end{theorem}

\begin{remark}
The situation for integrable $2$-plane distributions on compact (orientable) $3$-manifolds is more complicated. For instance, the second-factor distribution on $S^1\times S^2$ is not even $C^0$-approximable by everywhere twisted distributions. For this and further results, cf.\ \cite{EliashbergThurston}, \S2.4 and \S2.9.
\end{remark}


\chapter{Solutions on manifolds of dimension $\geq3$} \label{SIX}

This chapter contains the main results of the thesis, namely all the theorems about the pseudo-Riemannian prescribed scalar curvature problem on manifolds of dimension $n\geq3$. Metrics of index $\in\set{3,\dots,n-3}$ are quite easy to deal with, given our results from Chapter \ref{FIVE}.
\smallskip\\
For metrics of index $1$ or $2$, we prove that at least those functions are scalar curvatures which are positive somewhere (on each connected component of the manifold under consideration). In the proofs, we have to distinguish two cases. The first case consists of functions which are positive everywhere, the second case takes care of functions which are positive somewhere but zero somewhere else.
\smallskip\\
In the first case, we apply the same analytic technique as in the case of metrics of index $\in\set{3,\dots,n-3}$: the method of sub- and supersolutions. If applicable, this method gives the best results on the distribution version of the prescribed scalar curvature problem, because the prescribed distributions have to be perturbed only slightly.
\smallskip\\
In the second case, we employ the technique that Kazdan and Warner developed for the Riemannian prescribed scalar curvature problem: it uses the implicit function theorem for Banach spaces and a theorem about $L^p$ approximation of functions via pullback by diffeomorphisms (cf.\ Appendix \ref{D2}). The diffeomorphisms can be chosen from the identity component of the diffeomorphism group, but they are usually not close to the identity. That's why prescribed distributions have to be perturbed considerably in this approach.
\medskip\\
The analytic techniques employed in this chapter are not very elaborate. We have to solve a quite complicated partial differential equation (cf.\ Theorem \ref{THEEQUATION}), containing many terms about which we don't have much information in general. So our philosophy is to get by on the simplest criteria which imply that solutions exist. The present work contains only those results that these simple criteria yield.
\smallskip\\
To do better, one would have to do more analysis in order to prove better criteria for the existence of solutions, or one would have to do more geometry in order to arrange that the terms in the equation have nice properties. Maybe one would have to do both.
\smallskip\\
Some rough ideas about the ``more geometry'' approach are contained in Section \ref{SIXesc}. There we discuss the esc Conjecture which, if true, would solve completely the homotopy class version of the Lorentzian prescribed scalar curvature problem on manifolds of dimension $\geq4$.

\section{The method of sub- and supersolutions} \label{SIXsusu}

Recall the definition \ref{PDdef} of the second-order differential operator $\PD_{g,V,s}\colon C^\infty(M,\R_{>0})\to C^\infty(M,\R)$. We want to employ the method of sub- and supersolutions (cf.\ Theorem \ref{choquetbruhatleray}) in order to prove that the equation $\PD_{g,V,s}(f)=0$ has a solution $f\in C^\infty(M,\R_{>0})$. Since the symbol of the semilinear operator $\PD_{g,V,s}$ equals the symbol of $2\laplace_g$ (where $g$ is a Riemannian metric), $\PD_{g,V,s}$ is positively elliptic in the terminology of Appendix \ref{AppendixB}. Thus $f\in C^\infty(M,\R_{>0})$ is a supersolution of our elliptic equation if and only if $0 \geq \PD_{g,V,s}(f)$; and it is a subsolution if and only if $0 \leq \PD_{g,V,s}(f)$.
\smallskip\\
According to the philosophy mentioned above, we try to get by on the simplest criterion for the existence of sub- and supersolutions one can imagine: we discuss under which conditions our equation admits \emph{constant} sub- and supersolutions. No other functions will be used as sub- or supersolutions in the present thesis.

\subsection{Metrics of index $\in\set{3,\dots,n-3}$}

In the statements of the following lemmata, we identify positive real numbers with constant positive functions on a manifold $M$.

\begin{lemma}[supersolutions from twistedness] \label{supsol1}
Let $M$ be a compact manifold, let $g$ be a Riemannian metric on $M$, let $V$ be an everywhere twisted $q$-plane distribution on $M$, let $s\in C^\infty(M,\R)$. Then there is a number $c_+\in\R_{>0}$ such that every constant $c \geq c_+$ satisfies $0 > \PD_{g,V,s}(c)$.
\end{lemma}
\Proof
Let $n\define\dim(M)$, and let $H$ denote the $g$-orthogonal distribution of $V$. We have $2\leq q\leq n-1$ because $V$ is everywhere twisted.
\smallskip\\
Since $V$ is everywhere twisted and $M$ is compact, there exists an $\eps\in\R_{>0}$ such that $\abs{\Twist_V}_g^2$ is everywhere $\geq\eps$. Moreover, there exists a $C\in\R_{>0}$ such that the absolute value of each of the functions $\abs{\Twist_H}_g^2$, $\xi_{g,V}$, $\scal_g$, $s$ is everywhere $\leq C$. Finally, there exists a $c_+\in\R_{>0}$ such that for every $x\in\R$ with $x \geq c_+$, we have
\[
\frac{\eps}{5C}\cdot\frac{x(1+x^2)^2}{2} \geq \max\left\{\frac{(1+x^2)^2}{2x^3},\; \frac{(1+x^2)^2}{x},\; \frac{1+x^2}{x},\; x^{\frac{2q}{n-1}-1}(1+x^2)^{1-\frac{1}{n-1}}\right\} \;\;.
\]
(This is of course the main point of the lemma: considered as a function in $f$, the absolute value of the coefficient of $\abs{\Twist_V}_g^2$ in $\PD_{g,V,s}(f)$ grows faster than the absolute values of all the other zeroth-order coefficients as $f$ tends to infinity.)
\smallskip\\
Hence we obtain for every constant $c \geq c_+$:
\[ \begin{split}
\PD_{g,V,s}(c) &= \frac{(1+c^2)^2}{2c^3}\abs{\Twist_H}^2_g -\frac{c(1+c^2)^2}{2}\abs{\Twist_V}^2_g +\frac{(1+c^2)^2}{c}\xi_{g,V} +\frac{1+c^2}{c}\scal_g\\
&\mspace{20mu}-c^{\frac{2q}{n-1}-1}(1+c^2)^{1-\frac{1}{n-1}}s\\
&\leq -\frac{c(1+c^2)^2}{2}\abs{\Twist_V}^2_g +4C\cdot\frac{\eps}{5C}\cdot\frac{c(1+c^2)^2}{2}\\
&\leq -\eps\cdot\frac{c(1+c^2)^2}{2}\Big(1-\frac{4}{5}\Big)\\
&<0 \;\;,
\end{split} \]
as claimed.
\end{proof}

\begin{lemma}[subsolutions from twistedness] \label{subsol1}
Let $M$ be a compact manifold, let $g$ be a Riemannian metric on $M$, let $V$ be a $q$-plane distribution on $M$, let $s\in C^\infty(M,\R)$, and assume that the $g$-orthogonal distribution $H$ of $V$ is everywhere twisted. Then there is a number $c_-\in\R_{>0}$ such that every constant $c\in\R_{>0}$ with $c \leq c_-$ satisfies $0 < \PD_{g,V,s}(c)$.
\end{lemma}
\Proof
Let $n\define\dim(M)$. We have $1\leq q\leq n-2$ because $H$ is everywhere twisted.
\smallskip\\
Since $H$ is everywhere twisted and $M$ is compact, there exists an $\eps\in\R_{>0}$ such that $\abs{\Twist_H}_g^2$ is everywhere $\geq\eps$. Moreover, there exists a $C\in\R_{>0}$ such that the absolute value of each of the functions $\abs{\Twist_V}_g^2$, $\xi_{g,V}$, $\scal_g$, $s$ is everywhere $\leq C$. Finally, there exists a $c_-\in\R_{>0}$ such that for every $x\in\R$ with $0 < x \leq c_-$, we have
\[
\frac{\eps}{5C}\cdot\frac{(1+x^2)^2}{2x^3} \geq \max\left\{\frac{x(1+x^2)^2}{2},\; \frac{(1+x^2)^2}{x},\; \frac{1+x^2}{x},\; x^{\frac{2q}{n-1}-1}(1+x^2)^{1-\frac{1}{n-1}}\right\} \;\;.
\]
(That's the main point here: the coefficient of $\abs{\Twist_H}_g^2$ in $\PD_{g,V,s}(f)$ grows faster than the absolute values of all the other zeroth-order coefficients as $f$ tends to zero.)
\smallskip\\
Hence we obtain for every constant $c \leq c_-$:
\[ \begin{split}
\PD_{g,V,s}(c) &= \frac{(1+c^2)^2}{2c^3}\abs{\Twist_H}^2_g -\frac{c(1+c^2)^2}{2}\abs{\Twist_V}^2_g +\frac{(1+c^2)^2}{c}\xi_{g,V} +\frac{1+c^2}{c}\scal_g\\
&\mspace{20mu}-c^{\frac{2q}{n-1}-1}(1+c^2)^{1-\frac{1}{n-1}}s\\
&\geq \eps\cdot\frac{(1+c^2)^2}{2c^3} -4C\cdot\frac{\eps}{5C}\cdot\frac{(1+c^2)^2}{2c^3}\\
&>0 \;\;.\qedhere
\end{split} \]
\end{proof}

Now we can prove our main result about metrics of index $q\in\set{3,\dots,n-3}$. Recall that $\Metr_q(M)$ denotes the space of all index-$q$ metrics on $M$, and that $\Distr_q(M)$ denotes the space of all $q$-plane distributions on $M$.

\begin{theorem} \label{main1}
Let $M$ be a compact $n$-manifold, let $q\in\set{3,\dots,n-3}$, and let $s\in C^\infty(M,\R)$. Then every connected component of $\Metr_q(M)$ contains a metric with scalar curvature $s$.
\smallskip\\
Moreover, let $V$ be a $q$-plane distribution on $M$, let $H$ be an $(n-q)$-plane distribution which is complementary to $V$, let $\mathscr{V}\subseteq\Distr_q(M)$ be a $C^0$-neighbourhood of $V$, and let $\mathscr{H}\subseteq\Distr_{n-q}(M)$ be a $C^0$-neighbourhood of $H$. Then there is a pseudo-Riemannian metric of index $q$ on $M$ with scalar curvature $s$ which makes some element of $\mathscr{V}$ timelike and makes some element of $\mathscr{H}$ spacelike.
\end{theorem}
\Proof
Because the map $\TMC$ which we defined in Appendix \ref{AppendixC} is a bijection between the sets of connected components of $\Distr_q(M)$ resp.\ $\Metr_q(M)$, and because connected components are $C^0$-open, it suffices to prove the second statement. So we consider $V$, $H$, $\mathscr{V}$, $\mathscr{H}$ as in the statement of the theorem.
\smallskip\\
The set of all distributions which are complementary to $H$ is a $C^0$-open neighbourhood of $V$ in $\Distr_q(M)$. By Theorem \ref{FIVEMAIN}, its intersection with $\mathscr{V}$ contains an everywhere twisted distribution $V'$ (because $3\leq q\leq n-1$). The set of all distributions which are complementary to $V'$ is a $C^0$-open neighbourhood of $H$ in $\Distr_{n-q}(M)$. Its intersection with $\mathscr{H}$ contains an everywhere twisted distribution $H'$, again by Theorem \ref{FIVEMAIN} (because $3\leq n-q\leq n-1$). We will prove that there is a pseudo-Riemannian metric of index $q$ on $M$ with scalar curvature $s$ which makes $V'$ timelike and $H'$ spacelike.
\smallskip\\
We choose a Riemannian metric $g$ on $M$ which makes $V'$ and $H'$ orthogonal to each other. By Theorem \ref{THEEQUATION}, we are done if we can prove that the elliptic PDE $\PD_{g,V',s}(f) = 0$ has a solution $f\in C^\infty(M,\R_{>0})$. Lemma \ref{supsol1} tells us that it has a constant supersolution $f_+\in C^\infty(M,\R_{>0})$, and Lemma \ref{subsol1} says that every sufficiently small constant $f_-\in C^\infty(M,\R_{>0})$ is a subsolution. In particular, we may assume that $f_- < f_+$. Since our PDE has a form to which the method of sub- and supersolutions applies (cf.\ \ref{choquetbruhatleray} and \ref{susuexample}), we conclude that there is a function $f\in C^\infty(M,\R_{>0})$ with $f_-\leq f\leq f_+$ such that $\PD_{g,V',s}(f) = 0$.
\end{proof}

\begin{remark}
The proof of the preceding theorem shows that when $V$ and/or $H$ is everywhere twisted, then we can find a suitable metric which makes $V$ itself timelike and/or $H$ itself spacelike.
\smallskip\\
The discussion of $C^\infty$-denseness of everywhere twisted distributions in Section \ref{twistCinfty} makes it reasonable to assume that the theorem remains true in many (probably all) cases when we replace $C^0$-neighbourhoods by $C^\infty$-neighbourhoods. In fact, the results of that section imply already that we can replace $C^0$ by $C^\infty$ when $V$ and $H$ are integrable. It should be easy to analyse the general situation; I just haven't done that yet.
\end{remark}

\begin{theorem}[real-analytic version] \label{main1ra}
Let $M$ be a compact real-analytic $n$-manifold, let $3\leq q\leq n-3$, and let $s\colon M\to\R$ be real-analytic. Then every connected component of $\Metr_q(M)$ contains a real-analytic metric with scalar curvature $s$.
\smallskip\\
Moreover, let $V$ be a $q$-plane distribution on $M$, let $H$ be an $(n-q)$-plane distribution which is complementary to $V$, let $\mathscr{V}\subseteq\Distr_q(M)$ be a $C^0$-neighbourhood of $V$, and let $\mathscr{H}\subseteq\Distr_{n-q}(M)$ be a $C^0$-neighbourhood of $H$. Then there is a real-analytic pseudo-Riemannian metric of index $q$ on $M$ with scalar curvature $s$ which makes some element of $\mathscr{V}$ timelike and makes some element of $\mathscr{H}$ spacelike.
\end{theorem}
\Proof
The proof is the same as that of Theorem \ref{main1}, with the following additional observations: By Corollary \ref{radense} (applied to sections in the real-analytic Grassmann bundles $\Gr_q(TM)\to M$ and $\Gr_{n-q}(TM)\to M$), the distributions $V'$ and $H'$ can be chosen real-analytic. Again by \ref{radense} (applied to the real-analytic bundles $\Sym_0(V')\to M$ and $\Sym_0(H')\to M$), these real-analytic vector bundles admit real-analytic Riemannian metrics, and these metrics define a real-analytic Riemannian metric $g$ on the manifold $M$ which makes $V'$ orthogonal to $H'$. Thus the elliptic PDE $\PD_{g,V',s}=0$ has real-analytic coefficients. Its smooth solution $f$ is therefore in fact real-analytic; cf.\ Theorem \ref{hopf}. Hence the index-$q$ metric $h\define \change(g,f,K\compose f,V')$ which solves our problem (cf.\ Theorem \ref{THEEQUATION}) is real-analytic (because $K\in C^\infty(\R_{>0},\R)$ is real-analytic).
\end{proof}

\begin{remark}[lower regularity]
When the prescribed function $s$ is not smooth but only contained in some Hölder space $C^{k,\alpha}(M,\R)$ (with $k\in\N$ and $\alpha\in\oointerval{0}{1}$) or some Sobolev space $\Sob{p}{k}(M,\R)$ (with $k\in\N$ and $p\in\R_{>1}$), then still something can be said about existence of solution metrics of the prescribed scalar curvature problem. Namely, in the proof of \ref{main1}, the coefficients of our PDE $\PD_{g,V',s}=0$ are then contained in $C^{k,\alpha}(M,\R)$ resp.\ $\Sob{p}{k}(M,\R)$ (we can still choose smooth $V'$ and $g$).
\smallskip\\
The sub- and supersolution theorem in \cite{ChoquetBruhatLeray} tells us in the Hölder case that the PDE has a solution $f\in C^{k+2,\alpha}(M,\R_{>0})$, so our pseudo-Riemannian solution metric $\change(g,f,K\compose f,V')$ has regularity $C^{k+2,\alpha}$ as well.
\smallskip\\
In the Sobolev case, we can apply the sub- and supersolution theorem 6.5 from \cite{KazdanKramer}, provided $p>n$ and $s\in L^\infty(M,\R)$; the latter assumption is needed in the proofs of Lemma \ref{supsol1} and \ref{subsol1}. This yields an everywhere positive solution $f\in\Sob{p}{2}(M,\R)\subseteq C^1(M,\R)$ which, by elliptic regularity, is even contained in $\Sob{p}{k+2}(M,\R)$. So we get a solution metric of regularity $\Sob{p}{k+2}$.
\smallskip\\
Similar remarks apply
to all the other theorems below which are proved via the sub- and supersolution method.
\end{remark}

\subsection{Metrics of index $1$ or $2$}

The approach from the preceding subsection works for metrics of index $2$ or $n-2$ as well, but leads to weaker results.

\begin{theorem}[metrics of index $2$ in dimension $\geq5$, no restriction on $s$] \label{main2a}
Let $n\geq5$, let $M$ be a compact $n$-manifold, let $s\in C^\infty(M,\R)$, let $V$ be an orientable $2$-plane distribution on $M$, let $H$ be an $(n-2)$-plane distribution which is complementary to $V$ and admits a nowhere vanishing section. Let $\mathscr{V}\subseteq\Distr_2(M)$ be a $C^0$-neighbourhood of $V$, and let $\mathscr{H}\subseteq\Distr_{n-2}(M)$ be a $C^0$-neighbourhood of $H$. Then there is a pseudo-Riemannian metric of index $2$ on $M$ with scalar curvature $s$ which makes some element of $\mathscr{V}$ timelike and makes some element of $\mathscr{H}$ spacelike. There is a pseudo-Riemannian metric of index $n-2$ on $M$ with scalar curvature $s$ which makes some element of $\mathscr{V}$ spacelike and makes some element of $\mathscr{H}$ timelike.
\end{theorem}
\Proof
The proof is the same as that of Theorem \ref{main1}, except that we use Theorem \ref{FIVEMAINtwoplane} for the existence of the everywhere twisted distribution $V'$.
\end{proof}

\emph{Remark.} What would we have to do to get rid of the orientability assumption on $V$ and the assumption on $TM/V \cong H$ to split off a trivial line bundle? I.e., what would we have to do to generalise (at least the first paragraph of) Theorem \ref{main1} to metrics of index $2$ (and $n-2$), provided $n\geq5$? Since arbitrarily small constant subsolutions exist, it would suffice to arrange that there is a constant supersolution. First, we should make the zero set of $\Twist_V$ as small as possible. By a dimension count as in the proof of Proposition \ref{genericformalsolution}, we see that the twistedness relation has a formal solution over (an open neighbourhood of) the $(n-3)$-skeleton of $M$ (when we have represented $M$ as a simplicial or CW complex). Hence standard arguments from Chapter \ref{FIVE} show that there is a $2$-plane distribution $V$ which is twisted on a neighbourhood of the $(n-3)$-skeleton. When $f$ tends to $\infty$, the second-fastest growing coefficient in the zeroth-order terms of $\PD_{g,V,s}(f)$ is the coefficient of $\xi_{g,V}$. So we would get a constant supersolution if $\xi_{g,V}$ were negative on the zero set of $\Twist_V$. Assume we could prove that, given any $2$-plane distribution $V$ on $M$ and a tubular neighbourhood $U$ of the $(n-3)$-skeleton of $M$, there exists a Riemannian metric $g$ on $M$ with $\xi_{g,V}\restrict(M\without U) <0$; then the desired generalisation of the first paragraph of Theorem \ref{main1} would be true.

\medskip
In dimension $4$, the approach from the preceding subsection works at least for parallelisable manifolds:

\begin{theorem}[metrics of index $2$ in dimension $4$, no restriction on $s$] \label{main2b}
Let $M$ be a compact $4$-manifold, let $s\in C^\infty(M,\R)$, let $V,H$ be complementary $2$-plane distributions on $M$ which are trivial as vector bundles (so in particular $M$ is parallelisable). Let $\mathscr{V}\subseteq\Distr_2(M)$ be a $C^0$-neighbourhood of $V$, and let $\mathscr{H}\subseteq\Distr_2(M)$ be a $C^0$-neighbourhood of $H$. Then there is a pseudo-Riemannian metric of index $2$ on $M$ with scalar curvature $s$ which makes some element of $\mathscr{V}$ timelike and makes some element of $\mathscr{H}$ spacelike.
\end{theorem}
\Proof
The proof is the same as that of Theorem \ref{main1}, except that we use Theorem \ref{FIVEMAINtwoplane} for the existence of the everywhere twisted distributions $V'$ and $H'$.
\end{proof}

We have to use a different argument for Lorentzian metrics. In that case, Lemma \ref{subsol1} still gives us a subsolution in many cases, but Lemma \ref{supsol1} yields no supersolution anymore because line distributions are always integrable (i.e.\ \emph{nowhere} twisted instead of everywhere twisted). We will discuss in Section \ref{SIXesc} a strategy how to overcome that problem in general. For the moment, let us collect the fruits within easy reach.

\begin{lemma}[supersolution from a sign condition on $s$] \label{supsol2}
Let $M$ be a compact manifold, let $g$ be a Riemannian metric on $M$, let $V$ be a $q$-plane distribution on $M$. Then there exists a constant $s_0\in\R$ such that every function $s\in C^\infty(M,\R)$ with $s\geq s_0$ satisfies the inequality $0>\PD_{g,V,s}(1)$.
\end{lemma}
\Proof
By the definition of $\PD_{g,V,s}$, the function $u_{g,V}\define \PD_{g,V,s}(1)+ 2^{1-1/(n-1)}s \in C^\infty(M,\R)$ is independent of $s$. Since $M$ is compact, there is an $s_0\in\R$ such that $\PD_{g,V,s}(1) = u_{g,V} -2^{1-1/(n-1)}s$ is everywhere negative whenever $s\geq s_0$.
\end{proof}

We state the following theorem not only for metrics of index $1$ or $2$ because even for higher index, it gives us a bit more information than Theorem \ref{main1}.

\begin{theorem}[metrics of index $1$ or $2$ in dimension $\geq5$, everywhere positive $s$] \label{main3a}
Let $n\geq5$, let $M$ be a compact $n$-manifold, let $q\in\set{1,\dots,n-3}$, let $s\in C^\infty(M,\R)$ be everywhere positive. Then every connected component of $\Metr_q(M)$ contains a metric with scalar curvature $s$.
\smallskip\\
Moreover, let $V$ be a $q$-plane distribution on $M$, let $H$ be an $(n-q)$-plane distribution which is complementary to $V$, and let $\mathscr{H}\subseteq\Distr_{n-q}(M)$ be a $C^0$-neighbourhood of $H$. Then there is a pseudo-Riemannian metric of index $q$ on $M$ with scalar curvature $s$ which makes $V$ timelike and makes some element of $\mathscr{H}$ spacelike.
\end{theorem}
\Proof
As before, it suffices to prove the second statement. The set of all distributions which are complementary to $V$ is a $C^0$-neighbourhood of $H$ in $\Distr_{n-q}(M)$. By Theorem \ref{FIVEMAIN}, its intersection with $\mathscr{H}$ contains an everywhere twisted distribution $H'$. We choose a Riemannian metric $g$ which makes $V$ and $H'$ orthogonal. By Lemma \ref{supsol2}, there is a constant $s_0\in\R$ such that for every function $\tilde{s}\in C^\infty(M,\R)$ with $\tilde{s}\geq s_0$, the constant $1$ is a supersolution of the equation $\PD_{g,V,\tilde{s}}(f)=0$. On the other hand, every sufficiently small positive constant is a subsolution, by Lemma \ref{subsol1}; so we find a subsolution which is smaller than $1$. The method of sub- and supersolutions (cf.\ \ref{choquetbruhatleray}, \ref{susuexample}) shows that there exists an $f\in C^\infty(M,\R_{>0})$ with $\PD_{g,V,\tilde{s}}(f)=0$.
\smallskip\\
Hence Theorem \ref{THEEQUATION} proves that every function $\tilde{s}\in C^\infty(M,\R)$ with $\tilde{s}\geq s_0$ is the scalar curvature of some index-$q$ metric which makes $V$ timelike and $H'$ spacelike. Since our given function $s$ is everywhere positive, there is a (large) constant $c\in\R_{>0}$ with $cs\geq s_0$, and thus there exists an index-$q$ metric $h$ with scalar curvature $cs$ which makes $V$ timelike and $H'$ spacelike. The metric $ch$ has scalar curvature $s$ and makes $V$ timelike and $H'$ spacelike.
\end{proof}

\begin{theorem}[Lorentzian metrics in dimension $4$, everywhere positive $s$] \label{main3b}
Let $M$ be a compact connected orientable $4$-manifold which either has nonempty boundary, or is closed with $\sigma_M\not\equiv 2\mod4$.\footnote{Cf.\ Notation \ref{intersectionform}.} Let $s\in C^\infty(M,\R)$ be everywhere positive. Then every connected component of $\Metr_1(M)$ which consists of time-orientable metrics contains a metric with scalar curvature $s$.
\smallskip\\
Moreover, let $V$ be an orientable line distribution on $M$, let $H$ be a $3$-plane distribution which is complementary to $V$, and let $\mathscr{H}\subseteq\Distr_3(M)$ be a $C^0$-neighbourhood of $H$. Then there is a Lorentzian metric on $M$ with scalar curvature $s$ which makes $V$ timelike and makes some element of $\mathscr{H}$ spacelike.
\end{theorem}
\Proof
The proof is the same as that of Theorem \ref{main3a}, except that we apply Theorem \ref{FIVEMAINfourdim} instead of \ref{FIVEMAIN} to get an everywhere twisted $H'$.
\end{proof}

\begin{theorem}[Lorentzian metrics in dimension $3$, everywhere positive $s$] \label{main3c}
Let $M$ be a compact orientable $3$-manifold, let $s\in C^\infty(M,\R)$ be everywhere positive. Then every connected component of $\Metr_1(M)$ contains a metric with scalar curvature $s$.
\end{theorem}
\Proof
Let $\mathscr{C}$ be a connected component of $\Metr_1(M)$. The connected component $\SDC(\mathscr{C})$ of $\Distr_2(M)$ contains a contact structure $H'$, by Theorem \ref{threecontact}. We choose a Riemannian metric $g$ on $M$ and denote the $g$-orthogonal distribution of $H'$ by $V$. As in the proof of \ref{main3a}, we find now a Lorentzian metric $h$ with scalar curvature $s$ which makes $V$ timelike and $H'$ spacelike. By construction, $h$ is contained in $\mathscr{C}$.
\end{proof}

\emph{Remark.} In many cases, one can improve this $3$-dimensional result in such a way that one obtains a $C^0$-closeness statement as in the preceding theorems. Cf.\ Theorem 2.4.1 in \cite{EliashbergThurston} for the necessary $C^0$-approximation statement.

\begin{theorem}[metrics of index $2$ in dimension $4$, everywhere positive $s$] \label{main3d}
Let $M$ be a compact $4$-manifold, let $s\in C^\infty(M,\R)$ be everywhere positive, let $V$ be a $2$-plane distribution on $M$ which admits a nowhere vanishing section, let $H$ be an orientable $2$-plane distribution which is complementary to $V$, and let $\mathscr{H}\subseteq\Distr_2(M)$ be a $C^0$-neighbourhood of $H$. Then there is a pseudo-Riemannian metric of index $2$ on $M$ with scalar curvature $s$ which makes $V$ timelike and makes some element of $\mathscr{H}$ spacelike.
\end{theorem}
\Proof
The proof is the same as that of Theorem \ref{main3a}, except that we apply Theorem \ref{FIVEMAINtwoplane} instead of \ref{FIVEMAIN} to get an everywhere twisted $H'$.
\end{proof}

\begin{remark}
Recall that when a semi-Riemannian metric $h$ of index $q$ on an $n$-manifold $M$ has scalar curvature $s$ and makes a $q$-plane distribution $V$ timelike and a complementary distribution $H$ spacelike, then the metric $-h$ of index $n-q$ has scalar curvature $-s$, makes $V$ spacelike and $H$ timelike. We get therefore similar theorems as those above for metrics of index $n-1$ or $n-2$ and prescribed functions $s$ which are everywhere negative.
\end{remark}

\begin{remark}[real-analytic metrics] \label{ranaremark}
By the same arguments as in the proof of Theorem \ref{main1ra}, all the theorems above have real-analytic counterparts; i.e., if the manifold $M$ is equipped with a real-analytic atlas, if the function $s$ is real-analytic, and if the distribution $V$ in the theorems \ref{main3a}--\ref{main3d} is real-analytic, then we can find a real-analytic metric with the desired properties.
\end{remark}

\begin{remark}[boundary values] \label{boundaryvalue}
Consider the proofs of the theorems above. If the manifold $M$ has a boundary and a function $\varphi\in C^\infty(\mfbd M,\R_{>0})$ is given, then we can find a solution $f\in C^\infty(M,\R_{>0})$ of the equation $\PD_{g,V,s}(f)=0$ with boundary values $\varphi$: We choose our sub- and supersolutions so small resp.\ large that $\varphi$ lies between them (note that we can replace Lemma \ref{supsol2} by a statement where an arbitrary function instead of the constant $1$ is a supersolution) and then use Theorem \ref{choquetbruhatleray}.
\smallskip\\
A typical application would look like this: \emph{Let $n\geq4$, let $M$ be a compact connected $n$-manifold with nonempty boundary such that $M$ admits a line distribution $V$ which is transverse to the boundary; if $n=4$, assume that $M$ and $V$ are orientable. Let $s\in C^\infty(M,\R)$ be everywhere positive. Then we can find a Lorentzian metric on $M$ with scalar curvature $s$ which makes $V$ timelike and makes $\mfbd M$ spacelike (i.e., each vector in $T(\mfbd M)$ becomes spacelike).}
\smallskip\\
It takes only a few additional observations to verify this statement: In \ref{main3a} and \ref{main3b}, we can choose $H$ in such a way that its restriction to $\mfbd M$ is $T(\mfbd M)$. In the proofs, we choose the Riemannian metric $g$ (which makes $V$ orthogonal to some everywhere twisted distribution $H'$) in such a way that $H$ is spacelike with respect to the Lorentzian metric $\switch(g,V)$. (We can start from an arbitrary Riemannian metric $\tilde{g}$ which makes $V$ and $H'$ orthogonal and take $g = \stre(\tilde{g},\tilde{c},V)$ for a sufficiently large constant $\tilde{c}$. Then the lightcones of $\switch(g,V)$ lie very narrow around $V$, so $H$ is spacelike.)
\smallskip\\
Now we proceed as in the proofs of \ref{main3a} and \ref{main3b}, but choose a solution $f\in C^\infty(M,\R_{>0})$ of the equation $\PD_{g,V,s}(f)=0$ with $f\restrict\mfbd M \equiv 1$. The restriction of our Lorentzian solution metric $ch = c\,\change(g,f,K\compose f,V)$ to the boundary $\mfbd M$ is then conformal to $\stre(\switch(g,V),1,V) = \switch(g,V)$. Hence it makes the boundary spacelike, as claimed.
\end{remark}

\begin{remark}
At first sight, it looks a bit strange that the theorems above are true if $s$ is everywhere \emph{positive}: We know that everywhere \emph{negative} functions are easier to realise as scalar curvatures of \emph{Riemannian} metrics --- in fact, there are obstructions to the existence of Riemannian metrics with positive scalar curvature (cf.\ Appendix \ref{D1}) ---, and the proofs in the Riemannian as well as e.g.\ in the Lorentzian case employ the same general elliptic equation from Theorem \ref{THEEQUATION} (cf.\ also \ref{equationRiemann}), which contains no factor $(-1)^q$ or something like that.
\smallskip\\
But there is nothing wrong here. First, the easy sign in the Riemannian case (i.e.\ negative scalar curvature) is probably also ``easy'' in the Lorentzian case, in the sense that there are no obstructions to the existence of Lorentzian metrics with negative scalar curvature on manifolds of dimension $\geq4$ which admit a Lorentzian metric; cf.\ Section \ref{SIXesc}. It is just that our method of proof from above does not produce such metrics.
\smallskip\\
Second, the reason why the easier sign in the Lorentzian case (i.e.\ positive scalar curvature) is not easy in the Riemannian case is just that one cannot construct twisted $n$-plane distributions on an $n$-manifold. (Thus one has no subsolution, whereas an everywhere positive $s$ yields a supersolution in the same way as in the proof of Lemma \ref{supsol2}.) This makes the two problems very different, despite the similarity of the involved elliptic equations.
\end{remark}


\section{The Kazdan/Warner method} \label{SIXKW}

Our aim here is to find solutions of the equation $\PD_{g,V,s}(f)=0$ via the implicit function theorem for Banach spaces (cf.\ \ref{implicitfunction}) and the Kazdan/Warner approximation theorem \ref{KWapproximation1}. In order to do this, we have to interpret $\PD_{g,V,s}$ as an operator from a Sobolev space of functions to an $L^p$ space, and we have to find a function $f\colon M\to\R_{>0}$ such that the derivative of $\PD_{g,V,s}$ at the point $f$ is invertible.
\smallskip\\
This technique has been applied very successfully by Kazdan and Warner to the Riemannian prescribed curvature problem. In their solution, they proved also that the derivative of the relevant elliptic operator is invertible in every point $f$ from an open and \emph{dense} set of functions $\in C^2(M,\R_{>0})$; cf.\ \cite{KazdanWarner1}, \S4. We will not try to establish a similar perturbation theorem because our operator is considerably more complicated than the one that Kazdan and Warner discussed; in particular, ours is not formally self-adjoint.
\smallskip\\
According to the philosophy mentioned in the introduction to this chapter, we use only the simplest criterion for the invertibility of a linear elliptic operator, namely Theorem \ref{invertible}.

\subsection{The derivative}

\begin{remarkdefinition} \label{Sdef}
Let $n\in\N_{\geq2}$, let $(M,g)$ be a Riemannian $n$-manifold, let $V$ be a $q$-plane distribution on $M$. For each $f\in C^\infty(M,\R_{>0})$, there is a unique function $s\in C^\infty(M,\R)$ such that $\PD_{g,V,s}(f)=0$. We denote this function by $S(f)$. It is obviously given by
\[ \begin{split}
S(f) &= f^{-\mu}(1+f^2)^{-\nu}\bigg( 2\laplace_g(f) +a_{n,q}(f)\abs{df}^2_g +b_{n,q}(f)\abs{df}^2_{g,V}\\
&\mspace{150mu}+\frac{2(1+f^2)}{f^2}\eval{\divergence^V_g}{df}_{g,H} +2(1+f^2)\eval{\divergence^H_g}{df}_{g,V}\\
&\mspace{150mu}+\frac{(1+f^2)^2}{2f^3}\abs{\Twist_H}^2_g -\frac{f(1+f^2)^2}{2}\abs{\Twist_V}^2_g +\frac{(1+f^2)^2}{f}\xi_{g,V} +\frac{1+f^2}{f}\scal_g \bigg) \;\;,
\end{split} \]
where we used the abbreviations \fbox{$\mu\define \frac{2q}{n-1}-1$} and \fbox{$\nu\define \frac{n-2}{n-1}$}, and where $a_{n,q},b_{n,q}\in C^\infty(\R_{>0},\R)$ are the functions from Definition \ref{PDdef}.
\end{remarkdefinition}

Recall the basic facts about Sobolev spaces from Appendix \ref{sobolevappendix}, in particular the inclusion $\Sob{p}{2}(M,\R_{>0})\subseteq C^1(M,\R_{>0})$ for $p>\dim(M)$ (cf.\ \ref{sobolev}, \ref{compdef}). Cf.\ also Remark \ref{implicitremark} for an explanation of what is happening now. In the following, we identify positive real numbers with positive constant functions on $M$.

\begin{definition} \label{Phidefinition}
Let $n\in\N_{\geq2}$, let $(M,g)$ be a compact Riemannian $n$-manifold, let $V$ be a $q$-plane distribution on $M$, let $H$ denote the $g$-orthogonal distribution of $V$, let $p\in\R$ with $p>n$, let $c\in\R_{>0}$. We denote the affine subspace $c+\Sobzero{p}{2}(M,\R)$ of $\Sob{p}{2}(M,\R)$ by $A$, and we consider the open subset $\mathscr{N}^c\define A\cap\Sob{p}{2}(M,\R_{>0})$ of $A$. We define a map $\Phi^c\colon \mathscr{N}^c\times L^p(M,\R)\to L^p(M,\R)$ by
\[ \begin{split}
\Phi^c(f,s) &\define 2\laplace_g(f) +a_{n,q}(f)\abs{df}^2_g +b_{n,q}(f)\abs{df}^2_{g,V} +\frac{2(1+f^2)}{f^2}\eval{\divergence^V_g}{df}_{g,H} +2(1+f^2)\eval{\divergence^H_g}{df}_{g,V}\\
&\mspace{20mu}+\frac{(1+f^2)^2}{2f^3}\abs{\Twist_H}^2_g -\frac{f(1+f^2)^2}{2}\abs{\Twist_V}^2_g +\frac{(1+f^2)^2}{f}\xi_{g,V} +\frac{1+f^2}{f}\scal_g -f^\mu(1+f^2)^\nu s \;\;.
\end{split} \]
In other words, we extend the definition of the differential operator $\PD_{g,V,s}\colon C^\infty(M,\R_{>0})\to C^\infty(M,\R)$ to $\mathscr{N}^c\subseteq \Sob{p}{2}(M,\R_{>0})$, and we make the dependence on $s$ explicit, at the same time allowing arbitrary functions $s\in L^p(M,\R)$.
\smallskip\\
The map $\Phi^c$ is well-defined since $\laplace_g$ is well-defined as a linear map $\Sob{p}{2}(M,\R)\to L^p(M,\R)$ and since all the lower-order terms are well-defined via the inclusions $\Sob{p}{2}(M,\R)\subseteq C^1(M,\R)$ and $C^0(M,\R)\subseteq L^p(M,\R)$.
\smallskip\\
For every $s\in L^p(M,\R)$, we define the map $\Phi^c_s\colon \mathscr{N}^c\to L^p(M,\R)$ by $\Phi^c_s(f)\define \Phi^c(f,s)$.
\end{definition}

\begin{lemma} \label{Phideriv}
Let $n\in\N_{\geq2}$, let $(M,g)$ be a compact Riemannian $n$-manifold, let $V$ be a $q$-plane distribution on $M$, let $H$ denote the $g$-orthogonal distribution of $V$, let $p\in\R$ with $p>n$, let $c\in\R_{>0}$. Then $\Phi^c\colon \mathscr{N}^c\times L^p(M,\R)\to L^p(M,\R)$ is continuous. $\Phi^c_s\colon \mathscr{N}^c\to L^p(M,\R)$ is (\Frechet) differentiable for every $s\in L^p(M,\R)$, and its derivative $D_f\Phi^c_s\colon \Sobzero{p}{2}(M,\R) \to L^p(M,\R)$ in the point $f\in\mathscr{N}^c$ is given by
\[ \begin{split}
(D_f\Phi^c_s)(v) &= 2\laplace_g(v) +2a_{n,q}(f)\eval{df}{dv}_g +a_{n,q}'(f)\abs{df}_g^2v +2b_{n,q}(f)\eval{df}{dv}_{g,V} +b_{n,q}'(f)\abs{df}^2_{g,V}v\\
&\mspace{20mu}+\frac{2(1+f^2)}{f^2}\eval{\divergence^V_g}{dv}_{g,H} -\frac{4}{f^3}\eval{\divergence^V_g}{df}_{g,H}\,v +2(1+f^2)\eval{\divergence^H_g}{dv}_{g,V} +4f\eval{\divergence^H_g}{df}_{g,V}\,v\\
&\mspace{20mu}+\frac{(1+f^2)(f^2-3)}{2f^4}\abs{\Twist_H}^2_gv -\frac{(1+f^2)(1+5f^2)}{2}\abs{\Twist_V}^2_gv +\frac{(1+f^2)(3f^2-1)}{f^2}\xi_{g,V}\,v\\
&\mspace{20mu}-\frac{1-f^2}{f^2}\scal_gv -f^\mu(1+f^2)^\nu\frac{\mu+(\mu+2\nu)f^2}{f(1+f^2)}sv \;\;.
\end{split} \]
The map $\mathscr{N}^c\times L^p(M,\R)\to \Lin(\Sobzero{p}{2}(M,\R), L^p(M,\R))$ given by $(f,s)\mapsto D_f\Phi^c_s$ is continuous.
\end{lemma}
\Proof
Note that for every distribution $U$ on $M$, the map $\psi\colon \Sob{p}{2}(M,\R)\to L^p(M,\R)$ given by $u\mapsto \abs{du}_{g,U}^2$ is differentiable, and that its derivative is given by $(D_u\psi)(v) = 2\eval{du}{dv}_{g,U}$; this follows immediately from the existence of a constant $c_0\in\R_{>0}$ with $\norm{.}_{C^0}\leq c_0\norm{.}_{\Sob{p}{1}}$ (cf.\ \ref{sobolev}):
\begin{multline*}
\norm{\abs{d(u+v)}_{g,U}^2 -\abs{du}_{g,U}^2 -2\eval{du}{dv}_{g,U}}_{L^p} = \norm{\abs{dv}_{g,U}^2}_{L^p}\\
\leq \volume(M)^{1/p}\norm{\abs{dv}_{g,U}}_{C^0}^2
\leq c_0^2\volume(M)^{1/p}\norm{\abs{dv}_{g,U}}_{\Sob{p}{1}}^2
\leq c_0^2\volume(M)^{1/p}\norm{v}_{\Sob{p}{2}}^2 \;\;.
\end{multline*}
In particular, the map $\psi$ is continuous. Now the fact that every $k$th-order linear differential operator $C^\infty(M,\R)\to C^\infty(M,\R)$ induces a continuous linear map $\Sob{p}{k}(M,\R)\to L^p(M,\R)$, together with the continuity statements from \ref{multcont} and \ref{compdef}, implies that $\Phi^c$ is continuous.
\smallskip\\
Since a continuous linear map is in each point its own derivative and the product rule \ref{productrule} and the ``chain rule'' \ref{chainrule} hold, we see that $\Phi^c_s$ is differentiable at each $f\in \mathscr{N}^c$, and we see (after computing a few derivatives of functions $\R_{>0}\to\R$) that $D_f\Phi^c_s$ is given by the claimed formula.
\smallskip\\
By Lemma \ref{conti1} and Lemma \ref{conti2}, this formula implies that the map $(f,s)\mapsto D_f\Phi^c_s$ is continuous.
\end{proof}

For the following proposition, note that if the boundary of $M$ is nonempty, then the domain $\mathscr{N}^c$ of $\Phi^c$ contains precisely one constant function, namely $c$.

\begin{proposition} \label{invertibleprop}
Let $n\in\N_{\geq2}$, let $(M,g)$ be a compact Riemannian $n$-manifold, let $V$ be a $q$-plane distribution on $M$, let $H$ denote the $g$-orthogonal distribution of $V$, let $p\in\R$ with $p>n$.
\smallskip\\
If $V$ is everywhere twisted, then there is a number $c_+\in\R_{>0}$ such that for every constant $c\geq c_+$, the function $S(c)\in C^\infty(M,\R)$ is everywhere negative and the linear map $D_c\Phi^c_{S(c)}\colon \Sobzero{p}{2}(M,\R) \to L^p(M,\R)$ is bijective.
\smallskip\\
If $H$ is everywhere twisted, then there is a number $c_-\in\R_{>0}$ such that for every constant $c$ with $0< c\leq c_-$, the function $S(c)\in C^\infty(M,\R)$ is everywhere positive and the linear map $D_c\Phi^c_{S(c)}\colon \Sobzero{p}{2}(M,\R) \to L^p(M,\R)$ is bijective.
\end{proposition}
\Proof
For any constant $c>0$, the zeroth-order coefficient of the linear differential operator $D_c\Phi^c_{S(c)}$ is
\renewcommand{\FILL}{\mspace{20mu}}
\[ \begin{split}
&(D_c\Phi^c_{S(c)})(1)\\[1ex]
&\FILL= \frac{(1+c^2)(c^2-3)}{2c^4}\abs{\Twist_H}^2_g -\frac{(1+c^2)(1+5c^2)}{2}\abs{\Twist_V}^2_g +\frac{(1+c^2)(3c^2-1)}{c^2}\xi_{g,V} -\frac{1-c^2}{c^2}\scal_g\\
&\FILL\mspace{20mu}-c^\mu(1+c^2)^\nu\frac{\mu+(\mu+2\nu)c^2}{c(1+c^2)}S(c)\\
&\FILL= \frac{(1+c^2)(c^2-3)}{2c^4}\abs{\Twist_H}^2_g -\frac{(1+c^2)(1+5c^2)}{2}\abs{\Twist_V}^2_g +\frac{(1+c^2)(3c^2-1)}{c^2}\xi_{g,V} -\frac{1-c^2}{c^2}\scal_g\\
&\FILL\mspace{20mu}-\frac{\mu+(\mu+2\nu)c^2}{c(1+c^2)} \bigg(\frac{(1+c^2)^2}{2c^3}\abs{\Twist_H}^2_g -\frac{c(1+c^2)^2}{2}\abs{\Twist_V}^2_g +\frac{(1+c^2)^2}{c}\xi_{g,V} +\frac{1+c^2}{c}\scal_g \bigg)\\
&\FILL= \frac{1+c^2}{2c^4}\Big(c^2-3-\mu-(\mu+2\nu)c^2\Big)\abs{\Twist_H}^2_g -\frac{1+c^2}{2}\Big(1+5c^2-\mu-(\mu+2\nu)c^2\Big)\abs{\Twist_V}^2_g\\
&\FILL\mspace{20mu}+\frac{1+c^2}{c^2}\Big(3c^2-1-\mu-(\mu+2\nu)c^2\Big)\xi_{g,V} +\frac{1}{c^2}\Big(-1+c^2-\mu-(\mu+2\nu)c^2\Big)\scal_g \;\;;
\end{split} \]
cf.\ \ref{Sdef} and \ref{Phideriv}. Note that $0\leq\nu<1$ and $-1\leq\mu\leq \frac{2n}{n-1}-1=1+\frac{2}{n-1}\leq3$; hence $-3-\mu < 0$ and $-5+\mu+2\nu < 0$.
\smallskip\\
Assume that $V$ is everywhere twisted. Since $M$ is compact, there is then an $\eps\in\R_{>0}$ such that $\abs{\Twist_V}^2_g\geq\eps$, and there is a $C\in\R_{>0}$ such that each of the functions $\abs{\Twist_H}^2_g$, $\abs{\xi_{g,V}}$, $\abs{\scal_g}$ is $\leq C$. Since $-5+\mu+2\nu \neq 0$, the absolute value of the coefficient of $\abs{\Twist_V}^2_g$ in the expression above increases like $c^4$ as $c$ tends to $\infty$, whereas the absolute values of the coefficients of the other three functions increase at most like $c^2$ as $c\to\infty$. Thus $-5+\mu+2\nu < 0$ implies that there is a constant $\tilde{c}_+>0$ such that $(D_c\Phi^c_{S(c)})(1) < 0$ for all $c\geq\tilde{c}_+$.
\smallskip\\
A very similar argument (look at the formula in \ref{Sdef}) shows that there is a constant $\bar{c}_+>0$ such that $S(c)$ is everywhere negative for all $c\geq\bar{c}_+$. We define $c_+$ to be $\max\set{\tilde{c}_+,\bar{c}_+}$. For every $c\geq c_+$, the function $S(c)$ is everywhere negative; moreover, since the zeroth-order coefficient of the positively elliptic operator $D_c\Phi^c_{S(c)}\colon \Sobzero{p}{2}(M,\R)\to L^p(M,\R)$ is everywhere negative, $D_c\Phi^c_{S(c)}$ is in fact bijective; cf.\ Theorem \ref{invertible}. This completes the proof in the case when $V$ is everywhere twisted.
\smallskip\\
Now assume that $H$ is everywhere twisted. Since $M$ is compact, there is an $\eps\in\R_{>0}$ such that $\abs{\Twist_H}^2_g\geq\eps$, and there is a $C\in\R_{>0}$ such that each of the functions $\abs{\Twist_V}^2_g$, $\abs{\xi_{g,V}}$, $\abs{\scal_g}$ is $\leq C$. Since $-3-\mu \neq 0$, the absolute value of the coefficient of $\abs{\Twist_H}^2_g$ in the formula for $(D_c\Phi^c_{S(c)})(1)$ increases like $c^{-4}$ as $c$ tends to $0$, whereas the absolute values of the coefficients of the other three functions increase at most like $c^{-2}$ as $c\to0$. Thus $-3-\mu < 0$ implies that there is a constant $c_->0$ such that $(D_c\Phi^c_{S(c)})(1) < 0$ for all $c$ with $0<c\leq c_-$.
\smallskip\\
Again a similar argument involving the formula in \ref{Sdef} shows that there is a constant $\bar{c}_->0$ such that $S(c)$ is everywhere positive for all $c$ with $0<c\leq\bar{c}_-$. We define $c_-\define \min\set{\tilde{c}_-,\bar{c}_-}>0$ and get the desired statement mutatis mutandis as before.
\end{proof}

\subsection{The main theorems}

Taking into account what we have already obtained via the sub- and supersolution method, it is now easy to prove the following theorems.

\begin{theorem}[metrics of index $1$ or $2$ in dimension $\geq5$, somewhere positive $s$] \label{main4a}
Let $n\geq5$, let $M$ be a compact connected $n$-manifold, let $q\in\set{1,2}$, and let $s\in C^\infty(M,\R)$ be somewhere positive. Then every connected component of $\Metr_q(M)$ contains a metric with scalar curvature $s$.
\end{theorem}
\Proof
If $s$ is everywhere positive, then the statement of the theorem follows from \ref{main3a}. So we assume that $s$ is somewhere positive and somewhere $0$. Let $\mathscr{C}$ be a connected component of $\Metr_q(M)$. By Theorem \ref{FIVEMAIN}, some distribution $H$ in the connected component $\SDC(\mathscr{C})$ of $\Distr_{n-q}(M)$ is everywhere twisted (because $3\leq n-q\leq n-1$). We choose any Riemannian metric $g$ on $M$ and denote the $g$-orthogonal distribution of $H$ by $V$.
\smallskip\\
With respect to these data $V,H,g$, we consider the maps $S\colon C^\infty(M,\R_{>0})\to C^\infty(M,\R)$ from Definition \ref{Sdef} and $\Phi^c\colon \mathscr{N}^c\times L^p(M,\R)\to L^p(M,\R)$ from Definition \ref{Phidefinition}. By Proposition \ref{invertibleprop}, there exists a number $c_-\in\R_{>0}$ such that for every constant $c$ with $0< c\leq c_-$, the function $S(c)\in C^\infty(M,\R)$ is everywhere positive and the linear map $D_c\Phi^c_{S(c)}\colon \Sobzero{p}{2}(M,\R) \to L^p(M,\R)$ is bijective. We choose a constant $c$ with $0<c\leq c_-$. By Lemma \ref{Phideriv}, $\Phi^c_{\tilde{s}}$ is differentiable for every $\tilde{s}\in L^p(M,\R)$, and the map $\Phi^c$ and the map $(f,\tilde{s})\mapsto D_f\Phi^c_{\tilde{s}}$ are continuous.
\smallskip\\
Note that $\Phi^c(c,S(c))=0$ by the definitions of $S$ and $\Phi^c$. Hence the implicit function theorem \ref{implicitfunction} tells us that there exist an open neighbourhood $\mathscr{U}\subseteq L^p(M,\R)$ of $S(c)$ and a continuous function $U\colon \mathscr{U}\to \mathscr{N}^c$ such that $U(S(c))=c$ and $\Phi^c(U(\tilde{s}),\tilde{s})=0$ for all $\tilde{s}\in\mathscr{U}$.
\smallskip\\
Since our function $s\in C^\infty(M,\R)$ is somewhere positive and somewhere $0$, and since the function $S(c)\in C^\infty(M,\R)$ is everywhere positive, there is a constant $r>0$ such that $\inf(rs)\leq S(c)\leq \sup(rs)$. The Kazdan/Warner approximation theorem \ref{KWapproximation1} implies that there exists a diffeomorphism $\varphi\in\Diff^0(M)$ such that $rs\compose\varphi\in \mathscr{U}$. Thus there is a solution $f\in\Sob{p}{2}(M,\R_{>0})$ of the elliptic equation $\Phi^c_{rs\compose\varphi}(f)=0$. This equation has smooth coefficients, so $f$ is actually smooth, by elliptic regularity (cf.\ \ref{sobolevregularity}).
\smallskip\\
In other words, $f\in C^\infty(M,\R_{>0})$ solves the equation $\PD_{g,V,rs\compose\varphi}(f)=0$. Hence Theorem \ref{THEEQUATION} yields a pseudo-Riemannian metric $h$ on $M$ with scalar curvature $rs\compose\varphi$ which makes $V$ timelike and is thus contained in $\mathscr{C}$. Since $\varphi\in\Diff^0(M)$, the metric $(\varphi^{-1})^\ast(rh)$ lies in the same connected component of $\Metr_q(M)$ as $h$. It has scalar curvature $s$.
\end{proof}

\begin{theorem}[Lorentzian metrics in dimension $4$, somewhere positive $s$] \label{main4b}
Let $M$ be a compact connected orientable $4$-manifold which either has nonempty boundary, or is closed with $\sigma_M\not\equiv 2\mod4$. Let $s\in C^\infty(M,\R)$ be somewhere positive. Then every connected component of $\Metr_1(M)$ which consists of time-orientable metrics contains a metric with scalar curvature $s$.
\end{theorem}
\Proof
The proof is the same as that of Theorem \ref{main4a}, except that we apply Theorem \ref{FIVEMAINfourdim} instead of \ref{FIVEMAIN} in order to get an everywhere twisted $H$, and that we invoke \ref{main3b} instead of \ref{main3a} for everywhere positive $s$.
\end{proof}

\begin{theorem}[Lorentzian metrics in dimension $3$, somewhere positive $s$] \label{main4c}
Let $M$ be a compact connected orientable $3$-manifold, let $s\in C^\infty(M,\R)$ be somewhere positive. Then every connected component of $\Metr_1(M)$ contains a metric with scalar curvature $s$.
\end{theorem}
\Proof
The proof is the same as that of Theorem \ref{main4a}, except that we apply Theorem \ref{threecontact} instead of \ref{FIVEMAIN} in order to get an everywhere twisted $H$, and that we invoke \ref{main3c} instead of \ref{main3a} for everywhere positive $s$.
\end{proof}

\begin{theorem}[metrics of index $2$ in dimension $4$, somewhere positive $s$] \label{main4d}
Let $M$ be a compact connected $4$-manifold, let $s\in C^\infty(M,\R)$ be somewhere positive, let $\mathscr{C}$ be a connected component of $\Metr_2(M)$ consisting of space-orientable metrics, such that the elements of $\TDC(\mathscr{C})$ admit a nowhere vanishing section. Then $\mathscr{C}$ contains a metric with scalar curvature $s$.
\end{theorem}
\Proof
The proof is the same as that of Theorem \ref{main4a}, except that we apply Theorem \ref{FIVEMAINtwoplane} instead of \ref{FIVEMAIN} in order to find an everywhere twisted $H$ inside $\mathscr{C}$, and that we invoke \ref{main3d} instead of \ref{main3a} for everywhere positive $s$.
\end{proof}

\subsection{Product manifolds} \label{SIXwhatever}

Until now, all our solutions of the prescribed scalar curvature problem have been constructed via everywhere twisted distributions. But what about, say, the diffeotopy class problem mentioned in the introduction chapter? If we start with \emph{integrable} distributions $V$ and $H$, can we find solutions which, for some diffeomorphism $\varphi\in\Diff^0(M)$, make $\varphi^\ast V$ timelike and $\varphi^\ast H$ spacelike? Such solutions would admit timelike and spacelike foliations. For none of our solutions so far is it clear whether it admits both a timelike and a spacelike foliation.
\smallskip\\
Second, one might have got the misleading impression that positive functions are easier to realise as Lorentzian scalar curvatures than negative functions. But the real picture should be a bit different: In the Lorentzian problem, positive functions are easier to realise as scalar curvatures than negative functions as long as we restrict ourselves to choosing the background \emph{distributions} in a clever way. But \emph{negative} functions should be easier to realise when we restrict ourselves to choosing the background \emph{metric} nicely (while the background distributions might be integrable).
\smallskip\\
We will now investigate these issues in the simplest example, namely on a product manifold $M = B\times N$. Our background distributions in the proof of the following theorem will be the first-factor and second-factor distributions on $B\times N$. The proof employs the Kazdan/Warner method but is different than our proofs above: It arranges the invertibility of the relevant operator not by adjusting the sign of the zeroth-order coefficient but by a perturbation argument like in the Kazdan/Warner article \cite{KazdanWarner1}. However, our perturbation argument is comparatively trivial: it just uses the fact that the spectrum of the Laplacian is discrete.
\smallskip\\
For simplicity, we assume (here without essential loss of generality) that $B$ and $M$ have no boundary. Note that the following theorem shows in particular that in the Lorentzian case $q=1$ on product manifolds of dimension $\geq4$, every somewhere \emph{negative} function is a scalar curvature.

\begin{theorem} \label{main5}
Let $q,m\in\N_{\geq1}$, let $B$ be a closed connected $q$-manifold, let $N$ be a closed connected $m$-manifold, let $M$ be the product manifold $B\times N$. Let $V,H$ denote the first-factor resp.\ second-factor distribution on $M$, and let $s\in C^\infty(M,\R)$ be a somewhere negative function. If $m\geq3$, or if $m=2$ and $\chi(N)<0$, or if $q=m=2$ and $\chi(B)>0$, or if $q\geq3$ and $B$ admits a Riemannian metric with positive scalar curvature, then there exists a pseudo-Riemannian metric $h$ of index $q$ on $M$ with scalar curvature $s$, and there exists a diffeomorphism $\varphi\in\Diff^0(M)$ such that $\varphi^\ast(V)$ is $h$-timelike and $\varphi^\ast(H)$ is $h$-spacelike.
\end{theorem}
\Proof
If $m\geq3$, or if $m=2$ and $\chi(N)<0$, or if $q=m=2$ and $\chi(B)>0$, or if $q\geq3$ and $B$ admits a Riemannian metric with positive scalar curvature, then there exist Riemannian metrics $g_B$ on $B$ and $g_N$ on $N$ with \emph{constant} scalar curvatures $k_B$ resp.\ $k_N$ such that $k_N<0$ or $k_B>0$; this follows from the results on the Riemannian prescribed scalar curvature problem (cf.\ Appendix \ref{D1}). In the case $q>1$, we can arrange in addition that $k_B\neq0$.
\smallskip\\
We consider the product metric $g\define g_B\oplus g_N$ on $M$; it makes $V$ and $H$ orthogonal. Let $n\define q+m$.
\smallskip
In this situation, the functions $\eval{\divergence^V_g}{df}_{g,H}$, $\eval{\divergence^H_g}{df}_{g,V}$, $\abs{\Twist_V}_g^2$, $\abs{\Twist_H}_g^2$ vanish for every function $f\in C^\infty(M,\R)$, and we have $\scal_g = k_B+k_N$ and $\xi_{g,V} = -k_B$.
\smallskip\\
Let us prove this: $\Twist_V$ and $\Twist_H$ vanish because $V$ and $H$ are integrable, and the scalar curvature of a product metric is the sum of the scalar curvatures of the factors. Concerning the functions $\eval{\divergence^V_g}{df}_{g,H}$, $\eval{\divergence^H_g}{df}_{g,V}$, $\xi_{g,V}$, consider the formulae \ref{moreONformulae}, \ref{sigmatauframe}, \ref{qON}, \ref{xidef}. For each $(x,y)\in B\times N$, we can choose a $g_B$-ON frame on a neighbourhood of $x$ and a $g_N$-ON frame on a neighbourhood of $y$ in $N$. Clearly these two frames together define a $V$-adapted $g$-ON frame $(e_1,\dots,e_n)$ on a neighbourhood of $(x,y)$ in $M$ with the property that the corresponding ON Christoffel symbols $\Gamma^k_{ij}$ vanish except when $i,j,k:V$ or $i,j,k:H$. (Just observe that $[e_i,e_j]=0$ if $i:V$ and $j:H$, and that $g([e_i,e_j],e_k)=0$ if $i,j:V$ and $k:H$; analogously with the roles of $V$ and $H$ reversed. Thus $g([e_i,e_j],e_k)=0$ except when $i,j,k:V$ or $i,j,k:H$.) Hence we infer $0 = \sigma_{g,V} = \tau_{g,V} = \sigma_{g,H} = \tau_{g,H} = \lvert\divergence^V_g\rvert_{g,H}^2 = \lvert\divergence^H_g\rvert_{g,V}^2 = \eval{\divergence^V_g}{df}_{g,H} = \eval{\divergence^H_g}{df}_{g,V} = \qual^V_g$ from the cited formulae. This shows also $\xi_{g,V} = -\scal^{V,V}_g = -\fscal_{g,V} = -k_B$ (cf.\ \ref{fscallemma}), so our claim above was true.
\smallskip\\
Now we consider the map $S\colon C^\infty(M,\R_{>0})\to C^\infty(M,\R)$ from \ref{Sdef}. For a constant function $c>0$, we have
\[
S(c) = \rho(c)\Big(-(1+c^2)k_B +(k_B+k_N)\Big) = \rho(c)\Big(-c^2k_B +k_N\Big) \;\;,
\]
where $\rho\colon\R_{>0}\to\R_{>0}$ is a positive-valued function. Since $k_B>0$ or $k_N<0$, there exists a nonempty open interval $I\subseteq\R_{>0}$ such that $S(c)<0$ for all $c\in I$.
\smallskip\\
For some $p>n$, consider the operator $\Phi\colon \Sob{p}{2}(M,\R_{>0})\to L^p(M,\R)$ from Definition \ref{Phidefinition} (we omit the $c$ in the notation $\Phi^c$ because our manifold $M$ is closed, so $\mathscr{N}^c = \Sob{p}{2}(M,\R_{>0})$ and thus $\Phi^c$ does not depend on $c$). For a constant function $c>0$, the derivative of $\Phi_{S(c)}$ in the point $c$ is given by (cf.\ \ref{Phideriv}):
\[ \begin{split}
(D_c\Phi_{S(c)})(v) &= 2\laplace_g(v) -\frac{(1+c^2)(3c^2-1)}{c^2}k_Bv -\frac{1-c^2}{c^2}(k_B+k_N)v\\
&\mspace{20mu}-\frac{\mu+(\mu+2\nu)c^2}{c(1+c^2)}\Big(-\frac{(1+c^2)^2}{c}k_B +\frac{1+c^2}{c}(k_B+k_N)\Big)v\\
&= 2\laplace_g(v) +\lambda(c)v \;\;,
\end{split} \]
where
\[ \begin{split}
\lambda(c) = -\frac{(1+c^2)(3c^2-1)}{c^2}k_B -\frac{1-c^2}{c^2}(k_B+k_N) -\frac{\mu+(\mu+2\nu)c^2}{c^2}\Big(-(1+c^2)k_B +k_B+k_N\Big) \;\;.
\end{split} \]
We claim that the function $\lambda\colon\R_{>0}\to\R$ is not constant on the interval $I$. If $q=1$, then $k_B=0$ and thus $\lambda(c) = (-1+c^2-\mu-(\mu+2\nu)c^2)k_N/c^2$; thus $\lambda$ is not constant on $I$ because $k_N\neq0$ and $1+\mu = 2q/(n-1)\neq 0$. If $q>1$, then the leading term of the polynomial $c^2\lambda(c)$ is $(-3+\mu+2\nu)k_Bc^4$ since we chose $k_B\neq0$ and have $-3+\mu+2\nu = -2(n-q)/(n-1)\neq 0$; so again $\lambda$ is not constant on $I$.
\smallskip\\
Since the spectrum of the elliptic operator $2\laplace_g$ is discrete, we can find a constant $c\in I$ such that $D_c\Phi_{S(c)} = 2\laplace_g +\lambda(c)\colon \Sob{p}{2}(M,\R)\to L^p(M,\R)$ is invertible (cf.\ \cite{Taylor3}, \S13.7). Hence the implicit function theorem \ref{implicitfunction} gives us a neighbourhood $\mathscr{U}$ of $S(c)\in L^p(M,\R)$ and a function $U\colon \mathscr{U}\to \Sob{p}{2}(M,\R_{>0})$ such that $\Phi(U(\tilde{s}),\tilde{s})=0$ for all $s\in\mathscr{U}$.
\smallskip\\
Since $S(c)<0$ and our prescribed function $s\in C^\infty(M,\R)$ is somewhere negative, there is an $r\in\R_{>0}$ with $\inf(rs)\leq S(c)\leq\sup(rs)$. The Kazdan/Warner approximation theorem \ref{KWapproximation1} tells us that there exists a diffeomorphism $\varphi\in\Diff^0(M)$ such that $rs\compose\varphi\in\mathscr{U}$. The function $f\define U(rs\compose\varphi)\in \Sob{p}{2}(M,\R_{>0})$ solves the elliptic equation $\Phi_{rs\compose\varphi}(f)=0$ and is thus smooth by elliptic regularity. Since $\PD_{g,V,rs\compose\varphi}(f)=0$, we infer from Theorem \ref{THEEQUATION} that there is a pseudo-Riemannian metric $\tilde{h}$ of index $q$ on $M$ with scalar curvature $rs\compose\varphi$, such that $V$ is $\tilde{h}$-timelike and $H$ is $\tilde{h}$-spacelike.
\smallskip\\
The metric $h\define (\varphi^{-1})^\ast(r\tilde{h})$ has scalar curvature $s$; the distribution $(\varphi^{-1})^\ast(V)$ is $h$-timelike, and $(\varphi^{-1})^\ast(H)$ is $h$-spacelike.
\end{proof}


\section{The esc Conjecture} \label{SIXesc}

Recall the esc Conjecture \ref{esc} from the introduction chapter:

\begin{escconjecture}
Let $M$ be a compact manifold of dimension $\geq4$, and let $s\in C^\infty(M,\R)$. Then every connected component of the space of Lorentzian metrics on $M$ contains a metric with scalar curvature $s$.
\end{escconjecture}

The aim of this section is to give a very brief and rough outline of how one might try to prove this conjecture. We start with another conjecture:

\begin{conjecture} \label{chinegative}
Let $n\geq4$, let $M$ be an $n$-manifold, let $H$ be an integrable $(n-1)$-plane distribution on $M$. Then there is a Riemannian metric $g$ on $M$ such that the function $\chi_{g,\bot H}$ is everywhere negative.
\end{conjecture}

This conjecture implies the esc Conjecture in dimension $n\geq5$, and in many cases also in dimension $4$. We will first sketch why this is so; then we will discuss a possible strategy to prove \ref{chinegative}.
\medskip\\
In the Lorentzian case of the prescribed scalar curvature problem, we have to solve the elliptic equation (cf.\ \ref{lorentzequation})
\begin{equation} \begin{split} \label{LPDE}
0 &= 2\laplace_g(f) -\frac{4f^4+\alpha(n)}{f^3(1+f^2)}\abs{df}^2_g +\frac{\alpha(n)}{f^3}\abs{df}^2_{g,V} +\frac{2(1+f^2)}{f^2}\eval{\divergence^V_g}{df}_{g,H} +2(1+f^2)\eval{\divergence^H_g}{df}_{g,V}\\
&\mspace{20mu}+\frac{(1+f^2)^2}{2f^3}\abs{\Twist_H}^2_g +\frac{(1+f^2)^2}{f}\xi_{g,V} +\frac{1+f^2}{f}\scal_g -f\Big(\frac{1+f^2}{f^2}\Big)^{\alpha(n)}s \;\;.
\end{split} \end{equation}
We try to do this via the sub- and supersolution method, again with \emph{constant} sub- and supersolutions. The main problem is to find a supersolution, in particular when $s$ is everywhere negative.
\smallskip\\
First we invoke a strong theorem due to W.\ Thurston (cf.\ \cite{Thurston1976}) which says that each homotopy class of $(n-1)$-plane distributions on an $n$-manifold $M$ contains an integrable distribution. So we choose $H$ to be integrable. In this way, we get rid of the term $\frac{(1+f^2)^2}{2f^3}\abs{\Twist_H}^2_g$ which has the wrong sign as far as existence of supersolutions of \eqref{LPDE} is concerned.
\smallskip\\
If $f$ is constant, then the right hand side of \eqref{LPDE} is now
\[ \begin{split}
\frac{(1+f^2)^2}{f}\xi_{g,V} +\frac{1+f^2}{f}\scal_g -f\Big(\frac{1+f^2}{f^2}\Big)^{\frac{n-2}{n-1}}s &= \frac{1+f^2}{f}\Big((1+f^2)\xi_{g,V} +\scal_g -f^{\frac{2}{n-1}}(1+f^2)^{-\frac{1}{n-1}}s\Big) \;\;.
\end{split} \]
If we choose the constant $f>0$ very small, then the coefficient of $s$ is nearly zero, so we get rid of the problems with the possibly wrong sign of $s$. (At first sight, it might not look like a good idea to choose a very small supersolution, because we also need an even smaller subsolution. We will see in a moment that the idea is not so bad after all.)
\smallskip\\
Note that $\chi_{g,V} = \scal_g +\xi_{g,V}$. If Conjecture \ref{chinegative} is true, then we find, for a suitable metric $g$, a constant supersolution $f_+$ of our equation (in the strict sense that $\PD_{g,\bot_gH,s}(f_+)<0$). Now we need a subsolution.
\medskip\\
Observe that the right hand side of \eqref{LPDE} is, for fixed $f,g,s$, a function of the distribution $H$ (if $V$ is defined to be $\bot_gH$), and that this function $\Distr_{n-1}(M)\to C^\infty(M,\R)$ is continuous with respect to the $C^2$-topology on $\Distr_{n-1}(M)$ and the $C^0$-topology on $C^\infty(M,\R)$. (One has to verify that all the involved functions depend only on the $2$-jet of $H\in C^\infty(M\ot\Gr_{n-1}(TM))$. This can be checked in a similar way as we verified in Chapter \ref{FIVE} that $\Twist_H$ depends only on the $1$-jet of $H$. The continuity of the RHS of \eqref{LPDE} is then obvious.)
\smallskip\\
Thus there is a $C^2$-neighbourhood $\mathscr{U}\in\Distr_{n-1}(M)$ of our given integrable distribution $H$ such that, for all $H'\in\mathscr{U}$, the function $f_+$ is a supersolution of the elliptic equation $\PD_{g,\bot_gH',s}(f)=0$.
\smallskip\\
By the $C^\infty$-approximation results of Section \ref{twistCinfty}, $\mathscr{U}$ contains an everywhere twisted distribution $H_0$ if $n\geq5$ and in many cases also if $n=4$. The function $f_+$ is a supersolution of $\PD_{g,\bot_gH_0,s}(f)=0$. By Lemma \ref{subsol1}, every sufficiently small (and this means \emph{really} small in our case here because we have already quite a small supersolution) constant is a subsolution.
\smallskip\\
Hence we have a supersolution and a smaller subsolution. The sub- and supersolution method proves that there is a solution, so we have found a Lorentzian metric with scalar curvature $s$ in the given connected component of $\Metr_1(M)$.
\medskip\\
This shows how Conjecture \ref{chinegative} implies the esc Conjecture in dimensions $\geq5$ and in many cases also in dimension $4$. (We will not discuss the remaining four-dimensional cases here. But note that in the orientable time-orientable case, Theorem \ref{fourCinfty} reduces the problem in such a way that it can be handled inside one manifold chart.)
\smallskip\\
Let me conclude with a few vague ideas of how to prove Conjecture \ref{chinegative}. Using \ref{lineformulae}, \ref{fscallemma}, and the integrability of $H$, we can write the function $\chi_{g,V}$ as follows:
\[ \begin{split}
\chi_{g,V} &= \scal_g +\xi_{g,V}\\
&= \scal^{H,H}_g -2\qual^V_g -2\qual^H_g +\xi_{g,V}\\
&= \fscal_{g,H} -\lvert{\divergence^H_g}\rvert_{g,V}^2 +\sigma_{g,H} -2\partial_V\divergence_g(V) -2\sigma_{g,V} +2\divergence_g(\nabla_VV) +2\sigma_{g,V} -2\tau_{g,H}\\
&\mspace{20mu}+2\partial_V\divergence_g(V) +\divergence_g(V)^2 +\frac{\sigma_{g,H}+\tau_{g,H}}{2}\\
&= \fscal_{g,H} +2\divergence_g(\nabla_VV) \;\;.
\end{split} \]
The problem with $\chi_{g,V}$ curvature is, as we have seen in Section \ref{chisection}, that it behaves under stretching along $V$ and under stretching along $H$ in a way which makes it impossible to construct metrics with negative $\chi_{g,V}$ by a simple scaling argument.
\smallskip\\
The idea why Conjecture \ref{chinegative} should be true (in contrast to the analogous statement for \emph{positive} $\chi_{g,V}$, to which we can easily construct counterexamples) is of course the philosophy that negative curvature is easy to produce. For instance, the ($1$-parametric, relative) h-principle techniques that J.\ Lohkamp developed for Riemannian metrics with negative scalar or Ricci curvature (cf.\ \cite{Lohkamp1995}) can be used to prove the following statement:
\smallskip\\
\emph{For every manifold of dimension $n\geq4$ and every integrable $(n-1)$-plane distribution $H$ on $M$, there exists a Riemannian metric $g$ on $M$ such that $\fscal_{g,H}$ is everywhere negative.}
\smallskip\\
On the other hand, the term $2\divergence_g(\nabla_VV)$ has no preferred sign. For instance, its mean value on a closed manifold is zero (like every mean value of the divergence of a vector field).
\medskip\\
In spite of this, it seems that Lohkamp's results cannot be applied directly to produce a metric with $\chi_{g,V}<0$: the two terms $\fscal_{g,H}$ and $2\divergence_g(\nabla_VV)$ are just not independent enough (although $\fscal_{g,H}$ depends only on the metric along the leaves of the foliation while $2\divergence_g(\nabla_VV)$ depends also on the transverse part). Part of the problem is that there are no $C^1$-dense h-principles for the negative scalar curvature relation, only $C^0$-dense ones.
\smallskip\\
I think that for a proof of Conjecture \ref{chinegative}, one has to start from scratch and prove flexibility results for a certain second-order partial differential relation on the leaves of the foliation\footnote{I am not going to write down the relation here since we will not discuss it anyway. It is a partial differential relation for a pair $(g,a)$ where $g$ is a Riemannian metric on the leaves and $a$ is a real-valued function on the leaves. The relation depends on a given closed $1$-form $\alpha$ on the leaves, and this $1$-form is determined by the choice of a line bundle transverse to the foliation and the choice of a Riemannian metric on this line bundle. In a sufficiently general situation, there will be no particularly clever choice of these background data; so one must prove a certain $1$-parametric relative $h$-principle for the relation, without having knowledge about $\alpha$. This is essentially a local problem: one may assume that $g$ and $a$ live on (some relatively compact open subset of) $\R^{n-1}$.}; these results would be analogous to (but not directly deducible from) Lohkamp's flexibility results for negative scalar and Ricci curvature. So probably no easy proof is available.
\medskip\\
Of course, there might still be a totally different way to prove the esc Conjecture, e.g.\ by doing more analysis to obtain better solvability criteria for the elliptic equation. But even then, a proof of Conjecture \ref{chinegative} would yield additional insight into the geometry of the problem.


\chapter{The two-dimensional Lorentzian case} \label{SEVEN}

The aim of the present chapter is to prove the theorems \ref{introclosedtwo} and \ref{introopentwo}, which solve the plain problem for Lorentzian surfaces. The main idea is that our elliptic equation from Theorem \ref{THEEQUATION} becomes considerably simpler when we choose the background data $V$ and $g$ in such a way that the distribution $V$ is parallel with respect to the Riemannian metric $g$. Such a choice is possible on every compact $2$-manifold which admits a Lorentzian metric --- but not inside \emph{each} homotopy class of line distributions, so we cannot solve the homotopy class problem in this way.
\smallskip\\
On manifolds with nonempty boundary, we solve the resulting simple equation by a variational technique. On closed manifolds (i.e.\ on the torus and the Klein bottle), we use again the Kazdan/Warner method.

\section{Simplification of the problem}

\subsection{Removing the $\abs{df}_g^2$ term}

For Lorentzian metrics on surfaces, we can write our elliptic equation so that no squares of first derivatives appear. This is a particular feature of the $2$-dimensional case: By a suitable substitution $f=F\compose u$ (where $I\subseteq\R$ is an open interval, $F\in C^\infty(I,\R_{>0})$ and $u\in C^\infty(M,I)$), we can always get rid of the $\abs{df}_g^2$ term in the definition \ref{PDdef} of the operator $\PD_{g,V,s}$; but if $1\leq q\leq n-1$, then the $\abs{df}_{g,V}^2$ term vanishes if and only if $n=2$ (because its coefficient $b_{n,q}(f)$ vanishes if and only if $q=1$ and $n-1=q$).

\begin{proposition} \label{twodimequationA}
Let $(M,g)$ be a $2$-dimensional Riemannian manifold, let $V$ be a line distribution on $M$, let $s\in C^\infty(M,\R)$, and let $H$ denote the $g$-orthogonal distribution of $V$. If the elliptic PDE
\begin{equation} \begin{split} \label{PDEtwoa}
0 &= 2\laplace_g(w) +\frac{2}{\sin(w)^2}\eval{\divergence^V_g}{dw}_{g,H} +\frac{2}{\cos(w)^2}\eval{\divergence^H_g}{dw}_{g,V} +\frac{1}{\sin(w)\cos(w)}\xi_{g,V} +\frac{\cos(w)}{\sin(w)}\scal_g\\
&\mspace{20mu}-\sin(w)\cos(w)s
\end{split} \end{equation}
has a solution $w\in C^\infty(M,\oointerval{0}{\frac{\pi}{2}})$, then there is a Lorentzian metric $h$ on $M$ with scalar curvature $s$, such that $V$ is timelike with respect to $h$, and $H$ is $h$-orthogonal to $V$.
\end{proposition}
\Proof
Since the twistedness of a line distribution vanishes, our usual elliptic equation from Theorem \ref{THEEQUATION} (cf.\ also \ref{lorentzequation}) has in the $2$-dimensional case the form
\[ \begin{split}
0 &= 2\laplace_g(f) -\frac{4f}{1+f^2}\abs{df}^2_g +\frac{2(1+f^2)}{f^2}\eval{\divergence^V_g}{df}_{g,H} +2(1+f^2)\eval{\divergence^H_g}{df}_{g,V}\\
&\mspace{20mu}+\frac{(1+f^2)^2}{f}\xi_{g,V} +\frac{1+f^2}{f}\scal_g -fs \;\;.
\end{split} \]
If this equation has a solution $f\in C^\infty(M,\R_{>0})$, then there is a Lorentzian metric on $M$ with scalar curvature $s$ which makes $V$ timelike and $H$ orthogonal to $V$.
\smallskip\\
Now we assume that \eqref{PDEtwoa} has a solution $w\in C^\infty(M,\oointerval{0}{\frac{\pi}{2}})$ and consider the function $f\define\tan(w)\in C^\infty(M,\R_{>0})$. Since $\tan'(w)=\frac{1}{\cos(w)^2} = 1+\tan(w)^2$ and $\tan''(w) = \frac{2\sin(w)}{\cos(w)^3}$, we obtain from \eqref{PDEtwoa} (after multiplication by $\frac{1}{\cos(w)^2}$):
\[ \begin{split}
0 &= \frac{2}{\cos(w)^2}\laplace_g(w) +\frac{4\sin(w)}{\cos(w)^3}\abs{dw}_g^2 -\frac{4\sin(w)}{\cos(w)^3}\abs{dw}_g^2 +\frac{2}{\sin(w)^2\cos(w)^2}\eval{\divergence^V_g}{dw}_{g,H}\\
&\mspace{20mu}+\frac{2}{\cos(w)^4}\eval{\divergence^H_g}{dw}_{g,V} +\frac{1}{\sin(w)\cos(w)^3}\xi_{g,V} +\frac{1}{\sin(w)\cos(w)}\scal_g -\frac{\sin(w)}{\cos(w)}s\\
&= 2\laplace_g(\tan(w)) -\frac{4\tan(w)}{(1+\tan(w)^2)\cos(w)^4}\abs{dw}_g^2 +\frac{2(1+\tan(w)^2)}{\tan(w)^2\cos(w)^2}\eval{\divergence^V_g}{dw}_{g,H}\\
&\mspace{20mu}+\frac{2(1+\tan(w)^2)}{\cos(w)^2}\eval{\divergence^H_g}{dw}_{g,V} +\frac{(1+\tan(w)^2)^2}{\tan(w)}\xi_{g,V} +\frac{1+\tan(w)^2}{\tan(w)}\scal_g -\tan(w)s\\
&= 2\laplace_g(f) -\frac{4f}{1+f^2}\abs{df}^2_g +\frac{2(1+f^2)}{f^2}\eval{\divergence^V_g}{df}_{g,H} +2(1+f^2)\eval{\divergence^H_g}{df}_{g,V}\\
&\mspace{20mu}+\frac{(1+f^2)^2}{f}\xi_{g,V} +\frac{1+f^2}{f}\scal_g -fs \;\;.
\end{split}\]
This implies the statement of the proposition.
\end{proof}

\begin{remark}
The $\xi_{g,V}$ and $\scal_g$ terms from Equation \eqref{PDEtwoa} can be written in a more symmetric form (cf.\ Definition \ref{lineobjects}):
\[
\frac{1}{\sin(w)\cos(w)}\xi_{g,V} +\frac{\cos(w)}{\sin(w)}\scal_g = \frac{2\sin(w)}{\cos(w)}\Big(\partial_V\divergence_g(V) +\divergence_g(V)^2\Big) -\frac{2\cos(w)}{\sin(w)}\Big(\partial_H\divergence_g(H) +\divergence_g(H)^2\Big) \;.
\]
\end{remark}
\Proof
Since $V$ and $H$ are line distributions, the formulae \ref{lineformulae} imply $\frac{1}{2}(\sigma_{g,H}+\tau_{g,H}) = \sigma_{g,H} = \eval{\divergence^H_g}{\divergence^H_g}_{g,V} = \divergence_g(V)^2$ and thus
\[ \begin{split}
\xi_{g,V} &= 2\partial_V\divergence_g(V) +2\divergence_g(V)^2 \;\;,\\
\scal_g &= -2(\qual^V_g +\qual^H_g) = -2\Big(\partial_V\divergence_g(V) +\divergence_g(H)^2 +\partial_H\divergence_g(H) +\divergence_g(V)^2\Big) \;\;.
\end{split} \]
Taking $\frac{1}{\sin(w)\cos(w)} -\frac{\cos(w)}{\sin(w)} = \frac{1-\cos(w)^2}{\sin(w)\cos(w)} = \frac{\sin(w)}{\cos(w)}$ into account, we get the claimed equation.
\end{proof}

The main advantage of writing our elliptic equation in the form of Proposition \ref{twodimequationA} is that in this way, it has variational form if the line distribution $V$ is $g$-parallel.

\subsection{Parallel line distributions}

Recall that a vector field $X$ on a semi-Riemannian manifold $(M,g)$ is called \emph{parallel} if and only if $\nabla X=0$, where $\nabla\colon C^\infty(M\ot TM)\to C^\infty(M\ot T^\ast M\otimes TM)$ denotes the \LeviCivita\ connection with respect to $g$. A line distribution $V$ on $(M,g)$ is called \emph{parallel} if and only if for every $x\in M$, there exist an open neighbourhood $U$ of $x$ and a parallel vector field $X\in C^\infty(U\ot V)$ such that $V\restrict U = \R X$ (i.e., the restriction of $V$ to $U$ is pointwise the span of $X$).
\smallskip\\
Moreover, recall the following basic facts:

\begin{facts} \label{parallelfacts}
Every parallel vector field on a Riemannian manifold has constant length. A constant-length section in a parallel line distribution is parallel. A Riemannian manifold admits a nowhere vanishing parallel vector field if and only if it admits an orientable parallel line distribution. The orthogonal distribution of any parallel line distribution on a Riemannian $2$-manifold is parallel. If a Riemannian $2$-manifold admits a parallel line distribution, then it is flat.
\end{facts}
\Proof
Every parallel vector field $X$ on a Riemannian manifold $(M,g)$ has constant length because $\partial_vg(X,X) = 2g(X,\nabla_vX)=0$ for all $v\in TM$. A constant-length section in a parallel line distribution is thus locally a constant multiple of a parallel vector field, hence parallel. If a Riemannian manifold admits a nowhere vanishing parallel vector field $X$ then it admits an orientable parallel line distribution, namely the span of $X$. Conversely, if a Riemannian manifold admits an orientable parallel line distribution $V$, then it admits a nowhere vanishing parallel vector field, for instance any unit-length section in $V$.
\smallskip\\
Let $V$ be a parallel line distribution on a Riemannian $2$-manifold, and let $H$ be the $g$-orthogonal distribution of $V$. There exist local unit-length sections $e_0,e_1$ in $V$ resp.\ $H$. Being locally a constant multiple of a parallel vector field, $e_0$ is parallel. By the usual ON frame rules from Chapter \ref{TWO}, we have $g(\nabla_{e_1}e_1,e_1)=0$ and $g(\nabla_{e_1}e_1,e_0) = -g(\nabla_{e_1}e_0,e_1) = 0$, hence $\nabla_{e_1}e_1=0$. Moreover, $g(\nabla_{e_0}e_1,e_0) = -g(\nabla_{e_0}e_0,e_1) = 0$ and $g(\nabla_{e_0}e_1,e_1)=0$, hence $\nabla_{e_0}e_1=0$. Thus $\nabla e_1=0$, i.e., $e_1$ is parallel. Since this holds for all local unit-length sections $e_1$ in $H$, the line distribution $H$ is parallel.
\smallskip\\
If a Riemannian $2$-manifold admits a parallel line distribution, then it thus admits local ON frames consisting of parallel vector fields. This implies that the manifold is flat.
\end{proof}

\begin{examples} \label{parallelexamples}
The euclidean metric on $\R^2$ induces a flat Riemannian metric on the $2$-torus $T^2=\R^2/\Z^2$. Analogously, it induces a flat Riemannian metric on the Klein bottle $\klein=\R^2/\Gamma$; cf.\ Example \ref{Kleinbottle}. (In fact, there is a double cover map $T^2\to\klein$ which is a local isometry with respect to these metrics; but that's not important for us.)
\smallskip\\
With respect to these metrics, $T^2$ and $\klein$ admit canonical nonvanishing parallel vector fields $X_T$ and $X_\klein$, respectively: both are induced by the unit vector field $e_1$ on $\R^2$ (where $(e_1,e_2)$ is the standard basis of $\R^2$). The orientable line distribution $H$ in Example \ref{Kleinbottle} is spanned by $X_\klein$.
\end{examples}

\begin{proposition} \label{parallelexistence}
If $M$ is a compact connected $n$-manifold with nonempty boundary, then $M$ admits a Riemannian metric with a nonvanishing parallel vector field. If $M$ is a nonempty connected closed $2$-manifold, then the following statements are equivalent:
\begin{enumerate}
\item
$M$ admits a Lorentzian metric.
\item
$M$ admits a Riemannian metric with a parallel line distribution.
\item
$M$ admits a Riemannian metric with a nonvanishing parallel vector field.
\item
$M$ is diffeomorphic to either the $2$-torus $T^2$ or the Klein bottle $\klein$.
\end{enumerate}
\end{proposition}
\Proof
By a theorem of P.\ Percell (cf.\ \cite{Percell1981}), every compact connected $n$-manifold with nonempty boundary admits a Riemannian metric with a nonvanishing parallel vector field. Now let $M$ be a nonempty closed $2$-manifold.
\smallskip\\
The diffeomorphy classification of nonempty connected closed $2$-manifolds tells us that there are exactly two diffeomorphism types of such manifolds with vanishing Euler characteristic, namely $T^2$ and $\klein$; cf.\ e.g.\ \cite{Hirsch}, Theorem 9.3.11. By Proposition 5.37 in \cite{ONeill}, a nonempty connected closed manifold admits a Lorentzian metric if and only if\footnote{For this equivalence, cf.\ also Proposition \ref{propeleven} and Theorem \ref{baumtheorem}.} it admits a line distribution if and only if its Euler characteristic is zero. Thus (ii)$\implies$(iv)$\iff$(i).
\smallskip\\
The implication (iii)$\implies$(ii) is trivial and has already been mentioned. In \ref{parallelexamples}, we have seen that $T^2$ and $\klein$ admit Riemannian metrics with parallel vector fields. Hence (iv)$\implies$(iii), so the proof is complete.
\end{proof}

\begin{remark}
If a $2$-manifold $M$ admits a Riemannian metric with a parallel line distribution, then in general not every homotopy class of line distributions on $M$ will contain a line distribution which is parallel with respect to some Riemannian metric on $M$. Consider the case $M=T^2$, for instance: Example \ref{torusexample} shows that there is a canonical bijection between $\Z\times\Z$ and the set of homotopy classes of line distributions on $T^2$. But only the homotopy class corresponding to $(0,0)\in\Z^2$ contains a distribution which is parallel with respect to some Riemannian metric on $T^2$.
\smallskip\\
Namely, assume that the line distribution $V$ on $T^2$ is parallel with respect to $g$. Via the projection $\pr\colon \R^2 \to \R^2/\Z^2 = T^2$, we pull back $V$ and $g$ to $\R^2$. Since the pullback metric $\tilde{g}$ on $\R^2$ is complete and flat, it is affinely isometric to the euclidean metric; cf.\ \cite{Wolf}, Corollary 1.9.6. The pullback distribution $\tilde{V}$ is $\tilde{g}$-parallel. Thus there exists an affine isomorphism $A\colon\R^2\to\R^2$ such that $\tilde{V}=A^\ast(\tilde{V}_{0,0})$; here $\tilde{V}_{0,0}$ denotes the (euclidean-flat) first-factor distribution on $\R\times\R$, i.e.\ the $\pr$-pullback of the distribution $V_{0,0}$ from Example \ref{torusexample}.
\smallskip\\
We may assume that $A$ is contained in the identity component of the group of affine automorphisms of $\R^2$, since otherwise we compose it with a reflection along some axis which is parallel to $\tilde{V}_{0,0}$. Now any path from $A$ to the identity defines a path from $\tilde{V}$ to $\tilde{V}_{0,0}$ in the space of line distributions on $\R^2$. Since each distribution on the path projects down to a line distribution on $T^2$, we get a path from $V$ to $V_{0,0}$ in the space of line distributions on $T^2$. Hence $V$ is contained in the homotopy class corresponding to $(0,0)$.
\end{remark}

\begin{proposition} \label{twodimequationB}
Let $(M,g)$ be a Riemannian $2$-manifold, let $V$ be a parallel line distribution on $(M,g)$, and let $s\in C^\infty(M,\R)$. If the elliptic equation
\begin{equation} \label{PDEtwob}
0 = \laplace_g(u) -\frac{s}{2}\sin(u)
\end{equation}
has a solution $u\in C^\infty(M,\oointerval{0}{\pi})$, then $M$ admits a Lorentzian metric with scalar curvature $s$ which makes $V$ timelike.
\end{proposition}
\Proof
Let $(e_0,e_1)$ be any $V$-adapted local $g$-orthonormal frame, where without loss of generality $e_0$ is a local section in $V$. By \ref{parallelfacts}, $e_0$ and $e_1$ are parallel, so all ON Christoffel symbols with respect to the frame $(e_0,e_1)$ vanish. In particular (cf.\ the formulae in Subsection \ref{TWOTWO}), we have $\divergence^V_g=0$, $\divergence^H_g=0$, $\xi_{g,V}=0$, and $\scal_g=0$. Thus Proposition \ref{twodimequationA} says that if the PDE
\[
0 = 2\laplace_g(w) -\sin(w)\cos(w)s
\]
has a solution $w\in C^\infty(M,\oointerval{0}{\frac{\pi}{2}})$, then there is a Lorentzian metric on $M$ with scalar curvature $s$ which makes $V$ timelike. Now assume that \eqref{PDEtwob} has a solution $u\in C^\infty(M,\oointerval{0}{\pi})$. Then the function $w\define u/2\in C^\infty(M,\oointerval{0}{\frac{\pi}{2}})$ satisfies $2\laplace_g(w) -\sin(w)\cos(w)s = \laplace_g(u) -\frac{s}{2}\sin(u) = 0$. This implies the statement of the proposition.
\end{proof}

It remains to prove that Equation \eqref{PDEtwob} has solutions. We have to distinguish two cases.


\section{The closed case}

This is the only place in the whole thesis where we encounter an obstruction to realising certain functions as scalar curvatures of pseudo-Riemannian metrics (with prescribed index). In fact, the obstruction is well-known: it is the Gauss/Bonnet theorem for Lorentzian surfaces (cf.\ \cite{Avez1962} or \cite{Alty1995} for the general semi-Riemannian Gauss/Bonnet theorem, \cite{BirmanNomizu} for the $2$-dimensional Lorentzian case).

\begin{theorem}[Gauss/Bonnet for Lorentzian surfaces] \label{LorentzGaussBonnet}
Let $M$ be a closed $2$-manifold. Then every Lorentzian metric $h$ on $M$ satisfies the equation
\[
\int_{(M,h)}\scal_h = 0 \;\;.
\]
In particular, if a function $s\in C^\infty(M,\R)$ is the scalar curvature of some Lorentzian metric on $M$, then either $s$ is the constant $0$ or $s$ changes its sign.
\end{theorem}
\Proof
Using the formulae from Chapters \ref{TWO} and \ref{THREE}, we can give our own proof (which is not new but just a special case of Avez' classical ``reduce to the Riemannian case'' proof): We choose an $h$-timelike line distribution $V$ on $M$ and consider the Riemannian metric $g\define\switch(h,V)$. Then the densities of $g$ and $h$ are equal; that is, on every open oriented subset of $M$, the volume forms $\vol_g$ and $\vol_h$ coincide. To see this, just choose a local oriented $V$-adapted $g$-orthonormal frame $(e_0,e_1)$ of $TM$. Then the $2$-form $\vol_g$ is completely determined by the function $\vol_g(e_0,e_1)$, which is just the constant $1$. Since $(e_0,e_1)$ is also a local oriented $h$-orthonormal frame, the function $\vol_h(e_0,e_1)$ is the constant $1$ as well. Thus $\vol_g=\vol_h$.
\smallskip\\
Now Theorem \ref{switchformulae} and \ref{lineformulae} yield
\[ \begin{split}
\scal_h -\scal_g &= -2\scal^{V,V}_g +4\qual^V_g -2\lvert\divergence^V_g\rvert_{g,H}^2 +2\lvert\divergence^H_g\rvert_{g,V}^2 -2\sigma_{g,V} +2\sigma_{g,H}\\
&= 4\qual^V_g +4\lvert\divergence^H_g\rvert_{g,V}^2 -4\tau_{g,V} \;\;.
\end{split} \]
An application of Lemma \ref{adjointdivergence} (choose $u\equiv1$ there) shows that
\[
\int_{(M,h)}\scal_h -\int_{(M,g)}\scal_g = \int_{(M,g)}(\scal_h-\scal_g)
= 4\int_{(M,g)}\Big(\qual^V_g +\lvert\divergence^H_g\rvert_{g,V}^2 -\tau_{g,V}\Big)
= 0 \;\;.
\]
Since every closed manifold which admits a Lorentzian metric has zero Euler characteristic, the classical Riemannian Gauss/Bonnet theorem implies
\[
\int_{(M,g)}\scal_g = 4\pi\chi(M) = 0
\]
and thereby completes the proof.
\end{proof}

Our aim is now to prove that the necessary condition from the preceding theorem is also sufficient; i.e., we want to prove that every function $s\in C^\infty(M,\R)$ on $M=T^2$ or $M=\klein$ which changes its sign is the scalar curvature of some Lorentzian metric on $M$. (The constant function $0$ on $M$ is the scalar curvature of an obvious Lorentzian metric.) Once again, we employ the Kazdan/Warner method.

\begin{definition} \label{twoPhi}
Let $M$ be either the $2$-torus or the Klein bottle equipped with its standard flat Riemannian metric $g$, and let $p\in\R_{>1}$. We define the map $\Phi\colon \Sob{p}{2}(M,\R)\times L^p(M,\R)\to L^p(M,\R)$ by
\[
\Phi(u,s)\define \laplace_g(u) -\frac{s}{2}\sin(u) \;\;.
\]
This map is well-defined because $\sin(u)\in L^\infty(M,\R)$ (or, alternatively, because $\sin(u)\in \Sob{p}{2}(M,\R)$) and thus $s\sin(u)\in L^p(M,\R)$ for all $u\in\Sob{p}{2}(M,\R)$ and $s\in L^p(M,\R)$.
\smallskip\\
For every $s\in L^p(M,\R)$, we define the map $\Phi_s\colon \Sob{p}{2}(M,\R)\to L^p(M,\R)$ by $\Phi_s(u) = \Phi(u,s)$.
\end{definition}

\begin{lemma} \label{twoPhilemma}
Let $M$ be either the $2$-torus or the Klein bottle equipped with its standard flat Riemannian metric $g$, and let $p\in\R_{>1}$. Then $\Phi$ is continuous. $\Phi_s$ is (\Frechet) differentiable for every $s\in L^p(M,\R)$, and its derivative $D_u\Phi_s\colon \Sob{p}{2}(M,\R)\to L^p(M,\R)$ in the point $u$ is given by
\[
(D_u\Phi_s)(v) = \laplace_g(v) -\frac{s}{2}\cos(u)v \;\;.
\]
The map $\Sob{p}{2}(M,\R)\times L^p(M,\R)\to \Lin(\Sob{p}{2}(M,\R),L^p(M,\R))$ given by $(u,s)\mapsto D_u\Phi_s$ is continuous.
\end{lemma}
\Proof
The map $\laplace_g\colon \Sob{p}{2}(M,\R)\to L^p(M,\R)$ is continuous and linear, in particular differentiable with $D_u\laplace_g = \laplace_g$ for each $u\in\Sob{p}{2}(M,\R)$. Multiplication $\Sob{p}{2}(M,\R)\times L^p(M,\R)\to L^p(M,\R)$ is well-defined and continuous because of $2-\dim(M)/p>0$; cf.\ \ref{multcont}. For the same reason, the map $\sin\colon \Sob{p}{2}(M,\R)\to \Sob{p}{2}(M,\R)$ is well-defined and differentiable with $(D_u\sin)(v) = \cos(u)v$; cf.\ \ref{chainrule}. These facts imply that $\Phi$ is continuous, and that $\Phi_s$ is differentiable for every $s\in L^p(M,\R)$, with $D_u\Phi_s$ given by the claimed formula. The map $(u,s)\mapsto D_u\Phi_s$ is continuous as a consequence of Lemma \ref{conti1}.
\end{proof}

\begin{remark}
When we compare Definition \ref{twoPhi} and Lemma \ref{twoPhilemma} with the contents of Section \ref{SIXKW}, we see that we don't have to choose $p>\dim(M)$ here; $p>\dim(M)/2$ suffices. (This would still be true if we were discussing Equation \eqref{PDEtwoa} instead of \eqref{PDEtwob}.) The reason is that our two-dimensional equations do not contain squares of first derivatives. But $f\mapsto \abs{df}_{g,V}^2$ is only well-defined as a map $\Sob{p}{2}(M,\R)\to L^p(M,\R)$ if $1-\dim(M)/p>0$; and this map will always occur in the $\geq3$-dimensional pseudo-Riemannian case.
\end{remark}

\begin{lemma} \label{twoinvertible}
Let $M$ be either the $2$-torus or the Klein bottle equipped with its standard flat Riemannian metric $g$. Then there is a function $u\in C^\infty(M,(0,\pi))$ such that $\cot(u)\laplace_g(u)\in C^\infty(M,\R)$ is everywhere nonnegative, and is zero only on a set of measure $0$.
\end{lemma}
\emph{Remark}. It would suffice for our application of the lemma to know that $\cot(u)\laplace_g(u)$ is everywhere nonnegative and not identically zero.
\Proof
For $[x,y]\in T^2 = \R^2/\Z^2$, we define
\[
u([x,y]) \define \frac{\pi}{2}+\sin(2\pi x) \;\;.
\]
This yields obviously a well-defined function $u\in C^\infty(T^2,(0,\pi))$. For the standard flat metric $g$ on $T^2$, we have $\laplace_g(u) = \partial_x\partial_xu +\partial_y\partial_yu = -4\pi^2\sin(2\pi x)$. Hence $(\laplace_gu)([x,y])$ is negative if $0<x<\frac{1}{2}$, and is positive if $\frac{1}{2}<x<1$. Since $u([x,y])$ is contained in $\oointerval{\frac{\pi}{2}}{\pi}$ if $0<x<\frac{1}{2}$ and contained in $\oointerval{0}{\frac{\pi}{2}}$ if $\frac{1}{2}<x<1$, we see that $\cot(u([x,y]))$ is negative if $0<x<\frac{1}{2}$, and is positive if $\frac{1}{2}<x<1$. Hence $\cot(u)\laplace_g(u)$ is everywhere nonnegative, and is zero only on the union of two circles $\set{[0],[\frac{1}{2}]}\times(\R/\Z)$, i.e.\ on a set of measure zero.
\smallskip\\
This takes care of the case $M=T^2$. In the case of the Klein bottle, which is also a quotient of euclidean $\R^2$, and is a circle bundle over $S^1$ (cf.\ \ref{Kleinbottle}), exactly the same formula as in the torus case above defines a suitable function $u\in C^\infty(\klein,(0,\pi))$ (with $x$ denoting the coordinate in the base of the circle bundle). The proof is the same as before.
\end{proof}

Now we can prove our main theorem about Lorentzian metrics on $T^2$ and $\klein$.

\begin{theorem}[Lorentzian metrics on closed $2$-manifolds] \label{closedtwotheorem}
Let $M$ be either the torus $T^2$ or the Klein bottle $\klein$, let $s\in C^\infty(M,\R)$ be a function which is either identically zero or changes its sign. Then there is a Lorentzian metric on $M$ whose scalar curvature is $s$.
\end{theorem}
\Proof
Let $g$ be the standard flat metric on $M$, which admits a parallel line distribution $V$; cf.\ \ref{parallelexamples}. Since Equation \eqref{PDEtwob} has in the case $s\equiv0$ the obvious solution $u\equiv\frac{\pi}{2}$, it remains to prove that every sign-changing function on $M$ is the scalar curvature of a Lorentzian metric. We choose a number $p>2$.
\smallskip\\
By Lemma \ref{twoinvertible}, there is a function $u_0\in C^\infty(M,(0,\pi))$ such that $\cot(u_0)\laplace_g(u_0)\in C^\infty(M,\R)$ is everywhere nonnegative, and is zero only on a set of measure $0$. We consider the linear elliptic second-order differential operator $P\colon C^\infty(M,\R)\to C^\infty(M,\R)$ given by
\[
P(v) \define \laplace_g(v) -\laplace_g(u_0)\cot(u_0)v \;\;.
\]
Theorem \ref{invertible} tells us that the operator $\bar{P}\colon \Sob{p}{2}(M,\R)\to L^p(M,\R)$ induced by $P$ is bijective.
\smallskip\\
Consider the function
\[
s_0\define \frac{2\laplace_g(u_0)}{\sin(u_0)} \in C^\infty(M,\R)
\]
and the differentiable function $\Phi\colon \Sob{p}{2}(M,\R)\times L^p(M,\R)\to L^p(M,\R)$ from Definition \ref{twoPhi}. We have $\Phi(u_0,s_0)=0$ and, by Lemma \ref{twoPhilemma}, $D_{u_0}\Phi_{s_0} = \bar{P}$; moreover, the map $(u,s)\mapsto D_u\Phi_s$ is continuous.
\smallskip\\
Since $D_{u_0}\Phi_{s_0}$ is bijective, the implicit function theorem \ref{implicitfunction} implies that there exist an open neighbourhood $\mathscr{U}$ of $s_0$ in $L^p(M,\R)$ and a continuous function $U\colon \mathscr{U}\to \Sob{p}{2}(M,\R)$ such that $\Phi(U(s),s) = 0$ for all $s\in\mathscr{U}$. On a perhaps smaller neighbourhood $\mathscr{U}'$ of $s_0$ in $L^p(M,\R)$, the map $U$ is $\Sob{p}{2}(M,\oointerval{0}{\frac{\pi}{2}})$-valued\footnote{Cf.\ Definition \ref{compdef}.} (because $\Sob{p}{2}(M,\oointerval{0}{\frac{\pi}{2}})$ is open in $\Sob{p}{2}(M,\R)$ and $U$ is continuous).
\smallskip\\
Now let $s\in C^\infty(M,\R)$ be any function which is somewhere positive and somewhere negative. Then there is a constant $c\in\R_{>0}$ such that $\inf(cs)\leq s_0\leq \sup(cs)$. The Kazdan/Warner approximation theorem \ref{KWapproximation1} yields a diffeomorphism $\varphi\in\Diff^0(M)$ with $cs\compose\varphi \in \mathscr{U}'$. The function $u\define U(cs\compose\varphi)\in \Sob{p}{2}(M,\R)$ solves the smooth elliptic equation $\Phi(u,cs\compose\varphi)=0$ and is therefore smooth; cf.\ \ref{sobolevregularity}.
\smallskip\\
Hence Proposition \ref{twodimequationB} implies that $M$ admits a Lorentzian metric $h$ with scalar curvature $cs\compose\varphi$. The Lorentzian metric $(\varphi^{-1})^\ast(ch)$ has scalar curvature $s$. This completes the proof.
\end{proof}

\begin{remark}
The proof of the preceding theorem shows in fact that there exists \emph{one} connected component of the space of Lorentzian metrics on $M$ which contains for \emph{every} sign-changing or zero $s\in C^\infty(M,\R)$ a metric with scalar curvature $s$. In other words, we have solved the homotopy class version of the prescribed scalar curvature problem for at least one connected component of the space of Lorentzian metrics on $M$. For $M=T^2$, this is, in the notation from \ref{torusexample} and \ref{distrcompdef}, the connected component $\tmc(V_{0,0})$.
\end{remark}

\begin{remark}
The solution above of the homotopy class problem for one specific connected component of $\Distr_1(M)$ suggests the following approach to the general homotopy class problem for Lorentzian metrics on closed $2$-manifolds: One can describe explicitly a representative for each homotopy class of line distributions on $T^2$ or $\klein$; for instance, the distributions $V_{k,l}$ from Example \ref{torusexample} are such representatives in the torus case. This allows us to compute, for the standard flat metrics on $T^2$ and $\klein$, explicitly the functions $\xi_{g,V}$, $\eval{\divergence^V_g}{dw}_{g,H}$, $\eval{\divergence^H_g}{dw}_{g,V}$ which appear in the PDE \eqref{PDEtwoa}. Like in the proofs of \ref{closedtwotheorem} and \ref{twoinvertible}, it might be possible to specify in each case explicit points $\in C^\infty(M,\oointerval{0}{\frac{\pi}{2}})$ at which the derivative of the right hand side of \eqref{PDEtwoa} is invertible. Then the homotopy class problem would be solved completely. It might actually be easy to do this; I just haven't tried yet.
\end{remark}


\section{The open case}

Equation \eqref{PDEtwob} has variational form, i.e., it is the Euler/Lagrange equation of a functional. This suggests that we should use direct methods in the calculus of variations to prove that \eqref{PDEtwob} admits a solution; that is, we should try to show that the functional has a minimum. However, this works only on compact manifolds with \emph{nonempty} boundary since the \Poincare\ inequality turns out to be crucial in our analysis.

\smallskip
Because we are not really interested in \eqref{PDEtwob} as a boundary value problem, we solve it only for constant boundary values $\frac{\pi}{2}$. Except for the condition that the solution have values in $\oointerval{0}{\pi}$, this is a standard exercise in the calculus of variations. Since our functional has such a simple form, we can argue on a quite elementary level instead of appealing to strong theorems.

\begin{definition}
Let $(M,g)$ be a compact Riemannian $2$-manifold, and let $s\in L^1(M,\R)$. We define the functional $E_s\colon \Sobzero{2}{1}(M,\R) \to \R$ by
\[
E_s(v) \define \int_{(M,g)}\Big(\abs{dv}_g^2 +s\sin(v)\Big) \;\;.
\]
Note that $E_s$ is well-defined since $v\mapsto \sin\compose\,v$ is well-defined as a map $\Sobzero{2}{1}(M,\R)\subseteq L^2(M,\R)\to L^\infty(M,\R)$, and since $v\mapsto \abs{dv}_g$ is well-defined as a map $\Sobzero{2}{1}(M,\R)\to L^2(M,\R)$.
\end{definition}

\begin{lemma} \label{EulerLagrange}
Let $(M,g)$ be a compact Riemannian $2$-manifold, and let $s\in C^\infty(M,\R)$. The functional $E_s\colon \Sobzero{2}{1}(M,\R) \to \R$ is (\Frechet) differentiable, and its derivative $D_vE_s\colon \Sobzero{2}{1}(M,\R) \to \R$ in the point $v$ is given by
\[
(D_vE_s)(w) = \int_{(M,g)}\Big(2\eval{dv}{dw}_g +s\cos(v)w\Big) \;\;.
\]
If $v\in\Sobzero{2}{1}(M,\R)$ is a critical point of $E_s$, then $v\in C^\infty(M,\R)$ and $2\laplace_g(v) -s\cos(v) = 0$.
\end{lemma}
\Proof
Consider the functions $f_0,f_1\colon\R\to\R$ defined by
\begin{align*}
f_0(x) &\define \begin{cases} \frac{x-\sin(x)}{x^2} &\text{if $x\neq0$}\\ 0 &\text{if $x=0$} \end{cases} \;\;,
&f_1(x) &\define \begin{cases} \frac{1-\cos(x)}{x^2} &\text{if $x\neq0$}\\ \tfrac{1}{2} &\text{if $x=0$} \end{cases} \;\;.
\end{align*}
Since they are continuous and bounded, there exists a $c\in\R$ such that $f_i\compose w\in L^\infty(M,\R)$ and $\norm{f_i\compose w}_{L^\infty}\leq c$ for all $w\in L^2(M,\R)$ and $i\in\set{0,1}$. Using Hölder's inequality, we obtain for all $v,w\in\Sobzero{2}{1}(M,\R)$:
\[ \begin{split}
\norm{\sin(v+w)-\sin(v)-\cos(v)w}_{L^1} &= \norm{\sin(v)\cos(w)+\cos(v)\sin(w)-\sin(v)-\cos(v)w}_{L^1}\\
&= \norm{\sin(v)(\cos(w)-1) +\cos(v)(\sin(w)-w)}_{L^1}\\
&\leq \norm{1-\cos(w)}_{L^1} +\norm{w-\sin(w)}_{L^1}\\
&= \norm{f_1(w)\cdot w^2}_{L^1} +\norm{f_0(w)\cdot w^2}_{L^1}\\
&\leq \norm{f_1(w)}_{L^\infty}\cdot\norm{w}_{L^2}^2 +\norm{f_0(w)}_{L^\infty}\cdot\norm{w}_{L^2}^2\\
&\leq 2c\norm{w}_{\Sob{2}{1}}^2 \;\;.
\end{split} \]
Thus
\[ \begin{split}
&\mspace{20mu}\abs{E_s(v+w)-E_s(v)-\int_{(M,g)}\Big(2\eval{dv}{dw}_g +s\cos(v)w\Big)}\\
&\leq \norm{\abs{d(v+w)}_g^2 -\abs{dv}^2_g -2\eval{dv}{dw}_g}_{L^1} +\norm{s\Big(\sin(v+w)-\sin(v)-\cos(v)w\Big)}_{L^1}\\
&\leq \norm{dw}_{L^2}^2 +2c\norm{s}_{L^\infty}\cdot\norm{w}_{\Sob{2}{1}}^2\\
&\leq \Big(1+2c\norm{s}_{L^\infty}\Big)\norm{w}_{\Sob{2}{1}}^2 \;\;.
\end{split}\]
This shows that $E_s$ is differentiable in $v$ and that $D_vE_s$ is given by the claimed formula.
\smallskip\\
For all $v\in\Sob{2}{1}(M,\R)$ and $w\in\Sobzero{2}{1}(M,\R)$, we have $\int_{(M,g)}\eval{dv}{dw}_g = -\int_{(M,g)}\laplace_g(v)w$, where $\laplace_g(v)\in\Sob{2}{-1}(M,\R)$. If $v$ is a critical point of $E_s$, i.e.\ $D_vE_s = 0$, then $-2\laplace_g(v) +s\cos(v) = 0\in\Sob{2}{-1}(M,\R)$ (by definition of what it means for an element of $\Sob{2}{-1}(M,\R)$ to vanish). Since $-2\laplace_g(v) +s\cos(v) = 0$ is an elliptic equation with smooth coefficients, elliptic regularity tells us that $v\in C^\infty(M,\R)$.
\end{proof}

Recall that a sequence $(x_k)_{k\in\N}$ in some Banach space $X$ over $\R$ \emph{converges weakly} to $x\in X$ if and only if for every element $\varphi$ of the dual space $X^\ast$, the sequence $\varphi(x_k)$ converges (in $\R$) to $\varphi(x)$; and recall that a function $E\colon X\to\R$ is \emph{(sequentially) weakly lower semicontinuous} if and only if every sequence $(x_k)_{k\in\N}$ in $X$ which converges weakly to some $x\in X$ satisfies the inequality $E(x)\leq\liminf_{k\to\infty}E(x_k)$.

\begin{lemma} \label{wlsc}
Let $(M,g)$ be a compact Riemannian $2$-manifold, and let $s\in L^p(M,\R)$ for some $p\in\R_{>1}$. Then $E_s\colon \Sobzero{2}{1}(M,\R) \to\R$ is weakly lower semicontinuous.
\end{lemma}
\Proof
The map $E\colon\Sobzero{2}{1}(M,\R)\to\R$ given by $v\mapsto \norm{dv}_{L^2}^2$ is obviously convex (i.e., $E(tv+(1-t)w) \leq tE(v)+(1-t)E(w)$ for all $t\in[0,1]$ and $v,w\in\Sobzero{2}{1}(M,\R)$) and lower semicontinuous (it is even continuous). It is thus weakly lower semicontinuous; cf.\ e.g.\ Theorem 3.1.2 in \cite{Dacorogna}. Since the sum of two weakly lower semicontinuous functionals is weakly lower semicontinuous, it remains to prove that the map $\tilde{E}\colon \Sobzero{2}{1}(M,\R)\to\R$ given by $v\mapsto \int_{(M,g)}s\sin(v)$ is weakly lower semicontinuous.
\smallskip\\
Let $(v_k)_{k\in\N}$ be a sequence in $\Sobzero{2}{1}(M,\R)$ which converges weakly to $v\in\Sobzero{2}{1}(M,\R)$. We will show that $(\tilde{E}(v_k))_{k\in\N}$ converges to $\tilde{E}(v)$. Let $q\in\R_{>1}$ denote the unique number such that $\frac{1}{p}+\frac{1}{q} = 1$. As a consequence of the Rellich/Kondrakov theorem (cf.\ \ref{rellich}), $(v_k)_{k\in\N}$ converges in $L^q$ to $v$. Hölder's inequality yields
\[
\lvert\tilde{E}(v_k)-\tilde{E}(v)\rvert \leq \norm{s\sin(v_k)-s\sin(v)}_{L^1} \leq \norm{s}_{L^p}\cdot\norm{\sin(v_k)-\sin(v)}_{L^q} \leq \norm{s}_{L^p}\cdot\norm{v_k-v}_{L^q} \;\;;
\]
hence $(\tilde{E}(v_k))_{k\in\N}$ converges to $\tilde{E}(v)$, as claimed. Since this is true for every weakly convergent sequence in $\Sobzero{2}{1}(M,\R)$, the map $\tilde{E}$ is weakly (lower semi)continuous.
\end{proof}

Recall that a functional $E\colon X\to\R$ on some Banach space $X$ is \emph{coercive} if and only if there exist $\alpha\in\R_{>0}$ and $\beta\in\R$ such that $E(x)\geq \alpha\norm{x}+\beta$ for all $x\in X$.

\begin{lemma} \label{coercive}
Let $(M,g)$ be a compact connected Riemannian $2$-manifold with nonempty boundary, and let $s\in L^1(M,\R)$. Then $E_s\colon\Sobzero{2}{1}(M,\R)\to\R$ is coercive.
\end{lemma}
\Proof
Since $M$ is connected with nonempty boundary, the \Poincare\ inequality gives us a $\tilde{c}\in\R_{>0}$ such that $\norm{dv}_{L^2}\geq\tilde{c}\norm{v}_{L^2}$ for all $v\in\Sobzero{2}{1}(M,\R)$. Thus $\norm{v}_{\Sob{2}{1}} = \norm{v}_{L^2} +\norm{dv}_{L^2} \leq c\norm{dv}_{L^2}$ for all $v\in\Sobzero{2}{1}(M,\R)$, where $c\define\tilde{c}+1$. Let $\alpha\define 2c\in\R_{>0}$ and $\beta\define -1-\norm{s}_{L^1}$. We obtain for all $v\in\Sobzero{2}{1}$:
\[
E_s(v) = \norm{dv}_{L^2}^2 +\int_{(M,g)}s\sin(v) \geq 2\norm{dv}_{L^2} -1 -\norm{s\sin(v)}_{L^1} \geq 2c\norm{dv}_{\Sob{2}{1}} -1-\norm{s}_{L^1} = \alpha\norm{v}_{\Sob{2}{1}} +\beta \;\;.\qedhere
\]
\end{proof}

\emph{Remark.} Note that arbitrarily large constant functions spoil the coercivity of $E_s$ on \emph{closed} manifolds.

\begin{proposition} \label{sinproposition0}
Let $(M,g)$ be a compact connected $2$-dimensional Riemannian manifold with nonempty boundary, and let $s\in C^\infty(M,\R)$. Then the boundary value problem
\begin{equation} \label{sinequation0}
0 = \laplace_g(u) -\frac{s}{2}\sin(u), \mspace{60mu} u\restrict\mfbd M = \tfrac{\pi}{2}
\end{equation}
has a solution $u\in C^\infty(M,\R)$.
\end{proposition}
\Proof
By \ref{wlsc} and \ref{coercive}, the functional $E_s\colon \Sobzero{2}{1}(M,\R)\to\R$ is weakly lower semicontinuous and coercive. Since the Banach space $\Sobzero{2}{1}(M,\R)$ is reflexive (being in fact a Hilbert space), the standard Theorem 3.1.1 in \cite{Dacorogna} tells us that there is a $v\in\Sobzero{2}{1}(M,\R)$ where $E_s$ achieves its minimum. In particular, $v$ is a critical point of $E_s$, and thus we get $v\in C^\infty(M,\R)$ and $2\laplace_g(v)-s\cos(v)=0$, by Lemma \ref{EulerLagrange}.
\smallskip\\
The function $u\define v+\frac{\pi}{2}\in C^\infty(M,\R)$ has the properties $\laplace_g(u) -\frac{s}{2}\sin(u) = \laplace_g(v) -\frac{s}{2}\cos(v) = 0$ and $u\restrict\mfbd M = \frac{\pi}{2}$.
\end{proof}

However, the preceding result is not good enough for our problem because we need a solution $u$ whose values lie strictly between $0$ and $\pi$ (cf.\ Proposition \ref{twodimequationB}). The following $C^0$ a priori estimate tells us that every solution has this property if $\norm{s}_{L^2}$ is sufficiently small. This is the second place where the nonemptiness of the manifold boundary is important.

\begin{lemma} \label{apriori}
Let $(M,g)$ be a compact connected Riemannian $2$-manifold with nonempty boundary. Then there is a constant $c\in\R_{>0}$ such that for every $s\in L^2(M,\R)$ and every solution $u\in\Sob{2}{2}(M,\R)\subseteq C^0(M,\R)$ of the equation $\laplace_g(u) = s\,\sin(u)$ with boundary values $u\restrict\mfbd M = \frac{\pi}{2}$, the following inequality holds:
\[
\norm{u-\tfrac{\pi}{2}}_{C^0}\leq c\,\norm{s}_{L^2} \;\;.
\]
\end{lemma}
\Proof
By the Sobolev imbedding theorem, there is a constant $c_0>0$ such that $\norm{v}_{C^0} \leq c_0\norm{v}_{\Sob{2}{2}}$ for all $v\in\Sob{2}{2}(M,\R)$. The elliptic estimate for the operator $\laplace_g\colon\Sob{2}{2}(M,\R)\to L^2(M,\R)$ yields a constant $c_1>0$ such that $\norm{v}_{\Sob{2}{2}} \leq c_1(\norm{v}_{L^2} +\norm{\laplace_g(v)}_{L^2})$ for all $v\in\Sob{2}{2}(M,\R)$. The \Poincare\ inequality gives us a $c_2>0$ such that $\norm{v}_{L^2}\leq c_2\norm{dv}_{L^2}$ for all $v\in\Sobzero{2}{1}(M,\R)$. We claim that the proposition holds with the constant $c \define c_0c_1(c_2^2+1)$.
\smallskip\\
Let $u\in\Sob{2}{2}(M,\R)$ be a solution of $\laplace_g(u) = s\,\sin(u)$ which is $\frac{\pi}{2}$ on the boundary, and define $v\define u-\frac{\pi}{2}$. Since $v$ vanishes on the boundary, we have $\int_{(M,g)}\eval{dv}{dv}_g = -\int_{(M,g)}\laplace_g(v)v$. Thus Hölder's inequality yields
\[
\norm{dv}_{L^2}^2
= \int_{(M,g)}\!\!\!\!-\laplace_g(v)v
= \int_{(M,g)}\!\!\!\!-\laplace_g(u)v
= \int_{(M,g)}\!\!\!\!-s\sin(u)v
\leq \norm{sv}_{L^1}
\leq \norm{s}_{L^2}\cdot\norm{v}_{L^2}
\leq c_2\norm{s}_{L^2}\cdot\norm{dv}_{L^2} \;,
\]
i.e.\ $\norm{dv}_{L^2} \leq c_2\norm{s}_{L^2}$. We get
\[ \begin{split}
\norm{u-\tfrac{\pi}{2}}_{C^0} &= \norm{v}_{C^0} \leq c_0\norm{v}_{\Sob{2}{2}} \leq c_0c_1\big(\norm{v}_{L^2} +\norm{\laplace_g(v)}_{L^2}\big)\\
&\leq c_0c_1\big(c_2\norm{dv}_{L^2} +\norm{s\sin(u)}_{L^2}\big) \leq c_0c_1\big(c_2^2\norm{s}_{L^2} +\norm{s}_{L^2}\big) = c\norm{s}_{L^2} \;\;,
\end{split} \]
as claimed.
\end{proof}

\begin{corollary} \label{sinproposition}
Let $(M,g)$ be a connected compact $2$-dimensional Riemannian manifold with nonempty boundary. Then there is a $\delta\in\R_{>0}$ such that, for every function $s\in C^\infty(M,\R)$ with $\norm{s}_{L^2}\leq\delta$, the boundary value problem
\begin{equation} \label{sinequation}
0 = \laplace_g(u) -\frac{s}{2}\sin(u), \mspace{60mu} u\restrict\mfbd M = \tfrac{\pi}{2}
\end{equation}
has a solution $u\in C^\infty(M,(0,\pi))$.
\end{corollary}
\Proof
Let $c\in\R_{>0}$ be the constant from Lemma \ref{apriori}. We choose $\delta\in\R_{>0}$ so small that every $s\in C^\infty(M,\R)$ with $\norm{s}_{L^2}\leq\delta$ satisfies $\frac{c}{2}\norm{s}_{L^2} < \frac{\pi}{2}$. Now the statement follows from Proposition \ref{sinproposition0}.
\end{proof}

\renewcommand{\mod}{\mt{ mod }}
\begin{remark}
By a trick that Matthias Kurzke pointed out to me, the boundary value problem \eqref{sinequation} has a solution $u\in C^\infty(M,[0,\pi])$ for \emph{every} $s\in C^\infty(M,\R)$ (not necessarily satisfying the smallness condition $\norm{s}_{L^2}\leq\delta$).
\smallskip\\
Namely, let us look at the proof of the standard result (Theorem 3.1.1 in \cite{Dacorogna}) on which our proof of Proposition \ref{sinproposition0} is based. It employs a minimising sequence $(u_k)_{k\in\N}$ in $\Sobzero{2}{1}(M,\R)$ of the functional $E_s$ (recall that \emph{minimising} means that the sequence converges to the infimum of the functional $E_s$). We define a new sequence $(\tilde{u}_k)_{k\in\N}$ in $\Sobzero{2}{1}(M,\R)$ by
\[
\tilde{u}_k \define \mu\compose u_k \;\;,
\]
where the (sawtooth-shaped) function $\mu\colon\R\to[-\frac{\pi}{2},\frac{\pi}{2}]$ is defined as follows: Let $x\mod2\pi$ denote the unique number $y\in[-\frac{\pi}{2},\frac{3\pi}{2})$ such that a $k\in\Z$ exists with $x = 2k\pi +y$. Then
\[
\mu(x) \define \begin{cases} x\mod2\pi & \text{if $(x\mod2\pi)\in[-\tfrac{\pi}{2},\tfrac{\pi}{2}]$}\\
\pi-(x\mod2\pi) & \text{if $(x\mod2\pi)\in[\tfrac{\pi}{2},\tfrac{3\pi}{2}]$} \end{cases} \;\;.
\]
Since $\mu$ is a Lipschitz continuous function with $\mu(0)=0$ and left composition with such functions maps elements of $\Sobzero{2}{1}(M,\R)$ to $\Sobzero{2}{1}(M,\R)$, the elements of the sequence $(\tilde{u}_k)$ are indeed contained in $\Sobzero{2}{1}(M,\R)$.

\smallskip
Moreover, $E_s(u_k) = E_s(\tilde{u}_k)$ for all $k\in\N$ because $\sin(\mu(x))=\sin(x)$ for all $x\in\R$, and because the absolute value of the derivative of $\mu$ is almost everywhere $1$, so $\norm{d\tilde{u}_k}_{L^2}$ and $\norm{du_k}_{L^2}$ are equal. In particular, $(\tilde{u}_k)_{k\in\N}$ is a minimising sequence of $E_s$. Now the proof of the standard theorem shows that $(\tilde{u}_k)_{k\in\N}$ has a weakly convergent subsequence whose limit $u\in\Sobzero{2}{1}(M,\R)$ is a minimum of $E_s$.
\smallskip\\
Since weak convergence in $\Sob{2}{1}(M,\R)$ implies convergence in $L^2(M,\R)$, we conclude that the limit $u$ has almost everywhere values in $[0,\pi]$. Now the arguments from the proof of \ref{sinproposition0} show that $u$ is contained in $C^\infty(M,[0,\pi])$ and solves Equation \eqref{sinequation}, as claimed.
\medskip\\
But in spite of this result, it seems that we cannot avoid a smallness condition like $\norm{s}_{L^2} < \delta$ in order to get a solution with values in $\oointerval{0}{\pi}$.
\end{remark}

Now we come to our main result for the Lorentzian prescribed scalar curvature problem on open $2$-manifolds.

\begin{theorem} \label{twoopenmain}
Let $M$ be a connected compact $2$-manifold with nonempty boundary, let $s\in C^\infty(M,\R)$. Then there is a Lorentzian metric on $M$ with scalar curvature $s$.
\end{theorem}
\Proof
By Proposition \ref{parallelexistence}, $M$ admits a Riemannian metric $g$ with a parallel line distribution $V$. Corollary \ref{sinproposition} shows that there exist a (small) constant $c\in\R_{>0}$ and a function $u\in C^\infty(M,\oointerval{0}{\pi})$ with $0 = \laplace_g(u) -\frac{cs}{2}\sin(u)$. Now Proposition \ref{twodimequationB} tells us that $cs$ is the scalar curvature of some Lorentzian metric $h$ on $M$. Hence $s$ is the scalar curvature of the Lorentzian metric $ch$.
\end{proof}

\begin{remark}[lower regularity]
If the prescribed function $s$ is not smooth, then still something can be said about solutions of Equation \eqref{sinequation} and thus about solutions of the prescribed scalar curvature problem. As an example, let us consider the case where $s$ is contained in $L^2(M,\R)$ (but is not necessarily continuous); I leave it to the reader to state all sorts of similar results under various regularity assumptions on $s$.
\smallskip\\
If $s\in L^2(M,\R)$, then weak lower semicontinuity and coercivity of the functional $E_s$ are still satisfied, as we have proved. Hence $E_s$ assumes its minimum in some point $v\in\Sobzero{2}{1}(M,\R)$. As a modification of Lemma \ref{EulerLagrange}, we claim that $E_s$ is \Gateaux\ differentiable and that its \Gateaux\ derivative is given by the formula in \ref{EulerLagrange}. (If $s$ is contained in $L^\infty(M,\R)$, then the proof of \ref{EulerLagrange} goes through and shows that $E_s$ is even \Frechet\ differentiable.)
\smallskip\\
Namely, for all $t\in\R$ and all $w\in\Sobzero{2}{1}(M,\R)$ with norm $1$, we have (in the notation of \ref{EulerLagrange}):
\[
\norm{s\Big(\sin(v+tw)-\sin(v)-\cos(v)tw\Big)}_{L^1} \leq \norm{s}_{L^2}\cdot\norm{\sin(v+tw)-\sin(v)-\cos(v)tw}_{L^2}
\]
and
\[ \begin{split}
\norm{\sin(v+tw)-\sin(v)-\cos(v)tw}_{L^2} &\leq \norm{f_1(tw)t^2w^2}_{L^2} +\norm{f_0(tw)t^2w^2}_{L^2}\\
&\leq t^2\Big(\norm{f_1}_{L^\infty} +\norm{f_0}_{L^\infty}\Big)\norm{w^2}_{L^2} \;\;.
\end{split}\]
Now the remaining arguments from the proof of \ref{EulerLagrange} show that $E_s$ is \Gateaux\ differentiable with the derivative given by the same formula as before.
\smallskip\\
Since the \Gateaux\ derivative vanishes in the point $v$ where $E_s$ achieves its minimum, we obtain again the equation $2\laplace_g(v)-s\cos(v) = 0\in\Sob{2}{-1}(M,\R)$. Since $\cos(v)\in L^\infty(M,\R)$ and thus $s\cos(v)\in L^2(M,\R)$, elliptic regularity tells us that $v\in\Sobzero{2}{2}(M,\R)\subseteq C^0(M,\R)$. Lemma \ref{apriori} applies again and provides a $u\in\Sob{2}{2}(M,\R)\subseteq C^0(M,\R)$ with values in $\oointerval{0}{\pi}$ which satisfies the equation $\laplace_g(u) = \frac{s}{2}\sin(u)$, provided $\norm{s}_{L^2}$ is sufficiently small. The rest of our arguments above goes through and yields a Lorentzian metric on $M$ with scalar curvature $s$.
\medskip\\
To summarise: \emph{Let $M$ be a connected compact $2$-manifold with nonempty boundary, let $s\in L^2(M,\R)$. Then there is a Lorentzian metric on $M$ of Sobolev regularity $\Sob{2}{2}$ whose scalar curvature is $s$.}
\end{remark}

\begin{finalremark}
What about the \emph{homotopy class problem} on open $2$-manifolds? To solve that, one would have to find solutions of Equation \eqref{PDEtwoa}, without possessing any useful knowledge about $\divergence^V_g$, $\divergence^H_g$, $\xi_{g,V}$. However, for what it's worth, by the results of Kazdan and Warner (cf.\ Theorem \ref{KWopen}), one can choose the Riemannian background metric $g$ in such a way that $\scal_g$ is any function one wants. I have not yet thought about this homotopy class problem, so I add it to the long list of questions that the present work leaves open.
\end{finalremark}


\begin{appendix}


\chapter{Topological miscellanea}

\section{Background in differential topology} \label{diffbackground}

\subsection{Jets and $C^r$-topologies} \label{jets}

\begin{notation} \label{Crnotation}
Let $p\colon E\to M$ be a (smooth) fibre bundle, and let $r\in\N\cup\set{\infty}$. We denote the space of $C^r$ (i.e.\ $r$-times continuously differentiable) sections in the bundle $p$ by $C^r(M\text{$\xleftarrow{p}$}E)$, or, when the projection map is clear from the context, by $C^r(M\ot E)$. When $U$ is a subset of $M$, we write simply $C^r(U\ot E)$ instead of $C^r(U\ot E\restrict U)$.
\end{notation}

\begin{notation}[jet bundles]
For an introduction to the language of jets and jet bundles, cf.\ \cite{Saunders}, in particular Chapter 4. Our notation of the relevant objects is as follows.
\smallskip\\
If $p\colon E\to M$ is a (smooth) fibre bundle, then we denote the total space of its $r$th-order jet bundle by $J^rE$; the projection map $p$ is suppressed in this notation, but that won't cause any confusion. This manifold $J^rE$ is the total space of several (smooth) bundles: It is the total space of a fibre bundle $p^r\colon J^rE\to M$, whose fibre over $x\in M$ we denote by $J^r_xE$. If $r=1$ for instance, then there is a projection $p^1_0 \colon J^1E \to E$ which turns $J^1E$ into an affine bundle modelled on the vector bundle $p^\ast(T^\ast M)\otimes\ker(Tp)$ over $E$; here $\ker(Tp)$ is defined by the vector bundle morphism $Tp\colon TE\to TM$ and forms a sub vector bundle of $TE\to E$.\footnote{This vector bundle $\ker(Tp)\to E$ is usually called the \emph{vertical bundle (of $p$)} and denoted by $VE$ or $Vp$, but there are already too many $V$s floating around in this thesis where the letter stands for \emph{vertical} as well.} These bundle projections are related by $p^1 = p\compose p^1_0 \colon J^1E \to M$.
\smallskip\\
We will use M.\ Gromov's h-principle methods in Chapter \ref{FIVE}. Following Gromov, the manifold $J^1E$ is usually denoted by $E^{(1)}$ or $E^1$ in the literature on the h-principle (e.g.\ \cite{GromovPDR}, \cite{Spring}, \cite{Geigeshprinciple}, \cite{EliashbergMishachev}). We will not adopt that notation.
\smallskip\\
For $f\in C^r(M\ot E)$, we denote the $r$-jet prolongation of $f$ by $j^rf$, as usual; $j^rf$ is a continuous section in the bundle $p^r\colon J^rE\to M$. In the case $r=1$ for instance, we have $p^1_0\compose j^1f = f$. We denote the value of $j^rf$ in the point $x\in M$ by $j^r_xf \in J^r_xE$. (By definition of the jet bundle, every $\sigma\in J^rE$ has the form $j^r_xf$, where $U$ is an open neighbourhood of $x=p^r(\sigma)$ and $f\in C^r(U\ot E)$.)
\end{notation}

\begin{definition}[topologies on sets of fibre bundle sections] \label{sectiontopologies}
Let $p\colon E\to M$ be a (smooth) fibre bundle, and let $r\in\N\cup\set{\infty}$.
\smallskip\\
Recall that the \emph{compact-open $C^r$-topology} on the set $C^r(M\ot E)$ is defined as follows. There is a canonical inclusion $C^r(M\ot E)\hookrightarrow C^0(M\ot J^rE) \subseteq C^0(M,J^rE)$ which sends each $C^r$ section to its $r$-jet prolongation. The compact-open $C^r$-topology on $C^r(M\ot E)$ is the topology induced by the compact-open topology on $C^0(M,J^rE)$. (Note that $J^\infty E$ is an infinite-dimensional manifold in a certain sense; cf.\ \cite{Saunders}, Chapter 7. But this is not relevant here since all we need is the topology on $J^\infty E$, which is the inverse limit of an inverse system in the category of topological spaces and continuous maps, consisting of the spaces $J^rE$, where $r\in\N$, and the projections $J^sE\to J^rE$, where $r\leq s$.)
\smallskip\\
Alternatively, we can consider the trivial fibre bundle $M\times E\to M$ and the corresponding jet space $J^r(M\times E)$ (which is usually denoted by $J^r(M,E)$), and define the compact-open $C^r$-topology on $C^r(M\ot E)$ via the inclusion $C^r(M\ot E) \subseteq C^r(M,E) \hookrightarrow C^0(M,J^r(M,E))$. Both definitions are equivalent.
\medskip\\
Recall that the \emph{fine} (synonymously:\ \emph{Whitney}) \emph{$C^r$-topology} on the set $C^r(M\ot E)$ is defined as follows. For every section $f\in C^r(M\ot E)$ and every neighbourhood $\mathscr{N}$ of the image of its $r$-jet prolongation $(j^rf)(M)$ in $J^rE$, we consider the set $\mathscr{U}(f,\mathscr{N})$ consisting of all $g\in C^r(M\ot E)$ such that $(j^rg)(M)\subseteq\mathscr{N}$. The collection of all these sets $\mathscr{U}(f,\mathscr{N})$ is the basis of a topology on $C^r(M\ot E)$, namely the fine $C^r$-topology.

\smallskip
Let $s\in\N\cup\set{\infty}$ with $r\leq s$. The compact-open [resp.\ fine] $C^r$ topology on $C^s(M\ot E)\subseteq C^r(M\ot E)$ is the subspace topology induced by the compact-open [fine] $C^r$-topology on $C^r(M\ot E)$.
\end{definition}

Let us recall briefly some standard facts about compact-open topologies. We assume that the set $C^r(M\ot E)$ is equipped with the compact-open $C^r$-topology, unless explicitly stated otherwise.

\begin{facts} \label{topologyfacts}
Let $r,s\in\N\cup\set{\infty}$ with $r\geq s$, and let $E\to M$ be a (smooth) fibre bundle.
\smallskip\\
The fine $C^r$-topology on the set $C^r(M\ot E)$ is finer\footnote{This is a (non-strict) order relation, i.e., it includes the case of equality.} than the compact-open $C^r$ topology; both are equal if $M$ is compact. The compact-open [resp.\ fine] $C^r$-topology on $C^r(M\ot E)$ is finer than the subspace topology induced by the compact-open [resp.\ fine] $C^s$-topology on $C^s(M\ot E)$. The topological space $C^r(M\ot E)$ is second countable, and its topology is induced by a complete metric. $C^r(M\ot E)$ is a dense subset of $C^s(M\ot E)$. If $M$ is compact, then $C^r(M\ot E)$ is a \Frechet\ manifold (even a Banach manifold if $r<\infty$), in particular locally contractible. If $E_0\to M$ is an open subbundle of $E\to M$, then $C^r(M\ot E_0)$ is an open subspace of $C^r(M\ot E)$. In particular, if $U\to M$ is an open subbundle of the jet bundle $J^rE\to M$, then $C^0(M\ot U)$ is an open subspace of $C^0(M\ot J^rE)$, hence $\set{f\in C^r(M\ot E) \suchthat \forall x\in M\colon j^r_xf\in U_x}$ is an open subset of $C^r(M\ot E)$. If $E\to M$ is a vector bundle, then $C^r(M\ot E)$ is a \Frechet\ space; it is even a Banach space if $M$ is compact and $r<\infty$. The inclusion $C^r(M\ot E) \to C^s(M\ot E)$ is (continuous and) a homotopy equivalence; in particular, it induces a bijection $\pi_0(C^r(M\ot E)) \to \pi_0(C^s(M\ot E))$ between the sets of path components.
\end{facts}
\Proof
Cf.\ \cite{Hirsch}, \cite{Palais1966}.
\end{proof}

Note that we distinguish between compact-open and fine topologies, although the main results of this thesis are restricted to the case of \emph{compact} manifolds, where both coincide. We do so because whenever compactness is not essential (e.g.\ for the main theorems in Chapter 5), we will explain the general case, in order to point out clearly where compactness is indispensable. Namely, it is indispensable in Chapters 6 and 7, where we have to solve partial differential equations.

\begin{remark} \label{fineremark}
Let $E\to M$ be a vector bundle over a noncompact connected manifold $M$, let $r\in\N\cup\set{\infty}$. The vector space $V\define C^r(M\ot E)$ equipped with the \emph{fine} $C^r$-topology is not a topological vector space: addition $V\times V\to V$ is continuous, but scalar multiplication $\R\times V\to V$ is not. But $V$ is a topological module over the ring $C^r(M,\R)$ when the latter is equipped with the fine $C^r$-topology; i.e., scalar multiplication $C^r(M,\R)\times V\to V$ is continuous.
\smallskip\\
This failure of the fine $C^r$-topology is closely related to the reason why one usually doesn't consider connected components of the space $C^r(M\ot E)$ equipped with the fine $C^r$-topology, where $E\to M$ is a fibre bundle over a noncompact connected manifold. For example, the space of all Riemannian metrics on $M$ is a \emph{convex} subset of the vector space $C^\infty(M\ot\Sym(TM))$, but it has uncountably many connected components with respect to the fine $C^\infty$-topology. With respect to the compact-open $C^\infty$-topology, it has only one connected component because $C^\infty(M\ot\Sym(TM))$ is a topological vector space then (cf.\ \ref{topologyfacts}) and thus the obvious straight path between any two metrics is continuous.
\end{remark}

\subsection{Approximations}

\begin{theorem}[smooth approximation of fibre bundle sections] \label{smoothapproximation}
Let $r\in\N\cup\set{\infty}$; let $\pi\colon E\to M$ be a smooth fibre bundle; let $K$ be a closed subset of $M$; let $U$ be an open neighbourhood of $K$; let $h\in C^r(M\ot E)$ be a $C^r$-section whose restriction to $U$ is smooth; and let $\mathscr{N}$ be a neighbourhood of the image of $j^rh \colon M\to J^rE$. Then there is a continuous (with respect to the compact-open $C^r$-topology on $C^r(M\ot E)$) map $H\colon [0,1]\to C^r(M\ot E)$ with $H_0 = h$ such that the map $\bar{H}\colon [0,1]\times M\to E$ given by $(t,x)\mapsto H_t(x)$ is of class $C^r$ and, moreover, the following properties are satisfied:
\begin{itemize}
\item
$H_t\restrict K = h\restrict K$ for all $t\in[0,1]$.
\item
For each $t\in[0,1]$, the image of $j^rH_t \colon M\to J^rE$ is contained in $\mathscr{N}$.
\item
$H_1$ is smooth.
\end{itemize}
\end{theorem}
\Proof
Cf.\ \cite{Spring}, Theorem 1.1.
\end{proof}

We cannot expect that the preceding theorem generalises without modification from the smooth to the real-analytic category: Given a section $h$ which is real-analytic on some subset $K\subseteq M$ with nonempty interior, its desired real-analytic continuation $H_1$ (if it exists at all) is completely fixed on each connected component of $M$ which intersects the interior of $K$ nontrivially; in particular, we cannot arrange in general that $H_1$ lies in a given neighbourhood of $h$. We thus have to omit this ``rel $K$'' part of the theorem in order to find a real-analytic analogue.
\smallskip\\
The remaining parts are still true. This is proved, but not explicitly stated, in \cite{Shiga1964} (proof of Theorem 2 on p.~135). Probably there exists an explicit reference somewhere else, but I don't know one.

\begin{theorem}[real-analytic approximation of fibre bundle sections] \label{analyticapproximation}
Let $r\in\N\cup\set{\infty}$; let $\pi\colon E\to M$ be a real-analytic fibre bundle; let $h\in C^r(M\ot E)$ be a $C^r$-section; and let $\mathscr{N}$ be a neighbourhood of the image of $j^rh \colon M\to J^rE$. Then there is a continuous (with respect to the compact-open $C^r$-topology on $C^r(M\ot E)$) map $H\colon [0,1]\to C^r(M\ot E)$ with $H_0 = h$ such that the map $\bar{H}\colon [0,1]\times M\to E$ given by $(t,x)\mapsto H_t(x)$ is of class $C^r$ and, moreover, the following properties are satisfied:
\begin{itemize}
\item
For each $t\in[0,1]$, the image of $j^rH_t \colon M\to J^rE$ is contained in $\mathscr{N}$.
\item
$H_1$ is real-analytic.
\end{itemize}
\end{theorem}
\begin{proof}[\textsc{Remarks on the proof}]
The idea is to approximate the section $h$ by a real-analytic map $\in C^r(M,E)$ (which might not be a section in the bundle $\pi$); this is known to be possible. Then we have to deform the approximating map into a section, which can be done by a standard tubular neighbourhood technique (where the tubular neighbourhood is taken with respect to a real-analytic metric).
\end{proof}

\begin{corollary} \label{radense}
Let $r\in\N\cup\set{\infty}$, and let $E\to M$ be a smooth [real-analytic] fibre bundle. Then the set of smooth [real-analytic] sections in $E$ is dense in $C^r(M\ot E)$ with respect to the fine $C^r$-topology.
\end{corollary}
\Proof
We apply Theorem \ref{smoothapproximation} in the case $K=U=\leer$ resp.\ Theorem \ref{analyticapproximation}. Then the statement follows immediately from the definition of the fine $C^r$-topology.
\end{proof}

\begin{corollary}[real-analytic approximation of diffeomorphisms] \label{realanalyticd}
Let $M$ be a real-analytic manifold, and let $r\in\N_{\geq1}\cup\set{\infty}$. We equip the set $\Diff^{(r)}(M)$ of all $C^r$ diffeomorphisms $M\to M$ with the subspace topology induced by the fine $C^r$-topology on $C^r(M,M)$. Then the subset of real-analytic diffeomorphisms is dense in $\Diff^{(r)}(M)$.
\end{corollary}
\Proof
The set $\Diff^{(r)}(M)$ is open in $C^r(M,M)$ with respect to the fine $C^r$-topology; cf.\ \cite{Hirsch}, Theorem 2.1.7. Now we can apply Theorem \ref{analyticapproximation} or \cite{Hirsch}, Theorem 2.5.1.
\end{proof}

\begin{theorem} \label{realanalyticatlas}
Let $M$ be a smooth manifold. Then there is a real-analytic structure on $M$ (i.e., a real-analytic atlas which is compatible with the given smooth atlas). If $\mathscr{A}_0$ and $\mathscr{A}_1$ are two real-analytic structures on $M$, then there is a real-analytic diffeomorphism $(M,\mathscr{A}_0)\to(M,\mathscr{A}_1)$ (which is of course a smooth diffeomorphism $M\to M$).
\end{theorem}
\Proof
For existence, cf.\ e.g.\ the remarks on p.~66 of \cite{Hirsch}. For uniqueness, approximate the smooth diffeomorphism $\id_M\colon (M,\mathscr{A}_0)\to(M,\mathscr{A}_1)$ by a real-analytic one, employing \ref{realanalyticd}.
\end{proof}

\begin{lemma} \label{smoothpath}
Let $E\to M$ be a (smooth) fibre bundle. If the sections $s_0$ and $s_1$ are contained in the same path component of $C^\infty(M\ot E)$ (with respect to the compact-open $C^\infty$- or $C^0$-topology\footnote{This doesn't make a difference; cf.\ \ref{topologyfacts}.}), then there is a path $s\colon[0,1]\to C^\infty(M\ot E)$ from $s_0$ to $s_1$ which is smooth in the sense that the map $[0,1]\times M\to E$ given by $(t,x)\mapsto s(t)(x)$ is smooth.
\end{lemma}
\Proof
Via the obvious projection $[0,1]\times M\to M$, we pull back the bundle $E\to M$ to a fibre bundle $\tilde{E}\to [0,1]\times M$. Since $s_0$ and $s_1$ lie in the same path component of $C^0(M\ot E)$, there is a continuous section $\tilde{s}\in C^0([0,1]\times M\ot \tilde{E})$ such that $\tilde{s}\restrict(\set{0}\times M)=s_0$ and $\tilde{s}\restrict(\set{1}\times M)=s_1$. By Theorem \ref{smoothapproximation} (applied to the case $K=\set{0,1}\times M$), there is a smooth section $s\in C^\infty([0,1]\times M\ot \tilde{E})$ such that $s\restrict(\set{0,1}\times M) = \tilde{s}\restrict(\set{0,1}\times M)$. (Theorem \ref{smoothapproximation} assumes that $\tilde{s}$ is already smooth on a neighbourhood of $K$, but this can be arranged by extending $\tilde{E}$ to a pull-back bundle $\bar{E}$ over $[-1,2]\times M$ and then extending $\tilde{s}$ to a section $\bar{s}\in C^0([-1,2]\times M\ot\bar{E})$ such that $\bar{s}(t,.)=s_0(.)$ for all $t\in[-1,0]$ and $\bar{s}(t,.)=s_1(.)$ for all $t\in[1,2]$. Then we apply the theorem to $K=\set{-1,2}\times M$ and reparameterise the interval $[-1,2]$ to $[0,1]$.)

\smallskip
The section $s\in C^\infty([0,1]\times M\ot \tilde{E})$ can be identified with a path $s\colon[0,1]\to C^\infty(M\ot E)$, and this path has the desired property.
\end{proof}


\section{Background in obstruction theory} \label{obstructionappendix}

\begin{theorem}[sections of fibrations over CW complexes] \label{CWobstruction}
Let $n\in\N_{\geq1}$, let $(M,N)$ be a relative CW complex, and let $p\colon E\to M$ be a fibration\footnote{in the sense of homotopy theory; cf.\ e.g.\ \cite{Whitehead}, Definition I.7.1} with $(n-1)$-connected\footnote{Recall that a topological space $X$ is $k$-connected if and only if it is path connected and the group $\pi_i(X;x)$ is trivial for all $i\in\set{1,\dots,k}$ and every base point $x\in X$.} fibres.\footnote{I.e., we demand that $p^{-1}(x)$ is $(n-1)$-connected for each $x\in M$. Note that if $p\colon E\to M$ is a fibration and $x,y$ are contained in the same path component of $M$, then the fibres $p^{-1}(x)$ and $p^{-1}(y)$ are homotopy equivalent; cf.\ e.g.\ \cite{Whitehead}, Corollary IV.8.4.} If $n=1$, assume that the fundamental group of each fibre is abelian. Let $f\colon N\to E$ be a continuous map with $p\compose f = \incl_{N,M}$. Then there is a continuous map $\bar{f}\colon M_n\to E$ such that $\bar{f}\restrict N = f$ and $p\compose\bar{f} = \id_{M_n}$; here $M_n$ denotes the $n$-skeleton of the relative CW complex $(M,N)$.
\end{theorem}
\Proof
This is (a special case of) Theorem VI.6.1 in \cite{Whitehead}, for instance: in Whitehead's notation, take $X=E$, $B=K=M$, $L=N$, $\phi=\id_M$. (Whitehead assumes that $M$ is path connected, but that is not necessary since every CW complex is the topological sum of its path components, and the theorem holds obviously for a topological sum if it holds for each component. Whitehead's assumption in the case $n=1$ that the path connected fibres be $1$-simple is equivalent to the abelianness of their fundamental groups; cf.\ e.g.\ \cite{Spanier}, p.~384.)
\end{proof}

This result implies the theorem we will use in Chapter \ref{FIVE}:

\begin{theorem}[sections in smooth bundles] \label{smoothobstruction}
Let $n\in\N_{\geq2}$, let $M$ be an $n$-manifold, let $E\to M$ be a (smooth) fibre bundle with $(n-1)$-connected typical fibre $F$. Let $K$ be a closed subset of $M$, let $U\subseteq M$ be an open neighbourhood of $K$, and let $s_0\in C^\infty(U\ot E)$ be a section over $U$. Then there is a section $s\in C^\infty(M\ot E)$ whose restriction to $K$ is $s_0\restrict K$.
\end{theorem}
\Proof
There is a relative CW complex $(M,N)$ such that $K\subseteq N\subseteq U$. Since every fibre bundle over a manifold is a fibration (cf.\ \cite{Whitehead}, I.7.13), Theorem \ref{CWobstruction} gives us a continuous section $\bar{f}\in C^0(M\ot E)$ with $\bar{f}\restrict N = s_0\restrict N$. By Theorem \ref{smoothapproximation}, there is a $s\in C^\infty(M\ot E)$ with $s\restrict K = \bar{f}\restrict K = s_0\restrict K$.
\end{proof}

\emph{Remark.} It is easy to see that the theorem holds also for $n\leq1$, but we don't need that.


\section{Gromov's h-principle theorems} \label{hprinciple}

The proofs of our results in Chapter \ref{FIVE} rely heavily on the theorem which is known as Gromov's \emph{convex integration} technique or the \emph{h-principle for open ample partial differential relations}. The relevant definitions and the statement of the theorem (in the special case that we need) are contained in Subsection \ref{hprincipleONE} below. General references for the h-principle are \cite{EliashbergMishachev}, \cite{Geigeshprinciple}, \cite{Spring}, and the bible \cite{GromovPDR}.
\medskip\\
In its simplest version, the h-principle deals with the following problem: Given a fibre bundle $E\to M$ and an open subset $\mathscr{R}$ of $J^1E$ (i.e.\ an open \emph{first-order partial differential relation})\footnote{Note that we do not demand that $\mathscr{R}\to M$ is a sub fibre bundle of $J^1E\to M$; local triviality might not hold.}, is there a \emph{solution of $\mathscr{R}$}, i.e., is there a section $\sigma\in C^1(M\ot E)$ such that the image of $j^1\sigma\in C^0(M\ot J^1E)$ is contained in $\mathscr{R}$? (For instance, the existence problem for everywhere twisted distributions that we consider in Chapter 5 has this form.)
\smallskip\\
A necessary condition for the existence of a solution $\sigma\in C^1(M\ot E)$ is obviously the existence of a \emph{formal solution} of $\mathscr{R}$, i.e.\ of a section $\alpha\in C^0(M\ot J^1E)$ whose image is contained in $\mathscr{R}$. It is mainly a problem in algebraic topology to check whether this necessary condition is satisfied. Gromov's h-principle theorems tell us that, under certain conditions on $\mathscr{R}$, the existence of a formal solution implies existence of a solution. In other words, these theorems reduce global problems in \emph{differential} topology to simpler problems in \emph{algebraic} topology.

\subsection{The convex integration method} \label{hprincipleONE}

For an explanation of the name \emph{convex integration} method, let me refer you to the general references cited above which explain the idea of the proof of Theorem \ref{hprincipleone} below. We will just care about its statement (which I cite from \cite{Spring}: cf.\ Theorem 4.2 or the more general Theorem 8.12), and for simplicity, we consider only the case of a first-order relation.
\smallskip\\
First we have to define the notion of \emph{ampleness} of a first-order partial differential relation.

\begin{definition}[$J^1_{\bot W}E$] \label{perp}
Let $n,k\in\N$, let $M$ be an $n$-manifold, let $x\in M$, let $W$ be an $(n-1)$-dimensional sub vector space of $T_xM$, and let $E\to M$ be a fibre bundle with $k$-dimensional fibres.
\smallskip\\
We define $J^1_{\bot W}E$ to be the set of equivalence classes of sections $\sigma\in C^1(M\ot E)$ under the equivalence relation $\sim_x$, where by definition $\sigma_0\sim_x\sigma_1$ holds if and only if $\sigma_0(x)=\sigma_1(x)$ and the restrictions to $W$ of the derivatives $T_x\sigma_0\colon T_xM\to T_{\sigma_0(x)}E$ and $T_x\sigma_1\colon T_xM\to T_{\sigma_1(x)}E$ are equal.
\smallskip\\
Since $J^1_xE$ is defined analogously with $W$ replaced by $T_xM$, there is a canonical projection $p^1_{\bot W} \colon J^1_xE \to J^1_{\bot W}E$ which sends each equivalence class $[\sigma]\in J^1_xE$ to the equivalence class $[\sigma]\in J^1_{\bot W}E$. Moreover, there is a canonical projection $p^{\bot W}_0 \colon J^1_{\bot W}E \to E_x$ which sends each equivalence class $[\sigma]$ to $\sigma(x)$.
\end{definition}

\begin{facts} \label{perpfacts}
In the situation of the previous definition, both projections $p^1_{\bot W}$ and $p^{\bot W}_0$ admit the structure of an affine bundle in a natural way; \cite{Spring}, p.~90. We need only the affine structure for $p^1_{\bot W}$, which can be described as follows.
\smallskip\\
Recall that the bundle $p^1_0 = p^1_{\bot W}\compose p^{\bot W}_0\colon J^1_xE\to E_x$ is an affine bundle modelled on the vector bundle $\Lin(T_xM,T(E_x))\to E_x$, where $T_xM$ denotes the trivial vector bundle $T_xM\times E_x\to E_x$. Via $p^{\bot W}_0$, we can pull back this bundle to a vector bundle $(p^{\bot W}_0)^\ast\Lin(T_xM,T(E_x))$ over $J^1_{\bot W}E$. The affine bundle $p^1_{\bot W}$ is modelled on the sub vector bundle $\xi$ of $(p^{\bot W}_0)^\ast\Lin(T_xM,T(E_x))$ consisting of those $\gamma\in\Lin(T_xM,T(E_x))$ whose restriction to $W$ is zero.
\smallskip\\
In particular, the rank of the affine bundle $p^1_{\bot W}$ is equal to the rank of $T(E_x)$, that is, to the dimension $k$ of the fibres of $E$.
\smallskip\\
Moreover, we obtain in particular the following fact: Let $e\in J^1_{\bot W}E$, let $F\define(p^1_{\bot W})^{-1}(\set{e}) \subseteq J^1_xE$ denote the fibre of $p^1_{\bot W}$ over $e$, and let $b\define p^{\bot W}_0(e)\in E_x$. Then the affine space $F$ is modelled on the vector space $\xi_e = \set{\gamma\in\Lin(T_xM,T_b(E_x)) \suchthat \gamma\restrict W=0}$, and the affine space action $\xi_e\times F \to F$ is just the restriction of the affine space action $\Lin(T_xM,T_b(E_x))\times (p^1_0)^{-1}(\set{b})\to (p^1_0)^{-1}(\set{b})$ induced by the affine bundle structure of $J^1_xE\to E_x$.
\end{facts}

\begin{definition}[ample subset of an affine space] \label{ampledefone}
Let $A$ be an affine space modelled on a finite-dimensional real vector space. An open subset $S$ of $A$ is called \emph{ample} if and only if the convex hull of each connected component of $S$ is $A$.
\end{definition}

\begin{examples} \label{ampleexamples}
Let $A$ be an affine space modelled on a finite-dimensional real vector space. Then $\leer$ is an ample subset of $A$ because it has no connected component. $A$ is an ample subset of $A$. The complement of each affine subspace of $A$ with codimension $\geq2$ is an ample subset of $A$: it is connected, and its convex hull is $A$. The complement of each nonempty affine subspace of $A$ with codimension $1$ is \emph{not} ample: it has precisely two connected components, and the convex hull of each of these is an open half-space in $A$.
\end{examples}

\begin{definition}[ample first-order partial differential relation] \label{ampledeftwo}
Let $E\to M$ be a fibre bundle, and let $\mathscr{R}$ be an open subset of $J^1E$. Then $\mathscr{R}$ is called \emph{ample} if and only if for every $x\in M$, every codimension-$1$ sub vector space $W$ of $T_xM$, and every $e\in J^1_{\bot W}E$, the intersection of $\mathscr{R}$ with the fibre over $e$ of the affine bundle $p^1_{\bot W} \colon J^1_xE \to J^1_{\bot W}E$ is an ample subset of that fibre.
\end{definition}

\emph{Remark.} Ampleness is defined in slightly varying ways in the literature, sometimes explicitly involving the space $J^1_{\bot W}E$, sometimes not. The version above can be found in \cite{Spring}, \S6.1.1. In fact, given an $(n-1)$-plane distribution $W$ on $M$, Spring defines\footnote{using a different notation: $E^{(1)}$ instead of $J^1E$; $E^\bot_x$ instead of $J^1_{\bot W}E$, etc.} a space $J^1_{\bot W}E$ which fits into a sequence of bundle projections $J^1E\to J^1_{\bot W}E\to E\to M$. The fibres over $x$ of these bundles is what we considered in Definition \ref{perp} above --- and it is all we need to understand the statement of Gromov's theorem \ref{hprincipleone}, because ampleness is defined fibrewise, i.e., it does not depend on (say, derivatives of) a distribution $W$ on $M$ but only on a sub vector space of each tangent space $T_xM$. (This fact is obvious and, moreover, explicitly remarked on p.~132 in \cite{Spring}.)

\begin{definition}
We denote by $C^0(M\ot\mathscr{R})$ the set of all $\varphi\in C^0(M\ot J^1E)$ whose image is contained in $\mathscr{R}$. We say that $\varphi_0,\varphi_1$ are contained in the same path component of $C^0(M\ot\mathscr{R})$ if and only if there is a map $\Phi\colon [0,1]\to C^0(M\ot\mathscr{R})$ with $\Phi(0)=\varphi_0$ and $\Phi(1)=\varphi_1$ which is continuous with respect to the subspace topology on $C^0(M\ot\mathscr{R})$ induced by the compact-open $C^0$-topology on $C^0(M\ot J^1E)$.
\end{definition}

\begin{theorem}[Gromov's h-principle for ample relations; first-order $C^0$-dense non-relative non-parametric version] \label{hprincipleone}
Let $E\to M$ be a fibre bundle, let $\mathscr{R}\subseteq J^1E$ be open and ample, and let $\varphi\in C^0(M\ot\mathscr{R})$ (i.e., $\varphi$ is a formal solution of $\mathscr{R}$). Let $p^1_0\colon J^1E\to E$ denote the standard projection, let $h\define p^1_0\compose\varphi\in C^0(M\ot E)$, and let $\mathscr{N}$ be a neighbourhood of the image of $h$ in $E$. Then there is a section $\sigma\in C^\infty(M\ot E)$ such that $j^1\sigma\in C^0(M\ot\mathscr{R})$ (i.e., $\sigma$ is a solution of $\mathscr{R}$), such that the image of $\sigma$ is contained in $\mathscr{N}$ (i.e., $\sigma$ is close to $h$ with respect to the fine $C^0$-topology), and such that, moreover, $j^1\sigma$ is contained in the same path component of $C^0(M\ot\mathscr{R})$ as $\varphi$; the latter assertion implies in particular that $\sigma$ is contained in the same connected component of $C^0(M\ot E)$ as $h$.
\end{theorem}

\emph{Remarks.} The last claim of Theorem \ref{hprincipleone} follows from the preceding one because if a continuous map $\Phi\colon [0,1]\to C^0(M\ot J^1E)$ with $\Phi(0)=\varphi_0$ and $\Phi(1)=\varphi_1$ exists, then there is a continuous map $\Psi\colon [0,1]\to C^0(M\ot E)$ with $\Psi(0)=p^1_0\compose\varphi_0$ and $\Psi(1)=p^1_0\compose\varphi_1$: namely, composition with $p^1_0$ from the left defines a continuous (with respect to the compact-open topologies) map $\pi\colon C^0(M\ot J^1E)\to C^0(M\ot E)$, so the map $\Psi\define \pi\compose\Phi$ has the required property.
\smallskip\\
There is also a relative version of Theorem \ref{hprincipleone} (cf.\ Theorem 4.2 in \cite{Spring}; our Theorem \ref{hprincipleone} is the special case $K_0=\leer$). We do not state it here because we will not care about its consequences for the distribution version of the prescribed scalar curvature problem.

\subsection{The h-principle for diff-invariant relations} \label{openhprinciple}

Theorem \ref{hprincipletwo} below (usually called the \emph{h-principle for open diff-invariant relations} or the \emph{covering homotopy method}) is not applied in the present thesis, but we mention it occasionally. The \emph{diff-invariance} condition on a partial differential relation means roughly that the relation is invariant under a natural action of the diffeomorphism group. For a precise definition, let me refer you to \cite{Geigeshprinciple} (Definition 3.2) or \cite{EliashbergMishachev}, \S7.

\begin{theorem}[Gromov] \label{hprincipletwo}
Let $M$ be an open manifold, let $E\to M$ be a fibre bundle, let $r\in\N$, and let $\mathscr{R}$ be an open diff-invariant subset of $J^rE$. Then for every connected component $\mathscr{C}$ of $C^0(M\ot\mathscr{R})$, there is a section $\sigma\in C^r(M\ot E)$ with $j^r\sigma\in\mathscr{C}$.
\smallskip\\
In particular, for every $\varphi\in C^0(M\ot\mathscr{R})$ (i.e., for every \emph{formal solution} of $\mathscr{R}$), there is a section $\sigma\in C^r(M\ot E)$ with $j^r\sigma\in C^0(M\ot\mathscr{R})$ (i.e.\ a \emph{solution} of $\mathscr{R}$) such that $\sigma$ is contained in the same connected component of $C^r(M\ot E)$ as $p^r_0\compose\varphi$ (where $p^r_0\colon J^rE\to E$ denotes the standard projection).
\end{theorem}

So the difference to Theorem \ref{hprincipleone} is that ampleness is replaced by diff-invariance, and that $M$ must be open (and that we stated in \ref{hprincipleone} only the first-order version of the convex integration method).

\begin{remark}
Theorem \ref{hprincipletwo} contains only a weak form of Gromov's h-principle theorem for open diff-invariant relations. There is a ``parametric'' version of the theorem which says that
\[
j^r\colon \set{\sigma\in C^r(M\ot E)\suchthat j^r\sigma\in C^0(M\ot\mathscr{R})} \;\to\; C^0(M\ot\mathscr{R})
\]
is a homotopy equivalence; our version above is just the statement that this map induces a surjective map between the sets of connected components.
\smallskip\\
Cf.\ Chapter 3 in \cite{Geigeshprinciple} for a proof of that parametric h-principle (in particular, for a proof of Theorem \ref{hprincipletwo}), and cf.\ Theorem 7.2.4 in \cite{EliashbergMishachev} for a relative version. Gromov's original article \cite{Gromov1969english} contains the non-relative parametric version. Notice that there is in general no (global) $C^0$-dense version of this h-principle; cf.\ Chapter 7 in \cite{EliashbergMishachev} for further information.
\end{remark}


\section{Contact structures and even-contact structures} \label{CONTACT}

Since our discussion of everywhere twisted distributions (cf.\ Definition \ref{twistednessdef}) in Chapter \ref{FIVE} is related to contact and even-contact structures (cf.\ the remarks at the beginning of that chapter), we review here the definitions and some basic facts which explain the connection.

\subsection{Contact structures} \label{CONTACTONE}

We will discuss contact structures in the general --- that is, not necessarily coorientable --- form:

\begin{definition}
Let $\xi$ be a $2n$-plane distribution on a $(2n+1)$-manifold $M$. A $1$-form $\alpha\in\Omega^1(U)$ defined on some open subset $U$ of $M$ is a \emph{local contact form for $\xi$} if and only if $\xi\restrict U = \ker(\alpha)$ and the $(2n+1)$-form $\alpha\wedge(d\alpha)^n\in\Omega^{2n+1}(U)$ vanishes nowhere (the notation $^n$ refers to the wedge product). The distribution $\xi$ is a \emph{contact structure} if and only if for every $x\in M$ there is an open neighbourhood $U$ of $x$ and a local contact form $\alpha\in\Omega^1(U)$ for $\xi$.
\end{definition}

\begin{facts} \label{contactfacts}
Let $M$ be a $(2n+1)$-manifold.
\begin{enumerate}
\item
If $\alpha,\alpha'\in\Omega^1(U)$ are local contact forms for the contact structure $\xi$ on $M$, then there is a nowhere vanishing function $f\in C^\infty(U,\R)$ such that $\alpha' = f\alpha$; this follows from the condition $\ker(\alpha)=\ker(\alpha')$. Conversely, if $\alpha\in\Omega^1(U)$ is a local contact form for the contact structure $\xi$ and $f\in C^\infty(U,\R)$ vanishes nowhere, then $f\alpha$ is a local contact form for $\xi$; that's because $\ker(f\alpha) = \ker(\alpha)$ and $f\alpha\wedge(d(f\alpha))^n = f\alpha\wedge(f\,d\alpha+df\wedge\alpha)^n = f^{n+1}\alpha\wedge(d\alpha)^n$ (since $\alpha\wedge\alpha=0$).
\item
\emph{Existence of contact structures.} The $1$-form $\beta = dz +\sum_{j=1}^nx_jdy_j$ on $\R^{2n+1}$ ($(x_1,\dots,x_n,y_1,\dots,y_n,z)$ are the standard coordinates on $\R^{2n+1}$) is a contact form for the contact structure $\ker(\beta)$.
\smallskip\\
\emph{Local uniqueness of contact structures.} If, for $i\in\set{0,1}$, we have a $(2n+1)$-manifold $M_i$, a point $x_i\in M_i$, and a contact form $\alpha_i$ for the contact structure $\ker(\alpha_i)$ on $M_i$, then there are open neighbourhoods $U_0\subseteq M_0$ of $x_0$ and $U_1\subseteq M_1$ of $x_1$ and a diffeomorphism $\varphi\colon U_0\to U_1$ such that $\varphi(x_0)=x_1$ and $\alpha_0\restrict U_0 = \varphi^\ast(\alpha_1\restrict U_1)$. This is Darboux's theorem for contact structures; cf.\ e.g.\ \cite{Geigessurvey}, Theorem 2.24.
\item
If the contact structure $\xi$ on $M$ is \emph{coorientable} (synonymously: \emph{transversely orientable}), i.e.\ the line bundle $TM/\xi$ over $M$ is orientable and thus trivial, then there is a \emph{global} contact form for $\xi$, i.e.\ a local contact form for $\xi$ which is defined on the whole of $M$. To see this, we choose a vector field $X$ on $M$ which is everywhere transverse to $\xi$, and consider local contact forms $\alpha_i\in\Omega^1(U_i)$ (where $I$ is an index set) such that $(U_i)_{i\in I}$ is a covering of $M$. For each $i\in I$, we define a nowhere vanishing function $\lambda_i\in C^\infty(U_i,\R)$ by $\lambda_i\alpha_i(X)=1$. Then $\lambda_i\alpha_i = \lambda_j\alpha_j$ for all $i,j\in I$ since both $1$-forms have the same value on $X$ and vanish on $\xi$. Thus there is a $1$-form $\alpha\in\Omega^1(M)$ whose restriction to each $U_i$ is $\lambda_i\alpha_i$. It follows from the facts in paragraph (i) that $\alpha$ is a contact form for $\xi$.
\smallskip\\
Contact structures are sometimes defined to be kernels of \emph{global} contact forms. That definition yields precisely the coorientable contact structures in our general sense.
\item
If $n$ is even, then every contact structure on $M$ is orientable (as a vector bundle). Namely, every local contact form $\alpha\in\Omega^1(U)$ induces an orientation on the vector bundle $\xi\restrict U$ via the top-rank form $(d\alpha)^n\restrict(\xi\restrict U)$; this form vanishes nowhere since $\alpha\wedge(d\alpha)^n$ vanishes nowhere and $\xi=\ker(\alpha)$. Every two local contact forms $\alpha,\alpha'\in\Omega^1(U)$ induce the same orientation on $\xi\restrict U$ because there is a nowhere vanishing function $f\in C^\infty(U,\R)$ with $\alpha' = f\alpha$, and thus $(d\alpha')^n\restrict\xi = (fd\alpha +df\wedge\alpha)^n\restrict\xi = f^n(d\alpha)^n\restrict\xi$ (where we used $\alpha\restrict\xi = 0$); since $n$ is even, $f^n$ is everywhere positive.
\item
If $n$ is odd, then the existence of a contact structure $\xi$ on $M$ implies that $M$ is orientable. Namely, every local contact form $\alpha\in\Omega^1(U)$ induces an orientation on $U$ via the $(2n+1)$-form $\alpha\wedge(d\alpha)^n$. Every two local contact forms $\alpha,\alpha'\in\Omega^1(U)$ induce the same orientation on $U$ because there is a nowhere vanishing function $f\in C^\infty(U,\R)$ with $\alpha'\wedge(d\alpha')^n = f^{n+1}\alpha\wedge(d\alpha)^n$; since $n$ is odd, $f^{n+1}$ is everywhere positive.
\item
Let us summarise the results of (iv) and (v): If $\xi$ is a contact structure on $M$, then let $L$ be the line bundle $TM/\xi$ on $M$. If
\begin{itemize}
\item
$n$ is odd: then $M$ is orientable; $L$ is orientable if and only if $\xi$ is orientable.
\item
$n$ is even: then $\xi$ is orientable; $M$ is orientable if and only if $L$ is orientable.
\end{itemize}
\item
For every $n\in\N_{\geq1}$, there is a $(2n+1)$-manifold $M$ which admits a non-coorientable contact structure; take e.g.\ $M = \R^{n+1}\times\RP^n$: The $2n$-plane distribution $\ker(\sum_{j=0}^ny_jdx_j)$, where $(x_0,\dots,x_n)$ are the standard coordinates on $\R^{n+1}$ and $[y_0:\ldots:y_n]$ are the standard homogeneous coordinates on $\RP^n$, is a non-coorientable contact structure; cf.\ \cite{Geigessurvey}, Example 2.14 and Proposition 2.15.
\end{enumerate}
\end{facts}

If $M$ is an odd-dimensional manifold of dimension $\geq5$, then it is in general not easy to decide whether $M$ admits a contact structure; cf.\ H.\ Geiges' review \cite{Geiges2001} for results in this direction. As we have just seen, if a $3$-manifold admits a contact structure, then it is orientable. A theorem by J.\ Martinet (cf.\ \cite{Martinet1971}) tells us that the converse is also true. In fact, every homotopy class of $2$-plane distributions on an orientable $3$-manifold contains a contact structure. This follows in the open case from M.\ Gromov's h-principle theorems (from the covering homotopy method, to be precise). In the closed case, it has been proved by R.\ Lutz for homotopy classes of \emph{(co-)orientable} $2$-plane distributions (cf.\ \cite{Lutz1977}; \cite{Geigessurvey}, Section 3), and by Y.\ Eliashberg in the most general case (cf.\ \cite{Eliashberg1989}).

\begin{theorem}[existence of contact structures on orientable $3$-manifolds] \label{threecontact}
Let $M$ be an orientable $3$-manifold (with or without boundary). Then every homotopy class of $2$-plane distributions on $M$ contains a contact structure.
\end{theorem}
\begin{proof}[\textsc{Remarks on the proof.}]
For closed $M$ (this is the hard case), that is proved in \cite{Eliashberg1989}. In the open connected case, we can apply Gromov's h-principle Theorem \ref{hprincipletwo} since the contact condition, viewed as a partial differential relation on the total space $J^1\Gr_2(TM)$ of the $1$-jet bundle of the Grassmann bundle $\Gr_2(TM)\to M$, is obviously open and diff-invariant, and suitable formal solutions exist (as we reprove in Proposition \ref{formalsolution} and Remark \ref{threecontactremark}).
\end{proof}

\subsection{Even-contact structures} \label{CONTACTTWO}

\begin{definition}[even-contact structure]
Let $\xi$ be a $(2n+1)$-plane distribution on a $(2n+2)$-manifold $M$. A $1$-form $\alpha\in\Omega^1(U)$ defined on some open subset $U$ of $M$ is a \emph{local even-contact form for $\xi$} if and only if $\xi\restrict U = \ker(\alpha)$ and the $(2n+1)$-form $\alpha\wedge(d\alpha)^n\in\Omega^{2n+1}(U)$ vanishes nowhere. The distribution $\xi$ is an \emph{even-contact structure} if and only if for every $x\in M$ there is an open neighbourhood $U$ of $x$ and a local even-contact form $\alpha\in\Omega^1(U)$ for $\xi$. \footnote{The name \emph{even-contact form/structure} was introduced by V.\ L.\ Ginzburg, who investigated these objects in \cite{Ginzburg1992}, in 1992. In D.\ McDuff's article \cite{McDuff1987} (Example 2.6, Lemma 2.7, and \S7) from 1987, in which she proved the ampleness of the corresponding partial differential relation, even-contact forms are called \emph{non-degenerate $1$-forms}. It seems that not much work has been done on even-contact structures, in contrast to contact structures.}
\end{definition}

We collect some elementary facts about even-contact structures, analogous to the facts in \ref{contactfacts}.

\begin{facts} \label{evencontactfacts}
Let $M$ be a $(2n+2)$-manifold.
\begin{enumerate}
\item
If $\alpha,\alpha'\in\Omega^1(U)$ are local even-contact forms for the even-contact structure $\xi$ on $M$, then there is a nowhere vanishing function $f\in C^\infty(U,\R)$ such that $\alpha' = f\alpha$. Conversely, if $\alpha\in\Omega^1(U)$ is a local contact form for the contact structure $\xi$ and $f\in C^\infty(U,\R)$ vanishes nowhere, then $f\alpha$ is a local contact form for $\xi$.
\item
\emph{Existence.} Let $N$ be a $(2n+1)$-manifold, let $B\in\set{\R,S^1}$, let $\pr_N\colon N\times B\to N$ and $\pr_B\colon N\times B\to B$ denote the obvious projections. If $\xi$ is a contact structure on $N$, then the vector bundle $\eta\define (\pr_N^\ast\xi)\oplus(\pr_B^\ast TB)$ --- which is a sub vector bundle of $(\pr_N^\ast TM)\oplus(\pr_B^\ast TB) = T(N\times B)$ --- is an even-contact structure on $N\times B$.
\smallskip\\
Namely, by Darboux' theorem, $\xi$ is locally the kernel of a $1$-form $\alpha\in\Omega^1(U)$ (where $U$ is an open subset of $N$) which has in suitable local coordinates $(x_1,\dots,x_n,y_1,\dots,y_n,z_1)$ the form $dz_1 +\sum_{j=1}^nx_jdy_j$. With respect to local coordinates on $N\times B$ of the form $(x_1,\dots,x_n,y_1,\dots,y_n,z_1,z_2)$, the $1$-form $\pr_N^\ast\alpha$ is given by $dz_1 +\sum_{j=1}^nx_jdy_j$, which is an even-contact form. Its kernel is $\pr_N^\ast(\ker(\alpha))\oplus(\pr_B^\ast TB)$, i.e.\ $\eta$. Thus $\eta$ is an even-contact structure, as claimed.
\medskip\\
In particular, even-contact structures exist in all even dimensions.
\item
\emph{Local uniqueness.} If $\alpha$ is an even-contact form on $M$ and $x$ is a point in $M$, then there exist local coordinates $(x_1,\dots,x_n,y_1,\dots,y_n,z_1,z_2)$ on a neighbourhood $U$ of $x$ such that $\alpha\restrict U = dz_1 +\sum_{j=1}^nx_jdy_j$; cf.\ \cite{McDuff1987}, Proposition 7.2 and the remark preceding it.
\item
If an even-contact structure $\xi$ is coorientable, then there exists a \emph{global} even-contact form for $\xi$. The proof is the same as for the analogous statement for contact structures.
\item
If $n$ is even, then every even-contact structure $\xi$ on $M$ admits a nowhere vanishing $2n$-form $\beta\in C^\infty(M\ot\Lambda^{2n}(\xi^\ast))$.
\smallskip\\
Namely, there is an open cover $\mathscr{U}$ of $M$ such that every $U\in\mathscr{U}$ admits an even-contact form $\alpha_U$ for $\xi\restrict U$. Via a subordinate partition of unity $(\varphi_U)_{U\in\mathscr{U}}$ we define $\beta$ to be the global form $\sum_{U\in\mathscr{U}}\varphi_U\cdot(d\alpha_U)^n\restrict(\xi\restrict U)$. This $\beta$ vanishes indeed nowhere: For all $U,V\in\mathscr{U}$, there is a nowhere vanishing function $f_{UV}\in C^\infty(U\cap V,\R)$ with $\alpha_V\restrict(U\cap V) = f_{UV}\alpha_U\restrict(U\cap V)$; hence $(d\alpha_U)^n\restrict\xi = f_{UV}^n(d\alpha_V)^n\restrict\xi$ holds on $U\cap V$. For each $U\in\mathscr{U}$, we obtain
\[
\beta\restrict U = \sum_{V\in\mathscr{U}}\varphi_Vf_{UV}^n\cdot(d\alpha_U)^n\restrict\xi \;\;.
\]
This form vanishes nowhere since $\sum_{V\in\mathscr{U}}\varphi_Vf_{UV}^n$ is everywhere positive and $(d\alpha_U)^n\restrict\xi$ vanishes nowhere. (In order to verify the latter fact, choose locally a $\gamma\in\Omega^1(M)$ such that the top-rank form $\alpha\wedge\gamma\wedge(d\alpha_U)^n$ vanishes nowhere; that is possible e.g.\ by (iii): take $\gamma=dz_2$. Then extend a local frame $(v_1,\dots,v_{2n+1})$ of $\xi$ to a local frame of $TM$, and evaluate $\alpha\wedge\gamma\wedge(d\alpha_U)^n$ on this frame to obtain the statement.) Hence $\beta$ vanishes nowhere, as claimed.
\item
If $n$ is odd and $\xi$ is an even-contact structure on $M$, then $M$ admits a nowhere vanishing $(2n+1)$-form whose restriction to $\xi$ is zero. This is proved similarly as in (iv): we patch the local forms $\alpha_U\wedge(d\alpha_U)^n$ together via a partition of unity.
\end{enumerate}
\end{facts}

\subsection{Everywhere twistedness vs.\ contact and even-contact} \label{CONTACTTHREE}

\begin{proposition} \label{twistvswedge}
Let $n\in\N_{\geq1}$, let $M$ be an $n$-manifold, let $H$ be an $(n-1)$-plane distribution on $M$, and let $x\in M$. Then the following statements are equivalent:
\begin{enumerate}
\item
The twistedness of $H$ vanishes in $x$.
\item
There is an open neighbourhood $U$ of $x$ and a $1$-Form $\alpha\in\Omega^1(U)$ with $H\restrict U = \ker(\alpha)$ such that $\alpha\wedge d\alpha\in\Omega^3(U)$ vanishes in $x$.
\item
For every open neighbourhood $U$ of $x$ and every $1$-form $\alpha\in\Omega^1(U)$ with $H\restrict U = \ker(\alpha)$, the $3$-form $\alpha\wedge d\alpha$ vanishes in $x$.
\end{enumerate}
\end{proposition}
\Proof
(iii)$\implies$(ii): A sufficiently small open neighbourhood $U$ of $x$ admits a vector field $X$ which is transverse to $H$. Define $\alpha\in\Omega^1(U)$ by $\alpha(X)=1$ and $\alpha\restrict H = 0$. Then $H\restrict U = \ker(\alpha)$ and so, by statement (iii), $\alpha\wedge d\alpha$ vanishes in $x$.
\smallskip\\
(ii)$\implies$(i): We have to show that for all local sections $v,w$ in $H$, the value in $x$ of the Lie bracket $[v,w]$ lies in $H_x$; in other words, we must prove that $\alpha([v,w])$ vanishes in $x$ whenever $v,w$ are vector fields on $U$ with $\alpha(v)=\alpha(w)=0$. This is easy: From $d\alpha(v,w) = \partial_v(\alpha(w)) -\partial_w(\alpha(v)) -\alpha([v,w]) = -\alpha([v,w])$, we deduce
\[
(\alpha\wedge d\alpha)(v,w,[v,w]) = \alpha([v,w])d\alpha(v,w) = -\alpha([v,w])^2 \;\;.
\]
Since statement (ii) tells us that the leftmost function vanishes in $x$, so does $\alpha([v,w])$.
\smallskip\\
(i)$\implies$(iii): Let $U$ be an open neighbourhood of $x$, and let $\alpha\in\Omega^1(U)$ with $H\restrict U = \ker(\alpha)$. There is a vector field $X$ on $U$ such that $\alpha(X)=1$. In order to prove that a $3$-form $\beta$ on $U$ vanishes in $x$, it suffices to show that the function $\beta(v,w,X)$ vanishes in $x$ for all sections $v,w$ in $H\restrict U$.
\smallskip\\
If $v,w$ are sections in $H\restrict U$, then $\alpha(v)=\alpha(w)=0$ and, by statement (i), the function $\alpha([v,w])$ vanishes in $x$. As in the proof of (ii)$\implies$(i), we compute $(\alpha\wedge d\alpha)(v,w,X) = \alpha(X)d\alpha(v,w) = -\alpha([v,w])$. Therefore the $3$-form $\alpha\wedge d\alpha$ vanishes in $x$.
\end{proof}

\begin{corollary}
Contact structures on manifolds of dimension $\geq3$ are everywhere twisted. A $2$-plane distribution on a $3$-manifold is everywhere twisted if and only if it is contact. Even-contact structures on manifolds of dimension $\geq4$ are everywhere twisted. A $3$-plane distribution on a $4$-manifold is everywhere twisted if and only if it is even-contact.
\end{corollary}
\Proof
By the equivalence (i)$\iff$(iii) from the preceding proposition, an $(m-1)$-plane distribution on an $m$-manifold is everywhere twisted if and only if it is locally the kernel of a $1$-form $\alpha$ with nowhere vanishing $\alpha\wedge d\alpha$. All contact structures manifolds of dimension $2n+1\geq3$, as well as all even-contact structures on manifolds of dimension $2n+2\geq4$, have this property by definition. If $n=1$, then the contact and even-contact conditions are by definition equivalent to this property.
\end{proof}


\section{Bundles of Grassmannians} \label{Grassmannappendix}

This section reviews a few definitions and basic facts about Grassmann manifolds (also known as \emph{Grassmann varieties} or simply \emph{Grassmannians}) and bundles of Grassmann manifolds; they occur in Chapter \ref{FIVE} and Appendix \ref{AppendixC}.
\smallskip\\
The Grassmann manifold of $q$-dimensional sub vector spaces of $\R^n$ is often denoted by $G(q,n)$. Our notation here is $\Gr_q(E)$ for the $q$th Grassmannian on an arbitrary $n$-dimensional real vector space $E$. (This $\Gr_q$ can be regarded as a functor from the category of finite-dimensional vector spaces and vector space isomorphisms to the category of manifolds and diffeomorphisms, or from the category of vector bundles and fibrewise bijective vector bundle morphisms to the category of fibre bundles and fibrewise diffeomorphic bundle maps, but we don't need to adopt this functorial viewpoint.)
\medskip\\
Since the affine space in the following definition does not seem to have a standard name, I introduce the notation $\Compl(H)$ for it. Its affine structure plays an important role at several places in the thesis, so we verify the affine space axioms in detail.

\begin{definition}[the affine space $\Compl(H)$] \label{affdef}
Let $E$ be an $\R$-vector space, and let $H$ be a sub vector space of $E$ with finite codimension. We denote the set of all sub vector spaces of $E$ which are complementary to $H$ by $\Compl(H)$. We turn this set into an affine space modelled on the vector space $\Lin(E/H,H)$ by the following definition. Let $\pi\colon E\to E/H$ be the canonical projection. $\Lin(E/H,H)$ acts on $\Compl(H)$ via the map $+\colon \Lin(E/H,H)\times\Compl(H) \to \Compl(H)$ given by $(\lambda,V) \mapsto \lambda+V$, where $\lambda+V = \set{v+\lambda(\pi(v)) \suchthat v\in V}$.
\end{definition}
\begin{proof}[\textsc{Verification of the affine space axioms.}]
First we check that the map $+$ is well-defined, i.e.\ $\lambda+V\in\Compl(H)$ if $\lambda\in\Lin(E/H,H)$ and $V\in\Compl(H)$. To see this, note that $\incl_{V,E} +\lambda\compose\pi\restrict V \in \Lin(V,E)$ is injective: its composition with the projection $E = V\oplus H \to V$ is the identity on $V$. Hence $\dim(\lambda+V)=\dim(V)=\dim(E/H)$. Moreover, $(\lambda+V)\cap H = \set{0}$ since if $v+\lambda(\pi(v))\in H$ then $v=0$ and thus $v+\lambda(\pi(v))=0$. Because $E/H$ is finite-dimensional, $V$ and $H$ are complementary. This shows that $+$ is well-defined.
\smallskip\\
Obviously $(\lambda_0+\lambda_1)+V = \lambda_0+(\lambda_1+V)$ and $0_{\Lin(E/H,H)}+V = V$ hold for all $\lambda_0,\lambda_1\in\Lin(E/H,H)$ and $V\in\Compl(H)$, so $+$ is indeed an action of the vector space $\Lin(E/H,H)$ on the set $\Compl(H)$.
\smallskip\\
If $\lambda+V=V$ for some $\lambda\in\Lin(E/H,H)$ and $V\in\Compl(H)$, then for each $v\in V$, we have $\lambda(\pi(v))\in V$, hence $\lambda(\pi(v))=0$. This implies $\lambda=0$ since $\pi\restrict V\colon V\to E/H$ is surjective. Thus $+$ is a free action.
\smallskip\\
Let $V_0,V_1\in\Compl(H)$. We want to show that there is a $\lambda\in\Lin(E/H,H)$ with $\lambda+V_0 = V_1$. For $i\in\set{0,1}$, let $p_i\colon E = V_i\oplus H \to H$ denote the projection onto the second factor. Since $p_0-p_1\in\Lin(E,H)$ maps every element of $H$ to $0$, there is a $\lambda\in\Lin(E/H,H)$ such that $p_0-p_1 = \lambda\compose\pi$. We get \mbox{$\lambda+V = \set{v+\lambda(\pi(v)) \suchthat v\in V_0} = \set{v+p_0(v)-p_1(v) \suchthat v\in V_0} = \set{v-p_1(v) \suchthat v\in V_0}$}. This is a subset of $V_1$ by definition of $p_1$. On the other hand, every $v_1\in V_1$ has a unique decomposition $v_1 = v+p_1(v_1)$ with $v\in V_0$, and clearly $p_1(v_1) = -p_0(v)$ holds; thus $v_1\in \lambda+V$. To summarise, we have $\lambda+V_0 = V_1$. This proves that $+$ is a transitive action and completes our verification of the axioms.
\end{proof}

\begin{figure}[bth]
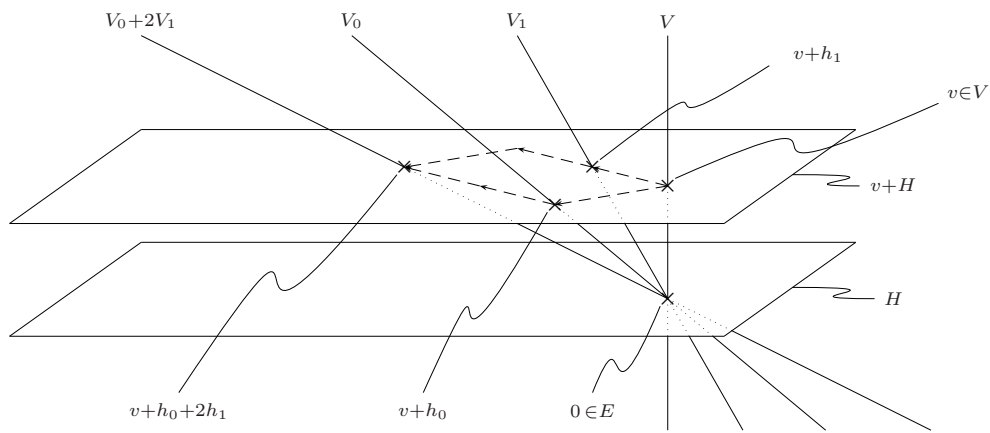
 \label{Complpicture}
\setlength{\unitlength}{1cm}
\vspace*{3.5cm}
\[
\psset{unit=0.5cm,nodesep=2pt,linewidth=0.3pt,origin={-5,0},dotsep=2pt}
\everypsbox{\scriptstyle}
\psline(1.5,-1)(5,1.5)
\psline(5,1.5)(-14,1.5)
\psline(-14,1.5)(-17.5,-1)
\psline(-17.5,-1)(1.5,-1)
\psline(1.5,2)(5,4.5)
\psline(5,4.5)(-14,4.5)
\psline(-14,4.5)(-17.5,2)
\psline(-17.5,2)(1.5,2)
\psline[linestyle=dashed,arrows=->](0,3)(-2,3.5)
\psline[linestyle=dashed,arrows=->](-2,3.5)(-4,4)
\psline[linestyle=dashed,arrows=->](0,3)(-3,2.5)
\psline[linestyle=dashed,arrows=->](-3,2.5)(-5,3)
\psline[linestyle=dashed,arrows=->](-5,3)(-7,3.5)
\psline[linestyle=dashed,arrows=->](-4,4)(-7,3.5)
\psline(0,-3.5)(0,-1)
\psline[linestyle=dotted](0,-1)(0,0)
\psline(0,0)(0,2)
\psline[linestyle=dotted](0,2)(0,3)
\psline(0,3)(0,7)
\psline(2,-3.5)(0.57143,-1)
\psline[linestyle=dotted](0.57143,-1)(0,0)
\psline(0,0)(-1.14286,2)
\psline[linestyle=dotted](-1.14286,2)(-2,3.5)
\psline(-2,3.5)(-4,7)
\psline(4.2,-3.5)(1.2,-1)
\psline[linestyle=dotted](1.2,-1)(0,0)
\psline(0,0)(-2.4,2)
\psline[linestyle=dotted](-2.4,2)(-3,2.5)
\psline(-3,2.5)(-8.4,7)
\psline(7,-3.5)(1.70588,-0.85294)
\psline[linestyle=dotted](1.70588,-0.85294)(0,0)
\psline(0,0)(-4,2)
\psline[linestyle=dotted](-4,2)(-7,3.5)
\psline(-7,3.5)(-14,7)
\psset{dotstyle=x,dotsize=5pt} \psdot(0,0) 
\psdot(0,3) 
\psdot(-2,3.5) 
\psdot(-3,2.5) 
\psdot(-7,3.5)
\rput[b](5,7.2){\Rnode{V}{V}}
\rput[b](1,7.2){\Rnode{V1}{V_1}}
\rput[b](-3.4,7.2){\Rnode{V0}{V_0}}
\rput[b](-9,7.2){\Rnode{V0plus2V1}{V_0+2V_1}}
\rput(11,0){\Rnode{H}{H}}
\rput(11,3){\Rnode{H}{v+H}}
\rput(13,5.5){\Rnode{v}{v\in V}}
\rput(9,6.5){\Rnode{h1}{v+h_1}}
\rput(-1.5,-3){\Rnode{h0}{v+h_0}}
\rput(-8,-3){\Rnode{resultv}{v+h_0+2h_1}}
\rput(3,-3){\Rnode{zero}{0\,\in E}}
\psbezier[linewidth=0.1pt](-0.2,-0.2)(-1.8,-4)(-1,0)(-2,-2.5)
\psbezier[linewidth=0.1pt](-3.2,2.3)(-6.8,-4)(-4,3)(-6.5,-2.5)
\psbezier[linewidth=0.1pt](-7.2,3.3)(-13,-4)(-8,5)(-13,-2.5)
\psbezier[linewidth=0.1pt](3.32,0.3)(7,0.3)(2.5,0)(5.5,0)
\psbezier[linewidth=0.1pt](3.32,3.3)(7,3.3)(2.5,3)(5.1,3)
\psbezier[linewidth=0.1pt](2.7,6.2)(-1.5,3.5)(2.5,7)(-1.8,3.7)
\psbezier[linewidth=0.1pt](0.2,3.2)(6.5,6)(0,2)(7.2,5.2)
\]
\vspace*{1.3cm}
\caption{The vector space structure of $\Compl_V(H)$ visualised (in the case $\dim E= 3$, $\dim V = 1$).}
\end{figure}

\begin{remark}
If we fix, in the situation of the preceding definition, a subspace $V\in\Compl(H)$, then this choice turns the affine space $\Compl(H)$ into a vector space with zero element $V$, since $V$ defines an affine isomorphism $\Lin(E/H,H)\to\Compl(H)$ by $\lambda\mapsto\lambda+V$. We denote this vector space by $\Compl_V(H)$ (cf.\ Figure \ref{Complpicture}). Note that the isomorphism $E/H\cong V$ induces a canonical vector space isomorphism $\Compl_V(H)\cong\Lin(V,H)$.
\end{remark}

\begin{definition}[Grassmann manifolds] \label{Grassmannatlas}
Let $n,q\in\N$, and let $E$ be an $n$-dimensional $\R$-vector space. We equip the set $\Gr_q(E)$ of all $q$-dimensional sub vector spaces of $E$ with the structure of a real-analytic manifold as follows. (In fact, we even equip it with the structure of an algebraic variety, but that's not important for us.) If $q>n$, then $\Gr_q(E)=\leer$ and we are done; so assume $q\leq n$. Let $\mathscr{A}$ be the set of all pairs $(U,\varphi)$, where $U=\Compl(H)$ for some $(n-q)$-dimensional subspace $H$ of $E$ and $\varphi\colon\Compl(H)\to\R^{q(n-q)}$ is an affine isomorphism (note that $\dim\Compl(H) = \dim\Lin(E/H,H) = q(n-q)$). This $\mathscr{A}$ is a real-analytic atlas on the set $\Gr_q(E)$. The set $\Gr_q(E)$ together with this real-analytic structure is called the \emph{$q$-th Grassmann manifold} (synonymously:\ \emph{Grassmannian}) \emph{of $E$}.
\end{definition}

\begin{remark} \label{grassmanntangent}
The \emph{universal vector bundle $\Uni_q(E)$ over $\Gr_q(E)$} is the sub vector bundle of the trivial vector bundle $\Gr_q(E)\times E\to\Gr_q(E)$ whose fibre over $V\in\Gr_q(E)$ is $\set{V}\times V \subseteq \set{V}\times E$. We define the vector bundle $\bot\Uni_q(E)$ over $\Gr_q(E)$ to be the quotient of the trivial bundle $\Gr_q(E)\times E\to \Gr_q(E)$ and its sub vector bundle $\Uni_q(E)$.
\smallskip\\
Then the tangent bundle of the manifold $\Gr_q(E)$ is canonically isomorphic to $\Lin(\Uni_q(E),\bot\Uni_q(E))$. In particular, for each $V\in\Gr_q(E)$, the tangent space $T_V\Gr_q(E)$ is canonically isomorphic to $\Lin(V,E/V)$.
\end{remark}

\begin{definition}[$\Compl(H)$ on the vector bundle level]
Let $\pi\colon E\to M$ be a smooth [resp.\ real-analytic] rank-$n$ vector bundle, and let $H$ be a smooth [real-analytic] rank-$q$ sub vector bundle of $E$. Then $\Lin(E/H,H)$ is a vector bundle. Let $\Pi \colon \Compl(H) \to M$ be the set-theoretic bundle whose fibre over $x\in M$ is the affine space $\Compl(H_x)$ modelled on the vector space $\Lin(E_x/H_x,H_x)$. We turn this set-theoretic bundle into a smooth [real-analytic] affine bundle modelled on the vector bundle $\Lin(E/H,H)$, as follows.
\smallskip\\
Since $H$ is a smooth [real-analytic] sub vector bundle of $E$, every point in $M$ has an open neighbourhood $U$ which admits smooth [real-analytic] vector bundle trivialisations $\phi\colon \pi^{-1}(U)\to U\times\R^n$ and $\phi_H\colon \pi^{-1}(U)\cap H \to U\times\R^q$ such that the diagram
\[
\setlength{\arraycolsep}{1cm}
\begin{array}{cc}
\Rnode{p00}{\pi^{-1}(U)\cap H}
& \Rnode{p01}{U\times\R^q} \\[1cm]
\Rnode{p10}{\pi^{-1}(U)}
& \Rnode{p11}{U\times\R^n}
\end{array}
\psset{arrows=->,nodesep=5pt,linewidth=0.3pt} \everypsbox{\scriptstyle}
\ncline{p00}{p01}\Aput{\phi_H}
\ncline{p10}{p11}\Aput{\phi}
\ncline{p00}{p10}\Bput{\incl}
\ncline{p01}{p11}\Aput{\id_U\times i}
\]
commutes (where $i\colon\R^q\to\R^n$ is the standard inclusion). Of course, $\phi_H$ is determined by $\phi$.
\smallskip\\
To every such $\phi$ we associate the map $\Psi_\phi \colon \Pi^{-1}(U) \to U\times\Compl(\R^q)$ (where $\Compl(\R^q)$ refers to complements of $\R^q\subseteq\R^n$) given by $V\mapsto \im(\phi\restrict V)$. Note that the projection of $\im(\phi\restrict V)\subseteq \set{\Pi(V)}\times\R^n$ onto the second factor is indeed a sub vector space of $\R^n$ which is complementary to $\R^q$; moreover, observe that $\Psi_\phi$ is bijective. The set $\mathscr{A}_{\Compl(H)}$ of all pairs $(U\times\Compl(\R^q), \Psi_\phi)$, where $\phi\colon \pi^{-1}(U)\to U\times\R^n$ is a vector bundle trivialisation as above, is a smooth [real-analytic] atlas on the set $\Compl(H)$; since $U\times\Compl(\R^q)$ is an open subset of the affine space $\R^{\dim M}\times\Compl(\R^q)$, it makes sense here to use the maps $\Psi_\phi$ as charts. With respect to this smooth [real-analytic] structure, $\Compl(H)$ is a smooth [real-analytic] affine bundle modelled on the vector bundle $\Lin(E/H,H)$; in fact, $\mathscr{A}_{\Compl(H)}$ is an affine bundle atlas.
\end{definition}

\begin{remark}
If, in the situation of the preceding definition, a smooth [resp.\ real-analytic] sub vector bundle $V$ of $E$ is given which is complementary to $H$, then $\Compl(H)$ becomes a smooth [real-analytic] vector bundle with zero section $V$. Moreover, there is then a canonical vector bundle isomorphism $\Compl(H)=\Lin(E/H,H)\cong\Lin(V,H)$, induced by the vector bundle isomorphism $E/H\cong V$.
\end{remark}

\begin{definition}
Let $E\to M$ be a smooth [real-analytic] vector bundle, and let $q\in\set{0,\dots,\rank(E)}$. Then there is a unique smooth [real-analytic] structure on the total space of the set-theoretic bundle $\Gr_q(E)\to M$ (whose fibre over $x\in M$ is $\Gr_q(E_x)$) such that for each open subset $U$ of $M$ and each smooth [real-analytic] rank-$q$ sub vector bundle $H$ of $E\restrict U$, the set $\Compl(H)$ is an open subset and a smooth [real-analytic] submanifold of $\Gr_q(E)$. We equip $\Gr_q(E)\to M$ with this structure, thereby turning it into a smooth [real-analytic] fibre bundle.
\end{definition}

\begin{remark} \label{equivsmoothdefs}
Let $E\to M$ be a smooth [resp.\ real-analytic] vector bundle, let $q\in\set{0,\dots,\rank(E)}$, let $k\in\N\cup\set{\infty}$, let $V$ be a set-theoretic section in $\Gr_q(E)\to M$ (i.e.\ a set-theoretic map $M\to\Gr_q(E)$ which assigns to each $x\in M$ an element of $\Gr_q(E_x)$). Then there are two definitions of the statement that \emph{$V$ is of class $C^k$} [or \emph{real-analytic}]: According to the first definition, $V$ is $C^k$ if and only if it is $C^k$ as a map between the manifolds $M$ and $\Gr_q(E)$. According to the second definition, $V$ is $C^k$ if and only if the set-theoretic sub vector bundle $V$ of $E$ is of class $C^k$, i.e., if and only if there is a sub vector bundle atlas of $V$ in $E$ all of whose transition maps are $C^k$. These two definitions are equivalent.
\end{remark}


\chapter{Analytic background} \label{AppendixB}

\section{Ellipticity and elliptic regularity}

\emph{Throughout this section, let $I$ be an open interval in $\R$.}
\medskip\\
Let $M$ be a manifold, let $k\in\N_{\geq1}$. Recall that a map $P\colon C^\infty(M,I)\to C^\infty(M,\R)$ or $P\colon C^k(M,I)\to C^0(M,\R)$ is a \emph{$k$th-order (partial) differential operator} if and only if there exists a map $\Phi\colon J^k(M,I)\to\R$ such that $P(u)(x)=\Phi(j^k_x(u))$ for all $x\in M$ and $u\in C^\infty(M,I)$. Such a $\Phi$ is then uniquely determined. A $k$th-order partial differential operator is \emph{quasilinear} if and only if the restriction of the defining map $\Phi$ to each fibre of the vector bundle $J^k(M,I)\to J^{k-1}(M,I)$ is linear.
\smallskip\\
We define the notion of \emph{ellipticity} only for quasilinear second-order differential operators $C^\infty(M,I)\to C^\infty(M,\R)$; that is all we need in this thesis.

\begin{definition}[elliptic, positively elliptic] \label{ellipticdef}
Let $M$ be a manifold, let $P\colon C^\infty(M,I)\to C^\infty(M,\R)$ be a quasilinear second-order differential operator induced by a map $\Phi\colon J^2(M,I)\to\R$. Recall that the vector bundle $J^2(M,I)\to J^1(M,I)$ is the pullback of the vector bundle $\Sym(TM)\to M$ (whose fibre over $x\in M$ consists of the symmetric bilinear forms $T_xM\times T_xM \to \R$) via the bundle projection $p^1\colon J^1(M,I)\to M$. So for each $b\in J^1(M,I)$, the restriction of $\Phi$ to the fibre over $b$ of the bundle $J^2(M,I)\to J^1(M,I)$ is a linear map $\Sym(T_{p^1(b)}M)\to\R$, i.e., it can be identified with an element $\beta_b$ of $\Sym(T^\ast_{p^1(b)}M)$.
\smallskip\\
The quasilinear operator $P$ is called \emph{elliptic} if and only if the symmetric bilinear form $\beta_b$ is definite for each $b\in J^1(M,I)$. The operator $P$ is called \emph{positively elliptic} if and only if $\beta_b$ is positive definite for each $b\in J^1(M,I)$.
\end{definition}

\begin{example} \label{ellipticexample}
Let $(M,g)$ be a semi-Riemannian $n$-manifold. Then $\laplace_g\colon C^\infty(M,\R)\to C^\infty(M,\R)$ is a second-order (quasi)linear differential operator induced by the following section $\beta$ in the $p^1$-pullback of the vector bundle $\Sym(T^\ast M)\to M$: for each $b\in J^1(M,\R)$, the bilinear form $\beta_b\colon T^\ast_{p^1(b)}M\times T^\ast_{p^1(b)}M \to \R$ is the bilinear form $\eval{.}{.}_g$ induced by $g$; cf.\ Notation \ref{musical}. This is obvious from the expression of $\laplace_g$ with respect to local coordinates: for each $u\in C^\infty(M,\R)$, we have $\laplace_g(u) = \sum_{i,j=1}^ng^{ij}\partial_i\partial_ju +\lot$ (where $\lot$ denotes terms which depend only on the $1$-jet of $u$).
\smallskip\\
Hence $\laplace_g$ is elliptic if and only if the index $q$ of $g$ is equal to $0$ or $n$. It is positively elliptic if and only if $q=0$, i.e., if and only if $g$ is Riemannian.
\end{example}

\begin{theorem}[elliptic regularity]
Let $M$ be a manifold, let $\Phi\colon J^2(M,I) \to \R$ be smooth, and let $P\colon C^2(M,I) \to C^0(M,\R)$ be the second-order differential operator induced by $\Phi$ via $P(u)(x)=\Phi(j^2_x(u))$. If $P$ is elliptic and $u\in C^2(M,I)$ is a solution of $P(u)=0$ which satisfies a local $C^{2,\alpha}$ Hölder condition\footnote{I.e., there is an $\alpha\in\oointerval{0}{1}$ for which $M$ admits an atlas $\mathscr{A}$ such that, for each chart $\varphi\colon U\to V\subseteq\R^n$ in the atlas, the function $u\compose\varphi^{-1}\colon V\to I$ is of Hölder class $C^{2,\alpha}$.}, then $u$ is smooth.
\end{theorem}

A proof of the following (as well as of the preceding) theorem can be found in the article \cite{Hopf1931} by E.\ Hopf.\footnote{A more general theorem can be found in \cite{Morrey} (Theorem 6.7.6).} Hopf stated it only for functions defined on a subset of euclidean space, but since it is obviously a local result, it generalises immediately to arbitrary manifolds. Note that a real-analytic structure on a manifold $M$ induces a real-analytic structure on the jet manifold $J^2(M,I)$.

\begin{theorem}[real-analytic version of elliptic regularity] \label{hopf}
Let $M$ be a real-analytic manifold, let $\Phi\colon J^2(M,I)\to\R$ be a real-analytic map, and let $P\colon C^2(M,I) \to C^0(M,\R)$ be the second-order differential operator induced by $\Phi$ via $P(u)(x)=\Phi(j^2_x(u))$. If $P$ is elliptic and $u\in C^2(M,I)$ is a solution of $P(u)=0$ which satisfies a local $C^{2,\alpha}$ Hölder condition, then $u$ is real-analytic.
\end{theorem}


\section{Sobolev spaces} \label{sobolevappendix}

\subsubsection{Generalities and basic inequalities}

Let $(M,g)$ be a compact Riemannian $n$-manifold, let $k\in\N$ and $p\in\R_{>1}$. Then we can consider the Sobolev space $\Sob{p}{k}(M,\R)$ (which is a Banach space); cf.\ e.g.\ \cite{Hebey} or Chapter 2 in \cite{Aubin1982} for its definition and basic properties. We denote the norm on $\Sob{p}{k}(M,\R)$ by $\norm{.}_{\Sob{p}{k}}$. The norm depends on the Riemannian metric $g$ on $M$, but the topology it induces is independent of $g$. (Even for fixed $g$, different definitions of $\norm{.}_{\Sob{p}{k}}$ are used in the literature, but all these definitions yield the same banachisable topology on $\Sob{p}{k}(M,\R)$.) So $\Sob{p}{k}(M,\R)$ is well-defined as a banachisable topological vector space for every compact manifold $M$ (without specification of a Riemannian metric). For example, $\Sob{p}{0}(M,\R) = L^p(M,\R)$.
\smallskip\\
The smooth functions $M\to\R$ form a dense subset of $\Sob{p}{k}(M,\R)$. The closure in $\Sob{p}{k}(M,\R)$ of the set of all smooth functions $f\colon M\to\R$ whose support is a compact subset of $M\without\mfbd M$ is denoted by $\Sobzero{p}{k}(M,\R)$. If $M$ is closed, then clearly $\Sobzero{p}{k}(M,\R) = \Sob{p}{k}(M,\R)$. If there is a continuous imbedding $\Sob{p}{k}(M,\R)\to C^0(M,\R)$ (the latter being equipped with the usual $C^0$ norm, i.e.\ with the compact-open topology), then obviously every $u\in\Sobzero{p}{k}(M,\R)$ vanishes on the boundary of $M$.

\begin{theorem}[the \Poincare\ inequality] \label{poincareinequality}
Let $(M,g)$ be a compact connected Riemannian manifold with nonempty boundary. Then there is a constant $c\in\R_{>0}$ such that the inequality
\[
\norm{u}_{L^2} \leq c\,\lVert\abs{du}_g\rVert_{L^2}
\]
holds for all $u\in \Sobzero{2}{1}(M,\R)$ (and thus in particular for all $u\in C^\infty(M,\R)$ with $u\restrict\mfbd M = 0$).
\end{theorem}
\emph{Remark.} Usually we write just $\norm{du}_{L^2}$ instead of $\lVert\abs{du}_g\rVert_{L^2}$. I just wanted to emphasise here that the $L^2$ norm on the right hand side of the inequality is the $L^2$ norm of a (continuous) function.

\bigskip
Let $M$ be a compact manifold, let $k\in\N$ and $p\in\R_{>1}$. Recall that every $k$th-order linear differential operator $P\colon C^\infty(M,\R) \to C^\infty(M,\R)$ induces a continuous linear map $\bar{P}\colon \Sob{p}{k}(M,\R)\to L^p(M,\R)$.

\begin{theorem}[the elliptic estimate] \label{ellipticestimate}
Let $M$ be a compact manifold, let $P\colon C^\infty(M,\R)\to C^\infty(M,\R)$ be an elliptic second-order differential operator. Then there is a constant $c\in\R_{>0}$ such that for every $u\in\Sob{2}{2}(M,\R)$, the induced operator $\bar{P}\colon \Sob{2}{2}(M,\R)\to L^2(M,\R)$ satisfies the inequality
\[
\norm{u}_{\Sob{2}{2}} \leq c\big(\norm{u}_{L^2} +\lVert\bar{P}(u)\rVert_{L^2}\big) \;\;.
\]
\end{theorem}
\Proof
Cf.\ e.g.\ Theorem III.5.2 in \cite{LawsonMichelsohn}.
\end{proof}

\subsubsection{Imbedding theorems}

\begin{theorem}[Sobolev imbedding into $C^r$] \label{sobolev}
Let $M$ be a compact $n$-manifold, let $k,r\in\N$ and $p\in\R_{>1}$ with $k>\frac{n}{p}+r$. Then there is a canonical continuous injective linear map $\Sob{p}{k}(M,\R)\hookrightarrow C^r(M,\R)$; here $C^r(M,\R)$ is equipped with the $C^r$ topology.
\end{theorem}
\Proof
Cf.\ \cite{Aubin1982}, Theorem 2.30, second part (b).
\end{proof}
\emph{Remark.} The theorem allows us to identify $\Sob{p}{k}(M,\R)$ with a sub vector space of $C^r(M,\R)$ if the condition $k>\frac{n}{p}+r$ is satisfied.

\begin{lemma}
Let $T\colon X\to Y$ be a compact operator between Banach spaces over $\R$. If a sequence $(x_m)_{m\in\N}$ converges weakly in $X$ to $x$, then $(T(x_m))_{m\in\N}$ converges in $Y$ to $T(x)$.
\end{lemma}
\Proof
Cf.\ e.g.\ Proposition 3.4.34 and 3.4.36 in \cite{Megginson}.
\end{proof}

\begin{theorem}[the $p=n$ special case of the Rellich/Kondrakov theorem] \label{rellich}
Let $M$ be a compact $n$-manifold, let $q\in\R_{>1}$. Then there is a canonical compact injective linear map $\Sob{n}{1}(M,\R)\hookrightarrow L^q(M,\R)$. If a sequence in $\Sob{n}{1}(M,\R)$ converges weakly in $\Sob{n}{1}(M,\R)$ to $u$, then it converges to $u$ in $L^q(M,\R)$.
\end{theorem}
\Proof
For the first statement, cf.\ e.g.\ \cite{Palais}, Theorem 9.1 (and Theorem 3.4.37 in \cite{Megginson}). The second statement follows from the preceding lemma.
\end{proof}

\subsubsection{Continuity and differentiability}

\begin{theorem} \label{multcont}
Let $(M,g)$ be a compact Riemannian $n$-manifold, let $k\in\N$ and $p\in\R_{>1}$ with $k>\frac{n}{p}$, and let $l\in\set{0,\dots,k}$. Then the product of an element of $\Sob{p}{k}(M,\R)$ with an element of $\Sob{p}{l}(M,\R)$ is well-defined as an element of $\Sob{p}{l}(M,\R)$; moreover, the multiplication map $\Sob{p}{k}(M,\R)\times \Sob{p}{l}(M,\R)\to \Sob{p}{l}(M,\R)$ is (bilinear and) continuous.
\end{theorem}
\Proof
Cf.\ e.g.\ \cite{Palais}, Corollary 9.7.
\end{proof}

\begin{theorem}[the Sobolev product rule] \label{productrule}
Let $M$ be a compact $n$-manifold, let $k\in\N$ and $p\in\R_{>1}$ with $k>\frac{n}{p}$. Let $U$ be an open subset of a Banach space $B$, let $F,G\colon U\to \Sob{p}{k}(M,\R)$ be (\Frechet) differentiable maps. Then their pointwise product $H\define F\cdot G \colon U\to \Sob{p}{k}(M,\R)$ is also differentiable; for each $u\in U$, the derivative $D_uH\colon B\to \Sob{p}{k}(M,\R)$ is given by
\[
(D_uH)(v) = (D_uF)(v)\cdot G(u) +F(u)\cdot(D_uG)(v) \;\;.
\]
\end{theorem}
\Proof
Cf.\ \cite{LangRFA}, Chapter XIII, \S3; take Theorem \ref{multcont} into account.
\end{proof}

\begin{definition} \label{compdef}
Let $M$ be a compact $n$-manifold, let $k\in\N$ and $p\in\R_{>1}$ with $k>\frac{n}{p}$, and let $I$ be an open subset of $\R$. By Theorem \ref{sobolev}, there is a continuous inclusion $i\colon \Sob{p}{k}(M,\R)\to C^0(M,\R)$. We define $\Sob{p}{k}(M,I)$ to be the $i$-preimage of $C^0(M,I)$. This $\Sob{p}{k}(M,I)$ is an open subset of $\Sob{p}{k}(M,\R)$ since $C^0(M,I)$ is an open subset of $C^0(M,\R)$.
\smallskip\\
Let $f\in C^\infty(I,\R)$. We define the map $\comp_f\colon \Sob{p}{k}(M,I) \to \Sob{p}{k}(M,\R)$ by $u\mapsto f\compose u$; it is well-defined and continuous by \cite{Palais}, Theorem 9.10. (That theorem deals only with the case $I=\R$, but the proof of the general case is analogous.)
\end{definition}

The following rule is standard (cf.\ e.g.\ \cite{KazdanWarner1}; also \cite{McDuffSalamonNew}, Proposition B.1.20):

\begin{theorem}[the Sobolev chain rule] \label{chainrule}
Let $M$ be a compact $n$-manifold, let $k\in\N$ and $p\in\R_{>1}$ with $k>\frac{n}{p}$, let $I$ be an open subset of $\R$, and let $f\in C^\infty(I,\R)$. Then $\comp_f\colon \Sob{p}{k}(M,I) \to \Sob{p}{k}(M,\R)$ is differentiable; its derivative $D_u\comp_f \colon \Sob{p}{k}(M,\R) \to \Sob{p}{k}(M,\R)$ in the point $u$ is given by
\[
(D_u\comp_f)(v) = (f'\compose u)\cdot v \;\;.
\]
\end{theorem}

In order to prove that certain derivatives are continuous, we use the following simple lemmata.

\begin{lemma} \label{conti1}
Let $M$ be a compact $n$-manifold, let $U$ be a topological space, let $X$ be a Banach space, let $p\in\R_{>1}$, let $k\in\N$ with $k>\frac{n}{p}$, let $L\in\Lin(X,\Sob{p}{k}(M,\R))$, let $F\colon U\to L^p(M,\R)$ be a continuous map. Then the map $m\colon U\to\Lin(X,L^p(M,\R))$ given by $m(u)(v) = F(u)\cdot L(v)$ is continuous.
\end{lemma}
\Proof
By \ref{multcont}, there is a constant $c\in\R_{>0}$ such that for all $u,u_0\in U$, we have
\[ \begin{split}
\norm{m(u)-m(u_0)}_{\Lin(\dots)} &= \sup\set{\norm{F(u)L(v)-F(u_0)L(v)}_{L^p} \suchthat \norm{v}_X\leq1}\\
&\leq c\,\sup\set{\norm{F(u)-F(u_0)}_{L^p}\cdot\norm{L(v)}_{\Sob{p}{k}} \suchthat \norm{v}_X\leq1}\\
&\leq c\,\norm{F(u)-F(u_0)}_{L^p}\cdot\norm{L}_{\Lin(\dots)} \;\;.
\end{split} \]
If $u$ tends to $u_0$ in $U$, then $m(u)$ tends to $m(u_0)$ since $F$ is continuous. This implies that $m$ is continuous.
\end{proof}

\begin{lemma} \label{conti2}
Let $(M,g)$ be a compact Riemannian $n$-manifold, let $V$ be a distribution on $M$, let $p\in\R_{>1}$, let $k\in\N$ with $k>\frac{n}{p}+1$. Then the map $\sigma\colon \Sob{p}{k}(M,\R)\to \Lin(\Sob{p}{k}(M,\R), L^p(M,\R))$ given by $\sigma(u)(v)\define \eval{du}{dv}_{g,V}$ is continuous.
\end{lemma}
\Proof
By \ref{multcont}, there is a constant $c\in\R_{>0}$ such that for all $u,u_0\in\Sob{p}{k}(M,\R)$, we have
\[ \begin{split}
\norm{\sigma(u)-\sigma(u_0)}_{\Lin(\dots)} &= \sup\set{\norm{\eval{du}{dv}_{g,V} -\eval{du_0}{dv}_{g,V}}_{L^p} \suchthat \norm{v}_{\Sob{p}{k}}\leq1}\\
&\leq \sup\set{\lVert{\abs{du -du_0}_g\cdot\abs{dv}_g\rVert}_{L^p} \suchthat \norm{v}_{\Sob{p}{k}}\leq1}\\
&\leq c\,\sup\set{\lVert{\abs{du -du_0}_g\rVert}_{L^p}\cdot \lVert{\abs{dv}_g\rVert}_{\Sob{p}{k-1}} \suchthat \norm{v}_{\Sob{p}{k}}\leq1}\\
&\leq c\lVert{\abs{du -du_0}_g\rVert}_{L^p}\\
&\leq c\norm{u-u_0}_{\Sob{p}{k}} \;\;.
\end{split} \]
If $u$ tends to $u_0$ in $\Sob{p}{k}(M,\R)$, then $\sigma(u)$ tends to $\sigma(u_0)$. This implies that $\sigma$ is continuous.
\end{proof}

\subsubsection{Invertibility of elliptic operators and the inverse function theorem}

\begin{proposition}
Let $M$ be a compact manifold, let $P\colon C^\infty(M,\R)\to C^\infty(M,\R)$ be a linear elliptic second-order differential operator. If the induced continuous linear operator $\bar{P}\colon \Sobzero{p}{2}(M,\R)\to L^p(M,\R)$ is injective, then it is bijective.
\end{proposition}
\Proof
Cf.\ e.g.\ the remarks and references in \cite{KazdanWarner1}, \S3.
\end{proof}

The \emph{zeroth-order coefficient} of a linear partial differential operator $P\colon C^\infty(M,\R) \to C^\infty(M,\R)$ is by definition the function $P(1)\in C^\infty(M,\R)$.

\begin{theorem} \label{invertible}
Let $M$ be a compact manifold, let $P\colon C^\infty(M,\R)\to C^\infty(M,\R)$ be a positively elliptic linear second-order differential operator whose zeroth-order coefficient is everywhere nonpositive and, moreover, is negative in at least one point of each boundaryless connected component of $M$. Then the induced linear operator $\bar{P}\colon \Sobzero{p}{2}(M,\R)\to L^p(M,\R)$ is bijective.
\end{theorem}
\Proof
It suffices to consider the case of a connected manifold $M$. By the preceding proposition, it just remains to show that for every $u\in \Sobzero{p}{2}(M,\R)$, the equation $\bar{P}(u)=0$ implies $u=0$. This equation tells us, by elliptic regularity (e.g.\ in the version of Corollary 8.11 in \cite{GilbargTrudinger}), that $u$ is smooth. From the sign condition on the zeroth-order coefficient and the strong maximum principle (cf.\ e.g.\ \cite{GilbargTrudinger}, Theorem 3.5), we infer that if $u$ assumes its maximum and minimum in the interior of $M$, then $u$ is constant. If $M$ is closed, then this must happen; so $u$ vanishes because the zeroth-order coefficient of $P$ is not identically zero. If the boundary of $M$ is nonempty, then maximum and minimum of $u$ are assumed on the boundary. Since $u\in \Sobzero{p}{2}(M,\R)\cap C^\infty(M,\R)$ vanishes on the boundary, it is identically zero. Thus $u=0$ in each case.
\end{proof}

We state the implicit function theorem for Banach spaces in the (restricted) form in which it is applied in Chapters \ref{SIX} and \ref{SEVEN}. Recall that if $X,Z$ are real Banach spaces, then $\Lin(X,Z)$ denotes the Banach space of all continuous linear maps $X\to Z$.

\begin{theorem}[the implicit function theorem for Banach spaces] \label{implicitfunction}
Let $X,S,Z$ be Banach spaces over $\R$, let $A$ be an affine space modelled on $X$, let $(u_0,s_0)\in A\times S$, let $\mathscr{N}\subseteq A$ be an open neighbourhood of $u_0$, let $\Phi\colon \mathscr{N}\times S \to Z$ be a continuous function with $\Phi(u_0,s_0)=0$ such that, for every fixed $s\in S$, the map $\Phi_s\colon \mathscr{N}\to Z$ given by $u\mapsto\Phi(u,s)$ is (\Frechet) differentiable\footnote{I.e., the partial derivative of $\Phi$ in the first component exists.}, such that the map $\mathscr{N}\times S\to \Lin(X,Z)$ given by $(u,s)\mapsto D_u\Phi_s$ is continuous at $(u_0,s_0)$, and such that the derivative $D_{u_0}\Phi_{s_0}\in \Lin(X,Z)$ is bijective. Then there exist an open neighbourhood $\mathscr{U}\subseteq S$ of $s_0$ and a continuous function $U\colon \mathscr{U}\to \mathscr{N}$ such that $U(s_0)=u_0$ and $\Phi(U(s),s)=0$ for all $s\in\mathscr{U}$.
\end{theorem}
\Proof
Cf.\ e.g.\ \cite{ZeidlerI}, Theorem 4.B (p.~150). (The version there assumes that $A$ is a Banach space, but of course the theorem holds for general affine spaces modelled on Banach spaces, or even for Banach manifolds.)
\end{proof}

\begin{remark} \label{implicitremark}
We apply this theorem in Chapters \ref{SIX} and \ref{SEVEN} in the following situation: Let $M$ be a compact $n$-manifold, let $p\in\R$ with $p>n$, let $f\in\Sob{p}{2}(M,\R_{>0})\subseteq C^1(M,\R_{>0})$ (cf.\ Definition \ref{compdef}). Then we consider the case
\begin{itemize}
\item
$X = \Sobzero{p}{2}(M,\R)$ (thus $X = \Sob{p}{2}(M,\R)$ if $M$ is closed);
\item
$A = f+\Sobzero{p}{2}(M,\R) \equiv \set{f+u \suchthat u\in\Sobzero{p}{2}(M,\R)} \subseteq \Sob{p}{2}(M,\R)$ (thus $A = \Sob{p}{2}(M,\R)$ if $M$ is closed);
\item
$S = Z = L^p(M,\R)$;
\item
$\mathscr{N} = A\cap\Sob{p}{2}(M,\R_{>0})$ (this is an open subset of $A$; moreover, $\mathscr{N} = \Sob{p}{2}(M,\R_{>0})$ if $M$ is closed).
\end{itemize}
\end{remark}

We conclude this section with a version of elliptic regularity for Sobolev spaces (a very restricted version which suffices for our needs). Let $M$ be a compact $n$-manifold, let $p>n$, let $I\subseteq\R$ be an open interval, let $P_1\colon C^\infty(M,I) \to C^\infty(M,\R)$ be a first-order differential operator. Then $P_1$ induces a map $\bar{P_1}\colon \Sob{p}{2}(M,I) \to L^p(M,\R)$: since $\Sob{p}{2}(M,I) \subseteq C^1(M,I)$ by Theorem \ref{sobolev}, we have $P_1(u)\in C^0(M,\R)\subseteq L^p(M,\R)$ for all $u\in\Sob{p}{2}(M,I)$.

\begin{theorem}[elliptic regularity for Sobolev spaces] \label{sobolevregularity}
Let $(M,g)$ be a compact Riemannian $n$-manifold, let $p>n$, let $I\subseteq\R$ be an open interval, let $P_1\colon C^\infty(M,I) \to C^\infty(M,\R)$ be a first-order differential operator, let $P$ denote the elliptic quasilinear second-order differential operator $\laplace_g +P_1\colon C^\infty(M,I) \to C^\infty(M,\R)$, and let $u\in\Sob{p}{2}(M,I)$ satisfy $\bar{P}(u)=0$, where $\bar{P} = \bar{\laplace_g}+\bar{P_1} \colon \Sob{p}{2}(M,I) \to L^p(M,\R)$ denotes the operator induced by $P$. Then $u\in C^\infty(M,I)$.
\end{theorem}
\Proof
Cf.\ e.g.\ \cite{Taylor3}, \S14, Proposition 4.9; take Theorem \ref{sobolev} into account.
\end{proof}


\section{The method of sub- and supersolutions}

For the elliptic equation that is relevant for the pseudo-Riemannian prescribed scalar curvature problem, the method of sub- and supersolution can be used to prove that solutions exist. There is a version of this method for Hölder spaces, due to Y.\ Choquet-Bruhat and J.\ Leray (cf.\ \cite{ChoquetBruhatLeray}), and a version for Sobolev spaces (cf.\ e.g.\ \cite{KazdanKramer}). Since we apply the method only to equations with smooth coefficients, it doesn't matter which version we use. The article \cite{ChoquetBruhatLeray} is a better reference in our situation because it deals with manifolds instead of subsets of euclidean space, and because it does not assume that the equation is defined for all $u\colon M\to\R$, but just e.g.\ for $u\colon M\to\R_{>0}$. (However, formulations which do not a priori have these features can of course easily be brought into this form as well, so the lack of these properties is just a technical inconvenience.)

\begin{definition}[sub-/supersolution]
Let $M$ be a manifold, let $I\subseteq\R$ be an open interval, and let $P\colon C^2(M,I) \to C^0(M,\R)$ be a positively elliptic (cf.\ Definition \ref{ellipticdef}) second-order differential operator. A function $f\in C^2(M,I)$ is a \emph{subsolution} [resp. \emph{supersolution}] \emph{of $P$} if and only if $0\leq P(f)$ [resp.\ $0\geq P(f)$].
\end{definition}

If a positively elliptic differential equation has a solution, then is has (for tautological reasons) a subsolution $f_-$ and a supersolution $f_+$ with $f_-\leq f_+$. Sub- and supersolution theorems say that, under certain conditions on the equation, the converse is true as well. We need only the following special case (note that the boundary $\mfbd M$ may be empty):

\begin{theorem}[Choquet-Bruhat/Leray] \label{choquetbruhatleray}
Let $(M,g)$ be a compact Riemannian manifold, let $I\subseteq\R$ be an open interval, let $\varphi\in C^\infty(\mfbd M,I)$, let $a\in C^\infty(T^\ast M\times I,\R)$. We consider the differential operator $P\colon C^\infty(M,I)\to C^\infty(M,\R)$ given by
\[
P(f) \define 2\laplace_g(f) +a(df,f) \;\;.
\]
Let $f_-\in C^\infty(M,I)$ be a subsolution of $P$, and let $f_+\in C^\infty(M,I)$ be a supersolution of $P$, such that $f_- < f_+$ and $f_-\restrict\mfbd M < \varphi < f_+\restrict\mfbd M$. Let $\pi\colon T^\ast M\to M$ denote the bundle projection. Assume that there is a $\mu\in\R_{>0}$ such that the following inequality holds for all $p\in T^\ast M$ and all $u\in[f_-(\pi(p)),f_+(\pi(p))]$:
\[
\abs{a(p,u)} \leq \mu\cdot(1+\abs{p}_g^2) \;\;.
\]
Then there exists a function $f\in C^\infty(M,I)$ such that $P(f)=0$ and $f_- \leq f \leq f_+$ and $f\restrict\mfbd M = \varphi$.
\end{theorem}
\Proof
This is a special case of \Theoreme\ 1 in \cite{ChoquetBruhatLeray} (except that, as a matter of taste, I assumed $a$ to be a function on $T^\ast M\times I$ instead of $TM\times I$, and I reversed the sign of $a$): with the notation there, we consider the case where the functions $\beta_0,\gamma_0$ are constant, where $m=2$, where the map $A$ is given by $A(x,p,u)=2p$, and where $a,\varphi$ are smooth. Then the uniform ellipticity condition $(1)_1$ of Choquet-Bruhat and Leray is automatically satisfied, and their condition $(1)_2$ is equivalent to the condition on $a$ in our theorem.\footnote{Choquet-Bruhat and Leray define the notions \emph{subsolution} and \emph{supersolution} by strict inequalities, but that makes no difference. If you don't believe me in this respect, note that all the sub- and supersolutions in this thesis are sub- and supersolutions in the strict sense of Choquet-Bruhat and Leray.}
\smallskip\\
Thus \Theoreme\ 1 in \cite{ChoquetBruhatLeray} yields a function $f\in C^{2,\alpha}(M,\R)$ with $P(f)=0$ and $f_- \leq f \leq f_+$ and $f\restrict\mfbd M = \varphi$. By elliptic regularity, $f$ is smooth in our special case here.
\end{proof}

We apply the sub- and supersolution method in Chapter \ref{SIX} as follows.

\begin{example} \label{susuexample}
Let $(M,g)$ be a compact Riemannian $n$-manifold, let $V$ be a $q$-plane distribution on $M$, let $H$ denote the $g$-orthogonal distribution of $V$, let $\pi\colon T^\ast M\to M$ denote the bundle projection, and let $s\in C^\infty(M,\R)$. With the functions $a_{n,q}$ and $b_{n,q}$ from \ref{PDdef}, we define a map $a\in C^\infty(T^\ast M\times \R_{>0},\R)$ by
\[ \begin{split}
a(p,u) &\define a_{n,q}(u)\abs{p}^2_g +b_{n,q}(u)\abs{p}^2_{g,V} +\frac{2(1+u^2)}{u^2}\Big\langle{(\divergence^V_g)(\pi(p))},\;{p}\Big\rangle_{g,H} +2(1+u^2)\Big\langle{(\divergence^H_g)((\pi(p))},\;{p}\Big\rangle_{g,V}\\
&\mspace{20mu}+\frac{(1+u^2)^2}{2u^3}\abs{\Twist_H}^2_g(\pi(p)) -\frac{u(1+u^2)^2}{2}\abs{\Twist_V}^2_g(\pi(p)) +\frac{(1+u^2)^2}{u}\xi_{g,V}(\pi(p))\\
&\mspace{20mu}+\frac{1+u^2}{u}\scal_g(\pi(p)) -u^{\frac{2q}{n-1}-1}(1+u^2)^{1-\frac{1}{n-1}}s(\pi(p)) \;\;.
\end{split} \]
The positively elliptic quasilinear second-order differential operator $P\colon C^\infty(M,\R_{>0})\to C^\infty(M,\R)$ given by
\[
P(f) \define 2\laplace_g(f) +a(df,f)
\]
is obviously just the operator $\PD_{g,V,s}$ from \ref{PDdef}.
\smallskip\\
If $P$ admits a subsolution $f_-\in C^\infty(M,\R_{>0})$ and a supersolution $f_+\in C^\infty(M,\R_{>0})$ with $f_- < f_+$, then the compactness of $M$, the continuity of all involved functions, and the estimate $\abs{p}_{g,V}\leq\abs{p}_g$ imply that there is a constant $\mu\in\R_{>0}$ such that
\[
\abs{a(p,u)} \leq \mu\cdot(1+\abs{p}_g^2)
\]
holds for all $p\in T^\ast M$ and all $u\in[\inf(f_-),\sup(f_+)]\subseteq \R_{>0}$. (Because several coefficient functions in the definition of $a$ tend to $\infty$ as $u$ tends to $0$ or $\infty$, there is usually no constant $\mu>0$ such that this inequality holds for \emph{all} $u\in\R_{>0}$. But that is no problem.)
\smallskip\\
We choose any boundary values $\varphi\in C^\infty(\mfbd M,\R_{>0})$ such that $f_-\restrict\mfbd M < \varphi < f_+\restrict\mfbd M$. (We will not be interested in boundary value problems in Chapter \ref{SIX}, so it makes no difference which $\varphi$ we choose.)
\smallskip\\
Now we can deduce from Theorem \ref{choquetbruhatleray} that there exists a solution $f\in C^\infty(M,\R)$ of the equation $P(f)=0$ such that $f_-\leq f\leq f_+$.
\end{example}


\chapter{The topology of metrics and distributions} \label{AppendixC}

This subsection summarises the relation between distributions and semi-Riemannian metrics, which forms the basis for our discussion of the distribution problem and the homotopy class problem in Chapter \ref{ONE}. Such a summary is not easy to find in the literature, although all these things are certainly well-known. For instance, \S40 in \cite{Steenrod} contains the main results in a coarse form, but not the facts from the exposition below about metrics which make a fixed distribution timelike, and distributions which are timelike with respect to a fixed metric.
\smallskip\\
Everything in this subsection generalises to arbitrary semi-Riemannian vector bundles, but we restrict ourselves to the case of tangent bundles, for otherwise the important points might be obscured by unnecessary generality.

\section{Distributions vs. semi-Riemannian metrics: first facts}

\begin{definition}
Let $n\in\N$ and $q\in\set{0,\dots,n}$.
\smallskip\\
\emph{Vector space level.} Let $W$ be an $n$-dimensional real vector space equipped with a symmetric nondegenerate bilinear form $g$ of index $q$. A sub vector space of $W$ is called \emph{maximally timelike} if and only if it is timelike (cf.\ \ref{distributionremarks}(ii)) and has dimension $q$. A sub vector space of $W$ is called \emph{maximally spacelike} if and only if it is spacelike and has dimension $n-q$. If $V$ is any sub vector space of $W$, then the \emph{$g$-orthogonal subspace of $V$}, denoted by $\bot_gV$, is the vector space formed by all $w\in W$ such that $g(w,v)=0$ for all $v\in V$.
\smallskip\\
\emph{Vector bundle level.} Let $(M,g)$ be an $n$-dimensional semi-Riemannian manifold of index $q$. A distribution on $M$ is called \emph{maximally timelike} if and only if it is timelike and has rank $q$. A distribution on $M$ is called \emph{maximally spacelike} if and only if it is spacelike and has rank $n-q$. If $V$ is any distribution on $M$, then the \emph{$g$-orthogonal distribution of $V$}, denoted by $\bot_gV$, is the distribution whose fibre over each point $x\in M$ is $\bot_gV_x$ (where $V_x$ denotes the fibre of $V$ over $x$); it is easy to check that this is indeed a smooth sub vector bundle of $TM\to M$.
\end{definition}

\begin{facts} \label{factsten}
We state the following obvious facts on the vector bundle level, but they hold already on the vector space level.
\smallskip\\
If $V$ is a maximally timelike [maximally spacelike] distribution on the semi-Riemannian manifold $(M,g)$, then the $g$-orthogonal distribution of $V$ is complementary to $V$ and maximally spacelike [maximally timelike]. (Note that if $g$ is pseudo-Riemannian and $V$ is an arbitrary distribution on $M$, then the $g$-orthogonal distribution of $V$ will in general not be complementary to $V$; cf.\ \ref{nondegfacts} in Chapter \ref{TWO}.)
\smallskip\\
If $V$ is a maximally timelike distribution on $(M,g)$ and $H$ is a maximally spacelike distribution on $(M,g)$, then $V$ and $H$ are complementary.
\smallskip\\
For every distribution $V$ on a manifold $M$, there is a distribution on $M$ which is complementary to $V$ (take the orthogonal distribution with respect to any Riemannian metric on $M$).
\end{facts}

\begin{proposition} \label{propeleven}
Let $M$ be a manifold, and let $V,H$ be complementary distributions on $M$, where $V$ has rank $q$. Then there is a semi-Riemannian metric (of index $q$) on $M$ which makes $V$ timelike, $H$ spacelike, and $V$ and $H$ orthogonal.
\end{proposition}
\Proof
We equip the vector bundles $V$ and $H$ with Riemannian metrics $g_V$ and $g_H$, respectively. The bilinear form $(-g_V)\oplus g_H$ on the vector bundle $TM = V\oplus H$ is a semi-Riemannian metric on $M$ which makes $V$ and $H$ orthogonal, $V$ timelike, and (hence) $H$ spacelike.
\end{proof}

This proposition explains why the assumption \emph{which admits a semi-Riemannian metric of index $q$} from the plain problem is superfluous in the distribution problem(s). The converse of the proposition is also true:

\begin{theorem}[existence of time-/spacelike distributions] \label{baumtheorem}
Let $(M,g)$ be a semi-Riemannian manifold. Then there is a $g$-orthogonal decomposition $TM = V\oplus H$, where $V$ is a maximally timelike and $H$ is a maximally spacelike distribution on $M$.
\end{theorem}
\Proof
Cf.\ \cite{Baum}, Satz 0.48. The main ingredients of the proof are: (A) Every fibre bundle with contractible fibres admits a section; this is a corollary of Theorem \ref{smoothobstruction}. (B) The space $\OO(q,n-q)/(\OO(q)\times\OO(n-q))$ is contractible since $\OO(q)\times\OO(n-q)$ is a maximally compact subgroup of $\OO(q,n-q)$.\footnote{Note that in \cite{Baum}, H.\ Baum uses a different convention for pseudo-orthogonal groups: what she calls $\OO(n,n-q)$ is our $\OO(q,n-q)$. The latter notation seems to be standard.} This argument can be varied slightly; cf.\ Remark \ref{altargument} below.
\end{proof}

This tells us that the situation in the distribution problem is no more special than the situation in the plain problem: Every solution of the plain problem is also a solution of the distribution problem, for suitable distributions $V$ and $H$.

\begin{proposition}[maximally time-/spacelike distributions are isomorphic] \label{isomorphy}
Let $(M,g)$ be a semi-Riemann\-ian manifold. Every two maximally timelike distributions on $M$ are isomorphic as (smooth) vector bundles over $M$.\footnote{Recall that the notion of morphism of vector bundles \emph{over $M$} is defined by a commutative diagram where the map between the base spaces is the identity, in contrast to a more general definition of morphism of vector bundles over possibly different base spaces; cf.\ \cite{Husemoller}, \S\S\ 2.3.1, 2.3.2.} Every two maximally spacelike distributions on $M$ are isomorphic as vector bundles over $M$.
\end{proposition}
\Proof
Choose a maximally spacelike distribution $H$ on $(M,g)$; this is possible by the previous theorem. (In fact, we don't need the theorem here. Since we assume that two maximally timelike distributions are given, we can just take the $g$-orthogonal distribution of one of them.) Every maximally timelike distribution is complementary to $H$, hence isomorphic to the quotient vector bundle $TM/H$. Isomorphy of maximally spacelike distributions is proved analogously.
\end{proof}

The preceding proposition is also a corollary of \ref{isomorphicVBs} and \ref{distri} below.

\begin{remark}
Different distributions of the same rank on a manifold $M$ will in general not be isomorphic as vector bundles over $M$. For example, if $V_0$ and $V_1$ are distributions on $M$, their first Stiefel/Whitney classes might be different. This happens if $V_0$ is orientable but $V_1$ is not orientable, for instance; cf.\ Examples \ref{Kleinbottle} and \ref{torusexample} below.
\end{remark}

\section[The fixed metric viewpoint]{Affine structures and convex subsets. The fixed metric viewpoint}

Recall from Appendix \ref{Grassmannappendix} the definition of the Grassmann bundle $\Gr_q(TM)\to M$ associated to the vector bundle $TM\to M $.

\begin{remark}[the topology of Grassmannians]
Let $W$ be an $n$-dimensional real vector space. Note that $\Gr_1(W)$ is diffeomorphic to the real-projective space $\RP^{n-1}$, and that $\Gr_q(W)$ is diffeomorphic to $\Gr_{n-q}(W)$ for all $q\in\set{0,\dots,n}$: for any scalar product $g$ on $W$, the map $V\mapsto \bot_gV$ is a diffeomorphism $\Gr_q(W)\to\Gr_{n-q}(W)$.
\smallskip\\
The $q(n-q)$-dimensional manifold $\Gr_q(W)$ is connected. However, it is not simply connected if $0<q<n$: the fundamental group of $\Gr_q(W)$ is isomorphic to $\Z_2$ if $1<q<n$ and $(n,q)\neq(2,1)$; it is isomorphic to $\Z$ if $(n,q)=(2,1)$ (cf.\ \cite{Whitehead}, Theorems 10.16, 10.12; take $\pi_1(\Gr_q(W)) \cong \pi_1(\Gr_{n-q}(W))$ and $\pi_0(\OO(q))\cong\Z_2$ into account). Many of the higher homotopy groups have been calculated as well by algebraic topologists.
\smallskip\\
Since too few homotopy groups of $\Gr_q(W)$ vanish, not every Grassmann bundle $\Gr_q(TM)\to M$ over an $n$-manifold $M$ admits a section; that is, not every $n$-manifold admits a $q$-plane distribution. And even if $\Gr_q(TM)\to M$ admits a section, the space of all sections will in general not be connected. One of the aims of this subsection is to investigate this phenomenon more closely.
\end{remark}

\begin{definition}[$\Time$ and $\Space$] \label{TimeSpacedef}
Let $n\in\N$ and $q\in\set{0,\dots,n}$.
\smallskip\\
\emph{Vector space level.} Let $W$ be an $n$-dimensional real vector space equipped with a symmetric nondegenerate bilinear form $g$ of index $q$. We define $\Time(g)$ [resp.\ $\Space(g)$] to be the set of all maximally timelike [maximally spacelike] sub vector spaces of $(W,g)$.
\smallskip\\
\emph{Vector bundle level.} Let $(M,g)$ be a semi-Riemannian $n$-manifold of index $q$. For each $x\in M$, we denote the restriction of $g$ to $T_xM$ by $g_x$. We define a (set-theoretic) fibre bundle $\Time(g)$ [resp.\ $\Space(g)$] over $M$ by declaring that its fibre over each $x\in M$ be $\Time(g_x)$ [resp.\ $\Space(g_x)$].
\end{definition}

\begin{facts}
In the vector space level situation of the previous definition, $\Time(g)$ and $\Space(g)$ are nonempty open subsets of the Grassmannians $\Gr_q(W)$ and $\Gr_{n-q}(W)$, respectively.
\smallskip\\
The total spaces of $\Time(g)\to M$ and $\Space(g)\to M$ are open subsets of the total spaces of the Grassmann bundles $\Gr_q(TM)\to M$ and $\Gr_{n-q}(TM)\to M$, respectively. The set-theoretic fibre bundles $\Time(g)\to M$ and $\Space(g)\to M$ inherit the structure of (locally trivial) smooth fibre bundles from these Grassmann bundles.
\smallskip\\
These statements are easy to verify from the definitions (cf.\ Appendix \ref{Grassmannappendix}) of the topologies on the Grassmannians resp.\ Grassmann bundles. Henceforth, we consider $\Time(g)\to M$ and $\Space(g)\to M$ as smooth bundles.
\end{facts}

\begin{definition}[affine structure on the set of complementary distributions]{\ } \label{affinedef}
\smallskip\\
\emph{Vector space level.} Let $W$ be an $n$-dimensional real vector space, and let $H$ be an $(n-q)$-dimensional sub vector space of $W$. We denote the set of all sub vector spaces of $W$ which are complementary to $H$ by $\Compl(H)$. Recall from Appendix \ref{Grassmannappendix} that this set is canonically equipped with the structure of an affine space modelled on the vector space $\Lin(W/H,H)$. Namely, the affine space operation $+ \colon \Lin(W/H,H)\times\Compl(H) \to \Compl(H)$ is given by $(\lambda,V)\mapsto \set{v+\lambda(\pi(v)) \suchthat v\in V}$, where $\pi\colon W\to W/H$ denotes the canonical projection.
\smallskip\\
\emph{Vector bundle level.} Let $M$ be a manifold, let $H$ be a distribution on $M$. We define a (set-theoretic) fibre bundle $\Compl(H)\to M$ by declaring that its fibre over $x\in M$ be $\Compl(H_x)$ (where $H_x$ is the fibre of $H$ over $x$). It is easy to check that the total space $\Compl(H)$ is an open subset of $\Gr_q(TM)$ and thus inherits a smooth structure, which turns $\Compl(H)\to M$ into a smooth affine bundle modelled on the vector bundle $\Lin(W/H,H)\to M$.
\end{definition}

\begin{lemma}
Let $W$ be a finite-dimensional real vector space equipped with a nondegenerate symmetric bilinear form. Let $v,w\in W$ be timelike [resp.\ spacelike] vectors such that $v-w$ is spacelike [timelike]. Then $tv+(1-t)w$ is timelike [spacelike] for every $t\in[0,1]$.
\end{lemma}
\Proof
We consider only the version for timelike $v,w$; the other version is proved analogously. Let $g$ denote our bilinear form. The function $f\colon[0,1]\to\R$ defined by $t\mapsto g(tv+(1-t)w,tv+(1-t)w)$ is negative in $0$ and $1$ since $v$ and $w$ are timelike. Its maximum is assumed in one of these points since the second derivative of $f$ is everywhere positive: we have $f(t) = g(t(v-w)+w,t(v-w)+w) = t^2g(v-w,v-w) +2tg(v-w,w) +g(w,w)$. Thus $f$ is everywhere negative.
\end{proof}

\begin{proposition}[convexity of $\Time$ and $\Space$] \label{convexityproposition}
{\ }\smallskip\\
\emph{Vector space level.} Let $W$ be a finite-dimensional real vector space equipped with a symmetric nondegenerate bilinear form $g$, and let $H$ be a maximally spacelike [resp.\ maximally timelike] sub vector space of $(W,g)$. Then $\Time(g)$ [resp.\ $\Space(g)$] is a nonempty open convex subset of the affine space $\Compl(H)$.
\smallskip\\
\emph{Vector bundle level.} Let $(M,g)$ be a semi-Riemannian manifold, and let $H$ be a maximally spacelike [resp.\ maximally timelike] distribution on $M$. Then $\Time(g)$ [resp.\ $\Space(g)$] is a (locally trivial) open convex subbundle, with nonempty fibres, of the affine bundle $\Compl(H)$.
\end{proposition}
\Proof
Nonemptyness and openness have already been mentioned. It suffices to prove convexity on the vector space level; so let $V_0$ and $V_1$ be maximally timelike subspaces of $(W,g)$ (the spacelike case works analogously). We have to show that each subspace $V_t\define (1-t)V_0+tV_1\in \Compl(H)$, where $t\in[0,1]$, is timelike.
\smallskip\\
Fix $t\in[0,1]$, and let $\pi\colon W\to W/H$ denote the canonical projection. With the map $\lambda\in\Lin(W/H,H)$ determined by $V_1 = \lambda+V_0 = \set{v_0+\lambda(\pi(v_0)) \suchthat v_0\in V_0}$, we have $V_t = t\lambda+V_0 = \set{v_0+t\lambda(\pi(v_0)) \suchthat v_0\in V_0}$, by definition of the affine structure on $\Compl(H)$. We have to prove that every nonzero element of $V_t$ is timelike, i.e.\ that $v_0+t\lambda(\pi(v_0))$ is timelike for every nonzero $v_0\in V_0$.
\smallskip\\
Consider such a $v_0$. The vector $v_1\define v_0+\lambda(\pi(v_0))\in V_1$ is timelike since $V_1$ is timelike ($v_1$ can't be zero since $\lambda(\pi(v_0))\in H$ and $V_0\cap H = \set{0}$). If $v_1=v_0$, then $v_0+t\lambda(\pi(v_0))=v_0$ is timelike. If $v_1-v_0\in H$ is nonzero and thus spacelike, then the preceding lemma implies that $v_0+t\lambda(\pi(v_0)) = tv_1+(1-t)v_0$ is timelike.
\end{proof}

Let me emphasise that the convex structures from the previous proposition depend on the choice of $H$, while the sets on which they are defined do not.

\begin{remark} \label{altargument}
Let $(M,g)$ be a semi-Riemannian manifold. The fibre bundle $\Time(g)$ over $M$ has contractible fibres: the fibre over $x$ is convex with respect to the affine structure defined by the choice of a maximally spacelike subspace of $T_xM$. Since every fibre bundle with contractible fibres admits a section, there is a maximally timelike distribution on $M$ (whose $g$-orthogonal distribution is maximally spacelike). This argument is a minor variation of the proof of Theorem \ref{baumtheorem}.
\end{remark}

\section{The fixed distribution viewpoint}

Now we change our viewpoint: instead of fixing a metric and considering the set of maximally timelike distributions, we fix a distribution and consider the set of metrics which make this distribution maximally timelike. Note the analogy between the following discussion and the treatment of the fixed metric viewpoint.

\begin{definition} \label{symqdef}
\emph{Vector space level.} Let $W$ be a finite-dimensional real vector space. We use the notation $\Sym(W)$ for the vector space of all symmetric bilinear forms on $W$. We denote its open subset consisting of all nondegenerate forms with index $q$ on $W$ by $\Sym_q(W)$.
\smallskip\\
\emph{Vector bundle level.} Let $E\to M$ be a vector bundle. We denote the vector bundle of all symmetric bilinear forms on $E$ by $\Sym(E)\to M$, and its (locally trivial) open sub fibre bundle consisting of all nondegenerate forms with index $q$ by $\Sym_q(E)\to M$.
\end{definition}

$\Sym_q$ can be turned into a (co)functor in an obvious way (on the vector space as well as on the vector bundle level), analogously to the definition of the functor $\Gr_q$. We will neither explain that in detail, nor will we use it.

\begin{facts}
Let $n\in\N$, let $W$ be an $n$-dimensional real vector space, and let $q\in\set{0,\dots,n}$.
\smallskip\\
$\Sym_q(W)$ is a manifold of dimension $n(n+1)/2$. There is a canonical diffeomorphism $\Sym_q(W) \to \Sym_{n-q}(W)$, given by $g\mapsto -g$. Note that $\Sym_0(W)$ and $\Sym_n(W)$ are nonempty convex subsets of $\Sym(W)$, hence contractible. $\Sym_q(W)$ is diffeomorphic to $\Sym_q(\R^n)$.
\smallskip\\
Let $\eval{.}{.}_q$ be the element of $\Sym_q(\R^n)$ which is, with respect to the standard basis $(e_1,\dots,e_n)$ of $\R^n$, given by the diagonal matrix $E_q=\diag(-1,\dots,-1,1,\dots,1)$, where exactly $q$ entries are equal to $-1$. Recall that $\OO(q,n-q)$ is by definition the (closed) subgroup of $\GL(n)$ consisting of all $\eval{.}{.}_q$-isometries, i.e.\ of all $A\in\GL(n)$ with $A^\top E_qA = E_q$. The map $\GL(n)\to \Sym_q(\R^n)$ which sends each $A\in\GL(n)$ to the bilinear form $(v,w)\mapsto \eval{Av}{Aw}_q$ (which is given by the matrix $A^\top E_qA$) induces a diffeomorphism $\GL(n)/\OO(q,n-q)\to\Sym_q(\R^n)$. Since every connected component of $\GL(n)$ intersects $\OO(q,n-q)$ nontrivially, this implies in particular that $\Sym_q(W)$ is connected. (Cf.\ \ref{steenrod} for further information on the homotopy type of $\Sym_q(W)$.)
\end{facts}
\Proof
That the map $a\colon\GL(n)\to\Sym_q(\R^n)$ given by $A\mapsto A^\top E_qA$ is well-defined and surjective follows from Sylvester's inertia theorem (cf.\ e.g.\ \cite{SchejaStorch2}, \S72). We see immediately from the definition of $\OO(q,n-q)$ that $a$ induces an injective map $\bar{a}\colon\GL(n)/\OO(q,n-q)\to\Sym_q(\R^n)$. This map is smooth and surjective since $a$ is smooth and surjective. For every $A\in\GL(n)$, the kernel of the derivative $T_Aa \colon T_A\GL(n) \to T_{a(A)}\Sym_q(\R^n)$ is equal to $T_A\OO(q,n-q)$. So the smooth map $\bar{a}$ is regular and thus a diffeomorphism. The other facts do not require a proof.
\end{proof}

\begin{definition}[$\Timifier$ and $\Spacifier$] \label{TimifierSpacifierdef}
Let $n\in\N$ and $q\in\set{0,\dots,n}$.
\smallskip\\
\emph{Vector space level.} Let $W$ be an $n$-dimensional real vector space, and let $V$ be a $q$-dimensional [resp.\ $(n-q)$-dimensional] sub vector space of $W$. We define $\Timifier(V)$ [resp.\ $\Spacifier(V)$] to be the (nonempty open) subset of $\Sym_q(W)$ consisting of all forms which make $V$ (maximally) timelike [spacelike].
\smallskip\\
\emph{Vector bundle level.} Let $M$ be an $n$-manifold, and let $V$ be a $q$-plane [resp.\ $(n-q)$-plane] distribution on $M$. We define $\Timifier(V)$ [resp.\ $\Spacifier(V)$] to be the (locally trivial) open sub fibre bundle of $\Sym_q(TM)$ whose fibre over $x$ is $\Timifier(V_x)$ [resp.\ $\Spacifier(V_x)$].
\end{definition}

\begin{remark}[$\Timifier$ and $\Spacifier$ are not convex]
Let $W$ be an $n$-dimensional real vector space, and let $V$ be a $q$-dimensional [resp. an $(n-q)$-dimensional] sub vector space of $W$. If $0<q<n$, then $\Timifier(V)$ [resp.\ $\Spacifier(V)$] is not convex as a subset of $\Sym(W)$. We check this in the $\Timifier$ case; the $\Spacifier$ case is analogous.
\smallskip\\
Before we start, note that for all $g_0,g_1\in\Timifier(V)$ and $t\in[0,1]$, the convex combination $(1-t)g_0+tg_1\in\Sym(W)$ makes $V$ timelike. The problem is that it can be degenerate or even nondegenerate with index greater than $q$. Now we show that this does always happen for some $g_0,g_1$ if $0<q<n$.
\smallskip\\
Without loss of generality, let $W=\R^n$ with basis $(e_1,\dots,e_n)$, and let $V=\spann\set{e_1,e_3,\dots,e_{q+1}}$; that is, $V=\spann\set{e_1}$ if $n=2$. Let $D$ be the diagonal $(n-2)\times(n-2)$-matrix $\diag(-1,\dots,-1,1,\dots,1)$, with $q-1$ entries equal to $-1$ and $n-q-1$ entries equal to $1$. Let $C_0 = \bigl(\begin{smallmatrix} -1&\phantom{-}2\\ \phantom{-}2&-1 \end{smallmatrix}\bigr)$ and $C_1 = \bigl(\begin{smallmatrix} -1&-2\\ -2&-1 \end{smallmatrix}\bigr)$. We consider the symmetric bilinear forms $g_0,g_1$ on $W$ given by the block matrices
\[
\begin{pmatrix} C_0&0\\ 0&D\end{pmatrix} \mspace{30mu}
\text{and} \mspace{30mu} \begin{pmatrix} C_1&0\\ 0&D\end{pmatrix} \;\;,
\]
respectively. These forms are nondegenerate with index $q$ since $C_0$ and $C_1$ are nondegenerate with index $1$ (this follows from $\det(C_i)<0$). It is easy to check that $g_0$ and $g_1$ make $V$ timelike.
\smallskip\\
However, the convex combination $\frac{1}{2}g_0+\frac{1}{2}g_1$ is given by the $n\times n$-matrix $\diag(-1,\dots,-1,1\dots,1)$ with $q+1$ entries equal to $1$, so it is nondegenerate with index $q+1$. (There are also convex combinations of $g_0$ and $g_1$ which are degenerate, for continuity reasons.) This proves the claimed nonconvexity.
\end{remark}

\begin{definition}[bundles of metrics] \label{metricbundledef}
Let $n\in\N$ and $q\in\set{0,\dots,n}$.
\smallskip\\
\emph{Vector space level.} Let $W$ be an $n$-dimensional real vector space. Consider the smooth vector bundle $p(\Sym^{(q)}) \colon \Sym^{(q)} \to \Gr_q(W)$ whose fibre over $V\in\Gr_q(W)$ is the vector space $\Sym(V)$ consisting of all symmetric bilinear forms on $V$. (The smooth structure of this bundle is induced by the universal vector bundle\footnote{cf.\ \ref{grassmanntangent}} over the Grassmannian $\Gr_q(W)$, from which $p(\Sym^{(q)})$ is obtained in a functorial way, via the cofunctor that assigns to each vector space $V$ the vector space $\Sym(V)$; cf.\ e.g.\ \cite{KolarMichorSlovak}, II.6.7.) This vector bundle has as an open convex subbundle the fibre bundle $p(\Sym_0^{(q)}) \colon \Sym_0^{(q)} \to \Gr_q(W)$ whose fibre over $V\in\Gr_q(W)$ is the set $\Sym_0(V)$ consisting of all positive definite bilinear forms on $V$.
\smallskip\\
If $V$ is a $q$-dimensional sub vector space of $W$, then we denote the restrictions of $p(\Sym^{(n-q)})$ and $p(\Sym_0^{(n-q)})$ to the open subspace $\Compl(V)\subseteq\Gr_{n-q}(W)$ by $p(\Sym^{\bot V}) \colon \Sym^{\bot V} \to \Compl(V)$ and $p(\Sym_0^{\bot V}) \colon \Sym_0^{\bot V} \to \Compl(V)$, respectively.
\medskip\\
\emph{Vector bundle level.} Let $M$ be an $n$-manifold. Consider the smooth vector bundle $p(\Sym^{(q)}) \colon \Sym^{(q)} \to \Gr_q(TM)$ whose fibre over $V\in\Gr_q(TM)$ is the vector space $\Sym(V)$. (The smooth structure of this bundle is induced by the universal vector bundle over the total space of the Grassmannian bundle $\Gr_q(TM)\to M$.) This vector bundle has as an open convex subbundle the fibre bundle $p(\Sym_0^{(q)}) \colon \Sym_0^{(q)} \to \Gr_q(TM)$ whose fibre over $V\in\Gr_q(TM)$ is the set $\Sym_0(V)$ consisting of all positive definite bilinear forms on $V$.
\smallskip\\
If $V$ is a $q$-plane distribution on $M$, then we denote the restrictions of $p(\Sym^{(n-q)})$ and $p(\Sym_0^{(n-q)})$ to $\Compl(V)\subseteq\Gr_{n-q}(TM)$ by $p(\Sym^{\bot V}) \colon \Sym^{\bot V} \to \Compl(V)$ and $p(\Sym_0^{\bot V}) \colon \Sym_0^{\bot V} \to \Compl(V)$, respectively.
\end{definition}

Recall that a \emph{deformation} of some smooth manifold $X$ onto some point $a\in X$ is a continuous map $r\colon [0,1]\times X \to X$ such that for all $t\in[0,1]$, the map $r_t\colon X\to X$ given by $x\mapsto r(t,x)$ is smooth and we have $r_0=\id_X$, $r_1\equiv a$, and $r(t,a)=a$. A \emph{bundle deformation} of some smooth fibre bundle $p\colon F\to M$ onto some section $s\in C^\infty(M\ot F)$ is a continuous map $r\colon [0,1]\times F\to F$ such that for all $t\in[0,1]$ and $x\in M$, the map $r_t\colon F\to F$ is smooth, $r_t(x)$ lies in the fibre over $x$, and we have $r_0=\id_F$, $r_1=s\compose p$, and $r_t\compose s = s$.

\begin{proposition}[$\Timifier$ and $\Spacifier$ are contractible] \label{timifier}
Let $n\in\N$ and $q\in\set{0,\dots,n}$.
\smallskip\\
\underline{\emph{Vector space level.}} Let $W$ be an $n$-dimensional real vector space, and let $V$ be a $q$-dimensional sub vector space of $W$. Then there is a canonical diffeomorphism $\Timifier(V) \to \Sym_0(V)\times\Sym_0^{\bot V}$; it assigns to each $g$ the pair $(-g\restrict V,\; g\restrict\bot_gV)$. The bundle $\Sym_0^{\bot V}\to\Compl(V)$ admits a section.
\smallskip\\
Every triple $(b,\xi,w)$, where $b\in\Sym_0(V)$ and $\xi\in C^\infty(\Compl(V)\ot\Sym_0^{\bot V})$ and $w\in\Compl(V)$, defines a deformation $r_{b,\xi,w} \colon [0,1]\times\Sym_0(V)\times\Sym_0^{\bot V} \to \Sym_0(V)\times\Sym_0^{\bot V}$ onto the point $(b,\xi(w))$, as follows: using the abbreviation $p$ for the projection $p(\Sym_0^{\bot V}) \colon \Sym_0^{\bot V} \to \Compl(V)$,
\[
r_{b,\xi,w}(t,\beta,\alpha) = \begin{cases}
\Big((1-t)\beta +tb,\; (1-2t)\alpha +2t\xi(p(\alpha))\Big) &\text{if $t\leq\frac{1}{2}$}\\
\Big((1-t)\beta +tb,\; \xi\big((2-2t)p(\alpha) +(2t-1)w\big)\Big) &\text{if $t\geq\frac{1}{2}$}
\end{cases}
\]
(this makes sense since $\Sym_0(V)$, $\Compl(V)$ and all fibres of $p$ are equipped with convex structures).
\smallskip\\
In particular, $\Timifier(V)$ is contractible.\footnote{It is in fact diffeomorphic to $\R^{n(n+1)/2}$, and the bundle $\Sym_0^{\bot V}\to\Compl(V)$ is trivial. But that's not important for us.} Analogously, $\Spacifier(V)$ is contractible. There is a canonical diffeomorphism $\Spacifier(V) \to \Sym_0(V)\times\Sym_0^{\bot V}$ which maps each $g$ to $(g\restrict V,\;-g\restrict\bot_gV)$.
\medskip\\
\underline{\emph{Vector bundle level.}} Let $M$ be an $n$-manifold, and let $V$ be a $q$-plane distribution on $M$. Consider the fibre bundles $\Sym_0(V)\to M$ and $\Sym_0^{\bot V} \to \Compl(V)\to M$ over $M$, and their product fibre bundle $P_V\colon \Sym_0(V)\times\Sym_0^{\bot V} \to M$. The bundles $\Timifier(V)\to M$ and $P_V\to M$ are isomorphic as fibre bundles over $M$; an isomorphism is given fibrewise by the map that we described in the vector space setting above. The bundles $\Sym_0(V)\to M$ and $\Sym_0^{\bot V}\to \Compl(V)$ and $\Compl(V)\to M$ admit sections. Every triple $(b,\xi,w)$, where $b\in C^\infty(M\ot\Sym_0(V))$ and $\xi\in C^\infty(\Compl(V)\ot\Sym_0^{\bot V})$ and $w\in C^\infty(M\ot\Compl(V))$, defines a bundle deformation of $P_V\to M$ onto the section $(b,\xi\compose w)\in C^\infty(M\ot P_V)$, by declaring that the restriction of this bundle deformation to the fibre over $x\in M$ be the deformation $r_{b(x),\xi\restrict\Compl(V_x),w(x)}$.
\smallskip\\
In particular, $\Timifier(V)\to M$ admits a bundle deformation. Analogously, $\Spacifier(V)\to M$ admits a bundle deformation. There is a canonical isomorphism of fibre bundles over $M$ between $P_V$ and $\Spacifier(V)$, which is given fibrewise by the map that we described in the vector space setting above.
\end{proposition}
\Proof
All the statements are obvious or require only routine verifications. We just remark that the bundles under consideration admit sections because they have contractible fibres, and that the inverse diffeomorphism $\Sym_0(V)\times\Sym_0^{\bot V} \to \Timifier(V)$ maps each $(g_0,g_1)$ to the unique metric whose restriction to $V$ is $-g_0$, whose restriction to $H\define p(\Sym_0^{\bot V})(g_1)\in\Compl(V)$ is $g_1$, and which makes $V$ and $H$ orthogonal.
\end{proof}

\section{Topologies}

Let $r\in\N\cup\set{\infty}$. The definitions and basic properties of the compact-open [resp.\ fine] $C^r$-topology on the space $C^r(M\ot F)$ of $C^r$ sections in a fibre bundle $F\to M$ are reviewed in Appendix \ref{diffbackground}. In our context here, where we shall consider the connected components of $C^r(M\ot F)$ for certain bundles $F\to M$, the \emph{fine} $C^r$ topologies are not appropriate; cf.\ Remark \ref{fineremark}. (The compact-open versus fine distinction is not important for the main results of this thesis since they refer to \emph{compact} manifolds $M$, for which each compact-open $C^r$ topology on $C^r(M\ot F)$ coincides with the respective fine $C^r$-topology. But I intend to provide some general background for the prescribed scalar curvature problem here, so we should not make unnecessary restrictive assumptions.) The value of $r$ makes no essential difference as far as connected components of $C^r(M\ot F)$ are concerned; cf.\ \ref{openfacts} below.

\begin{definition}[topology on the set of distributions] \label{distrdef}
Let $M$ be an $n$-manifold, and let $q\in\set{0,\dots,n}$. We denote the set of all $q$-plane distributions on $M$ by $\Distr_q(M)$. It can be identified with the set $C^\infty(M\ot \Gr_q(TM))$ of sections in the Grassmann bundle $\Gr_q(TM)\to M$; cf.\ Remark \ref{equivsmoothdefs}. (The point here is that smoothness of distributions --- which are by definition vector bundles --- and smoothness of Grassmann bundle sections are defined in different ways.) We equip $\Distr_q(M) = C^\infty(M\ot \Gr_q(TM))$ with the compact-open $C^\infty$-topology.
\end{definition}

\begin{definition}[topology on the set of metrics] \label{metrdef}
Let $M$ be an $n$-manifold, and let $q\in\set{0,\dots,n}$. We denote the set of all semi-Riemannian metrics of index $q$ on $M$ by $\Metr_q(M) = C^\infty(M\ot\Sym_q(TM))$ (cf.\ \ref{symqdef} above), and we equip it with the compact-open $C^\infty$-topology.
\end{definition}

The following statements are easy consequences of the general facts reviewed in Appendix \ref{diffbackground}; cf.\ \ref{topologyfacts}.

\begin{facts} \label{openfacts}
Let $M$ be an $n$-manifold and let $q\in\set{0,\dots,n}$.
\smallskip\\
$\Metr_q(M)$ is an open subset of the \Frechet\ space $C^\infty(M\ot\Sym(TM))$; it is a convex subset if $q\in\set{0,n}$. (For $q\in\set{1,\dots,n-1}$, the space $\Metr_q(M)$ is in general not even path connected. In fact, that's the main point of this whole appendix chapter. We will investigate a concrete example below; cf.\ \ref{torusexample}.)
\smallskip\\
Let $H$ be an $(n-q)$-plane distribution on $M$. Then $C^\infty(M\ot\Compl(H))$ is an open subset of $\Distr_q(M)$. Since $\Compl(H)\to M$ is an affine bundle modelled on the vector bundle $\Lin(TM/H,H)\to M$, the set $C^\infty(M\ot\Compl(H))$ is an affine space modelled on the \Frechet\ space $C^\infty(M\ot\Lin(TM/H,H))$, and the (pointwise defined) affine space operation is continuous.
\smallskip\\
The sets $C^\infty(M\ot\Timifier(H))$ and $C^\infty(M\ot\Spacifier(H))$ are open in $\Metr_{n-q}(M)$ and $\Metr_q(M)$, respectively. If $g$ is a semi-Riemannian metric with index $q$ [resp.\ $n-q$] on $M$ and $H$ is $g$-spacelike [$g$-timelike], then $C^\infty(M\ot\Time(g))$ [resp.\ $C^\infty(M\ot\Space(g))$] is a nonempty open convex subset of $C^\infty(M\ot\Compl(H))$.
\smallskip\\
For every fibre bundle $F\to M$ and each $r\in\N$, the inclusion map $C^\infty(M\ot F) \to C^r(M\ot F)$ is a homotopy equivalence with respect to the compact-open $C^r$ topology on the target space. In particular, the inclusion $C^\infty(M\ot\Sym_q(TM)) \to C^r(M\ot\Sym_q(TM))$ induces a bijection $\pi_0(\Metr_q(M)) \to \pi_0(C^r(M\ot\Sym_q(TM)))$, and the inclusion $C^\infty(M\ot\Gr_q(TM)) \to C^r(M\ot\Gr_q(TM))$ induces a bijection $\pi_0(\Distr_q(M)) \to \pi_0(C^r(M\ot\Gr_q(TM)))$, for each $r\in\N$. (This is important when we want to compute $\pi_0(\Metr_q(M))$ and $\pi_0(\Distr_q(M))$ with the methods of algebraic topology, because these methods apply to the compact-open $C^0$-topology. Cf.\ Example \ref{torusexample} below.)
\end{facts}

\begin{proposition} \label{isomorphicVBs}
Let $M$ be an $n$-manifold and let $q\in\set{0,\dots,n}$. If the $q$-plane distributions $V_0$ and $V_1$ lie in the same path component of $\Distr_q(M)$, then they are isomorphic as vector bundles over $M$.
\end{proposition}
\Proof
We consider the vector bundle $E=\pr^\ast(TM)$ over $[0,1]\times M$ which is the pullback of the tangent bundle $TM\to M$ by the obvious projection $\pr\colon [0,1]\times M\to M$. There is a path $\gamma\colon[0,1]\to\Distr_q(M)$ from $V_0$ to $V_1$ which is smooth in the sense of Lemma \ref{smoothpath}; that is, the section $\xi$ in the Grassmann bundle $\Gr_q(E)\to [0,1]\times M$ defined by $\xi(t,x) = \gamma(t)(x)$ is smooth and can thus be identified with a smooth rank-$q$ sub vector bundle $\xi$ of $E$ (by \ref{equivsmoothdefs}).
\smallskip\\
For $i\in\set{0,1}$, let $f_i\colon M\to[0,1]\times M$ be the inclusion $x\mapsto(i,x)$. These two maps are homotopic: the identity $[0,1]\times M \to [0,1]\times M$ is a homotopy from $f_0$ to $f_1$. Hence the vector bundles $f_0^\ast(\xi) \cong V_0$ and $f_1^\ast(\xi) \cong V_1$ are isomorphic; cf.\ e.g.\ \cite{Husemoller}, Theorem 3.4.7.\footnote{This theorem applies to $C^0$ vector bundles, but it holds for smooth vector bundles as well: one can either generalise the proof directly, or apply Theorem 4.3.5 in \cite{Hirsch}.}
\end{proof}

If $F\to M$ is a fibre bundle, then the following terminology is customary: The path component of $C^\infty(M\ot F)$ which contains a given section $s$ is called the \emph{homotopy class of $s$}. Two sections in $F\to M$ are called \emph{homotopic} if and only if they lie in the same homotopy class. We will usually adopt this terminology from the next subsection on.

\section[Path components]{Path components of $\Metr_q(M)$ versus path components of $\Distr_q(M)$ and $\Distr_{n-q}(M)$}

Now we come to the core of this subsection.

\begin{proposition}[path component of time-/spacelike distributions] \label{distri}
Let $(M,g)$ be a semi-Riemannian $n$-manifold of index $q$. The set $C^\infty(M\ot\Time(g))$ of all maximally timelike distributions on $M$ is a contractible subset of $\Distr_q(M)$. The set $C^\infty(M\ot\Space(g))$ of all maximally spacelike distributions on $M$ is a contractible subset of $\Distr_{n-q}(M)$.
\end{proposition}
\Proof
$(M,g)$ admits a maximally spacelike distribution $H$ and a maximally timelike distribution (by \ref{baumtheorem}). It follows from Proposition \ref{convexityproposition} that $C^\infty(M\ot\Time(g))$ is a nonempty convex subset of the affine space $C^\infty(M\ot\Compl(H))$. Since this affine space is modelled on a topological vector space, the affine space operation being continuous (cf.\ \ref{openfacts}), $C^\infty(M\ot\Time(g))$ is contractible. Analogously, $C^\infty(M\ot\Space(g))$ is contractible.
\end{proof}

\begin{definition}[time-/spacelike distribution component] \label{distrcompdef}
Let $M$ be an $n$-manifold.
\smallskip\\
Let $g$ be a semi-Riemannian metric on $M$ with index $q$. We call the unique path component of $\Distr_q(M)$ [resp.\ $\Distr_{n-q}(M)$] which contains all the maximally timelike [maximally spacelike] distributions on $(M,g)$ the \emph{$g$-timelike [$g$-spacelike] distribution component}. (We shall omit the prefix ``$g$-'' when the metric is clear from the context.) This path component is well-defined because the set of all maximally timelike [maximally spacelike] distributions is contractible by the preceding proposition, thus in particular nonempty and path connected.
\smallskip\\
We define the function $\tdc \colon \Metr_q(M) \to \pi_0(\Distr_q(M))$ by sending each $g$ to its timelike distribution component, and $\sdc \colon \Metr_q(M) \to \pi_0(\Distr_{n-q}(M))$ by sending each $g$ to its spacelike distribution component.
\end{definition}

Note that the timelike [resp.\ spacelike] distribution component contains usually distributions which are not timelike [spacelike]: it contains, for instance, all distributions which are complementary to a given spacelike [timelike] one, and clearly not all of those are timelike [spacelike], except in the trivial cases $q\in\set{0,n}$ or $M=\leer$.

\begin{fact} \label{surjfact}
Let $M$ be an $n$-manifold, and let $q\in\set{0,\dots,n}$. Then $\tdc \colon \Metr_q(M) \to \pi_0(\Distr_q(M))$ and $\sdc \colon \Metr_q(M) \to \pi_0(\Distr_{n-q}(M))$ are surjective. This follows immediately from the last fact in \ref{factsten} and Proposition \ref{propeleven}.
\end{fact}

\begin{proposition} \label{tdcsdc}
Let $n\in\N$ and $q\in\set{0,\dots,n}$, let $M$ be an $n$-manifold. If $g_0$ and $g_1$ are contained in the same path component of $\Metr_q(M)$, then $\tdc(g_0)=\tdc(g_1)$ and $\sdc(g_0)=\sdc(g_1)$.
\end{proposition}
\Proof
We prove only $\tdc(g_0)=\tdc(g_1)$; the proof of $\sdc(g_0)=\sdc(g_1)$ is analogous. If $g_0$ and $g_1$ are contained in the same path component of $\Metr_q(M)$, then there is a path $g_{\blank} \colon [0,1] \to \Metr_q(M)$ from $g_0$ to $g_1$ which is \emph{smooth} (in the sense of Lemma \ref{smoothpath}).
\smallskip\\
Let $\pr\colon [0,1]\times M \to M$ be the projection onto the second factor, and consider the pullback fibre bundle $\pr^\ast(\Gr_q(TM)) \to [0,1]\times M$ of the Grassmann bundle $\Gr_q(TM)\to M$. We define a subbundle of $\pr^\ast(\Gr_q(TM))$: to each $(t,x)\in [0,1]\times M$, we assign as fibre the set $\Time(g_t(x))$. It is easy to verify that the total space of this bundle is an open subset of $\pr^\ast(\Gr_q(TM))$ and thereby inherits a smooth structure which turns it into a (locally trivial) smooth fibre bundle.
\smallskip\\
Since the fibres of this bundle are contractible, the bundle admits a smooth section. This section defines a (smooth) path $V_{\blank}\colon [0,1]\to\Distr_q(M)$ with $V_0 \in C^\infty(M\ot\Time(g_0))$ and $V_1\in C^\infty(M\ot\Time(g_1))$. The path component $\tdc(g_0)\in\pi_0(\Distr_q(M))$ contains $V_0$, the path component $\tdc(g_1)$ contains $V_1$, and $V_0$ and $V_1$ lie in the same path component.
\end{proof}

\begin{definition} \label{TDCSDC}
Let $M$ be an $n$-manifold. We define a function $\TDC \colon \pi_0(\Metr_q(M)) \to \pi_0(\Distr_q(M))$ by sending the path component of each $g$ to $\tdc(g)$, and we define a function $\SDC \colon \pi_0(\Metr_q(M)) \to \pi_0(\Distr_{n-q}(M))$ by sending the path component of each $g$ to $\sdc(g)$. These functions are well-defined by the preceding proposition.
\end{definition}

$\TDC$ and $\SDC$ are surjective by \ref{surjfact}. In order to prove that they are bijections, we will now introduce reverse constructions to those above; i.e., instead of mapping metrics to path components of distributions, we map distributions to path components of metrics.

\begin{proposition}[reverse of \ref{distri}] \label{metri}
Let $M$ be an $n$-manifold, and let $q\in\set{0,\dots,n}$.
\smallskip\\
Let $V$ be a $q$-plane distribution on $M$. The set $C^\infty(M\ot\Timifier(V))$ of all semi-Riemannian metrics with index $q$ on $M$ which make $V$ (maximally) timelike is a contractible open subset of $\Metr_q(M)$.
\smallskip\\
Let $H$ be an $(n-q)$-plane distribution on $M$. The set $C^\infty(M\ot\Spacifier(H))$ of all semi-Riemannian metrics with index $q$ on $M$ which make $H$ (maximally) spacelike is a contractible open subset of $\Metr_q(M)$.
\end{proposition}
\Proof
We prove only the $\Timifier$ case; the $\Spacifier$ case is analogous. By Proposition \ref{timifier}, the bundles $\Sym_0(V)\to M$ and $\Sym_0^{\bot V}\to\Compl(V)$ and $\Compl(V)\to M$ admit sections, and every triple of such sections defines a bundle deformation $r$ of $\Timifier(V)$ onto a certain section $s\in C^\infty(M\ot\Timifier(V))$. The map $[0,1] \times C^\infty(M\ot\Timifier(V)) \to C^\infty(M\ot\Timifier(V))$ which maps each $(t,f)$ to the section given by $x\mapsto r(t,f(x))$ is then a deformation (in the topological sense, without any smoothness requirements) of $C^\infty(M\ot\Timifier(V))$ onto the point $s$.
\end{proof}

\begin{definition}[reverse of \ref{distrcompdef}] \label{metrcompdef}
Let $M$ be an $n$-manifold, and let $q\in\set{0,\dots,n}$.
\smallskip\\
Let $V$ be a $q$-plane distribution on $M$. We call the unique path component of $\Metr_q(M)$ that contains all metrics which make $V$ timelike the \emph{$V$-timifying metric component}.
\smallskip\\
Let $H$ be an $(n-q)$-plane distribution on $M$. We call the unique path component of $\Metr_q(M)$ that contains all metrics which make $H$ spacelike the \emph{$H$-spacifying metric component}.
\smallskip\\
We define the function $\tmc \colon \Distr_q(M) \to \pi_0(\Metr_q(M))$ by sending each $V$ to its timifying metric component, and $\smc \colon \Distr_{n-q}(M) \to \pi_0(\Metr_q(M))$ by sending each $H$ to its spacifying distribution component.
\end{definition}

Note that the timifying [spacifying] metric component contains usually metrics which are not timifying [spacifying]. It is clear that $\tmc$ and $\smc$ are surjective (cf.\ \ref{baumtheorem}).

\begin{proposition}[reverse of \ref{tdcsdc}] \label{tmcsmc}
Let $n\in\N$ and $q\in\set{0,\dots,n}$, let $M$ be an $n$-manifold. If $V_0$ and $V_1$ are contained in the same path component of $\Distr_q(M)$ [resp.\ $\Distr_{n-q}(M)$], then $\tmc(V_0)=\tmc(V_1)$ [resp.\ $\smc(V_0)=\smc(V_1)$].
\end{proposition}
\Proof
There is a path $V_{\blank}\colon[0,1]\to\Distr_q(M)$ which is smooth (in the sense of Lemma \ref{smoothpath}). Consider the vector bundle $\pr^\ast(\Sym(TM))$ over $[0,1]\times M$ which is the pullback of $\Sym(TM)\to M$ by the projection $\pr \colon [0,1]\times M \to M$. We define a (set-theoretic) sub fibre bundle of $\pr^\ast(\Sym(TM))$ by assigning to each $(t,x) \in [0,1]\times M$ the set $\Timifier(V_t(x))$. It is easy to verify that the total space of this fibre bundle is an open subset of the total space of $\pr^\ast(\Sym(TM))$ and thereby inherits a smooth structure which turns it into a (locally trivial) smooth fibre bundle.
\smallskip\\
Since the fibres of this bundle are contractible, the bundle admits a smooth section. This section defines a (smooth) path $g_{\blank}\colon [0,1]\to\Metr_q(M)$ with $g_0 \in C^\infty(M\ot\Timifier(V_0))$ and $g_1\in C^\infty(M\ot\Timifier(V_1))$. The path component $\tmc(V_0)\in\pi_0(\Metr_q(M))$ contains $g_0$, the path component $\tmc(V_1)$ contains $g_1$, and $g_0$ and $g_1$ lie in the same path component.
\end{proof}

\begin{definition}[reverse of \ref{TDCSDC}] \label{TMCSMC}
Let $M$ be an $n$-manifold, and let $q\in\set{0,\dots,n}$. We define a function $\TMC\colon \pi_0(\Distr_q(M)) \to \pi_0(\Metr_q(M))$ by sending the path component of each $V$ to $\tmc(V)$, and we define a function $\SMC\colon \pi_0(\Distr_{n-q}(M)) \to \pi_0(\Metr_q(M))$ by sending the path component of each $H$ to $\smc(H)$. These functions are well-defined by the preceding proposition.
\end{definition}

Now we can state one of the main results of this subsection:

\begin{theorem}
Let $M$ be an $n$-manifold, let $q\in\set{0,\dots,n}$. Then the function $\TDC \colon \pi_0(\Metr_q(M)) \to \pi_0(\Distr_q(M))$ is bijective with inverse $\TMC$, and the function $\SDC \colon \pi_0(\Metr_q(M)) \to \pi_0(\Distr_{n-q}(M))$ is bijective with inverse $\SMC$.\hfill\qed
\end{theorem}

\begin{remark} \label{steenrod}
If $n\in\N$ and $q\in\set{0,\dots,n}$, then the manifolds $\Gr_q(\R^n)$ and $\Sym_q(\R^n)$ are homotopy equivalent. More precisely, if $W$ is an $n$-dimensional real vector space and $g$ is a positive-definite bilinear form on $W$, then the inclusion $i_g\colon\Gr_q(W)\to\Sym_q(W)$ which maps $V$ to the unique bilinear form $h$ with $h\restrict V = -g\restrict V$, $h\restrict\bot_gV = g\restrict\bot_gV$, and $h(v,w)=0$ for all $v\in V$ and $w\in\bot_gV$, is a homotopy equivalence. This is Theorem 40.8 in \cite{Steenrod}.\footnote{Steenrod constructs an explicit deformation. That the inclusion $\Gr_q(\R^n) = \OO(n)/(\OO(q)\times\OO(n-q)) \to \GL(n)/\OO(q,n-q) = \Sym_q(\R^n)$ is a homotopy equivalence could also be proved by considering the inclusion induced chain map between the homotopy exact sequences of the fibrations $\OO(q)\times\OO(n-q) \to \OO(n)$ and $\OO(q,n-q) \to \GL(n)$: One applies the maximal compactness of the subgroup $\OO(n)$ in $\GL(n)$, the maximal compactness of $\OO(q)\times\OO(n-q)$ in $\OO(q,n-q)$, the $5$-lemma, and the fact that weak homotopy equivalences are homotopy equivalences.}
\smallskip\\
Analogous things hold on the bundle level: If $(M,g)$ is a Riemannian $n$-manifold, then there is a bundle deformation of $\Sym_q(TM)\to M$ onto its (pointwise defined) subbundle $i_g(\Gr_q(TM))\to M$; cf.\ Theorem 40.10 in \cite{Steenrod}. This yields, in particular, a bijection between $\pi_0(\Metr_q(M))$ and $\pi_0(\Distr_q(M))$.
\smallskip\\
However, Steenrod uses his Theorem 40.10 only to deduce our statements \ref{propeleven} and \ref{baumtheorem} from above. He does neither mention path components --- in particular does he not discuss the (in)dependence of the bijection $\pi_0(\Metr_q(M)) \to \pi_0(\Distr_q(M))$ on the choice of $g$ ---, nor does he introduce concepts like $\Time$ or $\Timifier$ to clarify the correspondence between metrics and distributions.
\end{remark}

Finally, let us discuss the concept of \emph{complementary distribution component}.

\begin{definition}[complementary distribution component] \label{cdcdef}
Let $M$ be an $n$-manifold, let $q\in\set{0,\dots,n}$. If $V\in\Distr_q(M)$, then the unique path component of $\Distr_{n-q}(M)$ which contains all distributions that are complementary to $V$ is called the \emph{complementary distribution component of $V$}. (Recall that the distributions which are complementary to $V$ form a topological affine space, in particular a path connected space.) We define $\cdc \colon \Distr_q(M) \to \pi_0(\Distr_{n-q}(M))$ by sending each $V$ to its complementary distribution component.
\end{definition}

\begin{proposition}
Let $M$ be an $n$-manifold. Then $\cdc = \SDC\compose\tmc = \TDC\compose\smc$.
\end{proposition}
\Proof
Let $V\in\Distr_q(M)$. Choose $H\in\Compl(V)\subseteq\cdc(V)$. By \ref{propeleven}, there exists a $g\in\Metr_q(M)$ which makes $V$ timelike and $H$ spacelike. Since $H\in\sdc(g)=\SDC(\tmc(V))$ and $H\in\cdc(V)$, we get $\cdc(V)=\SDC(\tmc(V))$. This holds for all $V$, hence $\cdc = \SDC\compose\tmc$. Similarly $\cdc = \TDC\compose\smc$.
\end{proof}

\begin{definition} \label{CDCdef}
If $V_0$ and $V_1$ lie in the same path component of $\Distr_q(M)$, then $\cdc(V_0) = \cdc(V_1)$ because of the preceding proposition and $\tmc(V_0) = \tmc(V_1)$. (Of course there is an alternative direct argument.) We define $\CDC \colon \pi_0(\Distr_q(M)) \to \pi_0(\Distr_{n-q}(M))$ to be the map which sends the path component of each $V\in\Distr_q(M)$ to its complementary distribution component.
\end{definition}

\begin{proposition}
Let $M$ be an $n$-manifold. Then $\CDC = \SDC\compose\TMC = \TDC\compose\SMC$.\hfill\qed
\end{proposition}

\section{Time-orientability and space-orientability}

\begin{definition}[time-/space-orientability]
A semi-Riemannian manifold $(M,g)$ is called \emph{time-orientable} [\emph{space-orientable}] if and only if there is a maximally timelike [maximally spacelike] distribution on $M$ which is orientable as a vector bundle. Hence, by Theorem \ref{baumtheorem} and Proposition \ref{isomorphy}, $(M,g)$ is time-orientable [space-orientable] if and only if \emph{all} maximally timelike [maximally spacelike] distributions on $M$ are orientable.
\smallskip\\
A \emph{time-orientation} [\emph{space-orientation}] of $(M,g)$ is a choice of orientation on one maximally timelike [maximally spacelike] distribution on $M$. Such a choice induces an orientation on each maximally timelike [maximally spacelike] distribution on $M$. Namely, this vector bundle orientation is induced pointwise (depending smoothly on the point\footnote{This smooth (or, equivalently, continuous) dependence is easy to check using the smooth fibre bundle morphism $\Gr^+_q(TM)\to \Gr_q(TM)$ which appears fibrewise in the following definition. We omit the details.}): for every $x\in M$, every orientation on one maximally timelike [maximally spacelike] subspace of $T_xM$ induces an orientation on all maximally timelike [maximally spacelike] subspaces of $T_xM$. To see how, we restrict to the timelike case for convenience; the spacelike case works analogously. Let $q=\ind(g)$.
\smallskip\\
Consider the smooth fibre bundle $p\colon \Gr^+_q(T_xM) \to \Gr_q(T_xM)$, where $\Gr_q(T_xM)$ is the Grassmannian (consisting of all $q$-dimensional sub vector spaces of $T_xM$), $\Gr^+_q(T_xM)$ is the oriented Grass\-mann\-ian (consisting of all oriented $q$-dimensional sub vector spaces of $T_xM$), and the projection is the map which forgets the orientation; each fibre has exactly two elements\footnote{We use a definition which yields two orientations on a $0$-dimensional vector space; cf.\ e.g.\ \cite{SchejaStorch2}, \S74.}. Moreover, consider a fixed maximally timelike [maximally spacelike] subspace $V$ of $T_xM$ and an orientation $o$ on $V$. We want to understand how $o$ defines a section of the fibre bundle $p_{\mathfrak{T}} \define p\restrict\Time(g)_x$. (Such a section assigns to each element of $\Time(g)_x$, i.e.\ to each maximally timelike subspace of $T_xM$, an orientation, and that's precisely what we want.)
\smallskip\\
(If the twofold covering map $p$ admitted a unique section $\gamma\colon \Gr_q(T_xM) \to \Gr^+_q(T_xM)$ with $\gamma(V)=(V,o)$, then we would take $\gamma\restrict\Time(g)_x$ as our desired section. However, $p$ does not admit any section at all if $0<q<n$.\footnote{Because this fact is not important for us, let me just remark that this is an easy application of the lifting theorem. In order to prove that the group homomorphism $\pi_1(p)\colon \pi_1(\Gr^+_q(T_xM);o) \to \pi_1(\Gr_q(T_xM);V)$ is not surjective (as it would have to be for a section to exist), write $T_xM$ as an internal direct sum of subspaces $V_0\oplus V_1\oplus H_0\oplus H_1$, where $V_0\oplus V_1=V$ and $\dim V_1 = \dim H_1 = 1$. A half-rotation in the $2$-dimensional subspace $V_1\oplus H_1$ defines an element of $\pi_1(\Gr_q(T_xM);V)$ which reverses the orientation of $V$ and is therefore not contained in the image of $\pi_1(p)$. (Use the fact that the canonical projection $S^1\to\RP^1$ induces a non-surjective homomorphism $\pi_1(S^1)\to \pi_1(\RP^1)$.)})
\smallskip\\
Since $\Time(g)_x$ is simply connected (cf.\ Proposition \ref{convexityproposition}), the fibre bundle $p_{\mathfrak{T}}$ admits a unique section $\gamma$ with $\gamma(V)=(V,o)$ by the lifting theorem\footnote{For example, with the notations from \cite{Bredon}, Theorem III.4.1, take $X = p_{\mathscr{U}}^{-1}(\mathscr{U}_x)\subseteq \Gr^+_q(T_xM)$, $W=Y=\mathscr{U}_x$, $f=\id_{Y}$, $p=p_{\mathscr{U}}$, $w_0=y_0=V$, $x_0=o$, $g=\gamma$.}. As explained, we use this section to define an orientation on each maximally timelike subspace of $T_xM$. This completes the present definition.
\end{definition}

\begin{remark}[alternative definitions]
Time-/space-orientability resp.\ -orientations can also be characterised by existence resp.\ choices of structure group reductions of the $\OO(q,n-q)$-principal bundle of $g$-orthonormal frames over $M$; cf.\ Satz 0.51 in \cite{Baum}. Or one identifies time-/space-orientations with sections in suitable orientation bundles; cf.\ \cite{ONeill}, pp.~240--242. In the Lorentzian case there is yet another description of time-orientability and -orientations by considering connected components of the set of timelike vectors in $TM$. It is clear that all these definitions are equivalent (in a sense which we do not spell out here because these things are not relevant in our context).
\end{remark}

Time- and space-orientability illustrate one aspect of the distribution problem: If we prescribe an orientable bundle $V$ [resp.\ $H$] in the (time [space]) distribution problem, then every solution metric is time-orientable [space-orientable]. If we prescribe a bundle $V$ [resp.\ $H$] which is not orientable, then every solution metric is not time-orientable [not space-orientable].

\begin{definition}
Let $M$ be an $n$-manifold, let $0\leq q\leq n$. A path component $C\in\Metr_q(M)$ is called \emph{time-orientable} [\emph{space-orientable}] if and only if one element of $C$ is time-orientable [space-orientable]. By \ref{isomorphicVBs} and \ref{tdcsdc}, this is equivalent to \emph{all} elements of $C$ being time-orientable [space-orientable].
\end{definition}

\section{Examples}

We illustrate the above mentioned issues of orientability and path components by some simple examples, with a preference for two-dimensional manifolds because they can be visualised easily and are nonetheless instructive.

\begin{example}[Klein bottle] \label{Kleinbottle}
Let $\klein$ be the Klein bottle $\R^2/\Gamma$, where $\Gamma$ is the discrete subgroup of the euclidean isometry group on $\R^2$ which is generated by the two isometries $(x,y)\mapsto(x+1,-y)$ and $(x,y)\mapsto(x,y+1)$. Let $g$ denote the flat Riemannian metric on $\klein$ which is induced by the euclidean metric on $\R^2$. The Klein bottle is a fibre bundle over $\R/\Z = S^1$ with typical fibre $\R/\Z$; the bundle projection $\tau$ is given by $[x,y]\mapsto[x]$. It is well-known that the manifold $\klein$ is not orientable and that the bundle $\tau$ is not orientable (these statements are equivalent since the base space $S^1$ is orientable).
\smallskip\\
We consider two line distributions on $\klein$. First, there is the fibre distribution $V$ on $\klein$ induced by the bundle structure: it assigns to each $[x,y]\in\klein$ the tangent space in $[x,y]$ to the fibre of $\tau$ over $x$. (The fibre is a one-dimensional submanifold of $\klein$, so the tangent space is a one-dimensional subspace of $T_{[x,y]}\klein$.) The second distribution $H$ is the $g$-orthogonal distribution of $V$.
\smallskip\\
The distribution $H$ is orientable, the distribution $V$ is not orientable. There is a Lorentzian metric $h$ on $M$ which makes $V$ timelike and $H$ spacelike and is thus space-orientable but not time-orientable. (With the notation introduced in \ref{switchdefinition}, we could take the Lorentzian metric $\switch(g,V)$, for instance, which is flat.) The metric $-h$ makes $V$ spacelike and $H$ timelike and is thus time-orientable but not space-orientable.
\smallskip\\
Note that the Klein bottle contains a Möbius strip as an open subset, which provides another example manifold with orientability properties very similar to those of the Klein bottle.
\end{example}

\begin{example}[product manifolds] \label{productexample}
For $i\in\set{1,2}$, let $M_i$ be an $n_i$-manifold. We define the \emph{first-factor distribution on $M_1\times M_2$} to be the $n_1$-plane distribution on $M_1\times M_2$ which is everywhere tangential to $M_1$, i.e.\ whose value in each point $(x_1,x_2)\in M_1\times M_2$ is the vector space $(T_{x_1}M_1)\oplus\set{0} \subseteq (T_{x_1}M_1)\oplus(T_{x_2}M_2) = T_{(x_1,x_2)}(M_1\times M_2)$. (In other words: When $p_1 \colon M_1\times M_2 \to M_1$ is the projection onto the first factor, then the first-factor distribution is the $p_1$-pullback of the tangent bundle $TM_1\to M_1$, if we interpret the pullback as a subbundle of the tangent bundle of $M_1\times M_2$.) We define the \emph{second-factor distribution on $M_1\times M_2$} to be the $n_2$-plane distribution on $M_1\times M_2$ which is everywhere tangential to $M_2$.
\smallskip\\
The first-factor [resp.\ second-factor] distribution on $M_1\times M_2$ is orientable if and only if the manifold $M_1$ [resp.\ $M_2$] is orientable. Since $M_1\times M_2$ is orientable if and only if $M_1$ and $M_2$ are orientable, we see that a pseudo-Riemannian manifold can be simultaneously not orientable, not time-orientable, and not space-orientable. (This phenomenon occurs already for certain metrics on the Klein bottle.)
\end{example}

Clearly, if a semi-Riemannian manifold has two of the three properties \emph{orientable}, \emph{time-orientable}, and \emph{space-orientable}, then it has also the third. Every combination of the properties which is not ruled out by this fact can occur. We have not yet seen an orientable semi-Riemannian manifold which is not time-orientable and not space-orientable, but the next example provides nonorientable line distributions on the (orientable) $2$-torus and thus closes this gap.
\medskip\\
Now we compute $\pi_0(\Distr_1(M))$ for the special case $M=T^2$; that is, we determine how many connected components the space of Lorentzian metrics on the $2$-torus has.

\begin{example}[the $2$-torus] \label{torusexample}
As we remarked above (cf.\ \ref{openfacts}), there is a canonical bijection between the set of path components of $\Distr_1(T^2) = C^\infty(T^2\ot \Gr_1(T(T^2)))$ and the set of path components of $C^0(T^2\ot \Gr_1(T(T^2)))$. The latter can be computed using standard tools of algebraic topology (recall that we refer to compact-open topologies here), as we will do now.
\smallskip\\
Since $T^2$ is parallelisable, the Grassmann bundle $\Gr_1(T(T^2))\to T^2$ is trivial; in fact, there is a canonical trivialisation. We can thus identify the space $C^0(T^2\ot \Gr_1(T(T^2)))$ with $C^0(T^2,\Gr_1(\R^2))$. We want to compute the set of its path components, i.e.\ the set $[T^2,\Gr_1(\R^2)]$ of (free) homotopy classes of continuous maps $T^2\to\Gr_1(\R^2)$. Because $\Gr_1(\R^2) = \RP^1 \cong S^1$ is an Eilenberg/Mac Lane space $K(\Z,1)$, there is a bijection between $[T^2,\Gr_1(\R^2)]$ and $H^1(T^2;\Z)$; cf.\ e.g.\ \cite{Spanier}, Theorem 8.1.10. This bijection is determined by the choice of a characteristic element in $H^1(S^1;\Z)$ (cf.\ \cite{Spanier}, 8.1.3), and there exists a canonical choice. From the Künneth theorem, we get a canonical group isomorphism $H^1(S^1\times S^1;\Z)\cong \Z\oplus\Z$. To sum up: there is a canonical bijection $\pi_0(\Distr_1(T^2))\cong \Z\times\Z$.
\smallskip\\
Let us rephrase this abstract result in more concrete terms. For each $(k,l)\in\Z\times\Z$, we consider the line distribution $V_{k,l}$ on $T^2 = \R^2/\Z^2$ which assigns to each $[x,y]\in\R^2/\Z^2$ the $1$-dimensional sub vector space of $\R^2 = T_{[x,y]}(T^2)$ that is spanned by the vector $\exp(k\pi ix +l\pi iy)\in\R^2$. This distribution is well-defined and smooth because $\exp(k\pi i(x+1) +l\pi i y) = (-1)^k\exp(k\pi ix +l\pi i y)$ and $\exp(k\pi ix +l\pi i(y+1)) = (-1)^l\exp(k\pi ix +l\pi i y)$ for all $x,y\in\R$. It is easy to see that the distribution $V_{k,l}$ is orientable if and only if the integers $k$ and $l$ are both even.
\smallskip\\
$V_{k,l}$ and $V_{k',l'}$ lie in the same path component of $\pi_0(\Distr_1(T^2))$, i.e. in the same homotopy class in $[T^2,\RP^1]$, if and only if $(k,l)=(k',l')$. Moreover, every homotopy class in $[T^2,\RP^1]$ contains one of the distributions $V_{k,l}$. These facts can be verified either in an elementary way (via $[T^2,\RP^1] = [S^1\times S^1,S^1]$, using only the well-known isomorphism $\pi_1(S^1)\cong\Z$), or by checking from the definitions that the map $\Z\times\Z \to \pi_0(\Distr_1(T^2))$ which we described above in an abstract way is induced by $(k,l)\mapsto V_{k,l}$. Here we do neither of both because the statement is not used anywhere in this thesis, except for illustrative purposes.
\end{example}

For more general manifolds $M$, the computation of $\pi_0(\Distr_q(M))$ is not so easy, not even in the case $q=1$. Let me refer you to Chapter VI in \cite{Whitehead} for relevant techniques and results in obstruction theory (cf.\ in particular Section 6, especially Theorem VI.6.13).

\begin{remark} \label{leafspace}
The above examples possess some properties which make them atypical; in particular, all the distributions we considered are \emph{integrable}, that is, they are tangent to some foliation of the manifold. (Subection \ref{twistsection} contains a discussion of (non)integrability of distributions. Note that all line distributions are integrable.) One can therefore consider the \emph{leaf space} of the corresponding foliation, i.e.\ the set of leaves equipped with the quotient topology of the topology on $M$.
\smallskip\\
To avoid possible confusion, we point out that the question in which path component of $\Distr_q(M)$ a particular integrable distribution $V$ lies is completely unrelated to the shape of the leaves, or the topology of the leaf space, of $V$. For example, consider for $t\in\R$ the distribution $V_t$ on the $2$-torus which assigns to each $[x,y]\in\R^2/\Z^2$ the $1$-dimensional subspace of $\R^2=T_{[x,y]}(\R^2/\Z^2)$ generated by the vector $(1,t)\in\R^2$. For irrational $t$, each leaf of the foliation defined by $V_t$ is dense in $T^2$. For rational $t$, each leaf is a closed submanifold of $T^2$, a so-called torus knot. The winding numbers of such a torus knot (which count how often the knot winds around each factor of $T^2 = S^1\times S^1$ before it closes) are not related to the path component of $\Distr_1(M)$ in which $V_t$ lies, simply because all $V_t$ are contained in the same path component: the map $\R\to\Distr_1(M)$ given by $t\mapsto V_t$ is a path which connects them.
\end{remark}


\chapter[The Riemannian problem]{The Riemannian prescribed scalar curvature problem} \label{D}

\section{Review of some classical results} \label{D1}

The aim of this section is to review some of the work that has been done on the Riemannian prescribed scalar curvature problem, focussing on theorems and proof ideas which can be compared to the results and methods of the present thesis and thus help to put the present work into context. Some of the results, in particular the Kazdan/Warner approximation theorem, are also applied directly in the thesis. Of course, I neither try to give a balanced overview of everything that has been done in the Riemannian case, nor do I care too much about putting events into their proper historical context. A more complete and historically accurate picture can be obtained, for instance, from the overview articles \cite{BerardBergery0}, \cite{RosenbergStolz1994}, \cite{RosenbergStolz2001}, \cite{Stolz1995}, \cite{Stolz2002}, from \S\S I, II in \cite{Kazdan}, and from Chapters 5, 6 in \cite{Aubin1998}.
\medskip\\
Since there are several conventions involving terms like \emph{conformal}, \emph{conformally equivalent}, \emph{pointwise conformal}, we should make clear which one we are going to use in the present thesis:
\begin{definition} \label{confdef}
Let $(M,g)$ be a semi-Riemannian manifold. A semi-Riemannian metric $\tilde{g}$ on $M$ is \emph{conformal to $g$} if and only if there is a function $\lambda\in C^\infty(M,\R_{>0})$ such that $\tilde{g} = \lambda g$. The \emph{conformal class of $g$} is the set of all metrics which are conformal to $g$.
\end{definition}
Note that several authors, including Kazdan and Warner, call conformal metrics \emph{pointwise conformal}; they call metrics $g,\tilde{g}$ \emph{conformal} if there is a function $\lambda\in C^\infty(M,\R_{>0})$ and a diffeomorphism $\varphi\colon M\to M$ such that $\tilde{g} = \lambda\varphi^\ast(g)$. We do not introduce a special term for metrics of this form. If at all, we are only interested in diffeomorphisms contained in $\Diff^0(M)$.

\subsection{The PDE} \label{D11}

Let $M$ be an $n$-manifold. The main approach to finding a Riemannian metric with prescribed scalar curvature $s\in C^\infty(M,\R)$ is to fix a Riemannian metric $g$ on $M$ and try to find a suitable metric in the conformal class of $g$. This approach is equivalent to finding a solution $\kappa\in C^\infty(M,\R_{>0})$ of the quasilinear elliptic PDE
\begin{subequations} \label{yamabe}
\begin{equation} \label{yamabe1}
0 = 2(n-1)\laplace_g(\kappa) -\frac{n(n-1)}{\kappa}\abs{d\kappa}_g^2 +\kappa\scal_g -\frac{1}{\kappa}s \;\;:
\end{equation}
the metric $\kappa^{-2}g$ has scalar curvature $s$ if and only if $\kappa$ solves this equation.
\smallskip\\
Using a substitution $\kappa=U\compose u$, we can get rid of the summand involving $\abs{d\kappa}_g^2$: If $n=2$, we take the diffeomorphism $U\in C^\infty(\R,\R_{>0})$ given by $U(x)=e^{-x}$. In this case, equation \eqref{yamabe1} has a solution $\kappa\in C^\infty(M,\R_{>0})$ if and only if the equation\footnote{When you compare \eqref{yamabe2} with the formulae given in \cite{KazdanWarner2a}, \cite{BerardBergery0}, \cite{Kazdan} etc., note that those refer to the Gaussian curvature $K_g = \frac{1}{2}\scal_g$.}
\begin{equation} \label{yamabe2}
0 = 2\laplace_g(u) -\scal_g +e^{2u}s
\end{equation}
has a solution $u\in C^\infty(M,\R)$. If $n\geq3$, we take the diffeomorphism $U\in C^\infty(\R_{>0},\R_{>0})$ given by $U(x)=x^{-2/(n-2)}$. In that case, equation \eqref{yamabe1} has a solution $\kappa\in C^\infty(M,\R_{>0})$ if and only if
\begin{equation} \label{yamabe3}
0 = {\textstyle\frac{4(n-1)}{n-2}}\laplace_g(u) -u\scal_g +u^{\frac{n+2}{n-2}}s
\end{equation}
\end{subequations}
has a solution $u\in C^\infty(M,\R_{>0})$.
\smallskip\\
However, \eqref{yamabe1} does not always have a solution for a fixed $g$. For instance, for each $g$ on a closed manifold of dimension $\geq3$, there is at most one constant $c\in\set{-1,0,1}$ such that \eqref{yamabe3} has a solution with $s=c$; cf.\ \cite{KazdanWarner1}, Theorem 6.2. (Namely, at most that element of $\set{-1,0,1}$ can occur which is the sign of the smallest eigenvalue of the linear elliptic operator $L_g\colon C^\infty(M,\R) \to C^\infty(M,\R)$ given by $u\mapsto -4(n-1)/(n-2)\laplace_g(u) +\scal_gu$.) In fact, \eqref{yamabe3} has indeed a solution with $s=c$ for this unique constant $c\in\set{-1,0,1}$, according to the solution of the famous \emph{Yamabe problem}, obtained by the work of H.\ Yamabe, N.\ Trudinger, T.\ Aubin, and R.\ Schoen (\cite{Yamabe1960}, \cite{Trudinger1968}, \cite{Aubin1976}, \cite{Schoen1984}; cf.\ also \cite{LeeParker}).
\smallskip\\
Equation \eqref{yamabe1} is \emph{almost} the Riemannian special case of the general (semi-Riemannian) equation \eqref{THEEQUATION} which we will use to solve the pseudo-Riemannian prescribed scalar curvature problem: Our equation \eqref{THEEQUATION} has a solution $f\in C^\infty(M,\R_{>0})$ if and only if \eqref{yamabe1} has a solution $\kappa\in C^\infty(M,\R_{>1})$.
\smallskip\\
However, in the pseudo-Riemannian case, it is possible to substitute away all the quadratic first-order terms in \eqref{THEEQUATION} \emph{only} if $n=2$; and even then, there will still remain nasty linear first-order terms, except if the given metric $g$ has very special properties. (It is precisely metrics with these properties that we will use to solve the plain problem for Lorentzian metrics on $2$-manifolds.)

\subsubsection{Introduction to the following subsections} \label{D1111}

The following theorems are the main results on the Riemannian prescribed scalar curvature problem (cf.\ \cite{KazdanWarner2a}, \cite{KazdanWarner2b}, \cite{KazdanWarner0}, \cite{KazdanWarner1}, \cite{KazdanWarner2}, \cite{KazdanWarner3}). Several proofs have been given for them. There is one proof which does \emph{not} employ the conformal deformation equation \eqref{yamabe}, namely that given in \cite{KazdanWarner2}; it does not seem to have a generalisation to the pseudo-Riemannian case. All other proofs that I know of are based on conformal deformation. Three analytic techniques are used to solve the corresponding elliptic equations \eqref{yamabe2} and \eqref{yamabe3}:
\begin{enumerate}
\item
the method of sub- and supersolutions;
\item
variational methods;
\item
an implicit function plus perturbation method.
\end{enumerate}
These are also the basic techniques that we will apply in this thesis. While (iii) is the most powerful method for the Riemannian problem, (i) is better suited to most cases of the pseudo-Riemannian problem.

\subsection{Closed $2$-manifolds} \label{D12}

\begin{theorem}[Kazdan/Warner] \label{KW2}
Let $M$ be a connected closed $2$-manifold.
\begin{itemize}
\item
If the Euler characteristic of $M$ is positive [resp.\ negative], then the set of functions $\in C^\infty(M,\R)$ which can be represented as the scalar curvature of a Riemannian metric consists exactly of those functions which are positive [negative] somewhere.
\item
If $M$ has Euler characteristic $0$, then the set of functions $\in C^\infty(M,\R)$ which can be represented as the scalar curvature of a Riemannian metric consists exactly of the constant $0$ and of those functions which change their sign (i.e., which are positive somewhere and negative somewhere else).
\end{itemize}
\end{theorem}
\begin{proof}[\textsc{References}]
The case $\chi(M)\leq0$ was proved in \cite{KazdanWarner2a}; by the method of sub- and supersolutions for $\chi(M)<0$ (Theorem 11.8), by variational methods for $\chi(M)=0$ (Theorem 6.3). The general case was solved (reproving the $\chi(M)\leq0$ case) in \cite{KazdanWarner1}, via the implicit function plus perturbation technique (Theorem 5.6).
\end{proof}

\begin{remark}
If $M$ is a connected closed $2$-manifold (not necessarily orientable), the Euler characteristic $\chi(M)\in\Z$ and the total scalar curvature of any Riemannian metric $g$ on $M$ are related by the Gauss/Bonnet theorem (cf.\ e.g.\ \cite{GallotHulinLafontaine}, \S3.111):
\[
4\pi\chi(M) = \int_{(M,g)}\scal_g \;\;.
\]
This implies immediately one half of Theorem \ref{KW2}: \emph{only} the specified functions can be represented as scalar curvatures. The other half of the theorem is the difficult one, of course.
\end{remark}

\subsection{Closed manifolds of dimension $\geq3$} \label{D13}

Before we state the main result, two preparatory theorems deserve to be mentioned separately.
\smallskip\\
Recall that the \emph{total scalar curvature} of a compact Riemannian manifold $(M,g)$ is the number
\[
\int_{(M,g)}\scal_g \;\;.
\]

\begin{theorem}[Aubin; \Eliasson] \label{eliasson}
Every nonempty closed manifold of dimension $\geq3$ admits a Riemannian metric of negative total scalar curvature.
\end{theorem}
\begin{proof}[\textsc{Remarks on the proof}]
A Riemannian metric with negative total scalar curvature can in general not be found in a given conformal class. The presumably easiest way to construct such a metric --- invented by L.\ \BerardBergery\ --- is sketched in \cite{Besse}, \S4.32. The original proof, discovered independently by T.\ Aubin and H.\ I.\ \Eliasson\ (\cite{Aubin1970}, \cite{Eliasson}), used a different approach: from a given Riemannian metric $g$ on $M$ and suitable functions $\lambda,\mu\in C^\infty(M,\R_{>0})$, $\varphi\in C^\infty(M,\R)$, they constructed a metric with negative total scalar curvature in the form
\[
\lambda g +\mu\,d\varphi\otimes d\varphi \;\;.
\]
This approach has a vague similarity to the constructions of the present thesis. A third (somewhat overpowered) proof could employ J.\ Lohkamp's existence results for metrics with negative Ricci curvature.
\end{proof}

Now we come to the main theorem.

\begin{trichotomy} \label{trichotomy}
Let $M$ be a nonempty connected closed manifold of dimension $\geq3$. Then exactly one of the following three alternatives holds:
\begin{enumerate} \renewcommand{\labelenumi}{(\Alph{enumi})}
\item[$(\mathscr{P})$]
Every function $\in C^\infty(M,\R)$ can be represented as the scalar curvature of a Riemannian metric.
\item[$(\mathscr{Z})$]
The set of functions $\in C^\infty(M,\R)$ which can be represented as the scalar curvature of a Riemannian metric consists exactly of the constant $0$ and of those functions which are negative somewhere.
\item[$(\mathscr{N})$]
The set of functions $\in C^\infty(M,\R)$ which can be represented as the scalar curvature of a Riemannian metric consists exactly of those functions which are negative somewhere.
\end{enumerate}
\end{trichotomy}
\begin{proof}[\textsc{References and remarks}]
Most parts of this theorem have been proved by Kazdan and Warner; cf.\ \cite{KazdanWarner1}, Theorem 6.4. The missing piece, namely that the manifold admits a Riemannian metric with constant positive scalar curvature if it admits a Riemannian metric with nonnegative but not identically zero scalar curvature, is provided by the solution of the Yamabe conjecture.
\end{proof}

\subsection{Open manifolds} \label{D14}

\begin{theorem}[Kazdan/Warner] \label{KWopen}
Let $M$ be a connected noncompact manifold of dimension $\geq2$ which is diffeomorphic to an open subset of a compact manifold. Then every function $\in C^\infty(M,\R)$ is the scalar curvature of some Riemannian metric on $M$.
\end{theorem}
\begin{proof}[\textsc{Idea of the proof, references}]
Cf.\ \cite{KazdanWarner2a}, \S6 in \cite{KazdanWarner0}, and Remark 5.7 in \cite{KazdanWarner1}. The idea of the proof is to extend the prescribed function $s\in C^\infty(M,\R)$ to a not necessarily continuous function $\bar{s}\in L^p(\bar{M},\R)$, where $\bar{M}\supset M$ is closed. If this is not possible for $s$ itself, then it is at least possible for some pullback $\varphi^\ast s$, where $\varphi\colon M\to M$ is a diffeomorphism. Then one tries to find a solution $\kappa\in\Sob{p}{2}(\bar{M},\R)$ (where $p$ is chosen so large that $\Sob{p}{2}(\bar{M},\R)\subseteq C^0(\bar{M},\R)$) of the PDE \eqref{yamabe1}. The restriction of this solution to $M$ will be smooth and thus solve the problem.
\end{proof}

\begin{corollary}
Let $n\in\N_{\geq2}$, and let $M$ be a connected compact $n$-manifold with nonempty boundary. Then every function $\in C^\infty(M,\R)$ is the scalar curvature of some Riemannian metric.
\end{corollary}
\Proof
Imbed $M$ into a closed manifold $\bar{M}$, extend the given function $s$ on $M$ to a smooth function $\bar{s}$ on $\bar{M}$, choose a point $x\in\bar{M}\without M$, apply the preceding theorem to the manifold $\tilde{M} \define \bar{M}\without\set{x}$ and the function $\tilde{s} \define \bar{s}\restrict\tilde{M}$ to get a metric $\tilde{g}$ on $\tilde{M}$ with scalar curvature $\tilde{s}$. The metric $\tilde{g}\restrict M$ has scalar curvature $s$.
\end{proof}

One can also say something about arbitrary noncompact manifolds, which are not necessarily imbeddable into a closed manifold. Some sample results:

\begin{theorem}[Gromov]
Every open manifold of dimension $\geq2$ admits a Riemannian metric of positive scalar curvature and a Riemannian metric of negative scalar curvature. (In fact, it even admits Riemannian metrics of positive resp.\ negative \emph{sectional} curvature.)
\end{theorem}
\Proof
This is an easy consequence of Gromov's h-principle theorem (cf.\ Appendix \ref{openhprinciple}) for diff-invariant open partial differential relations on open manifolds; cf.\ \cite{Gromov1969english}, Theorem 4.5.1.
\end{proof}

\begin{theorem}[Bland/Kalka]
Every manifold of dimension $\geq3$ admits a complete Riemannian metric of constant negative scalar curvature.
\end{theorem}
\Proof
Cf.\ the article \cite{BlandKalka1989} by J.\ Bland and M.\ Kalka; or the alternative proof by J.\ Lohkamp: Corollary 5.4 in \cite{Lohkamp1994a}.
\end{proof}

Unfortunately, this thesis contains no similar results for the pseudo-Riemannian problem, because completeness of pseudo-Riemannian metrics is a more subtle issue than completeness of Riemannian metrics; cf.\ Remark \ref{completenessremark}.

\subsection{Topological characterisation of the trichotomy phenomenon} \label{D15}

One would like to tell from topological invariants of a closed $n$-manifold $M$, where $n\geq3$, in which of the three classes $\mathscr{P}$, $\mathscr{Z}$, $\mathscr{N}$ from Theorem \ref{trichotomy} $M$ is contained. Each class is nonempty, in each dimension: $S^n\in\mathscr{P}$, $T^n\in\mathscr{Z}$, and $T^n\connsum T^n\in\mathscr{N}$, for example. Much work has been done on the characterisation of manifolds which admit Riemannian metrics of positive scalar curvature, i.e.\ the characterisation of $\mathscr{P}$. It is not necessary to repeat the well-known results here. The article \cite{RosenbergStolz2001} by J.\ Rosenberg and S.\ Stolz contains a recent overview with many references; cf.\ also \cite{GromovLawson2} by M.\ Gromov and H.\ B.\ Lawson. The class $\mathscr{Z}$ has been investigated by A.\ Futaki (\cite{Futaki1993}) and A.\ Dessai (\cite{Dessai2001}), for instance.
\smallskip\\
The main results of the present thesis show that the pseudo-Riemannian prescribed scalar curvature problem does in most cases not lead to similar topological questions.

\subsection{Real-analytic metrics} \label{D16}

Let us assume that the manifold $M$ is equipped with a real-analytic structure; that's no restriction since every smooth manifold admits such a structure (cf.\ \ref{realanalyticatlas}). If our prescribed function $s$ is real-analytic and there exists a smooth Riemannian metric with scalar curvature $s$, is there also a \emph{real-analytic} Riemannian metric with scalar curvature $s$?
\smallskip\\
The answer deserves to be mentioned here since we will deal with the analogous question in the pseudo-Riemannian case. The result for Riemannian metrics, Theorem \ref{realanalyticRiemann0} below, is easy to deduce from the Kazdan/Warner proofs and some standard (but deep) theorems and is therefore probably well-known, but I do not know a reference. Therefore we give a proof in Appendix \ref{D3}.

\begin{theorem} \label{realanalyticRiemann0}
Let $M$ be an open subset of some closed real-analytic manifold (of any dimension), and let $s\colon M\to\R$ be a real-analytic function which is the scalar curvature of some smooth Riemannian metric on $M$. Then there is a real-analytic Riemannian metric on $M$ with scalar curvature $s$.
\end{theorem}


\section{The Kazdan/Warner approximation theorem} \label{D2}

The following approximation theorem is an essential ingredient in Kazdan and Warner's solution of the Riemannian prescribed scalar curvature problem. We state it in a slightly sharper form than they did: instead of arbitrary diffeomorphisms, we use diffeomorphisms in $\Diff^0(M)$ with compact support. (Recall that a diffeomorphism $\varphi\colon M\to M$ has \emph{compact support} if and only if there is a compact subset $K$ of $M$ such that $\varphi\restrict(M\without K)$ is the identity.) The essential information for us lies not in the compact support, but in the $\Diff^0(M)$ part: we want to apply the theorem later to solve the pseudo-Riemannian homotopy or diffeotopy class problem (on compact manifolds).

\begin{theorem}[the Kazdan/Warner $L^p$ approximation theorem] \label{KWapproximation1}
Let $1\leq p<\infty$, let $M$ be a nonempty connected manifold of dimension $\geq2$ equipped with a Riemannian measure\footnote{A \emph{Riemannian measure} is a measure induced by a Riemannian metric.}, let $g\in L^p(M,\R)$ (where the $L^p$ space is defined with respect to the given measure), and let $f\in C^0(M,\R)\cap L^p(M,\R)$. Then the following statements are equivalent:
\begin{enumerate}
\item
For every open neighbourhood $\mathscr{U}$ of $g$ in $L^p(M,\R)$, there is a diffeomorphism $\varphi\in\Diff^0(M)$ with compact support such that $f\compose\varphi \in \mathscr{U}$;
\item
$\inf(f) \leq g \leq \sup(f)$ (where $\inf(f),\sup(f)\in[-\infty,\infty]$).
\end{enumerate}
\end{theorem}
\Proof
This is Theorem 2.1 in \cite{KazdanWarner1}, except for the claim that in the implication \textit{(ii)}$\implies$\textit{(i)}, we can not only find any diffeomorphism $\varphi\in\Diff(M)$ with $f\compose\varphi \in \mathscr{U}$, but even such a diffeomorphism in $\Diff^0(M)$, and with compact support. This claim is true since the proof of Kazdan and Warner remains valid for our stronger statement, word by word. We just have to take the following facts into account:
\begin{enumerate}
\item
Let $\tilde{\eps}\in\R_{>0}$, and let $x_1,x_2,y_1,y_2\in M$ with $x_1\neq x_2$ and $y_1\neq y_2$. Then there is a compactly supported $\varphi_1\in\Diff^0(M)$ with $\varphi_1(y_i)=x_i$ for $i\in\set{1,2}$, such that $\int_M\abs{f\compose\varphi_1-f}^p\leq\tilde{\eps}^p$.
\item
Let $M_1$ be a compact subset of $M$, let $k\in\N$, let $V_1,\dots,V_k$ be disjoint open subsets of $M_1$, and let $U_1,\dots,U_k$ be open subsets of $M_1$ with disjoint closures. Then there is a compactly supported $\varphi_2\in\Diff^0(M)$ whose restriction to $M\without M_1$ is the identity, such that $\varphi_2(U_i)\subseteq V_i$ for all $i\in\set{1,\dots,k}$.
\end{enumerate}
(Fact (ii) is applied twice in the Kazdan/Warner proof: once using precisely the above notation, once with $M_1$ replaced by $\bar{O_1}$, with $V_i$ replaced by $U_i$, and with $U_i$ replaced by $O_1\cap\Delta_i$.)
\smallskip\\
Both facts follow from elementary differential topology (so Kazdan and Warner did not bother to write down an argument, which would show immediately that all the diffeomorphisms are taken from $\Diff^0(M)$). Here is a proof sketch: For $j\in\set{0,1}$, the desired diffeotopy $\Phi_j\colon[0,1]\times M\to M$ with $\varphi_j=\Phi_j(1,.)$ can be constructed as the flow of a suitable vector field with compact support.\footnote{Cf.\ e.g.\ \S4.1 in \cite{Conlon}, in particular Theorem 4.1.11.}
\smallskip\\
For (i), we choose regular injective paths $w_i\colon[0,1]\to M$, $i\in\set{1,2}$, from $y_i$ to $x_i$, such that the images of $w_1$ and $w_2$ are disjoint (this is possible because $\dim M\geq2$). We choose neighbourhoods $N_i$ of $\im(w_i)$ which are disjoint and, moreover, so small that $2(\volume(N_1)^{1/p}+\volume(N_2)^{1/p})\max\set{\abs{f(x)} \suchthat x\in N_1\cup N_2} \leq \tilde{\eps}$. Now any vector field $X$ on $M$ with support in $N_1\cup N_2$ whose restriction to $\im(w_i)$ is given by $X\compose w_i = w_i'$ for $i\in\set{0,1}$ generates a suitable flow $\Phi_1$.
\smallskip\\
For (ii), we can choose smooth imbeddings $w_1,\dots,w_k \colon [0,1.5]\times D^{n-1}_1 \to M_1$ (where $n=\dim M$, and $D^{n-1}_r$ is the $0$-centered closed disk in $\R^{n-1}$ with radius $r$) such that the following conditions hold: (A) The images of $w_1,\dots,w_k$ are disjoint. (B) Each $U_i$ is a subset of $w_i([0.1,0.4]\times D^{n-1}_{0.9})$. (C) Each $V_i$ is a superset of $w_i([1,1.5]\times D^{n-1}_1)$. (These imbeddings are constructed by first choosing suitable imbedded paths which map $[0,1.5]$ to the interior of $M_1$, and then suitable tubular neighbourhoods of these imbedded paths. The corresponding bundles over $[0,1.5]$ are trivial since the base space is contractible.)
\smallskip\\
For each $x\in D^{n-1}_1$ and $i\in\set{1,\dots,k}$, let $w_i(x)\colon [0.1,1.4] \to M_1$ be the map given by $t\mapsto w_i(t,x)$. There is a vector field $X$ on $M$ whose support is a subset of $\bigcup_{i=1}^k\im(w_i)$, such that $X\compose w_i(x) = w_i(x)'$ for all $x\in D^{n-1}_{0.9}$ and $i\in\set{1,\dots,k}$. It generates a flow $\Phi_2$ with the desired properties.
\end{proof}

By dropping the compact support condition, we get a real-analytic analogue, Theorem \ref{realanalyticapproximation1} below. But first we need a lemma.

\begin{lemma} \label{KWapproxlemma}
Let $1\leq p<\infty$, let $M$ be a real-analytic manifold equipped with a measure, let $u\colon M\to\R$ be a not necessarily continuous function, let $\psi\colon M\to M$ be a smooth diffeomorphism such that $u\,\mathord{\compose}\,\psi\in L^p(M,\R)$, and let $\mathscr{U}$ be any neighbourhood of $u\compose\psi$ in $L^p(M,\R)$. Then there is a real-analytic diffeomorphism $\varphi\colon M\to M$ with $u\compose\varphi\in\mathscr{U}$; if $\psi\in\Diff^0(M)$, then we can choose $\varphi\in\Diff^0(M)$, too.
\end{lemma}
\Proof
There is an $\eps>0$ such that $\mathscr{U}$ contains an $L^p$-ball of radius $\eps$ around $u\compose\psi$. Since continuous functions with compact support are dense in $L^p(M,\R)$, the $L^p$-ball of radius $\frac{\eps}{2}$ around $u\compose\psi$ contains a continuous function, and by composing this function with $\psi^{-1}$ from the right, we get a $\tilde{u}\in C^0(M,\R)$ with $\norm{u\compose\psi-\tilde{u}\compose\psi}_{L^p}\leq\frac{\eps}{2}$.
\smallskip\\
We choose a (rapidly decreasing) continuous $\eta\in C^0(M,\R_{>0})$ with $\norm{\eta}_{L^p} <\frac{\eps}{2}$. The set $\mathscr{U}'$ of all $h\in C^0(M,\R)$ with $\abs{h-\tilde{u}\compose\psi}<\eta$ is a neighbourhood of $\tilde{u}\compose\psi$ with respect to the fine $C^0$-topology on $C^0(M,\R)$. Observe that $\norm{h-u\compose\psi}_{L^p}<\eps$ for all $h\in\mathscr{U}'$; hence $\mathscr{U}'\subseteq\mathscr{U}$.
\smallskip\\
Since the map $C^0(M,M)\to C^0(M,\R)$ given by $\varphi\mapsto u\compose\varphi$ is continuous with respect to the fine $C^0$-topologies on both sides (cf.\ \S2.4, Exercise 10(b) in \cite{Hirsch}), there is a fine $C^0$-neighbourhood $\mathscr{V}\subseteq C^0(M,M)$ of $\psi$ such that $u\compose\varphi \in \mathscr{U}'$ for all $\varphi\in\mathscr{V}$.
\smallskip\\
The set of diffeomorphisms is fine $C^\infty$-open in $C^\infty(M,M)$ (cf.\ \cite{Hirsch}, Theorem 2.1.7), and so is the set $\Diff^0(M)$. Thus there is a fine $C^\infty$-neighbourhood $\mathscr{V}'\subseteq\mathscr{V}\cap C^\infty(M,M)$ of $\psi$ all of whose elements are contained in $\Diff(M)$; if $\psi\in\Diff^0(M)$, then there is even such a neighbourhood $\subseteq\Diff^0(M)$. By the Grauert/Remmert/Morrey denseness theorem (cf.\ e.g.\ Theorem 2.5.1 in \cite{Hirsch}), $\mathscr{V}'$ contains a real-analytic map. This is the desired diffeomorphism $\varphi\in\Diff^0(M)$.
\end{proof}

\begin{theorem}[real-analytic Kazdan/Warner approximation] \label{realanalyticapproximation1}
Let $1\leq p<\infty$, let $M$ be a nonempty connected real-analytic manifold of dimension $\geq2$ equipped with a Riemannian measure, let $g\in L^p(M,\R)$, and let $f\in C^0(M,\R)\cap L^p(M,\R)$. Then the following statements are equivalent:
\begin{enumerate}
\item
For every open neighbourhood $\mathscr{U}$ of $g$ in $L^p(M,\R)$, there is a real-analytic diffeomorphism $\varphi\in\Diff^0(M)$ such that $f\compose\varphi \in \mathscr{U}$;
\item
$\inf(f) \leq g \leq \sup(f)$ (where $\inf(f),\sup(f)\in[-\infty,\infty]$).
\end{enumerate}
\end{theorem}
\Proof
The implication \textit{(i)}$\implies$\textit{(ii)} follows from Theorem 2.1 in \cite{KazdanWarner1}. In order to prove \textit{(ii)}$\implies$\textit{(i)}, we take any smooth diffeomorphism $\psi\in\Diff^0(M)$ such that $f\compose\psi\in\mathscr{U}$; this is possible by Theorem \ref{KWapproximation1}. Now we apply the preceding lemma and are done.
\end{proof}

We state another consequence of Lemma \ref{KWapproxlemma}, which will be applied in the proof of Theorem \ref{realanalyticRiemann0}.

\begin{corollary} \label{KWapprox2}
Let $1\leq p<\infty$, let $M$ be a connected real-analytic manifold equipped with a Riemannian measure, let $f\in L^p(M,\R)$ be zero on some nonempty open subset of $M$, and let $\eps\in\R_{>0}$. Then there is a real-analytic diffeomorphism $\varphi\colon M\to M$ such that $\norm{f\compose\varphi}_{L^p} <\eps$.
\end{corollary}
\Proof
Proposition 2.6 in \cite{KazdanWarner2a} says that a \emph{smooth} diffeomorphism $\psi\colon M\to M$ with $\norm{f\compose\psi}_{L^p} <\eps$ exists; now we apply Lemma \ref{KWapproxlemma}.
\end{proof}

\emph{Remark.} In the statement of the corollary, $\varphi$ can even be chosen from $\Diff^0(M)$. The diffeomorphism in the smooth analogon, i.e.\ Proposition 2.6 in \cite{KazdanWarner2a}, can be chosen from $\Diff^0(M)$, and with compact support.


\section[Real-analyticity]{The Riemannian real-analytic prescribed scalar curvature problem} \label{D3}

The aim of this chapter is to prove Theorem \ref{realanalyticRiemann0} (whose statement we repeat below as Theorem \ref{realanalyticRiemann1}, for convenience).
\medskip\\
We need the following theorem to deal with some special case in the proof of \ref{realanalyticRiemann1}. Recall that a manifold chart $(x_1,\dots,x_n)\colon U\to\R^n$ defined on some open subset $U$ of an $n$-manifold $M$ is called \emph{($g$-)harmonic} with respect to a Riemannian metric $g$ on $M$ if and only if $\laplace_g(x_i)=0$ for all $i\in\set{1,\dots,n}$.

\begin{theorem}[DeTurck/Kazdan] \label{deturckkazdan}
Let $(M,g)$ be a smooth\footnote{In fact, it would suffice to assume that $M$ and $g$ are $C^2$.} Einstein Riemannian manifold of dimension $n\geq3$. Then the set of all $g$-harmonic charts on $M$ is a real-analytic atlas on $M$ which is compatible with the given smooth atlas. The metric $g$ is real-analytic with respect to this atlas.
\end{theorem}
\Proof
Smooth harmonic charts exist around each point of $M$; cf.\ e.g.\ Lemma 1.2 in \cite{DeTurckKazdan1981}. Let $(x_1,\dots,x_n)\text{:}$ $U\to\tilde{U}$ and $y=(y_1,\dots,y_n)\colon V\to\tilde{V}$ be $g$-harmonic charts onto open subsets $\tilde{U},\tilde{V}$ of $\R^n$. We have to show that the smooth coordinate change $x_i\compose y^{-1} \colon y(U\cap V) \to \R$ is real-analytic for each $i\in\set{1,\dots,n}$ (we suppress restrictions of maps in our notation).
\smallskip\\
Since $g$ is Einstein, it is real-analytic in the $g$-harmonic coordinates $(y_1,\dots,y_n)$, by Theorem 5.2 in \cite{DeTurckKazdan1981}; that is, the metric $(y^{-1})^\ast g$ on $\tilde{V}\subseteq\R^n$ is real-analytic. The function $(y^{-1})^\ast x_i \equiv x_i\compose y^{-1}\in C^\infty(y(U\cap V),\R)$ satisfies the elliptic equation $\laplace_{(y^{-1})^\ast g}((y^{-1})^\ast x_i) = (y^{-1})^\ast(\laplace_gx_i) = 0$ with real-analytic coefficients and is thus real-analytic, by Theorem \ref{hopf}. This proves that the $g$-harmonic charts form a real-analytic atlas. We have already used the DeTurck/Kazdan theorem which says that $g$ is real-analytic with respect to each chart in this atlas.
\end{proof}

\begin{corollary} \label{deturckkazdan2}
Let $M$ be a real-analytic manifold of dimension $\geq3$. If $g$ is a smooth Einstein metric on $M$, then there is a smooth diffeomorphism $\varphi\colon M\to M$ such that $\varphi^\ast g$ is real-analytic. In particular, if $M$ admits a smooth Einstein metric with constant scalar curvature $\lambda$, then it admits also a real-analytic Einstein metric with constant scalar curvature $\lambda$.
\end{corollary}
\Proof
Let $\mathscr{A}$ be the given real-analytic structure on the smooth manifold $M$. By Theorem \ref{deturckkazdan}, there exists another real-analytic structure $\mathscr{A}'$ on the smooth manifold $M$ such that the smooth Einstein metric $g$ is real-analytic with respect to $\mathscr{A}'$. By the uniqueness of real-analytic structures compatible to a given smooth structure (cf.\ Theorem \ref{realanalyticatlas}), there is a smooth diffeomorphism $\varphi\in\Diff(M)$ which is a real-analytic diffeomorphism $(M,\mathscr{A})\to(M,\mathscr{A}')$. The metric $\varphi^\ast g$ is thus real-analytic with respect to $\mathscr{A}$. It is of course an Einstein metric with the same constant scalar curvature as $g$.
\end{proof}

Now we can prove Theorem \ref{realanalyticRiemann0}:

\begin{theorem} \label{realanalyticRiemann1}
Let $M$ be an open subset of some closed real-analytic manifold, and let $s\colon M\to\R$ be a real-analytic function which is the scalar curvature of some smooth Riemannian metric on $M$. Then there is a real-analytic Riemannian metric on $M$ with scalar curvature $s$.
\end{theorem}
\Proof
The case $\dim M < 2$ is trivial.
\smallskip\\
\emph{First case: $M$ is a closed $2$-manifold.} Then there is a real-analytic Riemannian metric on $M$ with constant (scalar) curvature, as we will prove now.
\smallskip\\
If the Euler characteristic $\chi(M)$ is nonnegative, i.e. $M\in\set{S^2,\RP^2,T^2,\klein}$, then the standard metric $\tilde{g}$ on $M$ is real-analytic with constant scalar curvature with respect to the standard real-analytic structure $\tilde{\mathscr{A}}$ on $M$. Since there is a real-analytic diffeomorphism $\varphi$ from our given real-analytic structure to $\tilde{\mathscr{A}}$ (cf.\ \ref{realanalyticatlas}), there exists a constant curvature metric on $M$ which is real-analytic with respect to our given structure, namely the metric $\varphi^\ast\tilde{g}$.
\smallskip\\
If $\chi(M)<0$, then we choose any real-analytic metric $\tilde{g}$ on $M$; this is always possible. The conformal class of $\tilde{g}$ contains a metric with negative constant scalar curvature; cf.\ e.g.\ Theorem 1 in \cite{Berger1971}. This metric is real-analytic since it has the form $e^{2u}\tilde{g}$, where $u\in C^\infty(M,\R)$ solves an elliptic equation (namely \eqref{yamabe2}) with real-analytic coefficients and is thus real-analytic; cf.\ Theorem \ref{hopf}.
\smallskip\\
Thus there is a real-analytic Riemannian metric on $M$ with constant (scalar) curvature, as claimed.
\smallskip\\
Now let $g$ be a real-analytic Riemannian metric on $M$ with constant scalar curvature $2k$, and let $2K\colon M\to\R$ be a real-analytic function which satisfies the necessary and sufficient condition for being the scalar curvature of some smooth Riemannian metric. We consider the proof of Theorem 5.1 in \cite{KazdanWarner1}. The operator $T_1\colon \Sob{p}{2}(M,\R)\to L^p(M,\R)$ there, which is given by $u\mapsto -e^{-2u}(\laplace_gu-k)$, has real-analytic coefficients. Since the diffeomorphism $\varphi\in\Diff(M)$ there can be chosen real-analytic (cf.\ our Theorem \ref{realanalyticapproximation1}), we can assume that $K\compose\varphi$ is real-analytic. By Theorem \ref{hopf}, the locally unique function $u$ with $T_1(u)=K\compose\varphi$ is real-analytic. The proof of Theorem 5.1 in \cite{KazdanWarner1} shows therefore that a real-analytic Riemannian metric with scalar curvature $2K$ exists.
\medskip\\
\emph{Second case: $M$ is a closed $n$-manifold, where $n\geq3$.} Let $s\colon M\to\R$ be real-analytic.
\smallskip\\
By Theorem \ref{eliasson}, there is a Riemannian metric on $M$ with negative total scalar curvature. Since this metric has a $C^2$ neighbourhood consisting of metrics with negative total scalar curvature, we can find a real-analytic Riemannian metric $\tilde{g}$ on $M$ with negative total scalar curvature; cf.\ \ref{analyticapproximation}. There is a smooth Riemannian metric $g$, conformal to $\tilde{g}$, with constant scalar curvature $-1$; cf.\ Theorem 4.1 in \cite{KazdanWarner0}. This metric is even real-analytic by Theorem \ref{hopf}, because it has the form $u^{4/(n-2)}\tilde{g}$, where $u\in C^\infty(M,\R_{>0})$ solves the elliptic equation \eqref{yamabe3}.
\smallskip\\
Lemma 2.15 in \cite{KazdanWarner0} yields an $L^p$-neighbourhood (where $p>n$) $\mathscr{U}$ of $-1\in L^p(M,\R)$ such that every real-analytic $\tilde{s}\in\mathscr{U}$ is the scalar curvature of some smooth Riemannian metric which is conformal to $g$ and thus (by the argument above involving \eqref{yamabe3} and \ref{hopf}) real-analytic. If $s$ is somewhere negative, then there exist a constant $\alpha\in\R_{>0}$ and a real-analytic diffeomorphism $\varphi\in\Diff(M)$ such that $\alpha s\compose\varphi \in \mathscr{U}$: choose $\alpha$ such that $\inf(\alpha s) \leq -1 \leq \sup(\alpha s)$ and use \ref{realanalyticapproximation1}. Hence $\alpha s\compose\varphi$ is the scalar curvature of some real-analytic Riemannian metric $g'$, so $s$ is the scalar curvature of the real-analytic metric $\alpha(\varphi^{-1})^\ast g'$.
\smallskip\\
This completes the proof if $M$ lies in the class $\mathscr{N}$ of the trichotomy theorem. Let us now assume that $M$ belongs to the class $\mathscr{P}$. Then the set of Riemannian metrics $g$ on $M$ with $\lambda_1(g)>0$ is nonempty and $C^2$-open and thus contains a real-analytic metric $g_1$; here $\lambda_1(g)\in\R$ denotes the smallest eigenvalue of the operator $L_g$ mentioned in Section \ref{D11}. There is also a real-analytic metric $g_0$ with $\lambda_1(g_0)<0$, as we proved above (negative total scalar curvature implies $\lambda_1(g_0)<0$). For some $t\in\oointerval{0}{1}$, the real-analytic metric $g_t\define(1-t)g_0+tg_1$ has $\lambda_1(g_t)=0$ (cf.\ the proof of Theorem 3.9 in \cite{KazdanWarner0}).
\smallskip\\
By Proposition 5.1 in \cite{KazdanWarner0}, there is a metric, conformal to $g_t$, with constant scalar curvature $0$. This metric is real-analytic by the standard argument from above. To finish the case of class $\mathscr{P}$, it remains to prove that $s$ is the scalar curvature of some analytic Riemannian metric if it is somewhere positive. The conformal class of $g_1$ contains a Riemannian metric $g$ of constant positive scalar curvature, by the Aubin/Schoen solution of the Yamabe problem. By our standard argument, $g$ is real-analytic. Now we go through the proof of Lemma 6.1 in \cite{KazdanWarner1}, which is analogous to the proof of Theorem 5.1 there that we considered above. We see in the same way as before (employing \ref{realanalyticapproximation1} and \ref{implicitfunction}) that $s$ is the scalar curvature of some real-analytic Riemannian metric.
\smallskip\\
It remains to consider the case where $M$ belongs to the class $\mathscr{Z}$. Since we have already dealt with all functions which are somewhere negative, we just have to realise the constant $0$ as the scalar curvature of some real-analytic metric. In order to do this, we choose a smooth metric $g$ with scalar curvature $0$. This metric is even Ricci-flat; cf.\ e.g.\ Lemma 5.2 in \cite{KazdanWarner3}. By Corollary \ref{deturckkazdan2}, there exists a real-analytic Ricci-flat metric on $M$. This completes the proof of the closed case.
\medskip\\
\emph{Third case: $M$ is an open manifold of dimension $\geq2$.} We go through the proof at the beginning of \S6 in \cite{KazdanWarner0} (or the proof sketched in Remark 5.7 in \cite{KazdanWarner1}). The smooth diffeomorphism $\varphi\in\Diff(M)$ there (with $s\compose\varphi\in L^p(M)$) can be chosen real-analytic by Lemma \ref{KWapproxlemma}. From our Corollary \ref{KWapprox2} (applied to a noncontinuous function), we deduce that the smooth diffeomorphism $\psi\in\Diff(M_1)$ in the Kazdan/Warner proof can be replaced by a real-analytic one, too. The rest of the proof consists of our standard arguments from above.
\smallskip\\
Now we have checked all cases, so the proof is complete.
\end{proof}


\end{appendix}


\cleardoublepage

\begin{thebibliography}{100}

\bibitem{Alty1995}
{\sc L.~J. Alty}, {\em The generalized {G}auss-{B}onnet-{C}hern theorem}, J.
  Math. Phys., 36 (1995), pp.~3094--3105.

\bibitem{Aubin1970}
{\sc T.~Aubin}, {\em M\'etriques riemanniennes et courbure}, J. Differential
  Geometry, 4 (1970), pp.~383--424.

\bibitem{Aubin1976}
\leavevmode\vrule height 2pt depth -1.6pt width 23pt, {\em \'{E}quations
  diff\'erentielles non lin\'eaires et probl\`eme de {Y}amabe concernant la
  courbure scalaire}, J. Math. Pures Appl. (9), 55 (1976), pp.~269--296.

\bibitem{Aubin1982}
\leavevmode\vrule height 2pt depth -1.6pt width 23pt, {\em Nonlinear analysis
  on manifolds. {M}onge-{A}mp\`ere equations}, vol.~252 of Grundlehren der
  Ma\-the\-matischen Wissenschaften [Fundamental Principles of Mathematical
  Sciences], Springer-Verlag, New York, 1982.

\bibitem{Aubin1998}
\leavevmode\vrule height 2pt depth -1.6pt width 23pt, {\em Some nonlinear
  problems in {R}iemannian geometry}, Springer Monographs in Mathematics,
  Springer-Verlag, Berlin, 1998.

\bibitem{Avez1962}
{\sc A.~Avez}, {\em Formule de {G}auss-{B}onnet-{C}hern en m\'etrique de
  signature quelconque.}, C. R. Acad. Sci. Paris, 255 (1962), pp.~2049--2051.

\bibitem{Baum}
{\sc H.~Baum}, {\em Spin-{S}trukturen und {D}irac-{O}peratoren \"uber
  pseudoriemannschen {M}annigfaltigkeiten}, vol.~41 of Teubner-Texte zur
  Mathematik [Teubner Texts in Mathematics], BSB B. G. Teubner
  Verlagsgesellschaft, Leipzig, 1981.

\bibitem{BeemEhrlichEasley}
{\sc J.~K. Beem, P.~E. Ehrlich, and K.~L. Easley}, {\em Global {L}orentzian
  geometry}, vol.~202 of Monographs and Textbooks in Pure and Applied
  Mathematics, Marcel Dekker Inc., New York, second~ed., 1996.

\bibitem{BerardBergery0}
{\sc L.~B{\'e}rard~Bergery}, {\em La courbure scalaire des vari\'et\'es
  riemanniennes}, in Bourbaki Seminar, Vol. 1979/80, vol.~842 of Lecture Notes
  in Math., Springer, Berlin, 1981, pp.~225--245.

\bibitem{Berger1971}
{\sc M.~S. Berger}, {\em Riemannian structures of prescribed {G}aussian
  curvature for compact {$2$}-manifolds}, J. Differential Geometry, 5 (1971),
  pp.~325--332.

\bibitem{Besse}
{\sc A.~L. Besse}, {\em Einstein manifolds}, vol.~10 of Ergebnisse der
  Mathematik und ihrer Grenzgebiete (3) [Results in Mathematics and Related
  Areas (3)], Springer-Verlag, Berlin, 1987.

\bibitem{BirmanNomizu}
{\sc G.~S. Birman and K.~Nomizu}, {\em The {G}auss-{B}onnet theorem for
  {$2$}-dimensional spacetimes}, Michigan Math. J., 31 (1984), pp.~77--81.

\bibitem{BlandKalka1989}
{\sc J.~Bland and M.~Kalka}, {\em Negative scalar curvature metrics on
  noncompact manifolds}, Trans. Amer. Math. Soc., 316 (1989), pp.~433--446.

\bibitem{Bourguignon1975}
{\sc J.-P. Bourguignon}, {\em Une stratification de l'espace des structures
  riemanniennes}, Compositio Math., 30 (1975), pp.~1--41.

\bibitem{Bredon}
{\sc G.~E. Bredon}, {\em Topology and geometry}, vol.~139 of Graduate Texts in
  Mathematics, Springer-Verlag, New York, 1997.
\newblock Corrected third printing of the 1993 original.

\bibitem{ChoquetBruhatLeray}
{\sc Y.~Choquet-Bruhat and J.~Leray}, {\em Sur le probl\`eme de {D}irichlet,
  quasilin\'eaire, d'ordre {$2$}}, C. R. Acad. Sci. Paris S\'er. A-B, 274
  (1972), pp.~A81--A85.

\bibitem{Conlon}
{\sc L.~Conlon}, {\em Differentiable manifolds}, Birkh\"auser Advanced Texts:
  Basler Lehrb\"ucher. [Birkh\"auser Advanced Texts: Basel Textbooks],
  Birkh\"auser Boston Inc., Boston, MA, second~ed., 2001.

\bibitem{Dacorogna}
{\sc B.~Dacorogna}, {\em Direct methods in the calculus of variations}, vol.~78
  of Applied Mathematical Sciences, Springer-Verlag, Berlin, 1989.

\bibitem{Dessai2001}
{\sc A.~Dessai}, {\em On the topology of scalar-flat manifolds}, Bull. London
  Math. Soc., 33 (2001), pp.~203--209.

\bibitem{DeTurck1983}
{\sc D.~M. DeTurck}, {\em The {C}auchy problem for {L}orentz metrics with
  prescribed {R}icci curvature}, Compositio Math., 48 (1983), pp.~327--349.

\bibitem{DeTurckKazdan1981}
{\sc D.~M. DeTurck and J.~L. Kazdan}, {\em Some regularity theorems in
  {R}iemannian geometry}, Ann. Sci. \'Ecole Norm. Sup. (4), 14 (1981),
  pp.~249--260.

\bibitem{Donaldson1987c}
{\sc S.~K. Donaldson}, {\em The orientation of {Y}ang-{M}ills moduli spaces and
  {$4$}-manifold topology}, J. Differential Geom., 26 (1987), pp.~397--428.

\bibitem{Eliashberg1989}
{\sc Y.~Eliashberg}, {\em Classification of overtwisted contact structures on
  {$3$}-manifolds}, Invent. Math., 98 (1989), pp.~623--637.

\bibitem{EliashbergMishachev}
{\sc Y.~Eliashberg and N.~Mishachev}, {\em Introduction to the
  {$h$}-principle}, vol.~48 of Graduate Studies in Mathematics, American
  Mathematical Society, Providence, RI, 2002.

\bibitem{EliashbergThurston}
{\sc Y.~M. Eliashberg and W.~P. Thurston}, {\em Confoliations}, vol.~13 of
  University Lecture Series, American Mathematical Society, Providence, RI,
  1998.

\bibitem{Eliasson}
{\sc H.~I. El{\'{\i}}asson}, {\em On variations of metrics}, Math. Scand., 29
  (1971), pp.~317--327 (1972).

\bibitem{FischerMarsden1975}
{\sc A.~E. Fischer and J.~E. Marsden}, {\em Linearization stability of
  nonlinear partial differential equations}, in Differential geometry (Proc.
  Sympos. Pure Math., Vol. XXVII, Part 2, Stanford Univ., Stanford, Calif.,
  1973), Amer. Math. Soc., Providence, R.I., 1975, pp.~219--263.

\bibitem{Futaki1993}
{\sc A.~Futaki}, {\em Scalar-flat closed manifolds not admitting positive
  scalar curvature metrics}, Invent. Math., 112 (1993), pp.~23--29.

\bibitem{GallotHulinLafontaine}
{\sc S.~Gallot, D.~Hulin, and J.~Lafontaine}, {\em Riemannian geometry},
  Universitext, Springer-Verlag, Berlin, second~ed., 1990.

\bibitem{Geiges2001}
{\sc H.~Geiges}, {\em Contact topology in dimension greater than three}, in
  European Congress of Mathematics, Vol. II (Barcelona, 2000), vol.~202 of
  Progr. Math., Birkh\"auser, Basel, 2001, pp.~535--545.

\bibitem{Geigessurvey}
\leavevmode\vrule height 2pt depth -1.6pt width 23pt, {\em Contact geometry},
  in Handbook of differential geometry, F.~J.~E. Dillen and L.~C.~A.
  Verstraelen, eds., vol.~2, 2003, pp.~1--86.
\newblock to appear.

\bibitem{Geigeshprinciple}
\leavevmode\vrule height 2pt depth -1.6pt width 23pt, {\em {$h$}-principles and
  flexibility in geometry}, Mem. Amer. Math. Soc., 164 (2003), pp.~viii+58.

\bibitem{Geroch1967}
{\sc R.~P. Geroch}, {\em Topology in general relativity}, J. Mathematical
  Phys., 8 (1967), pp.~782--786.

\bibitem{GilbargTrudinger}
{\sc D.~Gilbarg and N.~S. Trudinger}, {\em Elliptic partial differential
  equations of second order}, Classics in Mathematics, Springer-Verlag, Berlin,
  2001.
\newblock Reprint of the 1998 edition.

\bibitem{Ginzburg1992}
{\sc V.~L. Ginzburg}, {\em Calculation of contact and symplectic cobordism
  groups}, Topology, 31 (1992), pp.~767--773.

\bibitem{GompfStipsicz}
{\sc R.~E. Gompf and A.~I. Stipsicz}, {\em {$4$}-manifolds and {K}irby
  calculus}, vol.~20 of Graduate Studies in Mathematics, American Mathematical
  Society, Providence, RI, 1999.

\bibitem{Gromov1969english}
{\sc M.~Gromov}, {\em Stable mappings of foliations into manifolds}, Math.
  USSR, Izv., 3 (1969), pp.~671--694.

\bibitem{GromovPDR}
\leavevmode\vrule height 2pt depth -1.6pt width 23pt, {\em Partial differential
  relations}, vol.~9 of Ergebnisse der Mathematik und ihrer Grenzgebiete (3)
  [Results in Mathematics and Related Areas (3)], Springer-Verlag, Berlin,
  1986.

\bibitem{GromovLawson2}
{\sc M.~Gromov and H.~B. Lawson, Jr.}, {\em Positive scalar curvature and the
  {D}irac operator on complete {R}iemannian manifolds}, Inst. Hautes \'Etudes
  Sci. Publ. Math.,  (1983), pp.~83--196 (1984).

\bibitem{HawkingEllis}
{\sc S.~W. Hawking and G.~F.~R. Ellis}, {\em The large scale structure of
  space-time}, Cambridge University Press, London, 1973.
\newblock Cambridge Monographs on Mathematical Physics, No. 1.

\bibitem{Hebey}
{\sc E.~Hebey}, {\em Sobolev spaces on {R}iemannian manifolds}, vol.~1635 of
  Lecture Notes in Mathematics, Springer-Verlag, Berlin, 1996.

\bibitem{Hirsch}
{\sc M.~W. Hirsch}, {\em Differential topology}, vol.~33 of Graduate Texts in
  Mathematics, Springer-Verlag, New York, 1994.
\newblock Corrected reprint of the 1976 original.

\bibitem{HirzebruchHopf}
{\sc F.~Hirzebruch and H.~Hopf}, {\em Felder von {F}l\"achenelementen in
  4-dimensionalen {M}annigfaltigkeiten}, Math. Ann., 136 (1958), pp.~156--172.

\bibitem{Hopf1931}
{\sc E.~Hopf}, {\em {\"U}ber den funktionalen, insbesondere den analytischen
  {C}harakter der {L}\"osungen elliptischer {D}ifferentialgleichungen zweiter
  {O}rdnung}, Math. Z., 34 (1931), pp.~194--233.
\newblock Available in PDF format from {\tt http://gdz.sub.uni-goettingen.de}.

\bibitem{Husemoller}
{\sc D.~Husemoller}, {\em Fibre bundles}, vol.~20 of Graduate Texts in
  Mathematics, Springer-Verlag, New York, third~ed., 1994.

\bibitem{IrelandRosen}
{\sc K.~Ireland and M.~Rosen}, {\em A classical introduction to modern number
  theory}, vol.~84 of Graduate Texts in Mathematics, Springer-Verlag, New York,
  second~ed., 1990.

\bibitem{Kazdan}
{\sc J.~L. Kazdan}, {\em Prescribing the curvature of a {R}iemannian manifold},
  vol.~57 of CBMS Regional Conference Series in Mathematics, published for the
  Conference Board of the Mathematical Sciences, Washington, DC, 1985.

\bibitem{KazdanKramer}
{\sc J.~L. Kazdan and R.~J. Kramer}, {\em Invariant criteria for existence of
  solutions to second-order quasilinear elliptic equations}, Comm. Pure Appl.
  Math., 31 (1978), pp.~619--645.

\bibitem{KazdanWarner2a}
{\sc J.~L. Kazdan and F.~W. Warner}, {\em Curvature functions for compact
  {$2$}-manifolds}, Ann. of Math. (2), 99 (1974), pp.~14--47.

\bibitem{KazdanWarner2b}
\leavevmode\vrule height 2pt depth -1.6pt width 23pt, {\em Curvature functions
  for open {$2$}-manifolds}, Ann. of Math. (2), 99 (1974), pp.~203--219.

\bibitem{KazdanWarner0}
\leavevmode\vrule height 2pt depth -1.6pt width 23pt, {\em Scalar curvature and
  conformal deformation of {R}iemannian structure}, J. Differential Geometry,
  10 ({\dummysort{a}}1975), pp.~113--134.

\bibitem{KazdanWarner1}
\leavevmode\vrule height 2pt depth -1.6pt width 23pt, {\em Existence and
  conformal deformation of metrics with prescribed {G}aussian and scalar
  curvatures}, Ann. of Math. (2), 101 ({\dummysort{b}}1975), pp.~317--331.

\bibitem{KazdanWarner2}
\leavevmode\vrule height 2pt depth -1.6pt width 23pt, {\em A direct approach to
  the determination of {G}aussian and scalar curvature functions}, Invent.
  Math., 28 ({\dummysort{c}}1975), pp.~227--230.

\bibitem{KazdanWarner3}
\leavevmode\vrule height 2pt depth -1.6pt width 23pt, {\em Prescribing
  curvatures}, in Differential geometry (Proc. Sympos. Pure Math., Vol. XXVII,
  Stanford Univ., Stanford, Calif., 1973), Part 2, Amer. Math. Soc.,
  Providence, R.I., {\dummysort{d}}1975, pp.~309--319.

\bibitem{Klingenberg}
{\sc W.~P.~A. Klingenberg}, {\em Riemannian geometry}, vol.~1 of de Gruyter
  Studies in Mathematics, Walter de Gruyter \& Co., Berlin, second~ed., 1995.

\bibitem{KobayashiNomizu1}
{\sc S.~Kobayashi and K.~Nomizu}, {\em Foundations of differential geometry,
  {V}ol.\ {I}. {\rm Reprint of the 1963 original}}, Wiley Classics Library,
  John Wiley \& Sons Inc., New York, 1996.

\bibitem{KolarMichorSlovak}
{\sc I.~Kol{\'a}{\v{r}}, P.~W. Michor, and J.~Slov{\'a}k}, {\em Natural
  operations in differential geometry}, Springer-Verlag, Berlin, 1993.

\bibitem{KrieglMichor}
{\sc A.~Kriegl and P.~W. Michor}, {\em The convenient setting of global
  analysis}, vol.~53 of Mathematical Surveys and Monographs, American
  Mathematical Society, Providence, RI, 1997.

\bibitem{LangRFA}
{\sc S.~Lang}, {\em Real and functional analysis}, vol.~142 of Graduate Texts
  in Mathematics, Springer-Verlag, New York, third~ed., 1993.

\bibitem{Lang}
\leavevmode\vrule height 2pt depth -1.6pt width 23pt, {\em Differential and
  {R}iemannian manifolds}, vol.~160 of Graduate Texts in Mathematics,
  Springer-Verlag, New York, third~ed., 1995.

\bibitem{LawsonMichelsohn}
{\sc H.~B. Lawson, Jr. and M.-L. Michelsohn}, {\em Spin geometry}, vol.~38 of
  Princeton Mathematical Series, Princeton University Press, Princeton, NJ,
  1989.

\bibitem{LeBrunMcKenzie}
{\sc C.~LeBrun and M.~Wang}, eds., {\em Surveys in differential geometry:
  essays on {E}instein manifolds}, Surveys in Differential Geometry, VI,
  International Press, Boston, MA, 1999.

\bibitem{LeeParker}
{\sc J.~M. Lee and T.~H. Parker}, {\em The {Y}amabe problem}, Bull. Amer. Math.
  Soc. (N.S.), 17 (1987), pp.~37--91.

\bibitem{Lohkamp1994a}
{\sc J.~Lohkamp}, {\em Metrics of negative {R}icci curvature}, Ann. of Math.
  (2), 140 (1994), pp.~655--683.

\bibitem{Lohkamp1995}
\leavevmode\vrule height 2pt depth -1.6pt width 23pt, {\em Curvature
  {$h$}-principles}, Ann. of Math. (2), 142 (1995), pp.~457--498.

\bibitem{Lutz1977}
{\sc R.~Lutz}, {\em Structures de contact sur les fibr\'es principaux en
  cercles de dimension trois}, Ann. Inst. Fourier (Grenoble), 27 (1977),
  pp.~ix, 1--15.

\bibitem{Marathe1972}
{\sc K.~B. Marathe}, {\em A condition for paracompactness of a manifold}, J.
  Differential Geometry, 7 (1972), pp.~571--573.

\bibitem{Martinet1971}
{\sc J.~Martinet}, {\em Formes de contact sur les vari\'et\'es de dimension
  {$3$}}, in Proceedings of Liverpool Singularities Symposium, II (1969/1970),
  Berlin, 1971, Springer, pp.~142--163. Lecture Notes in Math., Vol. 209.

\bibitem{McDuff1987}
{\sc D.~McDuff}, {\em Applications of convex integration to symplectic and
  contact geometry}, Ann. Inst. Fourier (Grenoble), 37 (1987), pp.~107--133.
\newblock Available in PDF format from \texttt{http://www.numdam.org/}
  \texttt{item?id=AIF\_1987\_\_37\_1\_107\_0}.

\bibitem{McDuffSalamonNew}
{\sc D.~McDuff and D.~Salamon}, {\em {$J$}-holomorphic curves and symplectic
  topology}, vol.~52 of Colloquium Publications, American Mathematical Society,
  Providence, RI, 2004.

\bibitem{Megginson}
{\sc R.~E. Megginson}, {\em An introduction to {B}anach space theory}, vol.~183
  of Graduate Texts in Mathematics, Springer-Verlag, New York, 1998.

\bibitem{MilnorHusemoller}
{\sc J.~Milnor and D.~Husemoller}, {\em Symmetric bilinear forms},
  Springer-Verlag, New York, 1973.
\newblock Ergebnisse der Mathematik und ihrer Grenzgebiete, Band 73.

\bibitem{MilnorStasheff}
{\sc J.~W. Milnor and J.~D. Stasheff}, {\em Characteristic classes}, Princeton
  University Press, Princeton, N. J., 1974.
\newblock Annals of Mathematics Studies, No. 76.

\bibitem{Morrey}
{\sc C.~B. Morrey, Jr.}, {\em Multiple integrals in the calculus of
  variations}, Die Grundlehren der mathematischen Wissenschaften, Band 130,
  Springer-Verlag New York, Inc., New York, 1966.

\bibitem{ONeill}
{\sc B.~O'Neill}, {\em Semi-{R}iemannian geometry. {W}ith applications to
  relativity}, vol.~103 of Pure and Applied Mathematics, Academic Press Inc.
  [Harcourt Brace Jovanovich Publishers], New York, 1983.

\bibitem{Palais1966}
{\sc R.~S. Palais}, {\em Homotopy theory of infinite dimensional manifolds},
  Topology, 5 (1966), pp.~1--16.

\bibitem{Palais}
\leavevmode\vrule height 2pt depth -1.6pt width 23pt, {\em Foundations of
  global non-linear analysis}, W. A. Benjamin, Inc., New York-Amsterdam, 1968.

\bibitem{Percell1981}
{\sc P.~Percell}, {\em Parallel vector fields on manifolds with boundary}, J.
  Differential Geom., 16 (1981), pp.~101--104.

\bibitem{Petersen}
{\sc P.~Petersen}, {\em Riemannian geometry}, vol.~171 of Graduate Texts in
  Mathematics, Springer-Verlag, New York, 1998.

\bibitem{RosenbergStolz1994}
{\sc J.~Rosenberg and S.~Stolz}, {\em Manifolds of positive scalar curvature},
  in Algebraic topology and its applications, vol.~27 of Math. Sci. Res. Inst.
  Publ., Springer, New York, 1994, pp.~241--267.

\bibitem{RosenbergStolz2001}
\leavevmode\vrule height 2pt depth -1.6pt width 23pt, {\em Metrics of positive
  scalar curvature and connections with surgery}, in Surveys on surgery theory,
  Vol. 2, vol.~149 of Ann. of Math. Stud., Princeton Univ. Press, Princeton,
  NJ, 2001, pp.~353--386.

\bibitem{Saunders}
{\sc D.~J. Saunders}, {\em The geometry of jet bundles}, vol.~142 of London
  Mathematical Society Lecture Note Series, Cambridge University Press,
  Cambridge, 1989.

\bibitem{SchejaStorch2}
{\sc G.~Scheja and U.~Storch}, {\em Lehrbuch der {A}lgebra. {T}eil 2},
  Mathematische Leitf\"aden. [Mathematical Textbooks], B. G. Teubner,
  Stuttgart, 1988.

\bibitem{Schoen1984}
{\sc R.~Schoen}, {\em Conformal deformation of a {R}iemannian metric to
  constant scalar curvature}, J. Differential Geom., 20 (1984), pp.~479--495.

\bibitem{Seeley1964}
{\sc R.~T. Seeley}, {\em Extension of {$C\sp{\infty }$} functions defined in a
  half space}, Proc. Amer. Math. Soc., 15 (1964), pp.~625--626.

\bibitem{Shiga1964}
{\sc K.~Shiga}, {\em Some aspects of real-analytic manifolds and differentiable
  manifolds}, J. Math. Soc. Japan, 16 (1964), pp.~128--142.
\newblock Erratum: J. Math. Soc. Japan {17} (1965), 216--217.

\bibitem{Spanier}
{\sc E.~H. Spanier}, {\em Algebraic topology}, Springer-Verlag, New York, 19??.
\newblock Corrected reprint of the 1966 original.

\bibitem{Spanier1993}
\leavevmode\vrule height 2pt depth -1.6pt width 23pt, {\em Singular homology
  and cohomology with local coefficients and duality for manifolds}, Pacific J.
  Math., 160 (1993), pp.~165--200.

\bibitem{Spring}
{\sc D.~Spring}, {\em Convex integration theory. Solutions to the $h$-principle
  in geometry and topology}, vol.~92 of Monographs in Mathematics, Birkh\"auser
  Verlag, Basel, 1998.

\bibitem{Steenrod}
{\sc N.~Steenrod}, {\em The topology of fibre bundles}, Princeton Landmarks in
  Mathematics, Princeton University Press, Princeton, NJ, 1999.
\newblock Reprint of the 1957 edition, Princeton Paperbacks.

\bibitem{Stolz1995}
{\sc S.~Stolz}, {\em Positive scalar curvature metrics---existence and
  classification questions}, in Proceedings of the International Congress of
  Mathematicians, Vol.\ 1, 2 (Z\"urich, 1994), Basel, 1995, Birkh\"auser,
  pp.~625--636.

\bibitem{Stolz2002}
\leavevmode\vrule height 2pt depth -1.6pt width 23pt, {\em Manifolds of
  positive scalar curvature}, in Topology of high-dimensional manifolds, No. 1,
  2 (Trieste, 2001), vol.~9 of ICTP Lect. Notes, Abdus Salam Int. Cent.
  Theoret. Phys., Trieste, 2002, pp.~661--709.

\bibitem{Taylor3}
{\sc M.~E. Taylor}, {\em Partial differential equations {III}. Nonlinear
  equations}, vol.~117 of Applied Mathematical Sciences, Springer-Verlag, New
  York, 1997.
\newblock Corrected reprint of the 1996 original.

\bibitem{Thomas1967a}
{\sc E.~Thomas}, {\em Fields of tangent {$2$}-planes on even-dimensional
  manifolds}, Ann. of Math. (2), 86 (1967), pp.~349--361.

\bibitem{Thomas1967b}
\leavevmode\vrule height 2pt depth -1.6pt width 23pt, {\em Fields of tangent
  {$k$}-planes on manifolds}, Invent. Math., 3 (1967), pp.~334--347.

\bibitem{Thurston1976}
{\sc W.~P. Thurston}, {\em Existence of codimension-one foliations}, Ann. of
  Math. (2), 104 (1976), pp.~249--268.

\bibitem{Thurston}
\leavevmode\vrule height 2pt depth -1.6pt width 23pt, {\em Three-dimensional
  geometry and topology. {V}ol. 1}, vol.~35 of Princeton Mathematical Series,
  Princeton University Press, Princeton, NJ, 1997.
\newblock Edited by Silvio Levy.

\bibitem{Tipler1977}
{\sc F.~J. Tipler}, {\em Singularities and causality violation}, Ann. Physics,
  108 (1977), pp.~1--36.

\bibitem{Tod1999}
{\sc K.~P. Tod}, {\em General relativity}, in Surveys in differential geometry:
  essays on Einstein manifolds, Surv. Differ. Geom., VI, Int. Press, Boston,
  MA, 1999, pp.~329--364.

\bibitem{Trudinger1968}
{\sc N.~S. Trudinger}, {\em Remarks concerning the conformal deformation of
  {R}iemannian structures on compact manifolds}, Ann. Scuola Norm. Sup. Pisa
  (3), 22 (1968), pp.~265--274.

\bibitem{Wald}
{\sc R.~M. Wald}, {\em General relativity}, University of Chicago Press,
  Chicago, IL, 1984.

\bibitem{Whitehead}
{\sc G.~W. Whitehead}, {\em Elements of homotopy theory}, vol.~61 of Graduate
  Texts in Mathematics, Springer-Verlag, New York, 1978.

\bibitem{Wolf}
{\sc J.~A. Wolf}, {\em Spaces of constant curvature}, Publish or Perish Inc.,
  Houston, TX, fifth~ed., 1984.

\bibitem{Yamabe1960}
{\sc H.~Yamabe}, {\em On a deformation of {R}iemannian structures on compact
  manifolds}, Osaka Math. J., 12 (1960), pp.~21--37.

\bibitem{Yodzis1}
{\sc P.~Yodzis}, {\em Lorentz cobordism}, Comm. Math. Phys., 26 (1972),
  pp.~39--52.

\bibitem{Yodzis2}
\leavevmode\vrule height 2pt depth -1.6pt width 23pt, {\em Lorentz cobordism.
  {II}}, General Relativity and Gravitation, 4 (1973), pp.~299--307.

\bibitem{ZeidlerI}
{\sc E.~Zeidler}, {\em Nonlinear functional analysis and its applications {I}.
  {F}ixed-point theorems}, Springer-Verlag, New York, 1986.

\end{thebibliography}
\addtocontents{toc}{\protect\contentsline {chapter}{Bibliography}{\arabic{page}}}

\newcommand{\dummysort}[1]{}


\end{document}